%% file: piN_Collection_arxiv_v3.tex
\UseRawInputEncoding 

\documentclass{article}

\usepackage{amssymb}

\usepackage{proof}
\usepackage{latexsym}

\newcommand{\mS}{\mathbb{S}}
\newcommand{\mI}{\mathbb{I}}
\newcommand{\mT}{\mathbb{T}}
\newcommand{\mU}{\mathbb{U}}
\newcommand{\mV}{\mathbb{V}}
\newcommand{\mW}{\mathbb{W}}

\newcommand{\supp}{{\rm supp}}

\usepackage{proof}
\usepackage{latexsym}
\usepackage{amsfonts}
\usepackage{amssymb}

\input{header}

\title{
An ordinal analysis of $\Pi_{N}$-Collection
}
\author{Toshiyasu Arai
\\
Graduate School of Mathematical Sciences,
University of Tokyo
\\
3-8-1 Komaba, Meguro-ku,
Tokyo 153-8914, JAPAN
\\
tosarai@ms.u-tokyo.ac.jp}
\date{}
\begin{document}
\maketitle
\begin{abstract}
In this paper we give an ordinal analysis of a set theory with $\Pi_{N}$-Collection.
\end{abstract}

\section{Introduction}
Throughout in this paper $N$ denotes a fixed positive integer.
In this paper we give an ordinal analysis of a Kripke-Platek 
set theory with the axiom of Infinity and one of $\Pi_{N}$-Collection, denoted
by ${\sf KP}\ome+\Pi_{N}\mbox{-Collection}$.
Our proof is an extension of \cite{singlewfprf, singlestable}.
Since \cite{singlestable} has not yet appeared, some proofs are duplicated for the readers' conveniences.

In \cite{singlestable} we analyzed proof-theoretically
a set theory 
${\sf KP}\ell^{r}+(M\prec_{\Sigma_{1}}V)$ extending  ${\sf KP}\ell^{r}$ 
with an axiom stating that
`there exists a transitive set $M$ such that $M\prec_{\Sigma_{1}}V$'.
An ordinal analysis of an extension
${\sf KP}i+(M\prec_{\Sigma_{1}}V)$ is given in M. Rathjen\cite{RathjenAFML2}.
Our proof is an extension of \cite{KPPiN, singlestable}.
In \cite{KPPiN},
a set theory ${\sf KP}\Pi_{N}$ of $\Pi_{N}$-reflection is analyzed, while
\cite{KPPiN} 
is an extension of M. Rathjen's analysis in \cite{Rathjen94} 
for 
$\Pi_{3}$-reflection.

$\Sig^{1}_{N+2}\mbox{-DC+BI}$ [$\Sig^{1}_{N+2}\mbox{-AC+BI}$]
denotes a second order arithmetic obtained from $\mbox{ACA}_{0}+\mbox{BI}$
 by adding the axiom of $\Sig^{1}_{N+2}$-Dependent Choice
 [$\Sig^{1}_{N+2}$-Axiom of Choice], resp.
 It is easy to see that $\Sig^{1}_{N+2}\mbox{-DC+BI}$ is interpreted canonically to the set theory
${\sf KP}\ome+\Pi_{N}\mbox{-Collection}+(V=L)$ with the axiom $V=L$ of constructibility.
It is well known that $\Sig^{1}_{N+2}\mbox{-DC}_{0}$ implies $\Sig^{1}_{N+2}\mbox{-AC}$, which yields $\Del^{1}_{N+2}\mbox{-CA}$,
a fortiori $\Sig^{1}_{N+1}\mbox{-CA}$, cf.\,Lemma VII.6.6 of \cite{Simpson}.
Moreover it is known that $\Sig^{1}_{N+2}\mbox{-DC+BI}$ is 
$\Pi^{1}_{4}$-conservative over $\Sig^{1}_{N+2}\mbox{-AC+BI}$
[over $\Del^{1}_{N+2}\mbox{-CA+BI}$], resp., 
cf.\,Exercise VII.5.13 and
Theorem VII.6.16 of \cite{Simpson}.

Let $n$ be a positive integer.
We say that an ordinal $\alp$ is $n$-\textit{stable} if $L_{\alp}\prec_{\Sig_{n}}L$
for the constructible universe $L=\bigcup_{\alp}L_{\alp}$.
In general, a transitive and non-empty set $M$ is $n$-\textit{stable} if $M\prec_{\Sig_{n}}V$
for the universe $V$.
We see that
$(V,\in)\models {\sf KP}\ome+\Pi_{N}\mbox{-Collection}$ if
$V$ enjoys the $\Del_{0}(\{st_{i}\}_{0<i\leq N})$-collection, where
$st_{i}$ denotes the predicate
for the class $\{M\in V: M\prec_{\Sig_{i}}V\}$ of $i$-stable sets
in $V$.

We introduce an extension $S_{\mI_{N}}$ of ${\sf KP}\ome+\Pi_{N}\mbox{-Collection}$ in the language $\{\in\}\cup\{st_{i}\}_{0<i\leq N}$,
which codifies $\Sig(\{st_{i}\}_{0<i\leq N})$-reflection.
We aim to give an ordinal analysis of the theory $S_{\mI_{N}}$.

In the following theorems,
$\Ome$ denotes the least recursively regular ordinal $\ome_{1}^{CK}$, and $\psi_{\Ome}$ 
a collapsing function such that $\psi_{\Ome}(\alp)<\Ome$.
$\mI_{N}$ is an ordinal term denoting an ordinal such that
$L_{\mI_{N}}\models {\sf KP}\ome+\Pi_{N}\mbox{-Collection}+(V=L)$.

First we show the following Theorem \ref{thm:2}.

\begin{theorem}\label{thm:2}
Suppose $S_{\mI_{N}}\vdash\tht^{L_{\Ome}}$
for a $\Sig_{1}$-sentence $\tht$ in the language $\{\in\}$ of set theory.
Then 
$L_{\psi_{\Ome}(\veps_{\mI_{N}+1})}\models\tht$ holds.
\end{theorem}

It is not hard to see that the ordinal $\psi_{\Ome}(\veps_{\mI_{N}+1})$ is computable.
Let $<$ denote a computable well-ordering of type $\psi_{\Ome}(\veps_{\mI_{N}+1})$
on the set of natural numbers.
Conversely we show that $\Sig^{1}_{N+2}\mbox{{\rm -DC+BI}}$ proves that
each initial segment of $\psi_{\Ome}(\veps_{\mI_{N}+1})$ is well-founded.

\begin{theorem}\label{th:wf}
$\Sig^{1}_{N+2}\mbox{{\rm -DC+BI}}\vdash Wo[\alp]$
for {\rm each} $\alp<\psi_{\Ome}(\veps_{\mI_{N}+1})$.
\end{theorem}

For $T\supset{\rm ACA}_{0}$, $|T|$ denotes
the proof-theoretic ordinal of $T$, i.e., the supremum of order types
of computable well-orderings $\prec$ on the set of natural numbers for which
$T$ proves the fact that $\prec$ is a well-ordering.
Also let
$|{\sf KP}\ome+\Pi_{N}\mbox{-Collection}|_{\Sig_{1}^{\Ome}}$ denote the 
$\Sig_{1}^{\Ome}$-ordinal of ${\sf KP}\ome+\Pi_{N}\mbox{-Collection}$, i.e., the ordinal
$
\min\{\alp\leq\ome_{1}^{CK}: \fal \tht\in\Sig_{1}
\left(
{\sf KP}\ome+\Pi_{N}\mbox{-Collection}\vdash\tht^{L_{\Ome}} \Rarw L_{\alp}\models\tht
\right)\}
$.
For more on ordinal analysis see \cite{OA}.
We conclude the following Theorem \ref{th:main}, where
$\psi_{\Ome}(\veps_{\mI_{N}+1})$ denotes the order type of 
the initial segment $OT(\mI_{N})\cap\Ome$ of a notation system $OT(\mI_{N})$ of
ordinals.

\begin{theorem}\label{th:main}
$
|\Del^{1}_{N+2}\mbox{{\rm -CA+BI}}|
=
|\Sig^{1}_{N+2}\mbox{{\rm -AC+BI}}|=
|\Sig^{1}_{N+2}\mbox{{\rm -DC+BI}}|=
\\
|{\sf KP}\ome+\Pi_{N}\mbox{{\rm -Collection}}|_{\Sig_{1}^{\Ome}}
=
\psi_{\Ome}(\veps_{\mI_{N}+1})
$.
\end{theorem}

Let $\mathbf{Z}_{2}=\Sig^{1}_{\infty}\mbox{{\rm -DC}}$ be  the full second order arithmetic with
the Dependent Choice schema, and ${\sf ZFC}-{\sf Power}$ denote the set theory
{\sf ZFC} minus the power set axiom.
$\mathbf{Z}_{2}$ proves the ($\Pi^{1}_{1}$-)soundness of
$\Sigma^{1}_{N+2}\mbox{{\rm -DC}}+\mbox{BI}$,
and hence $\mathbf{Z}_{2}$ proves that
$(OT(\mI_{N}),<)$ is a well ordering for \textit{each} $N$.
$\mathbf{Z}_{2}$ is canonically interpreted in $({\sf ZFC}-{\sf Power})+(V=L)$, which is
$\Pi^{1}_{1}$-conservative over ${\sf ZFC}-{\sf Power}$.

Assume ${\sf ZFC}-{\sf Power}\vdash \tht$ for a sentence $\tht$.
$S_{\mI_{N}}$ subsumes $\Pi_{N}\mbox{-Collection}$ and 
$\Sig_{N}\mbox{-Separation}$, and there is an $N$ such that 
$S_{\mI_{N}}\vdash\tht$.
Therefore we conclude the following.

\begin{theorem}\label{th:mainZ}
\[
\psi_{\Omega}(\mathbb{I}_{\ome}):=\sup\{\psi_{\Omega}(\mathbb{I}_{N}) :
0<N<\ome\}
=|\mathbf{Z}_{2}|
= |{\sf ZFC}-{\sf Power}|_{\Sigma_{1}^{\Omega}}.\]
\end{theorem}

Let us mention the contents of this paper.
In the next section \ref{sec:9} a second order arithmetic $\Sig^{1}_{N+2}\mbox{-DC+BI}$ 
is interpreted to a set theory
${\sf KP}\ome+\Pi_{N}\mbox{-Collection}+(V=L)$, and
${\sf KP}\ome+\Pi_{N}\mbox{-Collection}$ is shown to be 
a subtheory of a set theory $S_{\mI_{N}}$.
In section \ref{sect:ordinalnotation} ordinals for our analysis of $\Pi_{N}$-Collection are introduced, 
and a computable notation system $OT(\mI_{N})$ is extracted.

Theorem \ref{thm:2} is proved in sections \ref{sec:operatorcont} and \ref{sec:proofonestep}.
In section \ref{sec:operatorcont} operator controlled derivations are introduced.
In section \ref{sec:proofonestep}, stable ordinals are removed
from derivations.
Although our proof of Theorem \ref{thm:2} is based on operator controlled derivations introduced by W. Buchholz\cite{Buchholz}, it is hard for us to give its sketch here.
See subsection \ref{subsec:preview} for an outline of the proof.

Theorem \ref{th:wf} is proved 
in sections \ref{sect:prop} and \ref{sec:distinguished}.
For $0\leq i\leq N$, we introduce \textit{$i$-maximal distinguished sets}, which are
$\Sig^{1}_{2+i}$-definable.
A $0$-maximal distinguished set is $\Sig^{1}_{2}$-definable as in \cite{singlewfprf}.
$\Sig^{1}_{N+2}$-(Dependent) Choice is needed to handle limits of $N$-stable ordinals.
Our proof of Theorem \ref{th:wf} is based on maximal distinguished class introduced
again by Buchholz\cite{Buchholz75}.
A sketch of the well-foundedness proof is outlined in subsection \ref{subsec:preview_wf}.

In the final section \ref{sect:consis} let us conclude some standard outcomes of
an ordinal analysis of the theory $\mathbf{Z}_{2}$.

IH denotes the Induction Hypothesis, MIH the Main IH, SIH the Subsidiary IH,
and SSIH the Sub-Subsidiary IH.

\section{$\Pi_{N}$-Collection}\label{sec:9}
In this section a second order arithmetic $\Sig^{1}_{N+2}\mbox{-DC+BI}$ is interpreted canonically to a set theory
${\sf KP}\ome+\Pi_{N}\mbox{-Collection}+(V=L)$, and
${\sf KP}\ome+\Pi_{N}\mbox{-Collection}$ is shown to be 
a subtheory of a set theory $S_{\mI_{N}}$.

For subsystems of second order arithmetic, we follow largely Simpson's monograph\cite{Simpson}.
The schema Bar Induction, {\rm BI} is denoted by {\rm TI} in \cite{Simpson}.
 {\rm BI} allows the transfinite induction schema for well-founded relations.
 
$\Sig^{1}_{N+2}\mbox{-AC+BI}$ denotes a second order arithmetic obtained from 
$\Pi^{1}_{1}\mbox{-CA}_{0}+\mbox{BI}$
 by adding the axiom $\Sig^{1}_{N+2}\mbox{-AC}$, 
 $\fal n\exi XF(n,X)\to \exi Y \fal n F(n,Y_{n})$
for each $\Pi^{1}_{N+1}$-formula $F(n,X)$, where $m\in Y_{n}\Lrarw (n,m)\in Y$
 for a bijective pairing function $(\cdot,\cdot)$.
$\Sig^{1}_{N+2}\mbox{-DC+BI}$ denotes a second order arithmetic obtained from 
$\Pi^{1}_{1}\mbox{-CA}_{0}+\mbox{BI}$
 by adding the axiom $\Sig^{1}_{N+2}\mbox{-DC}$ for each $\Pi^{1}_{N+1}$-formula $F(n,X,Y)$,
 $\fal n\fal X\exi YF(n,X,Y)\to \fal X_{0}\exi Y \fal n [Y_{0}=X_{0}\land F(n,Y_{n},Y_{n+1})]$.
 It is easy to see that the formulas $F$ can be $\Sig^{1}_{N+2}$ in the axioms.

The axioms of the set theory ${\sf KP}\ome+\Pi_{N}\mbox{-Collection}$ consists of those of ${\sf KP}\ome$
 (Kripke-Platek set theory with the Axiom of Infinity, cf.\cite{Ba, J3}) plus $\Pi_{N}\mbox{-Collection}$:
for each $\Pi_{N}$-formula $A(x,y)$ in the language of set theory, 
$\fal x\in a\exi yA(x,y)\to \exi b\fal x\in a\exi y\in b A(x,y)$.

$\Sig_{N}\mbox{-Separation}$ denotes the axiom $\exi y\fal x(x\in y\lrarw  x\in a \land \vphi(x))$
for each $\Sig_{N}$-formula $\vphi(x)$.
 $\Del_{N+1}\mbox{-Separation}$ denotes the axiom
 $\fal x\in a(\vphi(x)\lrarw \lnot \psi(x))\to \exi y\fal x(x\in y\lrarw  x\in a \land \vphi(x))$
for each $\Sig_{N+1}$-formulas $\vphi(x)$ and $\psi(x)$.

$\Sig_{N+1}\mbox{-Replacement}$ denotes the axiom stating that if
$\fal x\in a\exi !y\vphi(x,y)$, then there exists a function $f$ with its domain $\mathrm{dom}(f)=a$ such that
$\fal x\in a\,\vphi(x,f(x))$
for each $\Sig_{N+1}$-formula $\vphi(x,y)$.

\blem\label{lem:9Sep}
${\sf KP}\ome+\Pi_{N}\mbox{{\rm -Collection}}$ proves each of
$\Sig_{N}\mbox{{\rm -Separation}}$,
$\Del_{N+1}-\\
\mbox{{\rm Separation}}$ and
$\Sig_{N+1}\mbox{{\rm -Replacement}}$.
\elem
\bprf
We show that $\{x\in a:\vphi(x)\}$ exists as a set for each $\Sig_{i}$-formula 
$\vphi$ by (meta)induction on $i\leq N$.
The case $i=0$ follows from $\Del_{0}$-Separation.
Let $\vphi\equiv \exi y\, \tht(x,y)$ 
with a $\Pi_{i-1}$-matrix $\tht$. 
We have by logic $\fal x\in a\exi y(\exi z\tht(x,z)\rarw \tht(x,y))$.
By $\Pi_{i}$-Collection pick a set $b$ so that
$\fal x\in a\exi y\in b(\vphi(x)\rarw \tht(x,y))$.
In other words,
$\{x\in a:\vphi(x)\}=\{x\in a: \exi y\in b\, \tht(x,y)\}$.
If $i=1$, then $\exi c[c=\{x\in a: \exi y\in b\, \tht(x,y)\}]$ by $\Del_{0}$-Separation.
Let $2\leq i\leq N$.
By $\Pi_{i-2}$-Collection we obtain a $\Pi_{i-1}$-formula $\sig$ such that
$\exi y\in b\, \tht(x,y)\lrarw \sig(x)$.
By IH we obtain $\exi c[c=\{x\in a :\sig(x)\}]$.

$\Del_{N+1}\mbox{-Separation}$ follows from $\Sig_{N}\mbox{-Separation}$ as in \cite{Ba}, p.17, Theorem 4.5($\Del$ Separation),
and
$\Sig_{N+1}\mbox{-Replacement}$ follows from $\Del_{N+1}\mbox{-Separation}$ as in \cite{Ba}, p.17, Theorem 4.6($\Sig$ Replacement).
\eprf
\\

For a formula $A$ in the language of second order arithmetic let $A^{set}$ denote the formula obtained from $A$ by interpreting the first order variable $x$ as $x\in \ome$ and the second order variable $X$ as $X\subset\ome$.

The following is the Quantifier Theorem in p.125 of \cite{J3}, in which
${\sf KP}l^{r}$ is defined as a set theory for limits of admissible sets with restricted induction.
${\sf KP}l^{r}$ is a subtheory of ${\sf KP}\ome+\Pi_{N}\mbox{{\rm -Collection}}$.
$Ad(x)$ designates that $x$ is an admissible set.

\blem\label{lem:9Quant}
For each $\Sig^{1}_{N+1}$-formula $F(n,a,Y)$,
there exists a $\Sig_{N}$-formula $A_{\Sig}(d,n,a,Y)$ in the language of set theory so that for
\\
$F_{\Sig}(n,a,Y) :\Lrarw \exi d[Ad(d)\land Y\in d\land A_{\Sig}(d,n,a,Y)]$,
\[
{\sf KP}l^{r}\vdash n,a\in\ome\land Y\subset\ome\to \{F^{set}(n,a,Y)\lrarw F_{\Sig}(n,a,Y)\} 
.\]
\elem

For an ordinal $\alpha$, $L_{\alp}$ denotes the initial segment of G\"odel's constructible universe 
$L=\bigcup_{\alp}L_{\alp}$. 
$x\in L$ is a $\Sig_{1}$-formula.
$<_{L}$ denotes a canonical $\Del_{1}$ well ordering of $L$ such that
if $y<_{L}x\in L_{\alp}$, then $y\in L_{\alp}$, cf.\,p.162 of \cite{Ba}.
$V=L$ denotes the axiom of Constructibility.

\blem\label{lem:9scdset}
For each sentence $A$ in the language of second order arithmetic,
\[
\Sig^{1}_{N+2}{\rm -DC+BI}\vdash A
\Rarw 
{\sf KP}\ome+\Pi_{N}\mbox{{\rm -Collection}}+(V=L)\vdash A^{set}
.\]
\elem
\bprf
By the Quantifier Theorem \ref{lem:9Quant} $F^{set}(n,X,Y)$ is equivalent to a $\Pi_{N}$-formula $\vphi(n,X,Y)$ 
for a $\Pi^{1}_{N+1}$-formula $F(n,X,Y)$, $n\in\ome$ and $X\subset\ome$. 
It suffices to show for a $\Pi_{N}$-formula $\vphi(n,X,Y)$ that assuming
$\fal n\in\ome\fal X\subset\ome\exi Y\subset\ome \vphi(n,X,Y)$ and $X_{0}\subset\ome$,
there exists a function $f$ with its domain $\mathrm{dom}(f)=\ome$ such that
$\fal n\in\ome [f(0)=X_{0} \land \vphi(n,f(n),f(n+1))]$.
By induction on $k\in\ome$ using $V=L$ we see that there exists a unique family $(Y_{n})_{n<k}$ of subsets of $\ome$
such that
$\fal n<k[\vphi(n,Y_{n},Y_{n+1})\land \fal Z<_{L}Y_{n+1}\lnot\vphi(n,Y_{n},Z)]$, where
$\fal Z<_{L}Y\lnot\vphi(n,Y,Z)$ is equivalent to a $\Sig_{N}$-formula under
$\Pi_{N-1}$-Collection.
By $\Sig_{N+1}$-Replacement pick a function $g$ with $\mathrm{dom}(g)=\ome$ and $rng(g)\subset{}^{<\ome}{\cal P}(\ome)$ 
so that for any  $k\in\ome$ $g(k)$ is the unique sequence $(Y_{n})_{n<k}\in {}^{k}{\cal P}(\ome)$ with $Y_{0}=X_{0}$. 
Then the function $f(n)=\left(g(n+1)\right)(n)$ is a desired one.
\eprf
\\

\noindent
It is easy to see that
${\sf KP}\ome+\Pi_{N}\mbox{{\rm -Collection}}+(V=L)\vdash A \Rarw
{\sf KP}\ome+\Pi_{N}\mbox{{\rm -Collection}}\vdash A^{L}$ for any $A$, and
each $\Pi^{1}_{1}$-sentence $B$ on $\ome$ is absolute for $L$,
${\sf KP}\ome+\Pi_{N}\mbox{{\rm -Collection}}\vdash B\lrarw B^{L}$.

Next we show that ${\sf KP}\ome+\Pi_{N}\mbox{-Collection}$ is contained in a set theory 
$S_{\mI_{N}}$.
The language of the theory $S_{\mI_{N}}$ is $\{\in,\mathsf{M}_{0}\}\cup\{st_{i}\}_{0<i\leq N}$ with unary predicate constants $st_{i}$
and an individual constant $\mathsf{M}_{0}$.
$st_{i}(a)$ is intended to denote the fact that $a$ is an $i$-stable set
and $\mathsf{M}_{0}$ is intended to denote the least admissible set $L_{\ome_{1}^{CK}}$ above $L_{\ome}$.
The axioms of $S_{\mI_{N}}$ are obtained from those\footnote{In the axiom schemata $\Del_{0}$-Separation and $\Del_{0}$-Collection,
 $\Del_{0}$-formulas remain to mean a $\Del_{0}$-formula in which $st_{i}$ does not occur, 
 while the axiom of foundation may be applied to a formula in which $st_{i}$ may occur.} of ${\sf KP}\ome$ by adding the following axioms.
By a $\Del_{0}(\{st_{i}\}_{0<i<k})$-formula we mean a bounded formula in the language
$\mathcal{L}_{k}=\{\in,\mathsf{M}_{0}\}\cup\{st_{i}\}_{i<k}$.

\benu
\item
The axioms for the admissible set $\mathsf{M}_{0}$:
$\mathsf{M}_{0}\neq\emptyset$, $\fal x\in \mathsf{M}_{0}\fal y\in x(y\in \mathsf{M}_{0})$, and the axioms stating that
$(\mathsf{M}_{0},\in)\models{\sf KP}\ome$.

\item
$\Del_{0}(\{st_{i}\}_{0<i\leq N})\mbox{-collection}$:
\[
\fal x\in a\exi y\, \tht(x,y) \to \exi b
\fal x\in a \exi y\in b\, \tht(x,y)
\]
for each $\Del_{0}(\{st_{i}\}_{0<i\leq N})$-formula $\tht$ in which the predicates $st_{i}$ may occur.
Note that $\Sig_{1}(\{st_{i}\}_{0<i\leq N})\mbox{-collection}$ follows from this.

\item
\beqn\label{eq:LimN}
\fal a\exi b[a\in b \land st_{N}(b)]
\eeqn

\item
For each $i+1\leq N$:
\beqn\label{eq:stbl0}
st_{i+1}(a) \to M_{0}\in a \land\fal y\in a\fal z\in y(z\in a) \land lst_{i}(a)
\eeqn
where $lst_{i}(a) :\Lrarw \fal b\in a \exi c\in a \left(b\in c \land st_{i}(c)\right)$ and
$st_{0}(c):\Lrarw (0=0)$.

\item 
For $0<i\leq N$:
\beqn\label{eq:sucstable}
st_{i}(a)\land \vphi(u) \land u\in a \to \vphi^{a}(u)
\eeqn
for each $\Sig_{1}(\{st_{j}\}_{j<i})$-formula $\vphi\equiv(\exi x\,\tht)$ in the language 
$\mathcal{L}_{i}=\{\in,M_{0}\}\cup\{st_{j}\}_{j<i}$, where
$\vphi^{a}\equiv(\exi x\in a\,\tht)$.

\eenu

Note that if $lst_{i+1}(a)$ for a transitive set $a$, then $lst_{i}(a)$ holds.

\blem\label{lem:9setIM}
$S_{\mI_{N}}\vdash st_{i}(M) \land u\in M \to
[\vphi^{M}(u)\lrarw \vphi(u)]
$
for set-theoretic $\Sig_{i}$-formulas $\vphi$.
\elem
\bprf
Argue in $S_{\mI_{N}}$.
The case $i=1$ follows from the axiom (\ref{eq:sucstable}).
We show
\beqn\label{eq:9setIM2}
st_{k}(a)\land u\in a  \to \left[
\tht^{a}(u) \lrarw \exi b\in a\{st_{i}(b) \land u\in b \land \tht^{b}(u)\}
\right]
\eeqn
for $0\leq i<k\leq N+1$ and $\Pi_{1}(\{st_{j}\}_{j<i-1})$-formula $\tht(u)$, where
$a=V$,
$st_{N+1}(V):\Lrarw(0=0)$ and $\tht^{V}(u):\Lrarw\tht$ when $k=N+1$.

Assume $st_{k}(a)$ and $\tht^{a}(u)$ with $u\in a$.
By the axioms (\ref{eq:LimN}) and (\ref{eq:stbl0}) 
there exists a set $b\in a$ such that $st_{i}(b)$ and $u\in b$.
$\tht^{b}(u)$ follows logically.
Conversely assume $\tht^{b}(u)$ for 
$b\in a$ such that $st_{i}(b)$ and $u\in b$.
(\ref{eq:sucstable}) yields $\tht(u)$, a fortiori $\tht^{a}(u)$.
Thus (\ref{eq:9setIM2}) is shown.

Let $\vphi(u)\in\Sig_{1+n}(\{st_{j}\}_{j<i})$ and $st_{i+n}(a)$ with $u\in a$.
From (\ref{eq:9setIM2}) we see by (meta-)induction on $n$ that
there exists a $\Sig_{1}(\{st_{j}\}_{j<i+n})$-formula $\tht$ such that
$\vphi^{a}\lrarw\tht^{a}$ and $\vphi \lrarw\tht$.

Now we show $\vphi^{M}(u)\lrarw \vphi(u)$, where $0\leq n<N$,
$st_{1+n}(M)$, $\vphi\in\Sig_{1+n}$ and $u\in M$.
Suppose $\vphi^{M}(u)$.
Pick a $\Sig_{1}(\{st_{j}\}_{j<n})$-formula $\tht$ such that
$\vphi^{M}(u)\lrarw\tht^{M}(u)$ and $\vphi(u)\lrarw\tht(u)$.
$\tht(u)$ follows logically, and $\vphi(u)$ follows.
Conversely assume $\vphi(u)$. Then we obtain $\tht(u)$, and 
(\ref{eq:sucstable}) yields $\tht^{M}(u)$, and hence $\vphi^{M}(u)$.
\eprf

\blem\label{lem:9setI}
$S_{\mI_{N}}$ is an extension of ${\sf KP}\ome+\Pi_{N}\mbox{{\rm -Collection}}$.
Namely $S_{\mI_{N}}$ proves $\Pi_{N}$-Collection.
\elem
\bprf
Argue in $S_{\mI_{N}}$.
Let $A(x,y)$ be a $\Pi_{N}$-formula in the language of set theory.
We obtain by the axiom (\ref{eq:LimN})
and Lemma \ref{lem:9setIM}
\beqn\label{eq:9setI}
A(x,y)\lrarw\exi b(st_{N}(b) \land x,y\in b \land A^{b}(x,y))
\eeqn
Assume $\fal x\in a\exi yA(x,y)$. 
Then we obtain
$\fal x\in a\exi y\exi b(st_{N}(b) \land x,y\in b \land A^{b}(x,y))$ by (\ref{eq:9setI}).
 Since $st_{N}(b) \land x,y\in b \land A^{b}(x,y)$ is a 
 $\Sig_{1}(\{st_{i}\}_{0<i\leq N})$-formula,
pick a set $c$ such that 
$\fal x\in a\exi y\in c\exi b\in c(st_{N}(b) \land x,y\in b \land A^{b}(x,y))$
by $\Del_{0}(\{st_{i}\}_{0<i\leq N})$-Collection.
Again by (\ref{eq:9setI}) we obtain $\fal x\in a\exi y\in c A(x,y)$.
\eprf

\section{Ordinals for $\Pi_{N}$-Collection}\label{sect:ordinalnotation}

In this section up to subsection \ref{subsec:largecardinal}
we work in a set theory ${\sf ZFC}(\{St_{i}\}_{i})$,
where each $St_{i}\,(0<i\leq N)$ is a unary predicate symbol. 
Let $St_{0}$ denote the set of uncountable cardinals below $\mI_{N}$.
$\Ome$ and $\mI_{N}$ are strongly critical numbers with $\Ome<\mI_{N}$, i.e.,
non-zero ordinals closed under the binary Veblen function 
$\varphi\alpha\beta=\varphi_{\alpha}(\beta)$.
We assume that 
$St_{i+1}\subset St_{i}$ for $i<N$,
each $St_{i}$ is an unbounded class of ordinals below $\mI_{N}$ such that
the least element of $St_{i}$ is larger than $\Ome$, $\Ome<\min(\bigcup_{0<i\leq N} St_{i})$.
The predicate $St_{i}$ is identified with the class $\{\alp\in ON: \alp\in St_{i}\}$.
$\alp^{\dagger i}$ denotes the least ordinal$>\alp$ in the class $St_{i}$ when $\alp<\mI_{N}$.
$\alp^{\dagger i}:=\mI_{N}$ if $\alp\geq\mI_{N}$.
Put
$\alp^{\dagger}:=\alp^{\dagger 1}$.
Let $SSt_{i}:=\{\alp^{\dagger i}:\alp\in ON\}$ and $LSt_{i}=St_{i}\setminus SSt_{i}$.

$\Gamma_{a}$ denotes the $a$-th strongly critical number.
For ordinals $\alp$, $\veps(\alp)$ denotes the least epsilon number above $\alp$, and
$\Gam(\alp)$ the least strongly critical number above $\alp$.
For ordinals $\alp,\bet$, and $\gam$,
$\gam=\alp-\bet$ designates that $\alp=\bet+\gam$.
$\alpha\dot{+}\beta$ denotes the sum $\alpha+\beta$
when $\alpha+\beta$ equals to the commutative (natural) sum $\alpha\#\beta$, i.e., when
either $\alpha=0$ or $\alpha=\alpha_{0}+\omega^{\alpha_{1}}$ with
$\omega^{\alpha_{1}+1}>\beta$.

$u,v,w,x,y,z,\ldots$ range over sets in the universe,
$a,b,c,\alp,\bet,\gam,\del,\ldots$ range over ordinals$<\veps(\mI_{N})$, and
$\xi,\zeta,\eta,\ldots$ range over ordinals less than $\Gam(\mI_{N})$,
and ordinals$\leq\mI_{N}$ are denoted by $\pi,\kap,\rho,\sig,\tau,\lam,\ldots$.

Let $\mS\in St_{i}$ with $i>0$.
A `Mahlo degree' $m(\pi)$ of ordinals $\pi<\mS$
with higher reflections is defined to be a finite function 
$f:\mI_{N}\to\Gamma(\mI_{N})$.
Let $\Lam\leq\mI_{N}$ be a strongly critical number.
To denote ordinals$<\Gam(\mI_{N})$,
it is convenient for us to introduce an ordinal function 
$\tilde{\tht}_{b}(\xi;\Lam)<\Gam(\mI_{N})$,
$\xi<\Gam(\mI_{N})$ and $b<\Lam$ as in \cite{singlewfprf, singlestable}, 
which is a $b$-th iterate of the exponential
$\tilde{\tht}_{1}(\xi;\Lam)=\Lam^{\xi}$
with the base
$\Lam$.

\begin{definition}\label{df:Lam}
{\rm
Let $\Lam\leq\mI_{N}$ be a strongly critical number.
$\varphi_{b}(\xi)$ denotes the binary Veblen function on 
$(\mI_{N})^{\dagger 0}=\omega_{\mathbb{I}_{N}+1}$ with $\varphi_{0}(\xi)=\omega^{\xi}$, and
$\tilde{\varphi}_{b}(\xi;\Lam):=\varphi_{b}(\Lam\cdot \xi)$.

Let $b,\xi<(\mI_{N})^{\dagger 0}$.
$\theta_{b}(\xi)$ [$\tilde{\theta}_{b}(\xi;\Lam)$] denotes
a $b$-th iterate of $\varphi_{0}(\xi)=\omega^{\xi}$ [of $\tilde{\varphi}_{0}(\xi;\Lam)=\Lam^{\xi}$], resp.
Specifically ordinals
$\theta_{b}(\xi), \tilde{\theta}_{b}(\xi;\Lam)<(\mI_{N})^{\dagger 0}$ are defined by recursion on $b$ as follows.
$\theta_{0}(\xi)=\tilde{\theta}_{0}(\xi;\Lam)=\xi$, 
$\theta_{\omega^{b}}(\xi)=\varphi_{b}(\xi)$,
$\tilde{\theta}_{\omega^{b}}(\xi;\Lam)=\tilde{\varphi}_{b}(\xi;\Lam)$, and
$\theta_{c\dot{+}\omega^{b}}(\xi)=\theta_{c}(\theta_{\omega^{b}}(\xi))$,
$\tilde{\theta}_{c\dot{+}\omega^{b}}(\xi;\Lam)=
\tilde{\theta}_{c}(\tilde{\theta}_{\omega^{b}}(\xi;\Lam);\Lam)$.

}
\end{definition}

A finite set $SC(a)$ of strongly critical numbers is defined recursively as follows.
$SC(0)=\emptyset$, $SC(a)=\bigcup_{i\leq m}SC(a_{i})$ for $a=\ome^{a_{m}}\dot{+}\cdots\dot{+}\ome^{a_{0}}$,
and
$SC(a)=SC(b)\cup SC(c)$ for $a=\vphi_{b}(c)$ if $a$ is not strongly critical.
$SC(a)=\{a\}$ if $a$ is strongly critical.

Let $\Lam\leq\mI_{N}$ be a strongly critical number.
Let us define a normal form of non-zero ordinals $\xi<\Gamma(\Lam)$ and a finite set $SC_{\Lambda}(\xi)$
of strongly critical numbers$<\Lam$.
Let $\xi=\Lam^{\zeta}$.
If $\zeta<\Lam^{\zeta}$, then $\tilde{\theta}_{1}(\zeta;\Lam)$ is the normal form of $\xi$,
denoted by
$\xi=_{NF}\tilde{\theta}_{1}(\zeta;\Lam)$ and $SC_{\Lambda}(\xi)=SC_{\Lambda}(\zeta)\cup\{1\}$.
Assume $\zeta=\Lam^{\zeta}$, and let $b>0$ be the maximal ordinal such that
there exists an ordinal $\eta$ with
$\zeta=\tilde{\varphi}_{b}(\eta;\Lam)>\eta$.
Then $\xi=\tilde{\varphi}_{b}(\eta;\Lam)=_{NF}\tilde{\theta}_{\omega^{b}}(\eta;\Lam)$ and
$SC_{\Lambda}(\xi)=SC_{\Lambda}(b)\cup SC_{\Lambda}(\eta)$.

Let $\xi=\Lam^{\zeta_{m}}a_{m}+\cdots+\Lam^{\zeta_{0}}a_{0}$, where
$\zeta_{m}>\cdots>\zeta_{0}$ and $0<a_{0},\ldots,a_{m}<\Lam$.
Let $\Lam^{\zeta_{i}}=_{NF}\tilde{\theta}_{b_{i}}(\eta_{i};\Lam)$ with
$b_{i}=\omega^{c_{i}}$ for each $i$.
Then
$\xi=_{NF}\tilde{\theta}_{b_{m}}(\eta_{m};\Lam)\cdot a_{m}+\cdots+\tilde{\theta}_{b_{0}}(\eta_{0};\Lam)\cdot a_{0}$
and 
$SC_{\Lambda}(\xi):=\{a_{i}:i\leq m\}\cup\bigcup_{i\leq m}(SC_{\Lambda}(b_{i})\cup SC_{\Lambda}(\eta_{i}))$, where
$SC_{\Lambda}(0)=\emptyset$.

\begin{definition}\label{df:Lam2}
{\rm
Let $\xi<\Gamma(\Lam)$ be a non-zero ordinal with its normal form:
\begin{equation}\label{eq:CantornfLam}
\xi=\sum_{i\leq m}\tilde{\theta}_{b_{i}}(\xi_{i};\Lam)\cdot a_{i}=_{NF}
\tilde{\theta}_{b_{m}}(\xi_{m};\Lam)\cdot a_{m}+\cdots+\tilde{\theta}_{b_{0}}(\xi_{0};\Lam)\cdot a_{0}
\end{equation}
where
$\tilde{\theta}_{b_{i}}(\xi_{i};\Lam)>\xi_{i}$,
$\tilde{\theta}_{b_{m}}(\xi_{m};\Lam)>\cdots>\tilde{\theta}_{b_{0}}(\xi_{0};\Lam)$, 
$b_{i}=\omega^{c_{i}}<\Lam$, and
$0<a_{0},\ldots,a_{m}<\Lam$.
$\tilde{\theta}_{b_{0}}(\xi_{0};\Lam)$ is said to be the \textit{tail} of $\xi$, denoted 
$\tilde{\theta}_{b_{0}}(\xi_{0};\Lam)=tl(\xi)$, and
$\tilde{\theta}_{b_{m}}(\xi_{m};\Lam)$ the \textit{head} of $\xi$, denoted 
$\tilde{\theta}_{b_{m}}(\xi_{m};\Lam)=hd(\xi)$.

\begin{enumerate}
\item\label{df:Exp2.3}
 $\zeta$ is a \textit{segment} of $\xi$
 iff there exists an $n\, (0\leq n\leq m+1)$
 such that
 $\zeta=_{NF}\sum_{i\geq n}\tilde{\theta}_{b_{i}}(\xi_{i};\Lam)\cdot a_{i}=
 \tilde{\theta}_{b_{m}}(\xi_{m};\Lam)\cdot a_{m}+\cdots+\tilde{\theta}_{b_{n}}(\xi_{n};\Lam)\cdot a_{n}$
 for $\xi$ in (\ref{eq:CantornfLam}).

\item\label{df:thtminus}
Let $\zeta=_{NF}\tilde{\theta}_{b}(\xi;\Lam)$ with $\tilde{\theta}_{b}(\xi;\Lam)>\xi$ and $b=\omega^{b_{0}}$,
and $c$ be an ordinal.
An ordinal $\tilde{\theta}_{-c}(\zeta;\Lam)$ is defined recursively as follows.
If $b\geq c$, then $\tilde{\theta}_{-c}(\zeta;\Lam)=\tilde{\theta}_{b-c}(\xi;\Lam)$.
Let $c>b$.
If $\xi>0$, then
$\tilde{\theta}_{-c}(\zeta;\Lam)=\tilde{\theta}_{-(c-b)}(\tilde{\theta}_{b_{m}}(\xi_{m};\Lam);\Lam)$ for the head term 
$hd(\xi)=\tilde{\theta}_{b_{m}}(\xi_{m};\Lam)$ of 
$\xi$ in (\ref{eq:CantornfLam}).
If $\xi=0$, then let $\tilde{\theta}_{-c}(\zeta;\Lam)=0$.

\item\label{df:LamI}
Let $\xi<\Gamma(\mI_{N})$ be such that
$SC_{\mI_{N}}(\xi)\subset\Lam$ for a strongly critical number $\Lam<\mI_{N}$.
Then $\xi[\Lam : \mI_{N}]$ denotes an ordinal$<\Gamma(\Lam)$ obtained from $\xi$
by changing the base $\mI_{N}$ into $\Lam$.
This means that
$\xi[\Lam:\mI_{N}]$ is obtained from $\xi$ in (\ref{eq:CantornfLam})  
by replacing
$\tht_{b_{i}}(\xi_{i};\mI_{N})\cdot a_{i}$ by $\tht_{b_{i}}(\xi_{i}[\Lam:\mI_{N}];\Lam)\cdot a_{i}$.

\end{enumerate}
}
\end{definition}

\bprp\label{prp:thtcollapse}
Let $\xi,\zeta<\Gamma(\mI_{N})$ be such that
$SC_{\mI_{N}}(\xi,\zeta)\subset\Lam$ for a strongly critical number $\Lam<\mI_{N}$.
Then
$\xi<\zeta$ iff $\xi[\Lam:\mI_{N}]<\zeta[\Lam:\mI_{N}]$.
\eprp

\begin{definition}\label{df:Lam3}
{\rm

  \begin{enumerate}
  \item
A function $f:\Lam\to\Gamma(\Lam)$ with a \textit{finite} support
\\
${\rm supp}(f)=\{c<\Lam: f(c)\neq 0\}$ is said to be a \textit{finite function}
with \textit{base} $\Lam$
if 
$\fal i>0(a_{i}=1)$ and $a_{0}=1$ when $b_{0}>1$
in
$f(c)=_{NF}\tilde{\theta}_{b_{m}}(\xi_{m};\Lam)\cdot a_{m}+\cdots+\tilde{\theta}_{b_{0}}(\xi_{0};\Lam)\cdot a_{0}$
for any $c\in{\rm supp}(f)$.

It is identified with the finite function $f\!\upharpoonright\! {\rm supp}(f)$.
When $c\not\in {\rm supp}(f)$, let $f(c):=0$.
$f,g,h,\ldots$ range over finite functions.

Let
$\mathrm{fld}(f):=\bigcup\{\{c,f(c)\}: c\in {\rm supp}(f)\}$.

For an ordinal $c$, $f_{c}$ and $f^{c}$ are restrictions of $f$ to the domains
\\
${\rm supp}(f_{c})=\{d\in{\rm supp}(f): d< c\}$, ${\rm supp}(f^{c})=\{d\in{\rm supp}(f): d\geq c\}$, resp.
$g_{c}*f^{c}$ denotes the concatenated function such that
${\rm supp}(g_{c}*f^{c})={\rm supp}(g_{c})\cup {\rm supp}(f^{c})$, 
$(g_{c}*f^{c})(a)=g(a)$ for $a<c$, and
$(g_{c}*f^{c})(a)=f(a)$ for $a\geq c$.

For finite functions $f,g$,
$f\leq g:\Lrarw\forall c(f(c)\leq g(c))$.

\item\label{df:Exp2.5}
Let $f$ be a finite function and $c\leq \Lam$, $\xi<\Gam(\Lam)$ ordinals.
A relation $f<_{\Lam}^{c}\xi$ is defined by induction on the
cardinality of the finite set $\{d\in {\rm supp}(f): d>c\}$ as follows.
If $f^{c}=\emptyset$, then $f<_{\Lam}^{c}\xi$ holds.
Let $f^{c}\neq\emptyset$.
If $f^{c+1}=\emptyset$, then
$f<_{\Lam}^{c}\xi$ iff 
$f(c)< \xi$.
Otherwise for
$d=\min\{d> 0: c+d\in {\rm supp}(f)\}$,
$f<_{\Lam}^{c}\xi$ iff 
there exists a segment $\mu$ of $\xi$ such that
$f(c)< \mu$
and 
$f<_{\Lam}^{c+d} \tilde{\theta}_{-d}(tl(\mu);\Lam)$, where
$tl(\mu)$ is the tail of $\mu$ with base $\Lam$.

\end{enumerate}

}
\end{definition}

The following Proposition \ref{prp:tht4} is shown in \cite{singlewfprf}.

\bprp\label{prp:tht4}
\begin{enumerate}
\item\label{prp:tht4.1}
$\zeta\leq\xi<\Gamma(\Lam) \Rightarrow \tilde{\theta}_{-c}(\zeta;\Lam)\leq\tilde{\theta}_{-c}(\xi;\Lam)$.

\item\label{prp:tht4.2}
$\tilde{\theta}_{c}(\tilde{\theta}_{-c}(\zeta;\Lam);\Lam) \leq\zeta$ for $\zeta<\Gamma(\Lam)$.
\end{enumerate}
\eprp

Although the following Proposition \ref{prp:idless} is shown in \cite{singlestable},
let us reproduce its proof.

\bprp\label{prp:idless}
$f<_{\Lam}^{c}\xi\leq\zeta \Rightarrow f<_{\Lam}^{c}\zeta$.
\eprp
\bprf
By induction on the cardinality $n$
of the finite set
$\{d\in {\rm supp}(f): d> c\}=\{c+d_{1}<\cdots<c+d_{n}\}$ with $c<c+d_{1}$.
If $n=0$, then there is nothing to prove.
Let $n>0$.
We have $f(c)<\mu$, 
and $f<_{\Lam}^{c+d_{1}}\tilde{\theta}_{-d_{1}}(tl(\mu);\Lam)$ 
for a segment $\mu$ of $\xi$.
We show the existence of a segment $\lambda$ of $\zeta$ such that 
$\mu\leq\lambda$, 
and $\tilde{\theta}_{-d_{1}}(tl(\mu);\Lam) \leq \tilde{\theta}_{-d_{1}}(tl(\lambda);\Lam)$.
Then IH yields $f<_{\Lam}^{c+d_{1}}\tilde{\theta}_{-d_{1}}(tl(\lambda);\Lam)$, and $f<_{\Lam}^{c}\zeta$ follows.

If $\mu$ is a segment of $\zeta$, then $\lambda=\mu$ works.
Otherwise $\xi<\zeta$ and there exists a segment $\lambda$ of $\zeta$ such that $\mu<\lambda$,
and
$tl(\mu)<tl(\lambda)$.
We obtain $\tilde{\theta}_{-d_{1}}(tl(\mu);\Lam)\leq\tilde{\theta}_{-d_{1}}(tl(\lambda);\Lam)$
by Proposition \ref{prp:tht4}.\ref{prp:tht4.1}.
\eprf

\subsection{Skolem hulls and Mahlo classes}\label{subsec:psi}

In this subsection Skolem hulls $C_{a}(\beta)$, collapsing functions $\psi$ and
Mahlo classes $Mh^{a}_{i,c}(\xi)$ are introduced. 
$\psi$-functions are introduced in Buchholz\cite{Buchholz86}.

\begin{definition}
{\rm
Let $A\subset\mathbb{I}_{N}$ be a set, and $\alpha\leq\mathbb{I}_{N}$ a limit ordinal.
\[
\alpha\in M(A) :\Lrarw A\cap\alpha \mbox{ is stationary in }
\alpha
\Lrarw \mbox{ every club subset of } \alpha \mbox{ meets } A.
\]

}
\end{definition}

In the following Definition \ref{df:Cpsiregularsm}, 
$\varphi\alpha\beta=\varphi_{\alpha}(\beta)$ denotes the binary Veblen function on $(\mI_{N})^{\dagger 0}$.
For 
$a<\varepsilon(\mI_{N})$,
$\bet,c<\mI_{N}$, and
$\xi<\Gamma(\mI_{N})$, 
define simultaneously 
classes $C_{a}(\bet)\subset\Gamma(\mI_{N})$,
$Mh^{a}_{i,c}(\xi)\subset(\mathbb{I}_{N}+1)\,(i>0)$, and 
ordinals $\psi_{\mI_{N}}(a)\leq\mI_{N}$ and $\psi_{\kappa}^{f}(a)\leq\kappa$ by recursion on ordinals $a$ as follows.

\begin{definition}\label{df:Cpsiregularsm}
{\rm
Let $a<\varepsilon(\mI_{N})$,
$\bet,c<\mI_{N}$, and
$\xi<\Gamma(\mI_{N})$.

\benu

\item\label{df:Cpsiregularsm.1}
(Inductive definition of $C_{a}(\bet)$)

\benu
\item\label{df:Cpsiregularsm.10}
$\{0,\Omega,\mI_{N}\}\cup \beta\subset C_{a}(\beta)$, where 
$\Ome\in 
SSt_{0}$.

\item\label{df:Cpsiregularsm.11}
If $x, y \in C_{a}(\beta)$,
then $x+y\in C_{a}(\beta)$
and 
$\varphi xy\in C_{a}(\beta)$.

\item\label{df:Cpsiregularsm.1343}
Let $\alpha=\psi_{\pi}(b)$ with $\pi\in C_{a}(\beta)\cap SSt_{0}\cap\mI_{N}$, 
$b\in C_{a}(\beta)\cap a$ such that
$\{\pi,b\}\subset C_{b}(\alpha)$.
Then $\alpha\in C_{a}(\beta)$.

\item\label{df:Cpsiregularsm.LS}
Let $\alp=\psi_{\mI_{N}}(b)$ with $b\in C_{a}(\beta)\cap 
a$.
Then 
$\alp\in C_{a}(\beta)\cap(LSt_{N}\cup\{\mI_{N}\})$.

\item\label{df:Cpsiregularsm.12}
Let $\alp\in C_{a}(\beta)\cap\mI_{N}$.
Then $\alp^{\dagger i}\in C_{a}(\beta)\cap SSt_{i}$ for each $0< i\leq N$.

\item\label{df:Cpsiregularsm.1344}
Let $\alpha=\psi_{\pi}^{f}(b)$ 
with
$b<a$, and a finite function $f:\mI_{N}\to\Gamma(\mI_{N})$ such that
$\{\pi,b\}\cup \mathrm{fld}(f)\subset C_{a}(\beta)\cap
C_{b}(\alpha)$.
Then $\alpha\in C_{a}(\beta)$.

\eenu

\item\label{df:Cpsiregularsm.2}
 (Definitions of $Mh^{a}_{i,c}(\xi)$ and $Mh^{a}_{i,c}(f)$ for $0<i\leq N$)
\\
The classes $Mh^{a}_{i,c}(\xi)$ are defined for 
$c<\mI_{N}$, $a<\veps(\mI_{N})$ and $\xi<\Gamma(\mI_{N})$.
By main induction on ordinals $\pi<\mI_{N}$
with subsidiary induction on $c<\mI_{N}$ 
we define 
$\pi\in Mh^{a}_{i,c}(\xi)$ iff $\pi\in LSt_{i-1}$,
$\{a,c,\xi\}\subset C_{a}(\pi)$ and
 the following condition is met for any 
 finite functions $f,g:\mI_{N}\to\Gamma(\mI_{N})$
 such that $f<_{\mI_{N}}^{c}\xi$:
 \[
 \mathrm{fld}(f,g) \subset C_{a}(\pi) 
 \spand \pi\in Mh^{a}_{i,0}(g_{c})
 \Rarw \pi\in M(Mh^{a}_{i,0}(g_{c}*f^{c}))
\]
where $ \mathrm{fld}(f,g)= \mathrm{fld}(f)\cup \mathrm{fld}(g)$ and
\beqnarrs
Mh^{a}_{i,c}(f) & := & \bigcap\{Mh^{a}_{i,d}(f(d)): d\in {\rm supp}(f^{c})\}
\\
& = &
\bigcap\{Mh^{a}_{i,d}(f(d)): c\leq d\in {\rm supp}(f)\}.
\eeqnarrs
$Mh^{a}_{i,0}(g_{c})=\bigcap\{Mh^{a}_{i,d}(g(d)): d\in {\rm supp}(g_{c})\}
=\bigcap\{Mh^{a}_{i,d}(g(d)): c> d\in {\rm supp}(g)\}$.
When $f=\emptyset$ or $f^{c}=\emptyset$, let $Mh^{a}_{i,c}(\emptyset):=LSt_{i-1}$.

\item\label{df:Cpsiregularsm.3}
 (Definition of $\psi_{\pi}^{f}(a)$)
\\
 Let $a,\pi$ be ordinals,
and
 $f:\mI_{N}\to\Gamma(\mI_{N})$ a finite function.
Then $\psi_{i,\pi}^{f}(a)$ denotes the least ordinal $\kap<\pi$ such that
\beqn\label{eq:Psivec}
\kappa\in Mh^{a}_{i,0}(f) \spand
 C_{a}(\kappa)
 \cap\pi\subset\kappa
 \spand
   \{\pi,a\}\cup \mathrm{fld}(f)\subset C_{a}(\kappa)
\eeqn
if such a $\kap$ exists. Otherwise set $\psi_{i,\pi}^{f}(a)=\pi$.

\item
$\psi_{\Ome}(a):=
\min
(\{\Ome\}\cup
\{\bet : C_{a}(\bet)\cap\Ome\subset\bet
\}
)
$
and
\beqn\label{eq:psiI}
\psi_{\mI_{N}}(a):=
\min(
\{\mI_{N}\}\cup\{\kap\in LSt_{N} : C_{a}(\kap)\cap\mI_{N}\subset\kap 
\}
)
\eeqn

\item
For classes $A\subset\mathbb{I}_{N}$, let
$\alpha\in M^{a}_{i,c}(A)$ iff $\alpha\in A$ and
for any finite functions
$g:\mI_{N}\to\Gamma(\mI_{N})$
\begin{equation}\label{eq:Mca}
\alpha\in Mh_{i,0}^{a}(g_{c}) \spand 
\mathrm{fld}(g_{c})\subset C_{a}(\alpha) \Rarw
\alpha\in M \!\left( Mh_{i,0}^{a}(g_{c}) \cap A \right)
\end{equation}

\eenu
}

\edf

The following Propositions \ref{prp:Mhless}, \ref{prp:MMh} and \ref{prp:definability} are seen as in \cite{singlestable}.

\bprp\label{prp:Mhless}
Assume $\pi\in Mh^{a}_{i,c}(\zeta)$ and $\xi<\zeta$ with $\xi\in C_{a}(\pi)$.
Then 
$\pi\in Mh^{a}_{i,c}(\xi)\cap M^{a}_{i,c}(Mh^{a}_{i,c}(\xi))$.
\eprp
\bprf
Proposition \ref{prp:idless} yields
$\pi\in Mh^{a}_{i,c}(\xi)$. 
$\pi\in M^{a}_{i,c}(Mh^{a}_{i,c}(\xi))$ is seen from the function $f$ such that
$f<_{\mI_{N}}^{c}\zeta$ with ${\rm supp}(f)=\{c\}$ and $f(c)=\xi$.
\eprf

\bprp\label{prp:MMh}
Suppose $\pi\in Mh^{a}_{i,c}(\xi)$.
\begin{enumerate}
\item\label{prp:MMh.1}
Let $f<_{\mI_{N}}^{c}\xi$ with $\mathrm{fld}(f) \subset C_{a}(\pi)$.
Then
$\pi\in M_{i,c}^{a}(Mh^{a}_{i,c}(f^{c}))$.

\item\label{prp:MMh.2}
Let $\pi\in M^{a}_{i,d}(A)$ for $d>c$ and $A\subset\mathbb{I}_{N}$.
Then $\pi\in M_{i,c}^{a}(Mh^{a}_{i,c}(\xi)\cap A)$.
\end{enumerate}
\eprp
\bprf
\ref{prp:MMh}.\ref{prp:MMh.1}.
Let $g$ be a function such that $\pi\in Mh_{i,0}^{a}(g_{c})$ with $\mathrm{fld}(g_{c})\subset C_{a}(\pi)$.
We obtain
$\pi\in M \!\left( Mh_{i,0}^{a}(g_{c}) \cap Mh^{a}_{i,c}(f^{c}) \right)$
by Definition \ref{df:Cpsiregularsm}.\ref{df:Cpsiregularsm.2} of $\pi\in Mh^{a}_{i,c}(\xi)$.
\\
\ref{prp:MMh}.\ref{prp:MMh.2}.
Let $\pi\in M^{a}_{i,d}(A)$ for $d>c$. Then $\pi\in Mh^{a}_{i,c}(\xi)\cap A$.
Let $g$ be a function such that $\pi\in Mh_{i,0}^{a}(g_{c})$ with $\mathrm{fld}(g_{c})\subset C_{a}(\pi)$.
We obtain by (\ref{eq:Mca}) and $d>c$ with the function $g_{c}*h$,
$\pi\in M\! \left( Mh_{i,0}^{a}(g_{c}) \cap Mh^{a}_{i,c}(\xi)\cap A \right)$, where
${\rm supp}(h)=\{c\}$ and $h(c)=\xi$.
\eprf

\bprp\label{prp:definability}
Each of 
$x\in C_{a}(y)$,
$x\in Mh^{a}_{i,c}(f)$ and $x=\psi^{f}_{\kappa}(a)$
is a $\Delta_{1}(\{St_{i}\}_{0<i\leq N})$-predicate in ${\sf ZFC}(\{St_{i}\}_{0<i\leq N})$.
\eprp
\bprf
An inspection of Definition \ref{df:Cpsiregularsm} shows that
$x\in C_{a}(y)$, $\psi^{f}_{\kappa}(a)$ and $x\in Mh^{a}_{i,c}(f)$ 
are simultaneously defined by recursion on 
$a<\varepsilon(\mathbb{I}_{N})$, in which $x\in Mh^{a}_{i,c}(f)$ is defined by recursion on ordinals $x<\mathbb{I}_{N}$
with subsidiary recursion on $c<\mathbb{I}_{N}$.
\eprf

\subsection{A small large cardinal hypothesis}\label{subsec:largecardinal}

It is convenient for us to assume the existence of a small large cardinal in justification of 
Definition \ref{df:Cpsiregularsm}.
\textit{Shrewd cardinals} as well as $\mathcal{A}$\textit{-shrewd cardinals} 
 are introduced by M. Rathjen\cite{RathjenAFML2}.
 
\bdf\label{df:shrewd}
{\rm (Rathjen\cite{RathjenAFML2})\\
Let $\eta>0$.
A cardinal $\kappa$ is 
$\eta$-\textit{shrewd} iff for any $P\subset V_{\kappa}$, and a
set-theoretic formula $\varphi(x,y)$  if
$V_{\kappa+\eta}\models\varphi[P,\kappa]$, then there are $0<\kappa_{0},\eta_{0}<\kappa$ such that
$V_{\kappa_{0}+\eta_{0}}\models\varphi[P\cap V_{\kappa_{0}},\kappa_{0}]$.
For classes $\mathcal{A}$,
$\kappa$ is 
$\mathcal{A}$-$\eta$-\textit{shrewd} iff for any $P\subset V_{\kappa}$, and a formula 
$\varphi(x,y)$ in the language $\{\in,R\}$ with a unary predicate $R$
if
$(V_{\kappa+\eta};\mathcal{A})\models\varphi[P,\kappa]$, then there are $0<\kappa_{0},\eta_{0}<\kappa$ such that
$(V_{\kappa_{0}+\eta_{0}};\mathcal{A})\models\varphi[P\cap V_{\kappa_{0}},\kappa_{0}]$,
where $(V_{\alp};\mathcal{A})$ denotes the structure $(V_{\alp},\in;\mathcal{A}\cap V_{\alp})$,
and for the formulas $\vphi$ in the language 
$\{\in,R\}$, $R(t)$ is interpreted as $t\in\mathcal{A}\cap V_{\alp}$
in $(V_{\alp};\mathcal{A})\models \vphi$.
}
\edf
Obviously each $\mathcal{A}$-$\eta$-shrewd cardinal is $\eta$-shrewd.
We see easily that each $\eta$-shrewd cardinal is regular.
A cardinal $\kap$ is said to be \textit{$(<\eta)$-shrewd} 
[\textit{$\mathcal{A}$-$(<\eta)$-shrewd}] if
$\kap$ is $\del$-shrewd [$\mathcal{A}$-$\del$-shrewd] for every $\del<\eta$, resp.
%\bprf
%Let $\kap$ be an $\eta$-shrewd singular cardinal, and $f:\alp\to\kap$ be a cofinal map 
%from an $\alp<\kap$.
%We can assume that $\bet<f(\bet)$ and $f(\lam)=\sup_{\bet<\lam}f(\bet)$ 
%for limit ordinals $\lam<\alp$.
%Then $V_{\kap+\eta}\models\vphi[P,\kap]$ holds for 
%$P=\{(\bet,\gam)\in\kap\times\kap: \bet<\alp, \gam=f(\bet)\}$, where
%$\vphi(x,y)$ says that `$x$ is a graph of a cofinal map from an ordinal$<y$ to $y$'.
%By shrewdness pick $0<\kap_{0},\eta_{0}<\kap$ such that
% $V_{\kap_{0}+\eta_{0}}\models\vphi[P\cap V_{\kap_{0}},\kap_{0}]$.
%Let $\bet<\alp$ be the least ordinal such that $\kap_{0}\leq f(\bet)$.
%Then $\bet<\kap_{0}$, and $P\cap V_{\kap_{0}}$ is not a graph of any map.
%\eprf

On the other side subtle cardinals are introduced by R. Jensen and K. Kunen.
%A cardinal $\kap$ is subtle if for any sequence $(S_{\alp})_{\alp<\kap}$ with $S_{\alp}\subset\alp$,
%and club subset $C$ of $\kap$, there are $\bet<\del$ such that $\{\bet,\del\}\subset C$ and
%$S_{\del}\cap\bet=S_{\bet}$.
The following Lemma \ref{lem:2.7R} is shown in \cite{RathjenAFML2} by Rathjen.
\blem\label{lem:2.7R}{\rm (Lemma 2.7 of \cite{RathjenAFML2})}\\
Let $\pi$ be a subtle cardinal.
The set 
$\{\kap\in V_{\pi}: (V_{\pi};\mathcal{A}) \models \mbox{`$\kap$ is $\mathcal{A}$-shrewd'}\}$
of $\mathcal{A}$-shrewd cardinals in $(V_{\pi};\mathcal{A})$
is stationary in $\pi$ for each class $\mathcal{A}$.
\elem

\bdf\label{df:2shrewd}
{\rm
Let $\pi$ be a cardinal.
The classes $\mathcal{B}_{n}$ and $\mathcal{A}_{n}$ are defined recursively for $n<\ome$.
Let 
\beqnarrs
\mathcal{B}_{0} & = & \{\kap\in V_{\pi}: V_{\pi}\models\mbox{`$\kap$ is an uncountable cardinal'}\}
\\
\mathcal{A}_{n} & = & \{\la i,\sig\ra: i\leq n, \sig\in\mathcal{B}_{i}\}
\\
\mathcal{B}_{n+1} & = &
\{\kap\in V_{\pi}: (V_{\pi};\mathcal{A}_{n})\models\mbox{`$\kap$ is an 
$\mathcal{A}_{n}$-shrewd cardinal'}\}.
\eeqnarrs
We say that a cardinal $\kap\in V_{\pi}$ is $n$-\textit{shrewd} in $\pi$ iff 
$\kap\in\mathcal{B}_{n}$.
An $n$-shrewd carinal is an \textit{$n$-shrewd limit} iff the set of $n$-shrewd cardinals is cofinal in it.
}
\edf
$\mathcal{B}_{1}$ is the set of shrewd cardinals in $V_{\pi}$, and a $1$-shrewd cardinal
is a shrewd cardinal in $\pi$.
Each $\mathcal{A}_{n+1}$-shrewd cardinal is $\mathcal{A}_{n}$-shrewd, and
each $(n+1)$-shrewd cardinal is $n$-shrewd.

\blem\label{lem:2shrewd}
Let $\pi$ be a subtle cardinal.
\benu
\item\label{lem:2shrewd.1}
The set of $n$-shrewd cardinals in $\pi$ is stationary in $\pi$ for each $n<\ome$.
\item\label{lem:2shrewd.2}
Let $\kap$ be an $(n+1)$-shrewd cardinal in $\pi$.
If 
$(V_{\kappa+\eta};\mathcal{A}_{n})\models\varphi[P,\kappa]$
for $0<\eta<\pi$, $P\subset V_{\kappa}$ and a formula $\varphi(x,y)$ in $\{\in,R\}$,
then there are an $n$-shrewd limit $\kap_{0}<\kap$ 
and
$0<\eta_{0}<\kappa$ such that
$(V_{\kappa_{0}+\eta_{0}};\mathcal{A}_{n})\models\varphi[P\cap V_{\kappa_{0}},\kappa_{0}]$.
\eenu
\elem
\bprf
\ref{lem:2shrewd}.\ref{lem:2shrewd.1}.
From Lemma \ref{lem:2.7R} we see that
the set of 
$\mathcal{A}_{n-1}$-shrewd cardinals is stationary in
a subtle cardinal $\pi$.
\\
\ref{lem:2shrewd}.\ref{lem:2shrewd.2}.
Let $\kap$ be an $(n+1)$-shrewd cardinal in $\pi$.
Then $\kap$ is $n$-shrewd, and hence
$(V_{\kap+\eta};\mathcal{A}_{n})\models\exi x(x\in P)\land R(\la n,\kap\ra)$
for each $P=\{\alp\}\subset V_{\kap}$ with $\kap<\kap+\eta<\pi$.
Since $\kap$ is $\mathcal{A}_{n}$-shrewd, there are $0<\kap_{0},\eta_{0}<\kap$ such that
$(V_{\kap_{0}+\eta_{0}};\mathcal{A}_{n})\models\exi x(x\in P\cap V_{\kap_{0}}) \land R(\la n,\kap_{0}\ra)$.
This means that $\alp<\kap_{0}$ is $n$-shrewd.
Therefore
$\kap$ is an $n$-shrewd limit.

Suppose 
$(V_{\kappa+\eta};\mathcal{A}_{n})\models\varphi[P,\kappa]$
for $0<\eta<\pi$, $P\subset V_{\kappa}$ and a formula $\varphi(x,y)$ in $\{\in,R\}$.
Then $(V_{\kappa+\eta};\mathcal{A}_{n})\models\varphi[P,\kappa]\land R(\la n,\kap\ra)\land
\fal \alp<\kap\exi \sig<\kap(\sig>\alp\land R(\la n,\sig\ra))$.
Since $\kap$ is $\mathcal{A}_{n}$-shrewd, 
there are an $n$-shrewd limit $\kap_{0}<\kap$
and
$0<\eta_{0}<\kappa$ such that
$(V_{\kappa_{0}+\eta_{0}};\mathcal{A}_{n})\models\varphi[P\cap V_{\kappa_{0}},\kappa_{0}]$.
\eprf
\\

\noindent
In this subsection we work in an extension $T$ of ${\sf ZFC}$ by adding
the axiom stating that
there exists a regular cardinal $\mI_{N}$ in which the set of $N$-shrewd cardinals is stationary.
$\Ome$ denotes the least uncountable ordinal $\ome_{1}$,
For $0<i\leq N$,
$St_{i}=\mathcal{B}_{i}$ the class of $i$-shrewd cardinals in $V_{\mI_{N}}$.
$LSt_{i}$ denotes the class of 
$i$-shrewd limits in $V_{\mI_{N}}$.
Let $St_{N+1}=SSt_{N+1}=\{\mI_{N}\}$ with $\mI_{N}=\Ome^{\dagger (N+1)}$.
Also $St_{0}$ denotes the class of uncountable cardinals in $V_{\mI_{N}}$, and
$LSt_{0}$ the class of limit cardinals in $V_{\mI_{N}}$.
A \textit{successor $n$-shrewd cardinal} is an $n$-shrewd cardinal in $V_{\mI_{N}}$, but not in $LSt_{n}$.

\blem\label{lem:welldefinedness.2}
$T\vdash \forall a<\Gam(\mI_{N})[
\psi_{\mI_{N}}(a)<\mI_{N}]$.
\elem
\bprf
We see that the set
$C=\{\kappa<\mI_{N} : 
C_{a}(\kappa)\cap\mI_{N}\subset\kappa\}$
 is a club subset of the regular cardinal $\mI_{N}$.
This shows the existence of a $\kappa\in LSt_{N} \cap C$, and hence
$\psi_{\mI_{N}}(a)<\mI_{N}$ by the definition (\ref{eq:psiI}).
\eprf
\\

$\alp^{\dagger i^{(k)}}$ is defined by recursion on $k<\ome$ by
$\alp^{\dagger i^{(0)}}=\alp$ and
$\alp^{\dagger i^{(k+1)}}=(\alp^{\dagger i^{(k)}})^{\dagger i}$.

\bprp\label{prp:comparisonrud}
Let 
$a\in C_{a}(\psi_{\mI_{N}}(a))$, $b\in C_{b}(\psi_{\mI_{N}}(b))$,
$c\in C_{c}(\psi_{\Ome}(c))$ and $d\in C_{d}(\psi_{\Ome}(d))$.
\benu

\item\label{prp:comparisonrud.1}
$\psi_{\mI_{N}}(a)<\psi_{\mI_{N}}(b)$ iff 
$a<b$.

\item\label{prp:comparisonrud.2}
$\Ome^{\dagger N^{(k)}}<\psi_{\mI_{N}}(b)$ for every $k<\ome$.

\item\label{prp:comparisonrud.3}
Let $\alp=\psi_{\mI_{N}}(a)$ and $0<k<\ome$. 
Then $\alp^{\dagger N^{(k)}}<\psi_{\mI_{N}}(b)$ iff $\alp<\psi_{\mI_{N}}(b)$.
$\psi_{\mI_{N}}(b)<\alp^{\dagger N^{(k)}}$ iff $\psi_{\mI_{N}}(b)\leq\alp$.

\item\label{prp:comparisonrud.4}
$\psi_{\Ome}(c)<\psi_{\Ome}(d)$ iff 
$c<d$.

\item\label{prp:comparisonrud.5}
If $x<y$, then $\psi_{\mI_{N}}(x)\leq\psi_{\mI_{N}}(y)$.
\eenu
\eprp
\bprf
\ref{prp:comparisonrud}.\ref{prp:comparisonrud.2} and \ref{prp:comparisonrud}.\ref{prp:comparisonrud.3}. 
Let $\bet=\psi_{\mI_{N}}(b)$.
By the definition (\ref{eq:psiI}) and 
$\Ome\in C_{b}(\bet)\cap\mI_{N}\subset\bet$ we obtain $\Ome<\bet$.
Let $\alp\in\{\Ome,\psi_{\mI_{N}}(a)\}$.
If $\alp<\bet$, then $\bet\in LSt_{N}$ yields $\alp^{\dagger N^{(k)}}<\bet$.
\\
\ref{prp:comparisonrud}.\ref{prp:comparisonrud.5}.
We obtain $\psi_{\mI_{N}}(y)\in LSt_{N}$ and 
$C_{x}(\psi_{\mI_{N}}(y))\cap\mI_{N}\subset C_{y}(\psi_{\mI_{N}}(y))\cap\mI_{N}\subset\psi_{\mI_{N}}(y)$
by $x<y$ and Lemma \ref{lem:welldefinedness.2}.
Hence $\psi_{\mI_{N}}(x)\leq\psi_{\mI_{N}}(y)$.
\eprf

\subsection{$\psi$-functions}\label{subsec:psif}

In this subsection we work in ${\sf ZFC}(\{St_{i}\}_{0<i\leq N})$ with $St_{i}=\mathcal{B}_{i}$, and
show that $\psi_{i,\mS}^{f}(a)<\mS$ for $i$-shrewd cardinal $\mS$
in Lemma \ref{lem:limitcollapse.1},
and introduce an \textit{irreducibility} of finite functions in Definition \ref{df:irreducible}
using Lemma \ref{lem:stepdown},
which is needed to define a normal form in ordinal notations.

\blem\label{lem:welldefinedness.1}
Let $\mS$ be an $i$-shrewd cardinal with $0<i\leq N$,
$a<\veps(\mI_{N})$, $h:\mI_{N}\to\Gamma(\mI_{N})$ a finite function 
with
$\{a\}\cup SC(h)\subset C_{a}(\mS)$. Then
$\mS\in Mh^{a}_{i,0}(h)\cap M(Mh^{a}_{i,0}(h))$.
\elem
\bprf
By induction on $\xi<\Gamma(\mI_{N})$
we show $\mS\in Mh^{a}_{i,c}(\xi)$
for $\{a,c,\xi\}\subset C_{a}(\mS)$.

Let $\{a,c,\xi\}\cup \mathrm{fld}(f)\subset C_{a}(\mS)$ 
with $f<_{\mI_{N}}^{c}\xi$ and $a<\veps(\mI_{N})$.
We show $\mS\in M^{a}_{i,c}(Mh^{a}_{i,c}(f^{c}))$, 
which yields $\mS\in Mh^{a}_{i,c}(\xi)$.
IH yields $\mS\in Mh^{a}_{i,c}(f^{c})$ by Proposition \ref{prp:tht4}.\ref{prp:tht4.2},
$\tilde{\theta}_{-e}(\zeta;\mI_{N})\leq\zeta$.
By the definition (\ref{eq:Mca}) it suffices to show that
\[
\forall g
[
\mS\in Mh_{i,0}^{a}(g_{c}) \,\&\, SC(g_{c})\subset C_{a}(\mS) \Rightarrow
\mS\in M\left( Mh_{i,0}^{a}(g_{c}) \cap Mh^{a}_{i,c}(f^{c})\right)
]
.\]

Let $g:\mI_{N}\to\Gamma(\mI_{N})$ be a finite function such that $\mathrm{fld}(g_{c})\subset C_{a}(\mS)$ and
$\mS\in Mh_{i,0}^{a}(g_{c})$.
We have to show $\mS\in M(A\cap B)$ for $A=Mh_{i,0}^{a}(g_{c})\cap\mS$ and 
$B=Mh^{a}_{i,c}(f^{c})\cap\mS$.
Let $C$ be a club subset of $\mS$.

We have 
$\mS\in Mh_{i,0}^{a}(g_{c})\cap Mh^{a}_{i,c}(f^{c})$,
and
$\{a\}\cup \mathrm{fld}(g_{c},f^{c})\subset C_{a}(\mS)$.
Pick a $b<\mS$ so that $\{a\}\cup \mathrm{fld}(g_{c},f^{c})\subset C_{a}(b)$.
Since the cardinality of the set $ C_{a}(\mS)$ is equal to $\mS$, pick a bijection 
$F:\mS\to  C_{a}(\mS)$.
Each $\alpha<\Gamma(\mI_{N})$ with $\alpha\in C_{a}(\mS)$ is identified with its code, denoted
by $F^{-1}(\alpha)<\mS$.
Let $P$ be the class
$P=\{(\pi,d,\alpha)\in\mS^{3} : \pi\in Mh^{F(\alp)}_{i,F(d)}(F(\xi))\}$,
where 
$F(d)\in  C_{a}(\mS)\cap(c+1)$ and 
$F(\alpha)<\Gamma(\mI_{N})$ with 
$\{F(d), F(\alpha)\}\subset C_{a}(\pi)$.
For fixed $i$, $a$ and $c$, the set
$\{(d,\zeta)\in
\left( C_{a}(\mS)\cap(c+1)\right)
\times\Gamma(\mathbb{I}_{N})   : \mS\in Mh^{a}_{i,d}(\zeta)\}$
is defined from the classes $P$ and $\{St_{j}\}_{j<i}$ by recursion on ordinals $d\leq c$.

Let $\varphi$ be a formula in $\{\in\}\cup\{St_{j}\}_{j<i}$ 
such that $(V_{\mS+c^{\dagger i}};\{St_{j}\}_{j<i})\models\varphi[P, C,\mS,b]$ iff
$\mS\in Mh_{i,0}^{a}(g_{c})\cap Mh^{a}_{i,c}(f^{c})$ and $C$ is a club subset of $\mS$,
where $\{St_{j}\}_{j<i}=\mathcal{A}_{i-1}$.
Since $\mS$ is $i$-shrewd in $V_{\mI_{N}}$, pick $b<\mS_{0}<\eta<\mS$ such that
$(V_{\mS_{0}+\eta};\{St_{j}\}_{j<i})\models\varphi[P\cap \mS_{0},C\cap\mS_{0},\mS_{0},b]$.
We obtain $\mS_{0}\in A\cap B\cap C$.

Therefore $\mS\in Mh^{a}_{i,c}(\xi)$ is shown
for every
$\{c,\xi\}\subset C_{a}(\mS)$.
This yields $\mS\in Mh_{i,0}^{a}(h)$ for 
$\mathrm{fld}(h)\subset C_{a}(\mS)$.
$\mS\in M(Mh^{a}_{i,0}(h))$ follows from the $i$-shrewdness of $\mS$.
\eprf

\blem\label{lem:limitcollapse.1}
Let $\mS$ be 
an $i$-shrewd cardinal, $a$ an ordinal, and
$f:\mI_{N}\to\Gamma(\mI_{N})$ a finite function 
such that $\{a,\mS\}\cup \mathrm{fld}(f)\subset C_{a}(\mathbb{S})$.
Then
$\psi_{i,\mathbb{S}}^{f}(a)<\mathbb{S}$ holds.
\elem
\bprf
Suppose $\{a,\mS\}\cup \mathrm{fld}(f)\subset C_{a}(\mathbb{S})$.
By Lemma \ref{lem:welldefinedness.1} we obtain
$\mathbb{S}\in M(Mh^{a}_{i,0}(f))$.
The set
$C=\{\kappa<\mathbb{S}: 
 C_{a}(\kappa)\cap\mathbb{S}\subset\kappa, 
\{a,\mS\}\cup \mathrm{fld}(f)\subset C_{a}(\kappa)\}$
 is a club subset of 
 the regular cardinal $\mathbb{S}$,
 and $Mh^{a}_{i,0}(f)$ is stationary in $\mathbb{S}$.
This shows the existence of a $\kappa\in Mh_{i,0}^{a}(f)\cap C\cap\mathbb{S}$, and hence
$\psi_{i,\mathbb{S}}^{f}(a)<\mathbb{S}$ by the definition (\ref{eq:Psivec}).
\eprf

\bprp\label{prp:welldefinedness.suc}
Let $\alp$ be either $0$ or an $i$-shrewd cardinal for $0<i\leq N$
and $\mS=\alp^{\dagger i}$.
Assume $\{a,\mS\}\cup \mathrm{fld}(f)\subset C_{a}(\mathbb{S})$
for an ordinal $a$ and a finite function $f:\mI_{N}\to\Gamma(\mI_{N})$.
Then
$\alp^{\dagger j}<\psi^{f}_{i,\mS}(a)$ for every $j<i$, and $\psi^{f}_{i,\mS}(a)\in LSt_{i-1}\setm St_{i}$.
\eprp
\bprf
Let $\kap=\psi^{f}_{i,\mS}(a)<\mS$.
We obtain $\alp\in C_{a}(\kap)$ by $\alp^{\dagger i}=\mS\in C_{a}(\kap)$,
and $\alp^{\dagger j}\in C_{a}(\kap)\cap\mS$ for $\mS\in LSt_{j}$.
$\alp<\kap$ is seen from 
$\alp^{\dagger j}\in C_{a}(\kap)\cap\mS\subset\kap$
in the definition (\ref{eq:Psivec}).
\eprf
\\

The following Lemma \ref{lem:stepdown} and
Corollary \ref{cor:stepdown} are seen as in \cite{singlestable}.

\blem\label{lem:stepdown}
Assume $\mathbb{I}_{N}>\pi\in Mh^{a}_{i,d}(\xi)\cap Mh^{a}_{i,c}(\xi_{0})$, $\xi_{0}\neq 0$,
and $d<c$.
Moreover let $\xi_{1}\in C_{a}(\pi)$ for 
$\xi_{1}\leq\tilde{\theta}_{c-d}(\xi_{0};\mI_{N})$,
and $tl(\xi)\geq \xi_{1}$ when $\xi\neq 0$.
Then
$\pi\in Mh^{a}_{i,d}(\xi+\xi_{1})\cap M^{a}_{i,d}(Mh^{a}_{i,d}(\xi+\xi_{1}))$.
\elem
\bprf
$\pi\in M^{a}_{i,d}(Mh^{a}_{i,d}(\xi+\xi_{1}))$ follows from $\pi\in Mh^{a}_{i,d}(\xi+\xi_{1})$
 and $\pi\in Mh^{a}_{i,c}(\xi_{0})\subset M^{a}_{i,c}(Mh^{a}_{i,c}(\emptyset))$
 by Proposition \ref{prp:MMh}.\ref{prp:MMh.1}.

Let $f$ be a finite function such that 
$\mathrm{fld}(f)\subset C_{a}(\pi)$, and
$f<_{\mI_{N}}^{d}\xi+\xi_{1}$.
We show $\pi\in M^{a}_{i,d}(Mh^{a}_{i,d}(f^{d}))$ by main induction on 
the cardinality of the finite set $\{e\in {\rm supp}(f): e>d\}$
 with subsidiary induction on $\xi_{1}$.

First let $f<_{\mI_{N}}^{d}\mu$ for a segment $\mu$ of $\xi$. By Proposition \ref{prp:Mhless} we obtain
$\pi\in Mh^{a}_{i,d}(\mu)$ and $\pi\in M^{a}_{i,d}(Mh^{a}_{i,d}(f^{d}))$.

In what follows let $f(d)=\xi+\zeta$ with $\zeta<\xi_{1}$.
By SIH we obtain $\pi\in Mh^{a}_{i,d}(f(d))\cap M^{a}_{i,d}(Mh^{a}_{i,d}(f(d)))$.
If $\{e\in {\rm supp}(f): e>d\}=\emptyset$, then $Mh^{a}_{i,d}(f^{d})=Mh^{a}_{i,d}(f(d))$, and we are done.
Otherwise let $e=\min\{e\in {\rm supp}(f): e>d\}$.
By SIH we can assume $f<_{\mI_{N}}^{e}\tilde{\theta}_{-(e-d)}(tl(\xi_{1});\mI_{N})$.
We obtain $f<_{\mI_{N}}^{e}\tilde{\theta}_{-(e-d)}(\tilde{\theta}_{c-d}(\xi_{0};\mI_{N});\mI_{N})=\tilde{\theta}_{-e}(\tilde{\theta}_{c}(\xi_{0};\mI_{N});\mI_{N})$
by $\xi_{1}\leq\tilde{\theta}_{c-d}(\xi_{0};\mI_{N})$, Propositions \ref{prp:idless} and 
\ref{prp:tht4}.\ref{prp:tht4.1}.
We claim that $\pi\in M^{a}_{i,c_{0}}(Mh_{i,c_{0}}^{a}(f^{c_{0}}))$ for $c_{0}=\min\{c,e\}$.
If $c=e$, then the claim follows from the assumption $\pi\in Mh_{i,c}^{a}(\xi_{0})$ and $f<_{\mI_{N}}^{e}\xi_{0}$.
Let $e=c+e_{0}>c$. Then $\tilde{\theta}_{-e}(\tilde{\theta}_{c}(\xi_{0};\mI_{N});\mI_{N})=\tilde{\theta}_{-e_{0}}(hd(\xi_{0});\mI_{N})$, and
$f<_{\mI_{N}}^{c}\xi_{0}$ with $f(c)=0$ yields the claim.
Let $c=e+c_{1}>e$. Then $\tilde{\theta}_{-e}(\tilde{\theta}_{c}(\xi_{0};\mI_{N});\mI_{N})=\tilde{\theta}_{c_{1}}(\xi_{0};\mI_{N})$.
MIH yields the claim.

On the other hand we have
$Mh_{i,d}^{a}(f^{d})=Mh^{a}_{i,d}(f(d))\cap Mh_{i,c_{0}}^{a}(f^{c_{0}})$.
$\pi\in Mh^{a}_{i,d}(f(d))\cap M^{a}_{i,c_{0}}(Mh_{i,c_{0}}^{a}(f^{c_{0}}))$ with $d<c_{0}$ yields by 
Proposition \ref{prp:MMh}.\ref{prp:MMh.2},
$\pi\in M^{a}_{i,d}(Mh^{a}_{i,d}(f(d))\cap Mh_{i,c_{0}}^{a}(f^{c_{0}}))$, i.e.,
$\pi\in M^{a}_{i,d}(Mh^{a}_{i,d}(f^{d}))$.
\eprf

\bdf\label{df:nfform2}
{\rm
Let $f, g:\mI_{N}\to\Gamma(\mI_{N})$ be finite functions.

$ Mh^{a}_{i,0}(g)\prec Mh^{a}_{i,0}(f)$ iff the following holds:
\[
\forall\pi\in Mh^{a}_{i,0}(f)
\left(
\mathrm{fld}(g)\subset C_{a}(\pi
) \Rightarrow \pi\in M(Mh^{a}_{i,0}(g))
\right)
.
\]
}
\edf

\bcor\label{cor:stepdown}
Let $f,g:\mI_{N}\to\Gamma(\mI_{N})$
 be finite functions and $c\in{\rm supp}(f)$.
Assume  that 
there exists an ordinal
$d<c$ 
such that
$(d,c)\cap {\rm supp}(f)=(d,c)\cap {\rm supp}(g)=\emptyset$, 
$g_{d}\leq f_{d}$, 
$g(d)<f(d)+\tilde{\theta}_{c-d}(f(c);\mI_{N})\cdot\omega$,
and
$g<_{\mI_{N}}^{c}f(c)$.

Then
$Mh^{a}_{i,0}(g)\prec Mh^{a}_{i,0}(f)$ holds.
In particular if $\pi\in Mh^{a}_{i,0}(f)$ and
$\mathrm{fld}(g)\subset C_{a}(\pi)$, then
$\psi_{i,\pi}^{g}(a)<\pi$.
\ecor
\bprf
Let $\pi\in Mh^{a}_{i,0}(f)=\bigcap\{Mh^{a}_{i,e}(f(e)): e\in {\rm supp}(f)\}$ and
$\mathrm{fld}(g)\subset C_{a}(\pi)$.
Lemma \ref{lem:stepdown} with 
$\pi\in Mh^{a}_{i,d}(f(d))\cap Mh^{a}_{i,c}(f(c))$ yields
\\
$\pi\in Mh^{a}_{i,d}(g(d))\cap M^{a}_{i,c}(Mh^{a}_{i,c}(g^{c}))$.

On the other hand we have 
$\pi\in Mh^{a}_{i,0}(g_{d})=\bigcap\{Mh^{a}_{i,e}(g(e)): e\in {\rm supp}(g)\cap d\}\subset
\bigcap\{Mh^{a}_{i,e}(f(e)): e\in {\rm supp}(f)\cap d\}$ by Proposition \ref{prp:Mhless}.
Hence $\pi\in M(Mh^{a}_{i,0}(g))$.

Now suppose $\mathrm{fld}(g)\subset C_{a}(\pi)$.
The set
$C=\{\kappa<\pi:  C_{a}(\kappa)\cap\pi\subset\kappa, \{\pi,a\}\cup \mathrm{fld}(g)\subset C_{a}(\kappa)\}$
 is a club subset of the regular cardinal $\pi$,
 and $Mh^{a}_{i,0}(g)$ is stationary in $\pi$.
This shows the existence of a $\kappa\in Mh_{i,0}^{a}(g)\cap C\cap\pi$, and hence
$\psi_{i,\pi}^{g}(a)<\pi$ by the definition (\ref{eq:Psivec}).
\eprf

\bdf\label{df:irreducible}
{\rm

An \textit{irreducibility} of finite functions $f:\mI_{N}\to\Gamma(\mI_{N})$
 is defined by induction on the cardinality
$n$ of the finite set ${\rm supp}(f)$.
If $n\leq 1$, $f$ is defined to be irreducible.
Let $n\geq 2$ and $c<c+d$ be the largest two elements in ${\rm supp}(f)$, and let $g$ be 
a finite function
such that ${\rm supp}(g)={\rm supp}(f_{c})\cup\{c\}$, $g_{c}=f_{c}$ and
$g(c)=f(c)+\tilde{\theta}_{d}(f(c+d);\mI_{N})$.

Then $f$ is irreducible iff 
$tl(f(c))>\tilde{\theta}_{d}(f(c+d);\mI_{N})$ and
$g$ is irreducible.

}
\edf

\bdf\label{df:lx}
 {\rm 
 Let  $f,g:\mI_{N}\to\Gamma(\mI_{N})$ 
 be irreducible finite functions, and $b$ an ordinal.
Let us define a relation $f<^{b}_{lx}g$
by induction on the cardinality $\#\{e\in{\rm supp}(f)\cup{\rm supp}(g): e\geq b\}$ as follows.
$f<^{b}_{lx}g$ holds iff $f^{b}\neq g^{b}$ and
for the ordinal $c=\min\{c\geq b : f(c)\neq g(c)\}$,
one of the following conditions is met:

\begin{enumerate}

\item\label{df:lx.23}
$f(c)<g(c)$ and let $\mu$ be the shortest segment of $g(c)$ such that $f(c)<\mu$.
Then for any $c<c+d\in{\rm supp}(f)$,  
if $tl(\mu)\leq\tilde{\theta}_{d}(f(c+d);\mI_{N})$, then 
$f<_{lx}^{c+d}g$ holds.

\item\label{df:lx.24}
$f(c)>g(c)$ and let $\nu$ be the shortest segment of $f(c)$ such that $\nu>g(c)$.
Then there exist a $c<c+d\in {\rm supp}(g)$ such that
$f<_{lx}^{c+d}g$ and
$tl(\nu)\leq \tilde{\theta}_{d}(g(c+d);\mI_{N})$.

\end{enumerate}

}
\edf

In \cite{singlewfprf}
the following Proposition \ref{lem:psinucomparison} is shown.

\bprp\label{lem:psinucomparison}
Let $f:\mI_{N}\to\Gamma(\mI_{N})$.
If $f<^{0}_{lx}g$, then
$Mh^{a}_{i,0}(f)\prec Mh^{a}_{i,0}(g)$.
\eprp

\bprp\label{prp:psicomparison}
Let $f,g:\mI_{N}\to\Gamma(\mI_{N})$ be 
irreducible functions, and assume that
$\psi_{i,\pi}^{f}(b)<\pi$ and $\psi_{i,\kappa}^{g}(a)<\kappa$.

Then $\psi_{i,\pi}^{f}(b)<\psi_{i,\kappa}^{g}(a)$ iff one of the following cases holds:
\begin{enumerate}
\item\label{prp:psicomparison.0}
$\pi\leq \psi_{i,\kappa}^{g}(a)$.

\item\label{prp:psicomparison.1}
$b<a$, $\psi_{i,\pi}^{f}(b)<\kappa$, and 
$\mathrm{fld}(f)\cup\{\pi,b
\}\subset C_{a}(\psi_{i,\kappa}^{g}(a))$.

\item\label{prp:psicomparison.2}
$b>a$, and $\mathrm{fld}(g)\cup\{\kappa,a
\}\not\subset C_{b}(\psi_{i,\pi}^{f}(b))$.

\item\label{prp:psicomparison.25}
$b=a$, $\kappa<\pi$, and 
$\kap\not\in C_{b}(\psi_{i,\pi}^{f}(b))$.

\item\label{prp:psicomparison.3}
$b=a$, $\pi=\kappa$, $\mathrm{fld}(f)\subset C_{a}(\psi_{i,\kappa}^{g}(a))$, and
$f<^{0}_{lx}g$.

\item\label{prp:psicomparison.4}
$b=a$, $\pi=\kappa$, 
$\mathrm{fld}(g)\not\subset C_{b}(\psi_{i,\pi}^{f}(b))$.

\end{enumerate}

\eprp
\bprf
This is seen from Proposition \ref{lem:psinucomparison} as in \cite{KPPiN}.
\eprf

\subsection{A computable notation system for $\Pi_{N}$-collection}\label{subsec:decidable}

Although Propositions \ref{prp:comparisonrud}, \ref{prp:welldefinedness.suc}, and
\ref{prp:psicomparison} suffice for us to define a computable notation system for
$ C_{\veps(\mI_{N})}(0)$, we need a notation system closed under 
Mostowski collapsings
to remove stable ordinals from derivations as in \cite{singlestable}, 
cf.\,section \ref{sec:proofonestep}.
Two new constructors $\mI_{N}[\cdot]$ and $\mS^{\dagger\vec{i}}[\rho/\mS]$ are used to generate terms in $OT(\mI_{N})$.

\bdf\label{df:prec}
{\rm
$\rho\prec\sig$ denotes the transitive closure 
of the relation 
\\
$\{(\rho,\sig): \exi f,a(\rho=\psi_{\sig}^{f}(a))\}$.
Let $\rho\preceq\sig:\Lrarw \rho\prec\sig \lor \rho=\sig$.
}
\edf

Let $\mS\in SSt_{i}$ and $\rho\prec\mS$.
We define a set $M_{\rho}= C_{b}(\rho)$ from $\rho$ in (\ref{eq:M}) 
in such a way that $ C_{b}(\rho)\cap\mS\subset\rho$.
Then a Mostowski collapsing $M_{\rho}\ni\alp\mapsto\alp[\rho/\mS]$ in Definition \ref{df:notationsystem.2M} 
maps ordinal terms $\alp\in M_{\rho}$ to $\alp[\rho/\mS]<\mS$ isomorphically.
The transitive collapse $(M_{\rho})^{[\rho/\mS]}=\{\alp[\rho/\mS]: \alp\in M_{\rho}\}$
is an initial segment in $OT(\mI_{N})$ such that
$(M_{\rho})^{[\rho/\mS]}<\kap$ if $\rho<\kap\prec\mS$.
Note that both $\rho$ and $\kap$ can be interpreted as uncountable cardinals, and
the cardinality of the set $M_{\rho}$ is equal to $\rho$.

Let us define simultaneously the followings:
A set $OT(\mI_{N})$ of terms over constants $0,\Ome,\mI_{N}$ and
constructors $+,\vphi,\psi$, $\mI_{N}[*]$, $*^{\dagger i}\,(0<i\leq N)$, and 
$*_{0}[*_{1}/*_{2}]$.
Its subsets
$SSt_{i}$, $LSt_{i}$ with $St_{i}=SSt_{i}\cup LSt_{i}$, and sets 
$M_{\rho}\, (\rho\in\Psi)$,
finite sets 
$K_{X}(\alp)$ of subterms of $\alp$ for $X\subset OT(\mI_{N})$.
Let $SSt=\bigcup_{0<i\leq N}SSt_{i}$ and $LSt=\bigcup_{0<i\leq N}LSt_{i}$.
For each $\mS\in SSt$, there exists a unique $i$ such that $\mS\in SSt_{i}$.

For $i>0$, $\kap\in St_{i}$ is intended to designate that $\kap$ is an $i$-shrewd cardinal,
or $\kap$ is an $i$-stable ordinal.
$\kap\in SSt_{i}$ [$\kap\in LSt_{i}$] is intended to designate that $\kap$ is a 
successor $i$-stable ordinal [$\kap$ is a limit of $i$-stable ordinals], resp.
$\kap\in St_{0}$ is intended to designate that $\kap$ is an uncountable cardinal,
or $\kap$ is either a recursively regular ordinal or their limit.
We have
$St_{i}=SSt_{i}\cup LSt_{i}$ with $SSt_{i}\cap LSt_{i}=\emptyset$, and
$St_{i+1}\subset LSt_{i}$.
If $\mS\in SSt_{i}$, then the ordinal term
$\psi_{\mS}^{f}(a)$ in Definition \ref{df:notationsystem}.\ref{df:notationsystem.5} denotes
the ordinal $\psi_{i,\mS}^{f}(a)$ in (\ref{eq:Psivec}) of Definition \ref{df:Cpsiregularsm}.\ref{df:Cpsiregularsm.3}.

$\alp=_{NF}\alp_{m}+\cdots+\alp_{0}$ means that $\alp=\alp_{m}+\cdots+\alp_{0}$ 
with $\alp_{m}\geq\cdots\geq\alp_{0}$
and each $\alp_{i}$ is a non-zero additive principal number.
$\alp=_{NF}\vphi\bet\gam$ means that $\alp=\vphi\bet\gam$ and $\bet,\gam<\alp$.

Sets $SC(\alp)$ of strongly critical numbers are slightly modified as $SC(\Ome)=SC(\mI_{N})=\emptyset$. Specifically
$SC(0)=\emptyset$, $SC(\alp)=\bigcup_{i\leq m}SC(\alp_{i})$ for $\alp=_{NF}\alp_{m}+\cdots+\alp_{0}$,
and
$SC(a)=SC(b)\cup SC(c)$ for $a=_{NF}\vphi_{b}(c)$.
$SC(\Ome)=SC(\mI_{N})=\emptyset$.
$SC(a)=\{a\}$ if $a\not\in\{\Ome,\mI_{N}\}$ is strongly critical.

For $\alp=\psi_{\pi}^{f}(a)$, let $m(\alp)=f$.
$\mathrm{fld}(f)=\bigcup\{\{c,f(c)\}: c\in {\rm supp}(f)\}$.
Immediate subterms of terms are defined as follows.
$k(\alp_{m}+\cdots+\alp_{0})=\{\alp_{0},\ldots,\alp_{m}\}$,
$k(\vphi\alp\bet)=\{\alp,\bet\}$, 
$k(\psi_{\mI_{N}}(a))=\{\mI_{N},a\}$, and
\\
$k(\psi_{\sigma}^{f}(\alpha))=\{\sigma,\alpha\}\cup \mathrm{fld}(f)$.

Note that in the following Definition \ref{df:notationsystem}, e.g.,
there is no clause for constructing $\kap=\psi_{\mS}(a)$ from $a$
for $\mS\not\in SSt$.

\bdf\label{df:veci}
{\rm
\benu
\item
$\alp\in\Psi :\Lrarw \exi \kap,f,a(\alp=\psi_{\kap}^{f}(a))$
and
$\alp\in\Psi_{\mS} :\Lrarw \exi \kap\preceq\mS\exi f,a(\alp=\psi_{\kap}^{f}(a))$.

\item
For sequences $\vec{i}=(i_{0},i_{1},\ldots, i_{n})$ of numbers and ordinals $\alpha$,
let $\alpha^{\dagger\vec{i}}=(\cdots((\alpha^{\dagger i_{0}})^{\dagger i_{1}})\cdots)^{\dagger i_{n}}$.

\item
By $\vec{i}\leq i$ let us understand that
$\vec{i}=(i_{0}, i_{1}, \ldots, i_{n})$ is a non-empty and non-increasing sequence
of numbers such that $0<i_{n}\leq\cdots\leq i_{1}\leq i_{0}\leq i$.
\eenu
}
\edf

\bdf\label{df:p0}
{\rm
\benu

\item\label{df:p0.0}
Let $\alp\preceq\psi_{\mS}^{g}(b)$ for an $\mS\in SSt$ and a $g$ with
$b=\mathtt{ p}_{0}(\alp)$.
Then let
\beqn\label{eq:M}
M_{\alp}  :=   C_{b}(\alp)
\eeqn

\item\label{df:p0.1}
For $\alp\in\Psi$, an ordinal $\mathtt{ p}_{0}(\alp)$ is defined.

\benu
\item\label{df:p0.11}
If $\alp\preceq\psi_{\mS}^{g}(b)$,
then $\mathtt{ p}_{0}(\alp)=b$.

\item\label{df:p0.12}
If there are $\rho$ and $\bet\in M_{\rho}$ such that
$LSt_{i}\ni\rho\prec\mS\in SSt_{i+1}$ and $\alp=\bet[\rho/\mS]$,
then
$\mathtt{ p}_{0}(\alp)=\mathtt{ p}_{0}(\bet)$.

\item\label{df:p0.13}
$\mathtt{ p}_{0}(\alp)=0$ otherwise.
\eenu

\item\label{df:p0.2}
$\alp^{\dagger}:=\alp^{\dagger 1}$.

\eenu
}
\edf

\bdf\label{df:notationsystem}
{\rm
  (Definitions of $OT(\mI_{N})$
  and $K_{X}(\alp)$)
\\
Let 
$St_{i}=SSt_{i}\cup LSt_{i}\subset OT(\mI_{N})$ with $SSt_{i}\cap LSt_{i}=\emptyset$ and $St_{i+1}\subset LSt_{i}$.
For $\del,\alp\in OT(\mI_{N})$,
$K_{\del}(\alp)=K_{X}(\alp)$, where 
$X=\{\bet\in OT(\mI_{N}):\bet<\del\}$.

 \benu
 \item\label{df:notationsystem.1}
$\{0,\Ome, \mI_{N}\}\subset OT(\mI_{N})$ and
$0^{\dagger i }\in SSt_{i}$ for $0<i\leq N$.
Let $St_{N+1}=\{\mI_{N}\}$.
$m(\alp)=K_{X}(\alp)=\emptyset$
for
$\alp\in\{0,\mI_{N},\Ome\}\cup\{0^{\dagger i}: 0<i\leq N\}$.

 \item\label{df:notationsystem.2}

If $\alp=_{NF}\alp_{m}+\cdots+\alp_{0}\, (m>0)$ with $\{\alp_{i}:i\leq m\}\subset OT(\mI_{N})$,
then
$\alp\in OT(\mI_{N})$, and $m(\alp)=\emptyset$.

Let $\alp=_{NF}\vphi\bet\gam<\veps(\mI_{N})$ with $\{\bet,\gam\}\subset OT(\mI_{N})$.
Then
$\alp\in OT(\mI_{N})$ and $m(\alp)=\emptyset$.

In each case
$K_{X}(\alp)=K_{X}(k(\alp))$.

 \item\label{df:notationsystem.3}
Let $\alp=\psi_{\Ome}(a)$ with $a\in  OT(\mI_{N})$
and $K_{\alp}(a)<a$.
Then
$\alp\in OT(\mI_{N})$.

Let
$m(\alp)=\emptyset$.
$K_{X}(\alp)=\emptyset$ if $\alp\in X$.
$K_{X}(\alp)=\{a\}\cup K_{X}(a)$ if $\alp\not\in X$.

  \item\label{df:notationsystem.4}
 Let $\alp=\psi_{\mI_{N}}(a)$ with $a\in  OT(\mI_{N})$
such that $K_{\alp}(a)<a$.
Then
$\alp\in LSt_{N}$
and $\alp^{\dagger i}\in SSt_{i}$ for $0<i\leq N$.
For $\bet\in\{\alp,\alp^{\dagger i}\}$,
$m(\bet)=\emptyset$.
Also 
$K_{X}(\alp^{\dagger i})=\emptyset$ if $\alp^{\dagger i}\in X$.
$K_{X}(\alp^{\dagger i})=K_{X}(\alp)$ if $\alp^{\dagger i}\not\in X$.
$K_{X}(\alp)=\emptyset$ if $\alp\in X$.
$K_{X}(\alp)=\{a\}\cup K_{X}(a)$ if $\alp\not\in X$.

\item\label{df:notationsystem.5}
Let $\mT\in LSt_{k}\cup\{0\}$ and $\mS=\mT^{\dagger\vec{i}}\in SSt_{i+1}$
for a non-empty and non-increasing sequence of numbers
$\vec{i}=(i_{0}\geq i_{1}\geq\cdots\geq i_{n})$ such that
$i_{0}\leq k$ and $i_{n}=i+1$, cf.\,Proposition \ref{prp:Stclass}.
Let $\alp=\psi_{\mS}^{f}(a)$, 
where $\{a,\mS\}\subset OT(\mI_{N})$, and if $f\neq\emptyset$, then 
there are $\{d,\xi\}\subset OT(\mI_{N})$ such that
${\rm supp}(f)=\{d\}$, $0<f(d)=\xi<(\mI_{N})^{2}$, $d<\mI_{N}$.
If
$K_{\mS}(\{\mS,a\}\cup \mathrm{fld}(f))<a$ for
$\mathrm{fld}(f)=\{d,\xi\}$, 

\beqn\label{eq:notationsystem.5}
\mathrm{fld}(f)\subset C_{0}(SC(a))
\eeqn
and
\beqn\label{eq:notationsystem.55}
\forall b\left[
\mT=\psi_{\mI_{N}}(b) \Rarw b<a
\right]
\eeqn

Then
$\alp\in LSt_{i}$ and
$\alp^{\dagger j}\in SSt_{j}$ for $0<j\leq i$.

Let
$a=\mathtt{ p}_{0}(\alp)$,
$m(\alp)=f$.
$K_{X}(\alp)=\emptyset$ if $\alp\in X$.
$K_{X}(\alp)=\{a\}\cup K_{X}(\{a,\mS\}\cup \mathrm{fld}(f))$ if
$\alp\not\in X$.

$m(\alp^{\dagger j})=\emptyset$.
$K_{X}(\alp^{\dagger j})=\emptyset$ if $\alp^{\dagger j}\in X$.
$K_{X}(\alp^{\dagger j})=K_{X}(\alp)$ if $\alp^{\dagger j}\not\in X$.

 \item\label{df:notationsystem.6}

Let
$\{\pi,a,d\}\subset  OT(\mI_{N})$ with
$\pi\prec\mS\in SSt_{i+1}$,
$m(\pi)=f$,
 $d<c\in \supp(f)$,
and $(d,c)\cap \supp(f)=\emptyset$.

When $g\neq\emptyset$, let
$g$
be an irreducible finite function 
 such that 
$\mathrm{fld}(g)\subset OT(\mI_{N})$,
$g_{d}\leq f_{d}$, $(d,c)\cap \supp(g)=\emptyset$,
$g(d)<f(d)+\tilde{\theta}_{c-d}(f(c); \mI_{N})\cdot\omega$, 
and $g<_{\mI_{N}}^{c}f(c)$.

Then
$\alp=\psi_{\pi}^{g}(a)\in LSt_{i}$ and
$\alp^{\dagger j}\in SSt_{j}$ for $0<j\leq i$
 if 
$K_{\pi}(k(\alp))<a$, 
and 

\beqn\label{eq:notationsystem.6}
\mathrm{fld}(g)\cup\{\mathtt{p}_{0}(\alp)\}\subset M_{\alp}
\eeqn

Let $m(\alp)=g$. 
$K_{X}(\alp)=\emptyset$ if $\alp\in X$.
$K_{X}(\alp)=\{a\}\cup K_{X}(k(\alp))$ if $\alp\not\in X$.

$m(\alp^{\dagger j})=\emptyset$.
$K_{X}(\alp^{\dagger j})=\emptyset$ if $\alp^{\dagger j}\in X$.
$K_{X}(\alp^{\dagger j})=K_{X}(\alp)$ if $\alp^{\dagger j}\not\in X$.

\item\label{df:notationsystem.7}
Let $\mS\in SSt_{i}$ and $0<k\leq i$.
Then $\mS^{\dagger k}\in SSt_{k}$.

$m(\mS^{\dagger k})=\emptyset$.
$K_{X}(\mS^{\dagger k})=\emptyset$ if $\mS^{\dagger k}\in X$.
$K_{X}(\mS^{\dagger k})=K_{X}(\mS)$ if $\mS^{\dagger k}\not\in X$.

\item\label{df:notationsystem.8}
Let
$SSt_{i}^{M}=SSt_{i}\cup\{\alp[\rho/\mS]: \rho\prec\mS\in SSt^{M}, \alp\in M_{\rho}\cap SSt_{i}^{M}\}$
and
$SSt^{M}=\bigcup_{0<i\leq N} SSt_{i}^{M}$.
Also let
$LSt_{i}^{M}=LSt_{i}\cup\{\alp[\rho/\mS]: \rho\prec\mS\in SSt^{M}, \alp\in M_{\rho}\cap LSt_{i}^{M}\}$
and
$LSt^{M}=\bigcup_{0<i\leq N} LSt_{i}^{M}$.

Let $\rho\prec\mS\in SSt_{i+1}^{M}$ and $\vec{i}=(i_{0}\geq i_{1}\geq\cdots\geq i_{n})\,(n\geq 0)$ with $0<i_{n}\leq i_{0}\leq i+1$.
Then
$(\mS^{\dagger \vec{i}}[\rho/\mS]) \in SSt_{i_{n}}^{M}\subset OT(\mI_{N})$,
where a term $\mS^{\dagger \vec{i}}[\rho/\mS]$ is built from terms $\mS^{\dagger \vec{i}}$, $\rho$ and $\mS$ by the constructor $*_{0}[*_{1}/*_{2}]$.

 \item\label{df:notationsystem.9}
 Let $\alp=\bet[\rho/\mS]$ with
 $\mS<\bet\in M_{\rho}$, 
 $\rho\prec\mS$, and $\mS\in SSt^{M}$.
 Then $\alp\in OT(\mI_{N})\setm St$.

 \eenu

}
\edf

\bprp\label{prp:Stclass}
Let $\alp\in OT(\mI_{N})$.
\benu

\item\label{prp:Stclass.1}
$\alp\in LSt_{N}$ iff $\alp=\psi_{\mI_{N}}(a)$ for an $a$.
For $0<i<N$, $\alp\in LSt_{i}\cap\Psi$ iff there exists an $\mS\in SSt_{i+1}$ such that $\alp\prec\mS$.

\item\label{prp:Stclass.2}
$\bet\in SSt_{k}$ iff there exists 
an $\alp\in\{0\}\cup (LSt_{i}\cap\Psi)$ for an $k\leq i\leq N$ and
a non-empty and non-increasing sequence $\vec{i}=(i_{0}\geq i_{1}\geq\cdots\geq i_{n})$
of numbers such that $k=i_{n}>0$, $\alp\in LSt_{i}\Rarw i_{0}\leq i$ and
$\bet=\alp^{\dagger \vec{i}}$.

\item\label{prp:Stclass.3}
Let $\psi_{\mS}^{f}(a)\in OT(\mI_{N})$ with $\mS\in SSt$.
Let $\{(\mT_{m},\mS_{m},\vec{i}_{m})\}_{m\leq n}$ be a sequence such that
$\mT_{m}\in LSt\cap\Psi$, $\mS_{m}\in SSt$ and $\vec{i}_{m}$ sequences of numbers such that
$\mT_{0}=\psi_{\mI_{N}}(b)$,
$\mS_{m}=\mT_{m}^{\dagger\vec{i}_{m}}$ and $\mT_{m+1}\prec\mS_{m}\,(m<n)$, and
$\mS=\mS_{n}$.
Then $b<a$ holds.

\item\label{prp:Stclass.4}
$\alp\in SSt^{M}$ iff there exists a $\rho$ and an $\vec{i}$ such that
$\alp\in\{\rho^{\dagger\vec{i}},\mS^{\dagger\vec{i}}[\rho/\mS]\}$.
\eenu
\eprp
\bprf
\ref{prp:Stclass}.\ref{prp:Stclass.1} and \ref{prp:Stclass}.\ref{prp:Stclass.2}.
We see these
from Definitions \ref{df:notationsystem}.\ref{df:notationsystem.1},
 \ref{df:notationsystem}.\ref{df:notationsystem.4},
\ref{df:notationsystem}.\ref{df:notationsystem.5}, 
\ref{df:notationsystem}.\ref{df:notationsystem.6} and \ref{df:notationsystem}.\ref{df:notationsystem.7}.
\\
\ref{prp:Stclass}.\ref{prp:Stclass.3}.
Let $\mT_{m}=\psi_{\sig_{m}}^{g_{m}}(b_{m})$ and
$\mT_{m}\preceq\psi_{\mS_{m-1}}^{f_{m-1}}(a_{m-1})$ for $\psi_{\mS_{m-1}}^{f_{m-1}}=\psi_{\mI_{N}}$
and $a_{-1}=b$.
In general, if $\sig=\psi_{\tau}^{f}(c)\in C_{b}(\psi_{\sig}^{g}(b))$ with $\psi_{\sig}^{g}(b)<\sig$, then
$c<b$.
Hence 
$a_{m-1}<a_{m}$.
On the other we obtain
$b=a_{-1}<a_{0}$ by (\ref{eq:notationsystem.55})
 in Definition \ref{df:notationsystem}.\ref{df:notationsystem.5}.
Therefore $b=a_{-1}<a_{n}=a$.
\eprf
\\

Let $\mT\in St_{i}\cup\{0\}$ and $i\geq j$. 
If $\mT\in LSt_{i}\cup\{0\}$, then $\mT^{\dagger j}\in SSt_{j}\subset OT(\mI_{N})$.
Otherwise there is an $\mS\in LSt_{k}\cup\{0\}$, a non-increasing sequence
$\vec{i}=(i_{0}\geq\cdots\geq i_{n})$ such that $i_{n}=i$ and
$\mT=\mS^{\dagger\vec{i}}$, where $k\geq i_{0}$ when $\mS\neq 0$.
Then let $\mT^{\dagger j}:=\mS^{\dagger\vec{j}}$ for $\vec{j}=\vec{i}*(j)=(i_{0},\ldots,i_{n},j)$.

From Definition \ref{df:notationsystem} and Proposition \ref{prp:Stclass} we see that
for each $\alpha\in LSt_{i}$ with $i<N$, there is an $\mS\in SSt_{i+1}$ such that $\alpha\prec\mS$.
Also for each $\beta\in SSt_{i}$ there is a $\mT\in St_{k}$ with $k\geq i$ such that
$\beta=\mT^{\dagger i}$.

Sets $C_{\gam}(\del)$ are defined for $\{\gam,\del\}\subset OT(\mI_{N})$ in such a way that
$\alp\in C_{\gam}(\del)$ iff $K_{\del}(\alp)<\gam$
for
$\alp,\gam,\del\in OT(\mI_{N})$.
In particular $OT(\mI_{N})=C_{\veps(\mI_{N})}(0)$, and
$C_{\gam}(X)$ is closed under Mostowski collapsing $\alp\mapsto\alp[\rho/\mS]$
if $\gam\geq\mI_{N}$, and differs from sets defined in Definition \ref{df:Cpsiregularsm}.

We define terms $\alp[\rho/\mS]$, sets 
$K_{X}(\alp[\rho/\mS])$ and a relation $\bet<\gam$ on $OT(\mI_{N})$
recursively as follows.

\bdf\label{df:notationsystem.2M} 
{\rm
(Definitions of $\alp[\rho/\mS]$ and $K_{X}(\alp[\rho/\mS])$)
\\
Let $\rho\prec\mS\in SSt_{i+1}^{M}$.
We define a term $\alp[\rho/\mS]\in OT(\mI_{N})$ for $\alp\in M_{\rho}$
in such a way that
$\alp[\rho/\mS]=\alp$ iff $\alp<\rho$.
Moreover
$\alp[\rho/\mS]\in St$ iff either $\alp[\rho/\mS]=\alp\in St$ or $\alp[\rho/\mS]=\rho\in SSt$.

Also 
$K_{X}(\alp[\rho/\mS])$
is defined
recursively as follows.
The map $\alp\mapsto\alp[\rho/\mS]$ commutes with $\psi$, 
$\vphi$, $\mI_{N}[\cdot]$, and $+$.
$K_{X}(\alp[\rho/\mS])=\emptyset$ if $\alp[\rho/\mS]\in X$.
\benu
\item
$\alp[\rho/\mS]:=\alp$ when $\alp<\mS$.

In what follows assume $\alp\geq\mS$, $\alp[\rho/\mS]\geq\rho$ and 
$\alp[\rho/\mS]\not\in X$.

\item
$(\mS)[\rho/\mS]:=\rho$ and $(\mathbb{I}_{N})[\rho/\mS]:=\mathbb{I}_{N}[\rho]$.

For $\vec{i}\leq i+1$,
$(\mS^{\dagger \vec{i}})[\rho/\mS]:=(\mS^{\dagger \vec{i}}[\rho/\mS])\in SSt_{i+1}^{M}$, cf.\,Definition \ref{df:notationsystem}.\ref{df:notationsystem.8}.
Here
$\mS^{\dagger \vec{i}}[\rho/\mS]\neq\rho^{\dagger \vec{i}}$.

$K_{X}(\alp[\rho/\mS])=K_{X}(\rho)$ if 
$\alp[\rho/\mS]\in\{\mI_{N}[\rho],\mS^{\dagger \vec{i}}[\rho/\mS]\}$.

\item
Let $\alp=\psi_{\mI_{N}}(a)$. Then
$\alp[\rho/\mS]=
\psi_{\mI_{N}[\rho]}(a[\rho/\mS])$.

$K_{X}(\alp[\rho/\mS])=K_{X}(\{\rho, a[\rho/\mS]\})\cup \{a[\rho/\mS]\}$.

\item
Let $\alp=\psi_{\kap}^{f}(a)$. Then
$\alp[\rho/\mS]=
\psi_{\kap[\rho/\mS]}^{f[\rho/\mS]}(a[\rho/\mS])$, where $(f[\rho/\mS]):\mI_{N}[\rho]\to\Gamma(\mI_{N}[\rho])$,
${\rm supp}(f[\rho/\mS])=({\rm supp}(f))[\rho/\mS]=\{c[\rho/\mS]: c\in{\rm supp}(f)\}$ and
$(f[\rho/\mS])(c[\rho/\mS])=(f(c))[\rho/\mS]$
for $f:\mI_{N}[\rho]\to\Gamma(\mI_{N}[\rho])$ and $c\in{\rm supp}(f)$.

$K_{X}(\alp[\rho/\mS])=K_{X}(\{\kap[\rho/\mS], a[\rho/\mS]\}\cup SC(f[\rho/\mS]))\cup\{a[\rho/\mS]\}$.

$M_{\alp[\rho/\mS]}=C_{b[\rho/\mS]}(\alp[\rho/\mS])$ for
$b=\mathtt{p}_{0}(\alp)$ and $b[\rho/\mS]=\mathtt{p}_{0}(\alp[\rho/\mS])$.

\item
Let $\alp=\mI_{N}[\tau]\neq\mI_{N}$. Then
$\alp[\rho/\mS]=\mI_{N}[\tau[\rho/\mS]]$, where
$\mI_{N}[\tau]\in M_{\rho}$ iff $\tau\in M_{\rho}$.
$K_{X}(\alp[\rho/\mS])=K_{X}(\tau[\rho/\mS])$.

\item
Let $\alp=\psi_{\mI_{N}[\tau]}(a)$ for $\mI_{N}[\tau]\neq\mI_{N}$. Then
$\alp[\rho/\mS]=
\psi_{\mI_{N}[\tau[\rho/\mS]]}(a[\rho/\mS])$.

$K_{X}(\alp[\rho/\mS])=K_{X}(\{\tau[\rho/\mS], a[\rho/\mS]\})\cup \{a[\rho/\mS]\}$.

\item 
Let $\alp=\tau^{\dagger \vec{j}}$ with $\mS<\tau\in LSt^{M}$. Then
$\alp[\rho/\mS]=(\tau[\rho/\mS])^{\dagger \vec{j}}$,
where 
$\tau^{\dagger \vec{j}}\in M_{\rho}$
iff $\tau\in M_{\rho}$.
$K_{X}(\alp[\rho/\mS])=K_{X}(\tau[\rho/\mS])$.

\item
Let $\alp=\mT^{\dagger \vec{j}}[\tau/\mT]$, where $\tau\prec\mT\in SSt^{M}$.
Then $\alp[\rho/\mS]=\mT^{\dagger \vec{j}}_{1}[\tau_{1}/\mT_{1}]$, where
$\tau_{1}=\tau[\rho/\mS] \prec \mT_{1}=\mT[\rho/\mS]\in SSt^{M}$ and
$\mT^{\dagger \vec{j}}_{1}=(\mT_{1})^{\dagger \vec{j}}$.

$K_{X}(\alp[\rho/\mS])=K_{X}(\tau[\rho/\mS])$.

\item
Let $\alp=\vphi\bet\gam$. Then
$\alp[\rho/\mS]=\vphi(\bet[\rho/\mS])(\gam[\rho/\mS])$.

$K_{X}(\alp[\rho/\mS])=K_{X}(\bet[\rho/\mS], \gam[\rho/\mS])$.

\item
For $\alp=\alp_{m}+\cdots+\alp_{0}\, (m>0)$,
$\alp[\rho/\mS]=(\alp_{m}[\rho/\mS])+\cdots+(\alp_{0}[\rho/\mS])$.
$K_{X}(\alp[\rho/\mS])=\bigcup\{K_{X}(\alp_{i}[\rho/\mS])): i\leq m\}$.

\eenu

}
\edf

A relation $\alp<\bet$ for $\alp,\bet\in OT(\mI_{N})$ is defined according to
Lemmas \ref{lem:welldefinedness.2} and \ref{lem:limitcollapse.1},
Propositions \ref{prp:comparisonrud}, \ref{prp:welldefinedness.suc}, and
\ref{prp:psicomparison},
and
Corollary \ref{cor:stepdown}, provided that
$\alp\in C_{\gam}(\del)$ is replaced by $K_{\del}(\alp)<\gam$.
The relation enjoys
$\psi_{\kap}^{f}(a)<\kap$ according to Lemma \ref{lem:limitcollapse.1} and Corollary \ref{cor:stepdown}.
Moreover we obtain
$\mS^{\dagger i}<\psi_{\mS^{\dagger (i+1)}}^{g_{0}}(b_{0})<\mS^{\dagger (i+1)}$ for $i+1\leq N$,
and
$LSt_{N}\ni \tau_{0}=\psi_{\mI_{N}}(c_{0})<\psi_{\tau_{0}^{\dagger}}^{h_{0}}(d_{0})<\tau_{0}^{\dagger}<\mI_{N}$
by Proposition \ref{prp:welldefinedness.suc} and Lemma \ref{lem:welldefinedness.2}.
Hence if $\mS<\psi_{\mI_{N}}(c_{0})$, then
$\mS<\mS^{\dagger i}<\psi_{\mS^{\dagger(i+1)}}^{g_{0}}(b_{0})<\mS^{\dagger (i+1)}<
\tau_{0}=\psi_{\mI_{N}}(c_{0})<\psi_{\tau_{0}^{\dagger}}^{h_{0}}(d_{0})<\tau_{0}^{\dagger}<\mI_{N}$.
The Mostowski collapsing $\cdot[\rho/\mS]$ maps these inequalities isomorphically to
$\rho<\mS^{\dagger i}[\rho/\mS]<\psi_{\mS^{\dagger (i+1)}[\rho/\mS]}^{g}(b)<\mS^{\dagger (i+1)}[\rho/\mS]<\tau=\psi_{\mI_{N}[\rho]}(c)<
\psi_{\tau^{\dagger}}^{h}(d)<\tau^{\dagger}<\mI_{N}[\rho]<\rho^{\dagger 0}$, where
$b=b_{0}[\rho/\mS]$, etc.

\begin{definition}\label{df:precR}
{\rm
For terms $\pi,\kappa\in OT(\mI_{N})$,
a relation $\pi\prec^{R} \kap$ is defined recursively as follows.
 \benu
 \item\label{df:precR.1}
Let $\pi\prec\kap\preceq \mS\in SSt_{i+1}^{M}$, and
$\vec{i}\leq i+1$.
Then
each of
$\pi\prec^{R}\kap$, 
$\mS^{\dagger\vec{i}}[\pi/\mS]\prec^{R}\kap$ and $\mI_{N}[\pi]\prec^{R}\kap$
holds.
Moreover $\pi^{\dagger \vec{i}}\prec^{R}\kap$ holds provided that
$\pi^{\dagger \vec{i}}\not\in SSt$.

 \item\label{df:precR.2}
 $\tau\prec^{R}\pi\prec^{R}\kap\Rarw \tau\prec^{R}\kap$.
 \eenu
Let $\pi\preceq^{R}\kappa:\Lrarw \pi\prec^{R}\kappa \lor \pi=\kappa$.
For $\mS\in SSt$, let
$$
L(\mS):=\{\alp\in OT(\mI_{N}): \alp\prec^{R}\mS\}.
$$
}
\end{definition}

Note that $L(\mS)\cap SSt=\emptyset$, and
$SSt\ni\rho^{\vec{i}}\not\prec^{R}\mS$ for $LSt_{i}\ni\rho\prec\mS\in SSt_{i+1}$
and $\vec{i}\leq i$.
For $\mS\in SSt$ and a strongly critical number $\alpha$,
$\alpha\in L(\mS)$ iff $\alpha=\alpha_{0}[\rho/\mS]$ for a $\rho\prec\mS$ and an $\mS\leq\alpha_{0}\in M_{\rho}$.
For each strongly critical number $\Ome<\alpha<\mI_{N}$,
there exists a (unique) $\mS\in SSt$ such that $\alp\prec^{R}\mS$ iff
$\alpha\not\in St$.

Let $\bet\prec^{R}\mT$ and $\alp\prec^{R}\mS$ be such that $\bet<\alp$ and $\mT\neq\mS$.
Then either $L(\mT)<\alpha$ or $\beta<L(\mS)$ holds,
where $L(\mT)<\alpha\Lrarw\forall\beta\in L(\mT)(\beta<\alpha)$ and
$\beta<L(\mS)\Lrarw \forall\alpha\in L(\mS)(\beta<\alpha)$.

\bdf\label{df:lessthanpsi}
{\rm
Let $\bet,\alp\in OT(\mI_{N})\cap\mI_{N}$ be strongly critical numbers.
$\bet<\alp$ iff one of the following cases holds:

\benu

\item\label{df:lessthanpsi.1}
$\bet=\Ome$ and $\Ome\neq\alp\neq\psi_{\Ome}(a)$ for any $a$.

\item\label{df:lessthanpsi.21}
%SSt vs SSt
There are $\mS\in St_{k_{0}}$, $\mT\in St_{i_{0}}$ and $i,k\geq 1$ such that
$k_{0}\geq k$, $i_{0}\geq i$,
$\beta=\mT^{\dagger i}$, $\alpha=\mS^{\dagger k}$ and
either $i<k\spand \mT<\mS^{\dagger k}$ or
$i\geq k\spand \mT^{\dagger i}\leq \mS$.

\item\label{df:lessthanpsi.22}
%LSt vs SSt k=N

 \benu
 \item
 $\alpha\prec\mI_{N}$, and
there is $\mT\in St_{i_{0}}$ with $i_{0}\geq i$ such that
%\in SSt_{N}$.
$\beta=\mT^{\dagger i}$ and
$\mT<\alpha$.

 \item
 $\beta\prec\mI_{N}$,  and there is $\mS\in St_{k_{0}}$ with $k_{0}\geq k$ such that
$\alpha=\mS^{\dagger k}$, and
$\beta\leq\mS$.

 \eenu
 
\item\label{df:lessthanpsi.23}
%LSt vs SSt $k<N$
 \benu
 \item
There are $\mS\in St_{k_{0}}$ with $k_{0}>k$ and $\mT\in St_{i_{0}}$ with $i_{0}\geq i$ such that
$\alpha\prec\mS^{\dagger(k+1)}$, 
%\in SSt_{k+1}$.
$\beta=\mT^{\dagger i}$ and
%\in SSt_{i}$, and
%$\alpha\in LSt_{k}$ with
 either $i\leq k \spand \mT<\alpha$ or
$i>k \spand \beta\leq \mS$.

 \item
 There are $\mT\in St_{i_{0}}$ with $i_{0}>i$ and $\mS\in St_{k_{0}}$ with $k_{0}\geq k$ such that
$\beta\prec\mT^{\dagger(i+1)}$, 
%\in SSt_{k+1}$.
$\alpha=\mS^{\dagger k}$, and
%$\alpha\in LSt_{k}$ with
 either $k\leq i \spand \beta\leq\mS$ or
$k>i \spand \mT<\alpha$.

 \eenu

\item\label{df:lessthanpsi.24}
%LSt vs LSt $\spi$ vs $\psi$
$\bet=\psi_{\pi}^{f}(b)$, $\alp=\psi_{\kappa}^{g}(a)$ for some $\pi,\kappa,f,g,b,a$
for which
one of the following holds:

\begin{enumerate}
\item\label{df:lessthanpsi.240}
$\pi\leq \alp$.

\item\label{df:lessthanpsi.241}
$b<a$, $\bet<\kappa$, and 
$K_{\alp}(\mathrm{fld}(f)\cup\{\pi,b\})<a$

\item\label{df:lessthanpsi.242}
$b>a$ and $b\leq K_{\bet}(\mathrm{fld}(g)\cup\{\kappa,a\})$.

\item\label{df:lessthanpsi.243}
$b=a$, $\kappa<\pi$, and $b\leq K_{\bet}(\kap)$.

\item\label{df:lessthanpsi.244}
$b=a$, $\pi=\kappa$, $K_{\alp}(\mathrm{fld}(f))<a$, and
$f<^{0}_{lx}g$.

\item\label{df:lessthanpsi.245}
$b=a$, $\pi=\kappa$, and
$b\leq K_{\bet}(\mathrm{fld}(g))$.

 \eenu

\item\label{df:lessthanpsi.5}
%fake vs St
 \benu
 \item\label{df:lessthanpsi.5a}
$\alpha\in St$, and there are $\mT\in SSt$, $\rho\prec\mT$ and
$\mT\leq\beta_{0}\in M_{\rho}$ such that
$\beta=\beta_{0}[\rho/\mT]$,
and $\rho<\alpha$.

\item\label{df:lessthanpsi.5b}
$\beta\in St$, and there are $\mS\in SSt$, $\tau\prec\mS$ and $\mS\leq\alpha_{0}\in M_{\tau}$
such that
$\alpha=\alpha_{0}[\tau/\mS]$ and $\beta\leq\tau$.

 \eenu

\item\label{df:lessthanpsi.6}
%teisei fake vs fake
There are $\mT\neq\mS\in SSt$ such that
$\bet\prec^{R}\mT$, $\alp\prec^{R}\mS$ and one of the following holds:
 \benu
 \item $\mT<\alpha$.
 \item $\beta<\mS<\mT$.
 \eenu

\item\label{df:lessthanpsi.7}
%fake vs fake
There are $\mS\in SSt$, $\tau\prec\mS$, $\mS\leq\alpha_{0}\in M_{\tau}$,
$\rho\prec\mS$, $\mS\leq\beta_{0}\in M_{\rho}$ such that
$\alpha=\alpha_{0}[\tau/\mS]$, $\beta=\beta_{0}[\rho/\mS]$ and
either $\rho<\tau$ or $\rho=\tau \spand \beta_{0}<\alpha_{0}$.

\eenu
}
\edf

\bprp\label{prp:M_bnd}
Let $\rho\prec\mT$ with $\mT\in SSt$, and $\beta\in M_{\rho}$.
Then $\beta[\rho/\mT]<\mS$ for any $\rho<\mS\in St$.
\eprp
\bprf
This follows from 
Definition \ref{df:lessthanpsi}.\ref{df:lessthanpsi.5a}.
\eprf

\blem\label{lem:compT}
$(OT(\mI_{N}),<)$ is a computable linear order.
Specifically
each of $\alpha<\beta$ and $\alpha=\beta$ is decidable for $\alpha,\beta\in OT(\mI_{N})$, 
and $\alpha\in OT(\mI_{N})$ is decidable for terms $\alpha$
over symbols $\{0,\Ome,\mI_{N}, +,\varphi, \psi\}$, $\{ {*}^{\dagger i}:0<i\leq N\}$,
$\mI_{N}[*]$ and $*_{0}[*_{1}/*_{2}]$.

In particular the order type of the initial segment $\{\alpha\in OT(\mI_{N}): \alpha<\Omega\}$
is less than $\omega_{1}^{CK}$ if it is well-founded.
\elem

In what follows by ordinals we mean ordinal terms in $OT(\mI_{N})$.
$\ell\alp$ denotes the length of ordinal terms $\alp$, which means the number of
occurrences of symbols in $\alp$.

\bprp\label{prp:LS}
If $\mS\in St_{i+1}=SSt_{i+1}\cup LSt_{i+1}$ and $\alp<\mS$,
then $\alp^{\dagger i}<\mS$.
\eprp
\bprf
This is seen from Proposition \ref{prp:Stclass} and Definition \ref{df:lessthanpsi}.
\eprf

\bprp\label{prp:Erho}
$\{\mS\}\cup \mathrm{fld}(m(\rho))\cup\{\mathtt{p}_{0}(\rho)\}\subset M_{\rho}$ for $\rho\in\Psi_{\mS}$.
\eprp
\bprf
If $\rho=\psi_{\mS}^{f}(a)$ with an $\mS\in SSt$, then we obtain
$f=m(\rho)$, $a=\mathtt{ p}_{0}(\rho)$, 
$\{\mS\}\cup \mathrm{fld}(f)\cup\{\mathtt{p}_{0}(\rho)\}\subset C_{a}(\alp)=M_{\rho}$ by
Definition \ref{df:notationsystem}.\ref{df:notationsystem.5}.
Otherwise $\{\mS\}\cup \mathrm{fld}(m(\rho))\cup\{\mathtt{p}_{0}(\rho)\}\subset M_{\rho}$ follows from
 (\ref{eq:notationsystem.6})
in Definition \ref{df:notationsystem}.\ref{df:notationsystem.6}.
\eprf

\bprp\label{prp:Mostowski_base}
Let $\rho\in\Psi_{\mS}$ and $\mS\leq\alpha\in M_{\rho}$.
If $\alpha[\rho/\mS]\in C_{a}(\beta)$, then $\rho\in C_{a}(\beta)$.
\eprp
\bprf
This is seen by induction on $\ell\alpha$ from Definition \ref{df:notationsystem.2M}.
\eprf
\\

An ordinal term $\sig\in OT(\mI_{N})$ is said to be \textit{regular} if
either $\sig\in\{\mI_{N}\}\cup\{\sig\in OT(\mI_{N}): \exi\rho(\sig=\mI_{N}[\rho])\}$ or
$\psi_{\sig}^{f}(a)$ is in $OT(\mI_{N})$ for some $f$ and $a$.
$Reg$
denotes the set of regular terms.
Then
$Reg=SSt^{M}\cup\{\mI_{N}[\rho]:\exi\mS\in SSt^{M}(\rho\prec\mS)\}\cup\{\Ome,\mI_{N}\}$.
We see that for each $\alp\in\Psi$, there exists a $\kap\in Reg_{0}:=(Reg\setm\Psi)$ such that
$\alp\prec\kap$.
Such a $\kap$ is either in $\{\Ome,\mI_{N}\}$ or 
one of the form $\mI_{N}[\rho]$, $\rho^{\dagger\vec{i}}$
 or $\mS^{\dagger\vec{i}}[\rho/\mS]$ with a non-empty $\vec{i}$.

\bprp\label{prp:jumpover}
Let $\psi_{\pi}^{f}(a)<\psi_{\kap}^{g}(b)<\pi<\kap$ and $\pi\preceq\rho$ and $\kap\preceq\tau$
with $\{\rho,\tau\}\subset Reg_{0}$.
Then $\rho=\tau$.
\eprp
\bprf
From Definition \ref{df:lessthanpsi} we see that the only possible case is
Definition \ref{df:lessthanpsi}(\ref{df:lessthanpsi.241}).
\eprf

\blem\label{lem:Mostowskicollaps}
For $\rho\prec\mS$ and $\mS\in SSt$,
$\{\alpha[\rho/\mS]:\alpha\in M_{\rho}\}$ is a transitive collapse of $M_{\rho}$ in the following sense.
Let $\{\alpha,\beta,\gamma\}\subset M_{\rho}$.

\begin{enumerate}

\item\label{lem:Mostowskicollaps.1}
$\beta<\alpha\Leftrightarrow\beta[\rho/\mS]<\alpha[\rho/\mS]$.

\item\label{lem:Mostowskicollaps.15}
$\bet\prec^{R}\alp \Lrarw \bet[\rho/\mS]\prec^{R}\alp[\rho/\mS]$.

\item\label{lem:Mostowskicollaps.2}
$\mathbb{S}<\gam \Rightarrow 
\left(
K_{\gam}(\bet)<\alp
\Lrarw
K_{\gam[\rho/\mS]}(\bet[\rho/\mS])<\alp[\rho/\mS]
\right)$.

\item\label{lem:Mostowskicollaps.4}
$OT(\mI_{N})\cap\alpha[\rho/\mS]=\{\gamma[\rho/\mS]:\gamma\in M_{\rho}\cap\alpha\}$.

\eenu

\elem
\bprf
We show Lemmas \ref{lem:Mostowskicollaps}.\ref{lem:Mostowskicollaps.1}-
\ref{lem:Mostowskicollaps}.\ref{lem:Mostowskicollaps.2}
 simultaneously 
by induction on the sum $2^{\ell\alpha}+2^{\ell\beta}$ for
$\alpha,\beta\in M_{\rho}$.
We see easily that 
$\mathbb{S}>\Gamma(\mathbb{I}_{N}[\rho])>\alpha[\rho/\mS]>\rho$ when 
$\alpha>\mathbb{S}$.
Also $\alpha[\rho/\mS]\leq\alpha$.
\\
\ref{lem:Mostowskicollaps}.\ref{lem:Mostowskicollaps.15} and \ref{lem:Mostowskicollaps}.\ref{lem:Mostowskicollaps.2} are seen from IH.
\\
 \ref{lem:Mostowskicollaps}.\ref{lem:Mostowskicollaps.1}.
 Let $k(\psi^{g}_{\kappa}(a))=\mathrm{fld}(g)\cup\{\kap,a\}$.
Let $\mathbb{S}<\beta=\psi^{f}_{\pi}(b)<\psi^{g}_{\kappa}(a)=\alpha$ with
$k(\bet,\alp)\subset M_{\rho}$.
From IH with Definition \ref{df:lessthanpsi}
we see that 
$\beta[\rho/\mS]=\psi_{\pi[\rho/\mS]}^{f[\rho/\mS]}(b[\rho/\mS])<\psi_{\kap[\rho/\mS]}^{g[\rho/\mS]}(a[\rho/\mS])=\alpha[\rho/\mS]$.
Other cases are seen from IH.
\\
\ref{lem:Mostowskicollaps}.\ref{lem:Mostowskicollaps.2}.
Suppose $K_{\gam}(\bet)<\alp$
for $\mathbb{S}<\gam$.
Then $K_{\gam[\rho/\mS]}(\bet[\rho/\mS])<\alp[\rho/\mS]$
is seen from 
IH and Lemma \ref{lem:Mostowskicollaps}.\ref{lem:Mostowskicollaps.1}
 using the fact $\gamma[\rho/\mS]>\rho$.
\\

\noindent
\ref{lem:Mostowskicollaps}.\ref{lem:Mostowskicollaps.4}.
Let
$\beta\in OT(\mI_{N})\cap\alpha[\rho/\mS]$ for $\alpha\in M_{\rho}$.
We show by induction on $\ell\beta$ that there exists a 
$\gamma\in M_{\rho}$
such that $\beta=\gamma[\rho/\mS]$.
If $\beta<\rho$, then $\beta[\rho/\mS]=\beta$．
Also $\rho=\mathbb{S}[\rho/\mS]$ and $\mathbb{I}_{N}[\rho]=(\mathbb{I}_{N})[\rho/\mathbb{S}]$.
Let $\Gam(\mathbb{I}_{N}[\rho])>\alp[\rho/\mS]>\beta>\rho$.
We may assume $\mathbb{I}_{N}[\rho]>\beta>\rho$ by IH.

If $\bet=\mathbb{I}_{N}[\tau]$, then $\mathbb{I}_{N}[\tau]>\tau$.
Pick a $\kap\in M_{\rho}$ such that $\kap[\rho/\mS]=\tau$.
Then $\bet=(\mI_{N}[\kap])[\rho/\mS]$.

If $\bet=\tau^{\dagger \vec{i}}$, then $\tau^{\dagger \vec{i}}>\tau$.
Pick a $\kap\in M_{\rho}$ such that $\kap[\rho/\mS]=\tau$.
Then $\bet=(\kap^{\dagger \vec{i}})[\rho/\mS]$.

If $\bet=\mT_{1}^{\dagger\vec{j}}[\tau_{1}/\mT_{1}]$, then 
$\mT_{1}^{\dagger\vec{j}}[\tau_{1}/\mT_{1}]>\tau_{1}$.
Pick a $\tau\in M_{\rho}$ such that $\tau[\rho/\mS]=\tau_{1}$.
Then for $\tau\prec \mT\in SSt^{M}$, we obtain
$\bet=(\mT^{\dagger\vec{j}}[\tau/\mT])[\rho/\mS]$.

Finally let $\beta=\psi^{f}_{\pi}(b)$ with 
 $k(\beta)\subset C_{b}(\beta)$,
 $b<\Gam(\mI_{N}[\rho])$ and 
$f:\Lambda\to\Gamma(\Lambda)$ for 
$\pi\preceq\sig^{\dagger \vec{k}}$ with a $\vec{k}\neq\emptyset$.
We have $\bet\prec\sig^{\dagger \vec{k}}$, 
$\rho<\bet<\mI_{N}[\rho]$, and $\rho\prec\mS$.
By
Definition \ref{df:lessthanpsi}
we obtain $\sig\neq\mS$.
Suppose $\bet<\mS<\sig^{\dagger \vec{k}}$. Then $\alp<\rho$ by
Definition \ref{df:lessthanpsi}.
Hence we may assume $\sig^{\dagger \vec{k}}<\mS$.
Then we obtain 
$\rho<\sig^{\dagger k}<\mI_{N}[\rho]$.
Hence $\sig\prec^{R}\mI_{N}[\rho]$ or $\sig\prec^{R}\mS^{\dagger \vec{i}}[\rho/\mS]$ for an $\vec{i}$.
By IH with $\pi\leq\sig^{\dagger \vec{k}}$ there are 
$\{c,\kap,\lam\}\subset M_{\rho}$ and
$g:\lam\to\Gamma(\lam)$
such that $c[\rho/\mS]=b$, $\kap[\rho/\mS]=\pi$,
$\lam[\rho/\mS]=\Lam$,
$\mathrm{fld}(g)\subset M_{\rho}$,
$g[\rho/\mS]=f$ in the sense that
$({\rm supp}(g))[\rho/\mS]={\rm supp}(f)$ and
$(g(d))[\rho/\mS]=f(d[\rho/\mS])$ for every $d\in{\rm supp}(g)$.
Let $\gam=\psi_{\kap}^{g}(c)\in M_{\rho}$.
Then 
$\gam[\rho/\mS]=\psi_{\pi}^{f}(b)=\bet$ 
and $k(\gam)\subset C_{c}(\gam)$.

Other cases are seen from IH.
\eprf

\subsection{Sets $M_{\rho}$, trails and stepping-down}

In this subsection some facts on sets $M_{\rho}$, ordinal terms and finite functions are established.
These facts are needed in sections \ref{sec:operatorcont} and \ref{sec:proofonestep}.

\blem\label{prp:Hclosed}
\benu
\item\label{prp:Hclosed.1}
Let $\alp=\psi_{\Ome}(a)$ with $a\in C_{a}(\alp)$. 
Then $ C_{a}(\alp)\cap\Ome\subset\alp$.

\item\label{prp:Hclosed.2}
Let $\alp=\psi_{\mI_{N}}(a)$ with $a\in C_{a}(\alp)$. 
Then $ C_{a}(\alp)\cap\mI_{N}\subset\alp$.

\item\label{prp:Hclosed.3}
Let $\mS\in SSt$, and
$\alp=\psi_{\kap}^{f}(a)<\kap$ with 
$\kap\preceq\mS$ and 
$\{\kap,a\}\cup \mathrm{fld}(f)\subset C_{a}(\alp)$.
Then $ C_{a}(\alp)\cap\kap\subset\alp$.
\eenu
\elem
\bprf
We see $\bet\in C_{a}(\alp)\cap\Ome \Rarw \bet<\alp=\psi_{\Ome}(a)$
by induction on the lengths $\ell\bet$ of $\bet$.
Lemmas
\ref{prp:Hclosed}.\ref{prp:Hclosed.2} and \ref{prp:Hclosed}.\ref{prp:Hclosed.3}
are seen similarly using the fact
$\rho<\alp\Rarw \mI_{N}[\rho]<\alp$ for $\alp\in\{\psi_{\mI_{N}}(a),\psi_{\kap}^{f}(a)\}$.
\eprf

\bprp\label{prp:EmS}
Let $\mS\in SSt$, and
$\rho=\psi_{\kap}^{f}(a)<\kap$ with 
$\kap\preceq\mS$ and
$ C_{\gam}(\kap)\cap\mS\subset\kap$ for $\gam\leq a$.
Then $ C_{\gam}(\rho)\cap\mS\subset\rho$.
\eprp
\bprf
If $\kappa=\mathbb{S}$, then 
$ C_{\gamma}(\rho)\cap\mathbb{S}\subset C_{a}(\rho)\cap\mathbb{S}\subset\rho$
by $\gamma\leq a$ and Lemma \ref{prp:Hclosed}.\ref{prp:Hclosed.3}.
Let $\kappa=\psi_{\pi}^{g}(b)<\mathbb{S}$. 
We have $\kappa\in C_{a}(\rho)$ by (\ref{eq:Psivec}),
and hence $b<a$ by $\mS>\kappa>\rho$.
We obtain
$ C_{\gamma}(\rho)\cap\mathbb{S}\subset C_{\gamma}(\kappa)\cap\mathbb{S}\subset\kappa$.
$\gamma\leq a$ with Lemma \ref{prp:Hclosed}.\ref{prp:Hclosed.3}
yields
$ C_{\gamma}(\rho)\cap\mathbb{S}\subset C_{\gamma}(\rho)\cap\kappa\subset
 C_{a}(\rho)\cap\kappa\subset\rho$.
\eprf

\begin{definition}\label{df:divide}
{\rm (Mostowski uncollapsing)\\
Let $\alpha$ be an ordinal term and $\rho\prec\mS$ with
$\mS\in SSt$.
If there exists a $\beta\in M_{\rho}$ such that
$\alpha=\beta[\rho/\mS]$, then $\alpha[\rho/\mS]^{-1}:=\beta$. 
Otherwise $\alpha[\rho/\mS]^{-1}:=\alpha$.
Let 
 $X[\rho/\mS]^{-1}:=\{\alpha[\rho/\mS]^{-1}:\alpha\in X\}$ for a set $X$ of ordinal terms.

}
\end{definition}
We see that
ordinal terms $\rho$ and $\beta\in M_{\rho}$
 with $\rho\leq\alpha=\beta[\rho/\mS]< \Gam(\mathbb{I}_{N}[\rho])$
are uniquely determined from $\alp$, when
such $\beta$ and $\rho$ exist.

Following Buchholz\cite{Buchholz}, let us introduce operators $\calh_{\gam}$ on ordinals for ordinal analysis in sections \ref{sec:operatorcont} and \ref{sec:proofonestep}.
For a finite set $\Tht$ of ordinals, let $C_{\alp}(\Tht)$ denote the Skolem hull of $\{0,\Ome,\mI_{N}\}\cup\Tht$ 
under the functions $+, \vphi$, $\alp\mapsto\alp^{\dagger i}$ and
$\bet\mapsto\psi_{\pi}^{f}(\bet)\,(\bet<\alp)$ as in Definition \ref{df:Cpsiregularsm}.\ref{df:Cpsiregularsm.1}.

The set $C_{\alp}(\Tht)$ is closed under these functions in the sense that, e.g.,
if $\{\bet,\gam\}\subset C_{\alp}(\Tht)$, then $\bet+\gam\in C_{\alp}(\Tht)$.
However the set $C_{\alp}(\Tht)$ is not closed under these functions in a reverse direction.
This means, e.g., $\{\bet,\gam\}\not\subset C_{\alp}(\Tht)$ when $\bet+\gam\in C_{\alp}(\Tht)$.
For example $\Ome+\psi_{\Ome}(1)\in C_{0}(\{\Ome+\psi_{\Ome}(1)\})$, but
$\psi_{\Ome}(1)\not\in C_{0}(\{\Ome+\psi_{\Ome}(1)\})$.
It is convenient for us to have sets of ordinals closed under some functions in both directions.
This motivates the sets $\calh_{\gam}[\Tht]$ in the following Definition \ref{df:Buchholz4.3}.

We need sets of ordinals closed under $\mS^{\dagger}\mapsto\mS$ in subsection \ref{subsec:elimrfl}.
The set $\calh_{\gam}[\Tht]$ is better than $C_{\gam}(\Tht)$
unless we assume that the set $\Tht$ is closed under $\mS^{\dagger}\mapsto\mS$.

\bdf\label{df:Buchholz4.3}
{\rm
\benu
\item\label{df:Buchholz4.3.1}
We say that a set $\Tht$ of ordinals is \textit{$C$-bounded} if there are ordinals $\alp$ and $\bet$ such that $\Tht\subset C_{\alp}(\bet)$.
\item\label{df:Buchholz4.3.2}
For $C$-bounded sets $\Tht$ of ordinals, let
\[
\calh_{\gam}[\Tht]:=\bigcap\{C_{\alp}(\bet): \gam<\alp,\, \Tht\subset C_{\alp}(\bet)\}
\]
where
$\gam\in C_{\alp}(\bet) \Lrarw K_{\bet}(\gam)<\alp$
for $\{\gam,\alp,\bet\}\subset OT(\mI_{N})$.
\eenu
}
\edf

Each finite set $\Tht$ is $C$-bounded since $\Tht\subset C_{\alp}(\bet)$, where
$\bet=\max(\Tht)+1$ if $\Tht\neq\emptyset$, and $\bet=0$ else.
Also each ordinal $\bet=\{\gam:\gam<\bet\}\subset C_{\alp}(\bet)$ as well as
a set $C_{\gam}(\del)$  is $C$-bounded.

\bprp\label{prp:Mostowski_base_H}
Let $\rho\in\Psi_{\mS}$ and $\mS\leq\alpha\in M_{\rho}$.
If $\alpha[\rho/\mS]\in \mathcal{H}_{\gamma}[\Theta]$, then $\rho\in \mathcal{H}_{\gamma}[\Theta]$.
\eprp
\bprf
This seen from Proposition \ref{prp:Mostowski_base}.
\eprf

\bdf\label{df:cloesunder}
{\rm
The set $stm(\alp)$ of immediate subterms of $\alp\in OT(\mI_{N})$ is defined as follows.
$stm(\alp)=\emptyset$ if $\alp\in\{0,\Ome,\mI_{N}\}$.
$stm(\alp)=\{\alp_{i}:i\leq m\}$ if $\alp=_{NF}\alp_{m}+\cdots+\alp_{0}$.
$stm(\alp)=\{\bet,\gam\}$ if $\alp=_{NF}\vphi\bet\gam$.
$stm(\alp)=\{\sig,a\}\cup \mathrm{fld}(f)$ if $\alp=\psi_{\sig}^{f}(a)$.
$stm(\alp)=\{\mS\}$ if $\alp=\mS^{\dagger k}$.
$stm(\alp)=\{\tau\}$ if $\alp=\mI_{N}[\tau]$.
$stm(\alp)=\{\tau\}$ if $\alp=\mT^{\dagger \vec{i}}[\tau/\mT]$.

}
\edf

\bprp\label{prp:stm}
\benu
\item\label{prp:stm.1}
If $\bet\leq\alp\in C_{b}(\bet)$, then $stm(\alp)\subset C_{b}(\bet)$.

\item\label{prp:stm.2}
Let $\alp\not\in\Psi$. Then $stm(\alp)\subset\alp$, and $\alp\in C_{b}(\bet) \Lrarw stm(\alp)\subset C_{b}(\bet)$.

\item\label{prp:stm.3}
Let $\alp_{i}\in stm(\alp_{i+1})$ for $i<m$. Assume that $\alp_{0}<\alp_{k}$ for any $0<k\leq m$.
Then $\alp_{m}\in C_{b}(\bet) \Rarw \alp_{0}\in C_{b}(\bet)$.

\item\label{prp:stm.4}
Let $k_{\delta}(\alpha)=\max(\{0\}\cup K_{\delta}(\alpha))$.
If $\alpha\in C_{b}(\beta)$ and $k_{\delta}(\alpha)<b$, then $K_{\delta}(\alpha)\subset C_{b}(\beta)$ holds.
\eenu
\eprp
\bprf
\ref{prp:stm}.\ref{prp:stm.3}.
Let $\alp_{m}\in C_{b}(\bet)$. If $\alp_{0}<\bet$, then $\alp_{0}\in C_{b}(\bet)$.
Assume $\bet\leq\alp_{0}$. Then $\bet<\alp_{k}$ for any $0<k\leq m$ by the assumption.
From $\alp_{m}\in C_{b}(\bet)$ we see $\alp_{k-1}\in C_{b}(\bet)$ inductively on $0<k\leq m$.
\\
\ref{prp:stm}.\ref{prp:stm.4}. By induction on $\ell\alpha$.
If $\alpha\not\in\Psi$, then $\alpha_{2}<\alpha$ for each $\alpha_{2}\in stm(\alpha)$ and 
$K_{\delta}(\alpha)=\bigcup\{K_{\delta}(\alpha_{2}:\alpha_{2}\in stm(\alpha)\}$.
IH yields $K_{\delta}(\alpha)\subset C_{b}(\beta)$. Let $\alpha\in\Psi$ and
we may assume $\delta<\alpha=\psi_{\sigma}^{f}(\alpha_{1})$ for some $\sigma,\alpha_{1}$ and $f$.
Then $K_{\delta}(\alpha)=\{\alpha_{1}\}\cup\bigcup\{K_{\delta}(\alpha_{2}):\alpha_{2}\in stm(\alpha)\}$.
If $stm(\alpha)\subset C_{b}(\beta)$ and $\alpha_{1}<b$, then IH yields $K_{\delta}(\alpha)\subset C_{b}(\beta)$.
Let $\alpha<\beta$.
We have $stm(\alpha)\subset C_{\alpha_{1}}(\alpha)$, and $\alpha_{1}\leq k_{\delta}(\alpha)<b$.
We obtain $\alpha_{1}\in stm(\alpha)\subset C_{\alpha_{1}}(\alpha) \subset C_{b}(\beta)$.
Hence $\alpha_{1}\in C_{b}(\beta)$, and $K_{\delta}(\alpha_{2})\subset C_{b}(\beta)$ for $\alpha_{2}\in stm(\alpha)$ by IH.
$K_{\delta}(\alpha)=\{\alpha_{1}\}\cup K_{\delta}(\alpha)\subset C_{b}(\beta)$ follows.
\eprf

\bprp\label{prp:Buchholz4.3}
Let $\Tht$ be a $C$-bounded set of ordinals.

\benu
\item\label{prp:Buchholz4.3.0}
$\calh_{\gam}[C_{\alp}(\bet)]\subset C_{\alp}(\bet)$ if $\gam<\alp$.

\item\label{prp:Buchholz4.3.1}
$\{0,\Ome,\mI_{N}\}\cup\Tht\subset\calh_{\gam}[\Tht]$.

\item\label{prp:Buchholz4.3.3}
Let $\psi_{\sig}^{f}(a)\in OT(\mI_{N})$.
If $\{\sig,a\}\cup \mathrm{fld}(f)\subset \calh_{\gam}[\Tht]$ and $a\leq\gam$, then $\psi_{\sig}^{f}(a)\in\calh_{\gam}[\Tht]$.

\item\label{prp:Buchholz4.3.2}
If $\alp\in \calh_{\gam}[\Tht]$, then $stm(\alp)\cap\alp\subset\calh_{\gam}[\Tht]$.

\item\label{prp:Buchholz4.3.4}
Let $\alp\not\in\Psi$. Then 
$\alp\in\calh_{\gam}[\Tht]\Lrarw stm(\alp)\subset \calh_{\gam}[\Tht]$.

\item\label{prp:Buchholz4.3.5}
Let $\alp_{i}\in stm(\alp_{i+1})$ for $i<m$. Assume that $\alp_{0}<\alp_{k}$ for any $0<k\leq m$.
Then $\alp_{m}\in \calh_{\gam}[\Tht] \Rarw\alp_{0}\in \calh_{\gam}[\Tht]$.

\item\label{prp:Buchholz4.3.6}
Let $k_{\delta}(\Theta)=\max(\{0\}\cup \bigcup\{K_{\delta}(\alpha):\alpha\in\Theta\})$.
Then $K_{\delta}(\Theta)\subset\mathcal{H}_{\gamma}[\Theta]$ if $\gamma\geq k_{\delta}(\Theta)$.
\eenu

\eprp
\bprf
Let $\gam<b$ and $\Tht\subset C_{b}(\bet)$.
\\
\ref{prp:Buchholz4.3}.\ref{prp:Buchholz4.3.3}.
Let $\{\sig,a\}\cup \mathrm{fld}(f)\subset \calh_{\gam}[\Tht]$ and $a\leq\gam$. Then $\{\sig,a\}\cup \mathrm{fld}(f)\subset C_{b}(\bet)$ and
$a\leq\gam<b$. Hence $\psi_{\sig}^{f}(a)\in C_{b}(\bet)$, and
$\psi_{\sig}^{f}(a)\in\calh_{\gam}[\Tht]$.
\\
\ref{prp:Buchholz4.3}.\ref{prp:Buchholz4.3.2}. We show that
$\alp\in C_{b}(\bet)\Rarw stm(\alp)\cap\alp\subset C_{b}(\bet)$.
Let $\alp\in C_{b}(\bet)$ and $\del\in stm(\alp)\cap\alp$.
If $\alp<\bet$, then $\del<\bet$, and $\del\in C_{b}(\bet)$.
Otherwise we obtain $stm(\alp)\subset C_{b}(\bet)$ by Proposition \ref{prp:stm}.\ref{prp:stm.1}.
\\
\ref{prp:Buchholz4.3}.\ref{prp:Buchholz4.3.4}. 
This follows from Proposition \ref{prp:stm}.\ref{prp:stm.2}.
\\
\ref{prp:Buchholz4.3}.\ref{prp:Buchholz4.3.5}. 
This follows from Proposition \ref{prp:stm}.\ref{prp:stm.3}.
\\
\ref{prp:Buchholz4.3}.\ref{prp:Buchholz4.3.6}.
Let $\alpha\in\Theta$. We obtain $k_{\delta}(\alpha)\leq k_{\delta}(\Theta)\leq\gamma<b$.
Proposition \ref{prp:stm}.\ref{prp:stm.4} yields $K_{\delta}(\alpha)\subset C_{b}(\beta)$.
Hence $K_{\delta}(\Theta)\subset\mathcal{H}_{\gamma}[\Theta]$ follows.
\eprf

\blem\label{prp:EK2}
Let $\rho\in\Psi_{\mS}$ for an $\mS\in SSt$.
\benu
\item\label{prp:EK2.1}
$\mathcal{H}_{\gamma}[M_{\rho}]\subset M_{\rho}$ 
if $\gam<\mathtt{ p}_{0}(\rho)$.

\item\label{prp:EK2.2}
$M_{\rho}\cap\mS=\rho$ and $\rho\not\in M_{\rho}$.

\item\label{prp:EK2.3}
If $\sig<\rho$ and $\mathtt{ p}_{0}(\sig)\leq \mathtt{ p}_{0}(\rho)$, then
$M_{\sig}\subset M_{\rho}$.
\eenu
\elem
\bprf
Lemma \ref{prp:EK2}.\ref{prp:EK2.1} follows from Proposition \ref{prp:Buchholz4.3}.\ref{prp:Buchholz4.3.0}.
\eprf

\blem\label{lem:fTaut}
Let $\sigma\in\Psi_{\mS}$ and $\alpha\in M_{\sigma}$.
Assume  $\sigma\in C_{b}(\delta)$ for $\mathtt{p}_{0}(\sigma)\leq b$ and $\delta\leq\sigma$.
Then
$\alpha[\sigma/\mS]\in C_{b}(\delta) \Lrarw \alp\in C_{b}(\delta)$.

In particular, if $\tau\in\Psi_{\mT}$ with $\mT<\mS$ and 
$\sigma\in M_{\tau}$, then
$\alpha[\sigma/\mS]\in M_{\tau} \Lrarw \alp\in M_{\tau}$ holds.
\elem
\bprf
By induction on $\ell\alp$.
We may assume that $\alpha\geq\mS$ and $\alp=\psi_{\kap}^{f}(a)\in\Psi$ by Proposition \ref{prp:stm}.\ref{prp:stm.2}.

We obtain $\beta\in M_{\sigma}$ for each $\beta\in \{\kap,a\}\cup \mathrm{fld}(f)$.
IH yields $\beta[\sigma/\mS]\in C_{b}(\delta) \Lrarw \beta\in C_{b}(\delta)$.
Moreover if $\delta\leq\sigma\leq\alpha[\sigma/\mS]\in C_{b}(\delta)$, then $\forall\beta\in \{\kap,a\}\cup \mathrm{fld}(f)(\beta[\sig/\mS]\in C_{b}(\delta))$,
and if $\alp\in C_{b}(\delta)$, then $\{\kap,a\}\cup \mathrm{fld}(f)\subset C_{b}(\delta)$.

Let $M_{\sigma}=C_{d}(\sigma)$ with $d=\mathtt{p}_{0}(\sigma)\leq b$.
We have $a<d$ by $\sigma<\mS\leq\psi_{\kap}^{f}(a)\in C_{d}(\sigma)$.
We obtain $a[\sig/\mS]\leq a<b$.
Hence $\alp[\sig/\mS]\in C_{b}(\delta) \Lrarw \alp=\psi_{\kap}^{f}(a)\in C_{b}(\delta)$.

Let $M_{\tau}=C_{b}(\tau)$ with $b=\mathtt{p}_{0}(\tau)$.
We see $d<b$ from $\sigma\in C_{b}(\tau)$ and $\tau<\mT<\sigma$.
\eprf

\bdf\label{df:ES}
{\rm
For $\alp\in OT(\mI_{N})$ and $\mS\in SSt$, a finite set $E_{\mS}(\alp)$ of subterms of $\alp$ is defined recursively as follows.

\benu
\item 
$E_{\mS}(\alp)=E_{\mS}(SC(\alp)):=\bigcup\{E_{\mS}(\bet):\bet\in SC(\alp)\}$
if $\alp\not\in SC(\alp)$.

\item 
$E_{\mS}(\alp)=\bigcup\{E_{\mS}(\bet): \bet\in stm(\alp)\}$
if $\mS<\alp\in SC(\alp)$.

\item
$E_{\mS}(\mS)=\emptyset$.

\item
$E_{\mS}(\alp)=\{\alp\}$ if $\mS>\alp\in SC(\alp)$.

\eenu
}
\edf

\bprp\label{prp:EH_Pi11}
Let $\rho\in\Psi_{\mS}$ for an $\mS\in SSt$.
\benu
\item\label{prp:EH_Pi11.1}
$E_{\mS}(\alp)\subset \mS\cap(\alp+1)$.

\item\label{prp:EH_Pi11.5}
$E_{\mS}(\alpha)\subset\calh_{0}[\{\alp\}]$.

\item\label{prp:EH_Pi11.6}
Let $\alp\in M_{\rho}$.
Then $\alp[\rho/\mS]\in \calh_{\Gam(\rho^{\dagger})}[E_{\mS}(\alp)\cup\{\rho\}]$.

\item\label{prp:EH_Pi11.4}
Let $\alp\in M_{\rho}$. Then $\alp[\rho/\mS]\in \calh_{\Gam(\rho^{\dagger})}[\{\alp,\rho\}]$.

\item\label{prp:EH_Pi11.7}
Assume $\alp\in M_{\rho}$, $E_{\mS}(\alp)\subset\sig\in\Psi_{\mS}$ and ${\tt p}_{0}(\rho)\leq{\tt p}_{0}(\sig)$.
Then $\alp\in M_{\sig}$.

\item\label{prp:EH_Pi11.8}
If $\alp\in M_{\rho}$, then $E_{\mS}(\alp)\subset\rho$.
\eenu
\eprp
\bprf
\ref{prp:EH_Pi11}.\ref{prp:EH_Pi11.5}.
Let $a>0$.
We show $\alp\in C_{a}(b) \Rarw E_{\mS}(\alpha)\subset C_{a}(b)$ by induction on $\ell\alp$. 
Let $\alp\in C_{a}(b)$.
We may assume $\mI_{N}\neq\alp>\mS$.
Then $E_{\mS}(\alp)=\bigcup\{E_{\mS}(\bet):\bet\in stm(\alp)\}$.
If $\alp\not\in\Psi$, then $stm(\alp)\subset C_{a}(b)$ by Proposition \ref{prp:stm}.\ref{prp:stm.2}, and
IH yields $E_{\mS}(\alpha)\cup K_{\mS}(\alpha)\subset C_{a}(b)$.
Let $\alp=\psi_{\sig}^{f}(\alp_{1})\in C_{a}(b)$, and $\bet\in stm(\alp)$.
If $\bet\in C_{a}(b)$, then IH yields $E_{\mS}(\beta)\cup K_{\mS}(\beta)\subset C_{a}(b)$.
Assume $\bet\not\in C_{a}(b)$.
Then $E_{\mS}(\bet)\subset\mS<\alp<b$.
$E_{\mS}(\bet)\subset C_{a}(b)$ follows.
\\
\ref{prp:EH_Pi11}.\ref{prp:EH_Pi11.6}.
Let $a>\Gam(\rho^{\dagger})$.
By induction on $\ell\alp$ we show that
$E_{\mS}(\alp)\cup\{\rho\}\subset C_{a}(b) \Rarw \alp[\rho/\mS]\in C_{a}(b)$.
We may assume that
$\alp=\psi_{\sig}^{f}(\alp_{1})>\mS$.
Then $E_{\mS}(\alp)=\bigcup\{E_{\mS}(\bet):\bet\in stm(\alp)\}$, and IH yields $\bet[\rho/\mS]\in C_{a}(b)$ for each $\bet\in stm(\alp)$.
On the other hand we have $\alp_{1}[\rho/\mS]<\Gam(\rho^{\dagger})<a$.
$\alp[\rho/\mS]\in C_{a}(b)$ follows.
\\
\ref{prp:EH_Pi11}.\ref{prp:EH_Pi11.4}.
By Propositions \ref{prp:EH_Pi11}.\ref{prp:EH_Pi11.6} and \ref{prp:EH_Pi11}.\ref{prp:EH_Pi11.5}
we obtain $\alp[\rho/\mS]\in \calh_{\Gam(\rho^{+})}[E_{\mS}(\alp)\cup\{\rho\}]$ and $E_{\mS}(\alp)\subset\calh_{0}[\{\alp\}]$.
Hence $\alp[\rho/\mS]\in \calh_{\Gam(\rho^{\dagger})}[\{\alp,\rho\}]$.
\\
\ref{prp:EH_Pi11}.\ref{prp:EH_Pi11.7}.
By induction on $\ell\alp$. Let $b={\tt p}_{0}(\rho)\leq{\tt p}_{0}(\sig)=d$.
If $\alp<\mS$, then $E_{\mS}(\alp)=SC(\alp)\subset\sig$, and $\alp<\sig$.
We may assume $\mS<\alp=\psi_{\sig}^{f}(\alp_{1})\in M_{\rho}=C_{b}(\rho)$.
We obtain $\alp_{1}<b\leq d$, $E_{\mS}(\alp)=\bigcup\{E_{\mS}(\bet):\bet\in stm(\alp)\}$ and $stm(\alp)\subset C_{b}(\rho)$.
IH yields 
$stm(\alp) \subset M_{\sig}=C_{d}(\sig)$, and 
$\alp\in C_{d}(\sig)$ follows.
\\
\ref{prp:EH_Pi11}.\ref{prp:EH_Pi11.8}.
By induction on $\ell\alp$ using the fact $M_{\rho}\cap\mS=\rho$.
\eprf

\bprp\label{prp:inverse}
Let $\alpha\in M_{\rho}$ with $\rho\in\Psi_{\mS}$ and $\mS\in SSt$, and $\mathtt{p}_{0}(\rho)\leq a$.
If $E_{\mS}(\alpha)\cup\{\mS\}\subset C_{a}(\beta)$, then $\alpha\in C_{a}(\beta)$.
\eprp
\bprf
By induction on $\ell\alpha$.
We may assume that $\mS<\alpha\in\Psi$ by IH and Proposition \ref{prp:stm}.\ref{prp:stm.2}.
 and \ref{prp:Buchholz4.3}.\ref{prp:Buchholz4.3.4}.
Let $\alpha=\psi_{\sigma}^{f}(b)$.
 $\Theta\subset C_{a}(\beta)$ with $a>\gamma\geq\mathtt{p}_{0}(\rho)$.
We have $\bigcup \{E_{\mS}(\alpha_{0}):\alpha\in stm(\alpha)\}=E_{\mS}(\alpha)\subset C_{a}(\beta)$.
IH yields $stm(\alpha)\subset C_{a}(\beta)$.
On the other hand we have $\rho<\mS<\alpha\in M_{\rho}=C_{\mathtt{p}_{0}(\rho)}(\rho)$.
We obtain $b<\mathtt{p}_{0}(\rho)\leq a$, and 
$\alpha\in C_{a}(\beta)$ follows.
\eprf

\blem\label{lem:inverse}
Let $\alpha\in M_{\rho}$ with $\rho\in\Psi_{\mS}$ and $\mS\in SSt$, and $\max\{\mS,\mathtt{p}_{0}(\rho)\}\leq \gamma$.
Assume $\rho\in \mathcal{H}_{\gamma}[\Theta]$. Then
$\alpha[\rho/\mS]\in \mathcal{H}_{\gamma}[\Theta] \Lrarw \alpha\in \mathcal{H}_{\gamma}[\Theta]$.
\elem
\bprf
If $\{\alpha,\rho\}\subset \mathcal{H}_{\gamma}[\Theta]$, then $\alpha[\rho/\mS]\in \mathcal{H}_{\gamma}[\Theta]$
by Proposition \ref{prp:EH_Pi11}.\ref{prp:EH_Pi11.4}.

Assume $\rho\in\mathcal{H}_{\gamma}[\Theta]$.
Let $\Theta\subset C_{a}(\beta)$ with $a>\gamma\geq\mathtt{p}_{0}(\rho)$.
We have $\rho\in C_{a}(\beta)$.
It suffices to show that if $\alpha[\rho/\mS]\in C_{a}(\beta)$, then $\alpha\in C_{a}(\beta)$ by induction on $\ell\alpha$.
We may assume that $\mS\leq\alpha\in\Psi$ by IH and Proposition \ref{prp:stm}.\ref{prp:stm.2}.

If $\alpha=\mS$, then $\rho=\alpha[\rho/\mS]\in C_{a}(\beta)$.
By Proposition \ref{prp:Erho} we obtain $\mS\in M_{\rho}=C_{\mathtt{p}_{0}(\rho)}(\rho)$ with $\mathtt{p}_{0}(\rho)\leq \gamma<a$.
If $\rho\leq\beta$, then we obtain $\mS\in C_{a}(\beta)$.
Let $\beta<\rho$. Then $\rho\in C_{a}(\beta)$ yields $\mS\in C_{a}(\beta)$.

Next let $\alpha=\psi_{\sigma}^{f}(b)>\mS$. We have $\alpha[\rho/\mS]\in C_{a}(\beta)$.
If $\alpha[\rho/\mS]\geq\beta$, then we obtain $\{\delta[\rho/\mS]:\delta\in stm(\alpha)\}\subset C_{a}(\beta)$, and
$stm(\alpha)\subset C_{a}(\beta)$ by IH.
On the other hand we have $\rho<\mS<\alpha\in M_{\rho}=C_{\mathtt{p}_{0}(\rho)}(\rho)$, and $b<\mathtt{p}_{0}(\rho)< a$.
Hence $\alpha\in C_{a}(\beta)$.
Let $\alpha[\rho/\mS]<\beta$. We obtain $E_{\mS}(\alpha)\subset\rho\leq\alpha[\rho/\mS]<\beta$ by $\alpha\in M_{\rho}$ and 
Proposition \ref{prp:EH_Pi11.8}.\ref{prp:EH_Pi11.8}, and
$E_{\mS}(\alpha)\subset C_{a}(\beta)$.
On the other hand we have $\mS\in C_{a}(\beta)$ by the case $\alpha=\mS$.
Therefore $\alpha\in C_{a}(\beta)$ follows from Proposition \ref{prp:inverse}.
\eprf
\\

The trail in the following Proposition \ref{prp:trail} is needed  in Propositions \ref{prp:stblS1} and \ref{prp:g01}.

\bprp\label{prp:trail}
Let $\alp$ be a strongly critical number such that $\Ome<\alp<\mI_{N}$.
There exists a unique sequence
$(\alp_{n})_{n\leq m}$ such that 
either $\alp_{0}=\psi_{\mI_{N}}(a)$ for an $a$ or $\alpha_{0}=0$, $\alp_{m}=\alp$ and
each $\alp_{n+1}$ is one of the forms
$\psi_{\alp_{n}}^{f}(b)$, $\alp_{n}^{\dagger\vec{i}}$, $\mI_{N}[\alp_{n}]$, $\mS^{\dagger\vec{i}}[\alp_{n}/\mS]$
for some $f,b,\vec{i}$ and $\mS$.
The sequence $(\alp_{n})_{n\leq m}$ is said to be the {\rm trail} to $\alp$, and denoted by 
$\mathrm{trail}(\alp)$.

\benu
\item\label{prp:trail.1}
For a term $\alp_{n}$ in the trail to $\alp$, assume $\alp_{n}<\alp$.
Then
$\alp_{n}<\alp_{k}$ for $n<k\leq m$, 
$\alp\in\calh_{\gam}[\Tht] \Rarw \alp_{n}\in\calh_{\gam}[\Tht]$, and $E_{\mT}(\alp_{n})\subset E_{\mT}(\alp)$ for every $SSt\ni\mT\leq\alp_{n}$.

\item\label{prp:trail.2}
$\alp_{0}\leq\alp$, $\alp\in\calh_{\gam}[\Tht] \Rarw \alp_{0}\in\calh_{\gam}[\Tht]$, 
and $E_{\mT}(\alp_{0})\subset E_{\mT}(\alp)$ holds for every
$SSt\ni\mT\leq\alp_{0}$.
\eenu
\eprp
\bprf
This is seen by inspection of Definitions \ref{df:notationsystem} and \ref{df:notationsystem.2M}.
If $\alp_{n}>\alp_{n+1}$, then we would have $\alp_{n+1}\prec\alp_{n}$ and $\alp<\alp_{n}$ by
Definition \ref{df:lessthanpsi}.

In Proposition \ref{prp:trail}.\ref{prp:trail.1}, $\alp\in\calh_{\gam}[\Tht] \Rarw \alp_{n}\in\calh_{\gam}[\Tht]$ follows from
Proposition \ref{prp:Buchholz4.3}.\ref{prp:Buchholz4.3.5}.
\eprf

\bprp\label{prp:ESM}
Let $\rho\in\Psi_{\mS}$ with a successor stable ordinal $\mS$.
Assume $k_{\mS}(\alpha)=\max(\{0\}\cup K_{\mS}(\alpha))<\mathtt{p}_{0}(\rho)$ and
$E_{\mS}(\alp)\subset\rho$.
Then $\alp\in M_{\rho}=C_{\mathtt{p}_{0}(\rho)}(\rho)$.
\eprp
\bprf
By induction on $\ell\alp$. Note that
$k_{\mS}(\alpha)<\mathtt{p}_{0}(\rho)$ iff $\alpha\in C_{\mathtt{p}_{0}(\rho)}(\mS)$.
We have $\mS\in M_{\rho}$ by Proposition \ref{prp:Erho}.
By IH we may assume that $\mS<\alp<\mI_{N}$ and $\alp\in\Psi$ by Proposition \ref{prp:stm}.\ref{prp:stm.2}.
Let $\alp=\psi_{\sig}^{f}(a)$. For $stm(\alp)=\{\sig,a\}\cup \mathrm{fld}(f)$, we obtain $E_{\mS}(\alp)=E_{\mS}(stm(\alp))$ and
 $stm(\alp)\subset C_{\mathtt{p}_{0}(\rho)}(\mS)$ with $a<\mathtt{p}_{0}(\rho)$. 
IH yields $stm(\alp)\subset M_{\rho}$, and $\alp\in M_{\rho}$ by $a<\mathtt{p}_{0}(\rho)$.
\eprf
\\

The  axiom (\ref{eq:stbl0}) states that each $(i+1)$-stable ordinal is a limit of $i$-stable ordinals.
In proving the axiom (\ref{eq:stbl0})
in Lemma \ref{th:embedreg}
we need the fact that if both of an $(i+1)$-stable ordinal $\mT$ and an ordinal $\alp<\mT$ are `captured'
in $\calh_{\gam}[\Tht]$, then so is a successor $i$-stable ordinal $\mS$ such that $\alp<\mS<\mT$.
Or in other words, such an $\mS$ should be constructed from data included in ordinals $\mT$ and $\alp$.
The data we need are trails.

\bprp\label{prp:stblS1} 
Let $\mT\in SSt_{i+1}$ be a successor $(i+1)$-stable ordinal, and $\alp<\mT$ an ordinal.
Assume $\{\alp,\mT\}\subset\calh_{\gam}[\Tht]$.
Then there exists a successor $i$-stable ordinal $\alp<\mS<\mT$ such that
$\mS\in \calh_{\gam}[\Tht]$.
\eprp
\bprf
By induction on the lengths $\ell\alp$ of ordinal terms $\alp$.
By IH we may assume that $\Ome<\alp<\mI_{N}$ and $\alp\in SC(\alp)$.
Let $\mT=\mU^{\dagger(i+1)}$ with $\mU\in St\cup\{\Ome\}$. 
\\
\textbf{Case 1}. There exists a $k>0$ such that $\alp<\mU^{\dagger i^{(k)}}$, where
$\mU^{\dagger i^{(0)}}=\mU$ and $\mU^{\dagger i^{(k+1)}}=(\mU^{\dagger i^{(k)}})^{\dagger i}$:
Pick a $k>0$ such that $\alp<\mS=\mU^{\dagger i^{(k)}}$.
We obtain $\mS\in\calh_{\gam}[\Tht]$ by $\mT\in\calh_{\gam}[\Tht]$ and Proposition \ref{prp:Buchholz4.3}.\ref{prp:Buchholz4.3.4}.
$\mS<\mT$ is seen from $\mT\in LSt_{i}$.
\\
\textbf{Case 2}. Otherwise:
Then we see from Definition \ref{df:lessthanpsi} that there exists a $\rho\leq\alp$ such that $\rho\prec\mT$, and $\rho\in{\rm trail}(\alp)$ if $\rho<\alp$.
From $\alp\in\calh_{\gam}[\Tht]$ and Proposition \ref{prp:trail}.\ref{prp:trail.1} we see $\rho\in\calh_{\gam}[\Tht]$.
Pick a $k>0$ such that $\alp<\mS=\rho^{\dagger i^{(k)}}<\mT$.
We obtain $\mS\in\calh_{\gam}[\Tht]$ by Proposition \ref{prp:Buchholz4.3}.\ref{prp:Buchholz4.3.4}.
\eprf
\\

The following Definition \ref{df:special} is needed in subsection \ref{subsec:elimrfl}.

\bdf\label{df:special}
{\rm
Let
$s(f)=\max(\{0\}\cup{\rm supp}(f))$ for finite function $f$.

Let  
$\Lambda\leq\mI_{N}$ be a strongly critical number, which is a base for $\tilde{\tht}$-function.
Let $f:\Lambda\to\Gamma(\Lambda)$ be a non-empty and irreducible finite function.
Then $f$ is said to be a \textit{special} finite function with base $\Lambda$ if there exists an ordinal $\alpha$
such that $f(s(f))=\Lambda(\alpha+1)$.
For a special finite function $f$, $f^{\prime}$ denotes a finite function such that
${\rm supp}(f^{\prime})={\rm supp}(f)$,
$f^{\prime}(c)=f(c)$ for $c\neq s(f)$, and
$f^{\prime}(s(f))=\Lambda\cdot\alpha$ with $f(s(f))=\Lambda(\alpha+1)$.
}
\edf

A special function $h^{b}(g;a)$ is defined from ordinals $a,b$ and a finite function $g$
as in \cite{singlestable}.

The following Definition \ref{df:hstepdown} is mainly needed in subsection \ref{subsec:elimrfl}, but
also in Definition \ref{df:resolvent}, cf.\,the beginning of subsection \ref{subsec:elimrfl}.

\bdf\label{df:hstepdown}
{\rm
Let $f, g:\Lambda\to\Gamma(\Lambda)$ be special finite functions.

\begin{enumerate}
\item\label{df:hstepdown.1}
For 
ordinals $a\leq\Lambda$, $b\leq s(g)$, 
let us define a special finite function $h=h^{b}(g;a)$ 
as follows.
$s(h)=b$, and
$h_{b}=g_{b}$.
To define $h(b)$,
let $\{b=b_{0}<b_{1}<\cdots<b_{n}=s(g)\}=\{b,s(g)\}\cup\left((b,s(g))\cap {\rm supp}(g)\right)$.
Define recursively ordinals $\alpha_{i}$ by
$\alpha_{n}=\Lambda\cdot\alpha+a$ with $g(s(g))=\Lambda(\alpha+1)$.
$\alpha_{i}=g(b_{i})+\tilde{\theta}_{c_{i}}(\alpha_{i+1})$ for
$c_{i}=b_{i+1}-b_{i}$.
Finally let $h(b)=\alpha_{0}+\Lambda$.

Note that $h=h^{b}(g;a)=g$ when $b=s(g)$.

\item\label{df:hstepdown.4}
$f_{b}*g^{b}$ denotes a special function $h$ such that
${\rm supp}(h)={\rm supp}(f_{b})\cup{\rm supp}(g^{b})$,
$h^{\prime}(c)=f^{\prime}(c)$ for $c<b$, and
$h^{\prime}(c)=g^{\prime}(c)$ for $c\geq b$.

\end{enumerate}
}
\edf
The following Proposition \ref{prp:hstepdown} is seen as in \cite{singlestable}.

\bprp\label{prp:hstepdown}
Let 
$f, g:\Lambda\to\Gamma(\Lambda)$ be special finite functions such that
$f_{d}\leq g_{d}$ and
$f<_{\Lambda}^{d}g^{\prime}(d)$ for a $d\in{\rm supp}(g)$.

\begin{enumerate}

\item\label{prp:hstepdown.2}
For $b<d$ and $a<\Lambda$, 
$f_{b}\leq (h^{b}(g;a))_{b}$ and
$f<_{\Lambda}^{b}(h^{b}(g;a))^{\prime}(b)$.

\item\label{prp:hstepdown.3}
For $d\leq c<s(f)$ and $a\leq\Lambda$,
$h^{c}(f;a)<^{d}g^{\prime}(d)$.

\item\label{prp:hstepdown.4}
Let $b\leq e<d$, 
$a_{0}<a\leq \Lambda$, and
 $h=(h^{e}(g;a_{0}))*f^{e+1}$.
Then
$h_{b}=(h^{b}(g;a))_{b}$,
$f\leq h$, and
$h<_{\Lambda}^{b}(h^{b}(g;a))^{\prime}(b)$.

\end{enumerate}
\eprp
\bprf
\ref{prp:hstepdown}.\ref{prp:hstepdown.2}.
Let $h=h^{b}(g;a)$.
We have $f_{b}\leq g_{b}=h_{b}$.
Let $b+x\in{\rm supp}(f)\cap d\subset{\rm supp}(g)\cap d$.
Then $f(b+x)\leq  g(b+x)<\tilde{\theta}_{-x}(h^{\prime}(b))$ and
$g^{\prime}(d)< \tilde{\theta}_{-(d-b)}(h^{\prime}(b))$.
Proposition \ref{prp:idless} yields the proposition.
\\
\ref{prp:hstepdown}.\ref{prp:hstepdown.3}.
First we show $(h^{d}(f;a))(d)<g^{\prime}(d)$ when $c=d$.
Let 
$\{d=d_{0}<d_{1}<\cdots<d_{n}=s(f)\}=
\{d,s(f)\} \cup \left( (d, s(f))\cap{\rm supp}(f)\right)$.
Define recursively
$\alpha_{n}=f^{\prime}(d_{n})+a$ and
$\alpha_{i}=f(d_{i})+\tilde{\theta}_{d_{i+1}-d_{i}}(\alpha_{i+1})$.
Then $(h^{d}(f;a))(d)=\alpha_{0}+\Lambda$.
We have $f(d_{0})<^{d_{0}}g^{\prime}(d)$.
Let $\mu_{0}$ be the shortest part of $g^{\prime}(d)$ such that $f(d_{0})<\mu_{0}$,
and $\mu_{i+1}$ be the shortest part of $\tilde{\theta}_{-(d_{i+1}-d_{i})}(tl(\mu_{i}))$
such that
$f(d_{i+1})<\mu_{i+1}$.

We show by induction on $n-i$ that $\alpha_{i}<\mu_{i}$. 
We have $\alpha_{n}=f^{\prime}(d_{n})+a\leq f^{\prime}(d_{n})+\Lambda=f(d_{n})<\mu_{n}$.
For $i<n$,
$\tilde{\theta}_{d_{i+1}-d_{i}}(\alpha_{i+1})<
\tilde{\theta}_{d_{i+1}-d_{i}}(\tilde{\theta}_{-(d_{i+1}-d_{i})}(tl(\mu_{i})))\leq
tl(\mu_{i})$ by IH and Proposition \ref{prp:tht4}.\ref{prp:tht4.2}.
On the other hand we have $f(d_{i})<\mu_{i}$.
Hence
$\alpha_{i}=f(d_{i})+\tilde{\theta}_{d_{i+1}-d_{i}}(\alpha_{i+1})<\mu_{i}$.
We obtain $\alpha_{0}<\mu_{0}$,
and $\alpha_{0}+\Lambda<\mu_{0}\leq g^{\prime}(d)$
by $\Lambda<tl(\mu_{0})$.

Next let $d<c$. We have $(h^{c}(f;a))(d)=f(d)<g^{\prime}(d)$ and
$f(d_{i})<\mu_{i}$ for $d_{i}\in{\rm supp}(f)\cap(d,c)$.
Let $i< n$ be the least such that $c<d_{i+1}$.
If $c=d_{i}$, then $(h^{c}(f;a))(c)=\alpha_{i}+\Lambda<\mu_{i}$.
Let $d_{i}<c$. Then
$(h^{c}(f;a))(c)=\tilde{\theta}_{d_{i+1}-c}(\alpha_{i+1})+\Lambda$
and $\tilde{\theta}_{d_{i+1}-c}(\alpha_{i+1})\leq\tilde{\theta}_{d_{i+1}-d_{i}}(\alpha_{i+1})<
tl(\mu_{i})$. Hence $(h^{c}(f;a))(c)<tl(\mu_{i})$. 
\\
\ref{prp:hstepdown}.\ref{prp:hstepdown.4}.
Note that $h=(h^{e}(g;a_{0}))^{\prime}*f^{e+1}$.
We have $h_{b}=g_{b}=(h^{b}(g;a))_{b}$.
For $b+x\in{\rm supp}(g)\cap e$, $h(b+x)=(h^{e}(g;a_{0}))(b+x)=g(b+x)<\tilde{\theta}_{-x}((h^{b}(g;a))^{\prime}(b))$, and
$h(e)=(h^{e}(g;a_{0}))(e)<\tilde{\theta}_{-(e-b)}((h^{b}(g;a))^{\prime}(b))$ by $a_{0}<a$.
For $e<e+x\in{\rm supp}(f)\cap d$, we obtain
$f(e+x)\leq g(e+x)=h(e+x)<
\tilde{\theta}_{-(e+x-b)}((h^{b}(g;a))^{\prime}(b))$.
We obtain
$f(d+x)\leq h(d+x)<\tilde{\theta}_{-x}(g^{\prime}(d))\leq
\tilde{\theta}_{-(d+x-b)}((h^{b}(g;a))^{\prime}(b))$ for $d+x\in{\rm supp}(f)$.
Therefore
$h<_{\Lambda}^{b}(h^{b}(g;a))^{\prime}(b)$.

Finally $f(e)\leq g(e)\leq h(e)$ yields $f\leq h$.
\eprf

\section{Operator controlled derivations}\label{sec:operatorcont}
We prove Theorem \ref{thm:2} assuming that the notation system $OT(\mI_{N})$ is a well ordering.
Operator controlled derivations are introduced by Buchholz\cite{Buchholz}, which we follow.
In this section except otherwise stated,
$\alp,\bet,\gam,a,b,c,d,\ldots$ and $\xi,\zeta,\nu,\mu,\ldots$ range over ordinal terms in $OT(\mI_{N})$,
$f,g,h,\ldots$ range over finite functions.

\subsection{Classes of sentences}\label{subsec:classformula}

Following Buchholz\cite{Buchholz} let us introduce a language of ramified set theory $RS$.

\bdf
{\rm
$RS$\textit{-terms} and their \textit{levels} are inductively defined.
\benu
\item
For each $\alp\in OT(\mI_{N})\cap \mI_{N}$, ${\sf L}_{\alp}$ is an $RS$-term of level $\alp$.

\item
Let $\phi(x,y_{1},\ldots,y_{n})$ be a set-theoretic formula in the language $\{\in\}$, and $a_{1},\ldots,a_{n}$ 
$RS$-terms of levels$<\! \alp\in OT(\mI_{N})\cap\mI_{N}$.

Then 
$[x\in {\sf L}_{\alp}:\phi^{{\sf L}_{\alp}}(x,a_{1},\ldots,a_{n})]$ 
is an $RS$-term of level $\alp$.
\eenu
$|u|$ denotes the level of $RS$-terms $u$, and $Tm(\alp)$ the set of $RS$-terms of level$<\alp\in OT(\mI_{N})\cap(\mI_{N}+1)$.
$Tm=Tm(\mI_{N})$ is then the set of $RS$-terms, which are denoted by $u,v,w,\ldots$

Let us identify the individual constant $\mathsf{M}_{0}$ in the language of $S_{\mI_{N}}$
with the $RS$-term ${\sf L}_{\Ome}$.
}
\edf

\bdf\label{df:Mostowski_list}
{\rm
\benu

\item\label{df:Mostowski_list2}
Let $\rho_{i}\in\Psi_{\mS_{i}}$ be ordinals with successor stable ordinals $\mS_{i}$ for $i<n$.
For a list $\bar{\rho}=(\rho_{0},\rho_{1},\ldots,\rho_{n-1})$ of ordinals,
\\
$[\bar{\rho}]:=[\rho_{0}/\mS_{0}][\rho_{1}/\mS_{1}]\cdots[\rho_{n-1}/\mS_{n-1}]$ is said to be a \textit{collapsing list} denoting
the list of Mostowski collapsings
$[\rho_{i}/\mS_{i}]$ in Definition \ref{df:notationsystem.2M}.

$[\bar{\rho}]$ applies collapsings repeatedly to ordinals $\alpha$ as follows.
\\
$\alpha[\bar{\rho}]=(\cdots((\alpha[\rho_{0}/\mS_{0}])[\rho_{1}/\mS_{1}])\cdots)[\rho_{n-1}/\mS_{n-1}]$
for ordinals $\alpha$
 such that $\alpha\in M_{[\bar{\rho}]}:\Lrarw\forall i<n(
 (\cdots((\alpha[\rho_{0}/\mS_{0}])[\rho_{1}/\mS_{1}])\cdots)[\rho_{i-1}/\mS_{i-1}]\in M_{\rho_{i}})$.

The empty list $\bar{\rho}$ with $n=0$ is denoted by $\emptyset$.
Let $M_{[\emptyset]}=OT(\mI_{N})$ and $\alpha[\emptyset]=\alpha$.

\item\label{df:Mostowski_list3}
For an ordinal $\tau\in\Psi_{\mT}$ let
$[\bar{\rho}]*[\tau]:=[\rho_{0}/\mS_{0}][\rho_{1}/\mS_{1}]\cdots[\rho_{n-1}/\mS_{n-1}][\tau/\mT]$.
\eenu
}
\edf

\bdf
{\rm
\benu

\item
$RS$-\textit{literals} 
are 
$u\in v, u\not\in v$
and
$st_{i}^{[\bar{\rho}]}(u), \lnot st_{i}^{[\bar{\rho}]}(u)$
for $RS$-terms $u,v$, $0<i\leq N$ and collapsing lists $[\bar{\rho}]$.
When $\bar{\rho}=\emptyset$, we write $st_{i}(u)$ for $st_{i}^{[\emptyset]}(u)$.

$RS$-\textit{formulas} 
are constructed from $RS$-\textit{literals}
by propositional connectives $\lor,\land$, bounded quantifiers $\exi x\in u, \fal x\in u$ and
unbounded quantifiers $\exi x,\fal x$.
Unbounded quantifiers $\exists x,\forall x$ are denoted by $\exists x\in {\sf L}_{\mI_{N}},\forall x\in {\sf L}_{\mI_{N}}$, resp.

It is convenient for us not to restrict propositional connectives $\lor,\land$ to binary ones.
Specifically when $A_{i}$ are $RS$-formulas for $i<n<\ome$,
$A_{0}\lor\cdots\lor A_{n-1}$ and $A_{0}\land\cdots\land A_{n-1}$ are $RS$-formulas.
Even when $n=1$, $A_{0}\lor\cdots\lor A_{0}$ is understood to be different from the formula $A_{0}$.
For $\Gam=\{A_{i}:i<n\}$ we write
$\bigvee\Gam\equiv(A_{0}\lor\cdots\lor A_{n-1})$ and 
$\bigwedge\Gam\equiv(A_{0}\land\cdots\land A_{n-1})$.

\item
For $RS$-terms and $RS$-formulas $\iota$,
${\sf k}(\iota)$ denotes the set of ordinal terms $\alp$ such that the constant ${\sf L}_{\alp}$ occurs in $\iota$.
Specifically $\mathsf{k}(st_{i}^{[\bar{\rho}]}(u))
=\mathsf{k}(u)$.

Let $|\iota|=\max(\sfk(\iota)\cup\{0\})$ and $\sfk(n)=\emptyset$ and $|n|=0$ for natural numbers $n$.

\item
$\mathcal{L}_{i}=\{\in\}\cup\{st_{j}^{[\bar{\rho}]} :0<j<i, \,  [\bar{\rho}] \mbox{ {\rm is a collapsing list}}\}$.
For $RS$-formulas $A$, $A\in\mathcal{L}_{i}$ designates that $A$ is a formula over the language 
$\mathcal{L}_{i}$.

\item
$\Del_{0}(\mathcal{L}_{i})$-formulas, $\Sig_{1}(\mathcal{L}_{i})$-formulas and $\Sig(\mathcal{L}_{i})$-formulas are defined as in \cite{Ba}.
Specifically if $\psi$ is a $\Sig(\mathcal{L}_{i})$-formula, then so is the formula $\fal y\in z\, \psi$.
$\tht^{(u)}$ denotes a $\Del_{0}(\mathcal{L}_{i})$-formula obtained from a $\Sig(\mathcal{L}_{i})$-formula $\tht$
by restricting each unbounded existential quantifier to $u$.

\item 
For a $\Sig_{1}(\mathcal{L}_{i})$-formula $\psi(x_{1},\ldots,x_{m})$
and 
$u_{1},\ldots,u_{m}\in Tm(\kap)$ with $\kap\leq\mI_{N}$,
$\psi^{({\sf L}_{\kap})}(u_{1},\ldots,u_{m})$
 is a $\Sig_{1}(\mathcal{L}_{i}:\kap)$\textit{-formula}.
 $\Del_{0}(\mathcal{L}_{i}:\kap)$-formulas and $\Sig(\mathcal{L}_{i}:\kap)$-formulas are defined similarly

Let $\Sig_{1}(\kap)=\Sig_{1}(\mathcal{L}_{N+1}:\kap)$ and $\Delta_{0}(\kap)=\Delta_{0}(\mathcal{L}_{N+1}:\kap)$.
\item
For $\tht\equiv\psi^{({\sf L}_{\kap})}(u_{1},\ldots,u_{m})\in\Sig(\mathcal{L}_{i}:\kap)$ 
and $\lam<\kap$,
with $u_{1},\ldots,u_{m}\in Tm(\lam)$,
$\tht^{(\lam,\kap)}:\equiv \psi^{({\sf L}_{\lam})}(u_{1},\ldots,u_{m})$.

\eenu
}
\edf

In what follows we consider only \textit{sentences} without free variables.
Sentences are denoted $A,C$ possibly with indices.

For each sentence $A$, either a disjunction is assigned as $A\simeq\bigvee(A_{\iota})_{\iota\in J}$, or
a conjunction is assigned as $A\simeq\bigwedge(A_{\iota})_{\iota\in J}$.
In the former case $A$ is said to be a \textit{$\bigvee$-formula}, and in the latter
$A$ is a \textit{$\bigwedge$-formula}.
By $st_{i}(u)$ we understand that 
there is a \textit{successor} $i$-stable ordinal $\mS$ such that $\mathsf{L}_{\mS}=u$.
It is convenient for us, cf.\,Recapping \ref{lem:recapping}, to modify
the assignment of disjunctions and conjunctions to sentences from \cite{Buchholz} 
such that
if $A\simeq\bigvee(A_{\iota})_{\iota\in J}$ is a $\bigvee$-formula, then each 
$A_{\iota}$ is a $\bigwedge$-formula, and similarly for $\bigwedge$-formulas $A$.

\bdf\label{df:bigveewedge}
{\rm
If $A$ is a $\bigvee$-formula, then let $A^{\lor}:\equiv A$.
Otherwise let
$A^{\lor}:\equiv \left(\bigvee_{i\leq 0}B_{i}\right)$ with $B_{0}\equiv A$.
If $A$ is a $\bigwedge$-formula, then let $A^{\land}:\equiv A$.
Otherwise let
$A^{\land}:\equiv \left(\bigwedge_{i\leq 0}B_{i}\right)$ with $B_{0}\equiv A$.
}
\edf

\begin{definition}\label{df:assigndc}
{\rm

\benu
\item\label{df:assigndc.1}
For $v,u\in Tm(\mI_{N})$ with $|v|<|u|$, let
\[
(v\dot{\in} u) :\equiv
\left\{
\begin{array}{ll}
A(v) & \mbox{{\rm if }} u\equiv[x\in {\sf L}_{\alp}: A(x)]
\\
v\not\in {\sf L}_{0} & \mbox{{\rm if }} u\equiv {\sf L}_{\alp}
\end{array}
\right.
\]
and 
$(u=v):\equiv(\fal x\in u(x\in v)\land \fal x\in v(x\in u))$.

\item
When $A\simeq\bigvee(A_{\iota})_{\iota\in J}$,
let $\lnot A\simeq\bigwedge(A_{\iota})_{\iota\in J}$.

\item\label{df:assigndc0}
For $RS$-terms $v,u$, 
let 
$(v\in u)  \simeq  \bigvee(A_{w,0}\land A_{w,1}\land A_{w,2})_{w\in J}$, 
and
$(v\not\in u)  \simeq  \bigwedge(\lnot A_{w,0}\lor \lnot A_{w,1}\lor \lnot A_{w,2})_{w\in J}$,
where $J=Tm(|u|)$,
$A_{w,0}\equiv (w\dot{\in} u)^{\lor}$, $A_{w,1}\equiv(\forall x\in w(x\in v))^{\lor}$ and 
$A_{w,2}\equiv(\forall x\in v(x\in w))^{\lor}$.

\item
$(A_{0}\lor\cdots\lor A_{n-1})\simeq \bigvee(A_{\iota}^{\land})_{\iota\in J}$ 
and
$(A_{0}\land\cdots\land A_{n-1})\simeq \bigwedge(A_{\iota}^{\lor})_{\iota\in J}$ for 
$J:=n$.

\item\label{df:assigndc.3}
For $u\in Tm(\mI_{N})\cup\{{\sf L}_{\mI_{N}}\}$,
$\exists x\in u\, A(x)\simeq\bigvee(A_{v})_{v\in J}$
for 
$A_{v}:\equiv ((v\dot{\in} u)^{\lor} \land (A(v))^{\lor})$ and
$J=Tm(|u|)$, where
$Tm(|{\sf L}_{\mI_{N}}|)=Tm(\mI_{N})=Tm$ and 
$(v\dot{\in} {\sf L}_{\mI_{N}}):\equiv(v\not\in {\sf L}_{0})$.

\item\label{df:assigndc.4}
For a collapsing list $[\bar{\rho}]$, let
$st_{i}^{[\bar{\rho}]}(u)\simeq \bigvee( {\sf L}_{\mS[\bar{\rho}]}=u
)_{{\sf L}_{\mS[\bar{\rho}]}\in J_{i}}$ 
with 
$J_{i}=\{{\sf L}_{\mS[\bar{\rho}]}: \mS\in SSt_{i}\cap M_{[\bar{\rho}]},\, 
|u|\geq\mS[\bar{\rho}] \}$, where
$st_{i}^{[\bar{\rho}]}$ denotes a predicate symbol,
while
$SSt_{i}\subset OT(\mI_{N})$ in the definition of $J_{i}$.

\item\label{df:assigndc.5}
For $A\simeq\bigvee(A_{\iota})_{\iota\in J}$ let
$[\rho]J=\{\iota\in J: \sfk(\iota)\subset M_{\rho}\}$.
\eenu

}
\end{definition}

It is clear that
$\sfk(A_{\iota})\subset\calh_{0}[\sfk(A)\cup\sfk(\iota)]$.

The rank $\mathrm{rk}(\iota)$ of 
sentences or terms $\iota$ 
is defined as in \cite{Buchholz} so that the following Proposition \ref{lem:rank} holds.
For completeness let us reproduce it.

\bdf\label{df:rank}
{\rm

$\mathrm{rk}(\lnot A):=\mathrm{rk}(A)$.
$\mathrm{rk}(\mathsf{L}_{\alpha})=\omega\alpha$.
$\mathrm{rk}([x\in \mathsf{L}_{\alpha}: A(x)])=\max\{\omega\alpha+1, \mathrm{rk}(A(\mathsf{L}_{0}))+2\}$.
$\mathrm{rk}(v\in u)=\max\{\mathrm{rk}(v)+6,\mathrm{rk}(u)+1\}$.
$\mathrm{rk}(A_{0}\lor \cdots\lor A_{n-1})=\max(\{0\}\cup \{\mathrm{rk}(A_{i})+1:i<n\})$.
 $\mathrm{rk}(\exists x\in u\, A(x))=\max\{\mathrm{rk}(u), \mathrm{rk}(A(\mathsf{L}_{0}))+2\}$ for $u\in Tm(\mathbb{I}_{N})\cup\{\mathsf{L}_{\mathbb{I}_{N}}\}$.
And
$\rk(st_{i}^{[\bar{\rho}]}(u))=\rk(u)+5$.

For finite sets $\Delta$ of sentences, let
$\mathrm{rk}(\Delta)=\max(\{0\}\cup\{\mathrm{rk}(\delta):\delta\in\Delta\})$.
}
\edf

We see the following Proposition \ref{lem:rank} as in \cite{Buchholz}.

\bprp\label{lem:rank}
Let $A$ be a sentence with 
$A\simeq\bigvee(A_{\iota})_{\iota\in J}$ or $A\simeq\bigwedge(A_{\iota})_{\iota\in J}$.
\begin{enumerate}

\item\label{lem:rank0}
$\mbox{{\rm rk}}(A)<\mathbb{I}_{N}+\omega$.

\item\label{lem:rank1}
$\omega |u|\leq \mbox{{\rm rk}}(u)\in\{\omega |u|+i: i\in\omega\}$, and
$\ome |A|\leq \mbox{{\rm rk}}(A)\in\{\omega |A|+i  : i\in\omega\}$.

\item\label{lem:rank0.5}
$\rk(st_{i}^{[\bar{\rho}]}(u))\in\{\rk(u)+i:i<\ome\}$.

\item\label{lem:rank2}
$\forall\iota\in J(\mbox{{\rm rk}}(A_{\iota})<\mbox{{\rm rk}}(A))$.

\item\label{lem:rank3}
$\forall\iota\in J(\mathsf{k}(\iota)\subset \mathsf{k}(A_{\iota})\subset\mathsf{k}(A)\cup\mathsf{k}(\iota))$.

\item\label{lem:rank4}
If $\mathrm{rk}(A)=\omega \cdot\alpha$, then either $A\equiv (\exists x\in\mathsf{L}_{\alpha}B(x))$ or 
$A\equiv (\forall x\in\mathsf{L}_{\alpha}B(x))$.
\end{enumerate}

\eprp

\bdf\label{df:fml_collpase}
{\rm

Let $\rho\in\Psi_{\mS}$ for a successor stable ordinal $\mS$, and
$\sfk(\iota)\subset M_{\rho}$ for $RS$-terms and $RS$-formulas $\iota$.
Then
$\iota^{[\rho/\mS]}$ denotes the result of replacing each unbounded quantifier
$Qx$ by $Qx\in {\sf L}_{\mI_{N}[\rho]}$,
and each ordinal term $\alp\in {\sf k}(\iota)$
by $\alp[\rho/\mS]$ for the Mostowski collapse in
Definition \ref{df:notationsystem.2M}.
$\iota^{[\rho/\mathbb{S}]}$ is defined recursively as follows.

 \begin{enumerate}
 \item\label{df:fml_collpase1}
 $({\sf L}_{\alpha})^{[\rho/\mathbb{S}]}\equiv {\sf L}_{\alpha[\rho/\mathbb{S}]}$ with $\alpha\in M_{\rho}$.
 When $\{\alpha\}\cup\bigcup\{{\sf k}(u_{i}):i\leq n\}\subset M_{\rho}$,
 $\left([x\in {\sf L}_{\alpha}:\phi^{{\sf L}_{\alpha}}(x,u_{1},\ldots,u_{n})]\right)^{[\rho/\mathbb{S}]}$
 is defined to be the $RS$-term
 $[x\in {\sf L}_{\alpha[\rho/\mathbb{S}]}:\phi^{L_{\alpha[\rho/\mathbb{S}]]}}(x, (u_{1})^{[\rho/\mathbb{S}]},\ldots,(u_{n})^{[\rho/\mathbb{S}]})]$.
 
 \item\label{df:fml_collpase2}
 $(\lnot A)^{[\rho/\mathbb{S}]}\equiv\lnot A^{[\rho/\mathbb{S}]}$.
 $(u\in v)^{[\rho/\mathbb{S}]}\equiv\left(u^{[\rho/\mathbb{S}]}\in v^{[\rho/\mathbb{S}]}\right)$.
 \\
 $(A_{0}\lor \cdots\lor A_{n-1})^{[\rho/\mathbb{S}]}\equiv 
 \left(
 (A_{0})^{[\rho/\mathbb{S}]}\lor\cdots\lor(A_{n-1})^{[\rho/\mathbb{S}]}
 \right)$.
 $(\exists x\in u\, A)^{[\rho/\mathbb{S}]}\equiv (\exists x\in u^{[\rho/\mathbb{S}]} A^{[\rho/\mathbb{S}]})$.
 $(\exists x A)^{[\rho/\mathbb{S}]}\equiv(\exists x\in {\sf L}_{\mathbb{I}_{N}[\rho]} A^{[\rho/\mathbb{S}]})$.

\item\label{df:fml_collpase3}
 \[
 (st_{i}^{[\bar{\rho}]}(u))^{[\rho/\mS]}\equiv 
 \left\{
 \begin{array}{ll}
 st_{i}^{[\bar{\rho}]*[\rho]}(u^{[\rho/\mS]}) & \mbox{if } |u|\geq\mS
 \\
  st_{i}^{[\bar{\rho}]}(u) & \mbox{if } |u|<\mS
  \end{array}
  \right.
 \]
 \end{enumerate}

}
\edf

The following Propositions \ref{prp:levelcollaps}, \ref{lem:assigncollaps}
and \ref{lem:rank_M} are seen from Lemma \ref{lem:Mostowskicollaps}.

\bprp\label{prp:levelcollaps}
Let $\rho\in\Psi_{\mS}$.
\benu
\item\label{prp:levelcollaps1}
Let $v$ be an $RS$-term with $\sfk(v)\subset M_{\rho}$,
and $\alp=|v|$.
Then 
$v^{[\rho/\mS]}$ is an $RS$-term of level $\alp[\rho/\mS]$,
$\left| v^{[\rho/\mS]} \right|=\alp[\rho/\mS]$ and
$\sfk(v^{[\rho/\mS]})=\left(\sfk(v)\right)^{[\rho/\mS]}$.

\item\label{prp:levelcollaps2}
Let $\alp\leq\mI_{N}$ be such that $\alp\in M_{\rho}$. Then
$
\left(Tm(\alp)\right)^{[\rho/\mS]}
:= 
\{v^{[\rho/\mS]}: v\in Tm(\alp), \sfk(v)\subset M_{\rho}\}
= Tm(\alp[\rho/\mS])
$.

\item\label{prp:levelcollaps3}
Let $A$ be an $RS$-formula with $\sfk(A)\subset M_{\rho}$.
Then $A^{[\rho/\mS]}$ is an $RS$-formula such that
$\sfk(A^{[\rho/\mS]})\subset\{\alp[\rho/\mS]: \alp\in\sfk(A)\}\cup\{\mI_{N}[\rho]\}\subset\calh_{\mS}[\sfk(A)\cup\{\rho\}]$.
\eenu
\eprp
\bprf
\ref{prp:levelcollaps}.\ref{prp:levelcollaps1}.
We see easily that 
$v^{[\rho/\mathbb{S}]}$ is an $RS$-term of level $\alpha[\rho/\mathbb{S}]$.
\\
\ref{prp:levelcollaps}.\ref{prp:levelcollaps2}.
We see 
$\left(Tm(\alpha)\right)^{[\rho/\mathbb{S}]}\subset Tm(\alpha[\rho/\mathbb{S}])$ from 
Proposition \ref{prp:levelcollaps}.\ref{prp:levelcollaps1}.
Conversely let $u$ be an $RS$-term with ${\sf k}(u)=\{\beta_{i}: i<n\}$ 
and $\max\{\beta_{i}: i<n\}=|u|<\alpha[\rho/\mathbb{S}]$.
By 
Lemma \ref{lem:Mostowskicollaps}
there are ordinal terms $\gamma_{i}\in OT(\mI_{N})$ such that
$\gamma_{i}\in M_{\rho}$ and
$\gamma_{i}[\rho/\mathbb{S}]=\beta_{i}$.
Let $v$ be an $RS$-term obtained from $u$ by replacing each constant ${\sf L}_{\beta_{i}}$ by 
${\sf L}_{\gamma_{i}}$.
We obtain $v^{[\rho/\mathbb{S}]}\equiv u$, $v\in Tm(\alpha)$, and ${\sf k}(v)=\{\gamma_{i}: i<n\}\subset M_{\rho}$.
This means $v\in\left(Tm(\alpha)\right)^{[\rho/\mathbb{S}]}$.
\\
\ref{prp:levelcollaps}.\ref{prp:levelcollaps3}.
We see readily that $\sfk(A^{[\rho/\mS]})\subset\{\alp[\rho/\mS]: \alp\in\sfk(A)\}\cup\{\mI_{N}[\rho]\}$.
From this and Proposition \ref{prp:EH_Pi11}.\ref{prp:EH_Pi11.4}
$\sfk(A^{[\rho/\mS]})\subset\calh_{\mS}[\sfk(A)\cup\{\rho\}]$ follows.
\eprf

\bprp\label{lem:assigncollaps}
For $RS$-formulas $A$ with
$A\simeq \bigvee(A_{\iota})_{\iota\in J}$, assume 
$\sfk(A)\subset M_{\rho}$ with $\rho\in\Psi_{\mS}$.
Then
$A^{[\rho/\mS]}\simeq \bigvee\left((A_{\iota})^{[\rho/\mS]}\right)_{\iota\in [\rho]J}$
for $[\rho]J=\{\iota\in J : \sfk(\iota)\subset M_{\rho}\}$.
\eprp
\bprf
This is seen from Proposition \ref{prp:levelcollaps}.\ref{prp:levelcollaps2}.
Let $A$ be a formula $st_{i}^{[\bar{\rho}]}(u)$ such that $\mathsf{k}(st_{i}^{[\bar{\rho}]}(u))=\mathsf{k}(u)\subset M_{\rho}$.
Then $st_{i}^{[\bar{\rho}]}(u)\simeq \bigvee( {\sf L}_{\mS[\bar{\rho}]}=u)_{{\sf L}_{\mS[\bar{\rho}]}\in J_{i}}$ 
with 
$J_{i}=\{{\sf L}_{\mS[\bar{\rho}]}: \mS\in SSt_{i}\cap M_{[\bar{\rho}]},\, 
|u|\geq\mS[\bar{\rho}] \}$.
If $|u|<\mS$, then $|u|<\rho$, $(st_{i}^{[\bar{\rho}]}(u))^{[\rho/\mS]}\equiv st_{i}^{[\bar{\rho}]}(u)$, $[\rho]J_{i}=J_{i}$ and
$({\sf L}_{\mS[\bar{\rho}]}=u)^{[\rho/\mS]}\equiv ({\sf L}_{\mS[\bar{\rho}]}=u)$ by $u^{[\rho/\mS]}=u$.
Let $|u|\geq\mS$. We have
$(st_{i}^{[\bar{\rho}]}(u))^{[\rho/\mS]}\equiv st_{i}^{[\bar{\rho}]*[\rho]}(u^{[\rho/\mS]})\simeq \bigvee( {\sf L}_{\mS([\bar{\rho}]*[\rho])}=u^{[\rho/\mS]}
)_{{\sf L}_{\mS([\bar{\rho}]*[\rho])}\in I}$ 
with 
$I=\{{\sf L}_{\mS([\bar{\rho}]*[\rho])}: \mS\in SSt_{i}\cap M_{[\bar{\rho}]*[\rho]},\, 
|u^{[\rho/\mS]}|\geq\mS([\bar{\rho}]* [\rho])\}$.
For ${\sf L}_{\mS[\bar{\rho}]}\in J_{i}$, ${\sf L}_{\mS[\bar{\rho}]}\in [\rho]J_{i}$ iff $\mathsf{k}({\sf L}_{\mS[\bar{\rho}]})=\{\mS[\bar{\rho}]\}\subset M_{\rho}$
iff $\mS\in M_{[\bar{\rho}]*[\rho]}$,
and $( {\sf L}_{\mS[\bar{\rho}]}=u)^{[\rho/\mS]}\equiv ({\sf L}_{(\mS[\bar{\rho}])[\rho/\mS]}=u^{[\rho/\mS]})$.
On the other hand we have $\mS([\bar{\rho}]*[\rho])=(\mS[\bar{\rho}])[\rho/\mS]$, and $|u|\geq\mS[\bar{\rho}]$ iff 
$|u^{[\rho/\mS]}|=|u|[\rho/\mS] \geq(\mS[\bar{\rho}])[\rho/\mS]=\mS([\bar{\rho}]* [\rho])$ for $\mS\in M_{[\bar{\rho}]*[\rho]}$.
\eprf

\bprp\label{lem:rank_M}
Let $\sfk(\iota)\subset M_{\rho}$ with $\rho\in\Psi_{\mS}$.

Then
$\rk(\iota^{[\rho/\mS]})=\left(\rk(\iota)\right)[\rho/\mS]$.
\eprp
\bprf
We see that $\mbox{{\rm rk}}(\iota)\in M_{\rho}$
from Proposition \ref{lem:rank}.\ref{lem:rank1}.
The proposition is seen from 
the facts
$(\omega\alpha)[\rho/\mathbb{S}]=\omega(\alpha[\rho/\mathbb{S}])$ and $(\alpha+1)[\rho/\mathbb{S}]=\alpha[\rho/\mathbb{S}]+1$ when $\alpha\in M_{\rho}$.
\eprf

\subsection{A preview of elimination procedures of stable ordinals}\label{subsec:preview}

Let us explain briefly our elimination procedures of stable ordinals in this section and
section \ref{sec:proofonestep}.
We define two derivability relations
$(\mathcal{H}_{\gamma},\Theta; \mathtt{D})\vdash^{* a}_{c}\Gamma$
and 
$(\calh_{\gam},\Tht,\mathtt{Q})\vdash^{a}_{c,\vec{d},\xi,\Lambda}\Gam$
in subsections \ref{subsec:operatorcont} and \ref{subsec:capderivation}, resp.
In the former relation,
$c$ is a bound on ranks of the formulas for stability and of cut formulas as well as
stable ordinals collected in $\mathtt{D}$.

The r\^{o}le of the former calculus $\vdash^{* a}_{c}$ is twofold:
first finite proof figures are embedded in the calculus.
Here is the \textit{axiom} for the stability $L_{\mS}\prec_{\Sig_{1}}L$ of 
successor stable ordinals $\mS$
\begin{equation}\label{eq:axiom_stability}
\lnot B(u), \exists x\in \mathsf{L}_{\mS} B(x)
\end{equation}
where $B$ is a $\Delta_{0}$-formula 
possibly parameters from $L_{\mS}$.

Each formula in it is of rank$<\mI_{N}$, and
the cut rank $c$ in $\vdash^{* a}_{c}$ is lowered to $\mI_{N}+1$.
Then the derivation is collapsed down to a $\beta<\mI_{N}$ using the collapsing function
$\psi_{\mI_{N}}(\alp)$.
For Collapsing \ref{lem:Kcollpase.1} we need
a standard requirement (\ref{eq:controlder*1})
in operator controlled derivations.
Namely $\{\gamma,a,c\}
\cup \mathsf{k}(\Gamma)\cup\mathtt{D}
\subset\mathcal{H}_{\gamma}[
\Theta
]$ when $(\mathcal{H}_{\gamma},\Theta; \mathtt{D})\vdash^{* a}_{c}\Gamma$.
Moreover $\mS\in\mathtt{D}$ when $(\mathcal{H}_{\gamma},\Theta;\mathtt{D})\vdash^{*a}_{c}\lnot B(u), \exists x\in \mathsf{L}_{\mS} B(x)$.

The axiom (\ref{eq:axiom_stability}) says that
$\mS$ `reflects' $\Pi_{\mI_{N}}$-formulas in transfinite levels for a bigger ordinal $\mI_{N}>\mS$ 
such that $L=L_{\mI_{N}}$.
In order to analyze the reflections,
Mahlo classes $Mh^{a}_{i,c}(\xi)$ are introduced in Definition \ref{df:Cpsiregularsm}.\ref{df:Cpsiregularsm.2}.
$\pi\in Mh^{a}_{i,c}(\xi)$
reflects every fact 
$\pi\in Mh_{i,0}^{a}(g_{c})=\bigcap\{Mh^{a}_{i,d}(g(d)): c>d\in {\rm supp}(g)\}$ 
on the ordinals $\pi\in Mh_{i,c}^{a}(\xi)$ in lower level,
down to
`smaller' Mahlo classes $Mh_{i,c}^{a}(f)=\bigcap\{Mh^{a}_{i,d}(f(d)): c\leq d\in {\rm supp}(f)\}$.

This apparatus would suffice to analyze reflections in transfinite levels.
We need another for the axiom $L_{\mS}\prec_{\Sig_{1}}L$, i.e.,
a (formal) \textit{Mostowski collapsing}:
Assume that $B(u,v)$ with $v\in L_{\mS}$ for a $\Del_{0}$-formula $B$.
We need to find a substitute $u^{\prime}\in L_{\mS}$ for $u\in L$ such that
$B(u^{\prime},v)$.
For simplicity let us assume that $v=\bet<\mS$ and $u=\alp$ are ordinals.
We may assume that $\alp\geq\mS$.
Let $\rho<\mS$ be an ordinal, which is bigger than every ordinal$<\mS$
occurring in the `context' of $B(\alp,\bet)$.
This means that $\del<\rho$ holds for every ordinal $\del<\mS$ occurring in a `relevant' branch of a
derivation of $B(\alp,\bet)$.
Then we can define a Mostwosiki collapsing $\alp\mapsto\alp[\rho/\mS]$ for
ordinal terms $\alp$ such that
$\bet[\rho/\mS]=\bet$ for each relevant $\bet<\mS$ and $\mS[\rho/\mS]=\rho$, cf.\,Definition \ref{df:notationsystem.2M}.
Then we see that $B(\alp[\rho/\mS],\bet)$ holds.

Let $M_{\rho}$ be the set in (\ref{eq:M}) of Definition \ref{df:p0}.\ref{df:p0.0}.
Roughly $\alpha\in M_{\rho}$ if
every subterm $\bet<\mS$ of $\alp$ is smaller than $\rho$.
It is shown in Lemma \ref{prp:EK2}.\ref{prp:EK2.1} that
$\mathcal{H}_{\gamma}[M_{\rho}]\subset M_{\rho}$
if $\mathcal{H}_{\gamma}[\rho]\cap\mathbb{S}\subset\rho$.
Let $\calh_{\gam}[\Tht]\vdash^{a}_{c}\Gamma$, and
assume that $\{\gam,a,c\}\cup\sfk(\Gam)\subset\calh_{\gam}[\Tht]$.
Moreover let us assume that $\Tht\subset M_{\rho}$ holds.
Then we obtain
$\{\gam,a,c\}\cup\sfk(\Gam)\subset\calh_{\gam}[\Tht]\subset\calh_{\gam}[M_{\rho}]\subset M_{\rho}$.
This means that $\sfk(\Gam)\subset M_{\rho}$ holds
as long as $\Tht\subset M_{\rho}$ holds, i.e., as long as
we are concerned with branches for $\sfk(\iota)\subset M_{\rho}$
in, e.g., inferences $(\bigwedge)$: $A\simeq\bigwedge(A_{\iota})_{\iota\in J}$
\beqn\label{eq:preview1}
\infer[(\bigwedge)]{\calh_{\gam}[\Tht]\vdash^{a}_{c}\Gam,A}
{
\{
\calh_{\gam}[\Tht]\vdash^{a_{0}}_{c}\Gam,A, A_{\iota}
\}_{\iota\in J}
}
\leadsto
\infer[(\bigwedge)]{\calh_{\gam}[\Tht]\vdash^{a}_{c}\Gam,A}
{
\{
\calh_{\gam}[\Tht]\vdash^{a_{0}}_{c}\Gam,A, A_{\iota}
\}_{\iota\in J,\, \sfk(\iota)\subset M_{\rho}}
}
\eeqn
and dually $\sfk(\iota)\subset M_{\rho}$ for a minor formula $A_{\iota}$
of a $(\bigvee)$ with the main formula $A\simeq\bigvee(A_{\iota})_{\iota\in J}$,
\textit{provided that} 
$\mathcal{H}_{\gamma}[\rho]\cap\mathbb{S}\subset\rho$.
The proviso means that $\gam_{1}>\gam$ when $\rho=\psi_{\mS}^{f}(\gam_{1})$.
Such a $\rho$ is in $\calh_{\gam}[\Tht]$ only when $\rho\in\Tht$.

Let us try to replace the axiom (\ref{eq:axiom_stability})
by inferences for reflection of $\rho$ with $\Tht\subset M_{\rho}$:
If $B(u)^{[\rho/\mS]}$ holds, then $B(u)^{[\sig/\mS]}$ holds for some $\sig<\rho$.
{\small
\[
\hspace{-2mm}
\infer[({\rm rfl})]{\mathcal{H}_{\gamma}[\Theta\cup\{\rho\}]
\vdash
\Gamma^{[\rho/\mS]}
}
{
\mathcal{H}_{\gamma}[\Theta\cup\{\rho\}]
\vdash
\Gamma^{[\rho/\mS]}, B(u)^{[\rho/\mS]}
&
\{
\mathcal{H}_{\gamma}[\Theta\cup\{\rho,\sigma\}]
\vdash
\Gamma^{[\rho/\mS]}, \lnot B(u)^{[\sigma/\mathbb{S}]}
\}_{\Theta\subset M_{\sigma},\sig<\rho}
}
\]
}
In analyzing the inferences for reflections in transfinite levels,
formulas $\Gam^{[\rho/\mS]}$ are replaced by $\Gam^{[\sig/\mS]}$.
This means that $\alp[\sig/\mS]$ is substituted for each $\alp[\rho/\mS]$.
Namely a composition of uncollapsing and collapsing $\alp[\rho/\mS]\mapsto\alp\mapsto\alp[\sig/\mS]$
arises.
Hence we need $\alp\in M_{\sig}\subsetneq M_{\rho}$ for $\sig<\rho$.
However we have $\sig\not\in M_{\sig}$ although $\sig\in M_{\rho}$, and we cannot replace $[\rho/\mS]$ by $[\sig/\mS]$
in the upper part of $\Gamma^{[\rho/\mS]}, B(u)^{[\rho/\mS]}$.
The schema seems to be broken.

Instead of an explicit collapsing ${}^{[\rho/\mS]}$, formulas could put on \textit{caps} $\rho,\sig,\ldots$
in such a way that $\sfk(A^{(\sigma)})=\sfk(A)$ and $\mathrm{rk}(A^{(\sigma)})=\mathrm{rk}(A)$.
This means that the cap $\sig$ does not `occur' in a \textit{capped} formula $A^{(\sig)}$.
Caps are temporary ones. 
A capped formula $A^{(\rho)}$ might put on another cap $\sigma$ as $A^{(\sigma)}$ as long as
$\mathsf{k}(A)\subset M_{\rho}\cap M_{\sigma}$.

Each cap $\rho$ exceeds a \textit{threshold} $\gamma_{1}$ in the sense that
$\rho\not\in\calh_{\gam_{1}}[\rho]\cap\mS\subset\rho$.
Then every ordinal `occurring' in derivations is in the domain $M_{\rho}$
of the Mostowski collapsing $\alp\mapsto\alp[\rho/\mS]$.
The ordinal $\gam_{1}$ is a threshold, which means that
every ordinal occurring in the relevant branches for $\rho$ is in $C_{\gam_{1}}(0)$ and
the subscript $\gam<\gam_{1}$ in $\calh_{\gam}$, while
each $\rho\in\mathtt{Q}$ for a finite set $\mathtt{Q}$ of ordinals, exceeds $\gam_{1}$ in such a way that
${\tt p}_{0}(\rho)>\gam_{1}$ for the ordinal $\mathtt{p}_{0}(\rho)$ in Definition \ref{df:p0}.\ref{df:p0.1}.
This ensures us that $\calh_{\gam}[M_{\rho}]\subset M_{\rho}$.
In the end, inferences for reflections are removed 
by moving outside $C_{\gam_{1}}(0)$.

For a single stable ordinal,
a big enough (depending on a given finite proof figure) ordinal $\gamma_{1}$ suffices to be a threshold.
Now we have several stable ordinals $\mS,\mT,\ldots\in \mathrm{dom}(\mathtt{Q})$
for a \textit{finite} collection $\mathrm{dom}(\mathtt{Q})$ of stable ordinals, cf.\,Definition \ref{df:ffthreshold}.\ref{df:QJ.1}.
Axioms for stability and their children for reflections are eliminated
first for bigger $\mS>\mT$, and then smaller ones $\mT$.
Therefore we need an assignment $\mathrm{dom}(\mathtt{Q})\ni\mS\mapsto\gam_{\mS}^{\mathtt{Q}}$ for 
thresholds so that
$\gam_{\mS}^{\mathtt{Q}}<\gam_{\mT}^{\mathtt{Q}}$ if $\mS>\mT$ in Definition \ref{df:ffthreshold}.\ref{df:caphat.42}.
Furthermore each formula $A$ puts on a list $\bar{\rho}=(\rho_{0},\rho_{1},\ldots,\rho_{n-1})$ of ordinals $\bar{\rho}(\mS_{i})=\rho_{i}\in\Psi_{\mS_{i}}$ 
for each  $\mS_{i}\in\mathrm{dom}(\mathtt{Q})$, which we call a \textit{cap} over $\mathtt{Q}$ since a single formula might be related to several collapsing
processes $[\rho_{i}/\mS_{i}]$.
A capped formula is formally a pair $(A,\bar{\rho})$ of a formula $A$ and a cap $\bar{\rho}$, denoted by $A^{(\bar{\rho})}$.
We assume that $\mathsf{k}(A)\subset M_{\bar{\rho}}=\bigcap_{i<n}M_{\rho_{i}}$.

Let $B$ be a $\Delta_{0}$-formula
possibly parameters from $L_{\mS}$ such that $L\models \exists x\, B(x)$.
Let $\mT>\mS$ be a successor stable ordinal.
Then we obtain $L_{\mT}\models \exists x\, B(x)$ by $L_{\mT}\prec_{\Sig_{1}}L$,
and $L_{\mS}\prec_{\Sig_{1}}L_{\mT}$ yields $L_{\mS}\prec_{\Sigma_{1}}L$.
Therefore we may assume $u\in L_{\mT}$, and
$\mathrm{rk}(B(u))<\mT$ for such stable ordinals $\mT$ in the axiom (\ref{eq:axiom_stability}), 
cf.\,(\ref{eq:stbl_I_1}).

Let $L=\bigcup_{n<\omega}L_{\mS_{n}}$ be a universe which is a limit of stable ordinals $\mS_{0}<\mS_{1}<\cdots$.
We see easily that $\Sigma_{1}\mathrm{-Separation}$ holds in such universes $L$, and
it suffices to assume $\forall n<\omega(L_{\mS_{n}}\prec_{\Sigma_{1}}L_{\mS_{n+1}})$ for
the axiom (\ref{eq:axiom_stability}).
However this is too weak to have (\ref{eq:axiom_stability}) when 
$\mS_{\omega}=\sup_{n<\omega}\mS_{n}\in L$.
\\

In the second derivability relation
$(\calh_{\gam},\Tht,\mathtt{Q})\vdash^{a}_{c,\vec{d},\xi,\Lambda}\Gam$
of subsection \ref{subsec:capderivation}, 
an operator $\calh_{\gam}$ together with a finite set $\Tht$ of ordinals and
a finite family $\mathtt{Q}\subset\coprod_{\mS}\Psi_{\mS}$ controls
ordinals occurring in derivations, where
$\mathrm{dom}(\mathtt{Q})$ is a finite set of stable ordinals $\mS$ and
$\mathtt{Q}(\mS)$ is a finite set of ordinals $\rho\in\Psi_{\mS}$ for each $\mS\in \mathrm{dom}(\mathtt{Q})$.
Furthermore $\mathtt{Q}$ carries thresholds.
The first subscript $c$ is a bound on ranks of cut formulas, while 
the third one $\xi$ is a bound on stable ordinals.
The fourth $\Lambda$ is a bound of ordinals in $\mathrm{fld}(g)$ for finite functions $g$
to be concerned.

We have an inference $(i {\rm -rfl}_{\mS}(\rho))$ in the second calculus to derive the axiom (\ref{eq:axiom_stability}):
{\small
\[
\infer[(i {\rm -rfl}_{\mS}(\rho))]{
(\mathcal{H}_{\gamma},\Theta,\mathtt{Q})\vdash^{a}\Gamma
}
{
\{
(\mathcal{H}_{\gamma},\Theta,\mathtt{Q})\vdash^{a_{0}}\Gamma,\lnot\delta^{(\bar{\rho})}
\}_{\delta\in\Delta}
&
\{
(\mathcal{H}_{\gamma},\Theta,\mathtt{Q}^{\sigma})\vdash^{a_{0}}
\Delta^{(\bar{\rho}_{\sigma})},\Gamma
\}_{\sigma<\rho, \Theta\cap M_{\rho}\subset M_{\sigma}}
}
\]
}
where $\mS$ is a successor $i$-stable ordinal
$\mS\in \mathrm{dom}(\mathtt{Q})\cap \xi$, $\rho\in\mathtt{Q}(\mS)$,
$a_{0}<a$,
$\mathtt{Q}^{\sigma}=\mathtt{Q}\cup\{(\mS,\sigma)\}$,
 $\Delta$ is
a set of $\mathcal{L}_{i}$-formulas $\delta$, and
$\bar{\rho}_{\sigma}$ is obtained from the list $\bar{\rho}$ by replacing a component $\rho=\bar{\rho}(\mS)$ by $\sigma$.
The second subscript $\vec{d}=(d_{1},\ldots,d_{N})$ in
$(\mathcal{H}_{\gamma},\Theta,\mathtt{Q})\vdash^{a}_{c,\vec{d},\xi,\Lambda}\Gamma$ is a list of ordinals $d_{i}$ which is a bound of ranks of formulas $\delta\in\Delta$.
Thus $\mathrm{rk}(\Delta)=\max\{\mathrm{rk}(\delta):\delta\in \Delta\}<d_{i}$.

For the time being let us consider the simplest case $N=1$.
Then $\vec{d}$ is a single ordinal $d$.

In Lemma \ref{lem:capping} (Capping), the first calculus is embedded in the second in such a way that
if $
(\mathcal{H}_{b_{0}},\Theta;\mathtt{D}
)\vdash^{* a}_{c}\Gamma
$, then
$(\mathcal{H}_{\gamma},\Theta, \mathtt{Q})
\vdash^{c+2a}_{c,c,c,\Gamma(c)}
\Gamma^{(\bar{\rho})}$ for a $\gamma>b_{0}$, where each formula $A\in\Gamma$ puts on the same cap $\bar{\rho}$.

The axiom (\ref{eq:axiom_stability}) is derived in the second calculus as follows.
\[
\infer[({\rm rfl})]{\lnot B(u)^{(\bar{\rho})},(\exists x\in\mathsf{L}_{\mathbb{S}}B(x))^{(\bar{\rho})} }
 {
  \infer*[D_{0}]{\lnot B(u)^{(\bar{\rho})},B(u)^{(\bar{\rho})} }{}
  &
  \infer[(\bigvee)]{\lnot B(u)^{(\bar{\rho}_{\sigma})},(\exists x\in\mathsf{L}_{\mathbb{S}}B(x))^{(\bar{\rho})} }
  {
    \infer*[D_{1}]{ \lnot B(u)^{(\bar{\rho}_{\sigma})}, (B(u)^{[\sigma/\mathbb{S}]})^{(\bar{\rho})} }{}
   }
  }
  \]
where $D_{0}$ denotes a derivation of a tautology
$\lnot A^{(\bar{\rho})}, A^{(\bar{\rho})}$ with $A\equiv B(u)$, while $D_{1}$ is a derivation of another `tautology'
$\lnot A^{(\bar{\rho}_{\sigma})}, (A^{[\sigma/\mS]})^{(\bar{\rho})}$ assuming
$\mathsf{k}(A)\subset M_{\bar{\rho}}$, and $\sigma<\rho$ varies through ordinals such that $\mathsf{k}(A)\subset M_{\sigma}$.
Let $A\simeq\bigwedge(A_{\iota})_{\iota\in J}$.
The branches of inferences $(\bigwedge)$ with the capped major formula $A^{(\bar{\rho})}$ 
are restricted to $\iota\in J$ such that $\mathsf{k}(\iota)\subset M_{\bar{\rho}}$ in
(\ref{eq:preview1}),
cf.\,Tautology \ref{lem:tautology.cap}.\ref{lem:tautology.cap1}.

Since the stable ordinal $\mS$ is arbitrary, this means that the formula $A\simeq\bigwedge(A_{\iota})_{\iota\in J}$ is derived from
formulas $A_{\iota}$ such that $\mathsf{k}(\iota)\subset M_{\mathtt{Q}}=\bigcap\{ M_{\tau}: \tau\in\bigcup_{\mT}\mathtt{Q}(\mT)\}$.
$M_{\sigma}\subset M_{\rho}$ holds for $\{\sigma<\rho\}\subset\mathtt{Q}(\mS)$, cf.\,Definition \ref{df:ffthreshold}.\ref{df:QJ.2}, and
 this amounts to
$\mathsf{k}(\iota)\subset M_{\rho}$ for each least $\rho_{\mathtt{Q}(\mS)}=\min \mathtt{Q}(\mS)$.

Next let us examine a derivation of 
$(\mathcal{H}_{\gamma},\Theta,\mathtt{Q}^{\sigma})\vdash^{a}\lnot A^{(\bar{\rho}_{\sigma})}, (A^{[\sigma/\mS]})^{(\bar{\rho})}$, where 
$\sigma<\rho$ and
$\mathtt{Q}^{\sigma}(\mS)=\mathtt{Q}(\mS)\cup\{\sigma\}$.
We obtain 
$A^{[\sigma/\mS]}\simeq\bigwedge(A_{\iota}^{[\sigma/\mS]})_{\iota\in [\sigma]J}$ by Proposition \ref{lem:assigncollaps}.
Let $B\equiv A^{[\sigma/\mS]}$ with $B\simeq\bigwedge(B_{\nu})_{\nu\in I}$ and $B_{\nu}\equiv A_{\iota}^{[\sigma/\mS]}$
for $I=\{\iota^{[\sigma/\mS]}: \iota\in[\sigma]J\}$ and $\nu=\iota^{[\sigma/\mS]}$.
Let $\iota\in[\sigma]J$ be such that $|\iota|\geq\mS$.
Then $\nu=\iota^{[\sigma/\mS]}\neq\iota$, $\mS>|\iota^{[\sigma/\mS]}|\geq\sigma=\mS[\sigma/\mS]$ and 
$\mathsf{k}(\iota^{[\sigma/\mS]})\not\subset M_{\sigma}$ for $\sigma=\min\mathtt{Q}^{\sigma}(\mS)$.
Therefore the capped formula $(A^{[\sigma/\mS]})^{(\bar{\rho})}$ needs to be derived from formulas
$B_{\nu}$ with $\mathsf{k}(\nu)\not\subset M_{\sigma}$ when the cap $\rho$ is not the smallest one in $\mathtt{Q}^{\sigma}(\mS)$
even if $\rho$ is the least one in $\mathtt{Q}(\mS)$, i.e., $\rho=\rho_{\mathtt{Q}(\mS)}>\rho_{\mathtt{Q}^{\sigma}(\mS)}=\sigma$.

This leads us to the \textit{inverse} $\mathrm{inv}((A^{[\sigma/\mS]})^{(\bar{\rho})};\mathtt{Q})\equiv A$ of
a capped formula $(A^{[\sigma/\mS]})^{(\bar{\rho})}$ with respect to a finite family $\mathtt{Q}$ in Definition \ref{df:inverse},
which is obtained from $A^{[\sigma/\mS]}$ by an uncollapsing $[\sigma/\mS]^{-1}$ in Definition \ref{df:divide}.
$\mathrm{inv}((A^{[\sigma/\mS]})^{(\bar{\rho})};\mathtt{Q})\equiv A$ holds only when $\bar{\rho}(\mS)>\sigma\in\mathtt{Q}(\mS)$.
$\mathsf{k}(A^{(\bar{\rho})};\mathtt{Q})$ in Definition \ref{df:Qmin}.\ref{df:Qmin.1} is defined to be the set of ordinals
occurring in the inverse $\mathrm{inv}(A^{(\bar{\rho})};\mathtt{Q})$, and
subsets $[\bar{\rho}\mathtt{Q}]J\subset J$ are introduced in Definition \ref{df:Qmin}.\ref{df:Qmin.3}.
A standard requirement $\mathsf{k}(A)\subset\mathcal{H}_{\gamma}[\Theta]$ in controlled derivations is replaced by
$\mathsf{k}(A^{(\bar{\rho})};\mathtt{Q})\subset\mathcal{H}_{\gamma}[\Theta\cap M_{\mathtt{Q}}]$
when $(\mathcal{H}_{\gamma},\Theta,\mathtt{Q})\vdash\Gamma$ holds with $A^{(\bar{\rho})}\in\Gamma$.

A capped $\bigwedge$-formula $A^{(\bar{\rho})}$ is derived from formulas $A_{\iota}^{(\bar{\rho})}$ for $\iota\in[\bar{\rho}\mathtt{Q}]J$
in the context $\mathtt{Q}$, where $\mathsf{k}(A_{\iota}^{(\bar{\rho})};\mathtt{Q})\subset M_{\mathtt{Q}}$ if $\iota\in[\bar{\rho}\mathtt{Q}]J$.
An inference introducing a $\bigwedge$-formula $A\simeq\bigwedge(A_{\iota})_{\iota\in J}$ with a cap $\bar{\rho}$ then runs as follows:
\[
\infer[(\bigwedge)]{(\mathcal{H}_{\gamma},\Theta,\mathtt{Q})\vdash\Gamma,A^{(\bar{\rho})}}
{
\{(\mathcal{H}_{\gamma},\Theta\cup \mathsf{k}(A_{\iota}^{(\bar{\rho})};\mathtt{Q}),\mathtt{Q})\vdash\Gamma,A_{\iota}^{(\bar{\rho})}\}_{\iota\in[\bar{\rho}\mathtt{Q}]J}
}
\]

For $\mathsf{k}(\iota)\subset M_{\sigma}$, $\Theta_{\iota}=\Theta\cup\mathsf{k}(\iota)$, $d_{\iota}=\mathrm{rk}(A_{\iota})$ and $d=\mathrm{rk}(A)$ we obtain, cf.\,Tautology \ref{lem:tautology.cap}.\ref{lem:tautology.cap3}.
\[
\infer[(\bigwedge)]{(\mathcal{H}_{\gamma},\Theta,\mathtt{Q}^{\sigma})\vdash^{2d}\lnot A^{(\bar{\rho}_{\sigma})}, (A^{[\sigma/\mS]})^{(\bar{\rho})}}
{
 \infer[(\bigvee)]{(\mathcal{H}_{\gamma},\Theta_{\iota},\mathtt{Q}^{\sigma})\vdash^{2d_{\iota}+1}\lnot A^{(\bar{\rho}_{\sigma})}, (A_{\iota}^{[\sigma/\mS]})^{(\bar{\rho})}
 }
 {
 \infer*{(\mathcal{H}_{\gamma},\Theta_{\iota},\mathtt{Q}^{\sigma})\vdash^{2d_{\iota}}\lnot A_{\iota}^{(\bar{\rho}_{\sigma})}, (A_{\iota}^{[\sigma/\mS]})^{(\bar{\rho})}}{}
 }
}
\]
provided that $\mathrm{inv}((A^{[\sigma/\mS]})^{(\bar{\rho})};\mathtt{Q}^{\sigma})\equiv A\equiv \mathrm{inv}(A^{(\bar{\rho}_{\sigma})};\mathtt{Q}^{\sigma})$.

Next let us consider a cut-elimination procedure in the second calculus.
We need to lower the ranks $\mathrm{rk}(A^{(\bar{\rho})})=\mathrm{rk}(A)$ of capped formula $A^{(\bar{\rho})}$ down to
$\mathrm{rk}(A^{(\bar{\rho})})\leq\mS$ in eliminating the largest stable ordinal $\mS$ in the domain $\mathrm{dom}(\mathtt{Q})$.
Note that we may assume that $\mathrm{rk}(\Delta)<\mS$ for 
the set $\Delta$ in inferences $(\mathrm{rfl}_{\mT}(\tau))$ when $\mT<\mS$.

We have an inference $(\mathtt{D})$ in the first calculus $\vdash^{*}$ at which a stable ordinal $\mS$ is eliminated from the collection $\mathtt{D}$:
\[
\infer[(\mathtt{D})]{(\mathcal{H}_{\gamma},\Theta; \mathtt{D}
)\vdash^{* a}_{c}\Gamma
}
{
(\mathcal{H}_{\gamma},\Theta; \mathtt{D}\cup\{\mS\}
)\vdash^{* a_{0}}_{c}\Gamma
}
\]
where $\mS\in \mathcal{H}_{\gamma}[\Theta]\cap c$ and $a_{0}<a$.
Corresponding to it, there is an inference $(\mathrm{dom})$ in the second calculus $\vdash$:
\[
\infer[(\mathrm{dom})]{
(\mathcal{H}_{\gamma},\Theta, \mathtt{Q}
)\vdash\Gamma
}
{
(\mathcal{H}_{\gamma},\Theta, \mathtt{P}
)\vdash \Gamma(\mathtt{P})
}
\]
where $\mathrm{dom}(\mathtt{P})=\mathrm{dom}(\mathtt{Q})\cup\{\mS\}$ and the cap $\bar{\rho}$ of each $A^{(\bar{\rho})}\in\Gamma$
is extended to $\bar{\rho}*(\rho_{\mathtt{P}(\mS)})$ in $\Gamma(\mathtt{P})$.
An extension $\mathtt{P}=\mathtt{Q}[\mS]$ by a stable ordinal $\mS$ carries a threshold
$\gamma^{\mathtt{P}}_{\mS}$.

Due to this inference $(\mathrm{dom})$, an inversion lemma seems not to hold.
An inversion lemma for the calculus $\vdash$ would state that  'if
$(\mathcal{H}_{\gamma},\Theta,\mathtt{Q})\vdash^{a}\Gamma,A^{(\bar{\rho})}$ with $A\simeq\bigwedge(A_{\iota})_{\iota\in J}$, then
for each $\iota\in[\bar{\rho}\mathtt{Q}]J$,
$(\mathcal{H}_{\gamma},\Theta\cup\mathsf{k}(A_{\iota}^{(\bar{\rho})}; \mathtt{Q});\mathtt{Q})\vdash^{a}\Gamma,A_{\iota}^{(\bar{\rho})}$ holds'.
Consider the following trunk of a derivation:
\[
\infer[(\mathrm{dom})]{
(\mathcal{H}_{\gamma},\Theta,\mathtt{Q})\vdash\Gamma,A^{(\bar{\rho})}
}
{
 \infer[(\bigwedge)]{
 (\mathcal{H}_{\gamma},\Theta,\mathtt{P})\vdash\Gamma(\mathtt{P}),A^{(\bar{\rho})}(\mathtt{P})
 }
 {
 \{
  (\mathcal{H}_{\gamma},\Theta\cup\mathsf{k}(A_{\iota}^{(\bar{\rho})}(\mathtt{P}); \mathtt{P}),\mathtt{P})\vdash\Gamma(\mathtt{P}),A_{\iota}^{(\bar{\rho})}(\mathtt{P})
  \}_{\iota\in[\bar{\rho}(\mathtt{P})\mathtt{P}]J}
 }
}
\]
where $\bar{\rho}(\mathtt{P})=\bar{\rho}*(\rho_{\mathtt{P}(\mS)})$.
$[\bar{\rho}(\mathtt{P})\mathtt{P}]J\subsetneq[\bar{\rho}\mathtt{Q}]J$ may be the case.
For example let $\mathrm{dom}(\mathtt{Q})=\{\mT\}$ with $\mathtt{Q}(\mT)=\{\rho\}$,
and $\mathtt{P}(\mS)=\{\sigma\}$ for an ordinal $\sigma\in\Psi_{\mS}$ and $\mT>\mS$.
Even 
if $\mathsf{k}(\iota)\subset M_{\rho}$,
$\mathsf{k}(\iota)\not\subset M_{\sigma}$ is possible, e.g., when $\sigma\leq |\iota|<\mS<\rho$.

The lack of an inversion lemma is not a serious problem to prove a cut-elimination theorem \ref{lem:CE}
provided that inferences $(cut)$ are interchangeable with $(\mathrm{dom})$:
\[
\infer[(cut)]{(\mathcal{H}_{\gamma},\Theta, \mathtt{Q})\vdash\Gamma
}
{
 \infer*[D]{(\mathcal{H}_{\gamma},\Theta, \mathtt{Q})\vdash\Gamma,\lnot C^{(\bar{\rho})}}{}
&
\infer[(\mathrm{dom})]{
(\mathcal{H}_{\gamma},\Theta, \mathtt{Q})\vdash C^{(\bar{\rho})},\Gamma
 }
 {
 (\mathcal{H}_{\gamma},\Theta, \mathtt{P})\vdash C^{(\bar{\rho})}(\mathtt{P}),\Gamma(\mathtt{P})
 }
}
\]
turns to
\[
\infer[(\mathrm{dom})]{(\mathcal{H}_{\gamma},\Theta, \mathtt{Q})\vdash\Gamma
}
{
 \infer[(cut)]{(\mathcal{H}_{\gamma},\Theta, \mathtt{P})\vdash\Gamma(\mathtt{P})
 }
 {
 \infer*[D*\mathtt{P}]{(\mathcal{H}_{\gamma},\Theta, \mathtt{P})\vdash\Gamma(\mathtt{P}),\lnot C^{(\bar{\rho})}(\mathtt{P})}{}
&
(\mathcal{H}_{\gamma},\Theta, \mathtt{P})\vdash C^{(\bar{\rho})}(\mathtt{P}),\Gamma(\mathtt{P})
}
}
\]
where  a derivation $D*\mathtt{P}$ is obtained from the derivation $D$ by 
enlarging $\mathtt{Q}$ to $\mathtt{P}=\mathtt{Q}[\mS]$.

Consider inferences $(\mathrm{dom})$ in the derivation $D$:
\[
\infer[(\mathrm{dom})]{
(\mathcal{H}_{\gamma},\Theta, \mathtt{Q})\vdash \Phi
 }
 {
 (\mathcal{H}_{\gamma},\Theta, \mathtt{R})\vdash \Phi(\mathtt{R})
 }
 \]
where $\mathtt{R}=\mathtt{Q}[\mT]$ is an extension with a threshold $\gamma^{\mathtt{R}}_{\mT}$.
If $\mT\neq\mS$, then we may assume 
$\mathtt{P}[\mT]=(\mathtt{Q}[\mS])[\mT]=(\mathtt{Q}[\mT])[\mS]=\mathtt{R}[\mS]$, and the inference 
turns to
\[
\infer[(\mathrm{dom})]{
(\mathcal{H}_{\gamma},\Theta, \mathtt{P})\vdash \Phi(\mathtt{P})
 }
 {
 (\mathcal{H}_{\gamma},\Theta, \mathtt{P}[\mT])\vdash \Phi(\mathtt{P}[\mT])
 }
 \]
Let $\mT=\mS$. Since 
the threshold $\gamma^{\mathtt{R}}_{\mS}$ need not to be equal to a given ordinal
$\gamma^{\mathtt{P}}_{\mS}$, 
inferences $(\mathrm{dom})$ should be
\[
\infer[(\mathrm{dom})]{
(\mathcal{H}_{\gamma},\Theta, \mathtt{Q})\vdash \Phi
 }
 {
 \{
 (\mathcal{H}_{\gamma},\Theta, \mathtt{R})\vdash \Phi(\mathtt{R})
 \}_{\mathtt{R}}
 }
 \]
where $\mathtt{R}=\mathtt{Q}[\mT]$ varies through extensions with thresholds 
$\gamma^{\mathtt{R}}_{\mT}$ in the upper sequents.
Then we choose a branch for $\gamma^{\mathtt{R}}_{\mS}=\gamma^{\mathtt{P}}_{\mS}$ when
$\mT=\mS$, cf.\,\textbf{Case 1} in the proof of Lemma \ref{lem:prereduction_I_1}.

Here is another problem in adding a stable ordinal $\mT$ to derivations, 
cf.\,\textbf{Case 3} of Lemma \ref{lem:prereduction_I_1}.
Consider an inference $({\rm rfl}_{\mS}(\rho))$:
\[
\infer[({\rm rfl}_{\mS}(\rho))]{
(\mathcal{H}_{\gamma},\Theta,\mathtt{Q})\vdash\Gamma
}
{
\{
(\mathcal{H}_{\gamma},\Theta,\mathtt{Q})\vdash\Gamma,\lnot\delta^{(\bar{\rho})}
\}_{\delta\in\Delta}
&
\{
(\mathcal{H}_{\gamma},\Theta,\mathtt{Q}^{\sigma})\vdash
\Delta^{(\bar{\rho}_{\sigma})},\Gamma
\}_{\sigma<\rho, \Theta\cap M_{\rho}\subset M_{\sigma}}
}
\]
where $\mathrm{rk}(\Delta)<\mU$ for every $\mS<\mU\in\mathrm{dom}(\mathtt{Q})$, cf.\,(\ref{eq:stbl_I_1})
and (r1) in Definition \ref{df:controldercap}.
When we add a successor stable ordinal $\mS<\mT\leq\mathrm{rk}(\Delta)$, 
we have $\mathrm{rk}(\Delta)<\mU$ for every 
$\mT<\mU\in\mathrm{dom}(\mathtt{P})=\mathrm{dom}(\mathtt{Q})\cup\{\mT\}$, and
the inference turns to an inference $({\rm rfl}_{\mT}(\rho_{\mT}))$ for $\rho_{\mT}\in\mathtt{P}(\mT)$.
For this, there are upper sequents of inferences $({\rm rfl}_{\mS}(\rho))$ for each
extension $\mathtt{P}=\mathtt{Q}[\mT]$ by successor stable ordinals $\mT$ such that
$\mS\leq\mT\leq\mathrm{rk}(\Delta)$.
\\

Lemma  \ref{lem:recapping} (Recapping) is one of the main lemma, which states the following:
For finite functions $g$, let $s(g)=\max(\{0\}\cup \mathrm{supp}(g))$.
Then we see that $\mathrm{rk}(\Delta)<s(g)$ for sets $\Delta$ in inferences $({\rm rfl}_{\mS}(\rho))$.
Let $\mS=\max(\mathrm{dom}(\mathtt{Q})\cap SSt)<\xi$ be the largest successor stable ordinal 
in the collection $\mathrm{dom}(\mathtt{Q})$, and $\rho\in\mathtt{Q}(\mS)$.
Let $b$ be an ordinal such that $\mS+1\leq b<s(g)$ with $g=m(\rho)$ and $c\leq b\leq d$.
Assume $(\mathcal{H}_{\gamma},\Theta,\mathtt{Q})\vdash^{a}_{c,d,\xi,\Lambda}\Gamma$
and $\mathrm{rk}(\Gamma)<b$.
Let $\bar{\rho}$ be a cap such that its $\mS$-th entry is $\rho$, i.e., $\bar{\rho}(\mS)=\rho$, and
$\kappa<\rho$ be an ordinal such that $m(\kappa)\geq g_{1}$ for a finite function $g_{1}$
obtained from $g$ by `stepping-down' to $b$, i.e., $s(g_{1})=b$, cf.\,Definition \ref{df:hstepdown}.
Let $\bar{\tau}=\bar{\rho}[\kappa/\rho]$ be a cap such that $\bar{\tau}(\mS)=\kappa$, and
$\mathtt{R}=\mathtt{Q}^{(\kappa/\rho)}$ a finite family such that 
$\mathtt{R}(\mS)$ is obtained from $\mathtt{Q}(\mS)$ by replacing $\rho$ by $\kappa$.
Then $(\mathcal{H}_{\gamma},\Theta,\mathtt{R})\vdash^{a_{1}}_{c,d,\xi,\Lambda}\Gamma^{(\kappa/\rho)}$
holds for an ordinal $a_{1}$, where $\Gamma^{(\kappa/\rho)}$ is obtained from
$\Gamma$ by replacing each $A^{(\bar{\rho})}$ by $A^{(\bar{\tau})}$.

Note that $\mS=\max(\mathrm{dom}(\mathtt{Q})\cap SSt)$ is a condition holding only locally, and
$\mS<\max(\mathrm{dom}(\mathtt{P})\cap SSt)$ happens if $\mathtt{P}=\mathtt{Q}[\mV]$ is
an extension by an ordinal $\mV>\mS$. 
When $\xi$ is larger than a limit stable ordinal (non-projectible ordinal) above $\mS$,
there are unboundedly many such successor stable ordinals $\mV$.
Therefore each of conditions $\mS=\max(\mathrm{dom}(\mathtt{Q})\cap SSt)$ and $\mathrm{rk}(\Gamma)<b$
has to be relaxed.
Let $\Delta$ be a set in inferences $({\rm rfl}_{\mV}(\lambda))$ with $\mV\neq\mS$.
Recall that $\mathtt{Q}$ is extended to $\mathtt{P}=\mathtt{Q}[\mT]$ by ordinals $\mT$ such that
$\mV\leq\mT\leq\mathrm{rk}(\Delta)$.
We then see $\mT\neq\mS$ as follows. We may assume $\mV<\mS$.
Then we obtain $\mathrm{rk}(\Delta)<\mS$ by (r1).
In inferences $({\rm rfl}_{\mV}(\lambda))$ the $\mT$-th entry in caps in the left
upper sequents is replaced by smaller ordinals in the right upper sequents.
This means that the $\mT$-th entry is a main entry, and the $\mS$-th one $\rho$ is a side one
in $({\rm rfl}_{\mV}(\lambda))$ by $\mT\neq\mS$.
Inferences $({\rm rfl}_{\mV}(\lambda))$ are replaced by another $({\rm rfl}_{\mV}(\lambda))$
by recapping $\rho$ by $\kappa$.
The condition (r1) is crucial in \textbf{Case 1.2} and \textbf{Case 1.3.2} in the proof of 
Recapping \ref{lem:recapping}.

For $\mT\in St\cup\{0\}$, let $b=\mT+1$ if $\mT\in SSt\cup\{0\}$, and $b=\mT$ else.
Assume $(\mathcal{H}_{\gamma},\Theta,\mathtt{Q})\vdash^{a}_{b,d,\mT^{\dagger},\Lambda}\Gamma$.
In Lemma \ref{mlem:singlemainl_I_1} the second subscript $d$, which is a bound on
ranks of sets $\Delta$ in inferences $({\rm rfl}_{\mV}(\lambda))$, is shown to be lowered to $b$
by an iterated use of Recapping \ref{lem:recapping}.
It is proved by induction on an ordinal $on_{\Lambda}(g)<\Gamma(\Lambda)$
for finite functions $g:\mI_{N}\to\Gamma(\mI_{N})$ such that
ordinals in $\mathrm{fld}(g)$ are bounded by the fourth subscript $\Lambda$, cf.\,Definition \ref{df:on_f}.
Thus $\hat{a}>\Lambda$ for an ordinal height $\hat{a}$ of the resulting derivation, and
the fourth subscript becomes a larger one $\Lambda_{1}>\Lambda$.
$g$ is a function such that $m(\rho)\geq g$ for $\rho=\rho_{\mathtt{Q}(\mS)}=\min\mathtt{Q}(\mS)$ and
$\mS=\max(\mathrm{dom}(\mathtt{Q})\cap SSt)<\mT^{\dagger}$.
When $\mathtt{P}=\mathtt{Q}[\mV]$ is extended by an ordinal $\mV>\mS$, we need to consider another
function $k$ for $\tau\in\mathtt{P}(\mV)$ such that $k\leq m(\tau)$.
For a proof by induction to work, we introduce another thresholds $\delta^{\mathtt{Q}}_{\mS}$
so that $\delta^{\mathtt{Q}}_{\mS}=\delta^{\mathtt{P}}_{\mS}>\delta^{\mathtt{P}}_{\mV}$
if $\mS<\mV$, cf.\,Definition \ref{df:ffthreshold}.\ref{df:caphat.52_I_1}.
Lemma \ref{mlem:singlemainl_I_1} is shown by main induction on the thresholds 
$\delta^{\mathtt{Q}}_{\mS}$ with subsidiary induction on ordinals $on_{\Lambda}(g)$ (and with sub-subsidiary induction on $a$).
Therefore we have to specify a function $g$ from an ordinal $\rho$ or a given family $\mathtt{Q}$.
Although $g=m(\rho)$ is an intended one, a function $g$ is not uniquely determined 
from an ordinal $\rho$,
cf.\,(\ref{eq:KppiNlowerCase1b.12_I_1}) in \textbf{Case 1.3.2} of Recapping \ref{lem:recapping}.
For this \textit{guarded ordinals} are introduced in Definition \ref{df:guarded_ord_I_1}.
A guarded ordinal is a pair $(\rho,g)$ of an ordinal $\rho$ and a finite function $g$
such that $m(\rho)\geq g$.
Then for example, $\mathtt{Q}(\mS)$ is defined to be a finite set of guarded ordinals
$(\rho,g)$ with $\rho\prec\mS$.
\\

In Lemmas \ref{lem:main.2_L_I_1} and \ref{lem:main.2_S_I_1}
 we eliminate a maximal stable ordinal $\mS$ from derivations.
Lemma \ref{lem:main.2_L_I_1} for limit stable ordinals $\mS$ immediately follows
from Lemma \ref{mlem:singlemainl_I_1}.
When $\mS$ is a successor stable ordinal, we need another Lemma \ref{mlem:singlemainl_S_I_1}.
Suppose that
the ranks of $\Delta$ in the right upper sequent of
inferences $({\rm rfl}_{\mS}(\rho))$ are lowered to $\mathrm{rk}(\Delta)\leq\mS$.
Then the largest successor stable ordinal $\mS$ is removed
from derivations in the following way:
Let $(\mathcal{H}_{\gamma},\Theta,\mathtt{Q})\vdash_{\mS+1,\mS+1,\mS^{\dagger},\Lambda}\Gamma$.
Let
$\gamma_{1}=\gamma^{\mathtt{Q}}_{\mS}+\mS$
and
$\mathtt{Q}\restrict\mS$ denote the restriction of $\mathtt{Q}$ to $\mS$.
Then
$(\mathcal{H}_{\gamma_{1}},\Theta\cup\mathtt{Q}(\mS),\mathtt{Q}\restrict\mS)\vdash_{\mS,\mS,\mS,\Lambda}\Gamma^{[\mS]}$ holds,
where 
$\Gamma^{[\mS]}$ is obtained from $\Gamma$ by replacing each $A^{(\bar{\rho})}\in\Gamma$ by
a formula $B^{(\bar{\rho}\,\restrict\,\mS)}$.
$\bar{\rho}\restrict\mS$ is the restriction of $\bar{\rho}$ to $\mS$, and a cap over $\mathtt{Q}\restrict\mS$, while $B$ denotes either
the collapse $A^{[\rho_{\mathtt{Q}(\mS)}/\mS]}$ of $A$ or $A^{[\bar{\rho}(\mS)/\mS]}$ roughly speaking, cf.\,Definition \ref{df:restrict}.

In \textbf{Case 2} of the proof of Lemma \ref{mlem:singlemainl_S_I_1},
an inference $({\rm rfl}_{\mS}(\rho))$ is replaced by a series of $(cut)$'s as follows.
Pick an ordinal $\sigma\in \mathcal{H}_{\gamma_{1}}[\Theta\cup\mathtt{Q}(\mS)]\cap\rho$ such that $(\mathcal{H}_{\gamma},\Theta,\mathtt{Q})\vdash_{\mS+1,\mS+1,\mS^{\dagger},\Lambda}\Delta^{(\bar{\rho}_{\sigma})},\Gamma$ and $\mathrm{rk}(\Delta)\leq\mS$.
We obtain 
$(\mathcal{H}_{\gamma_{1}},\Theta\cup\mathtt{Q}(\mS),\mathtt{Q}\restrict\mS)\vdash_{\mS,\mS,\mS,\Lambda}(\Delta^{(\bar{\rho}_{\sigma})})^{[\mS]},\Gamma^{[\mS]}$.
Assume for simplicity that $\Delta=\{\exists x\in\mathsf{L}_{\mS}B(x)\}\subset\Sigma_{1}(\mS)$.
Then $(\Delta^{(\bar{\rho}_{\sigma})})^{[\mS]}=\{(\exists x\in\mathsf{L}_{\sigma}B(x))^{(\bar{\rho}\,\restrict\,\mS)}\}$ for $\sigma=\rho_{\mathtt{Q}^{\sigma}(\mS)}$.
On the other hand we have $(\mathcal{H}_{\gamma},\Theta,\mathtt{Q})\vdash_{\mS+1,\mS+1,\mS^{\dagger},\Lambda}\Gamma,
\lnot(\exists x\in\mathsf{L}_{\mS}B(x))^{(\bar{\rho})}$.
We obtain
$(\mathcal{H}_{\gamma_{1}},\Theta\cup\mathtt{Q}(\mS),\mathtt{Q}\restrict\mS)\vdash_{\mS,\mS,\mS,\Lambda}\Gamma^{[\mS]},
(\lnot\exists x\in\mathsf{L}_{\rho}B(x))^{(\bar{\rho}\,\restrict\,\mS)}$ for $\rho=\rho_{\mathtt{Q}(\mS)}$.
We have $\sigma<\rho$, and 
$\lnot\exists x\in\mathsf{L}_{\rho}B(x)$ yields $\lnot\exists x\in\mathsf{L}_{\sigma}B(x)$.
We obtain
$(\mathcal{H}_{\gamma_{1}},\Theta\cup\mathtt{Q}(\mS),\mathtt{Q}\restrict\mS)\vdash_{\mS,\mS,\mS,\Lambda}\Gamma^{[\mS]}$
by a $(cut)$ with $\mathrm{rk}(\exists x\in\mathsf{L}_{\sigma}B(x))<\mS$.
In \textbf{Case 1.2.1} of Lemma \ref{mlem:singlemainl_S_I_1}, an inference $(\bigvee)$
deriving $\exists x\in\mathsf{L}_{\mS}B(x)$ from $B(u^{[\sigma/\mathbb{S}]})$
turns to
\[
  \infer[(\bigvee)]{\Gamma^{[\mS]},(\exists x\in\mathsf{L}_{\rho}B(x))^{(\bar{\rho}\,\restrict\,\mS)} }
  {
   \Gamma^{[\mS]}, B(u^{[\sigma/\mathbb{S}]})^{(\bar{\rho}\,\restrict\,\mS)} 
   }
\]
where $\alpha[\sigma/\mathbb{S}]<\rho$ for any $\alpha\in M_{\sigma}$.

In Lemma \ref{lem:main.1_I_1} stable ordinals are removed from derivations as follows.
Suppose
$
(\mathcal{H}_{\gamma},\Theta, \mathtt{Q})
\vdash_{\xi,\xi,\xi,\Lambda}
\Gamma
$ and $\mathrm{rk}(\Gamma)<\mW^{\dagger}$,
where 
both of $\mW$ and $\xi$ are stable ordinals such that $\mW\in\mathrm{dom}(\mathtt{Q})\cap\mathcal{H}_{\gamma}[\Theta\cap M_{\mathtt{Q}}]$.
Then there is a set $\Theta\subset\Theta_{\mW^{\dagger}*}$  and a recapping $\mathtt{Q}_{\mW}$ 
of $\mathtt{Q}\restrict\mW^{\dagger}$ such that
$
(\mathcal{H}_{\gamma_{\mW}},\Theta_{\mW^{\dagger}*}, \mathtt{Q}_{\mW}
)
\vdash_{\mW^{\dagger},\mW^{\dagger},\mW^{\dagger},\Lambda_{1}}
\Gamma^{[\mW*]}
$ holds for an ordinal $\gamma_{\mW}$, where
in $\Gamma^{[\mW^{\dagger}*]}$, each $A^{(\bar{\rho})}\in\Gamma$ is replaced by 
a formula $B^{(\bar{\rho}\,\restrict\,\mW^{\dagger})}$, and $B$ is obtained from $A$ by applying collapsings $[\rho_{0}/\mT_{0}]$, $[\rho_{1}/\mT_{1}],\ldots$
repeatedly, where 
$\{\mT_{0}>\mT_{1}>\cdots\}=\{\mT\in\mathrm{dom}(\mathtt{Q}):\mT\geq\mW^{\dagger}\}$.

In proving Lemma \ref{lem:main.1_I_1} we need $\mT\in\mathrm{dom}(\mathtt{Q})$
for a stable ordinal $\mT$ such that $\mT\leq\mathrm{rk}(C)<\mT^{\dagger}$ when $C\in\Gamma$,
cf.\,(\ref{eq:controlder_cap_cover}), \textbf{Case 2.2} and \textbf{Case 3.2} in the proof of the lemma.
For this, we introduce \textit{stable rank} $\mathrm{srk}(\alpha)$ of ordinals $\alpha$ in
Definition \ref{df:pd_closed_cover}.
For $\alpha<\mI_{N}$, $\mS=\mathrm{srk}(\alpha)\in St\cup\{0\}$ is the ordinal such that
$\mS\leq\alpha<\mS^{\dagger}$.
Lemma \ref{lem:main.1_I_1} is then proved by main induction on a bound $\xi$ of stable ordinals
with subsidiary induction on $a$ using Lemmas \ref{lem:main.2_L_I_1} and \ref{lem:main.2_S_I_1}.

In particular 
$
(\mathcal{H}_{\gamma_{0}},\Theta_{0^{\dagger}*}, \mathtt{Q}_{0}
)
\vdash_{0^{\dagger},0^{\dagger},0^{\dagger},\Lambda_{1}}
\Gamma^{[0^{\dagger} *]}
$ follows for the least stable ordinal $\mS_{0}=0^{\dagger}$,
 where $\mathrm{rk}(\Gamma)<0^{\dagger}$
and $\mathtt{Q}_{0}$ denotes a family with 
$\mathrm{dom}(\mathtt{Q}_{0})\subset\{0\}$.
This means that every formula occurring in the derivation puts on the empty cap $\bar{\rho}$, and is uncapped.
By eliminating uncapped cut formulas and collapsing the derivations down to $<\Omega$,
we obtain Theorem  \ref{thm:2} for the case $N=1$.
\\

Finally let $N>1$.
The axiom for the stability $L_{\mS}\prec_{\Sig_{1}(\mathcal{L}_{i})}L$ of 
successor $i$-stable ordinals $\mS$ is stated as
$\lnot B(u), \exists x\in \mathsf{L}_{\mS} B(x)$, 
where $B$ is a $\Delta_{0}$-formula in the language $\mathcal{L}_{i}=\{\in\}\cup\{st_{j}:j<i\}$ 
possibly parameters from $L_{\mS}$ with a constraint:
$\forall\mT\in\mathtt{D}\cap \bigcup_{j\geq i}SSt_{j}(\mS<\mT \Rarw \mathrm{rk}(B(u))<\mT)$, 
since $L_{\mS}\prec_{\Sig_{1}(\mathcal{L}_{i})}L_{\mT}\prec_{\Sig_{1}(\mathcal{L}_{j})}L$, cf.\,(\ref{eq:stbl_I_1}).
The constraint is restated in the second calculus as:
$\forall\mT\in\mathrm{dom}(\mathtt{Q})\cap \bigcup_{j\geq i}SSt_{j}(\mS<\mT \Rarw \mathrm{rk}(\Delta)<\mT)$, cf.\,(r1) in Definition \ref{df:controldercap}.

When $N>1$, 
the second subscript in the second derivability relation
\\
$(\mathcal{H}_{\gamma},\Theta,\mathtt{Q})\vdash^{a}_{c,\vec{d},\xi,\Lambda}\Gamma$
is a vector $\vec{d}=(d_{1},d_{2},\ldots,d_{N})$, in which $d_{i}$ is a bound on ranks of sets
$\Delta$ of formulas in inferences $(i\mathrm{-rfl}_{\mS}(\rho))$ for $\mS\in SSt_{i}$.
There are several maximal stable ordinals
$\mS_{i}=\max(\mathrm{dom}(\mathtt{Q})\cap \bigcup_{j\geq i}SSt_{j})$ in the
finite collection $\mathrm{dom}(\mathtt{Q})$.
Let us assume $\mS_{N}<\cdots<\mS_{2}<\mS_{1}$.
Then $L_{\mS_{i}}\prec_{\Sigma_{1}(\mathcal{L}_{i})}L$ and 
there is no non-trivial upper bound on $d_{i}$ besides $d_{i}\leq\xi$.

Assume $(\mathcal{H}_{\gamma},\Theta,\mathtt{Q})\vdash^{a}_{b,\vec{d},\mT^{\dagger},\Lambda}\Gamma$.
In Lemma \ref{mlem:singlemainl_I_1} the ordinal $d_{i}$ is shown to be lowered to $b\in\{\mT,\mT+1\}$
provided that $d_{j}=b$ for any $j<i$, cf.\,(\ref{eq:recapping_PiN_0}) in Recapping \ref{lem:recapping}
and \textbf{Case 1.2.1} in the proof of lemma.
This means that $d_{1}$ is first lowered to $b$, and then $d_{2}$ lowered to $b$, and so on.
When we replace an ordinal $\rho\in\Psi_{\mS_{i}}$ by ordinals $\kappa$
 in Recapping \ref{lem:recapping}, there may occur formulas not in the language $\mathcal{L}_{i}$,
cf.\,\textbf{Case 1.2.2} in the proof of the lemma.

Now details follow.

\subsection{Operator controlled $*$-derivations}\label{subsec:operatorcont}
By a \textit{successor stable ordinal} we mean ordinals in $SSt=\bigcup_{0<i\leq N}SSt_{i}$,
and $\mS^{\dagger}:=\mS^{\dagger 1}$.

We define a derivability relation 
$(\mathcal{H}_{\gamma},\Theta; \mathtt{D})\vdash^{* a}_{c}\Gamma$
where $\mathtt{D}$ is a collection of 
stable ordinals, and
$c$ is a bound of ranks of the formulas as well as stable ordinals $\mS$ in the axioms $(i{\rm -stbl}(\mS))$, and one of ranks of cut formulas.

Let $\vec{i}=(i_{0}\geq i_{1}\geq\cdots\geq i_{n})$ be a (possibly empty) sequence
of numbers $i_{m}$ such that $N\geq i_{0}\geq i_{n}=i>0$.
Let $\beta\in (LSt_{k}\cap\Psi)\cup\{0\}$ with $k\geq i$, and
$\mS=\beta^{\dagger\vec{i}}\in St_{i}\cup\{0\}$.
Then let $\mS^{\dagger i}:=\beta^{\dagger\vec{j}}$ for the extended sequence
$\vec{j}=\vec{i}*(i)$.
$\mS^{\dagger i}$ is the least (successor) $i$-stable ordinal above $\mS$, cf.\,Proposition \ref{prp:Stclass}.

\bdf\label{df:pd_closed_cover}
{\rm

\benu
\item
The \textit{stable rank} of ordinals $\alpha$, denoted by $\mathrm{srk}(\alpha)$, is defined as follows.
Let $\mathrm{srk}(\alpha):=\mS$ for $\mS\in St\cup\{0\}$ with $\mS\leq\alpha<\mS^{\dagger}$
when $\alpha<\mI_{N}$.
$\mathrm{srk}(\alpha):=0$ if $\alpha\geq\mI_{N}$.

\item
Let $\mS=\mT^{\dagger \vec{i}}\in St$ with 
a $\mT\in (LSt\cap\Psi)\cup\{0\}$ and a non-increasing sequence 
$\vec{i}=(i_{0}\geq i_{1}\geq\cdots\geq i_{n})$ of numbers.
Then let 
\[
\mathrm{cl}(\mS)=\{\mT^{\dagger \vec{i}_{m}}: 0\leq m\leq n+1\}\cup
\bigcup\{\mathrm{cl}(\mU)\cap\mS: \mT\neq\mU\in\mathrm{trail}(\mT)\}
\]
for initial segments $\vec{i}_{m}=(i_{0}\geq i_{1}\geq\cdots\geq i_{m-1})$ of $\vec{i}$,
cf.\,Proposition \ref{prp:trail} for the trail $\mathrm{trail}(\alpha)$ to ordinals $\alpha$.

When $\alpha\not\in St$, 
let $\mathrm{cl}(\alpha)=\emptyset$.

\item
Let
$\mathrm{srk}(u):=\mathrm{srk}(\mathrm{rk}(u))$ for $RS$-terms $u$, and
$\mathrm{srk}(A):=\mathrm{srk}(\mathrm{rk}(A))$ for $RS$-formulas $A$.

\item
Let $\mathrm{srk}(\Gamma) :=\bigcup\{\mathrm{srk}(A) :A\in\Gamma\}$
for sequents $\Gamma$.

\item
A finite collection $\mathtt{D}$ of ordinals in $St\cup\{0\}$
is said to be \textit{predecessor closed}
if $\mS\in\mathtt{D}$,  then $\mathrm{cl}(\mS)\subset\mathtt{D}$.
\eenu
}
\edf

\bprp\label{prp:pd_closed_cover_I_1}
\benu
\item\label{prp:pd_closed_cover_I_1.1}
$\mathrm{srk}(\alpha)\in\mathcal{H}_{0}[\{\alpha\}]$, and
$\mathrm{cl}(\mS)\subset\mathcal{H}_{0}[\{\mS\}]$ for $\mS\in St\cup\{0\}$.

\item\label{prp:pd_closed_cover_I_1.3}
Let $\mS\in SSt$, and $\alpha\in M_{\rho}$ with $\rho\in\Psi_{\mS}$.
Then $\mathrm{srk}(\alpha[\rho/\mS])\in \mathrm{srk}(\alpha)\cup
(\mathrm{cl}(\mS)\cap\rho)\cup(\{\rho\}\cap St)$.

\item\label{prp:pd_closed_cover_I_1.31}
Let $\mS\in SSt$, $\rho\in\Psi_{\mS}$ and $\mathsf{k}(A)\subset M_{\rho}$ for a formula $A$.
Then
$\mathrm{srk}(A^{[\rho/\mS]})\in \{\mathrm{srk}(A)\}\cup(\mathrm{cl}(\mS)\cap\rho)\cup(\{\rho\}\cap St)$.

\item\label{prp:pd_closed_cover_I_1.4}
Let $A$ be a formula and $\mT$ a stable ordinal.
If $\mathrm{srk}(A)\subset\mT$, then $\mathsf{k}(A)<\mT$.
In particular when $A$ is bounded, $\mathrm{rk}(A)<\mT$ holds.

\eenu
\eprp
\bprf
\ref{prp:pd_closed_cover_I_1}.\ref{prp:pd_closed_cover_I_1.3}.
Let $\alpha\geq\mS\in SSt_{i}$.
We have $\rho\leq \alpha[\rho/\mS]<\rho^{\dagger}$ by Proposition \ref{prp:M_bnd}.
If $i>1$, then $\rho\in LSt_{i-1}\subset St$ and $\mathrm{srk}(\alpha[\rho/\mS])=\rho$.
Let $i=1$. Then $\rho^{\dagger}=\mS$.
Let $\mS=\beta^{\dagger\vec{i}}$,
where $\beta\in (LSt_{k}\cap\Psi)\cup\{0\}$ with $k\geq i$, and
$\vec{i}=(i_{0}\geq i_{1}\geq\cdots\geq i_{n})$ be a non-empty sequence
of numbers $i_{m}$ such that $k\geq i_{0}\geq i_{n}=1$.
Let $\mT=\beta^{\dagger\vec{j}}$ with $\vec{j}=(i_{0},\ldots,i_{n-1})$.
Then $\mT\in St_{1}\cup\{0\}$ and $\mT<\rho<\mT^{\dagger}=\mS$.
Hence $\mathrm{srk}(\alpha[\rho/\mS])=\mT\in\mathrm{cl}(\mS)\cap\rho$.
\\
\ref{prp:pd_closed_cover_I_1}.\ref{prp:pd_closed_cover_I_1.31}.
This is seen from Proposition \ref{prp:pd_closed_cover_I_1}.\ref{prp:pd_closed_cover_I_1.3}.
\eprf

\bdf\label{df:controlder*}
{\rm
Let $\Theta\subset \mI_{N}$ be a finite set of ordinals, $\mathtt{D}\subset St\cup\{0\}$ a predecessor closed set of ordinals,
and
$\gamma,a,c$ ordinals such that $\mathtt{D}\subset c$.

$(\mathcal{H}_{\gamma},\Theta;\mathtt{D})\vdash^{* a}_{c} \Gamma$ holds
for a set $\Gamma$ of formulas
if 

\begin{equation}\label{eq:controlder*1}
\{\gamma,a,c\}
\cup \mathsf{k}(\Gamma)\cup\mathtt{D}
\subset\mathcal{H}_{\gamma}[
\Theta
]
\end{equation}

\begin{equation}\label{eq:controlder*_cover}
\mathrm{srk}(\Gamma)
\subset \mathtt{D}
\end{equation}

and one of the following cases holds:

\begin{description}

\item[$(\bigvee)$]\footnote{The condition $|\iota|< a$ is absent in the inference $(\bigvee)$.}
There exist 
$A\simeq\bigvee\{A_{\iota}: \iota\in J\}$,
$a(\iota)<a$ and an $\iota\in J$ such that $A\in\Gamma$ and
$(\mathcal{H}_{\gamma},\Theta; \mathtt{D})\vdash^{* a(\iota)}_{c}\Gamma,
A_{\iota}$.

\item[$(\bigwedge)$]
There exist 
$A\simeq\bigwedge\{A_{\iota}: \iota\in J\}$, 
ordinals $a(\iota)<a$ 
for each $\iota\in J$ such that
$A\in\Gamma$ and
$(\mathcal{H}_{\gamma},\Theta\cup \mathsf{k}(\iota);  \mathtt{D}\cup\{\mathrm{srk}(A_{\iota})\}
)
\vdash^{* a(\iota)}_{c}\Gamma,
A_{\iota}$.

\item[$(cut)$]
There exist an ordinal $a_{0}<a$ and a formula $C$ 
such that 
$\mathrm{rk}(C)<c$, 
$(\mathcal{H}_{\gamma},\Theta;  \mathtt{D} )\vdash^{* a_{0}}_{c}\Gamma,\lnot C$
and
$(\mathcal{H}_{\gamma},\Theta ;\mathtt{D})\vdash^{* a_{0}}_{c}C, \Gamma$.

\item[$(\Sigma(St)\mbox{{\rm -rfl}})$]
There exist 
$a_{\ell}, a_{r}<a$ and a formula $C
\in\Sigma(\mathcal{L}_{N+1})$ 
such that 
$c>\mI_{N}$,
$(\mathcal{H}_{\gamma},\Theta ;\mathtt{D} )
\vdash^{* a_{\ell}}_{c}\Gamma,C$
and
$(\mathcal{H}_{\gamma},\Theta;  \mathtt{D})
\vdash^{* a_{r}}_{c}
\lnot \exists x\, C^{(x,\mI_{N})}, \Gamma$.

\item[$(\Sigma(\Ome)\mbox{{\rm -rfl}})$]
There exist ordinals
$a_{\ell}, a_{r}<a$ and a formula $C\in\Sigma(\mathcal{L}_{0}:\Ome)$ 
such that 
$(\mathcal{H}_{\gamma},\Theta ;\mathtt{D}
)\vdash^{* a_{\ell}}_{c}C, \Gamma$
and
$(\mathcal{H}_{\gamma},\Theta;  \mathtt{D}
)\vdash^{* a_{r}}_{c}
\lnot \exists x<\Omega\,C^{(x,\Omega)}, \Gamma$, where $c>\Omega$.

\item[$(\mathtt{D})$]
There exist a 
stable ordinal $\mS\in \mathcal{H}_{\gamma}[\Theta]\cap c$ such that $\mS\not\in\mathtt{D}$
and $\mathtt{D}\cup\{\mS\}$ is predecessor closed
for which
$(\mathcal{H}_{\gamma},\Theta;\mathtt{D}\cup\{\mS\})\vdash^{* a_{0}}_{c}\Gamma$ holds for an $a_{0}<a$.

\item[$(i{\rm -stbl}(\mS))$]
Let $0<i\leq N$
and $a>0$.
There exist a successor $i$-stable ordinal $\mS\in \mathtt{D}\cap SSt_{i}$ such that $\mS<c$,
a $\bigwedge$-formula 
$B({\sf L}_{0})\in\mathcal{L}_{i}$
with $\rk(B({\sf L}_{0}))<\mS$, and
a $u\in Tm(\mI_{N})$ such that $\mS\leq\rk(B(u))<c$,
$\{\lnot B(u),\exists x\in\mathsf{L}_{\mathbb{S}}\, B(x)\}\subset\Gamma$ and
\beqn\label{eq:stbl_I_1}
\forall\mT\in\mathtt{D}\cap \bigcup_{j\geq i}SSt_{j}(\mS<\mT \Rarw \mathrm{rk}(B(u))<\mT)
\eeqn

\end{description}
}
\edf

We will state some lemmas for the operator controlled derivations.
These can be shown as in \cite{Buchholz}.

\blem\label{lem:tautology*}{\rm (Tautology)}
$(\mathcal{H}_{0},\Theta\cup \mathsf{k}(A); \{\mathrm{srk}(A)\})
\vdash^{* 2d}_{0}
\lnot A, A$
 holds for
$d=\mathrm{rk}(A)$.
\elem
\bprf
By induction on $d$.
\eprf

\blem\label{lem:equality*}
{\rm (Equality)}
$(\mathcal{H}_{0},\Theta_{uv};\mathtt{D}_{uv})
\vdash^{* \omega(|u|\#|v|)\# 2d}_{0}
u\neq v,\lnot A(u), A(v)$
holds for $d=\mathrm{rk}(A(\mathsf{L}_{0}))$, $\Theta_{uv}=\Theta\cup \mathsf{k}(A,u,v)$
and $\mathtt{D}_{uv}=\mathrm{srk}(A(u),A(v))$, where
$\mathsf{k}(A,u,v)=\mathsf{k}(A)\cup\mathsf{k}(u)\cup\mathsf{k}(v)$ and
$\mathrm{srk}(A(u),A(v))=\{\mathrm{srk}(A(u)), \mathrm{srk}(A(v))\}$.

\elem
\bprf
By induction on $d$, cf.\,\cite{Buchholz,OA}.

Let $\Theta=\mathsf{k}(u,v,w)$ and $\mathtt{D}=\{\mathrm{srk}(A(\iota)):\iota\in\{u,v,w\}\}$.
First show that
$(\mathcal{H}_{\gamma},\Theta;\emptyset
)\vdash^{* \alp}_{0}
u\neq v,u\not\in w,v\in w$,
$(\mathcal{H}_{\gamma},\Theta;\emptyset
)\vdash^{* \alp}_{0}
u\neq v,u\neq w,v= w$
and
$(\mathcal{H}_{\gamma},\Theta;\emptyset
)\vdash^{* \alp}_{0}
u\neq v,w\not\in u, w\in v$
simultaneously by induction on the natural sum $|u|\#|v|\#|w|$, where
$\alp=\ome(|u|\#|v|\#|w|)$.
Then the lemma is seen by induction on $d=\rk(A(\mathsf{L}_{0}))$. 
\eprf

\blem\label{th:embedreg}{\rm (Embedding of Axioms)}
For each axiom $A$ in $S_{\mI_{N}}$ 
there is an $m<\omega$ such that
 $(\mathcal{H}_{0},\mathsf{k}(\Gamma)\cup\mathtt{D};\mathtt{D})
 \vdash^{* \mI_{N}\cdot 2+m}_{\mI_{N}+m} A,\Gamma$
holds for every sequent $\Gamma$ and 
every 
predecessor closed 
set $\mathtt{D}$ such that $\mathrm{srk}(\Gamma)\subset\mathtt{D}$.
\elem
\bprf
In the proof, let us suppress the operator $\calh_{0}$, and
write $\vdash^{*}$ for $\vdash^{* \mI_{N}+m}_{\mI_{N}+m}$ for an $m<\ome$.
We show  $(\mathcal{H}_{0},\emptyset;\emptyset) \vdash^{*} A$ for each axiom $A$ of $S_{\mI_{N}}$  in section \ref{sec:9}.
Adding $\mathsf{k}(\Gamma)$ to the first $\emptyset$, and sequents $\Gamma$ to $A$ in the resulting derivation,
we obtain  $(\mathcal{H}_{0},\mathsf{k}(\Gamma)\cup\mathtt{D};\mathtt{D})\vdash^{*} A,\Gamma$
provided that $\mathrm{srk}(\Gamma)\subset\mathtt{D}$ and $\mathtt{D}$ is predecessor closed.
\\
\textbf{the axiom (\ref{eq:sucstable})}:
Let $\vphi(y)\equiv(\exi x\,\tht(x,y))$ be a $\Sig_{1}(\{st_{j}\}_{j<i})$-formula such that
$\rk(\tht({\sf L}_{0},{\sf L}_{0}))<\ome$.
Also let $u,w$ be $RS$-terms, $\mS$ a successor $i$-stable ordinal, and 
$B(x)\equiv\tht(x,w)$.
We may assume that $B(u)$ is a $\bigwedge$-formula.

Assuming $\{\mathrm{srk}(\theta(\mathsf{L}_{0},w)),\mS\}\subset\mathtt{D}$, we show
\beqn\label{eq:embed.1}
\mathsf{k}(w)\cup\mathtt{D};\mathtt{D} \vdash^{*}
 w\not\in{\sf L}_{\mS},\lnot \exists x\, B(x),\exists x\in {\sf L}_{\mathbb{S}} B(x)
\eeqn
First assume $|w|<\mS$.
Then $\rk(B({\sf L}_{0}))=\rk(\tht({\sf L}_{0},w))<\mS$.

Assuming $\mathrm{srk}(\theta(\mathsf{L}_{0},u))\in\mathtt{D}$ and $\mathsf{k}(B(u))\subset\mathcal{H}_{0}[\Theta]$,
we show by induction on the ordinal $d_{u}=\mathrm{rk}(B(u))$
that
\beqn\label{eq:embed_I_1.1}
\Theta\cup\mathtt{D};\mathtt{D} \vdash^{* 2d_{u}+\omega}_{\mI_{N}}
 \lnot B(u),\exists x\in \mathsf{L}_{\mathbb{S}}B(x)
\eeqn
We may assume that $\mS\leq\mathrm{rk}(B(u))$ for otherwise we have $|u|<\mS$ and
$\Theta\cup\mathtt{D};\mathtt{D} \vdash^{* 2d_{u}}_{\mI_{N}} \lnot B(u), B(u)$ by Tautology \ref{lem:tautology*}, and
(\ref{eq:embed_I_1.1}) follows.

If $\forall\mT\in\mathtt{D}\cap \bigcup_{j\geq i}SSt_{j}(\mS<\mT \Rarw d_{u}<\mT)$, then we obtain (\ref{eq:stbl_I_1}) and
$\Theta\cup\mathtt{D};\mathtt{D}
\vdash^{* 1}_{\mI_{N}}\lnot B(u),\exists x\in \mathsf{L}_{\mathbb{S}}B(x)$ by an $(i{\rm -stbl}(\mS))$.

Otherwise let $\mV=\max\{\mT\in\mathtt{D}\cap \bigcup_{j\geq i}SSt_{j}:\mS<\mT\leq d_{u}\}$
and $\mV\in SSt_{k}$ with $k\geq i$.
We have $\forall\mT\in\mathtt{D}\cap \bigcup_{j\geq k}SSt_{j}(\mV<\mT \Rarw d_{u}<\mT)$, and
a $(k {\rm -stbl}(\mV))$ yields 
\beqn\label{eq:embed_I_1.2}
\Theta\cup\mathtt{D};\mathtt{D} \vdash^{* 1}_{\mI_{N}}
 \lnot B(u),\exists x\in \mathsf{L}_{\mathbb{V}}B(x)
\eeqn
Let $v\in Tm(\mV)$ and $\mathtt{D}_{v}=\mathtt{D}\cup\{\mathrm{srk}(\theta(\mathsf{L}_{0},v))\}$. 
We obtain $\mathrm{rk}(B(v))=d_{v}<\mV\leq\mathrm{rk}(B(u))=d_{u}$.
IH yields 
\[
\Theta\cup\mathsf{k}(v)\cup\mathtt{D};\mathtt{D}_{v} \vdash^{* 2d_{v}+\omega}_{\mI_{N}}
 \lnot B(v),\exists x\in \mathsf{L}_{\mathbb{S}}B(x)
\]

A $(\bigwedge)$ yields
\beqn\label{eq:embed_I_1.3}
\Theta\cup\mathtt{D}; \mathtt{D} \vdash^{* \mV}_{\mI_{N}}
 \lnot \exists x\in \mathsf{L}_{\mathbb{V}}B(x),\exists x\in \mathsf{L}_{\mathbb{S}}B(x)
\eeqn
A $(cut)$ with (\ref{eq:embed_I_1.2}) and (\ref{eq:embed_I_1.3}) yields (\ref{eq:embed_I_1.1}),
where
$\mathrm{rk}(\exists x\in \mathsf{L}_{\mathbb{V}}B(x))=\mV\leq d_{u}<\mI_{N}$.
We obtain (\ref{eq:embed.1}) by (\ref{eq:embed_I_1.1}).

Next assume $|w|\geq\mS$, and let $v\in Tm(\mS)$. Then $|v|<\mS$ and 
$(v\dot{\in}{\sf L}_{\mS})\equiv(v\not\in{\sf L}_{0})$.
We obtain by (\ref{eq:embed.1})
\[
\mathsf{k}(v)\cup\{\mS\}; \{\mathrm{srk}(\theta(\mathsf{L}_{0},v))\} \cup \mathrm{cl}(\mS)\vdash^{*}
 \lnot \exists x\, \tht(x,v),\exists x\in {\sf L}_{\mathbb{S}}\, \tht(x,v)
\]
Let $\mathrm{srk}(v)\cup\mathrm{srk}(w)\cup\{\mS\}\subset\mathtt{D}$. We obtain 
\[
\mathsf{k}(w,v)\cup\mathtt{D};\mathtt{D} \vdash^{* }
\lnot(v\dot{\in}{\sf L}_{\mS}), w\neq v, \lnot \exists x\, \tht(x,w),\exists x\in {\sf L}_{\mathbb{S}}\tht(x,w)
\]
by Equality \ref{lem:equality*} followed by $(cut)$'s with $|v|,|w|<\mI_{N}$ and $\rk(\exists x\,\tht(x,w))=\mI_{N}$.
Then a $(\bigvee)$ followed by a  $(\bigwedge)$ yields (\ref{eq:embed.1}), where
$(w\not\in{\sf L}_{\mS})\simeq\bigwedge(\lnot(v\dot{\in}{\sf L}_{\mS})\lor w\neq v)_{v\in Tm(\mS)}$.

Let $v$ be an $RS$-term with $|v|\geq\mS$, and 
$\{\mathrm{srk}(w=v), \mS\}\subset\mathtt{D}$.
We obtain by (\ref{eq:embed.1}) and Equality \ref{lem:equality*}
\[
\mathsf{k}(w,v)\cup\mathtt{D};\mathtt{D} \vdash^{* }
 {\sf L}_{\mS}\neq v, w\not\in v,\lnot \exists x\, \tht(x,w),\exists x\in v\, \tht(x,w)
\]
We have $\lnot st_{i}(v)\simeq\bigwedge({\sf L}_{\mS}\neq v)_{J}$ with $J=\{{\sf L}_{\mS}: |v|\geq\mS\in SSt_{i}\}$.
A $(\bigwedge)$ yields the axiom (\ref{eq:sucstable}) for 
$\mathrm{srk}(w=v)\in\mathtt{D}$
\[
\mathsf{k}(w,v)\cup\mathtt{D};\mathtt{D} \vdash^{*}
\lnot st_{i}(v),\lnot\vphi(w), w\not\in v, \vphi^{v}(w)
\]
\textbf{the axiom (\ref{eq:LimN})}:
Let $u$ be an $RS$-term
and $\bet=\alp^{\dagger N}$ for $\alp=|u|$.
Then $\bet\in\calh_{0}[\sfk(u)]$ and $\mathrm{cl}(\beta)\subset\mathtt{D}=\{\mathrm{srk}(u=u), \beta\}$.
We obtain 
$\mathsf{k}(u);\mathtt{D} \vdash^{*} u=u $ 
and 
$\mathsf{k}(u);\mathtt{D} \vdash^{*} {\sf L}_{\bet}={\sf L}_{\bet}$.
Hence
\[
\infer[(\bigwedge)]{
\emptyset;\emptyset \vdash^{* }
\fal x \exi y \left(x\in y \land st_{N}(y) \right)
 }
{
 \infer[(\mathtt{D})]{
 \mathsf{k}(u); \mathrm{srk}(u=u) \vdash^{* } \exi y \left(u\in y \land st_{N}(y)\right)}
 {
 \infer[(\bigvee)]{
\mathsf{k}(u); \mathrm{srk}(u=u)\cup\{\beta\} \vdash^{* } \exi y \left(u\in y \land st_{N}(y)\right)
   }
 {
  \infer[(\bigwedge)]{
\mathsf{k}(u);\mathtt{D} \vdash^{* } 
  u\in{\sf L}_{\bet} \land st_{N}({\sf L}_{\bet})
    }
  {
   \infer[(\bigwedge),(\bigvee)]{\mathsf{k}(u);\mathtt{D} \vdash^{* }u\in{\sf L}_{\bet} }
    {
    \mathsf{k}(u);\mathtt{D}\vdash^{*}u\dot{\in}\mathsf{L}_{\beta}
    &
     \mathsf{k}(u);\mathtt{D} \vdash^{*} u=u 
     }
        &
   \infer[(\bigvee)]{\mathsf{k}(u); \mathtt{D} \vdash^{*} st_{N}({\sf L}_{\bet})}
    {\mathsf{k}(u);\mathtt{D}  \vdash^{*} {\sf L}_{\bet}={\sf L}_{\bet}}
 }
}
}
}
\]
\textbf{the axiom (\ref{eq:stbl0})}:
Let $\mT\in SSt_{i+1}$
be a successor $(i+1)$-stable ordinal.
We obtain $\{\mT\} ;\mathrm{cl}(\mT)\vdash^{*}\tht({\sf L}_{\mT})$ for
$\tht(x)\equiv({\sf L}_{\Ome}\in x\land\fal y\in x\fal z\in y(z\in x))$ with ${\sf L}_{\Ome}\equiv \mathsf{M}_{0}$.

For a given $\alp<\mT$ pick a successor $i$-stable ordinal $\alp<\mS<\mT$ such that
$\mS\in\calh_{0}[\{\alp,\mT\}]$ by Proposition \ref{prp:stblS1}.

Let $|v|=\alp<\mT$ and $\mathtt{D}_{v}=\{\mathrm{srk}(\alpha)\}$.
We obtain $({\sf L}_{\mS}\dot{\in} {\sf L}_{\mT})\equiv ({\sf L}_{\mS}\not\in{\sf L}_{0})$, and
$\mathsf{k}(v);\mathtt{D}_{v}\vdash^{*}v=v$ and 
$\mathsf{k}(v)\cup\{\mT\}  ; \mathrm{cl}(\mS) \vdash^{*}{\sf L}_{\mS}={\sf L}_{\mS}$.
Hence 
$\mathsf{k}(v)\cup\{\mT\} ;\mathtt{D}_{v}\cup\mathrm{cl}(\mS) \vdash^{*} v\in{\sf L}_{\mS}\land st_{i}({\sf L}_{\mS})$, and
$\mathsf{k}(v)\cup\{\mT\};\mathtt{D}_{v}\cup\mathrm{cl}(\mT)\cup \mathrm{cl}(\mS)
\vdash^{*} \exi z\in{\sf L}_{\mT}(v\in z\land st_{i}(z))$ by a $(\bigvee)$.
We obtain
$\mathsf{k}(v)\cup\{\mT\};\mathtt{D}_{v}\cup\mathrm{cl}(\mT)
\vdash^{*} \exi z\in{\sf L}_{\mT}(v\in z\land st_{i}(z))$ 
by several $(\mathtt{D})$'s.
Let $w$ and $u$ be $RS$-terms, and $\mathtt{D}_{wu}=\{\mathrm{srk}(w=u)\}$.
Equality \ref{lem:equality*} yields 
$\mathsf{k}(w)\cup\{\mT\}; \mathtt{D}_{wu}\cup\mathrm{cl}(\mT)\vdash^{*} w\not\in{\sf L}_{\mT}, \exi z\in{\sf L}_{\mT}(w\in z\land st_{i}(z))$, and
$\mathsf{k}(w,u)\cup\{\mT\};\mathtt{D}_{wu}\cup\mathrm{cl}(\mT) \vdash^{*} 
u\neq{\sf L}_{\mT}, w\not\in u, \exi z\in u(w\in z\land st_{i}(z))$.
A $(\bigwedge)$ yields
$\mathsf{k}(w,u);\mathtt{D}_{wu} \vdash^{*} 
\lnot st_{i+1}(u), w\not\in u, \exi z\in u(w\in z \land st_{i}(z))$.

$\Del_{0}(\mathcal{L}_{N+1})$-Collection follows from an inference $(\Sig(St){\rm -rfl})$, and
the $\Del_{0}$-collection for the set $M_{0}={\sf L}_{\Ome}$ follows from an inference
$(\Sig(\Ome){\rm -rfl})$.
Other axioms in ${\sf KP}\ome$, i.e., axioms for pair, union, $\Del_{0}$-Separation and foundation
are seen
as in \cite{Buchholz, OA}.
\eprf

\blem\label{th:embedregthm}{\rm (Embedding)}
If $S_{\mI_{N}}\vdash \Gamma$ for sets $\Gamma$ of sentences, 
there are $m,k<\omega$ such that 
 $(\mathcal{H}_{0},\emptyset; \{0\} )\vdash_{\mI_{N}+m}^{* \mI_{N}\cdot 2+k}\Gamma$ 
holds.
\elem
\bprf
This follows from
Lemma \ref{th:embedreg}
as in \cite{Buchholz,OA}.
\eprf

\blem\label{prp:inversion*_I_1}{\rm (Weakening)}
Let $(\mathcal{H}_{\gamma},\Theta; \mathtt{D})\vdash^{* a}_{c}\Gamma$ with $c\geq\mI_{N}$.
\benu
\item\label{prp:inversion*_I_1.1}
Let $A$ be an $RS$-formula.
Let $\mT=\max(\{0\}\cup(\{\mathrm{srk}(A)\}\cap SSt))$,
$\mathrm{cl}(\mT)
=\{\mT_{1}>\cdots>\mT_{n}\}$ and
$\beta=\mT_{1}+\cdots+\mT_{n}$.

Then
$(\mathcal{H}_{\gamma},\Theta\cup\mathsf{k}(A); \mathtt{D}\cup \{\mathrm{srk}(A)\}
)\vdash^{* \beta+a}_{c}\Gamma,A$ holds.

\item\label{prp:inversion*_I_1.2}
$(\mathcal{H}_{\gamma},\Theta; \mathtt{D}\cup \{\mT\}) \vdash^{* \mT\cdot k+a}_{c}\Gamma$ holds for stable ordinals $\mT\in\mathcal{H}_{\gamma}[\Theta]$ such that $\mathtt{D}\cup \{\mT\}$ is predecessor closed,
where $k=1$ if $\mT\in SSt$, and $k=0$ else.
\eenu
\elem
\bprf
We obtain $\mathrm{srk}(A)\in \mathcal{H}_{\gamma}[\mathsf{k}(A)]$.
It suffices to show Lemma \ref{prp:inversion*_I_1}.\ref{prp:inversion*_I_1.2} by induction on $a$.
\\
\textbf{Case 1}.
First consider the axiom $(i{\rm -stbl}(\mS))$:
We have 
an $\mS\in\mathtt{D}\cap SSt_{i}\cap c$,
a formula 
$B(\mathsf{L}_{0})\in \mathcal{L}_{i}$ with $\rk(B(\mathsf{L}_{0}))<\mS$, and
a $u\in Tm(\mI_{N})$ such that $\mS\leq d_{u}=\mathrm{rk}(B(u))<c$,
$\{\lnot B(u),\exists x\in\mathsf{L}_{\mathbb{S}}B(x)\}\subset\Gamma$ and
$\forall\mU\in\mathtt{D}\cap \bigcup_{j\geq i}SSt_{j}(\mS<\mU\Rarw d_{u}<\mU)$, cf.\,(\ref{eq:stbl_I_1}).
We have $(\mathcal{H}_{\gamma},\Theta; \mathtt{D})\vdash^{* a}_{c}\Gamma$ for an $a>0$.
We may assume that $\mT$ is a $k$-successor stable ordinal in $SSt_{k}$ for a $k\geq i$
such that $\mS<\mT\leq d_{u}$.
We obtain $\forall\mU\in\mathtt{D}\cap \bigcup_{j\geq k}SSt_{j}(\mT<\mU\Rarw d_{u}<\mU)$ and by a 
$(k{\rm -stbl}(\mT))$
\beqn\label{eq:inversion*_I_1_1}
(\mathcal{H}_{\gamma},\Theta ;\mathtt{D}\cup\{\mT\})\vdash^{* 1}_{c}
 \lnot B(u),\exists x\in \mathsf{L}_{\mathbb{T}}B(x)
\eeqn
 
Let $v\in Tm(\mT)$ and $\mathrm{srk}(|v|)\leq d_{v}=\mathrm{rk}(B(v))<\mT\leq d_{u}<\mI_{N}\leq c$.
By (\ref{eq:embed_I_1.1}) we obtain
\[
(\mathcal{H}_{\gamma},\Theta\cup\mathsf{k}(v);\mathtt{D}\cup\mathrm{srk}(|v|))
 \vdash^{* 2d_{v}+\omega}_{c}
 \lnot B(v),\exists x\in \mathsf{L}_{\mathbb{S}}B(x)
 \]
A $(\bigwedge)$ yields
\beqn\label{eq:inversion*_I_1_2}
(\mathcal{H}_{\gamma},\Theta;\mathtt{D}\cup\{\mT\})
 \vdash^{* \mT}_{c}
 \lnot \exists x\in \mathsf{L}_{\mathbb{T}}B(x),\exists x\in \mathsf{L}_{\mathbb{S}}B(x)
\eeqn
A $(cut)$ with (\ref{eq:inversion*_I_1_1}) and (\ref{eq:inversion*_I_1_2}) yields
$(\mathcal{H}_{\gamma},\Theta;\mathtt{D}\cup\{\mT\}) \vdash^{* \mT+1}_{c}
 \lnot B(u),\exists x\in \mathsf{L}_{\mathbb{S}}B(x)$, where
 $\mathrm{rk}(\exists x\in \mathsf{L}_{\mathbb{T}} B(x))=\mT<\mI_{N}\leq c$ and $\mT+1\leq\mT+a$.
\\
\textbf{Case 2}.
Next consider the case when the last inference in $(\mathcal{H}_{\gamma},\Theta; \mathtt{D})\vdash^{* a}_{c}\Gamma$ is a $(\mathtt{D})$:
We have $(\mathcal{H}_{\gamma},\Theta; \mathtt{D}\cup\{\mS\})\vdash^{* a_{0}}_{c}\Gamma$ for an $a_{0}<a$ and
a stable ordinal $\mS\in\mathcal{H}_{\gamma}[\Theta]\cap c$ with $\mS\not\in\mathtt{D}$.
We obtain  $(\mathcal{H}_{\gamma},\Theta;  \mathtt{D}\cup\{\mS,\mT\})\vdash^{* \mT\cdot k+a_{0}}_{c}\Gamma$ by IH.
A $(\mathtt{D})$ yields $(\mathcal{H}_{\gamma},\Theta; \mathtt{D}\cup \{\mT\})\vdash^{* \mT\cdot k+a}_{c}\Gamma$ if $\mS\neq\mT$.
\eprf

\blem\label{lem:inversionreg*}{\rm (Inversion)}
Let  $A\simeq \bigwedge(A_{\iota})_{\iota\in J}$ and
$(\mathcal{H}_{\gamma},\Theta;\mathtt{D})\vdash^{* a}_{c}\Gamma,A$.
\benu
\item\label{lem:inversionreg*1}
Let $\mathrm{rk}(A)\geq\mI_{N}$, and $\iota\in J$ be such that $\mathsf{k}(\iota)\subset\mathcal{H}_{\gamma}[\Theta]$ and $\mathrm{srk}(A_{\iota})\in\mathtt{D}$.
Then 
$(\mathcal{H}_{\gamma},\Theta;\mathtt{D}
)
\vdash^{* a}_{c}
\Gamma,A_{\iota}$ holds.

\item\label{lem:inversionreg*2}
Let $A\equiv(\forall x\,B)\in\Pi_{1}(\mI_{N})$, and $\beta\in \mathcal{H}_{\gamma}[\Theta]\cap\mI_{N}$
with $\beta\in LSt_{N}$.
Then
$(\mathcal{H}_{\gamma},\Theta;\mathtt{D}\cup\{\beta\}
)
\vdash^{* a}_{c}
\Gamma,A^{(\beta,\mI_{N})}$ holds, where $A^{(\beta,\mI_{N})}\equiv(\forall x\in\mathsf{L}_{\beta} B)$.
\eenu
\elem
\bprf
By induction on $a$.
Since $\mathrm{rk}(A)\geq\mathbb{I}_{N}$, $A$ is not a major formula of any of $(i{\rm -stbl}(\mS))$.
We obtain $\beta\not\in SSt$ by $\beta\in LSt_{N}$.
\eprf

\blem\label{lem:reduction*}{\rm (Reduction)}
Let $(\mathcal{H}_{\gamma},\Theta;\mathtt{D})\vdash^{*a}_{c}\Gamma_{0},\lnot C$ and
$(\mathcal{H}_{\gamma},\Theta;\mathtt{D})\vdash^{*b}_{c}C,\Gamma_{1}$ for
$C\simeq\bigvee(C_{\iota})_{\iota\in J}$ with $\mathbb{I}_{N}\leq\mathrm{rk}(C)\leq c$, $c>\mI_{N}$ and $a\geq b$.

Then
$(\mathcal{H}_{\gamma},\Theta;\mathtt{D})\vdash^{* \mI_{N}\cdot b+a+b}_{c}\Gamma_{0},\Gamma_{1}$ holds.
\elem
\bprf 
By induction on $b$.
\\
\textbf{Case 1}. $b=0$. Then $(\mathcal{H}_{\gamma},\Theta;\mathtt{D})\vdash^{*b}_{c}C,\Gamma_{1}$ follows from a void $(\bigwedge)$
with a major formula in $\Gamma_{1}$. 
We have 
$(\mathcal{H}_{\gamma},\Theta;\mathtt{D})\vdash^{*0}_{c}\Gamma_{0},\Gamma_{1}$.
Let $b>0$.
\\
\textbf{Case 2}. 
$(\mathcal{H}_{\gamma},\Theta;\mathtt{D})\vdash^{*b}_{c}C,\Gamma_{1}$ follows from an $(i{\rm -stbl}(\mS))$. Since $\mathrm{rk}(C)\geq\mathbb{I}_{N}$, $C$ is not a major formula of any $(i{\rm -stbl}(\mS))$.
We obtain
$(\mathcal{H}_{\gamma},\Theta;\mathtt{D})\vdash^{*b}_{c}\Gamma_{0},\Gamma_{1}$.
\\
\textbf{Case 3}. $\mathrm{rk}(C)<c$: Then a $(cut)$ yields the lemma by $b\leq a<a+b$.
In what follows assume $\mathrm{rk}(C)=c>\mI_{N}$.
\\
\textbf{Case 4}. 
The last inference in $(\mathcal{H}_{\gamma},\Theta; \mathtt{D})\vdash^{*b}_{c}C,\Gamma_{1}$ is a $(\bigwedge)$ with a
major formula $A\in\Gamma_{1}$:
Let $A\simeq\bigwedge(A_{\iota})_{\iota\in J}$.
For each $\iota\in J$ we have $(\mathcal{H}_{\gamma},\Theta\cup\mathsf{k}(\iota);\mathtt{D}\cup\{\mathrm{srk}(A_{\iota})\}
)
\vdash^{* b(\iota)}_{c}C,\Gamma_{1},A_{\iota}$.

On the other hand we have $\mathrm{srk}(A_{\iota})<\mI_{N}$, and
$(\mathcal{H}_{\gamma},\Theta\cup\mathsf{k}(\iota);\mathtt{D}\cup\{\mathrm{srk}(A_{\iota})\})
\vdash^{* \mI_{N}+a}_{c}\Gamma_{0},\lnot C$ by Weakening \ref{prp:inversion*_I_1}.\ref{prp:inversion*_I_1.1}.
IH yields 
$(\mathcal{H}_{\gamma},\Theta\cup\mathsf{k}(\iota);\mathtt{D}\cup\{\mathrm{srk}(A_{\iota})\})
\vdash^{* \mI_{N}\cdot (b(\iota)+1)+a+b(\iota)}_{c}\Gamma_{0},\Gamma_{1},A_{\iota}$.
A $(\bigwedge)$ yields
$(\mathcal{H}_{\gamma},\Theta;\mathtt{D})\vdash^{* \mI_{N}\cdot b+a+b}_{c}\Gamma_{0},\Gamma_{1}$.
\\
\textbf{Case 5}. 
The last inference in $(\mathcal{H}_{\gamma},\Theta;\mathtt{D})\vdash^{*b}_{c}C,\Gamma_{1}$ is a $(\bigvee)$ with the
major formula $C$.
We have $(\mathcal{H}_{\gamma},\Theta;\mathtt{D})\vdash^{*b_{0}}_{c}C_{\iota},C,\Gamma_{1}$ for $\iota\in J$ and $b_{0}<b$.
IH yields $(\mathcal{H}_{\gamma},\Theta;\mathtt{D})\vdash^{* \mI_{N}\cdot b_{0}+a+b_{0}}_{c}C_{\iota},\Gamma_{0},\Gamma_{1}$.
Assuming $\mathsf{k}(\iota)\subset\mathsf{k}(C_{\iota})$, we obtain $\mathsf{k}(\iota)\subset\mathcal{H}_{\gamma}[\Theta]$ by (\ref{eq:controlder*1}), and
$\mathrm{srk}(C_{\iota})\in\mathtt{D}$ by (\ref{eq:controlder*_cover}).
$(\mathcal{H}_{\gamma},\Theta;\mathtt{D})\vdash^{* a}_{c}\Gamma_{0},\lnot C_{\iota}$ 
follows by 
Inversion \ref{lem:inversionreg*}.\ref{lem:inversionreg*1}.
We obtain $\mathrm{rk}(C_{\iota})<\mathrm{rk}(C)$ by Proposition  \ref{lem:rank}.\ref{lem:rank2}, and
$(\mathcal{H}_{\gamma},\Theta;\mathtt{D})\vdash^{*\mI_{N}\cdot b+a+b}_{c}\Gamma_{0},\Gamma_{1}$ by a $(cut)$.
\\
\textbf{Case 6}. 
The last inference in $(\mathcal{H}_{\gamma},\Theta;\mathtt{D})\vdash^{*b}_{c}C,\Gamma_{1}$ is a $(\mathtt{D})$:
We have $(\mathcal{H}_{\gamma},\Theta;\mathtt{D}\cup\{\mS\})\vdash^{*b_{0}}_{c}C,\Gamma_{1}$ for a $b_{0}<b$ and 
$\mS\in\mathcal{H}_{\gamma}[\Theta]$.
We obtain
$(\mathcal{H}_{\gamma},\Theta;\mathtt{D}\cup\{\mS\})\vdash^{*\mI_{N}+a}_{c}\Gamma_{0},\lnot C$ by Weakening \ref{prp:inversion*_I_1}.\ref{prp:inversion*_I_1.2}.

IH yields
$(\mathcal{H}_{\gamma},\Theta;\mathtt{D}\cup\{\mS\})\vdash^{*\mI_{N}\cdot (b_{0}+1)+a+b_{0}}_{c}\Gamma_{0},\Gamma_{1}$, and
$(\mathcal{H}_{\gamma},\Theta;\mathtt{D})\vdash^{*\mI_{N}\cdot b+a+b}_{c}\Gamma_{0},\Gamma_{1}$ by a $(\mathtt{D})$.

Other cases are seen from IH.
\eprf

\blem\label{lem:predcereg*}{\rm (Cut-elimination)}\\
Suppose
$(\mathcal{H}_{\gamma},\Theta;\mathtt{D})\vdash^{* a}_{\mathbb{I}_{N}+1+m}\Gamma$ for $m<\omega$.
Then $(\mathcal{H}_{\gamma},\Theta;\mathtt{D})\vdash^{*\, \omega_{m}(\mI_{N}+a)}_{\mathbb{I}_{N}+1}\Gamma$ holds.
\elem
\bprf
By main induction on $m$ with subsidiary induction on $a$.
It suffices to show $(\mathcal{H}_{\gamma},\Theta;\mathtt{D})\vdash^{*\, \omega^{\mI_{N}+a}}_{\mathbb{I}_{N}+1+m}\Gamma$ assuming
$(\mathcal{H}_{\gamma},\Theta;\mathtt{D})\vdash^{* a}_{\mathbb{I}_{N}+1+m+1}\Gamma$
by induction on $a$.

Let $(\mathcal{H}_{\gamma},\Theta;\mathtt{D})\vdash^{* a_{0}}_{\mathbb{I}_{N}+1+m+1}\Gamma,\lnot C$ and $(\mathcal{H}_{\gamma},\Theta;\mathtt{D})\vdash^{* a_{0}}_{\mathbb{I}_{N}+1+m+1}\Gamma,C$ for $a_{0}<a$ and $\mathrm{rk}(C)<\mI_{N}+1+m+1$.
IH yields
$(\mathcal{H}_{\gamma},\Theta;\mathtt{D})\vdash^{*\, \omega^{\mI_{N}+a_{0}}}_{\mathbb{I}_{N}+1+m}\Gamma,\lnot C$ and
$(\mathcal{H}_{\gamma},\Theta;\mathtt{D})\vdash^{* \,\omega^{\mI_{N}+a_{0}}}_{\mathbb{I}_{N}+1+m}\Gamma,C$.
Reduction \ref{lem:reduction*} then yields
$(\mathcal{H}_{\gamma},\Theta;\mathtt{D})\vdash^{*a_{1}}_{\mathbb{I}_{N}+1+m}\Gamma$
for $a_{1}=\mI_{N}\cdot  \omega^{\mI_{N}+a_{0}}+ \omega^{\mI_{N}+a_{0}}+ \omega^{\mI_{N}+a_{0}}<\omega^{\mI_{N}+a}$.
\eprf

\blem\label{lem:persistency}{\rm ($\Sigma$-persistency)}
Let  
$A\in\Sigma(\mathcal{L}_{N+1})$,
$(\mathcal{H}_{\gamma},\Theta;\mathtt{D})\vdash^{* a}_{c}
\Gamma,A^{(\alpha,\mathbb{I}_{N})}$ and $\alpha<\beta<\mI_{N}$ for $\beta\in LSt_{N}$.
Assume
 that
$\beta\in\mathcal{H}_{\gamma}[\Theta]$.
Then 
$(\mathcal{H}_{\gamma},\Theta;\mathtt{D}\cup\{\beta\})
\vdash^{* a}_{c}
\Gamma,A^{(\beta,\mathbb{I}_{N})}$.

\elem
\bprf
This is seen by induction on $a$.
We obtain $(\mathcal{H}_{\gamma},\Theta;\mathtt{D}\cup\{\beta\})\vdash^{* a}_{c}
\Gamma,A^{(\alpha,\mathbb{I}_{N})}$ by Weakening \ref{prp:inversion*_I_1}.\ref{prp:inversion*_I_1.2}
with $\mathrm{cl}(\beta)=\{\beta\}$.
(\ref{eq:controlder*1}) follows from the assumption.
\eprf

\blem\label{lem:Kcollpase.1}{\rm (Collapsing)}
Let
$
(\mathcal{H}_{\gamma},\Theta;\mathtt{D}
)\vdash^{* a}_{\mI_{N}+1,\gam_{0}}\Gamma
$ for $\Gamma\subset\Sigma(\mathcal{L}_{N+1})$.
Assume
$\Tht\subset C_{\gam+1}(\psi_{\mI_{N}}(\gam+1))$ and
$\hat{a}:=\gam+\omega^{a}<\gam_{0}$.

Then 
$(\mathcal{H}_{\hat{a}},\Theta;\mathtt{D}\cup\{\beta\})
\vdash^{* \beta}_{\beta}
\Gamma^{(\beta,\mI_{N})}$ holds for $\beta=\psi_{\mI_{N}}(\hat{a})$.
\elem
\bprf
By induction on $a$ as in \cite{Buchholz}.

We have
$\{\gamma, a\}\subset\mathcal{H}_{\gamma}[\Theta]$ by (\ref{eq:controlder*1}).
We obtain
$\beta\in\mathcal{H}_{\hat{a}}[\Theta]$, and
$\Tht\subset C_{\gam+1}(\psi_{\mI_{N}}(\gam+1))\cap\mI_{N}=\psi_{\mI_{N}}(\gam+1)\subset\beta$ by the assumption.

On the other hand we have 
$\sfk(\Gam)\cup\mathtt{D}\subset\mathcal{H}_{\gamma}[\Theta]$ by (\ref{eq:controlder*1}).
We obtain
\beqn\label{eq:Kcollpase.1}
\sfk(\Gam)\cup\mathtt{D}\subset\psi_{\mI_{N}}(\gam+1)
\subset\beta
\eeqn
\textbf{Case 1}. 
The last inference is an $(i {\rm -stbl}(\mS))$:
We have an $\mS\in\mathtt{D}$,
a $\bigwedge$-formula
$B(\mathsf{L}_{0})\in\Delta_{0}(\mathbb{S})$, and
a term $u\in Tm(\mathbb{I}_{N})$
such that
$\mathsf{k}(B(u))\cup\{\mS\}\subset \mathcal{H}_{\gamma}[\Theta]\cap\mathbb{I}_{N}$ by (\ref{eq:controlder*1}).
The assumption $\Theta\subset
C_{\gam+1}(\psi_{\mI_{N}}(\gam+1))$ yields
$\mathsf{k}(B(u))\cup\{\mS\}\subset
C_{\gam+1}(\psi_{\mI_{N}}(\gam+1))\cap\mathbb{I}_{N}=\psi_{\mathbb{I}_{N}}(\gamma+1)$, and
$\max\{\mathrm{rk}(B(u)),\mS\}<\psi_{\mathbb{I}_{N}}(\gamma+1)\leq\beta$ by Proposition \ref{lem:rank}.\ref{lem:rank1}.
\\
\textbf{Case 2}.
The case when the last inference is a $(\Sigma(St)\mbox{{\rm -rfl}})$:
We have ordinals
$a_{\ell}, a_{r}<a$ and a formula $C\in\Sigma(\mathcal{L}_{N+1})$ such that
$(\mathcal{H}_{\gamma},\Theta;\mathtt{D}
)\vdash^{* a_{\ell}}_{\mI_{N}+1}\Gamma,C$
and
$(\mathcal{H}_{\gamma},\Theta;\mathtt{D}
)\vdash^{* a_{r}}_{\mI_{N}+1}
\lnot \exists x\,C^{(x,\mI_{N})} , \Gamma$.

Let $\beta_{\ell}=\psi_{\mI_{N}}(\widehat{a_{\ell}})\in LSt_{N}$ with 
$\widehat{a_{\ell}}=\gamma+\omega^{a_{\ell}}$.
$\beta_{\ell}<\beta$ follows from $a_{\ell}\in\calh_{\gam}[\Tht]$.
IH with $\Sigma$-persistency \ref{lem:persistency} yields
$(\mathcal{H}_{\hat{a}},\Theta;\mathtt{D}\cup\{\beta,\beta_{\ell}\}
)\vdash^{* \beta_{\ell}}_{\beta}\Gamma^{(\beta,\mI_{N})},
C^{(\beta_{\ell},\mI_{N})}$, where
$\beta_{\ell}\in\mathcal{H}_{\widehat{a_{\ell}}}[\Theta]$ and
$\beta_{\ell}\in\mathcal{H}_{\widehat{a_{\ell}}}[\Theta]$.

On the other hand we have
$(\mathcal{H}_{\widehat{a_{\ell}}},\Theta\cup \{\beta_{\ell}\};\mathtt{D}\cup\{\beta_{\ell}\}
)\vdash^{* a_{r}}_{\mathbb{I}_{N}+1}
\lnot C^{(\beta_{\ell},\mathbb{I}_{N})},\Gamma$ by Weakening \ref{prp:inversion*_I_1}.\ref{prp:inversion*_I_1.2} and
 Inversion \ref{lem:inversionreg*}\ref{lem:inversionreg*1}, and
$(\mathcal{H}_{\widehat{a_{\ell}}},\Theta;\mathtt{D}\cup\{\beta_{\ell}\}
)\vdash^{* a_{r}}_{\mathbb{I}_{N}+1}
\lnot C^{(\beta_{\ell},\mathbb{I}_{N})},\Gamma$.

Let
$\beta_{r}=\psi_{\mathbb{I}_{N}}(\widehat{a_{r}})<\beta$ with 
$\widehat{a_{r}}=\widehat{a_{\ell}}+\omega^{a_{r}}=\gamma+\omega^{a_{\ell}}+\omega^{a_{r}}<\gamma+\omega^{a}=\hat{a}$.
IH yields
$(\mathcal{H}_{\hat{a}},\Theta;\mathtt{D}\cup\{\beta,\beta_{\ell}\})
\vdash^{* \beta_{r}}_{\beta}
\lnot C^{(\beta_{\ell},\mathbb{I}_{N})},\Gamma^{(\beta,\mathbb{I}_{N})}$.
A $(cut)$ yields
$(\mathcal{H}_{\hat{a}},\Theta;\mathtt{D}\cup\{\beta,\beta_{\ell}\})
\vdash^{* \beta_{r}+1}_{\beta}
\Gamma^{(\beta,\mathbb{I}_{N})}$ for
$\mathrm{rk}(C^{(\beta_{\ell},\mathbb{I}_{N})})<\beta$.
We obtain
$(\mathcal{H}_{\hat{a}},\Theta;\mathtt{D}\cup\{\beta\})
\vdash^{* \beta}_{\beta}
\Gamma^{(\beta,\mathbb{I}_{N})}$ by a $(\mathtt{D})$.
\\
\textbf{Case 3}.
The last inference is a 
$(cut)$:
There exist $a_{0}<a$ and a $\bigvee$-formula $C$
such that $\mathrm{rk}(C)\leq\mathbb{I}_{N}$,
$(\mathcal{H}_{\gamma},\Theta;\mathtt{D})
\vdash^{* a_{0}}_{\mathbb{I}_{N}+1}
\Gamma,\lnot C$
and
$(\mathcal{H}_{\gamma},\Theta;\mathtt{D})
\vdash^{* a_{0}}_{\mathbb{I}_{N}+1}C,\Gamma$.
We obtain $C\in\Sigma_{1}(\mathbb{I}_{N})\cup \Delta_{0}(\mathbb{I}_{N})$ by Proposition \ref{lem:rank}.\ref{lem:rank3}.
IH yields
$(\mathcal{H}_{\hat{a}},\Theta;\mathtt{D}\cup\{\beta,\beta_{0}\})
\vdash^{* \beta_{0}}_{\beta}C^{(\beta_{0},\mathbb{I}_{N})},\Gamma_{1}^{(\beta,\mathbb{I}_{N})}$
for $\beta_{0}=\psi_{\mathbb{I}_{N}}(\widehat{a_{0}})\in\mathcal{H}_{\widehat{a_{0}}}[\Theta]$ with $\widehat{a_{0}}=\gamma+\omega^{a_{0}}$.
On the other, we obtain
$(\mathcal{H}_{\widehat{a_{0}}},\Theta;\mathtt{D}\cup\{\beta_{0}\})
)\vdash^{* a_{0}}_{\mathbb{I}_{N}+1}
\lnot C^{(\beta_{0},\mathbb{I}_{N})},\Gamma_{0}$
by Weakening \ref{prp:inversion*_I_1}.\ref{prp:inversion*_I_1.2} and
Inversion \ref{lem:inversionreg*}.\ref{lem:inversionreg*2} when $C\in\Sigma_{1}(\mathbb{I}_{N})$.
Let $\beta_{1}=\psi_{\mathbb{I}_{N}}(\widehat{a_{1}})<\beta$ with $\widehat{a_{1}}=\widehat{a_{0}}+\omega^{a_{0}}=\gamma+\omega^{a_{0}}\cdot 2<\hat{a}$.
IH yields
$(\mathcal{H}_{\hat{a}},\Theta;\mathtt{D}\cup\{\beta,\beta_{0}\})
)\vdash^{* \beta_{1}}_{\beta}
\lnot C^{(\beta_{0},\mathbb{I}_{N})},\Gamma_{0}^{(\beta,\mathbb{I}_{N})}$.
We obtain
$(\mathcal{H}_{\hat{a}},\Theta;\mathtt{D}\cup\{\beta,\beta_{0}\})
)\vdash^{* \beta_{1}+1}_{\beta}
\Gamma^{(\beta,\mathbb{I}_{N})}$
by a $(cut)$ with $\mathrm{rk}(C^{(\beta_{0},\mathbb{I}_{N})})<\beta$.
A $(\mathtt{D})$ yields
$(\mathcal{H}_{\hat{a}},\Theta;\mathtt{D}\cup\{\beta\})
)\vdash^{* \beta}_{\beta}
\Gamma^{(\beta,\mathbb{I}_{N})}$.
\\
\textbf{Case 4}.
The last inference is a $(\bigwedge)$:
A formula $A\in\Sigma(\mathbb{I}_{N})$ with $A\simeq\bigwedge\left(A_{\iota}\right)_{\iota\in J}$ 
is introduced in $\Gamma$.
For every
$\iota\in J$ there exists an $a(\iota)<a$ 
such that
$(\mathcal{H}_{\gamma},\Theta\cup\mathsf{k}(\iota);\mathtt{D}\cup\{\mathrm{srk}(A_{\iota})\}
)
\vdash^{* a(\iota)}_{\mathbb{I}_{N}+1}\Gamma, A_{\iota}$.
We see $\mathsf{k}(\iota)\subset\psi_{\mathbb{I}_{N}}(\gamma+1)$ as follows.
For example let $A\equiv(\forall x\in u\, B(x))$. Then 
$u\in Tm(\mathbb{I}_{N})$ by 
$A\in\Sigma(\mathbb{I}_{N})$, and $J=Tm(|u|)$. We obtain 
$\mathsf{k}(u)\subset\mathcal{H}_{\gamma}[\Theta]\cap\mathbb{I}_{N}\subset\psi_{\mathbb{I}_{N}}(\gamma+1)$.
Hence  for $\iota\in J$, $|\iota|<|u|<\psi_{\mathbb{I}_{N}}(\gamma+1)$.

IH yields 
$(\mathcal{H}_{\hat{a}},\Theta\cup \mathsf{k}(\iota);\mathtt{D}\cup\{\beta\}\cup\mathrm{srk}(A_{\iota})
)
\vdash^{* \beta(\iota)}_{\beta}
\Gamma^{(\beta,\mI_{N})},
A_{\iota}^{(\beta,\mI_{N})}$
for each $\iota\in J$, where
$\beta(\iota)=\psi_{\mathbb{I}_{N}}(\widehat{a(\iota)})$ with
$\widehat{a(\iota)}=\gamma+\omega^{a(\iota)}$.
We obtain
$(\mathcal{H}_{\hat{a}},\Theta;\mathtt{D}\cup\{\beta\})\vdash^{* \beta}_{\beta}
\Gamma^{(\beta,\mI_{N})}$ by a $(\bigwedge)$ and $\beta(\iota)<\beta$.
\\
\textbf{Case 5}.
The last inference is a $(\bigvee)$:
A $\bigvee$-formula with
$A\in\Gamma$ is introduced.
Let $A\simeq\bigvee\left(A_{\iota}\right)_{\iota\in J}$.
There are an $\iota\in J$, an ordinal
 $a(\iota)<a$
such that
$(\mathcal{H}_{\gamma},\Theta;\mathtt{D})\vdash^{* a(\iota)}_{\mathbb{I}_{N}+1}
\Gamma,A_{\iota}$.
We may assume $\mathsf{k}(\iota)\subset\mathsf{k}(A_{\iota})$.
We obtain 
by (\ref{eq:controlder*1}),
$\mathsf{k}(\iota)\subset\mathcal{H}_{\gamma}[\Theta]\cap\mathbb{I}_{N}
\subset C_{\gamma}(\psi_{\mathbb{I}_{N}}(\gamma+1))\cap\mathbb{I}_{N}\subset
\psi_{\mathbb{I}_{N}}(\gamma+1)\subset\beta$.
IH yields $(\mathcal{H}_{\hat{a}},\Theta;\mathtt{D}\cup\{\beta\})
\vdash^{*\beta(\iota)}_{\beta}
\Gamma^{(\beta,\mathbb{I}_{N})},
A_{\iota}^{(\beta,\mathbb{I}_{N})}$
for $\beta(\iota)=\psi_{\mathbb{I}_{N}}(\widehat{a(\iota)})$ with
$\widehat{a(\iota)}=\gamma+\omega^{a(\iota)}$.
$(\mathcal{H}_{\hat{a}+},\Theta;\mathtt{D}\cup\{\beta\})
\vdash^{*\beta}_{\beta}
\Gamma^{(\beta,\mathbb{I}_{N})}$ 
follows from a $(\bigvee)$.
\\
\textbf{Case 6}.
The last inference is a $(\mathtt{D})$:
We have $(\mathcal{H}_{\gamma},\Theta;\mathtt{D}\cup\{\mS\})\vdash^{* a_{0}}_{\mathbb{I}_{N}+1}
\Gamma$ for $a_{0}<a$ and $\mS\in\mathcal{H}_{\gamma}[\Theta]$.
By the assumption we obtain $\mS<\psi_{\mI_{N}}(\gamma+1)\leq\beta$.
IH yields $(\mathcal{H}_{\hat{a}},\Theta;\mathtt{D}\cup\{\beta,\mS\})\vdash^{* \beta_{0}}_{\beta}
\Gamma^{(\beta,\mI_{N})}$, where 
$\beta_{0}=\psi_{\mathbb{I}_{N}}(\widehat{a_{0}})\in\mathcal{H}_{\hat{a}}[\Theta]$ with $\widehat{a_{0}}=\gamma+\omega^{a_{0}}$.
We obtain $(\mathcal{H}_{\hat{a}},\Theta;\mathtt{D}\cup\{\beta\})\vdash^{* \beta}_{\beta}
\Gamma^{(\beta,\mI_{N})}$
by a $(\mathtt{D})$.

Other cases are seen from IH.
\eprf

\subsection{Operator controlled derivations with caps}\label{subsec:capderivation}

Let $\Lambda_{0}=\beta=\psi_{\mI_{N}}(b_{0})$ be the ordinal in Collapsing \ref{lem:Kcollpase.1}
with $b_{0}=\gamma+\omega^{a}$ and $\{\gamma,a\}\subset C_{0}(0)$.
$\Lambda_{0}$ is a bound on the ordinals $\mS$ and $s(g)$ in inferences $(i {\rm -rfl}_{\mS}(g,x,f,\Delta))$ of Definition \ref{df:controldercap}.
The base of the $\tilde{\theta}$-function
$\tilde{\theta}_{b}(\xi)=\tilde{\theta}_{b}(\xi;\Lambda)$ in Definition \ref{df:Lam} and of finite functions
is the ordinal $\mI_{N}$.

\bdf\label{df:guarded_ord_I_1}
{\rm
By a \textit{guarded ordinal} we mean a pair $(\rho,g)$ of an ordinal $\rho\in\Psi_{\mS}$ 
with a successor stable ordinal $\mS\in SSt$ and a special finite function $g$ with base $\mI_{1}$
such that $g\leq m(\rho):\Lrarw \forall i((m(\rho))(i)\leq g(i))$.
$\Psi_{\mS}^{\mathtt{g}}$ denotes the set of such pairs $(\rho,g)$ with $\rho\in\Psi_{\mS}$.

Let $(\rho,g)_{0}:=\rho$ and $(\rho,g)_{1}:=g$. 
For sets $X$ of guarded ordinals let $(X)_{0}=\{\kappa: \exists g((\kappa,g)\in X)\}$.

For (unguarded) ordinals $\alpha$, let $(\alpha)_{0}:=\alpha$ and $(\alpha)_{1}=\emptyset$.
For guarded and/or unguarded ordinals $\alpha,\beta$, let
$\alpha<\beta:\Lrarw (\alpha)_{0}<(\beta)_{0}$.
}
\edf

A guarded ordinal $(\rho,g)$ is denoted by $\rho$ when no confusion likely occurs.

\bdf\label{df:ffthreshold}
{\rm
$\Lambda_{0}$ is the ordinal $\beta$ in Collapsing \ref{lem:Kcollpase.1} such that
$\Lambda_{0}=\psi_{\mI_{N}}(b_{0})$ with $b_{0}=\gamma+\omega^{a}$,
and $\{\gamma,a\}\subset C_{0}(0)$.

A triple $(\mathtt{Q},\gamma_{\cdot}^{\mathtt{Q}}, \delta^{\mathtt{Q}}_{\cdot})$ 
is said to be a \textit{finite family with thresholds}
if the following conditions are met.

\benu
 
\item\label{df:QJ}
$\mathtt{Q}$ 
is a finite function $\mathtt{Q}\subset\coprod_{\mS}\Psi_{\mS}$ for which the following hold:
 \benu
 \item\label{df:QJ.1}
 Its domain $\mathrm{dom}(\mathtt{Q})$ is a predecessor closed set of 
 ordinals $\mS\in (St\cup\{0\})\cap\Lambda_{0}$.

 \item\label{df:QJ.2}
For each $\mS\in\mathrm{dom}(\mathtt{Q})$,
$\mathtt{Q}(\mS)$ is a finite set of guarded ordinals in $\Psi^{\mathtt{g}}_{\mS}$
 such that
$f=g$ if $(\kappa,f),(\kappa,g)\in\mathtt{Q}(\mS)$,
$\forall\{\tau<\kappa\}\subset(\mathtt{Q}(\mS))_{0}(M_{\tau}\subset M_{\kappa})$,
$\fal\mU<\mS(\kappa\in M_{\mathtt{Q}(\mU)})$ with 
$M_{\mathtt{Q}(\mU)}=\bigcap\{M_{\rho}:\rho\in (\mathtt{Q}(\mU))_{0}\}$,
%cf. Tautology \ref{lem:tautology.cap}.\ref{lem:tautology.cap3} 
and $m(\kappa)$ is a special finite function
with base $\mI_{N}$
for each $(\kappa,g)\in\mathtt{Q}(\mS)$.

For $\mS\not\in SSt$, $\mathtt{Q}(\mS)=\emptyset$.

\item\label{df:QJ.3}
$(\mathtt{Q}(\mS))_{0}\cap St\subset\mathrm{dom}(\mathtt{Q})$ for every 
$\mS\in\mathrm{dom}(\mathtt{Q})$.

 \eenu

Let
$\rho_{\mathtt{Q}(\mS)}:=\min(\mathtt{Q}(\mS))_{0}$ for $\mS\in\mathrm{dom}(\mathtt{Q})\cap SSt$,
$\bigcup\mathtt{Q}=\bigcup\{(\mathtt{Q}(\mS))_{0}:\mS\in \mathrm{dom}(\mathtt{Q})\}$,
$M_{\mathtt{Q}}=\bigcap\{M_{\rho}:\rho\in\bigcup\mathtt{Q}\}$
and
$(\mathtt{Q})_{1}=\bigcup\{(\mathtt{Q}(\mS))_{1}:\mS\in\mathrm{dom}(\mathtt{Q})\}$.

\item\label{df:caphat.42}
$\gamma_{\cdot}^{\mathtt{ Q}}$ is a \textit{threshold function} $\mathrm{dom}(\mathtt{Q})\ni\mS\mapsto\gam_{\mS}^{\mathtt{Q}}$ such that
$\gamma^{\mathtt{Q}}_{\mS}+\mS<\gamma^{\mathtt{Q}}_{\mT}$ and
$\gamma^{\mathtt{Q}}_{\mS}\in M_{\mathtt{Q}(\mT)}$ for every $\{\mT<\mS\}\subset \mathrm{dom}(\mathtt{Q})$,
$\gamma^{\mathtt{Q}}_{\mS}=\alpha_{\mS}+\beta_{\mS}$ for an $\alpha_{\mS}$ and 
an additive principal number $\mI_{N}>\beta_{\mS}$,
and the following 
(\ref{eq:controlder_cap4}) is met.

\beqn\label{eq:controlder_cap4}
\forall \mS\in \mathrm{dom}(\mathtt{Q})\forall \rho\in (\mathtt{Q}(\mS))_{0}
\left[
\gamma^{\mathtt{Q}}_{\mS}<\mathtt{p}_{0}(\rho)\leq b(\rho)<\gamma^{\mathtt{Q}}_{\mS}+\mS
\right]
\eeqn
$b(\rho)=a$ for $\rho=\psi_{\sig}^{f}(a)$.

 \item\label{df:caphat.52_I_1}
$\delta_{\cdot}^{\mathtt{ Q}}$ is a function 
$\mathrm{dom}(\mathtt{Q})\ni\mS\mapsto\delta_{\mS}^{\mathtt{Q}}$ such that
$\delta^{\mathtt{Q}}_{\mS}+\omega\leq \delta^{\mathtt{Q}}_{\mT}$ and
$1<\delta_{\mS}^{\mathtt{Q}}<\omega^{\Lambda_{0}+1}$ for every 
$\{\mT<\mS\}\subset \mathrm{dom}(\mathtt{Q})$, and
\begin{equation}
\label{eq:controlder_cap00_I_1}
\forall\mS\in \mathrm{dom}(\mathtt{Q})
\left(
\delta^{\mathtt{Q}}_{\mS}\in C_{b_{0}+1}(\psi_{\mI_{N}}(b_{0}+1))
\right)
\end{equation}
where $\Lambda_{0}=\psi_{\mI_{N}}(b_{0})$.

\eenu

For a set $\Theta$ of ordinals, let
$\Theta_{\mathtt{Q}}=\Theta\cup
\{\delta^{\mathtt{Q}}_{\mS}: \mS\in\mathrm{dom}(\mathtt{Q})\}\cup
\bigcup\{SC_{\mI_{N}}(\mathrm{fld}(g)) : g\in (\mathtt{Q})_{1}\}$.

A triple $(\mathtt{Q}, \gamma^{\mathtt{Q}}_{\cdot}, \delta^{\mathtt{Q}}_{\cdot})$ is simply denoted by $\mathtt{ Q}$ 
when the threshold functions $\gamma^{\mathtt{Q}}_{\cdot}, \delta^{\mathtt{Q}}_{\cdot}$
are irrelevant.

Let
$\Lambda^{\mathtt{Q}}_{\mS}:=\psi_{\mI_{N}}(b_{0}+\delta^{\mathtt{Q}}_{\mS})$,
$\Lambda^{(n)}_{\mS}:=\psi_{\mI_{N}}(b_{0}+\delta^{\mathtt{Q}}_{\mS}+n)$ for $0\leq n<\omega$
and $0\neq\mS\in\mathrm{dom}(\mathtt{Q})$.
Also let $\Lambda^{\mathtt{Q}}_{\Lambda_{0}}:=\Lambda^{(n)}_{\Lambda_{0}}:=\Lambda_{1}=\psi_{\mI_{N}}(b_{0}+1)$
and
$\Lambda^{\mathtt{Q}}_{0}:=\Lambda^{(n)}_{0}:=\psi_{\mI_{N}}(b_{0}+\omega^{\Lambda_{0}+1})$
 for $0\leq n<\omega$.
}
\edf
Let $\{\mT<\mS\}\subset \mathrm{dom}(\mathtt{Q})$. Then
we have $1<\delta^{\mathtt{Q}}_{\mS}$ and 
$\delta^{\mathtt{Q}}_{\mS}+n<\delta^{\mathtt{Q}}_{\mT}$.
By (\ref{eq:controlder_cap00_I_1}) we obtain 
\beqn\label{eq:delta_Lambda}
\Gamma(\Lambda_{\Lambda_{0}}^{(n)})<\Lambda^{\mathtt{Q}}_{\mS}
<
\psi_{\mI_{N}}(b_{0}+\delta^{\mathtt{Q}}_{\mS}+1+n)=\Lambda^{(1+n)}_{\mS}
<\psi_{\mI_{N}}(b_{0}+\delta^{\mathtt{Q}}_{\mT})=\Lambda^{\mathtt{Q}}_{\mT}
\eeqn

\bdf\label{df:extend_family}
{\rm

Let $\mathtt{Q}$ be a finite family with thresholds $\gamma_{\cdot}^{\mathtt{Q}}, \delta_{\cdot}^{\mathtt{Q}}$.
\benu

\item\label{df:extend_family2}
The threshold function $\gam^{\mathtt{Q}}_{\cdot}$ with an ordinal $\gamma$ is said to have \textit{gaps} $\eta$ if
for every $\{\alpha<\beta\}\subset\{0,\gamma\}\cup\{\gam^{\mathtt{Q}}_{\mS}:\mS\in \mathrm{dom}(\mathtt{Q})\}$,
$\alpha+ \eta\leq \beta$ holds.

\item\label{df:extend_family3}
Let $\mS\in St\cup\{0\}$ be an ordinal such that $\mathrm{dom}(\mathtt{Q})\cup\{\mS\}$ is
predecessor closed.
$\mathtt{P}=\mathtt{Q}[\mS]$
denotes a triple $(\mathtt{P},\gamma^{\mathtt{P}}_{\cdot},\delta^{\mathtt{P}}_{\cdot})$ enjoying the following:

\benu
\item
$\mathrm{dom}(\mathtt{P})=\mathrm{dom}(\mathtt{ Q})\cup\{\mS\}$, 
$\mathtt{P}(\mT)=\mathtt{ Q}(\mT)$ for $\mT\neq\mS$. 
$\mathtt{P}(\mS)=\mathtt{Q}(\mS)\cup X$
for a finite set $X\subset\Psi^{\mathtt{g}}_{\mS}$ when $\mS$ is a successor stable ordinal.
$\mathtt{P}(\mS)=\mathtt{Q}(\mS)$ if $\mS\in LSt\cup\{0\}$.

\item
$\gam^{\mathtt{P}}_{\cdot}$ extends $\gam^{\mathtt{ Q}}_{\cdot}$ in such a way that
$\gam^{\mathtt{P}}_{\mT}=\gam^{\mathtt{ Q}}_{\mT}$ for $\mT\in \mathrm{dom}(\mathtt{ Q})$,
$\gam^{\mathtt{P}}_{\mS}>\gam^{\mathtt{Q}}_{\mT}+\mT$ and
$\gam^{\mathtt{Q}}_{\mT}\in M_{\mathtt{P}(\mS)}$ for every $\mS<\mT\in \mathrm{dom}(\mathtt{ Q})$,
$\gam^{\mathtt{ Q}}_{\mU}>\gam^{\mathtt{P}}_{\mS}+\mS$ and 
$\gam^{\mathtt{P}}_{\mS}\in M_{\mathtt{Q}(\mU)}$ for every $\mS>\mU\in\mathrm{dom}(\mathtt{ Q})$, and
$\forall \rho\in (\mathtt{P}(\mS))_{0}
\left[
\gamma^{\mathtt{P}}_{\mS}<\mathtt{p}_{0}(\rho)\leq b(\rho)<\gamma^{\mathtt{P}}_{\mS}+\mS
\right]
$, cf.\,(\ref{eq:controlder_cap4}).

\item
$\delta^{\mathtt{P}}_{\cdot}$ extends $\delta^{\mathtt{ Q}}_{\cdot}$ in such a way that
$\delta^{\mathtt{P}}_{\mT}=\delta^{\mathtt{ Q}}_{\mT}$ for $\mT\in \mathrm{dom}(\mathtt{ Q})$,
$\delta^{\mathtt{P}}_{\mS}>1$,
$\delta^{\mathtt{P}}_{\mS}\geq\delta^{\mathtt{ Q}}_{\mT}+\omega$ for every $\mS<\mT\in \mathrm{dom}(\mathtt{Q})$,
$\delta^{\mathtt{ Q}}_{\mU}\geq\delta^{\mathtt{P}}_{\mS}+\omega$ for every $\mS>\mU\in\mathrm{dom}(\mathtt{Q})$, and
$\delta^{\mathtt{P}}_{\mS}\in C_{b_{0}+1}(\psi_{\mI_{N}}(b_{0}+1))$, cf.\,(\ref{eq:controlder_cap00_I_1}).

\eenu

 By convention, let $\mathtt{Q}[\mS]:=\mathtt{Q}$ when $\mS\in\mathrm{dom}(\mathtt{Q})$.
 An extension $\mathtt{Q}[\{\mS_{1},\ldots,\mS_{n}\}]$ of $\mathtt{Q}$ by a finite set $\{\mS_{1}>\cdots>\mS_{n}\}$ of successor stable ordinals
 is defined recursively by $\mathtt{Q}[\{\mS_{1},\ldots,\mS_{n}\}]=(\mathtt{Q}[\{\mS_{1},\ldots,\mS_{n-1}\}])[\mS_{n}]$.
 
\item\label{df:extend_family4}
For an ordinal $\xi$, $\mathtt{R}=\mathtt{Q}\restrict\xi$ denotes the restriction of $\mathtt{Q}$ to $\xi$ defined as follows.

\benu
\item
$\mathrm{dom}(\mathtt{R})=\mathrm{dom}(\mathtt{ Q})\cap\xi$, 
$\mathtt{R}(\mT)=\mathtt{ Q}(\mT)$ for $\mT<\xi$.

\item
$\gam^{\mathtt{R}}_{\cdot}$ [$\delta^{\mathtt{R}}_{\cdot}$] is the restriction of $\gam^{\mathtt{ Q}}_{\cdot}$ [of $\delta^{\mathtt{Q}}_{\cdot}$] to $\xi$ in such a way that
$\gam^{\mathtt{R}}_{\mT}=\gam^{\mathtt{ Q}}_{\mT}$ 
[$\delta^{\mathtt{R}}_{\mT}=\delta^{\mathtt{Q}}_{\mT}$] for $\mT\in \mathrm{dom}(\mathtt{Q})\cap\xi$,
resp.

\eenu

\end{enumerate}
}
\edf

\bdf\label{df:good_extension}
{\rm
Let $\mathtt{P}=\mathtt{Q}[\mS]$ be an extension of $\mathtt{Q}$ by a stable ordinal $\mS$ with thresholds
$\gamma^{\mathtt{P}}_{\cdot},\delta^{\mathtt{P}}_{\cdot}$.

\benu
\item\label{df:good_extension.1}
The extension $(\mathtt{P},\gam^{\mathtt{P}}_{\cdot},\delta^{\mathtt{P}}_{\cdot})$ is said to be a \textit{good} extension of $\mathtt{Q}$ for an ordinal $\gamma$ and a set $\Theta$ of ordinals
 if the following conditions are met when $\mS\not\in\mathrm{dom}(\mathtt{Q})$:

\benu

\item
$\mS\in\mathcal{H}_{\gamma}[\Theta(\mathtt{Q})]$ and
$\gamma\leq \gamma^{\mathtt{P}}_{\mS}$.

\item
$\mathtt{P}(\mS)=\{(\rho_{\mathtt{P}(\mS)},g_{\mS})\}\subset\Psi^{\mathtt{g}}_{\mS}$ such that 
$\Theta(\mathtt{Q})\cup\{\delta^{\mathtt{P}}_{\mS}\}\cup\mathrm{fld}(g_{\mS})
\subset M_{\mathtt{P}}=M_{\mathtt{Q}}\cap M_{\rho_{\mathtt{P}(\mS)}}$,
where $\Theta(\mathtt{Q})=\Theta\cap M_{\mathtt{Q}}$, cf.\,Definition \ref{df:Qmin}.\ref{df:Qmin.4}.

\eenu

We say that $\mathtt{P}$ is a good extension by the ordinals $\gamma^{\mathtt{P}}_{\mS}, \delta^{\mathtt{P}}_{\mS}$
 and 
 $(\rho_{\mathtt{P}(\mS)},g_{\mS})$.

When $\mS\in\mathrm{dom}(\mathtt{Q})$, $\mathtt{Q}[\mS]=\mathtt{Q}$ is defined to be a good extension of $\mathtt{Q}$
(for any $\gamma,\Theta$).

\item\label{df:good_extension.11}
Let $\mathtt{P}=\mathtt{Q}[\mS]$ be a good extension
of $\mathtt{Q}$ for $\gamma$ and $\Theta$ 
 by ordinals $\gamma^{\mathtt{P}}_{\mS}, \delta^{\mathtt{P}}_{\mS}$
and a guarded ordinal  $(\rho_{\mathtt{P}(\mS)},g_{\mS})$.
Let $\Theta_{\gamma,\mS}=\Theta\cup\{\gamma^{\mathtt{Q}}_{\mU} :\mS<\mU\in\mathrm{dom}(\mathtt{Q})\}$ and
$\Theta_{\delta,\mS}=\Theta\cup\{\delta^{\mathtt{Q}}_{\mU}: \mS<\mU\in\mathrm{dom}(\mathtt{Q})\}$.

Then $\mathtt{P}$ is \textit{very good} if
$\gamma^{\mathtt{P}}_{\mS}\in\mathcal{H}_{\gamma^{\mathtt{P}}_{\mS}}[\Theta_{\gamma,\mS}(\mathtt{Q}\restrict\mS)]$, 
$\{\delta^{\mathtt{P}}_{\mS}\}\cup SC_{\mI_{N}}(\mathrm{fld}(g_{\mS}))\subset
\mathcal{H}_{\gamma}[\Theta_{\delta,\mS}(\mathtt{Q})]$ and
$\rho_{\mathtt{P}(\mS)}\in\mathcal{H}_{\gamma^{\mathtt{P}}_{\mS}+\mS}[\Theta_{\gamma,\mS}(\mathtt{Q}\restrict\mS)]$, cf.\,\textbf{Case 1} in the proof of Capping \ref{lem:capping} and
\textbf{Case 1.2} in Lemma \ref{lem:main.1_I_1}.

\eenu

}
\edf

\bdf\label{df:caphat}
{\rm
Let $\mathtt{Q}$ be a finite family.

\begin{enumerate}

\item\label{df:caphat.0}
A function $\bar{\rho}\in\prod_{\mS\in\mathrm{dom}(\mathtt{Q})\cap SSt}\mathtt{Q}(\mS)$ is said to be a \textit{cap} over $\mathtt{Q}$.
For each $\mS\in\mathrm{dom}(\mathtt{Q})\cap SSt$, $\bar{\rho}(\mS)\in\mathtt{Q}(\mS)\subset\Psi_{\mS}^{\mathtt{g}}$ holds.
$\rho\in\bar{\rho}$ designates that $(\bar{\rho}(\mS))_{0}=\rho$ for an $\mS\in\mathrm{dom}(\mathtt{Q})$.

Let $M_{\bar{\rho}}:=\bigcap\{M_{\bar{\rho}(\mS)}: \mS\in\mathrm{dom}(\mathtt{Q})\cap SSt\}$, and
$\bar{\rho}_{\mathtt{Q}}$ denotes the least cap such that 
$(\bar{\rho}_{\mathtt{Q}}(\mS))_{0}=\rho_{\mathtt{Q}(\mS)}$ for each $\mS\in\mathrm{dom}(\mathtt{Q})$.
The cap $\bar{\rho}_{\mathtt{Q}}$ over the family $\mathtt{Q}$ such that
$\mathrm{dom}(\mathtt{Q})\cap SSt=\emptyset$, is
the empty function $\bar{\rho}_{\mathtt{Q}}=\emptyset$. Then let
$M_{\emptyset}=OT(\mI_{N})$.

\item\label{df:caphat.11}
Let $\mS\in\mathrm{dom}(\mathtt{Q})$ and $(\sigma,f)\in\Psi_{\mS}^{\mathtt{g}}$ be a guarded ordinal.
Then $\mathtt{Q}^{\sigma}$ denotes an extension of $\mathtt{Q}$ such that
$\mathtt{Q}^{\sigma}(\mS)=\mathtt{Q}(\mS)\cup\{(\sigma,f)\}$.
$\mathrm{dom}(\mathtt{Q}^{\sigma})=\mathrm{dom}(\mathtt{Q})\cup\{\sigma\}$ if $\sigma\in LSt$.
In this case $\mathtt{Q}^{\sigma}(\sigma)=\emptyset$ and
threshold functions $\gamma^{\mathtt{Q}^{\sigma}}_{\cdot}$ and $\delta^{\mathtt{Q}^{\sigma}}_{\cdot}$
extend $\gamma^{\mathtt{Q}}_{\cdot}$ and $\delta^{\mathtt{Q}}_{\cdot}$, resp.,
and
ordinals $\gamma^{\mathtt{Q}^{\sigma}}_{\sigma}$ and $\delta^{\mathtt{Q}^{\sigma}}_{\sigma}$
are attached as in Definition \ref{df:extend_family}.\ref{df:extend_family3}.

Otherwise $\mathrm{dom}(\mathtt{Q}^{\sigma})=\mathrm{dom}(\mathtt{Q})$.

For a cap $\bar{\rho}$ over $\mathtt{Q}$, $\bar{\rho}[\sigma/\rho]$ denotes a cap $\bar{\tau}$ over $\mathtt{Q}^{\sigma}$
such that $\bar{\tau}(\mT)=\bar{\rho}(\mT)$ for $\mT\neq\mS$.
$\bar{\tau}(\mS)$ is defined as follows.
$\bar{\tau}(\mS)=(\sigma,f)$ if $(\bar{\rho}(\mS))_{0}=\rho$, and $\bar{\tau}(\mS)=\bar{\rho}(\mS)$ else.

\item\label{df:caphat.12}
Let $\bar{\rho}$  be  a cap over $\mathtt{Q}$, and $\mathtt{P}=\mathtt{Q}[\mS]$ be 
an extension with 
$\mS\not\in\mathrm{dom}(\mathtt{Q})$.
For $(\rho,g)\in\mathtt{P}(\mS)$, $\bar{\tau}=\bar{\rho}*(\rho,g)$ denotes a cap over $\mathtt{P}$ such that
$\bar{\tau}(\mS)=(\rho,g)$, and $\bar{\tau}(\mT)=\bar{\rho}(\mT)$ for 
$\mT\in \mathrm{dom}(\mathtt{Q})$.

\item\label{df:caphat.1}
By a \textit{capped formula} over $\mathtt{Q}$ we mean a pair $(A,\bar{\rho})$ of $RS$-sentence $A$
and a cap $\bar{\rho}$ over $\mathtt{Q}$
 such that
$\mathsf{k}(A)\subset M_{\bar{\rho}}$.
Such a pair is denoted by $A^{(\bar{\rho})}$.
A capped formula $A^{(\bar{\rho})}$ is said to be a $\Sigma(\pi)$-formula if
$A\in\Sigma(\pi)$.

Let $\sfk(A^{(\bar{\rho})}):=\sfk(A)$, $\mathrm{rk}(A^{(\bar{\rho})}):=\mathrm{rk}(A)$
and $\mathrm{srk}(A^{(\bar{\rho})}):=\mathrm{srk}(A)$.

\item\label{df:caphat.2}
A \textit{sequent} is a finite set of capped formulas, denoted by
$\Gamma_{0}^{(\bar{\rho}_{0})},\ldots,\Gamma_{m}^{(\bar{\rho}_{m})}$,
where each formula in the set $\Gamma_{i}^{(\bar{\rho}_{i})}$ puts on the cap 
$\bar{\rho}_{i}$.
When we write $\Gamma^{(\bar{\rho})}$, we tacitly assume that
$\mathsf{k}(\Gamma)\subset M_{\bar{\rho}}$.

\item\label{df:caphat.13}
Let $\Gamma$ be a sequent.
Let $\mathrm{rk}(\Gamma)=\max(\{0\}\cup\{\mathrm{rk}(A^{(\bar{\rho})}): A^{(\bar{\rho})}\in\Gamma\})$.
$\Gamma(\rho)=\{A^{(\bar{\rho})}: A^{(\bar{\rho})}\in\Gamma,\, \rho\in\bar{\rho}\}$ denotes the set of capped formulas with a cap $\rho\in\bar{\rho}$ in $\Gamma$.

\item\label{df:good_extension.3}
Let $\Gamma$ be a set of capped formulas over $\mathtt{Q}$ such that 
$\mathsf{k}(\Gamma)\subset M_{\mathtt{Q}}\cap M_{\rho}$,
where
$\{(\rho,g)\}=\mathtt{P}(\mS)$.
Then let
$\Gamma(\mathtt{P})=\{A^{(\bar{\rho}*(\rho,g))}: A^{(\bar{\rho})}\in\Gamma\}$.
\eenu
}
\edf

\bdf\label{df:inverse}
{\rm
\benu
\item\label{df:inverse1}
Let $\sigma\in\Psi_{\mS}$, and $A$ an uncapped formula.
Then a formula $A[\sigma/\mS]^{-1}$  is obtained from $A$ by replacing each $\alpha\in\mathsf{k}(A)$ by 
the uncollapsing $\alpha[\sigma/\mS]^{-1}$ of $\alpha$, cf.\,Definition \ref{df:divide},
and each $st_{i}^{[\bar{\rho}]*[\sigma]}(u)$ is replaced by  $st_{i}^{[\bar{\rho}]}(u[\sigma/\mS]^{-1})$ if
$|u[\sigma/\mS]^{-1}|\geq\mS$.

\item\label{df:inverse2}
Let $\mathtt{Q}$ be a finite family, and $A^{(\bar{\rho})}$ a capped formula for $\mathtt{Q}$.

If there are an $\mS\in\mathrm{dom}(\mathtt{Q})$ and a 
$\sigma\in(\mathtt{Q}(\mS))_{0}\cap(\bar{\rho}(\mS))_{0}$
such that
$A[\sigma/\mS]^{-1}\not\equiv A$
with $\mathsf{k}(A[\sigma/\mS]^{-1})\subset  M_{\rho_{\mathtt{Q}(\mS)}}$,
then
the \textit{inverse} $\mathrm{inv}(A^{(\bar{\rho})};\mathtt{Q})$ of $A^{(\bar{\rho})}$ for $\mathtt{Q}$ is defined to be the formula
$(A[\sigma/\mS]^{-1})^{(\bar{\rho})}$.

Otherwise let 
$\mathrm{inv}(A^{(\bar{\rho})};\mathtt{Q}):\equiv A^{(\bar{\rho})}$.

\eenu
By definition, the inverse $\mathrm{inv}(A^{(\bar{\rho})};\mathtt{Q})$ of $A^{(\bar{\rho})}$ puts on the same cap $(\bar{\rho})$ as $A^{(\bar{\rho})}$.
Let $\mathrm{inv}^{(\mathtt{u})}(A^{(\bar{\rho})};\mathtt{Q}):\equiv B$ when $\mathrm{inv}(A^{(\bar{\rho})};\mathtt{Q})\equiv B^{(\bar{\rho})}$.

}
\edf

Let $A\equiv B^{[\sigma/\mS]}\not\equiv B$.
Then $\mathrm{inv}(A^{(\bar{\rho})};\mathtt{Q})\equiv B^{(\bar{\rho})}$ holds only when 
$\sigma\in(\mathtt{Q}(\mS))_{0}\cap(\bar{\rho}(\mS))_{0}$.

\bdf\label{df:Qmin}
{\rm
Let $\mathtt{Q}$ be a finite family, $A^{(\bar{\rho})}$ a capped formula for $\mathtt{Q}$.

\benu

\item\label{df:Qmin.1}
$\mathsf{k}(A^{(\bar{\rho})};\mathtt{Q})  := 
\mathsf{k}(\mathrm{inv}(A^{(\bar{\rho})};\mathtt{Q}) )$.

\item\label{df:Qmin.3}
For $\iota\in J$ with either $A\simeq\bigvee(A_{\iota})_{\iota\in J}$ or $A\simeq\bigwedge(A_{\iota})_{\iota\in J}$,
let
\[
\iota\in[\bar{\rho}\mathtt{Q}]J
 :\Lrarw 
 \mathsf{k}(A_{\iota})\subset M_{\bar{\rho}}
\spand
\mathsf{k}(A_{\iota}^{(\bar{\rho})};\mathtt{Q})\subset 
M_{\mathtt{Q}}
.\]

\item\label{df:Qmin.4}
For finite sets $\Theta$ of ordinals, let
$\Theta(\mathtt{Q}):=
\Theta\cap M_{\mathtt{Q}}
$.

\eenu
}
\edf

\bdf\label{df:resolvent}
{\rm 
Let $\mathtt{Q}$ be a finite family with thresholds $\gam^{\mathtt{Q}}_{\cdot}, \delta^{\mathtt{Q}}_{\cdot}$, $\mS\in\mathrm{dom}(\mathtt{Q})$ and $\rho\in(\mathtt{Q}(\mS))_{0}$.
$H_{\rho}(f,\mathtt{Q},\Theta)$
 denotes the \textit{resolvent class} defined by
$(\kappa,k)\in H_{\rho}(f,\mathtt{Q},\Theta)$ iff
$\kappa\in \Psi_{\mS}\cap\rho$,
$m(\kap)$ is special such that
$f=k\leq m(\kappa)$,
$\mathtt{Q}(\mS)\cup\{(\kappa,f)\}$ enjoys the condition in Definition \ref{df:ffthreshold}.\ref{df:QJ.2},
$\gam^{\mathtt{Q}}_{\mS}< b(\kap)<\gam^{\mathtt{Q}}_{\mS}+\mS$,
$\gam^{\mathtt{Q}}_{\mS}<\mathtt{p}_{0}(\kappa)\leq\mathtt{p}_{0}(\rho)$ and
$
(\Theta(\mathtt{Q})\cap M_{\rho})
\cup ((\mathtt{Q}(\mS))_{0}\cap\rho)
\cup\{\mathtt{p}_{0}(\rho)\}
\subset M_{\kappa}$.

We write $\kappa\in H_{\rho}(f,\mathtt{Q},\Theta)$ for $(\kappa,f)\in H_{\rho}(f,\mathtt{Q},\Theta)$.
}
\edf

Let $\Gamma$ be a sequent, 
$\Theta$ a finite set of ordinals,
and
$\mathtt{Q}$ a finite family 
with thresholds.
We define another derivability  relation $(\mathcal{H}_{\gamma},\Theta, \mathtt{Q}
)\vdash^{a}_{c,\vec{d},\xi,\Lambda}\Gamma$, where 
$c$ is a bound on ranks of cut formulas, 
$\vec{d}=(d_{1},d_{2},\ldots,d_{N})$ with
a bound $d_{i}$ on ranks of reflected formulas in $\mathcal{L}_{i}$,
$\xi$ a bound on ordinals in $\mathrm{dom}(\mathtt{Q})$,
and $\Lambda$ a bound on ordinals in $\mathrm{fld}(f)$ for finite functions.

\bdf\label{df:controldercap}
{\rm
Let $\Lambda_{0}=\beta=\psi_{\mI_{N}}(b_{0})$ 
denote a fixed ordinal in Collapsing \ref{lem:Kcollpase.1} such that 
$b_{0}\in C_{b_{0}+1}(\psi_{\mI_{N}}(b_{0}+1))$.

Let $(\mathtt{Q}, \gamma^{\mathtt{Q}}_{\cdot}, \delta^{\mathtt{Q}}_{\cdot})$ be a finite family with thresholds, 
$\Theta\subset \Lambda_{0}$ a finite set of ordinals, and
$a,c,d_{1},d_{2},\ldots,d_{N},\xi,\Lambda$ (unguarded) ordinals with a stable ordinal $\xi\in St$ 
such that $\mathrm{dom}(\mathtt{Q})\subset  \xi$, 
$\Lambda>\Lambda_{0}$ a strongly critical number,
$c\leq d_{1}\leq d_{2}\leq\cdots\leq d_{N}\leq\xi\leq \Lambda_{0}$ and $a<\Lambda$.

Let
$\Gamma=\bigcup\{\Gamma_{\bar{\rho}}^{(\bar{\rho})}: \bar{\rho}\in\prod_{\mS\in\mathrm{dom}(\mathtt{Q})\cap SSt}\mathtt{Q}(\mS)\}$ 
be a set of bounded formulas 
such that $\mathrm{rk}(\Gamma)<\Lambda_{0}$, 
$\mathsf{k}(\Gamma_{\bar{\rho}})\subset M_{\bar{\rho}}$ 
for each cap $ \bar{\rho}\in\prod_{\mS\in\mathrm{dom}(\mathtt{Q})\cap SSt} \mathtt{Q}(\mS)$ over $\mathtt{Q}$.

Let $\vec{d}=(d_{1},d_{2},\ldots,d_{N})$.
$(\mathcal{H}_{\gamma},\Tht, \mathtt{ Q})\vdash^{a}_{c,\vec{d},\xi,\Lambda} \Gamma$ holds
if 
each of 
the following (\ref{eq:controlder_cap_cover}), 
(\ref{eq:controlder_cap2}), (\ref{eq:controlder_cap22}), 
(\ref{eq:fld}) and (\ref{eq:controlder_cap0})
holds, 
and one of the following cases $(\Sigma(\Ome)\mbox{{\rm -rfl}})$, $(\mathrm{dom})$, 
$(\bigvee)$, $(\bigwedge)$, $(cut)$, 
 and 
$(i {\rm -rfl}_{\mS}(g,x,f,\Delta))$ holds:

\begin{equation}
\label{eq:controlder_cap_cover}
\mathrm{srk}(\Gamma)\subset\mathrm{dom}(\mathtt{Q})
\end{equation}

\begin{equation}
\label{eq:controlder_cap2}
\forall\bar{\rho}\in\prod_{\mS\in\mathrm{dom}(\mathtt{Q})\cap SSt}\mathtt{Q}(\mS)
\forall A\in\Gamma_{\bar{\rho}}
\left(
\mathsf{k}(A^{(\bar{\rho})};\mathtt{Q})
 \subset
\mathcal{H}_{\gamma}[
\Theta(\mathtt{Q})
]
\right)
\end{equation}

\beqnarr
&&
\{
a,c,\xi, \Lambda\}\cup \vec{d}
\cup
(\mathrm{dom}(\mathtt{Q}) \cap SSt)
\cup
\{\delta^{\mathtt{Q}}_{\mS}:\mS\in\mathrm{dom}(\mathtt{Q})\}
\nonumber
\\
& \cup &
\bigcup\{\mathrm{fld}(g)
: g\in
(\mathtt{Q})_{1}\}
\subset
\mathcal{H}_{\gamma}[\Theta(\mathtt{Q})]
\label{eq:controlder_cap22}
\eeqnarr

\begin{equation}
\label{eq:fld}
\forall g\in(\mathtt{Q})_{1}
\left[
SC_{\mI_{N}}(\mathrm{fld}(g))\subset \Lambda
\right] 
\end{equation}

\begin{equation}
\label{eq:controlder_cap0}
\forall\mS\in \mathrm{dom}(\mathtt{Q})
(\gamma\leq\gamma^{\mathtt{Q}}_{\mS}) 
\spand 
\gamma^{\mathtt{Q}}_{0}\in\mathcal{H}_{\gamma^{\mathtt{Q}}_{0}}[\Theta(\mathtt{Q})]
\end{equation}

\begin{description}

\item[$(\Sigma(\Ome)\mbox{{\rm -rfl}})$]
There exist ordinals
$a_{\ell}, a_{r}<a$ and an uncapped formula $C\in\Sigma(\mathcal{L}_{0}:\Ome)$
such that 
$c> \Omega$,
$(\mathcal{H}_{\gamma},\Tht, \mathtt{ Q}
)\vdash^{a_{\ell}}_{c,\vec{d},\xi,\Lambda}\Gamma,C^{(\bar{\rho}_{\mathtt{Q}})}$
and
$(\mathcal{H}_{\gamma},\Tht, \mathtt{ Q}
)\vdash^{a_{r}}_{c,\vec{d},\xi,\Lambda}
(\lnot \exists x<\Omega\,C^{(x,\Ome)})^{(\bar{\rho}_{\mathtt{Q}})}, \Gamma$
for the least cap $\bar{\rho}_{\mathtt{Q}}$.

\item[$(\mathrm{dom})$]
There are a stable ordinal $\mS\not\in \mathrm{dom}(\mathtt{Q})$ 
such that $\mS\in\mathcal{H}_{\gamma}[\Theta(\mathtt{Q})]\cap\xi$ and 
$ \mathrm{dom}(\mathtt{Q})\cup\{\mS\}$ is predecessor closed, and
  an ordinal $a_{0}<a$ for which
$(\mathcal{H}_{\gamma},\Theta_{\mathtt{P}}, \mathtt{P}
)\vdash^{a_{0}}_{c,\vec{d},\xi,\Lambda}\Gamma(\mathtt{P})$ holds,
where
$\mathtt{P}=\mathtt{Q}[\mS]$ varies through good extensions of $\mathtt{Q}$ for $\gamma,\Theta$
by some ordinals 
$\gamma^{\mathtt{P}}_{\mS}, \delta^{\mathtt{P}}_{\mS}$ and $(\rho_{\mathtt{P}(\mS)},g_{\mS})$,
cf.\,Definition \ref{df:good_extension}.\ref{df:good_extension.1},
and one of them is very good, cf.\,Definition \ref{df:good_extension}.\ref{df:good_extension.11}.

$\Gamma(\mathtt{P})$ denotes the sequent in 
Definition \ref{df:caphat}.\ref{df:good_extension.3}.

\item[$(\bigvee)$]
There exist  $A\simeq\bigvee(A_{\iota})_{\iota\in J}$, a cap $\bar{\rho}$, 
$a_{0}\in\mathcal{H}_{\gamma}[\Theta(\mathtt{Q})]\cap a$,
and an 
$\iota\in [\bar{\rho}\mathtt{Q}]J$
such that
$A^{(\bar{\rho})}\in\Gamma$,  $A_{\iota}\simeq\bigwedge(B_{\nu})_{\nu\in J_{\iota}}$, 
$\mathsf{k}(A_{\iota}^{(\bar{\rho})};\mathtt{Q}) \subset\mathcal{H}_{\gamma}[\Theta(\mathtt{Q})]$, and
$\mathrm{srk}(A_{\iota})\in\mathrm{dom}(\mathtt{Q})$
for which
$(\mathcal{H}_{\gamma},\Theta_{\mathtt{Q}_{\nu}}\cup\mathsf{k}(B_{\nu}^{(\bar{\rho})};\mathtt{Q}),\mathtt{Q}_{\nu})\vdash^{a_{0}}_{c,\vec{d},\xi,\Lambda}
(\Gamma\cup\{B_{\nu}^{(\bar{\rho})}\})(\mathtt{Q}_{\nu})$ holds for each $\nu\in[\bar{\rho}\mathtt{Q}]J_{\iota}$ and
each good extension
$\mathtt{Q}_{\nu}=\mathtt{Q}[\mathrm{srk}(B_{\nu})]$ of $\mathtt{Q}$ for $\gamma$ and $\Theta$,
where one of the extensions is very good.
cf.\,Definition \ref{df:assigndc}.

\item[$(\bigwedge)$]
There exist 
$A\simeq\bigwedge(A_{\iota})_{ \iota\in J}$,  a cap 
$\bar{\rho}$ and
ordinals $a(\iota)<a$ such that
$\Gamma=\Gamma_{0}\cup\{A^{(\bar{\rho})}\}$ for which
$(\mathcal{H}_{\gamma},\Theta_{\mathtt{Q}_{\iota}}\cup \mathsf{k}(A_{\iota}^{(\bar{\rho})};\mathtt{Q}), \mathtt{Q}_{\iota}
)
\vdash^{a(\iota)}_{c,\vec{d},\xi,\Lambda}
(\Gamma_{0}\cup\{A_{\iota}^{(\bar{\rho})}\})(\mathtt{Q}_{\iota})$ holds
for each $\iota\in[\bar{\rho}\mathtt{Q}]J$ and each good extension 
$\mathtt{Q}_{\iota}=\mathtt{Q}[\mathrm{srk}(A_{\iota})]$ of $\mathtt{Q}$
for $\gamma$  and $\Theta$, where one of the extensions is very good.

\item[$(cut)$] 
There exist
an ordinal $a_{0}<a$, a capped 
$\bigvee$-formula $C^{(\bar{\rho}_{\mathtt{Q}})}$ with the least cap $\bar{\rho}_{\mathtt{Q}}$
such that $C\simeq\bigvee(C_{\iota})_{\iota\in J}$ and $\mathrm{rk}(C)<c$
for which
$(\mathcal{H}_{\gamma},\Tht, \mathtt{Q})\vdash^{a_{0}}_{c,\vec{d},\xi,\Lambda}C^{(\bar{\rho}_{\mathtt{Q}})},\Gamma$ holds, and
$(\mathcal{H}_{\gamma},\Theta_{\mathtt{Q}_{\iota}}\cup\mathsf{k}(C_{\iota}^{(\bar{\rho}_{\mathtt{Q}})};\mathtt{Q}),
\mathtt{Q}_{\iota})\vdash^{a_{0}}_{c,\vec{d},\xi,\Lambda}
(\Gamma\cup\{\lnot C_{\iota}^{(\bar{\rho}_{\mathtt{Q}})}\})(\mathtt{Q}_{\iota})$ holds for each $\iota\in[\bar{\rho}_{\mathtt{Q}}\mathtt{Q}]J$
and each good extension $\mathtt{Q}_{\iota}=\mathtt{Q}[\mathrm{srk}(C_{\iota})]$ of $\mathtt{Q}$
for $\gamma$ and $\Theta$, where one of the extensions is very good.

\item[$(i{\rm -rfl}_{\mS}(g,x,f,\Delta))$] 
There exist a successor $i$-stable ordinal
$c<\mS\in \mathrm{dom}(\mathtt{Q})\cap SSt_{i}$,
an ordinal $x\in \mathrm{supp}(g)$ with $g=(\bar{\rho_{\mathtt{Q}}}(\mS))_{1}$,
a special finite function $f:\mI_{N}\to\Gamma(\mI_{N})$ with base $\mI_{N}$ such that 
$\min\{s(g), s(f)\}\geq\mS+1$, and
a set $\Delta\subset\mathcal{L}_{i}$
 of $\bigvee$-formulas $\del$ such that

\benu
\item[(r1)]
$\mS\leq \mathrm{rk}(\Delta)<d_{i}$,
$\Delta\subset\bigvee(x):=\{\delta: \mathrm{rk}(\delta)<x,\, \delta \mbox{ {\rm is a}}
\bigvee\mbox{{\rm -formula}}\}$
and
$\forall\mU\in\mathrm{dom}(\mathtt{Q})\cap \bigcup_{j\geq i}SSt_{j}(\mS<\mU \Rarw \mathrm{rk}(\Delta)<\mU)$, cf.\,(\ref{eq:stbl_I_1}).

\item[(r2)]
$\mathrm{fld}(f)\subset \mathcal{H}_{\gamma}[\Theta(\mathtt{Q})]$,
$SC_{\mI_{N}}(\mathrm{fld}(f))\subset\Lambda$
and 
$
f_{x}\leq g_{x} \,\&\, f<_{\mI_{N}}^{x} g^{\prime}(x)$, 
cf.\,Definition \ref{df:notationsystem}.\ref{df:notationsystem.6}.

\eenu
Furthermore there is an ordinal $a_{0}\in\mathcal{H}_{\gamma}[\Theta(\mathtt{Q})]\cap a$ for which the following holds:

Let $\mT$ be a successor $j$-stable ordinal 
such that $j\geq i$, $\mS\leq\mT\leq\mathrm{rk}(\Delta)$ and $\min\{s(g), s(f)\}\geq\mT+1$.

Let $\mathtt{P}=\mathtt{Q}[\mT]$ be a good extension of $\mathtt{Q}$ for $\gamma,\Theta$ by ordinals
$\gamma^{\mathtt{P}}_{\mT},\delta^{\mathtt{P}}_{\mT}$ and 
$\bar{\rho}_{\mathtt{P}}(\mT)=(\rho_{\mathtt{P}(\mT)},g)$, where 
$(\bar{\rho_{\mathtt{Q}}}(\mS))_{1}=g=(\bar{\rho_{\mathtt{P}}}(\mT))_{1}$.

For each such extension $\mathtt{P}=\mathtt{Q}[\mT]$, the following 
 (r3) and (r4) hold:

\benu

 \item[(r3)]
 Let $\delta\simeq\bigvee(\delta_{\iota})_{\iota\in J_{\delta}}$ for $\delta\in\Delta$.
 $(\mathcal{H}_{\gamma},\Theta_{\mathtt{P}_{\iota}}\cup\mathsf{k}(\delta_{\iota}^{(\bar{\rho}_{\mathtt{P}})};\mathtt{P}), \mathtt{P}_{\iota}
 )\vdash^{a_{0}}_{c,\vec{d},\xi,\Lambda}
 (\Gamma\cup\{ \lnot\delta_{\iota}^{(\bar{\rho}_{\mathtt{P}})}\})(\mathtt{P}_{\iota})$ holds for 
each $\delta\in\Delta$, each $\iota\in[\bar{\rho}_{\mathtt{P}}\mathtt{P}]J_{\delta}$, and each good extension
$\mathtt{P}_{\iota}=\mathtt{P}[\mathrm{srk}(\delta_{\iota})]$ of $\mathtt{P}$ for $\gamma$ and
$\Theta$, where one of the extensions is very good.

\item[(r4)]
$
(\mathcal{H}_{\gamma},\Theta_{\mathtt{P}},\mathtt{P}^{\sigma}
)\vdash^{a_{0}}_{c,\vec{d},\xi,\Lambda}\Delta^{(\bar{\rho}_{\sigma})},\Gamma(\mathtt{P})$ holds
for each $\sigma\in H_{\rho_{\mathtt{P}(\mT)}}(f,\mathtt{P},\Theta)$,
where 
$\mathtt{P}^{\sigma}=\mathtt{P}\cup\{(\mT,(\sigma,f))\}$ denotes
 a finite family with thresholds obtained from $\mathtt{P}$ by adding $(\sigma,f)$ 
 to $\mathtt{P}(\mT)$, cf.\,Definition \ref{df:caphat}.\ref{df:caphat.11} ,
 and
 $\bar{\rho}_{\sigma}=\bar{\rho}_{\mathtt{P}^{\sigma}}=\bar{\rho}_{\mathtt{P}}[(\sigma,f)/(\rho_{\mathtt{P}(\mT)},g)]$,
 where $\bar{\rho}_{\mathtt{P}}(\mT)=(\rho_{\mathtt{P}(\mT)},g)$ and $\bar{\rho}_{\sigma}(\mT)=(\sigma,f)$.

\eenu

\end{description}

}
\edf

We write $\mathtt{Q}^{\sigma}=\mathtt{Q}\cup\{(\mS,\sigma)\}$ for 
$\mathtt{Q}\cup\{(\mS,(\sigma,f))\}$, and
 $\bar{\rho}_{\sigma}=\bar{\rho}_{\mathtt{Q}^{\sigma}}=\bar{\rho}_{\mathtt{Q}}[\sigma/\rho]$
 for $\bar{\rho}_{\mathtt{Q}}[(\sigma,f)/(\rho,g)]$.

Let $(\mathcal{H}_{\gamma},\Theta,\mathtt{Q})\vdash^{a}_{c,\vec{d},\xi,\Lambda}\Gamma$ for a set $\Gamma$ of bounded formulas.
We obtain
$\mathrm{rk}(\Gamma)<\xi$
by (\ref{eq:controlder_cap_cover}) and Proposition \ref{prp:pd_closed_cover_I_1}.\ref{prp:pd_closed_cover_I_1.4}.

An inversion lemma for the calculus $\vdash$ seems not to hold due to the inference $(\mathrm{dom})$, cf.\,subsection \ref{subsec:preview}.
In this reason, the major formula $A^{(\bar{\rho})}$ need not to be in the upper sequents $\Gamma_{0}\cup\{A_{\iota}^{(\bar{\rho})}\}$
in inferences $(\bigwedge)$.
Furthermore
we perform in advance an inversion step
in  inferences $(\bigvee)$, $(cut)$ and (r3) of $(i {\rm -rfl}_{\mS}(g,x,f,\Delta))$.

\blem\label{lem:controlder_cap0}
Let 
$(\mathcal{H}_{\gamma},\Theta,\mathtt{Q})\vdash^{a}_{c,\vec{d},\xi,\Lambda}\Gamma$.
Then
$\mathcal{H}_{\gamma}[\Theta(\mathtt{Q})]\subset M_{\rho}\cap \bigcap\{M_{\mathtt{Q}(\mT)}:\mS>\mT\in\mathrm{dom}(\mathtt{Q})\}$ for
$\rho\in(\mathtt{Q}(\mS))_{0}$.
\elem
\bprf
Let $\xi>\mS>\mT\in\mathrm{dom}(\mathtt{Q})$.
We have 
$\gamma\leq\gamma^{\mathtt{Q}}_{\mS}< \min\{\mathtt{p}_{0}(\rho),\gamma^{\mathtt{Q}}_{\mT}\}$ by
(\ref{eq:controlder_cap0}) and Definition \ref{df:ffthreshold}.\ref{df:QJ.2}.
This yields
$\mathcal{H}_{\gamma}[M_{\rho}]\subset M_{\rho}$ and $\mathcal{H}_{\gamma}[M_{\mathtt{Q}(\mT)}]\subset M_{\mathtt{Q}(\mT)}$
by Lemma \ref{prp:EK2}.\ref{prp:EK2.1}.
On the other hand we have $\Theta(\mathtt{Q})=\Theta\cap M_{\mathtt{Q}}\subset M_{\rho}\cap M_{\mathtt{Q}(\mT)}$.
\eprf

\blem\label{lem:guard_change}(Cf.\,(\ref{eq:KppiNlowerCase1b.12_I_1}).)
Let $g,k$ be special finite functions with base $\mI_{N}$ such that $g\leq k$, and
$(\rho,k)\in\mathtt{Q}(\mS)$.
$\mathtt{Q}^{[g/k]}=\mathtt{Q}^{[(\rho,g)/(\rho,k)]}$ denotes a finite family obtained from $\mathtt{Q}$ by replacing
$(\rho,k)\in\mathtt{Q}(\mS)$ by $(\rho,g)\in\mathtt{Q}^{[g/k]}(\mS)$.
$\bar{\rho}[g/k]$ denotes the cap over $\mathtt{Q}^{[g/k]}$ 
obtained from a cap $\bar{\rho}$ over $\mathtt{Q}$ by replacing 
a guarded ordinal $\bar{\rho}(\mS)=(\rho,k)$ by $(\rho,g)$.
Specifically $(\bar{\rho}[g/k])(\mS)=(\rho,k)$.
Let $(A^{(\bar{\rho})})^{[g/k]}:\equiv A^{(\bar{\rho}[g/k])}$ for a capped formula $A^{(\bar{\rho})}$,
and $\Gamma^{[g/k]}:=\{(A^{(\bar{\rho})})^{[g/k]}: A^{(\bar{\rho})}\in\Gamma\}$
for a sequent $\Gamma$.

Let $(\mathcal{H}_{\gamma},\Theta,\mathtt{Q}^{[g/k]})\vdash^{a}_{c,\vec{d},\xi,\Lambda}\Gamma^{[g/k]}$.
Then 
$(\mathcal{H}_{\gamma},\Theta,\mathtt{Q})\vdash^{a}_{c,\vec{d},\xi,\Lambda}\Gamma$ holds.
\elem
\bprf 
By induction on $a$.
Consider the case when the last inference in $(\mathcal{H}_{\gamma},\Theta,\mathtt{Q})\vdash^{a}_{c,\vec{d},\xi,\Lambda}\Gamma^{[g/k]}$ is an $(i{\rm -rfl}_{\mS}(g,x,f,\Delta))$.
We have $f_{x}\leq g_{x} \,\&\, f<_{\mI_{N}}^{x} g^{\prime}(x)$ by (r2).
We obtain $f_{x}\leq k_{x} \,\&\, f<_{\mI_{N}}^{x} k^{\prime}(x)$ by $g\leq k$ and
Proposition \ref{prp:idless}.
IH followed by an $(i {\rm -rfl}_{\mS}(k,x,f,\Delta))$ yields the lemma
\eprf

\blem\label{lem:tautology.cap}
{\rm (Tautology)}
Let $\mathtt{Q}$ be a finite family with thresholds 
and $\Theta$ a finite set of ordinals 
such that $\gamma\in\mathcal{H}_{\gamma}[\Theta(\mathtt{Q})]$, $\mathrm{dom}(\mathtt{Q})\subset\xi$ and 
(\ref{eq:controlder_cap0}) is enjoyed.

Let  $\sigma=\rho_{\mathtt{Q}(\mS)}$ for an $\mS\in \mathrm{dom}(\mathtt{Q})$, 
and
$A$ be a formula such that $\mathsf{k}(A)\subset  \mathcal{H}_{\gamma}[\Theta(\mathtt{Q})]$
and $\mathrm{srk}(A)\in\mathrm{dom}(\mathtt{Q})$.
Let $d=\mathrm{rk}(A)$ and $\vec{0}=(0,\ldots,0)$.

\benu
\item\label{lem:tautology.cap1}
$(\mathcal{H}_{\gamma},\Theta,\mathtt{Q}
)\vdash^{2d}_{0,\vec{0},\xi,\Lambda}
\lnot A^{(\bar{\rho}_{\mathtt{Q}})}, A^{(\bar{\rho}_{\mathtt{Q}})}$ holds.

$(\mathcal{H}_{\gamma},\Theta\cup\mathsf{k}(A_{\iota}),\mathtt{Q}_{\iota}
)\vdash^{2d_{\iota}+1}_{0,\vec{0},\xi,\Lambda}
(\{\lnot A_{\iota}^{(\bar{\rho}_{\mathtt{Q}})}, A^{(\bar{\rho}_{\mathtt{Q}})}\})(\mathtt{Q}_{\iota})$ holds for each $\iota\in [\bar{\rho}_{\mathtt{Q}}\mathtt{Q}]J$ when
$A\simeq\bigvee(A_{\iota})_{\iota\in J}$ and $d_{\iota}=\mathrm{rk}(A_{\iota})$,
where $\mathtt{Q}_{\iota}=\mathtt{Q}[\mathrm{srk}(A_{\iota})]$ is a good extension of $\mathtt{Q}$.

\item\label{lem:tautology.cap3}
Let $\rho\in(\mathtt{Q}(\mS))_{0}$ be such that $\sigma<\rho$, and 
$\bar{\rho}=\bar{\rho}_{\mathtt{Q}}[\rho/\sigma]$.
Then
$(\mathcal{H}_{\gamma},\Theta,\mathtt{Q}
)\vdash^{2d}_{0,\vec{0},\xi,\Lambda}
\lnot A^{(\bar{\rho}_{\mathtt{Q}})}, (A^{[\sigma/\mathbb{S}]})^{(\bar{\rho})}$ holds.

\eenu
\elem
\bprf
By induction on $d=\mathrm{rk}(A)$. 
By $\mathsf{k}(A)\subset  \mathcal{H}_{\gamma}[\Theta(\mathtt{Q})]$, we obtain
$\mathsf{k}(A)\subset M_{\mathtt{Q}}$ by Lemma \ref{lem:controlder_cap0}.
In the proof let us write $\vdash^{a}_{c}$ for $\vdash^{a}_{c,\vec{d},\xi,\Lambda}$.
\\
\ref{lem:tautology.cap}.\ref{lem:tautology.cap1}.
We have 
$\mathrm{inv}^{(\mathtt{u})}(A^{(\bar{\rho}_{\mathtt{Q}})};\mathtt{Q})\equiv A$ for the least cap $\bar{\rho}_{\mathtt{Q}}$
since there is no $\tau\in(\mathtt{Q}(\mT))_{0}$ such that $\tau<\bar{\rho}_{\mathtt{Q}}(\mT)=\rho_{\mathtt{Q}(\mT)}$.
Hence
$\mathsf{k}(A^{(\bar{\rho}_{\mathtt{Q}})};\mathtt{Q})=\mathsf{k}(A)\subset  \mathcal{H}_{\gamma}[\Theta(\mathtt{Q})]\cap M_{\mathtt{Q}}$.
We obtain
$d\in\mathcal{H}_{\gamma}[\Theta(\mathtt{Q})]$ by Proposition \ref{lem:rank}.\ref{lem:rank1} for (\ref{eq:controlder_cap22}).

Let $A\simeq\bigvee(A_{\iota})_{\iota\in J}$ and $A_{\iota}\simeq\bigwedge(B_{\nu})_{\nu\in J_{\iota}}$ for each $\iota\in J$.
Let $\iota\in[\bar{\rho}_{\mathtt{Q}}\mathtt{Q}]J$ and $\nu\in[\bar{\rho}_{\mathtt{Q}}\mathtt{Q}]J_{\iota}$.
Let $d_{\iota}=\mathrm{rk}(A_{\iota})$ and $d_{\nu}=\mathrm{rk}(B_{\nu})$.
We have $d_{\nu}<d_{\iota}<d$.
We obtain
$\mathsf{k}(B_{\nu})\subset M_{\mathtt{Q}}$, and
 $
\mathsf{k}(B_{\nu})\subset\mathcal{H}_{\gamma}[(\Theta\cup\mathsf{k}(B_{\nu}))(\mathtt{Q})]$.

IH yields
$(\mathcal{H}_{\gamma},\Theta\cup \mathsf{k}(B_{\nu},A_{\iota}),\mathtt{Q}_{\iota\nu}
)\vdash^{2d_{\iota}}_{0}
(\{\lnot A_{\iota}^{(\bar{\rho}_{\mathtt{Q}})}, B_{\nu}^{(\bar{\rho}_{\mathtt{Q}})}\})(\mathtt{Q}_{\iota\nu})$
for good extensions $\mathtt{Q}_{\iota\nu}=\mathtt{Q}_{\iota}[\mathrm{srk}(B_{\nu})]$.
A $(\bigvee)$ yields $(\mathcal{H}_{\gamma},\Theta\cup \mathsf{k}(A_{\iota}),\mathtt{Q}_{\iota}
)\vdash^{2d_{\iota}+1}_{0}
(\{\lnot A_{\iota}^{(\bar{\rho}_{\mathtt{Q}})}, A^{(\bar{\rho}_{\mathtt{Q}})}\})(\mathtt{Q}_{\iota})$, and
$(\mathcal{H}_{\gamma},\Theta,\mathtt{Q}
)\vdash^{2d}_{0}
\lnot A^{(\bar{\rho}_{\mathtt{Q}})}, A^{(\bar{\rho}_{\mathtt{Q}})}$ by a $(\bigwedge)$.
\\
\ref{lem:tautology.cap}.\ref{lem:tautology.cap3}.
We have $\sigma\in\mathrm{dom}(\mathtt{Q})$ when $\sigma\in LSt$ by 
Definition \ref{df:ffthreshold}.\ref{df:QJ.3}.
By Proposition \ref{prp:pd_closed_cover_I_1}.\ref{prp:pd_closed_cover_I_1.31}
we have $\mathrm{srk}(A^{[\sigma/\mS]})\in\{\mathrm{srk}(A)\}\cup
(\mathrm{cl}(\mS)\cap\sigma)\cup(\{\sigma\}\cap St)$
and $\{\mathrm{srk}(A)\}\cup(\mathrm{cl}(\mS)\cap\sigma)\cup(\{\sigma\}\cap St)\subset\mathrm{dom}(\mathtt{Q})$ by the assumption.
Hence (\ref{eq:controlder_cap_cover}) is enjoyed
in $(\mathcal{H}_{\gamma},\Theta,\mathtt{Q}
)\vdash^{2d}_{0}
\lnot A^{(\bar{\rho}_{\mathtt{Q}})}, (A^{[\sigma/\mathbb{S}]})^{(\bar{\rho})}$.
Note that $\mathrm{cl}(\sigma)=\{\sigma\}\cup(\mathrm{cl}(\mS)\cap\sigma)$ when $\sigma\in\Psi_{\mS}\cap LSt$.

Let $B\equiv A^{[\sigma/\mathbb{S}]}$. We have $\mathrm{rk}(B)<\mS$.
We claim that $\mathrm{inv}^{(\mathtt{u})}(B^{(\bar{\rho})};\mathtt{Q})\equiv A$.
We have $\mathrm{inv}^{(\mathtt{u})}(B^{(\bar{\rho})};\mathtt{Q})\equiv B[\sigma/\mS]^{-1}\equiv A$ 
if $B\not\equiv A$.
Otherwise we have $\mathrm{inv}^{(\mathtt{u})}(B^{(\bar{\rho})};\mathtt{Q})\equiv A^{[\sigma/\mathbb{S}]}\equiv A$ since 
$B[\sigma/\mS]^{-1}\equiv B\equiv A$ and $(\bar{\rho}(\mT))_{0}=\rho_{\mathtt{Q}(\mT)}$ for $\mT\neq\mS$.
In each case we obtain $\mathsf{k}(B^{(\bar{\rho})};\mathtt{Q})=\mathsf{k}(A)\subset\mathcal{H}_{\gamma}[\Theta(\mathtt{Q})] \cap M_{\mathtt{Q}}$.
Hence (\ref{eq:controlder_cap2}) 
and (\ref{eq:controlder_cap22}) 
are enjoyed.

Let $A\simeq\bigvee(A_{\iota})_{\iota\in J}$ or  $A\simeq\bigwedge(A_{\iota})_{\iota\in J}$.
By Proposition \ref{lem:assigncollaps} we obtain
 $A^{[\sigma/\mathbb{S}]}\simeq\bigvee(A^{[\sigma/\mathbb{S}]}_{\iota})_{\iota\in [\sigma]J}$ or  
 $A^{[\sigma/\mathbb{S}]}\simeq\bigwedge(A^{[\sigma/\mathbb{S}]}_{\iota})_{\iota\in [\sigma]J}$.
 Let $I=\{\iota^{[\sigma/\mathbb{S}]}: \iota\in[\sigma]J\}$.

We claim that for $\nu=\iota^{[\sigma/\mathbb{S}]}$ and $B_{\nu}\equiv A_{\iota}^{[\sigma/\mathbb{S}]}$
\beqn\label{eq:tautology.cap3}
\mathsf{k}(A_{\iota}^{(\bar{\rho}_{\mathtt{Q}})};\mathtt{Q})=\mathsf{k}(B_{\nu}^{(\bar{\rho})};\mathtt{Q}) \mbox{ and }
\iota\in [\bar{\rho}_{\mathtt{Q}}\mathtt{Q}]J \Lrarw \nu=\iota^{[\sigma/\mathbb{S}]}\in [\bar{\rho} \mathtt{Q}]I
\eeqn
where
 $\iota\in [\bar{\rho}_{\mathtt{Q}}\mathtt{Q}]J$ 
 iff $\mathsf{k}(A_{\iota}^{(\bar{\rho}_{\mathtt{Q}})};\mathtt{Q})=\mathsf{k}(A_{\iota})\subset M_{\mathtt{Q}}$, while
 $\iota^{[\sigma/\mathbb{S}]}\in [\bar{\rho} \mathtt{Q}]I$ iff  
 $\mathsf{k}(A_{\iota})\subset M_{\sigma}$ and
 $\mathsf{k}(B_{\nu}^{(\bar{\rho})};\mathtt{Q})=\mathsf{k}(A_{\iota})\subset M_{\mathtt{Q}}$.
 Hence the claim (\ref{eq:tautology.cap3}) follows.

For each $\iota\in [\bar{\rho}_{\mathtt{Q}}\mathtt{Q}]J$, we have $d_{\iota}=\mathrm{rk}(A_{\iota})<d$ by Proposition \ref{lem:rank}.\ref{lem:rank2}.
IH yields
$(\mathcal{H}_{\gamma},\Theta\cup \mathsf{k}(A_{\iota},C_{\mu}),\mathtt{Q}
)\vdash^{2d_{\iota}}_{0}
\lnot C_{\mu}^{(\bar{\rho}_{\mathtt{Q}})}, B_{\nu}^{(\bar{\rho})}$ if $A_{\iota}\simeq \bigvee(C_{\mu})_{\mu\in J_{\iota}}$.
When $A_{\iota}\simeq \bigwedge(C_{\mu})_{\mu\in J_{\iota}}$, we have
$(\mathcal{H}_{\gamma},\Theta\cup \mathsf{k}(A_{\iota},C_{\mu}),\mathtt{Q}
)\vdash^{2d_{\iota}}_{0}
\lnot A_{\iota}^{(\bar{\rho}_{\mathtt{Q}})}, (C_{\mu}^{[\sigma/\mS]})^{(\bar{\rho})}$.
In each case
$(\mathcal{H}_{\gamma},\Theta\cup \mathsf{k}(A_{\iota}),\mathtt{Q}
)\vdash^{2d_{\iota}+1}_{0}
\lnot A_{\iota}^{(\bar{\rho}_{\mathtt{Q}})}, B_{\nu}^{(\bar{\rho})}$.

Let $A\simeq\bigvee_{\iota\in J}\bigwedge_{\mu\in J_{\iota}}C_{\mu}$ with $C_{\mu}\simeq\bigvee_{\alpha\in J_{\mu}}D_{\alpha}$, and $d_{\mu}=\mathrm{rk}(C_{\mu})$.
IH yields 
$(\mathcal{H}_{\gamma},\Theta\cup  \mathsf{k}(D_{\alpha}, C_{\mu}),\mathtt{Q}
)\vdash^{2d_{\mu}}_{0}
\lnot D_{\alpha}^{(\bar{\rho}_{\mathtt{Q}})}, (C_{\mu}^{[\sigma/\mathbb{S}]})^{(\bar{\rho})}$ for
$\iota\in [\bar{\rho}_{\mathtt{Q}}\mathtt{Q}]J$, $\mu\in[\bar{\rho}_{\mathtt{Q}}\mathtt{Q}]J_{\iota}$ and $\alpha\in[\bar{\rho}_{\mathtt{Q}}\mathtt{Q}]J_{\mu}$.
We see
$\mathsf{k}(C_{\mu})=\mathsf{k}(C_{\mu}^{(\bar{\rho}_{\mathtt{Q}})};\mathtt{Q})=\mathsf{k}((C_{\mu}^{[\sigma/\mathbb{S}]})^{(\bar{\rho})} ;\mathtt{Q})$
as in (\ref{eq:tautology.cap3}).

 We obtain
 $(\mathcal{H}_{\gamma},\Theta\cup  \mathsf{k}(D_{\alpha})\cup \mathsf{k}((C_{\mu}^{[\sigma/\mathbb{S}]})^{(\bar{\rho})} ;\mathtt{Q}),\mathtt{Q}
)\vdash^{2d_{\mu}}_{0}
\lnot D_{\alpha}^{(\bar{\rho}_{\mathtt{Q}})}, (C_{\mu}^{[\sigma/\mathbb{S}]})^{(\bar{\rho})}$.
Let $\Theta_{0}=\Theta\cup  \mathsf{k}(A_{\iota})\cup \mathsf{k}((C_{\mu}^{[\sigma/\mathbb{S}]})^{(\bar{\rho})};\mathtt{Q})$.
We have 
$\mathsf{k}(C_{\mu}^{(\bar{\rho}_{\mathtt{Q}})};\mathtt{Q})=\mathsf{k}(C_{\mu})\subset
\mathcal{H}_{\gamma}[\Theta_{0}(\mathtt{Q})]\cap M_{\mathtt{Q}}$.
On the other hand we have
$\mathsf{k}(A_{\iota})=\mathsf{k}((A_{\iota}^{[\sigma/\mathbb{S}]})^{(\bar{\rho})};\mathtt{Q})\subset
\mathcal{H}_{\gamma}[\Theta_{0}(\mathtt{Q})]\cap M_{\mathtt{Q}}$. 
Hence
\[
\infer[(\bigwedge)]{
  (\mathcal{H}_{\gamma},\Theta,\mathtt{Q})\vdash^{2d+1}_{0,\xi}\lnot A^{(\bar{\rho}_{\mathtt{Q}})}, B^{(\bar{\rho})}
  }
{
\infer[(\bigvee)]{
  (\mathcal{H}_{\gamma},\Theta\cup  \mathsf{k}(A_{\iota}),\mathtt{Q})\vdash^{2d}_{0,\xi}\lnot A_{\iota}^{(\bar{\rho}_{\mathtt{Q}})}, B^{(\bar{\rho})}
 }
{
 \infer[(\bigvee)]{
(\mathcal{H}_{\gamma},\Theta\cup  \mathsf{k}(A_{\iota})\cup \mathsf{k}((C_{\mu}^{[\sigma/\mathbb{S}]})^{(\bar{\rho})} ;\mathtt{Q}),\mathtt{Q}
)\vdash^{2d_{\iota}}_{0,\xi}
\lnot A_{\iota}^{(\bar{\rho}_{\mathtt{Q}})}, (C_{\mu}^{[\sigma/\mathbb{S}]})^{(\bar{\rho})}
 }
 {
  (\mathcal{H}_{\gamma},\Theta\cup  \mathsf{k}(A_{\iota},D_{\alpha})\cup \mathsf{k}((C_{\mu}^{[\sigma/\mathbb{S}]})^{(\bar{\rho})} ;\mathtt{Q}),\mathtt{Q}
)\vdash^{2d_{\mu}}_{0,\xi}
\lnot D_{\alpha}^{(\bar{\rho}_{\mathtt{Q}})}, (C_{\mu}^{[\sigma/\mathbb{S}]})^{(\bar{\rho})}
  }
 }
}
\]
where $\iota$ ranges over $\iota\in [\bar{\rho}_{\mathtt{Q}}\mathtt{Q}]J$ by
(\ref{eq:tautology.cap3}).
The case $A\simeq\bigwedge(A_{\iota})_{\iota\in J}$ is similarly seen.
\eprf

\section{Elimination of stable ordinals}
\label{sec:proofonestep}

\subsection{Capping and recapping}\label{subsec:capping}
In this subsection the relation $\vdash^{*}$ is embedded in $\vdash$ by putting caps on formulas,
and then caps are changed to smaller caps.

\blem\label{lem:capping}{\rm (Capping)}
Let
$\Gamma\subset\Delta_{0}(\mathcal{L}_{N+1})$ be a set of uncapped formulas
with $\rk(\Gamma)<\Lambda_{0}$, where 
$\Lambda_{0}=\beta=\psi_{\mI_{N}}(b_{0})$ 
denotes a fixed ordinal in Collapsing \ref{lem:Kcollpase.1} such that 
$b_{0}\in C_{b_{0}+1}(\psi_{\mI_{N}}(b_{0}+1))$.

Let
$
(\mathcal{H}_{b_{0}},\Theta;\mathtt{D}
)\vdash^{* a}_{\Lambda_{0}}\Gamma
$, where $0\in\mathtt{D}$,
$\mathsf{k}(\Gamma)\cup\mathtt{D}\subset\Lambda_{0}$, cf.\,(\ref{eq:Kcollpase.1}),
$\Theta\subset C_{b_{0}+1}(\psi_{\mI_{N}}(b_{0}+1))$ 
and
$a\leq\Lambda_{0}<\Lambda_{1}=\psi_{\mI_{N}}(b_{0}+1)<\mI_{N}\leq b_{0}$.

Let $\delta^{\mathtt{Q}}_{\cdot}$ be a function
 $\mathtt{D}\ni\mS\mapsto\delta^{\mathtt{Q}}_{\mS}$
such that
$\delta^{\mathtt{Q}}_{\mT}+\omega^{1+a}\leq\delta^{\mathtt{Q}}_{\mS}$
and $\delta^{\mathtt{Q}}_{\mS}\in\mathcal{H}_{\gamma}[\Theta_{\delta,\mS}]$ 
for every 
$\{\mS<\mT\}\subset\mathtt{D}$ and 
$\omega^{a}\leq\delta^{\mathtt{Q}}_{\mS}\in C_{b_{0}+1}(\psi_{\mI_{N}}(b_{0}+1))\cap\omega^{\Lambda_{0}+1}$
for $0\neq\mS\in\mathtt{D}$, where
$\Theta_{\delta,\mS}=\Theta\cup\{\delta^{\mathtt{Q}}_{\mU}: \mS<\mU\in\mathtt{D}\}$.
Let $\delta^{\mathtt{Q}_{0}}_{0}=\omega^{\Lambda_{0}+1}$.
For $\mS\in\mathtt{D}$
let $\Lambda^{\mathtt{Q}}_{\mS}=\psi_{\mI_{N}}(b_{0}+\delta^{\mathtt{Q}}_{\mS})$.

$g_{0}:\mI_{N}\to\Gamma(\mI_{N})$ 
denotes a special finite function with base $\mI_{N}$
 such that
$\mathrm{supp}(g_{0})=\{\Lambda_{0}\}$ with $g_{0}(\Lambda_{0})=\mI_{N}\cdot 3$.

Let $\gamma=b_{0}+\omega^{\Lambda_{0}+1}$, and
$\Theta_{\gamma,\mS}=\Theta\cup\{\gamma^{\mathtt{Q}}_{\mU}: \mS<\mU\in\mathtt{D}\}$.
Let $\gamma^{\mathtt{Q}}_{\cdot}$ be a threshold function
 $\mathtt{D}\ni\mS\mapsto\gamma^{\mathtt{Q}}_{\mS}$
such that $\gamma\leq\gamma^{\mathtt{Q}}_{\mS}\in\mathcal{H}_{\gamma}[\Theta_{\gamma,\mS}]$,
$\gamma^{\mathtt{Q}}_{\mS}=\alpha_{\mS}+\beta_{\mS}$ for an $\alpha_{\mS}$ and 
an additive principal number $\mI_{N}>\beta_{\mS}\geq \omega^{\Lambda_{0}(1+a)}$, 
$\gamma^{\mathtt{Q}}_{\cdot}$ together with $\gamma$ 
has gaps $\omega^{\Lambda_{0}(1+a)}$,
and
$\mathrm{fld}(g_{0})=\{\Lambda_{0},\mI_{N}\cdot 3\}\subset C_{0}(SC(\gamma^{\mathtt{Q}}_{\mS}))$,
cf.\,(\ref{eq:notationsystem.5}),
Let $\gamma^{\mathtt{Q}}_{0}=\alpha_{0}+\beta_{0}\in\mathcal{H}_{\gamma}[\Theta]$ with 
$\alpha_{0}=\max(\{\gamma\}\cup\mathsf{k}_{0}(\Theta))$ and 
$\beta_{0}=\omega^{\Lambda_{0}\cdot 2}$,
where $k_{\mS}(\Theta)=\max(\{0\}\cup \bigcup\{K_{\mS}(\alpha):\alpha\in\Theta\})$.
Assume $\mathsf{k}_{\mS}(\Theta)<\gamma^{\mathtt{Q}}_{\mS}$.

Let
$\rho_{\mS}=\psi_{\mathbb{S}}^{g_{0}}(\gamma^{\mathtt{Q}}_{\mS}+1)$, and
$\mathtt{Q}$ be a finite family with thresholds $\gamma^{\mathtt{Q}}_{\cdot}$ such that
$\mathrm{dom}(\mathtt{Q})= \mathtt{D}\subset\Lambda_{0}$, cf.\,(\ref{eq:Kcollpase.1}) and
$\mathtt{Q}(\mS)=\{(\rho_{\mS},g_{0})\}$ for each 
$\mS\in \mathrm{dom}(\mathtt{Q})\cap SSt= \mathtt{D}\cap SSt$.
Let $\bar{\rho}\in\prod_{\mS\in\mathtt{D}\cap SSt}\mathtt{Q}(\mS)$ be the cap with 
$\bar{\rho}(\mS)=(\rho_{\mS},g_{0})$ for every $\mS\in\mathtt{D}\cap SSt$.
Assume 
$\Theta_{\mathtt{Q}}=\Theta\cup\{\delta^{\mathtt{Q}}_{\mS}:\mS\in\mathrm{dom}(\mathtt{Q})\}
\subset M_{\bar{\rho}}=\bigcap_{\mS\in \mathtt{D}\cap SSt}M_{\rho_{\mS}}$.

Then
$(\mathcal{H}_{\gamma},\Theta_{\mathtt{Q}}, \mathtt{Q})
\vdash^{\Lambda_{0}+2a}_{\Lambda_{0},\vec{d},\Lambda_{0},\Lambda_{1}}
\Gamma^{(\bar{\rho})}$ holds for $c=\xi=\Lambda_{0}$, $\vec{d}=(\Lambda_{0},\ldots,\Lambda_{0})$ and
$\Gamma^{(\bar{\rho})}=\{A^{(\bar{\rho})}:A\in\Gamma\}$.

\elem
\bprf
By induction on $a$.
Let us 
write $\vdash^{a}_{\Lambda_{0}}$ for $\vdash^{a}_{\Lambda_{0}, \vec{d}, \Lambda_{0}, \Lambda_{1}}$ in the proof.

For each $\mS\in\mathtt{D}$ we have $\rho_{\mS}=\rho_{\mathtt{Q}(\mS)}$ and
$\mathsf{k}(A^{(\bar{\rho})};\mathtt{Q})=\mathsf{k}(A)$.
On the other hand we have $\Theta_{\mathtt{Q}}\subset M_{\mathtt{Q}}$ by the assumption.
Let $A\in\Gamma$.
We obtain $\mathsf{k}(A)\subset\mathcal{H}_{\gamma}[\Theta]\subset M_{\mathtt{Q}}$ by (\ref{eq:controlder*1}) and 
$\Theta\subset M_{\mathtt{Q}}$, where $\Theta_{\mathtt{Q}}=\Theta_{\mathtt{Q}}\cap M_{\mathtt{Q}}$.
(\ref{eq:controlder_cap2}) and (\ref{eq:controlder_cap22}) in 
$(\mathcal{H}_{\gamma}, \Theta, \mathtt{ Q})
\vdash^{\Lambda_{0}+2a}_{\Lambda_{0}}
\Gamma^{(\bar{\rho})}$ 
follow from (\ref{eq:controlder*1}) in $(\mathcal{H}_{b_{0}},\Theta; \mathtt{D}
)\vdash^{*a}_{\Lambda_{0}}\Gamma
$.
We have $\mathrm{dom}(\mathtt{Q})=\mathtt{D}\subset\Lambda_{0}$.
Let $\mS\in \mathrm{dom}(\mathtt{Q})$. Then by (\ref{eq:controlder*1}) we obtain $\mS\in\mathcal{H}_{\gam}[\Tht]$.

We have $SC_{\mI_{N}}(\mathrm{fld}(m(\rho_{\mS})))=\{\Lambda_{0},3\}\subset
\psi_{\mI_{N}}(b_{0}+1)=\Lambda_{1}=\Lambda^{\mathtt{Q}}_{\Lambda_{0}}$.

For $\mT<\mS$,
$\rho_{\mS}\in M_{\rho_{\mT}}$ 
follows from $\{\mS,\Lambda_{0},\gamma^{\mathtt{Q}}_{\mS}\}\subset\calh_{\gamma}[\Tht]\subset M_{\rho_{\mT}}$
and $\gamma\leq \gamma^{\mathtt{Q}}_{\mS}<\gamma_{\mS}^{\mathtt{ Q}}+\mS<\gamma_{\mT}^{\mathtt{Q}}<\mathtt{p}_{0}(\rho_{\mT})$ by $\omega^{\Lambda_{0}(1+a)}\geq\Lambda_{0}>\mS$.
We obtain
$\Gamma(\Lambda^{\mathtt{Q}}_{\mS})<\psi_{\mI_{N}}(b_{0}+\delta^{\mathtt{Q}}_{\mS}+3N-2)<
\psi_{\mI_{N}}(b_{0}+\delta^{\mathtt{Q}}_{\mT})=\Lambda^{\mathtt{Q}}_{\mT}$
by $\delta^{\mathtt{Q}}_{\mS}+3N-2<\delta^{\mathtt{Q}}_{\mT}$, 
$\{b_{0},\delta^{\mathtt{Q}}_{\mS}\}\subset C_{b_{0}+1}(\psi_{\mI_{N}}(b_{0}+1))$.
\\
\textbf{Case 1}. First consider the case when the last inference is a $(\mathtt{D})$:  
There exists a stable ordinal
$\mS\in\mathcal{H}_{\gamma}[\Theta]$ such that $\mS\not\in \mathtt{D}$,  
\[
\infer[(\mathtt{D})]{(\calh_{b_{0}},\Tht;\mathtt{D})\vdash^{* a}_{\Lambda_{0}}\Gam}
{
(\mathcal{H}_{b_{0}},\Tht; \mathtt{D}\cup\{\mS\})\vdash^{* a_{0}}_{\Lambda_{0}}
\Gamma
}
\]
for an $a_{0}<a$. We have $\mS<\Lambda_{0}$, cf.\,\textbf{Case 6} in the proof of Collapsing \ref{lem:Kcollpase.1}.

IH followed by a $(\mathrm{dom})$ yields the lemma.
\[
\infer[(\mathrm{dom})]{(\calh_{\gam},\Theta_{\mathtt{Q}}, \mathtt{Q})\vdash^{\Lambda_{0}+2a}_{\Lambda_{0}}\Gamma
}
{
\{
(\mathcal{H}_{\gamma},\Theta_{\mathtt{Q}}\cup\{\delta^{\mathtt{P}}_{\mS}\}, \mathtt{P})\vdash^{\Lambda_{0}+2a_{0}}_{\Lambda_{0}}
\Gamma(\mathtt{P})
\}_{\mathtt{P}}
}
\]
where $\mathtt{P}$ varies through good extensions $\mathtt{Q}[\mS]$ of $\mathtt{Q}$ for $\gamma,\Theta$.

In what follows we show the existence of a very good extension $\mathtt{P}$ of $\mathtt{Q}$ for $\gamma,\Theta$ by ordinals $\gamma^{\mathtt{P}}_{\mS}$, $\delta^{\mathtt{P}}_{\mS}$ and
$(\rho_{\mS},g_{0})$.
Let $b$ be an ordinal such that $\mT=\psi_{\mI_{N}}(b)$ if $\mS=\mT^{\dagger^{\vec{i}}}$ for 
a sequence $\vec{i}=(i_{0}\geq i_{1}\geq\cdots\geq i_{n})$ with $\mT\in LSt_{k}\cap\Psi$ and 
$k\geq i_{0}\geq i_{n}=i$.
Otherwise let $b=0$.
We see $b\leq b_{0}<\gamma$ as follows, cf.\,(\ref{eq:notationsystem.55}). 
Let $\mS=(\psi_{\mI_{N}}(b))^{\dagger \vec{i}}$.
We obtain $\psi_{\mI_{N}}(b)\in\mathcal{H}_{b_{0}}[\Theta]$ by 
$\mS=(\psi_{\mI_{N}}(b))^{\dagger \vec{i}}\in\mathcal{H}_{b_{0}}[\Theta]$.
On the other hand we have $\Theta\subset C_{b_{0}}(\psi_{\mI_{N}}(b_{0}))$.
Hence $\psi_{\mI_{N}}(b)\in C_{b_{0}+1}(\psi_{\mI_{N}}(b_{0}+1))\cap\mI_{N}=\psi_{\mI_{N}}(b_{0}+1)$, and
$b<b_{0}+1$.

Let $\mathtt{P}=\mathtt{ Q}\cup\{(\mS,(\rho_{\mS},g_{0}))\}$ with $\mathtt{P}(\mS)=\{(\rho_{\mS},g_{0})\}$.
Let $\gam^{\mathtt{P}}_{\mT}=\gam^{\mathtt{Q}}_{\mT}$ and
$\delta^{\mathtt{P}}_{\mT}=\delta^{\mathtt{Q}}_{\mT}$ for $\mT\in \mathtt{D}$.
Ordinals $\delta^{\mathtt{P}}_{\mS}$, $\gamma^{\mathtt{P}}_{\mS}$ and $\rho_{\mS}$ are defined
as follows.

First let
\[
\delta^{\mathtt{P}}_{\mS} = 
\max(\{0\}\cup \{\delta^{\mathtt{Q}}_{\mU}:\mS<\mU\in\mathtt{D}\})+\omega^{1+a_{0}}.
\]
By the assumption we have 
$\delta^{\mathtt{P}}_{\mS}\in\mathcal{H}_{\gamma}[\Theta_{\delta,\mS}]\subset M_{\mathtt{Q}}\cap M_{\rho_{\mS}}$,
$\omega^{a}\leq\delta^{\mathtt{Q}}_{\mU}<\omega^{\Lambda_{0}+1}$,
$\{a_{0}, \delta^{\mathtt{Q}}_{\mU}\}\subset C_{b_{0}+1}(\psi_{\mI_{N}}(b_{0}+1))$
and $\delta^{\mathtt{Q}}_{\mU}+\omega^{1+a}\leq\delta^{\mathtt{Q}}_{\mV}$ for $\{\mV<\mU\}\subset\mathtt{D}$.
We obtain $\omega^{1+a_{0}}\leq\delta^{\mathtt{P}}_{\mS}\in 
C_{b_{0}+1}(\psi_{\mI_{N}}(b_{0}+1))\cap\omega^{\Lambda_{0}+1}$ with $a_{0}<a\leq\Lambda_{0}$,
and $\delta^{\mathtt{Q}}_{\mU}+\omega^{1+a_{0}}\leq\delta^{\mathtt{P}}_{\mS}$ for
$\mS<\mU\in\mathtt{D}$.
Let $\mS>\mV\in\mathtt{D}$.
If there is a $\mS<\mU\in\mathtt{D}$, then 
$\delta^{\mathtt{P}}_{\mS}+\omega^{1+a_{0}}=\delta^{\mathtt{Q}}_{\mU}+\omega^{1+a_{0}}\cdot 2<
\delta^{\mathtt{Q}}_{\mV}$  for the least such $\mU$ by $\omega^{1+a_{0}}\cdot 2<\omega^{1+a}$.
Otherwise $\delta^{\mathtt{P}}_{\mS}+\omega^{1+a_{0}}=\omega^{1+a_{0}}\cdot 2<\omega^{1+a}\leq
\delta^{\mathtt{Q}}_{\mV}$.

Second let
\beqnarr
\alpha_{\mS} & = & \max(\{\gamma\}\cup k_{\mS}(\Theta)\cup\{\gamma^{\mathtt{Q}}_{\mU}:\mS<\mU\in\mathtt{D}\})
\nonumber
\\
\delta_{\mS} & = & \Lambda_{0}\cdot(1+a_{0})+\max(\{0\}\cup E_{\mS}(\Theta))
\nonumber
\\
\beta_{\mS} & = & \omega^{\delta_{\mS}}
\nonumber
\\
\gamma^{\mathtt{P}}_{\mS} & = & \alpha_{\mS}+\beta_{\mS}
\label{eq:capping_case2_I_1}
\eeqnarr
Let $\rho_{\mS}=\psi_{\mathbb{S}}^{g_{0}}(\gamma^{\mathtt{P}}_{\mS}+1)$.
We obtain $\{\Lambda_{0},\mI_{N}\cdot 3\}\subset C_{0}(SC(\gamma^{\mathtt{P}}_{\mS}))$, cf.\,(\ref{eq:notationsystem.5}) and
$b_{0}\leq\gamma<\gamma^{\mathtt{P}}_{\mS}$, cf.\,(\ref{eq:notationsystem.55}).
Also $\mathrm{fld}(g_{0})\subset\mathcal{H}_{\gamma}[\Theta]$.

We have  $\{\gamma,\Lambda_{0},\mI_{N},a_{0}\}\subset\mathcal{H}_{\gamma}[\Theta]$ and
$\gamma<\gamma^{\mathtt{Q}}_{\mU}<\gamma^{\mathtt{P}}_{\mS}$.
On the other hand we have
$E_{\mS}(\Theta)\subset\mathcal{H}_{\gamma}[\Theta]$ by Proposition \ref{prp:EH_Pi11}.\ref{prp:EH_Pi11.5}, and
$k_{\mS}(\Theta)\in\mathcal{H}_{\gamma^{\mathtt{P}}_{\mS}}[\Theta]$
by
Proposition \ref{prp:Buchholz4.3}.\ref{prp:Buchholz4.3.6}.
We obtain
$\gamma^{\mathtt{P}}_{\mS}\in\mathcal{H}_{\gamma}[\Theta_{\gamma,\mS}]$ and
$\rho_{\mS}\in\mathcal{H}_{\gamma^{\mathtt{P}}_{\mS}+\mS}[\Theta_{\gamma,\mS}]$.

Next we show $\Tht\subset M_{\rho_{\mS}}$.
We obtain $\max (\{0\}\cup E_{\mS}(\Tht))\in C_{\gamma^{\mathtt{P}}_{\mS}}(\rho_{\mS})\cap\mS=\rho_{\mS}$,
and hence
$E_{\mS}(\Tht)\subset\rho_{\mS}$.
We obtain $\Theta\cap\mS\subset\rho_{\mS}$.
On the other hand we have $\mathtt{p}_{0}(\rho_{\mS})=\gamma^{\mathtt{P}}_{\mS}+1>\gamma^{\mathtt{P}}_{\mS}> k_{\mS}(\Theta)$.
Proposition \ref{prp:ESM} yields $\Theta\cap M_{\mathtt{Q}}=\Theta\subset C_{\mathtt{p}_{0}(\rho_{\mS})}(\rho_{\mS})=M_{\rho_{\mS}}$.

Third we show that $\gamma^{\mathtt{P}}_{\cdot}$ is a threshold function.
We have
$\gamma^{\mathtt{Q}}_{\mU}+\mU<\gamma^{\mathtt{Q}}_{\mT}$ for $\{\mT<\mU\}\subset\mathtt{D}$, 
and
$\gamma^{\mathtt{Q}}_{\mU}+\mU<\gamma^{\mathtt{P}}_{\mS}$ for $\mS<\mU\in\mathtt{D}\subset\Lambda_{0}$
by the definition (\ref{eq:capping_case2_I_1}).
Let $\mathtt{D}\ni\mT<\mS$.
We show $\gamma^{\mathtt{P}}_{\mS}+\mS<\gamma^{\mathtt{Q}}_{\mT}$.
We have $\beta_{\mS}+\omega^{\Lambda_{0}(1+a_{0})}<\omega^{\Lambda_{0}(1+a)}$.

If $\alpha_{\mS}=\gamma$, then 
$\gamma+\beta_{\mS}+\mS<\gamma+\omega^{\Lambda_{0}(1+a)}\leq \gamma^{\mathtt{Q}}_{\mT}$.
If $\alpha_{\mS}=\gamma^{\mathtt{Q}}_{\mU}$ for a $\mS<\mU\in\mathtt{D}$, then
$\gamma^{\mathtt{Q}}_{\mU}+\beta_{\mS}+\mS<\gamma^{\mathtt{Q}}_{\mU}+\omega^{\Lambda_{0}(1+a)}\leq\gamma^{\mathtt{Q}}_{\mT}$.

Let $\alpha_{\mS}=k_{\mS}(\Theta)$.
We have $\Theta\subset M_{\rho_{\mT}}=C_{\gamma^{\mathtt{Q}}_{\mT}}(\rho_{\mT})\subset C_{\gamma^{\mathtt{Q}}_{\mT}}(\mT)$.
Hence $k_{\mS}(\Theta)\leq k_{\mT}(\Theta)<\gamma^{\mathtt{Q}}_{\mT}=\alpha_{\mT}+\beta_{\mT}$
with an additive principal $\beta_{\mT}\geq \omega^{\Lambda_{0}(1+a)}$.
Let $d<\beta_{\mT}$ be an ordinal such that $k_{\mS}(\Theta)<\alpha_{\mT}+d$.
We obtain
$k_{\mS}(\Theta)+\beta_{\mS}+\mS<\gamma^{\mathtt{Q}}_{\mT}$
by $d+\beta_{\mS}+\mS<\beta_{\mT}$.
\\
\textbf{Case 2}. 
Second consider the case when the last inference is an $(i{\rm -stbl}(\mS))$: 
We have $\mS\in\mathtt{D}\cap SSt_{i}\cap \Lambda_{0}$,
a $\bigwedge$-formula
$B(\mathsf{L}_{0})\in\mathcal{L}_{i}$
with $\mathrm{rk}(B(\mathsf{L}_{0}))<\mS$, and
a term $u\in Tm(\mathbb{I}_{N})$
such that $\mS\leq \mathrm{rk}(B(u))<\Lambda_{0}$,
$\{\lnot B(u),\exists x\in\mathsf{L}_{\mathbb{S}}B(x)\}\subset\Gamma$, 
and $\mathsf{k}(B(u))\subset\mathcal{H}_{\gamma}[\Theta]$.
We have $\mathrm{dom}(\mathtt{Q})\subset\Lambda_{0}$ and
$\{\lnot B(u)^{(\bar{\rho})},(\exists x\in\mathsf{L}_{\mS}B(x))^{(\bar{\rho})}\}\subset\Gamma^{(\bar{\rho})}$.
Also $\forall\mU\in\mathtt{D}\cap \bigcup_{j\geq i}SSt_{j}(\mS<\mU \Rarw \mathrm{rk}(B(u))<\mU)$, cf.\,(\ref{eq:stbl_I_1}).
We have $(\mathcal{H}_{b_{0}},\Theta;\mathtt{D}
)\vdash^{* a}_{\Lambda_{0}}\Gamma$ for an $a>0$.

Let $d=\mathrm{rk}(B(u))<\Lambda_{0}$.
Let $\mT$ be a successor $i$-stable ordinal 
such that $\mS\leq\mT\leq d$.
Let $\mathtt{P}=\mathtt{Q}[\mT]$ be a good extension of $\mathtt{Q}$ for $\gamma,\Theta$ by ordinals
$\gamma^{\mathtt{P}}_{\mT},\delta^{\mathtt{P}}_{\mT}$ and 
$\bar{\rho}_{\mathtt{P}}(\mT)=\{(\rho_{\mT},g_{0})\}$.
We see that a very good extension exists as in \textbf{Case 1} for $a_{0}=0<a$.

Let $\bar{\tau}=\bar{\rho}_{\mathtt{P}}=\bar{\rho}*\{(\mT,\bar{\rho}_{\mathtt{P}}(\mT))\}$ and
$B(u)\simeq\bigwedge(A_{\iota})_{\iota\in I}$.
We obtain for $\iota\in[\bar{\tau}\mathtt{P}]I$ 
\beqn\label{eq:capping.1_I_1}
(\mathcal{H}_{\gamma},\Theta\cup\mathsf{k}(B(u))\cup\mathsf{k}(A_{\iota}),
\mathtt{P}[\mathrm{srk}(A_{\iota})])\vdash^{2d_{\iota}+1}_{0}  
(\{\lnot B(u)^{(\bar{\tau})}, A_{\iota}^{(\bar{\tau})}\})(\mathtt{P}[\mathrm{srk}(A_{\iota})])
\eeqn
by Tautology \ref{lem:tautology.cap}.\ref{lem:tautology.cap1}, where $d_{\iota}=\mathrm{rk}(A_{\iota})$.

Let $f$ be a special finite function with base $\mI_{N}$
such that $\mathrm{supp}(f)=\{\Lambda_{0}\}$ and
$f(\Lambda_{0})=\mI_{N}$.
We obtain $\min\{s(g_{0}), s(f)\}=\Lambda_{0}\geq\mT+1$,
$f_{\Lambda_{0}}=(g_{0})_{\Lambda_{0}}=\emptyset$ and $f<_{\Lambda_{N}}^{\Lambda_{0}}(g_{0})^{\prime}(\Lambda_{0})$
by $f(\Lambda_{0})=\mI_{N}<\mI_{N}\cdot 2=(g_{0})^{\prime}(\Lambda_{0})$.
Let $\sigma\in H_{\rho_{\mT}}(f,\mathtt{P},\Theta)$ with
$SC_{\mI_{N}}(\mathrm{fld}(f))\subset\Lambda_{1}$.

We obtain
$
(\mathcal{H}_{\gamma},\Theta\cup\mathsf{k}(B(u)),\mathtt{P}^{\sigma})\vdash^{2d}_{0}\lnot B(u)^{(\bar{\tau}_{\sigma})}, 
(B(u)^{[\sigma/\mT]})^{(\bar{\tau})}
$  by Tautology \ref{lem:tautology.cap}.\ref{lem:tautology.cap3},
where 
$B(u^{[\sigma/\mathbb{T}]})\equiv B(u)^{[\sigma/\mT]}$.
We see
$\mathrm{inv}((B(u^{[\sigma/\mathbb{T}]}))^{(\bar{\tau})};\mathtt{P}^{\sigma})=B(u)^{(\bar{\tau})}$
for $u^{[\sigma/\mathbb{T}]}\in Tm(\mathbb{T})$.
$\mathsf{k}((B(u^{[\sigma/\mathbb{T}]}))^{(\bar{\tau})};\mathtt{P}^{\sigma})=\mathsf{k}(B(u))\subset
\mathcal{H}_{\gamma}[(\Theta\cup\mathsf{k}(B(u)))\cap M_{\mathtt{Q}^{\sigma}}]\cap M_{\mathtt{P}^{\sigma}}$,
and $u^{[\sigma/\mathbb{T}]}\in[\bar{\tau}\mathtt{P}^{\sigma}]J$ follow.
A $(\bigvee)$ yields 
$(\mathcal{H}_{\gamma},\Theta\cup\mathsf{k}(B(u)),\mathtt{P}^{\sigma})\vdash^{2d+1}_{0}\lnot B(u)^{(\bar{\tau}_{\sigma})}, 
((B(u)^{[\sigma/\mT]})^{\lor})^{(\bar{\tau})}
$.

On the other hand we have
$(\mathcal{H}_{\gamma},\Theta\cup\mathsf{k}(B(u)),\mathtt{P}^{\sigma})\vdash^{0}_{0}\lnot B(u)^{(\bar{\tau}_{\sigma})}, 
(u^{[\sigma/\mT]}\dot{\in}\mathsf{L}_{\mT})^{(\bar{\tau})}$.
We obtain
$(\mathcal{H}_{\gamma},\Theta\cup\mathsf{k}(B(u)),\mathtt{P}^{\sigma})\vdash^{1}_{0}\lnot B(u)^{(\bar{\tau}_{\sigma})}, 
((u^{[\sigma/\mT]}\dot{\in}\mathsf{L}_{\mT})^{\lor})^{(\bar{\tau})}$ by a $(\bigvee)$
with $(u^{[\sigma/\mT]}\dot{\in}\mathsf{L}_{\mT})\equiv\bigwedge\emptyset$.
A $(\bigvee)$ yields
\beqn\label{eq:capping.2_I_1}
(\mathcal{H}_{\gamma},\Theta\cup\mathsf{k}(B(u)),\mathtt{P}^{\sigma})\vdash^{2d+2}_{\Lambda_{0}}\lnot B(u)^{(\bar{\tau}_{\sigma})},
(\exists x\in\mathsf{L}_{\mathbb{T}}B(x))^{(\bar{\tau})}
\eeqn
where $\Xi(\mathtt{P})=\{(\exists x\in\mathsf{L}_{\mathbb{T}}B(x))^{(\bar{\tau})}\}$
for $\Xi=\{(\exists x\in\mathsf{L}_{\mS} B(x))^{(\bar{\rho})}\}\subset\Gamma^{(\bar{\rho})}$.

For $\bar{\rho}=\bar{\rho}_{\mathtt{Q}}$,
an inference $(i{\rm -rfl}_{\mS}(g_{0},\Lambda_{0},f,\Delta))$ with $\Delta=\{B(u)\}$,
$\mathrm{rk}(B(u))<\Lambda_{0}\in\mathrm{supp}(g_{0})$, (\ref{eq:capping.1_I_1}) and (\ref{eq:capping.2_I_1})
yields
$
(\mathcal{H}_{\gamma},\Theta\cup\mathsf{k}(B(u)), \mathtt{ Q})
\vdash^{2d+2}_{\Lambda_{0}}
\lnot B(u)^{(\bar{\rho})}, (\exists x\in\mathsf{L}_{\mathbb{S}}B(x))^{(\bar{\rho})}
$, where $2d+2<\Lambda_{0}$ and $\{\lnot B(u)^{(\bar{\rho})},(\exists x\in\mathsf{L}_{\mS}B(x))^{(\bar{\rho})}\}\subset\Gamma^{(\bar{\rho})}$.
\\
\textbf{Case 3}.
Third the last inference introduces a bounded $\bigwedge$-formula $A$:
Let
$A\simeq\bigwedge\left(A_{\iota}\right)_{\iota\in J}$.
For every
$\iota\in J$,
$(\mathcal{H}_{b_{0}},\Theta\cup\sfk(A_{\iota}) ;\mathtt{D}\cup\{\mathrm{srk}(A_{\iota})\}
)
\vdash^{* a(\iota)}_{\Lambda_{0}}\Gamma, A_{\iota}$ holds for an $a(\iota)<a$.

Let  
$\iota\in [\bar{\rho}\mathtt{Q}]J$. Then 
$\sfk(A_{\iota})\subset C_{b_{0}+1}(\psi_{\mI_{N}}(b_{0}+1))\cap M_{\mathtt{Q}}$.
Hence the assumption is enjoyed for $\Theta_{\iota}$.
We define an extension $\mathtt{Q}_{\iota}=\mathtt{Q}[\mathrm{srk}(A_{\iota})]$ of $\mathtt{Q}$ by 
the ordinal $\mathrm{srk}(A_{\iota})$.

Ordinals $\delta^{\mathtt{Q}_{\iota}}_{\mS}$,
$\gamma^{\mathtt{Q}_{\iota}}_{\mS}$ and $\rho_{\mS}=\psi_{\mS}^{g_{0}}(\gamma^{\mathtt{Q}_{\iota}}_{\mS}+1)$
are defined for $\mS\in \{\mathrm{srk}(A_{\iota})\}\setminus\mathrm{dom}(\mathtt{Q})$
as follows.
Let
$\delta^{\mathtt{Q}_{\iota}}_{\mS}=\max(\{0\}\cup\{\delta^{\mathtt{Q}_{\iota}}_{\mU}: \mS<\mU\in\mathtt{D}\})+\omega^{1+a(\iota)}$ and
$\alpha_{\mS} = \max(\{\gamma\}\cup k_{\mS}(\Theta)\cup\{\gamma^{\mathtt{Q}_{\iota}}_{\mU}:\mS<\mU\in\mathtt{D}\})$.
Then let $\delta_{\mS}= \Lambda_{0}(1+a_{0})+ \max(\{0\}\cup E_{\mS}(\Theta))$, $\beta_{\mS} = \omega^{\delta_{\mS}}$ and
$\gamma^{\mathtt{Q}_{\iota}}_{\mS} = \alpha_{\mS}+\beta_{\mS}$.

We see that 
$\delta^{\mathtt{Q}_{\iota}}_{\mS}\in C_{b_{0}+1}(\psi_{\mI_{N}}(b_{0}+1))\cap\omega^{\Lambda_{0}+1}$,
$\gamma^{\mathtt{Q}_{\iota}}_{\mS}\in\mathcal{H}_{\gamma^{\mathtt{Q}_{\iota}}_{\mS}}[\Theta]$ and
$\Theta\subset M_{\rho_{\mS}}$ as in \textbf{Case 1}.
Moreover we have $\omega^{1+a(\iota)}\leq \delta^{\mathtt{Q}_{\iota}}_{\mS}$
and $\delta^{\mathtt{Q}_{\iota}}_{\mS}+\omega^{1+a(\iota)}\leq \delta^{\mathtt{Q}}_{\mT}$
for $\mS>\mT\in\mathrm{dom}(\mathtt{Q})$.
Hence $\gamma^{\mathtt{Q}_{\iota}}_{\cdot}$ with $\gamma$ has gaps $\omega^{\Lambda_{0}(1+a_{0})}$.

IH yields 
$(\mathcal{H}_{\gamma},\Theta_{\mathtt{Q}_{\iota}}\cup\sfk(A_{\iota}),\mathtt{Q}_{\iota}
)
\vdash^{\Lambda_{0}+2a(\iota)}_{\Lambda_{0}}
(\Gamma^{(\bar{\rho})}\cup\{\left(A_{\iota}\right)^{(\bar{\rho})}\})(\mathtt{Q}_{\iota})$.
$(\mathcal{H}_{\gamma},\Theta, \mathtt{ Q})\vdash^{\Lambda_{0}+2a}_{\Lambda_{0}}
\Gamma^{(\bar{\rho})}$ follows by a $(\bigwedge)$.
\\
\textbf{Case 4}.
Fourth the last inference introduces a $\bigvee$-formula $A\in\Gamma$:
Let $A\simeq\bigvee(A_{\iota})_{\iota\in J}$.
$(\mathcal{H}_{b_{0}},\Theta ;\mathtt{D}
)
\vdash^{* a(\iota)}_{\Lambda_{0}}\Gamma, A_{\iota}$ holds for an $a(\iota)<a$ and an $\iota\in J$.
We have $\mathsf{k}(A_{\iota})\subset\mathcal{H}_{\gamma}[\Theta(\mathtt{Q})]$ by (\ref{eq:controlder*1}),
$\mathrm{srk}(A_{\iota})\in\mathtt{D}=\mathrm{dom}(\mathtt{Q})$ by (\ref{eq:controlder*_cover}), and 
$\Theta\subset C_{b_{0}+1}(\psi_{\mI_{N}}(b_{0}+1))\cap M_{\mathtt{Q}}$.
Let $A_{\iota}\simeq\bigwedge(B_{\nu})_{\nu\in I}$ and $\nu\in[\bar{\rho}\mathtt{Q}]I$.
Let $\mathtt{Q}_{\nu}=\mathtt{Q}[\mathrm{srk}(B_{\nu})]$ be an extension of $\mathtt{Q}$ by $\mathrm{srk}(B_{\nu})$
defined as in \textbf{Case 3}.
By IH  we obtain
$(\mathcal{H}_{\gamma},\Theta_{\mathtt{Q}}, \mathtt{Q}
)
\vdash^{\Lambda_{0}+2a(\iota)}_{\Lambda_{0}}\Gamma^{(\bar{\rho})}, A_{\iota}^{(\bar{\rho})}$, and
$(\mathcal{H}_{\gamma},\Theta_{\mathtt{Q}_{\nu}}, \mathtt{Q}_{\nu}
)
\vdash^{\Lambda_{0}+2a(\iota)}_{\Lambda_{0}}
(\Gamma^{(\bar{\rho})}\cup\{A_{\iota}^{(\bar{\rho})}\})(\mathtt{Q}_{\nu})$.
We obtain
$(\mathcal{H}_{\gamma},\Theta_{\mathtt{Q}_{\nu}}\cup\mathsf{k}(B_{\nu}), \mathtt{Q}_{\nu}
)
\vdash^{2d_{\nu}+1}_{\Lambda_{0}}(\{\lnot A_{\iota}^{(\bar{\rho})}, B_{\nu}^{(\bar{\rho})}\})(\mathtt{Q}_{\nu})$
by Tautology \ref{lem:tautology.cap}.\ref{lem:tautology.cap1}, where $d_{\nu}=\mathrm{rk}(B_{\nu})<\mathrm{rk}(A_{\iota})<\mathrm{rk}(A)<\Lambda_{0}$.
A $(cut)$ yields
$(\mathcal{H}_{\gamma},\Theta_{\mathtt{Q}_{\nu}}\cup\mathsf{k}(B_{\nu}), \mathtt{Q}_{\nu}
)
\vdash^{\Lambda_{0}+2a(\iota)+1}_{\Lambda_{0}}
(\Gamma^{(\bar{\rho})}\cup\{B_{\nu}^{(\bar{\rho})}\})(\mathtt{Q}_{\nu})$.
We obtain
$(\mathcal{H}_{\gamma},\Theta, \mathtt{Q})\vdash^{\Lambda_{0}+2a}_{\Lambda_{0}}\Gamma^{(\bar{\rho})}$ by a $(\bigvee)$.

The lemma follows from IH when the last inference is 
one of $(cut)$ and  $(\Sigma(\Ome)\mbox{{\rm -rfl}})$.
Each cut formula puts on the cap $(\bar{\rho})$.
\eprf

\blem\label{lem:prereduction_I_1}
Let $(\mathcal{H}_{\gamma},\Theta, \mathtt{Q})\vdash^{a}_{c,\vec{d},\xi,\Lambda}\Gamma$ for a stable ordinal $\xi$.

\benu
\item\label{lem:prereduction_I_1.1}
Let $\mT$ be a stable ordinal $\mT\not\in \mathrm{dom}(\mathtt{Q})$ such that $\mT\in\mathcal{H}_{\gamma}[\Theta(\mathtt{Q})]\cap\xi$,
and $\mathtt{P}=\mathtt{Q}[\mT]$ a good extension of $\mathtt{Q}$ for $\gamma,\Theta$ by
 ordinals $\gamma^{\mathtt{P}}_{\mT}, \delta^{\mathtt{P}}_{\mT}$ and $\bar{\rho}_{\mathtt{P}}(\mT)$.
Then
\begin{equation}\label{eq:prereduction_I_1}
(\mathcal{H}_{\gamma},\Theta_{1}, \mathtt{P})\vdash^{a}_{c,\vec{d},\xi,\Lambda}\Gamma_{1}(\mathtt{P})
\end{equation}
holds for $\Theta_{1}=\Theta_{\mathtt{P}}$ and $\Gamma_{1}=\Gamma$.

\item\label{lem:prereduction_I_1.2}
Let $A^{(\bar{\rho})}$ be a capped bounded formula over $\mathtt{Q}$ such that 
$\mathrm{srk}(A)\in\mathcal{H}_{\gamma}[\Theta(\mathtt{Q})]\cap\xi$.
let $\mathtt{P}=\mathtt{Q}[\mathrm{srk}(A)]$ be a good extension of $\mathtt{Q}$ for $\gamma,\Theta$.
Then (\ref{eq:prereduction_I_1}) holds for $\Theta_{1}=\Theta_{\mathtt{P}}\cup\mathsf{k}(A^{(\bar{\rho})};\mathtt{Q})$ and $\Gamma_{1}=\Gamma\cup\{A^{(\bar{\rho})}\}$.

\eenu
\elem
\bprf
Each is shown by induction on $a$. 
We obtain $\mathrm{rk}(A)<\xi$ by $\mathrm{srk}(A)<\xi$ for bounded formulas $A$ by
Proposition \ref{prp:pd_closed_cover_I_1}.\ref{prp:pd_closed_cover_I_1.4}.
We may assume 
$\mT=\mathrm{srk}(A)\not\in\mathrm{dom}(\mathtt{Q})$ 
in Lemma \ref{lem:prereduction_I_1}.\ref{lem:prereduction_I_1.2}.
Let $\bar{\rho}(\mathtt{P})=\bar{\rho}*\bar{\kappa}$ for $\bar{\kappa}=((\rho_{\mathtt{P}(\mT)},g_{\mT}))$.

We have $\mathrm{inv}^{(\mathtt{u})}(B^{(\bar{\rho})};\mathtt{Q})\equiv \mathrm{inv}^{(\mathtt{u})}(B^{(\bar{\rho}(\mathtt{P}))};\mathtt{P})$.
By $\Theta(\mathtt{Q})\subset M_{\mathtt{P}}$, we obtain $\mathcal{H}_{\gamma}[\Theta(\mathtt{Q})]=\mathcal{H}_{\gamma}[\Theta(\mathtt{P})]$.
Hence each of
(\ref{eq:controlder_cap2}) and
 (\ref{eq:controlder_cap22}) is enjoyed  in (\ref{eq:prereduction_I_1}).
\\
\textbf{Case 1}. The last inference is a $(\mathrm{dom})$:
We have an ordinal $a_{0}<a$, and a stable ordinal $\mU\in\mathcal{H}_{\gamma}[\Theta(\mathtt{Q})]\cap\xi$
such that $\mU\not\in\mathrm{dom}(\mathtt{Q})$.
For each good extension $\mathtt{R}=\mathtt{Q}[\mU]$ of $\mathtt{Q}$ for $\gamma,\Theta$ 
by ordinals $\gamma^{\mathtt{R}}_{\mU}, \delta^{\mathtt{R}}_{\mU},\bar{\rho}_{\mathtt{R}}(\mU)$, we have
$(\mathcal{H}_{\gamma},\Theta_{\mathtt{R}},\mathtt{R})\vdash^{a_{0}}_{c,\vec{d},\xi,\Lambda}\Gamma(\mathtt{R})$.

First let $\mU\in\mathrm{dom}(\mathtt{P})$. This means $\mU=\mT$.
Let $\gamma^{\mathtt{R}}_{\mU}=\gamma^{\mathtt{P}}_{\mU}$, $\delta^{\mathtt{R}}_{\mU}=\delta^{\mathtt{P}}_{\mU}$ and $\bar{\rho}_{\mathtt{R}}(\mU)=\bar{\rho}_{\mathtt{P}}(\mU)$.
Then
$\mathtt{R}=\mathtt{P}_{0}=\mathtt{Q}[\mU]$.
Hence $(\mathcal{H}_{\gamma},\Theta_{\mathtt{P}_{0}},\mathtt{P}_{0})\vdash^{a_{0}}_{c,\vec{d},\xi,\Lambda}\Gamma(\mathtt{P}_{0})$.
Consider Lemma \ref{lem:prereduction_I_1}.\ref{lem:prereduction_I_1.2}.
We see that $\mathtt{P}=\mathtt{Q}[\mathrm{srk}(A)]$ is a good extension of $\mathtt{Q}[\mU]$ for $\gamma,\Theta$.
IH yields $(\mathcal{H}_{\gamma},\Theta_{1},\mathtt{P})\vdash^{a_{0}}_{c,\vec{d},\xi,\Lambda}\Gamma_{1}(\mathtt{P})$.

Next let $\mU\not\in\mathrm{dom}(\mathtt{P})$, i.e., $\mU\neq\mT$.
The extension $\mathtt{R}[\mT]$ is a good extension of $\mathtt{R}$ for $\gamma,\Theta$.
By IH we obtain
$(\mathcal{H}_{\gamma},(\Theta_{1})_{\mathtt{R}[\mT]},\mathtt{R}[\mT])\vdash^{a_{0}}_{c,\vec{d},\xi,\Lambda}\Gamma_{1}(\mathtt{R}[\mT])$ for each extension
$\mathtt{R}$ by the ordinals $\gamma^{\mathtt{R}}_{\mU}, \delta^{\mathtt{R}}_{\mU}, \bar{\rho}_{\mathtt{R}}(\mU)$.
On the other hand we have $\mathtt{P}[\mU]=\mathtt{R}[\mT]$.
Specifically the extension $\mathtt{P}[\mU]$ of $\mathtt{P}$ by the ordinals $\gamma^{\mathtt{R}}_{\mU}, \delta^{\mathtt{R}}_{\mU}, \bar{\rho}_{\mathtt{R}}(\mU)$
is equal to one $\mathtt{R}[\mT]$ of $\mathtt{R}$ by the ordinals $\gamma^{\mathtt{P}}_{\mT}, \delta^{\mathtt{P}}_{\mT}, \bar{\rho}_{\mathtt{P}}(\mU)$.
We obtain $(\mathcal{H}_{\gamma},(\Theta_{1})_{\mathtt{P}[\mU]},\mathtt{P}[\mU])\vdash^{a_{0}}_{c,\vec{d},\xi,\Lambda}\Gamma_{1}(\mathtt{P}[\mU])$, 
where 
if $\mathtt{R}=\mathtt{Q}[\mU]$ is a very good extension of $\mathtt{Q}$, then so is $\mathtt{P}[\mU]$.
A $(\mathrm{dom})$ yields (\ref{eq:prereduction_I_1}).
\\
\textbf{Case 2}. The last inference is a $(\bigwedge)$ introducing a $B^{(\bar{\rho})}\in\Gamma$: 
We see that it suffices to prune some branches as follows.
Let $B\simeq\bigwedge(B_{\iota})_{\iota\in J}$.
For each $\iota\in[\bar{\rho}\mathtt{Q}]J$ we have an
ordinal $a(\iota)<a$ such that
$(\mathcal{H}_{\gamma},\Theta_{\mathtt{Q}_{\iota}}\cup\mathsf{k}(B_{\iota}^{(\bar{\rho})};\mathtt{Q}) , \mathtt{Q}_{\iota}
)
\vdash^{a(\iota)}_{c, \vec{d},\xi,\Lambda}(\Gamma\cup\{B_{\iota}^{(\bar{\rho})}\})(\mathtt{Q}_{\iota})$ holds for each good extension $\mathtt{Q}_{\iota}=\mathtt{Q}[\mathrm{srk}(B_{\iota})]$ of $\mathtt{Q}$
for $\gamma,\Theta$, where 
$\mathsf{k}(B_{\iota}^{(\bar{\rho})};\mathtt{Q})=
 \mathsf{k}(B_{\iota}^{(\bar{\rho}(\mathtt{P}))};\mathtt{P})$.

Let $\iota\in[\bar{\rho}(\mathtt{P})\mathtt{P}]J$.
We obtain $M_{\mathtt{P}}\subset M_{\mathtt{Q}}$, $\mathsf{k}(B_{\iota})\subset M_{\mathtt{P}}$ and $\iota\in[\bar{\rho}\mathtt{Q}]J$.
Let 
$\mathtt{P}_{\iota}=\mathtt{P}[\mathrm{srk}(B_{\iota})]$ be a good extension of $\mathtt{P}$ for $\gamma,\Theta_{1}$.
Let $\mathtt{Q}_{\iota}=\mathtt{Q}[\mathrm{srk}(B_{\iota})]$ denote an extension of $\mathtt{Q}$ defined by
$\gamma^{\mathtt{Q}_{\iota}}_{\mU}=\gamma^{\mathtt{P}_{\iota}}_{\mU}$ and
$\delta^{\mathtt{Q}_{\iota}}_{\mU}=\delta^{\mathtt{P}_{\iota}}_{\mU}$ for $\mU\in \mathrm{srk}(B_{\iota})$.
Then $\mathtt{Q}_{\iota}$ is a good extension of $\mathtt{Q}$ for $\gamma,\Theta$.

IH yields 
$(\mathcal{H}_{\gamma},(\Theta_{1})_{\mathtt{Q}_{\iota}[\mT]}\cup \mathsf{k}(B_{\iota}^{(\bar{\rho})};\mathtt{Q}), \mathtt{Q}_{\iota}[\mT]
)
\vdash^{a(\iota)}_{c,\vec{d},\xi,\Lambda}
(\Gamma_{1}\cup\{B_{\iota}^{(\bar{\rho})}\})(\mathtt{Q}_{\iota}[\mT])$, where
$\mathtt{P}_{\iota}=\mathtt{P}[\mathrm{srk}(B_{\iota})]=\mathtt{Q}_{\iota}[\mT]$.
We obtain
$(\mathcal{H}_{\gamma},(\Theta_{1})_{\mathtt{P}_{\iota}}\cup
\mathsf{k}(B_{\iota}^{(\bar{\rho}(\mathtt{P}))};\mathtt{P}), \mathtt{P}_{\iota}
)
\vdash^{a(\iota)}_{c,\vec{d},\xi,\Lambda}
(\Gamma_{1}\cup\{B_{\iota^{(\bar{\rho})}}\})(\mathtt{P}_{\iota})$ for each $\iota\in[\bar{\rho}(\mathtt{P})\mathtt{P}]J$.
A $(\bigwedge)$ yields (\ref{eq:prereduction_I_1}).
\\
\textbf{Case 3}. The last inference is an $(i {\rm -rfl}_{\mU}(g,x,f,\Delta))$:
We have
$\xi>\mU\in \mathrm{dom}(\mathtt{Q})\cap SSt_{i}$ such that 
$\bar{\rho}_{\mathtt{Q}}(\mU)=(\rho_{\mathtt{Q}(\mU)},g)$,
a special finite function $f$, ordinals $x\in {\rm supp}(g)$, 
$a_{0}\in\mathcal{H}_{\gamma}[\Theta(\mathtt{Q})]\cap a$ and a set $\Delta$ of formulas.

If $\mT\in SSt_{k}$ with a $k\geq i$ and $\mU<\mT\leq\mathrm{rk}(\Delta)$, then
let $\mW=\mT$. Otherwise let $\mW=\mU$ and $k=i$.
We obtain $\mW<\mS\in \mathrm{dom}(\mathtt{Q})\cap \bigcup_{j\geq k}SSt_{j} \Rarw \mathrm{rk}(\Delta)<\mS$ for (r1).

Let $\mV$ be a successor $j$-stable ordinal 
such that $j\geq k$, $\mW\leq\mV\leq\mathrm{rk}(\Delta)$ and $\min\{s(g), s(f)\}\geq\mV+1$. Let
$\mathtt{R}=\mathtt{Q}[\mV]$ be a good extension of $\mathtt{Q}$ for $\gamma,\Theta$ by ordinals
$\gamma^{\mathtt{R}}_{\mV},\delta^{\mathtt{R}}_{\mV}$ and 
$\bar{\rho}_{\mathtt{R}}(\mV)=(\rho_{\mathtt{R}(\mV)},g)$.
Then 
 $(\mathcal{H}_{\gamma},\Theta_{\mathtt{R}_{\iota}}\cup\mathsf{k}(\delta_{\iota}^{(\bar{\rho}_{\mathtt{R}})};\mathtt{R}), \mathtt{R}_{\iota}
 )\vdash^{a_{0}}_{c,\vec{d},\xi,\Lambda}
 (\Gamma\cup\{\lnot\delta_{\iota}^{(\bar{\rho}_{\mathtt{Q}})}\})(\mathtt{R}_{\iota})$ holds for each 
 $\delta\in\Delta$ and $\iota\in[\bar{\rho}_{\mathtt{R}}\mathtt{R}]J_{\delta}$ with $\delta\simeq\bigvee(\delta_{\iota})_{\iota\in J_{\delta}}$,
 and each good extension $\mathtt{R}_{\iota}=\mathtt{R}[\mathrm{srk}(\delta_{\iota})]$.
On the other hand we have
$
(\mathcal{H}_{\gamma},\Theta_{\mathtt{R}},\mathtt{R}^{\sigma}
)\vdash^{a_{0}}_{c,\vec{d},\xi,\Lambda}\Delta^{(\bar{\rho}_{\mathtt{R}^{\sigma}})},\Gamma(\mathtt{R})$ for each $\sigma\in H_{\rho_{\mathtt{R}}(\mT)}(f,\mathtt{R},\Theta)$, where $\mathrm{rk}(\Delta)<d_{i}\leq d_{k}$
for the $i$-th element $d_{i}$ in $\vec{d}$.

Given a good extension $\mathtt{S}=\mathtt{P}[\mV]$ of $\mathtt{P}=\mathtt{Q}[\mT]$, let
$\mathtt{R}=\mathtt{Q}[\mV]$ be a good extension such that
$\mathtt{R}[\mT]=\mathtt{S}=\mathtt{P}[\mV]$.
Let $\delta\in\Delta$, $\iota\in[\bar{\rho}(\mathtt{S})\mathtt{S}]J_{\delta}$, and $\mathtt{S}_{\iota}=\mathtt{S}[\mathrm{srk}(\delta_{\iota}))$
 be a good extension of $\mathtt{S}$.
Let $\mathtt{R}_{\iota}$ be a good extension of $\mathtt{R}$ such that $\mathtt{R}_{\iota}[\mT]=\mathtt{S}_{\iota}$.
We obtain
 $(\mathcal{H}_{\gamma},(\Theta_{1})_{\mathtt{S}_{\iota}}\cup \mathsf{k}(\delta_{\iota}^{(\bar{\rho}_{\mathtt{R}})};\mathtt{R}), \mathtt{S}_{\iota}
 )\vdash^{a_{0}}_{c,\vec{d},\xi,\Lambda} (\Gamma_{1}\cup\{ \lnot\delta_{\iota}^{(\bar{\rho}_{\mathtt{R}})}\})(\mathtt{S}_{\iota})$ by IH.
 On the other hand we have $H_{\rho_{\mathtt{R}}(\mT)}(f,\mathtt{R},\Theta)=H_{\rho_{\mathtt{S}}(\mT)}(f,\mathtt{P},\Theta)$ by 
 $\Theta\cap M_{\mathtt{Q}}=\Theta\cap M_{\mathtt{R}}=\Theta\cap M_{\mathtt{P}}=\Theta\cap M_{\mathtt{S}}$.
 We obtain
 $
(\mathcal{H}_{\gamma},(\Theta_{1})_{\mathtt{S}},\mathtt{S}^{\sigma}
)\vdash^{a_{0}}_{c,\vec{d},\xi,\Lambda}\Delta^{(\bar{\rho}_{\mathtt{S}^{\sigma}})},\Gamma(\mathtt{S})$ 
for each $\sigma\in H_{\tau}(f,\mathtt{S},\Theta)$ by IH.
A $(k {\rm -rfl}_{\mW}(g,x,f,\Delta)$ yields 
(\ref{eq:prereduction_I_1}).

Other cases are seen from IH.
\eprf

\blem\label{lem:preprereduction}
Let $(\mathcal{H}_{\gamma},\Theta, \mathtt{Q})\vdash^{a}_{c,\vec{d},\xi,\Lambda}\Gamma$
and $\mS\in\mathrm{dom}(\mathtt{Q})$.
Let $(\sigma,f)\in\Psi_{\mS}^{\mathtt{g}}$ be a guarded ordinal such that 
$\Theta(\mathtt{Q})\subset M_{\sigma}$, 
$\mathrm{fld}(f)\subset\mathcal{H}_{\gamma}[\Theta(\mathtt{Q})$, 
$SC_{\mI_{N}}(\mathrm{fld}(f))\subset\Lambda$ and
$\mathtt{Q}^{\sigma}$ is a finite family with $\mathtt{Q}^{\sigma}(\mS)=\mathtt{Q}(\mS)\cup\{(\sigma,f)\}$.
Then
$(\mathcal{H}_{\gamma},\Theta, \mathtt{Q}^{\sigma})\vdash^{a}_{c,\vec{d},\xi,\Lambda}\Gamma$ holds.
\elem
\bprf
By induction on $a$. We may assume $\sigma\not\in(\mathtt{Q}(\mS))_{0}$.
By the assumptions each of (\ref{eq:controlder_cap_cover}), (\ref{eq:controlder_cap22}), 
(\ref{eq:fld}) and (\ref{eq:controlder_cap0}) holds in
$(\mathcal{H}_{\gamma},\Theta, \mathtt{Q}^{\sigma})\vdash^{a}_{c,\vec{d},\xi,\Lambda}\Gamma$.
Let $A^{(\bar{\rho})}\in\Gamma$.
We show 
$\mathsf{k}(A^{(\bar{\rho})};\mathtt{Q}^{\sigma})\subset\mathcal{H}_{\gamma}[\Theta(\mathtt{Q^{\sigma}})]$
for (\ref{eq:controlder_cap2}).
We have $\Theta(\mathtt{Q})=\Theta(\mathtt{Q}^{\sigma})$ by $\Theta(\mathtt{Q})\subset M_{\sigma}$,
and it suffices to show 
$\mathrm{inv}^{(\mathtt{u})}(A^{(\bar{\rho})};\mathtt{Q}^{\sigma})\equiv
\mathrm{inv}^{(\mathtt{u})}((A^{(\bar{\rho})};\mathtt{Q})$.
Suppose $\mathrm{inv}^{(\mathtt{u})}(A^{(\bar{\rho})};\mathtt{Q}^{\sigma})\not\equiv
\mathrm{inv}^{(\mathtt{u})}(A^{(\bar{\rho})};\mathtt{Q})$.
Then 
$\mathrm{inv}^{(\mathtt{u})}(A^{(\bar{\rho})};\mathtt{Q})\equiv A\equiv B^{[\sigma/\mS]}\not\equiv B\equiv \mathrm{inv}^{(\mathtt{u})}(A^{(\bar{\rho})};\mathtt{Q}^{\sigma})$ 
and $\sigma<(\bar{\rho}(\mS))_{0}$.
There is an ordinal $\alpha$ such that $\sigma\leq\alpha\in\mathsf{k}(B^{[\sigma/\mS]})\subset\mS$.
We would have
$\sigma\leq\alpha\in\mathsf{k}(B^{[\sigma/\mS]})=\mathsf{k}(A^{(\bar{\rho})};\mathtt{Q})\subset\mathcal{H}_{\gamma}[\Theta(\mathtt{Q})]\cap\mS\subset M_{\sigma}\cap\mS=\sigma$.
\eprf

\blem\label{lem:predcereg.Sa}{\rm (Reduction)}
Let 
$\bar{\rho}=\bar{\rho}_{\mathtt{Q}}$,
$\Pi^{(\bar{\rho})}=\{C_{m}^{(\bar{\rho})}: 1\leq m\leq n\}$ and $C_{m}\simeq\bigvee(C_{m,\iota})_{\iota\in J_{m}}$ for each $m\leq n$.

For each $m\leq n$, suppose that there exists an $a_{m}\in\mathcal{H}_{\gamma}[\Theta(\mathtt{Q})]$ such that
$(\mathcal{H}_{\gamma},\Theta_{\iota},\mathtt{Q}_{\iota})\vdash^{a_{m}}_{c,\vec{d},\xi,\Lambda}(\Gamma_{m}\cup\{\lnot C_{m,\iota}^{(\bar{\rho})}\})(\mathtt{Q}_{\iota})$ holds
for each $\iota\in[\bar{\rho}\mathtt{Q}]J_{m}$ and each good extension $\mathtt{Q}_{\iota}=\mathtt{Q}[\mathrm{srk}(C_{m,\iota})]$, 
where $\Theta_{\iota}=\Theta_{\mathtt{Q}_{\iota}}\cup\mathsf{k}(C_{m,\iota}^{(\bar{\rho})};\mathtt{Q})$.
Moreover suppose that
$(\mathcal{H}_{\gamma},\Theta, \mathtt{Q})\vdash^{a_{0}}_{c,\vec{d},\xi,\Lambda}\Pi^{(\bar{\rho})},\Gamma_{0}$.

 Let
$\mathrm{rk}(\Pi)\leq 
c+b<\mI_{N}$ with
$b\in\mathcal{H}_{\gamma}[\Theta(\mathtt{Q})]$, 
and
$\Gamma=\Gamma_{0}\cup\Gamma_{1}\cup\cdots\cup\Gamma_{n}$.

Then
\begin{equation}\label{eq:predcereg.Sa}
(\mathcal{H}_{\gamma},\Theta,\mathtt{Q})
\vdash^{\varphi_{b}(a)}_{c,\vec{d},\xi,\Lambda}\Gamma
\end{equation}
holds for the natural sum $a=a_{0}\# a_{1}\#\cdots\# a_{n}$.
\elem
\bprf
By main induction on $b$ with subsidiary induction on the natural sum $a$. 
By $b\in\mathcal{H}_{\gamma}[\Theta(\mathtt{Q})]$,
(\ref{eq:controlder_cap2}) is enjoyed in (\ref{eq:predcereg.Sa}).
\\
\textbf{Case 0}. $a_{m}=0$ for an $m\leq n$:
Since each of $\lnot C_{m,\iota}$ and $C_{m}$ is a $\bigvee$-formula,
the last inference is a void $(\bigwedge)$ with a major formula in $\Gamma_{m}$.
$(\mathcal{H}_{\gamma},\Theta,\mathtt{Q})\vdash^{0}_{c,\vec{d},\xi,\Lambda}\Gamma$ follows.

In what follows assume $a_{m}>0$ for each $m\leq n$.
\\
\textbf{Case 1}. $(\mathcal{H}_{\gamma},\Theta,\mathtt{Q})\vdash^{a_{0}}_{c,\vec{d},\xi,\Lambda}\Pi^{(\bar{\rho})},\Gamma_{0}$
follows by an inference other than a $(\bigvee)$ introducing a $C_{m}^{(\bar{\rho})}$: 
Let $a^{\prime}=a_{0}^{\prime}\#a_{1}\#\cdots\# a_{n}$.
 \\
\textbf{Case 1.1}. The last inference in $(\mathcal{H}_{\gamma},\Theta, \mathtt{Q})\vdash^{a_{0}}_{c,\vec{d},\xi,\Lambda}\Pi^{(\bar{\rho})},\Gamma_{0}$ is a $(\mathrm{dom})$:
For an $a_{0}^{\prime}<a_{0}$ and a $\mT\not\in \mathrm{dom}(\mathtt{Q})$ with 
$\mT\in\mathcal{H}_{\gamma}[\Theta(\mathtt{Q})]\cap\xi$,
we have 
\[
\infer[(\mathrm{dom})]{
(\mathcal{H}_{\gamma},\Theta,\mathtt{Q})\vdash^{a_{0}}_{c,\vec{d},\xi,\Lambda}\Pi^{(\bar{\rho})},\Gamma_{0}
}
{
(\mathcal{H}_{\gamma},\Theta_{\mathtt{P}},\mathtt{P})\vdash^{a_{0}^{\prime}}_{c,\vec{d},\xi,\Lambda}\Pi^{(\bar{\rho})}(\mathtt{P}),\Gamma_{0}(\mathtt{P})
}
\]
for each good extension $\mathtt{P}=\mathtt{Q}[\mT]$.

Given a good extension $\mathtt{P}$, pick a good extension $\mathtt{Q}_{\iota}=\mathtt{Q}[\mathrm{srk}(C_{m,\iota})]$
such that $\mathtt{P}_{\iota}=\mathtt{Q}_{\iota}[\mT]=\mathtt{P}[\mathrm{srk}(C_{m,\iota})]$.

By Lemma \ref{lem:prereduction_I_1}.\ref{lem:prereduction_I_1.1} we obtain
$(\mathcal{H}_{\gamma},(\Theta_{\iota})_{\mathtt{P}_{\iota}},\mathtt{P}_{\iota})\vdash^{a_{m}}_{c,\vec{d},\xi,\Lambda}(\Gamma_{m}\cup\{\lnot C_{m,\iota}^{(\bar{\rho})}\})(\mathtt{P}_{\iota})$,
and $\mathsf{k}(C_{m,\iota})=\mathsf{k}(C_{m,\iota}^{(\bar{\rho})};\mathtt{Q})=\mathsf{k}(C_{m,\iota}^{(\bar{\rho}(\mathtt{P}))};\mathtt{P})$
for $\iota\in[\bar{\rho}(\mathtt{P})\mathtt{P}]J_{m}$.

SIH yields
$(\mathcal{H}_{\gamma},\Theta_{\mathtt{P}},\mathtt{P})
\vdash^{\varphi_{b}(a^{\prime})}_{c,\vec{d},\xi,\Lambda}\Gamma(\mathtt{P})$.
We obtain (\ref{eq:predcereg.Sa}) by a $(\mathrm{dom})$. 
 \\
\textbf{Case 1.2}. The last inference is a $(\bigwedge)$:
Let $A^{(\bar{\tau})}\in\Gamma_{0}$ be the major formula such that $A\simeq\bigwedge(A_{\nu})_{\nu\in I}$, and
$\Gamma_{0}^{\prime}\cup\{A^{(\bar{\tau})}\}=\Gamma_{0}$.
For each $\nu\in[\bar{\tau}\mathtt{Q}]I$, there is an ordinal $a_{0}(\nu)<a$
we have 
\[
\infer[(\bigwedge)]{
(\mathcal{H}_{\gamma},\Theta,\mathtt{Q})\vdash^{a_{0}}_{c,\vec{d},\xi,\Lambda}\Pi^{(\bar{\rho})},\Gamma_{0}
}
{
(\mathcal{H}_{\gamma},\Theta_{\nu},\mathtt{P}_{\nu})\vdash^{a_{0}(\nu)}_{c,\vec{d},\xi,\Lambda}\Pi^{(\bar{\rho})}(\mathtt{P}_{\nu}),
\Gamma^{\prime}_{0}(\mathtt{P}_{\nu}), A_{\nu}^{(\bar{\tau})}(\mathtt{P}_{\nu})
}
\]
for each good extension $\mathtt{P}_{\nu}=\mathtt{Q}[\mathrm{srk}(A_{\nu})]$, where
$\Theta_{\nu}=\Theta_{\mathtt{P}_{\nu}}\cup\mathsf{k}(A_{\nu}^{(\bar{\tau})};\mathtt{Q})$.

Given a good extension $\mathtt{P}_{\nu}$, pick a good extension $\mathtt{Q}_{\iota}=\mathtt{Q}[\mathrm{srk}(C_{m,\iota})]$
such that $\mathtt{P}_{\iota\nu}=\mathtt{Q}_{\iota}[\mathrm{srk}(A_{\nu})]=\mathtt{P}_{\nu}[\mathrm{srk}(C_{m,\iota})]$.
Let $\Theta_{\iota\nu}=\Theta_{\mathtt{P}_{\iota\nu}}\cup\mathsf{k}(C_{m,\iota}^{(\bar{\rho})};\mathtt{Q})\cup\mathsf{k}(A_{\nu}^{(\bar{\tau})};\mathtt{Q})$.

By Lemma \ref{lem:prereduction_I_1}.\ref{lem:prereduction_I_1.2} we obtain
$(\mathcal{H}_{\gamma},\Theta_{\iota\nu},\mathtt{P}_{\iota\nu})\vdash^{a_{m}}_{c,\vec{d},\xi,\Lambda}(\Gamma_{m}\cup\{\lnot C_{m,\iota}^{(\bar{\rho})}, A_{\nu}^{(\bar{\tau})}\})(\mathtt{P}_{\iota\nu})$,
and $\mathsf{k}(C_{m,\iota})=\mathsf{k}(C_{m,\iota}^{(\bar{\rho})};\mathtt{Q})=\mathsf{k}(C_{m,\iota}^{(\bar{\rho}(\mathtt{P}))};\mathtt{P})$
for $\iota\in[\bar{\rho}(\mathtt{P})\mathtt{P}]J_{m}\subset[\bar{\rho}\mathtt{Q}]J_{m}$.

SIH yields
$(\mathcal{H}_{\gamma},\Theta_{\mathtt{P}_{\nu}},\mathtt{P}_{\nu})
\vdash^{\varphi_{b}(a^{\prime}(\nu))}_{c,\vec{d},\xi,\Lambda}\Gamma(\mathtt{P}_{\nu})$ for $a^{\prime}(\nu)=a_{0}(\nu)\#a_{1}\#\cdots\# a_{n}$.
We obtain (\ref{eq:predcereg.Sa}) by a $(\bigwedge)$.
 \\
\textbf{Case 1.3}. The last inference is an $(i {\rm -rfl}_{\mT}(g,x,f,\Delta))$:
 Then $\rho_{\mathtt{Q}(\mT)}=(\bar{\rho}(\mT))_{0}$ for a 
 $\mT\in \mathrm{dom}(\mathtt{Q})\cap SSt_{i}\subset\xi$.
 We have an ordinal $a_{0}^{\prime}<a_{0}$, 
a set $\Delta$ such that $\mathrm{rk}(\Delta)<d_{i}$.
Let $\mT\leq\mU\leq \mathrm{rk}(\Delta)$ with $j\geq i$, $\mU\in SSt_{j}$ and $\min\{s(g),s(f)\}\geq\mU+1$,
and 
$\mathtt{P}=\mathtt{Q}[\mU]$ be a good extension of $\mathtt{Q}$ for $\gamma,\Theta$.

We have
$(\mathcal{H}_{\gamma},\Theta_{\iota},
\mathtt{P}_{\iota}
)\vdash^{a_{0}^{\prime}}_{c,\vec{d},\xi,\Lambda} (\Pi^{(\bar{\rho})}\cup\Gamma_{0}\cup\{\lnot\delta_{\iota}^{(\bar{\rho})}\})(\mathtt{P}_{\iota})$ for each $\delta\in\Delta$ and
$\iota\in[\bar{\rho}\mathtt{P}]J$ with $\delta\simeq\bigvee(\delta_{\iota})_{\iota\in J}$ and each good extension 
$\mathtt{P}_{\iota}=\mathtt{P}[\mathrm{srk}(\delta_{\iota})]$, where
$\Theta_{\iota}=\Theta_{\mathtt{P}_{\iota}}\cup\mathsf{k}(\delta_{\iota}^{(\bar{\rho})};\mathtt{P})$.
We obtain
$(\mathcal{H}_{\gamma},\Theta_{\iota},
\mathtt{P}_{\iota}
)\vdash^{\varphi_{b}(a^{\prime})}_{c,\vec{d},\xi,\Lambda}
(\Gamma\cup\{\lnot\delta_{\iota}^{(\bar{\rho}_{\mathtt{P}})}\})(\mathtt{P}_{\iota})$ by SIH 
and Lemma \ref{lem:prereduction_I_1}.\ref{lem:prereduction_I_1.1} as in \textbf{Case 1}.

On the other hand we have 
 $
(\mathcal{H}_{\gamma},\Theta_{\mathtt{P}},
\mathtt{P}^{\kappa}
)\vdash^{a_{0}^{\prime}}_{c,\vec{d},\xi,\Lambda}\Delta^{(\bar{\rho}_{\kappa})}, (\Pi^{(\bar{\rho})}\cup\Gamma_{0})(\mathtt{P})$ for each 
$\kappa\in H_{\rho_{\mathtt{P}(\mU)}}(f,\mathtt{P},\Theta)$.
$
(\mathcal{H}_{\gamma},\Theta_{\mathtt{P}},
\mathtt{P}^{\kappa}
)\vdash^{\varphi_{b}(a^{\prime})}_{c,\vec{d},\xi,\Lambda}\Delta^{(\bar{\rho}_{\kappa})}, \Gamma(\mathtt{P})$ follows
by SIH, Lemmas \ref{lem:prereduction_I_1}.\ref{lem:prereduction_I_1.1} and \ref{lem:preprereduction}.
$(\mathcal{H}_{\gamma},\Theta,\mathtt{Q})\vdash^{\varphi_{b}(a^{\prime})+1}_{c,\vec{d},\xi,\Lambda}\Gamma$
follows by an $(i {\rm -rfl}_{\mT}(g,x,f,\Delta))$.
\\
\textbf{Case 2}.
In what follows we assume that  the last inference in $(\mathcal{H}_{\gamma},\Theta, \mathtt{Q})\vdash^{a_{0}}_{c,\vec{d},\xi}\Pi^{(\bar{\rho})},\Gamma_{0}$
 is a $(\bigvee)$ with a major formula $C_{m}^{(\bar{\rho})}$. 
 We may assume that $m=n$.
  \begin{equation}\label{eq:reduction_case23}
 \infer[(\bigvee)]{
 (\mathcal{H}_{\gamma},\Theta, \mathtt{Q})\vdash^{a_{0}}_{c,\vec{d},\xi,\Lambda}\Pi^{(\bar{\rho})},\Gamma_{0}
 }
 {
 (\mathcal{H}_{\gamma},\Theta_{\nu}, \mathtt{Q}_{\nu}
 )\vdash^{a^{\prime}_{0}}_{c,\vec{d},\xi,\Lambda}
 (\Pi^{(\bar{\rho})}\cup\{B_{\nu}^{(\bar{\rho})}\}\cup\Gamma_{0})(\mathtt{Q}_{\nu})
 }
  \end{equation}
for each $\nu\in [\bar{\rho}\mathtt{Q}] J_{n, \iota_{0}}$ and
each good extension $\mathtt{Q}_{\nu}=\mathtt{Q}[\mathrm{srk}(B_{\nu})]$ of $\mathtt{Q}$,
where $a_{0}^{\prime}\in\mathcal{H}_{\gamma}[\Theta(\mathtt{Q})]\cap a_{0}$,
$\iota_{0}\in [\bar{\rho}\mathtt{Q}] J_{n}$, 
$\Theta_{\nu}= \Theta_{\mathtt{Q}_{\nu}}\cup\mathsf{k}(B_{\nu}^{(\bar{\rho})};\mathtt{Q})$,
$\mathsf{k}(C_{n,\iota_{0}})=\mathsf{k}(C_{n,\iota_{0}}^{(\bar{\rho})};\mathtt{Q})\subset \mathcal{H}_{\gamma}[\Theta(\mathtt{Q})]$, $\mathrm{srk}(C_{n,\iota_{0}})\in\mathrm{dom}(\mathtt{Q})$,
 and 
$C_{n,\iota_{0}}\simeq\bigwedge(B_{\nu})_{\nu\in J_{n,\iota_{0}}}$.
We have
$\{\delta^{\mathtt{Q}}_{\mS}: \mS\in\mathrm{dom}(\mathtt{Q})\}\cup
\bigcup\{SC_{\mI_{N}}(\mathrm{fld}(g)): g\in(\mathtt{Q})_{1}\}
\subset\mathcal{H}_{\gamma}[\Theta(\mathtt{Q})]$ by (\ref{eq:controlder_cap22}).
 We obtain $\mathtt{Q}=\mathtt{Q}[\mathrm{srk}(C_{n,\iota_{0}})]=\mathtt{Q}_{\iota_{0}}$,
 $\mathcal{H}_{\gamma}[\Theta_{\iota_{0}}(\mathtt{Q}_{\iota_{0}})]=\mathcal{H}_{\gamma}[\Theta(\mathtt{Q})]$,
 and
\begin{equation}\label{eq:reduction_left}
(\mathcal{H}_{\gamma},\Theta,\mathtt{Q})\vdash^{a_{n}}_{c,\vec{d},\xi,\Lambda}
\Gamma_{n}, \lnot C_{n,\iota_{0}}^{(\bar{\rho})}
\end{equation}

SIH yields for $a^{\prime}=a^{\prime}_{0}\# a_{1}\#\cdots\#a_{n}$ 
\begin{equation}\label{eq:reduction_right}
\forall \nu\in[\bar{\rho}\mathtt{Q}]J_{n,\iota_{0}}
\left[
(\mathcal{H}_{\gamma},\Theta_{\nu},\mathtt{Q}_{\nu}
)\vdash^{\varphi_{b}(a^{\prime})}_{c,\vec{d},\xi,\Lambda}
(\Gamma\cup\{B_{\nu}^{(\bar{\rho})}\})(\mathtt{Q}_{\nu})
\right]
\end{equation}
We have $\mathrm{rk}(C_{n,\iota_{0}})<\mathrm{rk}(C_{n})\leq\mathrm{rk}(\Pi)\leq c+b$.
If $\mathrm{rk}(C_{n,\iota_{0}})<c$, then we obtain (\ref{eq:predcereg.Sa})
by a $(cut)$ with (\ref{eq:reduction_right}) and (\ref{eq:reduction_left}), where
$\max\{\varphi_{b}(a^{\prime}), a_{n}\}<\varphi_{b}(a)$ by $a_{n}<a_{0}\# a_{n}\leq a\leq \varphi_{b}(a)$.
Let
$c+b_{1}=\mathrm{rk}(C_{n,\iota_{0}})<c+b$.
MIH with (\ref{eq:reduction_right}) and (\ref{eq:reduction_left}) 
yields (\ref{eq:predcereg.Sa}), where $b_{1}\in\mathcal{H}_{\gamma}[\Theta(\mathtt{Q})]\cap M_{\mathtt{Q}}$
and $\varphi_{b_{1}}(\varphi_{b}(a^{\prime})\# a_{n})<\varphi_{b}(a)$.
\eprf

\blem\label{lem:CE}{\rm (Cut-elimination)}
Let
$(\mathcal{H}_{\gamma},\Theta,\mathtt{Q} )\vdash^{a}_{c+b,\vec{d},\xi,\Lambda}\Gamma$, where
$c+b<\mI_{N}$ and $c\in\mathcal{H}_{\gamma}[\Theta(\mathtt{Q})]$.
Then
$(\mathcal{H}_{\gamma},\Theta,\mathtt{Q})\vdash^{\varphi_{b}(a)}_{c,\vec{d},\xi,\Lambda}\Gamma$
holds.
\elem
\bprf
By induction on $a$. We may assume $b>0$. Let $\bar{\rho}=\bar{\rho}_{\mathtt{Q}}$.
Consider the case when the last inference is a $(cut)$.
We have an $a_{0}<a$ and a $\bigvee$-formula $C^{(\bar{\rho})}$ with $C\simeq\bigvee(C_{\iota})_{\iota\in J}$  such that
$\mathrm{rk}(C)<c+b$ and
$(\mathcal{H}_{\gamma},\Theta,\mathtt{Q} )\vdash^{a_{0}}_{c+b,\vec{d},\xi,\Lambda}C^{(\bar{\rho})},\Gamma$.
On the other hand we have
$(\mathcal{H}_{\gamma},\Theta_{\iota},\mathtt{Q}_{\iota} )\vdash^{a_{0}}_{c+b,\vec{d},\xi,\Lambda}
(\Gamma\cup\{\lnot C_{\iota}^{(\bar{\rho})}\})(\mathtt{Q}_{\iota})$
for each $\iota\in[\bar{\rho}\mathtt{Q}]J$ and each good extension 
$\mathtt{Q}_{\iota}=\mathtt{Q}[\mathrm{srk}(C_{\iota})]$, where 
$\Theta_{\iota}=\Theta_{\mathtt{Q}_{\iota}}\cup\mathsf{k}(C_{\iota}^{(\bar{\rho})};\mathtt{Q})$.

IH yields
$(\mathcal{H}_{\gamma},\Theta,\mathtt{ Q} )\vdash^{\varphi_{b}(a_{0})}_{c,\vec{d},\xi,\Lambda}
C^{(\bar{\rho})},\Gamma$ and
$(\mathcal{H}_{\gamma},\Theta_{\iota},\mathtt{Q}_{\iota} )\vdash^{\varphi_{b}(a_{0})}_{c,\vec{d},\xi,\Lambda}
(\Gamma\cup\{\lnot C_{\iota}^{(\bar{\rho})}\})(\mathtt{Q}_{\iota})$.
Let $c+b_{1}=\max\{c,\mathrm{rk}(C)\}<c+b$.
Then $b_{1}\in
\mathcal{H}_{\gamma}[\Theta(\mathtt{Q})]$ by (\ref{eq:controlder_cap2}).
We obtain
$(\mathcal{H}_{\gamma},\Theta,\mathtt{Q})\vdash^{\varphi_{b_{1}}(\varphi_{b}(a_{0})\cdot 2)}_{c,\vec{d},\xi,\Lambda}
\Gamma$ by
Reduction \ref{lem:predcereg.Sa}, where
$\varphi_{b_{1}}(\varphi_{b}(a_{0})\cdot 2)<\varphi_{b}(a)$.
\eprf

\bdf\label{df:recapping_I_1}
{\rm 
Let $\{\mS_{1}<\cdots<\mS_{m}\}\subset\mathrm{dom}(\mathtt{Q})$ be 
a non-empty list of successor stable ordinals.
Let  $\vec{\rho}=(\rho_{1},\ldots,\rho_{m})$ be ordinals for which
there is a special finite function $g$ such that $(\rho_{i},g)\in\mathtt{Q}(\mS_{i})$ for any $i\leq m$.

A cap $\bar{\tau}$ over $\mathtt{Q}$ is said to be \textit{good} for $\vec{\rho}$ iff
$\forall i((\bar{\tau}(\mS_{i}))_{0}=\rho_{i})$ if $\exists i((\bar{\tau}(\mS_{i}))_{0}=\rho_{i})$.

Let  $\vec{\kappa}=(\kappa_{1},\ldots,\kappa_{m})$  and
$\vec{\lambda}=(\lambda_{1},\ldots,\lambda_{m})$ be lists of ordinals such that
$\kappa_{k}\leq\lambda_{k}\leq\rho_{k}$ and there is a special finite function $g_{1}$ such that
$(\kappa_{k},g_{1})\in\Psi_{\mS_{k}}^{\mathtt{g}}$ for any $k$.

Let $A^{(\bar{\tau})}$ be a formula over $\mathtt{Q}$ such that
$\mathsf{k}(A^{(\bar{\tau})};\mathtt{Q})\subset M_{\kappa_{k}}\cap M_{\lambda_{k}}$ for each $k$,
and the cap $\bar{\tau}$ is good for $\vec{\rho}$.
Then $(A^{(\bar{\tau})})^{(\vec{\kappa}/\vec{\rho})}$
and $(A^{(\bar{\tau})})^{[\vec{\lambda}/\vec{\rho}]}$ denote formulas defined as follows.

\benu

\item
$(A^{(\bar{\tau})})^{(\vec{\kappa}/\vec{\rho})}\equiv 
A^{(\bar{\tau}_{\kappa})}$, where
$\bar{\tau}_{\kappa}(\mS_{k})=(\kappa_{k},g_{1})$ for each $k$ such that
$(\bar{\tau}(\mS_{k}))_{0}=\rho_{k}$, and $\bar{\tau}_{\kappa}(\mT)=\bar{\tau}(\mT)$ else.

\item
$(A^{(\bar{\tau})})^{[\vec{\lambda}/\vec{\rho}]}\equiv(B^{[\lambda_{k}/\mS_{k}]})^{(\bar{\tau})}$ 
if there is a $k$ such that
$A\equiv B^{[\rho_{k}/\mS_{k}]}\not\equiv B\equiv\mathrm{inv}^{(\mathtt{u})}(A^{(\bar{\tau})};\mathtt{Q})$ and 
$\rho_{k}<(\bar{\tau}(\mS_{k}))_{0}$.

Otherwise let $(A^{(\bar{\tau})})^{[\vec{\lambda}/\vec{\rho}]}\equiv A^{(\bar{\tau})}$.

\eenu

Let $\Gamma$ be a set of formulas $A^{(\bar{\tau})}$ such that
$\mathsf{k}(A^{(\bar{\tau})};\mathtt{Q})\subset M_{\kappa_{k}}\cap M_{\lambda_{k}}$ for each $k$,
and the cap $\bar{\tau}$ is good for $\vec{\rho}$.
Then let 
$\Gamma^{(\vec{\kappa}/\vec{\rho})}=
\{(A^{(\bar{\tau})})^{(\vec{\kappa}/\vec{\rho})}: A^{(\bar{\tau})}\in\Gamma\}$,
$\Gamma^{[\vec{\lambda}/\vec{\rho}]}=
\{(A^{(\bar{\tau})})^{[\vec{\lambda}/\vec{\rho}]}: A^{(\bar{\tau})}\in\Gamma\}$
and
 $\Gamma(\vec{\rho})=\{A^{(\bar{\tau})}\in\Gamma: \exists i((\bar{\tau}(\mS_{i}))_{0}=\rho_{i})\}$.

}
\edf

\blem\label{lem:recapping}{\rm (Recapping)}
Let 
$(\mathtt{Q}, \gam^{\mathtt{Q}}_{\cdot}, \delta^{\mathtt{Q}}_{\cdot})$ be
a finite family with thresholds.
Let 
$\{\mS_{1}<\cdots<\mS_{m}\}\subset\mathrm{dom}(\mathtt{Q})$ be a non-empty list of 
successor stable ordinals such that $\mS_{k}\in SSt_{i_{k}}$ with
$i=i_{1}\leq i_{2}\leq\cdots\leq i_{m}$, and $\vec{\rho}=(\rho_{1},\ldots,\rho_{m})$
 ordinals such that there is a special finite function $g$ such that
$\forall i((\rho_{i},g)\in\mathtt{Q}(\mS_{i}))$.

Let
$
(\mathcal{H}_{\gamma},\Theta, \mathtt{Q})
\vdash^{a}_{c,\vec{d},\xi,\Lambda}\Pi,\Phi, \Gamma
$, where $a<\Lambda$, $\vec{d}=(d_{1},d_{2},\ldots,d_{N})$ with
$d_{1}\leq d_{2}\leq\cdots\leq d_{N}\leq\xi<\Lambda$,
$\Pi(\vec{\rho})=\emptyset$, $\Phi(\vec{\rho})=\Phi$, 
$\Gamma(\vec{\rho})=\Gamma$ and 
the least cap $\bar{\rho}_{\mathtt{Q}}$
is good for $\vec{\rho}$.
Let $b$ be an ordinal such that $\mS_{m}+1\leq c\leq b\in\mathcal{H}_{\gamma}[\Theta(\mathtt{Q})]$,
$b<s(g)\leq \Lambda_{0}$, $b\leq d_{i_{m}}$ and $b^{\dagger}\geq\xi$.
Assume
\beqn\label{eq:recapping_PiN_0}
\forall j<i_{m} \left[\exists \mV\in\mathrm{dom}(\mathtt{Q})\cap SSt_{j}(\mS_{m}<\mV) \Rarw d_{j}\leq b\right] 
\eeqn
and
\beqn\label{eq:recapping_PiN}
\mS_{m}=\max(\mathrm{dom}(\mathtt{Q})\cap \bigcup_{j\geq i_{m}}SSt_{j})
\Rarw \Gamma \subset\bigvee(b)\cap \mathcal{L}_{i_{m}}
\eeqn

Let $g_{1}=h^{b}(g;\varphi_{d_{i_{m}}}(a))$.
Let 
$\vec{\kappa}=(\kappa_{1},\ldots,\kappa_{m})$, $\vec{\theta}=(\theta_{1},\ldots,\theta_{m})$ and
$\vec{\lambda}=(\lambda_{1},\ldots,\lambda_{m})$ be lists of ordinals such that
$\{\kappa_{k},\theta_{k}\}\subset H_{\rho_{k}}(g_{1},\mathtt{Q},\Theta)$, 
$\kappa_{k}\leq\theta_{k}\leq\lambda_{k}\leq\rho_{k}$, 
$\Theta(\mathtt{Q})\subset M_{\lambda_{k}}$ and
$(\mathtt{Q}(\mS_{k}))_{0}\cap\rho_{k}\subset\lambda_{k}$.

Let $h_{k}=g_{1}$ if $\lambda_{k}=\kappa_{k}$, and $h_{k}=g$ else.
Let $\mathtt{Q}^{[\vec{\kappa},\vec{\theta},\vec{\lambda}/\vec{\rho}]}$ be a finite family obtained from $\mathtt{Q}$
by replacing the guarded ordinal $(\rho_{k},g)$ in $\mathtt{Q}(\mS_{k})$
by $(\kappa_{k},g_{1})$, $(\theta_{k},g_{1})$ and $(\lambda_{k},h_{k})$, i.e.,
$\mathtt{Q}^{[\vec{\kappa},\vec{\theta},\vec{\lambda}/\vec{\rho}]}(\mS_{k})=
(\mathtt{Q}(\mS_{k})\cup\{ (\kappa_{k},g_{1}),(\theta_{k},g_{1}), (\lambda_{k},h_{k})\})\setminus\{(\rho_{k},g)\}$.
Then
\begin{equation}\label{eq:1.1.f_I_1}
(\mathcal{H}_{\gamma},\Theta, \mathtt{Q}^{[\vec{\kappa},\vec{\theta}, \vec{\lambda}/\vec{\rho}]}
)
\vdash^{\varphi_{d_{i_{m}}}(a)}_{c,\vec{d},\xi,\Lambda}
\Pi^{[\vec{\lambda}/\vec{\rho}]}, \Phi^{(\vec{\theta}/\vec{\rho})}, \Gamma^{(\vec{\kappa}/\vec{\rho})}
\end{equation}
holds.

\elem
\bprf
We show the lemma by induction on $a$.
Let $\mathtt{R}=\mathtt{Q}^{[\vec{\kappa},\vec{\theta},\vec{\lambda}/\vec{\rho}]}$.

We have $\varphi_{d_{i_{m}}}(a)<\Lambda<\mI_{N}$ and $g_{1}:\mI_{N}\to\Gamma(\mI_{N})$.
Also $SC_{\mI_{N}}(\mathrm{fld}(g))\subset\Lambda$ 
by (\ref{eq:fld}), and
$SC_{\mI_{N}}(\mathrm{fld}(g_{1}))\subset\Lambda$ follows.
We have $\mathrm{rk}(\Pi\cup\Phi\cup\Gamma)<\xi$ by (\ref{eq:controlder_cap_cover}).

Let
$\{\kappa_{k},\theta_{k}\}\subset H_{\rho_{k}}(g_{1},\mathtt{Q},\Theta)$.
By (\ref{eq:controlder_cap0}), Definition \ref{df:resolvent} and the assumption we have 
$\gamma\leq\gamma^{\mathtt{Q}}_{\mS}<\min\{\mathtt{p}_{0}(\kappa_{k}),\mathtt{p}_{0}(\theta_{k})\}$,
$\max\{\mathtt{p}_{0}(\kappa_{k}),\mathtt{p}_{0}(\theta_{k})\}\leq\mathtt{p}_{0}(\rho_{k})$,
$\Theta(\Theta)=\Theta(\mathtt{Q})\cap M_{\rho_{k}}
\subset M_{\kappa_{k}}\cap M_{\theta_{k}}\cap M_{\lambda_{k}}$ and $\mathtt{Q}(\mS_{k})\cap\rho_{k}\subset\kappa_{k}\subset\theta_{k}$.
We obtain $\Theta(\mathtt{Q})=\Theta(\mathtt{R})$.
On the other hand we have
$\{b,d_{i_{m}},a\}\cup\mathrm{fld}(g) \subset\mathcal{H}_{\gamma}[\Theta(\mathtt{Q})]$
by the assumption and (\ref{eq:controlder_cap22}).
Hence 
$\mathrm{fld}(g_{1})=\mathrm{fld}(h^{b}(g;\varphi_{d_{i_{m}}}(a)))\subset \mathcal{H}_{\gamma}[\Theta(\mathtt{Q})]\subset
M_{\kappa_{k}}\cap M_{\theta_{k}}$, cf.\,(\ref{eq:notationsystem.6}).

If $\rho_{k}=\rho_{\mathtt{Q}(\mS_{k})}$, then
$\rho_{\mathtt{R}(\mS_{k})}=\kappa_{k}$.
Otherwise $\rho_{\mathtt{R}(\mS_{k})}=\rho_{\mathtt{Q}(\mS_{k})}$.
In each case we obtain 
$(\mathtt{Q}(\mS_{k}))_{0}\cap\rho_{k}=(\mathtt{R}(\mS_{k}))_{0}\cap\kappa_{k}\subset
(\mathtt{R}(\mS_{k}))_{0}\cap\theta_{k}$, where
$(\mathtt{R}(\mS_{k}))_{0}\cap\theta_{k}=((\mathtt{R}(\mS_{k}))_{0}\cap\kappa_{k})\cup\{\kappa_{k}\}$
if $\kappa_{k}<\theta_{k}$, and
$(\mathtt{R}(\mS_{k}))_{0}\cap\theta_{k}=(\mathtt{R}(\mS_{k}))_{0}\cap\kappa_{k}$ else.

Let $A^{(\bar{\rho})}\in\Gamma$.
We obtain
$\mathrm{inv}^{(\mathtt{u})}((A^{(\bar{\rho})})^{(\vec{\kappa}/\vec{\rho})};\mathtt{R})\equiv\mathrm{inv}^{(\mathtt{u})}(A^{(\bar{\rho})};\mathtt{Q})$, which yields
$\mathsf{k}((A^{(\bar{\rho})})^{(\vec{\kappa}/\vec{\rho})} ; \mathtt{R})=\mathsf{k}(A^{(\bar{\rho})};\mathtt{Q})$.
We obtain $\mathcal{H}_{\gamma}[\Theta(\mathtt{R})]\subset M_{\kappa_{i}}$ by
Lemma \ref{prp:EK2}.\ref{prp:EK2.1}.
$\Theta(\mathtt{Q})=\Theta\cap M_{\mathtt{Q}}=\Theta\cap M_{\mathtt{R}}=\Theta(\mathtt{R})$
yields
$\mathsf{k}((A^{(\bar{\rho})})^{(\vec{\kappa}/\vec{\rho})};\mathtt{R})=\mathsf{k}(A^{(\bar{\rho})};\mathtt{Q})\subset\mathcal{H}_{\gamma}[\Theta(\mathtt{Q})]=\mathcal{H}_{\gamma}[\Theta(\mathtt{R})]$ by (\ref{eq:controlder_cap22}).

Let $A^{(\bar{\rho})}\in\Phi$.
We claim
$\mathrm{inv}^{(\mathtt{u})}((A^{(\bar{\rho})})^{(\vec{\theta}/\vec{\rho})};\mathtt{R})\equiv
\mathrm{inv}^{(\mathtt{u})}(A^{(\bar{\rho})};\mathtt{Q})$.
Let $A\equiv B^{[\tau/\mS_{k}]}\not\equiv B$.
If $\tau\in(\mathtt{Q}(\mS_{k}))_{0}\cap\rho_{k}$,
then $\mathrm{inv}^{(\mathtt{u})}(A^{(\bar{\rho})};\mathtt{Q})\equiv B$,
$\tau\in(\mathtt{R}(\mS_{k}))_{0}\cap\theta_{k}$
and
$\mathrm{inv}^{(\mathtt{u})}((A^{(\bar{\rho})})^{(\vec{\theta}/\vec{\rho})};\mathtt{R})\equiv B$.
Conversely let $\tau\in(\mathtt{R}(\mS_{k}))_{0}\cap\theta_{k}$.
Then $\mathrm{inv}^{(\mathtt{u})}((A^{(\bar{\rho})})^{(\vec{\theta}/\vec{\rho})};\mathtt{R})\equiv B$.
If $\tau\neq\kappa_{k}$, then $\tau\in(\mathtt{Q}(\mS_{k}))_{0}\cap\rho_{k}$ and
$\mathrm{inv}^{(\mathtt{u})}(A^{(\bar{\rho})};\mathtt{Q})\equiv B$.
Suppose $\tau=\kappa_{k}<\theta_{k}$. We have $\kappa_{k}\not\in(\mathtt{Q}(\mS_{k}))_{0}\cap\rho_{k}\subset\kappa_{k}$.
We would have 
$\mathrm{inv}^{(\mathtt{u})}(A^{(\bar{\rho})};\mathtt{Q})\equiv A\equiv B^{[\kappa_{k}/\mS_{k}]}$
and
$\mathsf{k}(B^{[\kappa_{k}/\mS_{k}]})=\mathsf{k}(A^{(\bar{\rho})};\mathtt{Q})\subset
\mathcal{H}_{\gamma}[\Theta(\mathtt{Q})]\subset M_{\kappa_{k}}$ by Lemma \ref{prp:EK2}.\ref{prp:EK2.1} and (\ref{eq:controlder_cap22}).
This is not the case since $M_{\kappa_{k}}\cap\mS_{k}\subset\kappa_{k}$, $B^{[\kappa_{k}/\mS_{k}]}\not\equiv B$,
 and
there is an $\alpha\in\mathsf{k}(B^{[\kappa_{k}/\mS_{k}]})$ such that 
$\kappa_{k}\leq\alpha=\beta[\kappa_{k}/\mS_{k}]<\mS_{k}$ for a $\beta\geq\mS_{k}$.
We obtain 
$\mathsf{k}((A^{(\bar{\rho})})^{(\vec{\theta}/\vec{\rho})};\mathtt{R})=
\mathsf{k}(A^{(\bar{\rho})};\mathtt{Q})\subset\mathcal{H}_{\gamma}[\Theta(\mathtt{Q})]=\mathcal{H}_{\gamma}[\Theta(\mathtt{R})]$ by Lemma \ref{prp:EK2}.\ref{prp:EK2.1} and (\ref{eq:controlder_cap22}).

Let $A^{(\bar{\tau})}\in\Pi$.
Then $(\bar{\tau}(\mS_{k}))_{0}\neq\rho_{k}$ for every $k$.
If $\rho_{k}
<(\bar{\tau}(\mS_{k}))_{0}$, then
$\lambda_{k}
\leq\rho_{k}<(\bar{\tau}(\mS_{k}))_{0}$ holds.
$\mathrm{inv}^{(\mathtt{u})}((A^{(\bar{\tau})})^{[\vec{\lambda}/\vec{\rho}]};\mathtt{R})\equiv\mathrm{inv}^{(\mathtt{u})}(A^{(\bar{\tau})};\mathtt{Q})$ follows.
We obtain
$\mathsf{k}((A^{(\bar{\tau})})^{[\vec{\lambda}/\vec{\rho}]};\mathtt{R})=\mathsf{k}(A^{(\bar{\tau})};\mathtt{Q})\subset
\mathcal{H}_{\gamma}[\Theta(\mathtt{Q})]
=\mathcal{H}_{\gamma}[\Theta(\mathtt{R})]$.
Therefore each of  (\ref{eq:controlder_cap2}) and (\ref{eq:controlder_cap22}) 
is enjoyed in (\ref{eq:1.1.f_I_1}).
\\
\textbf{Case 1}. 
First consider the case when the last inference is a $(j_{0} {\rm -rfl}_{\mV}(h,x,f,\Delta))$:
Then 
$\mV\in \mathrm{dom}(\mathtt{Q})\cap SSt_{j_{0}}\subset\xi$.
Let $\bar{\rho}=\bar{\rho}_{\mathtt{Q}}$. 
Then $\bar{\rho}(\mV)=(\rho_{\mathtt{Q}(\mV)},h)$.
We have  an ordinal $a_{0}<a$ and
a set
$\Delta\subset\mathcal{L}_{j_{0}}$
such that $\mT\leq\mathrm{rk}(\Delta)<d_{j_{0}}$,
$\forall\mU\in\mathrm{dom}(\mathtt{Q})\cap \bigcup_{j\geq j_{0}}SSt_{j}(\mV<\mU \Rarw \mathrm{rk}(\Delta)<\mU)$.
Also we have an ordinal $x\in\mathrm{supp}(h)$ and 
a special finite function $f$ such that $\min\{s(h), s(f)\}\geq\mV+1$,
$\mathrm{fld}(f)\subset \mathcal{H}_{\gamma}[\Theta(\mathtt{Q})]$,
$SC_{\mI_{N}}(\mathrm{fld}(f))\subset\Lambda$
and 
$f_{x}\leq h_{x} \,\&\, f<_{\mI_{N}}^{x} h^{\prime}(x)$.

Let $\mT$ be a successor $j$-stable ordinal 
such that $j\geq j_{0}$, $\mV\leq\mT\leq\mathrm{rk}(\Delta)$ and $\min\{s(h), s(f)\}\geq\mT+1$.
Let $\mathtt{P}=\mathtt{Q}[\mT]$ be a good extension of $\mathtt{Q}$ for $\gamma,\Theta$ 
such that
$(\bar{\rho_{\mathtt{Q}}}(\mV))_{1}=h=(\bar{\rho_{\mathtt{P}}}(\mT))_{1}$.
Let $\rho_{\mT}=\rho_{\mathtt{P}(\mT)}$.

We then have
$
(\mathcal{H}_{\gamma},\Theta_{\mathtt{P}_{\iota}}\cup\mathsf{k}(\delta_{\iota}^{(\bar{\rho}_{\mathtt{P}})};\mathtt{P}),
\mathtt{P}_{\iota}
)
\vdash^{a_{0}}_{c,\vec{d},\xi,\Lambda} (\Pi\cup\Phi\cup\Gamma\cup\{\lnot\delta_{\iota}^{(\bar{\rho})}\})(\mathtt{P}_{\iota})
$ 
for each $\delta\in\Delta$, $\iota\in[\bar{\rho}_{\mathtt{P}}\mathtt{P}]J$ and
each good extension $\mathtt{P}_{\iota}=\mathtt{P}[\mathrm{srk}(\delta_{\iota})]$.

On the other hand we have
$
(\mathcal{H}_{\gamma},\Theta_{\mathtt{P}},
\mathtt{P}^{\sigma}
)
\vdash^{a_{0}}_{c,\vec{d},\xi,\Lambda}
\Delta^{(\bar{\rho}_{\sigma})}, (\Pi\cup\Phi\cup\Gamma)(\mathtt{P})$
 for $\sigma\in H_{\rho_{\mT}}(f,\mathtt{P},\Theta)$, 
$\bar{\rho}_{\sigma}=\bar{\rho}_{\mathtt{P}}[\sigma/\rho_{\mT}]$ and
$\mathtt{P}^{\sigma}=\mathtt{P}\cup\{(\mT,(\sigma,f))\}$.
Note that $\bar{\rho}(\mathtt{P}_{\iota})$ as well as $\bar{\rho}_{\sigma}$ is good for $\vec{\rho}$.

Given a good extension $\mathtt{S}=\mathtt{R}[\mT]$, let $\mathtt{P}=\mathtt{Q}[\mT]$
be a good extension such that $\mathtt{P}^{[\vec{\kappa},\vec{\theta},\vec{\lambda}/\vec{\rho}]}=\mathtt{S}$.
\\
\textbf{Case 1.1}. 
$(\bar{\rho}(\mS_{1}))_{0}<\rho_{1}$:
Since $\bar{\rho}=\bar{\rho}_{\mathtt{Q}}$ is good for $\vec{\rho}$, this means
$(\bar{\rho}(\mS_{k}))_{0}<\rho_{k}$ for every $k$.
We have $(\bar{\rho}(\mS_{k}))_{0}<\kappa_{k}<\rho_{k}$.
Moreover we have $\rho_{\mathtt{Q}(\mV)}=\rho_{\mathtt{R}(\mV)}$.

$
(\mathcal{H}_{\gamma},\Theta_{\mathtt{S}_{\iota}}\cup\mathsf{k}(\delta_{\iota}^{(\bar{\rho})};\mathtt{Q}),
\mathtt{S}_{\iota}
)
\vdash^{\varphi_{d_{i_{m}}}(a_{0})}_{c,\vec{d},\xi,\Lambda} (\Pi^{[\vec{\lambda}/\vec{\rho}]}
\cup\Phi^{(\vec{\theta}/\vec{\rho})}\cup\Gamma^{(\vec{\kappa}/\vec{\rho})}\cup\{\lnot(\delta_{\iota}^{(\bar{\rho})})^{[\vec{\lambda}/\vec{\rho}]}\})(\mathtt{S}_{\iota})
$ for $\iota\in[\bar{\rho}_{\mathtt{S}}\mathtt{S}]J\subset [\bar{\rho}_{\mathtt{P}}\mathtt{P}]J$
and
$
(\mathcal{H}_{\gamma},\Theta_{\mathtt{S}},
\mathtt{S}^{\sigma}
)
\vdash^{\varphi_{d_{i_{m}}}(a_{0})}_{c,\vec{d},\xi,\Lambda}
(\Delta^{(\bar{\rho}_{\sigma})})^{[\vec{\lambda}/\vec{\rho}]}, 
(\Pi^{[\vec{\lambda}/\vec{\rho}]}\cup\Phi^{(\vec{\theta}/\vec{\rho})}\cup
\Gamma^{(\vec{\kappa}/\vec{\rho})})(\mathtt{P})
$
follow from IH and $\Theta\cap M_{\rho_{k}}\subset M_{\kappa_{k}}$.

 (\ref{eq:1.1.f_I_1}) follows by a $(j_{0} {\rm -rfl}_{\mV}(h,x,f,\Delta^{[\vec{\lambda}/\vec{\rho}]}))$.
In what follows assume $(\bar{\rho}(\mS_{k}))_{0}=\rho_{k}$ for every $k$.
We then have $(\Delta^{(\bar{\rho}_{\sigma})})(\vec{\rho})=\Delta^{(\bar{\rho}_{\sigma})}$.
\\
\textbf{Case 1.2}. $\mV\neq\mS_{k}$ for every $k$: 
Let $\mV\leq\mT\leq\mathrm{rk}(\Delta)$ and $\mT\in SSt_{j}$ with $j\geq j_{0}$.
Suppose $\mT=\mS_{k}$ for a $k\leq m$.
Then $i_{k}=j\geq j_{0}$ and $\mV<\mS_{k}\in\mathrm{dom}(\mathtt{Q})\cap SSt_{i_{k}}$.
We would have
$\mS_{k}=\mT\leq\mathrm{rk}(\Delta)<\mS_{k}$ by (r1).
Hence $\mT\neq\mS_{k}$ for every $k\leq m$.
\\
\textbf{Case 1.2.1}. $j_{0}\leq i_{m}$:
Then $\Delta\subset\mathcal{L}_{j_{0}}\subset\mathcal{L}_{i_{m}}$.
We see that the assumption (\ref{eq:recapping_PiN}) holds for $\Delta$ as follows.
Assuming
$\mS_{m}=\max(\mathrm{dom}(\mathtt{Q})\cap \bigcup_{j\geq i_{m}}SSt_{j})$, we show
$\mathrm{rk}(\Delta)<b$.
If $\mS_{m}<\mV$, then $j_{0}<i_{m}$ and $\mathrm{rk}(\Delta)<d_{j_{0}}\leq b$ by
(\ref{eq:recapping_PiN_0}).
Let $\mV<\mS_{m}$.
$\mV<\mS_{m}\in\mathrm{dom}(\mathtt{Q})\cap SSt_{i_{m}}$ yields
$\mathrm{rk}(\Delta)<\mS_{m}\leq b$ by (r1) and $j_{0}\leq i_{m}$.
Let us add $\Delta$ to $\Gamma$.
\\
$
(\mathcal{H}_{\gamma},\Theta_{\mathtt{S}_{\iota}}\cup\mathsf{k}((\delta_{\iota}^{(\bar{\rho})})^{(\vec{\kappa}/\vec{\rho})};\mathtt{S}),
\mathtt{S}_{\iota}
)
\vdash^{\varphi_{d_{i_{m}}}(a_{0})}_{c,\vec{d},\xi,\Lambda} (\Pi^{[\vec{\lambda}/\vec{\rho}]}
\cup\Phi^{(\vec{\theta}/\vec{\rho})}\cup\Gamma^{(\vec{\kappa}/\vec{\rho})}\cup
\{\lnot(\delta_{\iota}^{(\bar{\rho})})^{(\vec{\kappa}/\vec{\rho})}\})(\mathtt{S}_{\iota})
$ for $\iota\in[(\bar{\rho}[\vec{\kappa}/\vec{\rho}])\mathtt{S}]J$
follows from IH.
On the other hand we have
$
(\mathcal{H}_{\gamma},\Theta_{\mathtt{S}},
\mathtt{S}^{\sigma}
)
\vdash^{\varphi_{d_{i_{m}}}(a_{0})}_{c,\vec{d},\xi,\Lambda}
(\Delta^{(\bar{\rho}_{\sigma})})^{(\vec{\kappa}/\vec{\rho})},
(\Pi^{[\vec{\lambda}/\vec{\rho}]}\cup\Phi^{(\vec{\theta}/\vec{\rho})}\cup\Gamma^{(\vec{\kappa}/\vec{\rho})})(\mathtt{S})
$
 and 
$\mathsf{k}((\delta_{\iota}^{(\bar{\rho})})^{(\vec{\kappa}/\vec{\rho})};\mathtt{S})\subset M_{\kappa_{i}}$.
 (\ref{eq:1.1.f_I_1}) follows by a $(j_{0} {\rm -rfl}_{\mV}(h,x,f,\Delta^{(\vec{\kappa}/\vec{\rho})}))$.
\\
\textbf{Case 1.2.2}. $i_{m}<j_{0}$:
Then ordinals $\vec{\theta}$ play a r\^{o}le.
It may be the case 
$\Delta\not\subset\mathcal{L}_{i_{m}}$ or $\mathrm{rk}(\Delta)\geq b$, and
$\Delta$ does not enjoy (\ref{eq:recapping_PiN}).
In the upper sequents let us add $\Delta$ to $\Phi$.
\\
$
(\mathcal{H}_{\gamma},\Theta_{\mathtt{S}_{\iota}}\cup\mathsf{k}((\delta_{\iota}^{(\bar{\rho})})^{(\vec{\theta}/\vec{\rho})};\mathtt{S}),
\mathtt{S}_{\iota}
)
\vdash^{\varphi_{d_{i_{m}}}(a_{0})}_{c,\vec{d},\xi,\Lambda} (\Pi^{[\vec{\lambda}/\vec{\rho}]}
\cup\Phi^{(\vec{\theta}/\vec{\rho})}\cup\Gamma^{(\vec{\kappa}/\vec{\rho})}\cup
\{\lnot(\delta_{\iota}^{(\bar{\rho})})^{(\vec{\theta}/\vec{\rho})}\})(\mathtt{S}_{\iota})
$ for $\iota\in[(\bar{\rho}[\vec{\theta}/\vec{\rho}])\mathtt{S}]J$
and
$
(\mathcal{H}_{\gamma},\Theta_{\mathtt{S}},
\mathtt{S}^{\sigma}
)
\vdash^{\varphi_{d_{i_{m}}}(a_{0})}_{c,\vec{d},\xi,\Lambda}
(\Delta^{(\bar{\rho}_{\sigma})})^{(\vec{\theta}/\vec{\rho})},
(\Pi^{[\vec{\lambda}/\vec{\rho}]}\cup\Phi^{(\vec{\theta}/\vec{\rho})}\cup\Gamma^{(\vec{\kappa}/\vec{\rho})})(\mathtt{S})
$ follow from IH.
On the other hand we have
$\mathsf{k}((\delta_{\iota}^{(\bar{\rho})})^{(\vec{\theta}/\vec{\rho})};\mathtt{S})\subset M_{\theta_{k}}$.
 (\ref{eq:1.1.f_I_1}) follows by a $(j_{0} {\rm -rfl}_{\mV}(h,x,f,\Delta^{(\vec{\theta}/\vec{\rho})}))$.
\\
\textbf{Case 1.3}.
$\mV=\mS_{k}$ for a $k$:
Then $\mV\in SSt_{i_{k}}$ with $j_{0}=i_{k}\geq i$, $h=g$, $\rho_{\mathtt{R}(\mS_{k})}=\kappa_{k}$, and 
$\bar{\rho}[\vec{\kappa}/\vec{\rho}]=\bar{\rho}_{\mathtt{R}}$.
We have 
$\Delta\subset\bigvee(x)\cap\mathcal{L}_{i_{k}}$, $b_{0}=\mathrm{rk}(\Delta)<d_{i_{k}}$,
and $b_{0}<\mS_{k+1}$ by (r1) if $k<m$.

Let $\mT$ be a successor $j$-stable ordinal 
such that $j\geq i_{k}$, $\mS_{k}\leq\mT\leq\mathrm{rk}(\Delta)$ and $\min\{s(g), s(f)\}\geq\mT+1$,
and $\mathtt{P}=\mathtt{Q}[\mT]$ a good extension such that
$\mathtt{P}(\mT)=\{(\rho_{\mT},g)\}$.
Then the least cap $\bar{\rho}_{\mathtt{P}}$ is good for $\vec{\rho}*(\rho_{\mT})$
with respect to the list $\{\mS_{1}<\cdots<\mS_{k}<\mT<\mS_{k+1}<\cdots<\mS_{m}\}$ if 
$\mS_{k}<\mT$.
Let $\kappa_{\mT}=\kappa_{k}$, $\theta_{\mT}=\theta_{k}$ and $\lambda_{\mT}=\lambda_{k}$ 
if $\mT=\mS_{k}$.
Otherwise let $\kappa_{\mT}\in H_{\rho_{\mT}}(g_{1},\mathtt{P},\Theta)$ and
$\kappa_{\mT}=\theta_{\mT}\leq \lambda_{\mT}=\rho_{\mT}$.
We have
\begin{equation}\label{eq:Case1left_I_1}
(\mathcal{H}_{\gamma},\Theta_{\mathtt{P}_{\iota}}\cup\mathsf{k}(\delta_{\iota}^{(\bar{\rho}_{\mathtt{P}})};\mathtt{P}),
\mathtt{P}_{\iota}
)
\vdash^{a_{0}}_{c,\vec{d}, \xi,\Lambda}
(\Pi\cup\Phi\cup\Gamma\cup\{\lnot\delta_{\iota}^{(\bar{\rho})}\})(\mathtt{P}_{\iota})
\end{equation}
for each $\delta\in\Delta$, $\iota\in[\bar{\rho}_{\mathtt{P}}\mathtt{P}]J$ and each good extension $\mathtt{P}_{\iota}=\mathtt{P}[\mathrm{srk}(\delta_{\iota})]$.

We have $(\lnot\delta_{\iota})\in\bigvee(b_{0})$,
$b_{0}< x\leq s(g)$ and $b_{0}<d_{i_{k}}$.
On the other hand we have for $\sigma\in H_{\rho_{\mT}}(f,\mathtt{P},\Theta)$ and
$\mathtt{P}^{\sigma}=\mathtt{P}\cup\{(\mT,(\sigma,f))\}$
\begin{equation}\label{eq:KppiNlowerCase1b_I_1}
(\mathcal{H}_{\gamma},\Theta_{\mathtt{P}},
\mathtt{P}^{\sigma}
)
\vdash^{a_{0}}_{c,\vec{d},\xi,\Lambda}
\Delta^{(\bar{\rho}_{\sigma})}, (\Pi\cup\Phi\cup\Gamma)(\mathtt{P})
\end{equation}
where $\bar{\rho}_{\sigma}=\bar{\rho}_{\mathtt{P}}[\sigma/\rho_{\mT}]$ and
$\bar{\rho}_{\sigma}=\bar{\rho}_{\mathtt{P}^{\sigma}}$ is good for $\vec{\rho}*(\sigma)$.

We have $\Theta(\mathtt{P})\cap M_{\rho_{\mT}}
\subset M_{\sigma}$ for $\sigma=\rho_{\mathtt{P}^{\sigma}(\mT)}$.
$f$ is a finite function such that $x\in\mathrm{supp}(g)$ and
\beqn\label{eq:recapping_Pi11_r2_I_1}
f_{x}\leq g_{x} \,\&\, f<_{\mI_{N}}^{x}g^{\prime}(x) \,\&\, 
\mathrm{fld}(f)\subset\mathcal{H}_{\gamma}[\Theta(\mathtt{Q})]
\spand SC_{\mI_{N}}(\mathrm{fld}(f))\subset\Lambda
\eeqn
Let $\mathtt{T}=
(\mathtt{P}^{[\kappa_{\mT}, \theta_{\mT},\lambda_{\mT}/\rho_{\mT}]})
^{[\vec{\kappa},\vec{\theta},\vec{\lambda}/\vec{\rho}]}$.
If $\mT=\mS_{k}$, then $\mathtt{T}=\mathtt{R}$.
Otherwise $\mathtt{T}=
\mathtt{P}^{[\vec{\kappa},\kappa_{\mT},\vec{\theta},\theta_{\mT},\vec{\lambda},\lambda_{\mT}/\vec{\rho},\rho_{\mT}]}=
\mathtt{S}^{[\kappa_{\mT}, \theta_{\mT},\lambda_{\mT}/\rho_{\mT}]}$.

Let
$\Pi_{1}=\Pi(\mathtt{P})^{[\vec{\lambda}/\vec{\rho}]}$,
$\Phi_{1}=(\Phi(\mathtt{P})^{(\kappa_{\mT}/\rho_{\mT})})^{(\vec{\theta}/\vec{\rho})}$ and
$\Gamma_{1}=(\Gamma(\mathtt{P})^{(\kappa_{\mT}/\rho_{\mT})})^{(\vec{\kappa}/\vec{\rho})}$.
\\
\textbf{Case 1.3.1}. $b_{0}< b$: 
Let $\delta\in\Delta$, and $\iota\in J$ for $\delta\simeq\bigvee(\delta_{\iota})_{\iota\in J}$.
We have $\mathrm{rk}(\lnot\delta_{\iota})<\mathrm{rk}(\lnot\delta)\leq b$.
Hence the assumption (\ref{eq:recapping_PiN}) is enjoyed for $\Gamma\cup\{\lnot\delta_{\iota}^{(\bar{\rho})}\}$.
Let
$\iota\in[\bar{\rho}_{\mathtt{T}}\mathtt{T}]J$ and $n\leq m$.
We obtain
$\mathsf{k}(\delta_{\iota}^{(\bar{\rho}_{\mathtt{P}})};\mathtt{P})=\mathsf{k}(\delta_{\iota})=
\mathsf{k}(\delta_{\iota}^{(\bar{\rho}_{\mathtt{T}})};\mathtt{T})\subset 
M_{\kappa_{n}}\cap M_{\kappa_{\mT}}\subset M_{\theta_{n}} \cap M_{\theta_{\mT}}$, 
and $\iota\in[\bar{\rho}_{\mathtt{P}}\mathtt{P}]J$
by $\rho_{\mathtt{Q}(\mS_{n})}=\rho_{n}$,
$\rho_{\mathtt{R}(\mS_{n})}=\kappa_{n}$, $\rho_{\mathtt{T}(\mT)}=\kappa_{\mT}$
 and $\mathsf{k}(\delta_{\iota})\subset M_{\mathtt{T}}\subset M_{\mathtt{P}}$.
This yields 
$\{\kappa_{n},\theta_{n}\}\subset H_{\rho_{n}}(g_{1},\mathtt{P},\Theta
\cup\mathsf{k}(\delta_{\iota}^{(\bar{\rho}_{\mathtt{T}})};\mathtt{T}))$ and
 $\{\kappa_{\mT},\theta_{\mT}\}\subset H_{\rho_{\mT}}(g_{1},\mathtt{P},\Theta\cup\mathsf{k}(\delta_{\iota}^{(\bar{\rho}_{\mathtt{T}})};\mathtt{T}))$.
By IH with (\ref{eq:Case1left_I_1}) 
we obtain
\begin{equation}\label{eq:Case1cpi11_I_1}
(\mathcal{H}_{\gamma},\Theta_{\mathtt{T}_{\iota}}\cup\mathsf{k}(\delta_{\iota}^{(\bar{\rho}_{\mathtt{T}})};\mathtt{T}),
\mathtt{T}_{\iota}
)
\vdash^{\varphi_{d_{i_{m}}}(a_{0})}_{c,\vec{d},\xi,\Lambda}
(\Pi_{1}\cup
\Phi_{1}\cup
\Gamma_{1}\cup
\{\lnot\delta_{\iota}^{(\bar{\rho}_{\mathtt{T}})}\})(\mathtt{T}_{\iota})
\end{equation}
where $\delta_{\iota}^{(\bar{\rho}_{\mathtt{T}})}\equiv 
((\delta_{\iota}^{(\bar{\rho}_{\mathtt{P}})})^{(\kappa_{\mT}/ \rho_{\mT})})^{(\vec{\kappa}/\vec{\rho})}$.

Let $x_{1}=\min\{b,x\}>b_{0}$.
The following (\ref{eq:zigpi11.111_I_1}) is seen from Propositions \ref{prp:idless} and \ref{prp:hstepdown}.\ref{prp:hstepdown.2}.
\begin{equation}\label{eq:zigpi11.111_I_1}
f_{x_{1}}\leq (g_{1})_{x_{1}} \,\&\,
f<_{\mI_{N}}^{x_{1}}g_{1}^{\prime}(x_{1}) 
\end{equation}

Let $\sigma\in H_{\kappa_{\mT}}(f,\mathtt{P},\Theta)$.
Then $\Theta\cap M_{\rho_{\mT}}=\Theta\cap M_{\kappa_{\mT}}
\subset M_{\sigma}$ and $\sigma\in H_{\rho_{\mT}}(f,\mathtt{P},\Theta)$ by 
$\sigma<\kappa_{\mT}<\rho_{\mT}$.
By IH with (\ref{eq:KppiNlowerCase1b_I_1}) we obtain 
\begin{equation}\label{eq:1.1.fl_I_1}
(\mathcal{H}_{\gamma},\Theta_{\mathtt{T}},
\mathtt{T}\cup\{(\mT,(\sigma,f))\}
)
\vdash^{\varphi_{d_{i_{m}}}(a_{0})}_{c,\vec{d},\xi,\Lambda}
\Delta^{(\bar{\rho}_{\sigma}[\vec{\kappa}/\vec{\rho}])}, 
\Pi_{1},
\Phi_{1},
 \Gamma_{1}
\end{equation}
(\ref{eq:1.1.f_I_1})
follows by an inference $(j_{0} {\rm -rfl}_{\mV}(g_{1},x_{1},f,\Delta^{(\vec{\kappa}/\vec{\rho})}))$ with
(\ref{eq:zigpi11.111_I_1}), (\ref{eq:Case1cpi11_I_1}) and (\ref{eq:1.1.fl_I_1}).
\\
\textbf{Case 1.3.2}. 
$b\leq b_{0}$: In this case we need the assumption (\ref{eq:recapping_PiN}).
Let $\delta^{(\bar{\rho}_{\sigma})}\in\Delta^{(\bar{\rho}_{\sigma})}$ with $\bar{\rho}_{\sigma}=\bar{\rho}_{\mathtt{P}^{\sigma}}$.
By (\ref{eq:controlder_cap2}) we obtain $\mathrm{inv}^{(\mathtt{u})}(\delta^{(\bar{\rho}_{\sigma})};\mathtt{P}^{\sigma})\equiv \delta$,
and $\mathsf{k}(\delta)\subset \mathcal{H}_{\gamma}[\Theta(\mathtt{Q})]$. This yields
$b_{0}=\mathrm{rk}(\Delta)\in\mathcal{H}_{\gamma}[\Theta(\mathtt{P})]$.
We have $b\leq b_{0}< x\leq s(g)$ and $b_{0}<d_{i_{k}}$.
Let $g_{0}=h^{b_{0}}(g; \varphi_{d_{i_{m}}}(a_{0}))*f^{b_{0}+1}$, and
$\sigma\in H_{\rho_{\mT}}(g_{0},\mathtt{Q},\Theta)$.
We have 
$g_{0}\leq m(\sigma)$ 
by Definition \ref{df:resolvent}, 
$SC_{\mI_{N}}(\mathrm{fld}(g_{0}))\subset\Lambda$,  
and $\sigma\in H_{\rho_{\mT}}(h^{b_{0}}(g; \varphi_{d_{i_{m}}}(a_{0})),\mathtt{P},\Theta)$.
We see $\mathrm{fld}(g_{0})\subset\mathcal{H}_{\gamma}[\Theta(\mathtt{S})]$
from $\{b_{0},a_{0},d_{i_{m}}\}\cup \mathrm{fld}(g)\cup \mathrm{fld}(f)\subset\mathcal{H}_{\gamma}[\Theta(\mathtt{S})]$.

Suppose 
that there is a $\mV=\mS_{k}<\mU\in\mathrm{dom}(\mathtt{Q})\cap \bigcup_{j\geq i_{k}}SSt_{j}$.
Then we would have $b_{0}=\mathrm{rk}(\Delta)<
\mU<\xi\leq b^{\dagger}$ 
by (r1), and $b_{0}<b$.
Hence 
$\mS_{k}=\max(\mathrm{dom}(\mathtt{Q})\cap \bigcup_{j\geq i_{k}}SSt_{j})$.
This means that
$\mS_{k}=\mS_{m}$ and $i_{k}=i_{m}$.
By the assumption (\ref{eq:recapping_PiN}) we have $\Gamma\subset\bigvee(b)\cap \mathcal{L}_{i_{m}}$
with $b\leq b_{0}$. 
We obtain
$\Gamma\cup\{\lnot\delta_{\iota}^{(\bar{\rho})}\}\subset\bigvee(b_{0})\cap\mathcal{L}_{i_{m}}$.

Let $\mT$ be a successor $j$-stable ordinal 
such that $j\geq i_{m}$, $\mS_{m}\leq\mT\leq\mathrm{rk}(\Gamma)<b\leq b_{0}=\mathrm{rk}(\Delta)$ and 
$\min\{s(g), s(f)\}\geq\mT+1$,
and $\mathtt{P}=\mathtt{Q}[\mT]$ a good extension such that
$\mathtt{P}(\mT)=\{(\rho_{\mT},g)\}$.

Let $\mathtt{U}=\mathtt{T}^{[\sigma/\kappa_{\mT}]}=\mathtt{S}^{[\sigma/\kappa_{\mT}]}$ and
$\sigma\in H_{\kappa_{\mT}}(g_{0},\mathtt{Q},\Theta)$.
Then $\sigma<\kappa_{\mT}\leq\theta_{\mT}\leq\lambda_{\mT}\leq\rho_{\mT}$.
Let $\Gamma_{1}^{(\sigma/\kappa_{\mT})}=
(\Gamma(\mathtt{P})^{(\sigma/\rho_{\mT})})^{(\vec{\kappa}/\vec{\rho})}$.

As in \textbf{Case 1.3.1}, we obtain 
$(\mathcal{H}_{\gamma},\Theta_{\mathtt{U}_{\iota}}\cup\mathsf{k}(\delta_{\iota}),
\mathtt{U}_{\iota}
)
\vdash^{\varphi_{d_{i_{m}}}(a_{0})}_{c,\vec{d},\xi,\Lambda}
(\Pi_{1}\cup\Phi_{1}\cup
\Gamma_{1}^{(\sigma/\kappa_{\mT})}\cup
\{\lnot\delta_{\iota}^{(\bar{\rho}_{\mathtt{U}})}\})(\mathtt{U}_{\iota})$ for 
$\iota\in[\bar{\rho}_{\mathtt{U}}\mathtt{U}]J$
by $\sigma=\rho_{\mathtt{U}(\mS_{m})}$, IH and (\ref{eq:Case1left_I_1}),
where $\bar{\rho}_{\mathtt{U}}=\bar{\rho}_{\mathtt{S}}[\sigma/\kappa_{\mT}]$.

For $\mathtt{S}(\mT)\cup\{(\sigma,g_{0})\}=\mathtt{U}(\mT)\cup\{(\kappa_{\mT},g_{1})\}$, we have
$\sigma=\rho_{\mathtt{U}(\mT)}=\rho_{\mathtt{U}(\mT)\cup\{(\kappa_{\mT},g_{1})\}}=\rho_{\mathtt{S}(\mT)\cup\{(\sigma,g_{0})\}}$.
Hence (\ref{eq:controlder_cap22}) holds by adding the guarded ordinal $(\kappa_{\mT},g_{1})$ to $\mathtt{U}(\mT)$.

Let $C^{(\bar{\tau})}$ be a formula occurring in a derivation witnessing the fact
$(\mathcal{H}_{\gamma},\Theta_{\mathtt{U}_{\iota}}\cup\mathsf{k}(\delta_{\iota}),
\mathtt{U}_{\iota}
)
\vdash^{\varphi_{d_{i_{m}}}(a_{0})}_{c,\vec{d},\xi,\Lambda}
(\Pi_{1}\cup\Phi_{1}\cup
\Gamma_{1}^{(\sigma/\kappa_{\mT})}\cup
\{\lnot\delta_{\iota}^{(\bar{\rho}_{\mathtt{U}})}\})(\mathtt{U}_{\iota})$.
$\bar{\tau}$ is a cap over $\mathtt{U}_{\iota}$.
Then $(\bar{\tau}(\mT))_{0}\not\in\{\rho_{\mT},\kappa_{\mT}\}$ and
$\mathrm{inv}(C^{(\bar{\tau})};\mathtt{U})\equiv\mathrm{inv}(C^{(\bar{\tau})};\mathtt{S}^{\sigma})$,
where $\mathtt{S}^{\sigma}=\mathtt{S}\cup\{(\mT,(\sigma,g_{0}))\}$.
We obtain 
\begin{equation}\label{eq:L4.10case1rxi_I_1}
(\mathcal{H}_{\gamma},\Theta_{\mathtt{S}_{\iota}}\cup\mathsf{k}(\delta_{\iota}),
\mathtt{S}^{\sigma}_{\iota}
)
\vdash^{\varphi_{d_{i_{m}}}(a_{0})}_{c,\vec{d},\xi,\Lambda}
(\Pi_{1}\cup\Phi_{1}\cup
\Gamma_{1}^{(\sigma/\kappa_{\mT})}\cup
\{\lnot\delta_{\iota}^{(\bar{\rho}_{\mathtt{U}})}\})(\mathtt{S}_{\iota})
\end{equation}

We have
$f\leq g_{0}\leq m(\sigma)$.
Hence $\sigma\in H_{\rho_{\mT}}(f,\mathtt{P},\Theta)$.
We obtain by IH and (\ref{eq:KppiNlowerCase1b_I_1})
\[
(\mathcal{H}_{\gamma},\Theta_{\mathtt{S}},
\mathtt{S}\cup\{(\mT,(\sigma,f))\}
)
\vdash^{\varphi_{d_{i_{m}}}(a_{0})}_{c,\vec{d},\xi,\Lambda}
\Delta^{(\bar{\rho}_{\sigma,f})}, 
\Pi_{1}, \Phi_{1}, \Gamma_{1}
\]
where $\bar{\rho}_{\sigma,f}=(\bar{\rho}_{\sigma}[\vec{\kappa}/\vec{\rho}])[(\sigma,f)/(\sigma,g_{0})]$.
Lemma \ref{lem:guard_change} with $f\leq g_{0}$ yields
\begin{equation}\label{eq:KppiNlowerCase1b.12_I_1}
(\mathcal{H}_{\gamma},\Theta_{\mathtt{S}},
\mathtt{S}^{\sigma}
)
\vdash^{\varphi_{d_{i_{m}}}(a_{0})}_{c,\vec{d},\xi,\Lambda}
\Delta^{(\bar{\rho}_{\sigma,g_{0}})}, 
\Pi_{1}, \Phi_{1}, \Gamma_{1}
\end{equation}
where $\mathtt{S}^{\sigma}=\mathtt{S}\cup\{(\mT,(\sigma,g_{0}))\}$,
$\bar{\rho}_{\sigma,f}(\mT)=(\sigma,f)$ and $\bar{\rho}_{\sigma,g_{0}}(\mT)=(\sigma,g_{0})$
with $\bar{\rho}_{\sigma,g_{0}}=\bar{\rho}_{\sigma}[\vec{\kappa}/\vec{\rho}]=\bar{\rho}_{\mathtt{U}}$.

We have $\sigma=\rho_{\mathtt{S}(\mT)\cup\{(\sigma,g_{0})\}}$, and
$\mathrm{rk}(\delta)\leq b_{0}<d_{i_{m}}$.
Let $\#(\Delta)=n$ be the number of formulas in the set $\Delta$.
Reduction \ref{lem:predcereg.Sa} with (\ref{eq:L4.10case1rxi_I_1}) and (\ref{eq:KppiNlowerCase1b.12_I_1})
yields for every $\sigma\in H_{\kappa_{\mT}}(g_{0},\mathtt{S},\Theta) $
\begin{equation}\label{eq:L4.10case134_I_1}
(\mathcal{H}_{\gamma},\Theta_{\mathtt{S}},
\mathtt{S}^{\sigma}
)
\vdash^{a_{1}}_{c,\vec{d},\xi,\Lambda}
\Gamma_{1}^{(\sigma/\kappa_{\mT})},\Pi_{1}, \Phi_{1},\Gamma_{1}
\end{equation}
where
$2b\leq a_{1}=\varphi_{b_{0}}(\varphi_{d_{i_{m}}}(a_{0})\cdot(n+1))<\varphi_{d_{i_{m}}}(a)$
 by $b\leq b_{0}<d_{i_{m}}$ and $a_{0}<a$.

On the other, Tautology \ref{lem:tautology.cap}.\ref{lem:tautology.cap1}  yields
for each $A^{(\bar{\rho})}\in\Gamma(\mathtt{P})$ with $A\simeq\bigvee(A_{\iota})_{\iota\in J}$ and
$\iota\in[\bar{\rho}_{\kappa}\mathtt{S}]J$
\begin{equation}\label{eq:L4.10case1.1b_I_1}
(\mathcal{H}_{\gamma},\Theta_{\mathtt{S}_{\iota}}\cup\mathsf{k}( 
((A_{\iota}^{(\bar{\rho})})^{(\kappa_{\mT}/\rho_{\mT})})^{(\vec{\kappa}/\vec{\rho})}
;\mathtt{S}), 
\mathtt{S}_{\iota}
)\vdash^{2b}_{c,\vec{d},\xi,\Lambda}
(\Gamma_{1}\cup
\{\lnot ((A_{\iota}^{(\bar{\rho})})^{(\kappa_{\mT}/\rho_{\mT})})^{(\vec{\kappa}/\vec{\rho})}
\})
(\mathtt{S}_{\iota})
\end{equation}
where 
$\mathtt{S}_{\iota}=\mathtt{S}[\mathrm{srk}(A_{\iota})]$ is a good extension of $\mathtt{S}$.

For $g_{1}=h^{b}(g; \varphi_{d_{i_{m}}}(a))$, we obtain
$(g_{0})_{b}=g_{b}=(g_{1})_{b}$ and $g_{0}<_{\mI_{N}}^{b}g_{1}^{\prime}(b)$
by Proposition \ref{prp:hstepdown}.\ref{prp:hstepdown.4}.

Each of (\ref{eq:L4.10case1.1b_I_1}) and (\ref{eq:L4.10case134_I_1}) holds for
every successor $j$-stable ordinal  $\mT$ such that $j\geq i_{m}$,
$\mS_{m}\leq\mT\leq\mathrm{rk}(\Gamma)$ and 
$\min\{s(g), s(f)\}\geq\mT+1$, and each good extension $\mathtt{P}=\mathtt{Q}[\mT]$
such that $\mathtt{P}(\mT)=\{(\rho_{\mT},g)\}$.

By an inference rule 
$(i_{m}{\rm -rfl}_{\mS_{m}}(g_{1},b,g_{0},\Gamma^{(\vec{\kappa}/\vec{\rho})}))$ 
with 
its resolvent class 
$H_{\kappa_{m}}(g_{0},\mathtt{R},\Theta)$,
we conclude (\ref{eq:1.1.f_I_1}) by (\ref{eq:L4.10case1.1b_I_1}) and (\ref{eq:L4.10case134_I_1})
 with 
$\kappa_{m}=\rho_{\mathtt{R}(\mS_{m})}$,
$\mathrm{rk}(\Gamma)<b\leq d_{i_{m}}$ and $\Gamma\subset\bigvee(b)\cap\mathcal{L}_{i_{m}}$.
\\
\textbf{Case 2}. 
Second consider the case when the last inference $(\bigvee)$ introduces a $\bigvee$-formula
$A^{(\bar{\rho})}\in\Phi\cup\Gamma$:
Let $A\simeq\bigvee\left(A_{\iota}\right)_{\iota\in J}$.
We have an $\iota\in [\bar{\rho}\mathtt{Q}]J$ such that 
$\mathsf{k}(A_{\iota}^{(\bar{\rho})};\mathtt{Q})\subset\mathcal{H}_{\gamma}[\Theta(\mathtt{Q})]
=\mathcal{H}_{\gamma}[\Theta(\mathtt{R})]$, 
and
for $A_{\iota}\simeq\bigwedge(B_{\nu})_{\nu\in I}$ and $\nu\in[\bar{\rho}\mathtt{Q}]I$,
$
(\mathcal{H}_{\gamma},\Theta_{\mathtt{Q}_{\nu}}\cup\mathsf{k}(B_{\nu}^{(\bar{\rho})};\mathtt{Q}),\mathtt{Q}_{\nu}
)
\vdash^{a_{0}}_{c,\vec{d},\xi,\Lambda}
(\Pi\cup\Gamma\cup\{B_{\nu}^{(\bar{\rho})}\})(\mathtt{Q}_{\nu})
$ holds for each good extension $\mathtt{Q}_{\nu}= \mathtt{Q}[\mathrm{srk}(B_{\nu})]$.

Let $\bar{\rho}_{\theta}=\bar{\rho}[\vec{\theta}/\vec{\rho}]$ and $\nu\in[\bar{\rho}_{\theta}\mathtt{R}]I$
when $A^{(\bar{\rho})}\in\Phi$,
while $\bar{\rho}_{\kappa}=\bar{\rho}[\vec{\kappa}/\vec{\rho}]$ and $\nu\in[\bar{\rho}_{\kappa}\mathtt{R}]I$
when $A^{(\bar{\rho})}\in\Gamma$.
We see
$\mathrm{inv}^{(\mathtt{u})}(B_{\nu}^{(\bar{\rho})};\mathtt{Q})=\mathrm{inv}^{(\mathtt{u})}(B_{\nu}^{(\bar{\rho}_{\theta})};\mathtt{R})$
and
$\mathrm{inv}^{(\mathtt{u})}(B_{\nu}^{(\bar{\rho})};\mathtt{Q})=\mathrm{inv}^{(\mathtt{u})}(B_{\nu}^{(\bar{\rho}_{\kappa})};\mathtt{R})$, resp.
from $(\mathtt{Q}(\mS_{k}))_{0}\cap\rho_{k}=(\mathtt{R}(\mS_{k}))_{0}\cap\kappa_{k}\subset
(\mathtt{R}(\mS_{k}))_{0}\cap\theta_{k}$.
This yields
$\nu\in[\bar{\rho}\mathtt{Q}]I$,
$\mathsf{k}(B_{\nu}^{(\bar{\rho})};\mathtt{Q})=\mathsf{k}(B_{\nu}^{(\bar{\rho}_{\theta})};\mathtt{R})\subset M_{\theta_{k}}$ and 
$\mathsf{k}(B_{\nu}^{(\bar{\rho})};\mathtt{Q})=\mathsf{k}(B_{\nu}^{(\bar{\rho}_{\kappa})};\mathtt{R})\subset M_{\kappa_{k}}$, resp.
We obtain
$\{\kappa_{k},\theta_{k}\}\subset H_{\rho_{k}}(g_{1},\mathtt{Q},\Theta\cup\mathsf{k}(B_{\nu}^{(\bar{\rho})};\mathtt{Q}))$.

Likewise we see that
$\mathsf{k}(A_{\iota}^{(\bar{\rho}_{\theta})};\mathtt{R})=\mathsf{k}(A_{\iota}^{(\bar{\rho})};\mathtt{Q})
\subset\mathcal{H}_{\gamma}[\Theta(\mathtt{R})]\cap M_{\mathtt{R}}$
when $A^{(\bar{\rho})}\in\Phi$,
and
$\mathsf{k}(A_{\iota}^{(\bar{\rho}_{\kappa})};\mathtt{R})=\mathsf{k}(A_{\iota}^{(\bar{\rho})};\mathtt{Q})
\subset\mathcal{H}_{\gamma}[\Theta(\mathtt{R})]\cap M_{\mathtt{R}}$
when $A^{(\bar{\rho})}\in\Gamma$.
When $A^{(\bar{\rho})}\in\Phi$,
IH 
yields
$
(\mathcal{H}_{\gamma},\Theta_{\mathtt{R}_{\iota}}\cup\mathsf{k}(B_{\nu}^{(\bar{\rho}_{\theta})};\mathtt{R}),
\mathtt{R}_{\iota}
)\vdash^{\varphi_{d_{i_{m}}}(a_{0})}_{c,\vec{d},\xi,\Lambda}
(\Pi^{[\vec{\lambda}/\vec{\rho}]}\cup\Phi^{(\vec{\theta}/\vec{\rho})}\cup
\Gamma^{(\vec{\kappa}/\vec{\rho})}\cup\{B_{\nu}^{(\bar{\rho}_{\theta})}\})(\mathtt{R}_{\iota})
$ for $\nu\in[\bar{\rho}_{\theta}\mathtt{R}]I$.
Let $A^{(\bar{\rho})}\in\Gamma$. We obtain by IH,
$
(\mathcal{H}_{\gamma},\Theta_{\mathtt{R}_{\iota}}\cup\mathsf{k}(B_{\nu}^{(\bar{\rho}_{\kappa})};\mathtt{R}),
\mathtt{R}_{\iota}
)\vdash^{\varphi_{d_{i_{m}}}(a_{0})}_{c,\vec{d},\xi,\Lambda}
(\Pi^{[\vec{\lambda}/\vec{\rho}]}\cup\Phi^{(\vec{\theta}/\vec{\rho})}\cup\Gamma^{(\vec{\kappa}/\vec{\rho})}\cup\{B_{\nu}^{(\bar{\rho}_{\kappa})}\})(\mathtt{R}_{\iota})
$ for $\nu\in[\bar{\rho}_{\kappa}\mathtt{R}]I$,
where $\Gamma\cup\{B_{\nu}\}\subset\bigvee(b)\cap\mathcal{L}_{i_{m}}$ when 
$\mS_{m}=\max(\mathrm{dom}(\mathtt{Q})\cap \bigcup_{j\geq i_{m}}SSt_{j})$.
We obtain (\ref{eq:1.1.f_I_1})
by a $(\bigvee)$ in each case.
\\
\textbf{Case 3}.
Third consider the case when the last inference $(\bigvee)$ introduces a $\bigvee$-formula
$A^{(\bar{\tau})}\in\Pi$ with $\tau_{k}=(\bar{\tau}(\mS_{k}))_{0}\neq\rho_{k}$ for each $k$:
Although in (\ref{eq:singlemainl_I_1_3}) below, we have $\vec{\kappa}=\vec{\theta}=\vec{\lambda}$,
it may be the case $\kappa_{k}<\theta_{k}\leq\lambda_{k}$ in \textbf{Case 1.2.2}.

Let $A\simeq\bigvee\left(A_{\iota}\right)_{\iota\in J}$.
We have an $\iota\in [\bar{\tau}\mathtt{Q}]J$ such that 
$\mathsf{k}(A_{\iota}^{(\bar{\tau})};\mathtt{Q})\subset
\mathcal{H}_{\gamma}[\Theta(\mathtt{Q})]\cap M_{\mathtt{Q}}=
\mathcal{H}_{\gamma}[\Theta(\mathtt{R})]\cap M_{\mathtt{R}}$, 
and
for $A_{\iota}\simeq\bigwedge(B_{\nu})_{\nu\in I}$ and $\nu\in[\bar{\tau}\mathtt{Q}]I$,
$
(\mathcal{H}_{\gamma},\Theta_{\mathtt{Q}_{\nu}}\cup\mathsf{k}(B_{\nu}^{(\bar{\tau})};\mathtt{Q}), \mathtt{Q}_{\nu}
)
\vdash^{a_{0}}_{c,\vec{d},\xi,\Lambda}
(\{B_{\nu}^{(\bar{\tau})}\}\cup\Pi\cup\Phi\cup\Gamma)(\mathtt{Q}_{\nu})
$ holds for each good extension $\mathtt{Q}_{\nu}=\mathtt{Q}[\mathrm{srk}(B_{\nu})]$.
We claim that 
\begin{equation}\label{eq:clm_Case3_Recapping_I_1}
\mathrm{inv}(A_{\iota}^{(\bar{\tau})};\mathtt{Q})\equiv 
\mathrm{inv}((A_{\iota}^{(\bar{\tau})})^{[\vec{\lambda}/\vec{\rho}]};\mathtt{R})
\end{equation}
We have $\rho_{k}\neq\tau_{k}=(\bar{\tau}(\mS_{k}))_{0}$.
Let
$A_{\iota}\equiv C^{[\sigma/\mS_{k}]}\not\equiv C$ for a $C$ and a 
$\sigma\in(\mathtt{Q}(\mS_{k}))_{0}\cap\tau_{k}$.
Then $\mathrm{inv}^{(\mathtt{u})}(A_{\iota}^{(\bar{\tau})};\mathtt{Q})\equiv C$ and 
$\mathsf{k}(C)\subset \mathcal{H}_{\gamma}[\Theta(\mathtt{Q})]\cap M_{\mathtt{Q}}\subset M_{\mathtt{R}}$.
If $\sigma\neq\rho_{k}$, then 
$(A_{\iota}^{(\bar{\tau})})^{[\vec{\lambda}/\vec{\rho}]}\equiv A_{\iota}^{(\bar{\tau})}$ and
$\mathrm{inv}^{(\mathtt{u})}((A_{\iota}^{(\bar{\tau})})^{[\vec{\lambda}/\vec{\rho}]};\mathtt{R})\equiv C$
by $\sigma\in(\mathtt{R}(\mS_{k}))_{0}\cap\tau_{k}$.
Let $\sigma=\rho_{k}<\tau_{k}$. Then 
$(A_{\iota}^{(\bar{\tau})})^{[\vec{\lambda}/\vec{\rho}]}\equiv 
C^{[\lambda_{k}/\mS_{k}]}$ and
$\mathrm{inv}^{(\mathtt{u})}((A_{\iota}^{(\bar{\tau})})^{[\vec{\lambda}/\vec{\rho}]};\mathtt{R})\equiv C$
by $\lambda_{k}\in(\mathtt{R}(\mS_{k}))_{0}\cap\tau_{k}$.

Conversely let $(A_{\iota}^{(\bar{\tau})})^{[\vec{\lambda}/\vec{\rho}]}\equiv 
C^{[\lambda_{k}/\mS_{k}]}\not\equiv C$ and
$\mathrm{inv}^{(\mathtt{u})}((A_{\iota}^{(\bar{\tau})})^{[\vec{\lambda}/\vec{\rho}]};\mathtt{R})\equiv C$
with $\lambda_{k}<\tau_{k}$.
Then $A_{\iota}\equiv C^{[\rho_{k}/\mS_{k}]}\not\equiv C$.
On the other hand we have $\lambda_{k}<\tau_{k}$ by
$(\mathtt{Q}(\mS_{k})_{0}\cap\rho_{k}\subset\lambda_{k}$.
We obtain $\mathrm{inv}^{(\mathtt{u})}(A_{\iota}^{(\bar{\tau})};\mathtt{Q})\equiv C$.
(\ref{eq:clm_Case3_Recapping_I_1}) is shown, and
this yields
$\iota\in[\bar{\tau}\mathtt{R}]J$
and
$\mathsf{k}((A_{\iota}^{(\bar{\tau})})^{[\vec{\lambda}/\vec{\rho}]};\mathtt{R})\subset
\mathcal{H}_{\gamma}[\Theta(\mathtt{R})]\cap M_{\mathtt{R}}$.

Let $\nu\in[\bar{\tau}\mathtt{R}]I$. 
We see as in (\ref{eq:clm_Case3_Recapping_I_1}) that 
$\mathrm{inv}^{(\mathtt{u})}(B_{\nu}^{(\bar{\tau})};\mathtt{Q})\equiv
\mathrm{inv}^{(\mathtt{u})}((B_{\nu}^{(\bar{\tau})})^{[\vec{\lambda}/\vec{\rho}]};\mathtt{R})$,
$\nu\in [\bar{\tau}\mathtt{Q}]I$, 
$\mathsf{k}(B_{\nu}^{(\bar{\tau})};\mathtt{Q})=\mathsf{k}((B_{\nu}^{(\bar{\tau})})^{[\vec{\lambda}/\vec{\rho}]};\mathtt{R})\subset M_{\kappa_{k}}$, and
$\{\kappa_{k},\theta_{k}\}\subset H_{\rho_{k}}(g_{1},\mathtt{Q},\Theta\cup\mathsf{k}(B_{\nu}^{(\bar{\tau})};\mathtt{Q}))$.

IH with (\ref{eq:clm_Case3_Recapping_I_1}) followed by a $(\bigvee)$ yields (\ref{eq:1.1.f_I_1}).
\\
\textbf{Case 4}.
Fourth consider the case when the last inference is a $(cut)$ with a cut formula
$C^{(\bar{\rho})}$ with $\bar{\rho}=\bar{\rho}_{\mathtt{Q}}$ and $\mathrm{rk}(C)<c\leq b$:
Let $C\simeq\bigvee(C_{\iota})_{\iota\in J}$. We have an ordinal $a_{0}<a$ such that
$
(\mathcal{H}_{\gamma},\Theta, \mathtt{Q}
)
\vdash^{a_{0}}_{c,\vec{d},\xi,\Lambda}
C^{(\bar{\rho})},\Pi,\Phi,\Gamma
$, and
$
(\mathcal{H}_{\gamma},\Theta_{\mathtt{Q}_{\iota}}\cup\mathsf{k}(C_{\iota}), \mathtt{Q}_{\iota}
)
\vdash^{a_{0}}_{c,\vec{d},\xi,\Lambda}
(\{\lnot C_{\iota}^{(\bar{\rho})}\}\cup\Pi\cup\Phi\cup\Gamma)(\mathtt{Q}_{\iota})
$ holds for each $\iota\in[\bar{\rho}\mathtt{Q}]J$ and
each good extension $\mathtt{Q}_{\iota}=\mathtt{Q}[\mathrm{srk}(C_{\iota})]$.
First let $C\in\mathcal{L}_{i_{m}}$. We have $\mathrm{rk}(C)<b$.
IH yields
$
(\mathcal{H}_{\gamma},\Theta, \mathtt{R}
)
\vdash^{\varphi_{d_{i_{m}}}(a_{0})}_{c,\vec{d},\xi,\Lambda}
C^{(\bar{\rho}_{\kappa})},\Pi^{[\vec{\lambda}/\vec{\rho}]},\Phi^{(\vec{\theta}/\vec{\rho})},\Gamma^{(\vec{\kappa}/\vec{\rho})}
$, and
$
(\mathcal{H}_{\gamma},\Theta_{\mathtt{R}_{\iota}}\cup\mathsf{k}(C_{\iota}), \mathtt{R}_{\iota}
)
\vdash^{\varphi_{d_{i_{m}}}(a_{0})}_{c,\vec{d},\xi,\Lambda}
(\{\lnot C_{\iota}^{(\bar{\rho}_{\kappa})}\}\cup\Pi^{[\vec{\lambda}/\vec{\rho}]}\cup
\Phi^{(\vec{\theta}/\vec{\rho})}\cup \Gamma^{(\vec{\kappa}/\vec{\rho})})(\mathtt{R}_{\iota})
$, where $\bar{\rho}_{\kappa}=\bar{\rho}[\vec{\kappa}/\vec{\rho}]$.
Otherwise we obtain by IH,
$
(\mathcal{H}_{\gamma},\Theta, \mathtt{R}
)
\vdash^{\varphi_{d_{i_{m}}}(a_{0})}_{c,\vec{d},\xi,\Lambda}
C^{(\bar{\rho}_{\theta})},\Pi^{[\vec{\lambda}/\vec{\rho}]},\Phi^{(\vec{\theta}/\vec{\rho})},\Gamma^{(\vec{\kappa}/\vec{\rho})}
$, and
$
(\mathcal{H}_{\gamma},\Theta_{\mathtt{R}_{\iota}}\cup\mathsf{k}(C_{\iota}), \mathtt{R}_{\iota}
)
\vdash^{\varphi_{d_{i_{m}}}(a_{0})}_{c,\vec{d},\xi,\Lambda}
(\{\lnot C_{\iota}^{(\bar{\rho}_{\theta})}\}\cup\Pi^{[\vec{\lambda}/\vec{\rho}]}\cup
\Phi^{(\vec{\theta}/\vec{\rho})}\cup \Gamma^{(\vec{\kappa}/\vec{\rho})})(\mathtt{R}_{\iota})
$, where $\bar{\rho}_{\theta}=\bar{\rho}[\vec{\theta}/\vec{\rho}]$.
In each case we obtain (\ref{eq:1.1.f_I_1}) by a $(cut)$.

Other cases are seen from IH. 
\eprf

\blem\label{lem:exist_H}
Let $\mathtt{Q}$ be a finite family with thresholds $\gamma^{\mathtt{Q}}_{\cdot}$,
$(\rho,g)\in\mathtt{Q}(\mS)$ with
$\rho=\rho_{\mathtt{Q}(\mS)}$,
$\gamma_{1}=\gamma^{\mathtt{Q}}_{\mS}+\mS$,
$\Theta_{\mS}=\Theta\cup(\mathtt{Q}(\mS))_{0}$ and $\mathtt{R}=\mathtt{Q}\restrict\mS$.

Let 
$f:\mI_{N}\to\Gamma(\mI_{N})$ be a special finite function with base $\mI_{N}$
 for which one of the following holds:

 \benu
 
  \item\label{lem:exist_H1}
  $g=m(\rho)$,
  $\mathrm{fld}(f)\cup\{\mathtt{p}_{0}(\rho)\}\subset\mathcal{H}_{\gamma}[\Theta(\mathtt{Q})]$ for a $\gamma\leq\gamma^{\mathtt{Q}}_{\mS}$, and if $f\neq\emptyset$, then
  there exists $x\in\mathrm{supp}(g)$ and $y<x$ such that
  $(y,x)\cap\mathrm{supp}(g)=(y,x)\cap\mathrm{supp}(f)=\emptyset$, $f_{y}\leq g_{y}$,
$f(y)<g(y)+\tilde{\theta}_{x-y}(g(x))\cdot\omega$ and $f<_{\mI_{N}}^{x}g(x)$, 
cf.\,Definition \ref{df:notationsystem}.\ref{df:notationsystem.6}.

 \item\label{lem:exist_H15}
  $\mathrm{fld}(f)\cup\{\mathtt{p}_{0}(\rho)\}\subset\mathcal{H}_{\gamma}[\Theta(\mathtt{Q})]$ for a $\gamma\leq\gamma^{\mathtt{Q}}_{\mS}$, and
$f_{x}\leq g_{x}$ and $f<_{\mI_{N}}^{x}g^{\prime}(x)$ for an $x\in\mathrm{supp}(g)$
  
   \item\label{lem:exist_H2}
   $f=h^{b}(g;a)$ for some $\{b,a\}\cup \mathrm{fld}(g)\cup\{\mathtt{p}_{0}(\rho)\}\subset\mathcal{H}_{\gamma}[\Theta(\mathtt{Q})]$ such that
   $b<s(g)$ and $a<\mI_{N}$.
   \eenu
 Then there exists a $\sigma$ such that $\sigma\in H_{\rho}(f,\mathtt{Q},\Theta)$,
 $m(\sigma)=f$, $\mathtt{p}_{0}(\sigma)=\mathtt{p}_{0}(\rho)$ and 
 $\sigma\in
 \mathcal{H}_{\gamma_{1}}[\Theta_{\mS}(\mathtt{R})]$.

\elem
\bprf
\ref{lem:exist_H}.\ref{lem:exist_H1}.
Let $X=\Theta(\mathtt{Q})\cup((\mathtt{Q}(\mS))_{0}\cap\rho)
\cup\mathrm{fld}(f)\cup\{\mathtt{p}_{0}(\rho)\}$,
$E_{\mS}(X)=\bigcup\{E_{\mS}(\delta):\delta\in X\}$ and
$\alpha=\max(\{1\}\cup E_{\mS}(X))$.
Let $\sigma=\psi_{\rho}^{f}(\nu+\alpha)$ for $\nu=b(\rho)$.
We have $\mathtt{p}_{0}(\rho)=\mathtt{p}_{0}(\sigma)$, $m(\sigma)=f$, and 
$\gamma_{\mS}^{\mathtt{Q}}<\nu<\gamma_{\mS}^{\mathtt{Q}}+\mathbb{S}=\gamma_{1}$
by (\ref{eq:controlder_cap4}) in Definition \ref{df:ffthreshold}.\ref{df:caphat.42}.
We obtain
$\nu+\alpha<\gamma_{\mS}^{\mathtt{Q}}+\mathbb{S}=\gamma_{1}$ by $\alpha<\mathbb{S}$.

By $\mathrm{fld}(f)\cup\{\mathtt{p}_{0}(\rho)\}\subset\mathcal{H}_{\gamma}[\Theta(\mathtt{Q})]$ with
$\gamma\leq \gamma^{\mathtt{Q}}_{\mS}<\mathtt{p}_{0}(\rho)$ and
$\Theta(\mathtt{Q})\subset M_{\rho}$
we obtain $\mathrm{fld}(f)\cup\{\mathtt{p}_{0}(\rho)\}\subset M_{\rho}$.
Proposition \ref{prp:EH_Pi11}.\ref{prp:EH_Pi11.8} yields $E_{\mS}(X)\subset\rho$,  and $0<\alpha<\rho$.
Hence $\{\rho,\nu+\alpha\}\subset C_{\nu+\alpha}(\rho)$ by $\rho\in OT(\mI_{N})$, and $\alpha<\sigma$.
On the other hand we have $X\subset M_{\rho}=C_{\mathtt{p}_{0}(\rho)}(\rho)\subset C_{\mathtt{p}_{0}(\rho)}(\mS)$.
Proposition \ref{prp:ESM} yields
$X=\Theta(\mathtt{Q})\cup((\mathtt{Q}(\mS))_{0}\cap\rho)
\cup\mathrm{fld}(f)\cup\{\mathtt{p}_{0}(\rho)\}\subset M_{\sigma}$, cf.\,(\ref{eq:notationsystem.6}) in 
Definition \ref{df:notationsystem}.\ref{df:notationsystem.6}, and
$\mathrm{fld}(f)\subset C_{\nu+\alpha}(\rho)$ by $\mathtt{p}_{0}(\rho)\leq\nu<\nu+\alpha$.
Hence $\sigma\in OT(\mI_{N})$ by Definition \ref{df:notationsystem}.\ref{df:notationsystem.6}, and $\sigma\in H_{\rho}(f,\mathtt{Q},\Theta)$.

We have
$\mathrm{fld}(f)\subset\mathcal{H}_{\gamma}[\Theta(\mathtt{Q})]$ 
with $\gamma<\gamma_{1}$.
We obtain
$E_{\mS}(X)\subset \mathcal{H}_{\gamma_{1}}[\Theta_{\mS}(\mathtt{R})]$ by Proposition \ref{prp:EH_Pi11}.\ref{prp:EH_Pi11.5}.
Let $\{\rho\}\cup\Theta(\mathtt{Q})\subset\Theta_{\mS}(\mathtt{R})\subset C_{\delta}(\beta)$
for $\nu=b(\rho)<\gamma_{1}<\delta$.
Then $\alpha\in C_{\delta}(\beta)$ by $\alpha\in \mathcal{H}_{\gamma_{1}}[\Theta_{\mS}(\mathtt{R})]$.
If $\beta\leq\rho$, then $\nu\in C_{\delta}(\beta)$ by $\rho\in C_{\delta}(\beta)$.
Let $\rho\leq\beta$. Then $\nu\in C_{\delta}(\beta)$ by $\nu<\delta$ and $\nu\in C_{\nu}(\rho)$.
In each case we obtain $\nu\in C_{\delta}(\beta)$, and
$\nu\in\mathcal{H}_{\gamma_{1}}[\Theta_{\mS}(\mathtt{R})]$.
Hence $\sigma\in
\mathcal{H}_{\gamma_{1}}[\Theta_{\mS}(\mathtt{R})]$.

Finally let $\tau=\rho_{\mathtt{Q}(\mU)}$ with $\mU<\mS$.
We obtain $\sigma\in\mathcal{H}_{\gamma_{1}}[\Theta(\mathtt{Q})]\subset M_{\tau}$ by
$\Theta(\mathtt{Q})\subset M_{\tau}$ and
$\gamma_{1}=\gamma^{\mathtt{Q}}_{\mS}+\mS\leq \gamma^{\mathtt{Q}}_{\mU}<\mathtt{p}_{0}(\tau)$.
\\
\ref{lem:exist_H}.\ref{lem:exist_H15}.
Let $k=m(\rho)\geq g$.
If $f_{x}\leq g_{x}$ and $f<_{\mI_{N}}^{x}g^{\prime}(x)$ for an $x\in\mathrm{supp}(g)$, then
$f_{x}\leq k_{x}$ and $f<_{\mI_{N}}^{x}k^{\prime}(x)$ by Proposition \ref{prp:idless}, and we may assume $g=m(\rho)$.
Let
$y=\max(\{0\}\cup\{y\in \mathrm{supp}(g): y<x\})$. 
It is clear that $(y,x)\cap\mathrm{supp}(g)=(y,x)\cap\mathrm{supp}(f)=\emptyset$, $f_{y}\leq g_{y}$,
$f(y)<g(y)+\tilde{\theta}_{x-y}(g(x))\cdot\omega$ and $f<_{\mI_{N}}^{x}g(x)$.
Thus Lemma \ref{lem:exist_H}.\ref{lem:exist_H1} applies.
\\
\ref{lem:exist_H}.\ref{lem:exist_H2}.
Let $\{b=b_{0}<b_{1}<\cdots<b_{n}=s(g)\}$ with $n>0$.
Define recursively
ordinals $\alpha_{i}$ for $i\leq n$ as follows.
$\alpha_{n}=\mI_{N}\cdot\alpha+a$ with $g(s(g))=\mI_{N}(\alpha+1)$.
$\alpha_{i}=g(b_{i})+\tilde{\theta}_{b_{i+1}-b_{i}}(\alpha_{i+1})$.
Then
special functions $h_{i}\, (i\leq n)$ are defined by
$s(h_{i})=b_{i}$, $(h_{i})_{b_{i}}=g_{b_{i}}$,
$h_{n}=g$ and
$h_{i}(b_{i})=\alpha_{i}+\mI_{N}$ for $i<n$.
Then $h_{0}=h^{b}(g;a)=f$.

For $i<n$ we obtain $b_{i+1}\in\mathrm{supp}(h_{i+1})$, 
$(b_{i},b_{i+1})\cap\mathrm{supp}(h_{i})=(b_{i},b_{i+1})\cap\mathrm{supp}(h_{i+1})=\emptyset$,
$(h_{i})_{b_{i}}=(h_{i+1})_{b_{i}}=g_{b_{i}}$, and $h_{i}<_{\mI_{N}}^{b_{i+1}}h_{i+1}$.
By $\alpha_{i+1}<h_{i+1}(b_{i+1})$,
$h_{i}(b_{i})=g(b_{i})+\tilde{\theta}_{b_{i+1}-b_{i}}(\alpha_{i+1})+\mI_{N}<h_{i+1}(b_{i})+\tilde{\theta}_{b_{i+1}-b_{i}}(h_{i+1}(b_{i+1}))$.

We have $\mathrm{fld}(h_{i})\subset\mathcal{H}_{0}[\{b,a,\mI_{N}\}\cup\mathrm{fld}(g)]$.
We see
$\mathrm{fld}(h_{i})\cup\{\mathtt{p}_{0}(\rho)\}\subset\mathcal{H}_{\gamma}[\Theta(\mathtt{Q})]$
from the assumption $\{b,a,\mI_{N},\mathtt{p}_{0}(\rho)\}\cup\mathrm{fld}(g)\subset\mathcal{H}_{\gamma}[\Theta(\mathtt{Q})]$.
By Lemma \ref{lem:exist_H}.\ref{lem:exist_H1} pick ordinals $\{\sigma_{i}\}_{i\leq n}$ recursively as follows.
Let $\sigma_{n}=\rho$, and for $i<n$, $\sigma_{i}$ is an ordinal such that $\mathtt{p}_{0}(\sigma_{i})=\mathtt{p}_{0}(\rho)$,
$m(\sigma_{i})=h_{i}$, $\sigma_{i}\in OT(\mI_{N})$,
$\sigma_{i}\in H_{\sigma_{i+1}}(h_{i},\mathtt{Q},\Theta)$ and 
$\{\sigma_{i},b(\sigma_{i})\}
\subset\mathcal{H}_{\gamma_{1}}[\Theta(\mathtt{Q})]$.
We obtain $\sigma_{0}\in H_{\sigma_{1}}(f,\mathtt{Q},\Theta)$, $\sigma_{1}\leq\rho$ and 
$\mathtt{p}_{0}(\rho)=\mathtt{p}_{0}(\sigma_{1})=\mathtt{p}_{0}(\sigma_{0})$.
Therefore $\sigma_{0}\in H_{\rho}(f,\mathtt{Q},\Theta)$ is a desired one.
\eprf

 \bdf\label{df:on_f}{\rm (Cf.\,Definition \ref{df:hstepdown}.)\\
Let $\Lambda\leq\mI_{N}$ be a strongly critical number, and
$f:\mI_{N}\to\Gamma(\mI_{N})$ a special finite function such that
$SC_{\mI_{N}}(f)=\bigcup\{\{c\}\cup SC_{\mI_{N}}(f(c)): c\in\mathrm{supp}(f)\}\subset\Lambda$.
Then define
$on_{\Lambda}(f)=\alpha[\Lambda:\mI_{N}]$ for $\alpha=(h^{0}(f;\mI_{N}))^{\prime}(0)$ if $s(f)>0$,
where $h^{0}(f;\mI_{N})$ is a finite function in Definition \ref{df:hstepdown}.\ref{df:hstepdown.1}, and
$\alpha[\Lambda:\mI_{N}]$
denotes an ordinal obtained from $\alpha$ by changing the base $\mI_{N}$ to $\Lambda$, 
cf.\,Definition \ref{df:Lam2}.\ref{df:LamI}.
Otherwise let $on_{\Lambda}(f)=on_{\Lambda}(\emptyset)=0$.

}
\edf

Note that $on_{\Lambda}(f)\in\mathcal{H}_{\gamma}[\Theta]$ if $\{\Lambda\}\cup\mathrm{fld}(f)\subset\mathcal{H}_{\gamma}[\Theta]$, 
$on_{\Lambda}(f)<\Gamma(\Lambda)$ and $\Lambda<on_{\Lambda}(f)$ unless  $on_{\Lambda}(f)=0$.

\bprp\label{prp:hstepdown_on}
Let $\Lambda<\mI_{N}$ be a strongly critical number, and
$f,g:\mI_{N}\to\Gamma(\mI_{N})$ special finite functions such that $SC_{\mI_{N}}(f)\cup SC_{\mI_{N}}(g)\subset\Lambda$.
Assume
$f_{d}\leq g_{d}$ and $f<_{\mI_{N}}^{d}g^{\prime}(d)$ for a $d\in\mathrm{supp}(g)$. 
\benu

\item\label{prp:hstepdown_on_1}
$on_{\Lambda}(f)<on_{\Lambda}(g)$.

\item\label{prp:hstepdown_on_2}
Let $f_{1}=h^{b}(f;a)$ with $b<s(f)$ and $a<\Lambda$. 
Then $on_{\Lambda}(f_{1})<on_{\Lambda}(f)$.
\eenu
\eprp
\bprf
By Proposition \ref{prp:thtcollapse} it suffices to show
$(h^{0}(f;\mI_{N}))(0)<(h^{0}(g;\mI_{N}))(0)$ in Proposition
\ref{prp:hstepdown_on}.\ref{prp:hstepdown_on_1},
and $(h^{0}(f_{1};\mI_{N}))(0)<(h^{0}(f;\mI_{N}))(0)$
in Proposition \ref{prp:hstepdown_on}.\ref{prp:hstepdown_on_2}.
\\
\ref{prp:hstepdown_on}.\ref{prp:hstepdown_on_1}.
By Proposition \ref{prp:hstepdown}.\ref{prp:hstepdown.3} we obtain
$(h^{d}(f;\mI_{N}))(d)<g(d)\leq(h^{d}(g;\mI_{N}))(d)$, which yields $(h^{0}(f;\mI_{N}))(0)<(h^{0}(g;\mI_{N}))(0)$.
\\
\ref{prp:hstepdown_on}.\ref{prp:hstepdown_on_2}.
$(h^{0}(f_{1};\mI_{N}))(0)<(h^{0}(f;\mI_{N}))(0)$ is seen from $a<\Lambda$.
\eprf

\bdf\label{df:main_lemmas}
{\rm
Let $\mathtt{Q}$ be a finite family with thresholds $\gamma^{\mathtt{Q}}_{\cdot},\delta^{\mathtt{Q}}_{\cdot}$, $\Phi$ a sequent and $b,\alpha_{0}$ ordinals.
\benu
\item\label{df:main_lemmas.1}
Given $1\leq i\leq N$, an ordinal $d_{i}$, and strongly critical numbers $\Lambda<\Lambda^{+}$,
let us introduce ordinals $\mS,\kappa$, a finite family $\mathtt{Q}^{(\kappa/\rho)}$, 
a sequent $\Phi^{(\kappa/\rho)}$ and an ordinal 
$\alpha(\mathtt{Q})=\alpha(\mathtt{Q};\mS)$ as follows.

If $\mathrm{dom}(\mathtt{Q})\cap SSt_{i}=\emptyset$, then let 
$\alpha(\mathtt{Q})=\Lambda^{+}\cdot \omega^{\Lambda_{0}+1}$, $\mS=\kappa=\rho=0$,
$\mathtt{Q}^{(\kappa/\rho)}=\mathtt{Q}$ 
and $\Phi^{(\kappa/\rho)}=\Phi$.
In this case let $on_{\Lambda}(g)=0$ with $g=\emptyset$, and $\delta^{\mathtt{Q}}_{\mS}=\omega^{\Lambda_{0}+1}$.

Let $\mathrm{dom}(\mathtt{Q})\cap SSt_{i}\neq\emptyset$.
Then let $\mS=\max(\mathrm{dom}(\mathtt{Q})\cap SSt_{i})$ and
$\rho=\rho_{\mathtt{Q}(\mS)}$ with $(\rho,g)\in\mathtt{Q}(\mS)$.
If $s(g)\leq b$, then let $\kappa=\rho$.
Otherwise pick an ordinal $\kappa$ such that
$\kappa\in H_{\rho}(g_{1},\mathtt{Q},\Theta)$ with
$m(\kappa)\geq g_{1}=h^{b}(g;\varphi_{d_{i}}(\alpha_{0}))$.
Then let
$\mathtt{Q}^{(\kappa/\rho)}=(\mathtt{Q}\cup\{(\mS,(\kappa,g_{1})\})\setminus\{(\mS,(\rho,g))\}$
and
$\alpha(\mathtt{Q})=\Lambda^{+}\cdot \delta^{\mathtt{Q}}_{\mS}+on_{\Lambda}(g)$.
Let $\Phi=\Pi\cup\Gamma$ with $\Pi(\rho)=\emptyset$ and $\Gamma(\rho)=\Gamma$ for
$\vec{\rho}=(\rho)$.
$\Phi^{(\kappa/\rho)}$ denotes a sequent
$\Pi^{[\kappa/\rho]}\cup\Gamma^{(\kappa/\rho)}$
 in Definition \ref{df:recapping_I_1} for $\vec{\kappa}=\vec{\lambda}=(\kappa)$ and 
$\{\mS_{1}<\cdots<\mS_{m}\}=\{\mS\}$ with $m=1$.

\item\label{df:main_lemmas.2}
Given ordinals $\vec{d}=(d_{1},\ldots,d_{N})$, and strongly critical numbers 
$\Lambda^{(0)}<\Lambda^{(1)}<\Lambda^{(2)}<\cdots<\Lambda^{(3N-3)}<\Lambda^{(3N-2)}<\Lambda^{(3N-1)}$,
let us introduce ordinals $\mS_{i},\kappa_{i},\rho_{i}$, a finite family $\mathtt{Q}_{i}$,
a sequent $\Phi_{i}$
and ordinals
$\alpha_{i}(\mathtt{Q})$ and $\alpha_{i}$ recursively for $i\leq N$ as follows.
Let $\Phi_{0}=\Phi$ and $\mathtt{Q}_{0}=\mathtt{Q}$.

If $\mathrm{dom}(\mathtt{Q})\cap SSt_{i}=\emptyset$, then let 
$\alpha_{i}(\mathtt{Q})=\Lambda^{(3i-2)}\cdot \omega^{\Lambda_{0}+1}$, $\mS_{i}=\kappa_{i}=\rho_{i}=0$,
$\mathtt{Q}_{i}=\mathtt{Q}_{i-1}$ 
and $\Phi_{i}=\Phi_{i-1}$.

Let $\mathrm{dom}(\mathtt{Q})\cap SSt_{i}\neq\emptyset$.
Then let $\mS_{i}=\max(\mathrm{dom}(\mathtt{Q})\cap SSt_{i})$ and
$\rho_{i}=\rho_{\mathtt{Q}(\mS_{i})}$ with $(\rho_{i},g_{i})\in\mathtt{Q}(\mS_{i})$.
If $s(g_{i})\leq b$, then let $\kappa_{i}=\rho_{i}$.
Otherwise pick an ordinal $\kappa_{i}$ such that
$\kappa_{i}\in H_{\rho_{i}}(k_{i},\mathtt{Q},\Theta)$ with
$m(\kappa_{i})\geq k_{i}=h^{b}(g_{i};\varphi_{d_{i}}(\alpha_{i-1}))$.
Then let
$\mathtt{Q}_{i}=\mathtt{Q}_{i-1}^{(\kappa_{i}/\rho_{i})}=(\mathtt{Q}_{i-1}\cup\{(\mS_{i},(\kappa_{i},k_{i})\})\setminus\{(\rho_{i},g_{i})\}$,
$\alpha_{i}(\mathtt{Q})=\Lambda^{(3i-2)}\cdot \delta^{\mathtt{Q}}_{\mS_{i}}+on_{\Lambda^{(3i-3)}}(g_{i})$ and
$\Phi_{i}=\Phi_{i-1}^{(\kappa_{i}/\rho_{i})}$.

In each case let
$\alpha_{i}=\varphi_{\alpha_{i}(\mathtt{Q})}(\alpha_{i-1})$.
Let $\Phi^{(\vec{\kappa}/\vec{\rho})}=\Phi_{N}$.

\item\label{df:main_lemmas.3}
For an ordinal $d$, let $\vec{d}=(d,\ldots,d)$ denote the vector whose $i$-th entry is $d_{i}=d$
for every $i$.
\eenu
}
\edf

\blem\label{mlem:singlemainl_I_1}
Let $\mathtt{Q}$ be a finite family with thresholds $\gamma^{\mathtt{Q}}_{\cdot},\delta^{\mathtt{Q}}_{\cdot}$
such that 
$\mT\in \mathrm{dom}(\mathtt{Q})$, and
$\Lambda^{++}>\Lambda^{+}>\Lambda>\Lambda_{0}$ strongly critical numbers such that
$\{\Lambda^{++},\Lambda^{+},\Lambda\}\subset\mathcal{H}_{\gamma}[\Theta(\mathtt{Q})]$.
$b$ denotes an ordinal such that
$b=\mT+1$ if $\mT\in SSt\cup\{0\}$, 
and $b=\mT$ else.
Let $1\leq i\leq N$ be a number, and
$\vec{d}=(d_{1},\ldots,d_{N})=(b,\ldots b,d_{i},d_{i+1},\ldots,d_{N})$ be such that
$b\leq d_{i}\leq\cdots\leq d_{N}\leq\mT^{\dagger}$ and $\forall j<i(d_{j}=b)$.

Let
$
(\mathcal{H}_{\gamma},\Theta, \mathtt{Q}
)
\vdash^{a}_{b,\vec{d},\mT^{\dagger},\Lambda}
\Phi
$ for 
$a<\Lambda$.

Let $\mS,\kappa,\rho$ be ordinals, $\mathtt{Q}^{(\kappa/\rho)}$ a finite family,
$\Phi^{(\kappa/\rho)}$ a sequent and $\alpha(\mathtt{Q})=\alpha(\mathtt{Q};\mS)$ an ordinal 
defined from ordinals $b$, $\alpha_{0}=a$,
the number $i$, the ordinal $d_{i}$, and strongly critical numbers $\Lambda<\Lambda^{+}$
in Definition \ref{df:main_lemmas}.\ref{df:main_lemmas.1}.

Let $\vec{e}=(b,\ldots b,b,d_{i+1},\ldots,d_{N})$ be obtained from $\vec{d}$ by replacing the $i$-th entry
$d_{i}$ by $b$.
Then 
\begin{equation}\label{eq:singlemainl_1_I_1}
(\mathcal{H}_{\gamma},\Theta, \mathtt{Q}^{(\kappa/\rho)}
)
\vdash^{\hat{a}}_{b,\vec{e},\mT^{\dagger},\Lambda^{++}}
\Phi^{(\kappa/\rho)}
\end{equation}
holds for $\hat{a}=\varphi_{\alpha(\mathtt{Q})}(a)$.

\elem
\bprf
By main induction on $\alpha(\mathtt{Q})$ with subsidiary induction on $a$.

We have $\mU<b$ for every $\mU\in\mathrm{dom}(\mathtt{Q})\cap SSt$, and
$SC_{\mI_{N}}(\mathrm{fld}(g))\subset\Lambda$
by (\ref{eq:fld}) for finite functions $g=(\rho,g)_{1}$ with base $\mI_{N}$.

On the other hand we have 
$\Theta(\mathtt{Q})=\Theta(\mathtt{Q}^{(\kappa/\rho)})$,
$\{a,\Lambda^{+},\Lambda,\Lambda_{0},\delta^{\mathtt{Q}}_{\mS}\}\cup\mathrm{fld}(g)\subset
\mathcal{H}_{\gamma}[\Theta(\mathtt{Q})]$
by (\ref{eq:controlder_cap22}), and $\hat{a}\in\mathcal{H}_{\gamma}[\Theta(\mathtt{Q})]$.
(\ref{eq:controlder_cap22}) in (\ref{eq:singlemainl_1_I_1}) follows.
Furthermore $a<\Lambda$ and $SC_{\mI_{N}}(\mathrm{fld}(g))\subset\Lambda$ by (\ref{eq:fld}).
We obtain $on_{\Lambda}(g)<\Gamma(\Lambda)\leq\Lambda^{+}$ and
$\delta^{\mathtt{Q}}_{\mS}\leq\omega^{\Lambda_{0}+1}<\Lambda<\Lambda^{+}\leq\alpha(\mathtt{Q})<\Gamma(\Lambda^{+})=\Lambda^{++}$.
Hence $\hat{a}=\varphi_{\alpha(\mathtt{Q})}(a)<\Lambda^{++}$.

We have
$a<\Lambda$, $b^{\dagger}\geq\mT^{\dagger}$ and
$d_{j}=b\leq d_{i}$ for any $j<i$.
\\
\textbf{Case 1}. 
$s(g)> b$: Then  $\mathrm{dom}(\mathtt{Q})\cap SSt_{i}\neq\emptyset$.
Let $\mS=\max(\mathrm{dom}(\mathtt{Q})\cap SSt_{i})$,
$(\rho,g)\in\mathtt{Q}(\mS)$ with $\rho=\rho_{\mathtt{Q}(\mS)}$,
$g_{1}=h^{b}(g;\varphi_{d_{i}}(a))$ and
$\kappa\in H_{\rho}(g_{1},\mathtt{Q},\Theta)$.
Recapping \ref{lem:recapping} with $\vec{\kappa}=\vec{\theta}=\vec{\lambda}=(\kappa)$
yields
$(\mathcal{H}_{\gamma},\Theta, \mathtt{Q}^{(\kappa/\rho)}
)
\vdash^{\varphi_{d_{i}}(a)}_{b,\vec{d},\mT^{\dagger},\Lambda}
\Phi^{(\kappa/\rho)}$.
On the other hand we have $\kappa=\rho_{\mathtt{Q}^{(\kappa/\rho)}(\mS)}$,
$(\kappa,g_{1})\in\mathtt{Q}^{(\kappa/\rho)}(\mS)$ and 
$SC_{\mI_{N}}(\mathrm{fld}(g_{1}))\subset\Lambda$ by $a<\Lambda$.
Proposition \ref{prp:hstepdown_on}.\ref{prp:hstepdown_on_2} yields
$on_{\Lambda}(g_{1})<on_{\Lambda}(g)$. 
We obtain by $d_{i}\leq\mT^{\dagger}<\Lambda_{0}<\Lambda^{+}$,
$\varphi_{\Lambda^{+}\cdot\delta^{\mathtt{Q}}_{\mS}+on_{\Lambda}(g_{1})}(\varphi_{d_{i}}(a))<\varphi_{\Lambda^{+}\cdot\delta^{\mathtt{Q}}_{\mS}+on_{\Lambda}(g)}(a)$.
MIH yields (\ref{eq:singlemainl_1_I_1}).

In what follows assume either $\mathrm{dom}(\mathtt{Q})\cap SSt_{i}=\emptyset$
or $s(g)\leq b \spand \kappa=\rho$.
Then $\mathtt{Q}^{(\kappa/\rho)}=\mathtt{Q}$ and $\Phi^{(\kappa/\rho)}=\Phi$.
\\
\textbf{Case 2}. The last inference in $
(\mathcal{H}_{\gamma},\Theta, \mathtt{Q}
)
\vdash^{a}_{b,\vec{d},\mT^{\dagger},\Lambda}
\Phi
$ is a $(\mathrm{dom})$: We have an ordinal $a_{0}<a$, and a stable ordinal $\mU\in\mathcal{H}_{\gamma}[\Theta(\mathtt{Q})]\cap\mT^{\dagger}$
such that $\mU\not\in\mathrm{dom}(\mathtt{Q})$.
For each good extension $\mathtt{P}=\mathtt{Q}[\mU]$ of $\mathtt{Q}$ for $\gamma,\Theta$ 
by ordinals $\gamma^{\mathtt{P}}_{\mU}, \delta^{\mathtt{P}}_{\mU},\rho_{\mathtt{P}(\mU)}$, we have
$(\mathcal{H}_{\gamma},\Theta_{\mathtt{P}},\mathtt{P})\vdash^{a_{0}}_{b,\vec{d},\mT^{\dagger},\Lambda}\Gamma(\mathtt{P})$.

Let $\mS_{0}\in\{0\}\cup SSt_{i}$ be the ordinal for the lower sequent, and 
$\mS_{1}\in\{0\}\cup SSt_{i}$ for the upper sequents.
If $\mathrm{dom}(\mathtt{P})\cap SSt_{i}=\emptyset$, then let $\mS_{0}=\mS_{1}=0$.
If $\mV<\mU\in SSt_{i}$ for every $\mV\in\mathrm{dom}(\mathtt{Q})\cap SSt_{i}$, then 
$\mS_{1}=\mU>\mS_{0}$.
Otherwise
$\mS_{0}=\max(\mathrm{dom}(\mathtt{Q})\cap SSt_{i})=\max(\mathrm{dom}(\mathtt{P})\cap SSt_{i})=\mS_{1}$.

If $\mS_{1}=\mS_{0}$, then $\alpha(\mathtt{P})=\alpha(\mathtt{Q})$, and
SIH followed by a $(\mathrm{dom})$ yields (\ref{eq:singlemainl_1_I_1}).
Otherwise $\mS_{0}<\mS_{1}=\mU\in SSt_{i}\cap\mT^{\dagger}$ holds.
We show $\alpha(\mathtt{P};\mU)<\alpha(\mathtt{Q};\mS_{0})$.
We obtain 
$\delta^{\mathtt{Q}}_{0}=\omega^{\Lambda_{0}+1}>\delta^{\mathtt{P}}_{\mU}$ if $\mS_{0}=0$, and
$\delta^{\mathtt{Q}}_{\mS_{0}}=\delta^{\mathtt{P}}_{\mS_{0}}>\delta^{\mathtt{P}}_{\mU}$
for $\{\mS_{0}<\mU\}\subset\mathrm{dom}(\mathtt{P})$.
On the other, $\Lambda^{+}>on_{\Lambda}(k)$ for $\mathtt{P}(\mU)=\{(\rho_{\mathtt{P}(\mU)},k)\}$.
Hence
$\alpha(\mathtt{Q};\mS_{0})=\Lambda^{+}\cdot\delta^{\mathtt{Q}}_{\mS_{0}}+on_{\Lambda}(g)>
\Lambda^{+}\cdot\delta^{\mathtt{P}}_{\mU}+on_{\Lambda}(k)=\alpha(\mathtt{P};\mU)$.
MIH followed by a $(\mathrm{dom})$ yields (\ref{eq:singlemainl_1_I_1}).
\\
\textbf{Case 3}. The last inference is a $(j_{0} {\rm -rfl}_{\mU}(k,x,f,\Delta))$:
We have $\mU\in\mathrm{dom}(\mathtt{Q})\cap SSt_{j_{0}}$.
Let $\bar{\rho}=\bar{\rho}_{\mathtt{Q}}$.
We have  an ordinal $a_{0}<a$ and
a set
$\Delta\subset\mathcal{L}_{j_{0}}$ such that $\mU\leq\mathrm{rk}(\Delta)<\min\{s(k),d_{j_{0}}\}$.

Let $\mV$ be a successor $j$-stable ordinal 
such that $j\geq j_{0}$, $\mU\leq\mV\leq\mathrm{rk}(\Delta)$ and $\min\{s(k), s(f)\}\geq\mV+1$,
and $\mathtt{P}=\mathtt{Q}[\mV]$ a good extension such that
$\mathtt{P}(\mV)=\{(\rho_{\mV},k)\}$ with $\rho_{\mV}=\rho_{\mathtt{P}(\mV)}$ and
$(\rho_{\mathtt{Q}(\mU)},k)\in\mathtt{Q}(\mU)$.

Then
$
(\mathcal{H}_{\gamma},\Theta_{\mathtt{P}_{\iota}}\cup\mathsf{k}(\delta_{\iota}^{(\bar{\rho}_{\mathtt{P}})};\mathtt{P}),
\mathtt{P}_{\iota}
)
\vdash^{a_{0}}_{b,\vec{d},\mT^{\dagger},\Lambda} (\Phi(\mathtt{P})\cup\{\lnot\delta_{\iota}^{(\bar{\rho}_{\mathtt{P}})}\})(\mathtt{P}_{\iota})
$ holds
for each $\delta\in\Delta$, $\iota\in[\bar{\rho}_{\mathtt{P}}\mathtt{P}]J$ and
each good extension $\mathtt{P}_{\iota}=\mathtt{P}[\mathrm{srk}(\delta_{\iota})]$,
where $\delta\simeq\bigvee(\delta_{\iota})_{\iota\in J}$.

On the other hand we have
$
(\mathcal{H}_{\gamma},\Theta_{\mathtt{P}},
\mathtt{P}^{\sigma}
)
\vdash^{a_{0}}_{b,\vec{d},\mT^{\dagger},\Lambda}
\Delta^{(\bar{\rho}_{\sigma})}, \Phi(\mathtt{P})
$  for $\sigma\in H_{\rho_{\mV}}(f,\mathtt{Q},\Theta)$, 
$\bar{\rho}_{\sigma}=\bar{\rho}_{\mathtt{P}}[\sigma/\rho_{\mV}]$ and
$\mathtt{P}^{\sigma}=\mathtt{P}\cup\{(\mV,(\sigma,f))\}$.

$f$ is a finite function such that $x\in\mathrm{supp}(k)$,
$f_{x}\leq k_{x}$,  $f<_{\mI_{N}}^{x}k^{\prime}(x)$, $SC_{\mI_{N}}(\mathrm{fld}(f))\subset \Lambda$ and
$\mathrm{fld}(f)\subset\mathcal{H}_{\gamma}[\Theta(\mathtt{Q})]$.

Suppose $\mU\leq\mV=\mS$ for a $\mV$. Then $j_{0}\leq j=i$.
If $\mU<\mS$,
then $\mV\leq\mathrm{rk}(\Delta)<\mS<b$ by (r1).
Hence $\mV\neq\mS$ if $\mU\neq\mS$.

Let $\mU\in SSt_{i}$. Then $\mathrm{dom}(\mathtt{Q})\cap SSt_{i}\neq\emptyset$ and
$\mU\leq\mS=\max(\mathrm{dom}(\mathtt{Q})\cap SSt_{i})$.
If $\mU<\mS$, then $\mathrm{rk}(\Delta)<\mS<b$.
If $\mU=\mS$, then $k=(\bar{\rho}_{\mathtt{Q}}(\mS))_{1}=g$ and
$\mathrm{rk}(\Delta)<s(g)\leq b$.
In each case we obtain $\mathrm{rk}(\Delta)<b$.

Let $\widehat{a_{0}}=\varphi_{\alpha(\mathtt{Q})}(a_{0})$.
\\
\textbf{Case 3.1}. $\mU\neq\mS$:
We have $\rho_{\mathtt{P}_{\iota}(\mS)}=\rho_{\mathtt{Q}(\mS)}=\rho$
 with $(\rho,g)\in\mathtt{P}_{\iota}(\mS)$, and 
 $\delta^{\mathtt{P}_{\iota}}_{\mS}=\delta^{\mathtt{Q}}_{\mS}$.
 Hence $\alpha(\mathtt{P}_{\iota})=\alpha(\mathtt{P}_{\iota};\mS)=\alpha(\mathtt{Q};\mS)=\alpha(\mathtt{Q})$.
 SIH yields
\beqn\label{eq:singlemainl_I_1_1}
(\mathcal{H}_{\gamma},\Theta_{\mathtt{P}_{\iota}}\cup\mathsf{k}(\delta_{\iota}^{(\bar{\rho}_{\mathtt{P}})};\mathtt{P}),
\mathtt{P}_{\iota}
)
\vdash^{\widehat{a_{0}}}_{b,\vec{e},\mT^{\dagger},\Lambda^{++}} (\Phi(\mathtt{P})\cup\{\lnot\delta_{\iota}^{(\bar{\rho}_{\mathtt{P}})}\})(\mathtt{P}_{\iota})
\eeqn
and
\beqn\label{eq:singlemainl_I_1_2}
(\mathcal{H}_{\gamma},\Theta_{\mathtt{P}},
\mathtt{P}^{\sigma}
)
\vdash^{\widehat{a_{0}}}_{b,\vec{e},\mT^{\dagger},\Lambda^{++}}
\Delta^{(\bar{\rho}_{\kappa})}, \Phi(\mathtt{P})
\eeqn
for each $\sigma\in H_{\rho_{\mV}}(f,\mathtt{Q},\Theta)$.

We have $\mathrm{rk}(\Delta)<b$ if $j_{0}\leq i$, and $\mathrm{rk}(\Delta)<d_{j_{0}}$ if $j_{0}>i$.
A $(j_{0} {\rm -rfl}_{\mU}(k,x,f,\Delta))$ with (\ref{eq:singlemainl_I_1_1}) and (\ref{eq:singlemainl_I_1_2})
yields (\ref{eq:singlemainl_1_I_1}).
\\
\textbf{Case 3.2}. $\mU=\mS$:
We obtain $\mS\in\mathrm{dom}(\mathtt{Q})\cap SSt_{i}$,
$\Delta\subset\mathcal{L}_{i}$,
$\rho_{\mV}=\rho_{\mathtt{P}(\mV)}=\rho_{\mathtt{Q}(\mS)}=\rho$,
$k=(\bar{\rho}_{\mathtt{Q}}(\mS))_{1}=g$ and
$\mathrm{rk}(\Delta)<b$.
SIH yields (\ref{eq:singlemainl_I_1_1}).

First let $s(f)\leq b$.
SIH yields (\ref{eq:singlemainl_I_1_2}).
We conclude (\ref{eq:singlemainl_1_I_1}) by an $(i {\rm -rfl}_{\mS}(g,x,f,\Delta))$.

Next let $s(f)>b$.
We have $\rho_{\mathtt{P}(\mS)}=\rho_{\mathtt{Q}(\mS)}=\rho$
 with $(\rho,g)\in\mathtt{P}(\mS)$, and 
 $\delta^{\mathtt{P}^{\sigma}}_{\mS}=\delta^{\mathtt{Q}}_{\mS}$, while
 $(\sigma,f)\in\mathtt{P}^{\sigma}(\mS)$ with 
$\rho_{\mathtt{P}^{\sigma}(\mS)}=\sigma$.
We obtain $on_{\Lambda}(f)<on_{\Lambda}(g)$
by Proposition \ref{prp:hstepdown_on}.\ref{prp:hstepdown_on_1}.
This yields 
$\alpha(\mathtt{P}^{\sigma})=\alpha(\mathtt{P}^{\sigma};\mS)=\Lambda^{+}\cdot\delta^{\mathtt{Q}}_{\mS}+on_{\Lambda}(f)<
\Lambda^{+}\cdot\delta^{\mathtt{Q}}_{\mS}+on_{\Lambda}(g)=\alpha(\mathtt{Q};\mS)$ and
$\varphi_{\alpha(\mathtt{P}^{\sigma})}(a_{0})<\varphi_{\alpha(\mathtt{Q})}(a_{0})=\widehat{a_{0}}$.

Let $f_{1}=h^{b}(f;\varphi_{d_{i}}(a_{0}))$.
We obtain $\mathrm{fld}(f_{1})\subset\mathcal{H}_{\gamma}[\Theta(\mathtt{Q})]$,
$SC_{\mI_{N}}(\mathrm{fld}(f_{1}))\subset\Lambda$,
$(f_{1})_{x}=f_{x}\leq g_{x}$ and $f_{1}<_{\mI_{N}}^{x}g^{\prime}(x)$
by Proposition \ref{prp:hstepdown}.\ref{prp:hstepdown.3}.

Let $\kappa\in H_{\sigma}(f_{1},\mathtt{P}^{\sigma},\Theta)$.
We obtain by SIH (or by MIH) for  $\bar{\rho}_{\mathtt{P}^{\sigma}}(\mS)=(\sigma,f)$
%HERE we need to replace $\sigma$ by $\kappa$, not to retain $\sigma$.
\beqn\label{eq:singlemainl_I_1_3}
(\mathcal{H}_{\gamma},\Theta_{\mathtt{P}},
\mathtt{P}^{\kappa}
)
\vdash^{\widehat{a_{0}}}_{b,\vec{e},\mT^{\dagger},\Lambda^{++}}
\Delta^{(\bar{\rho}_{\kappa})}, \Phi(\mathtt{P})
\eeqn
where 
$\mathtt{P}^{\kappa}=(\mathtt{P}^{\sigma})^{(\kappa/\sigma)}$,
$\Delta^{(\bar{\rho}_{\kappa})}=(\Delta^{(\bar{\rho}_{\sigma})})^{(\kappa/\sigma)}$
and $\varphi_{\alpha(\mathtt{P}^{\sigma})}(a_{0})<\widehat{a_{0}}$.
Note that $\Phi^{(\kappa/\sigma)}=\Phi$ since $\Phi$ is a set of
capped formulas over $\mathtt{Q}$ with $\sigma\not\in (\mathtt{Q}(\mS))_{0}$,
and $\sigma\not\in\mathcal{H}_{\gamma}[\Theta(\mathtt{Q})]\subset M_{\sigma}$.

Given a $\kappa\in H_{\rho}(f_{1},\mathtt{Q},\Theta)$,
pick a $\sigma\in H_{\rho}(f,\mathtt{Q},\Theta)$ such that $\mathtt{p}_{0}(\sigma)=\mathtt{p}_{0}(\rho)$
by Lemma \ref{lem:exist_H}.\ref{lem:exist_H15}.
Then we obtain
$\kappa\in H_{\sigma}(f_{1},\mathtt{P}^{\sigma},\Theta)$.

We have (\ref{eq:singlemainl_I_1_1}) 
for each successor $j$-stable ordinal $\mV$
such that $j\geq i$, $\mS\leq\mV\leq\mathrm{rk}(\Delta)$ and $\min\{s(g), s(f)\}\geq\mV+1$,
and for each good extension $\mathtt{P}=\mathtt{Q}[\mV]$ such that
$\mathtt{P}(\mV)=\{(\rho_{\mV},g)\}$.
On the other hand we have (\ref{eq:singlemainl_I_1_3}) for each
$\kappa\in H_{\rho}(f_{1},\mathtt{Q},\Theta)$.
We conclude (\ref{eq:singlemainl_1_I_1}) by an $(i {\rm -rfl}_{\mS}(g,x,f_{1},\Delta))$
with $\mathrm{rk}(\Delta)<b$.

Other cases are seen from SIH.
\eprf

\bcor\label{cor:singlemainl_I_1}
Let $\mathtt{Q}$ be a finite family with thresholds $\gamma^{\mathtt{Q}}_{\cdot},\delta^{\mathtt{Q}}_{\cdot}$
such that 
$\mT\in \mathrm{dom}(\mathtt{Q})$, and
$\Lambda^{(n)}>\Lambda_{0}$ strongly critical numbers such that
$\{
\Lambda^{(0)}<\Lambda^{(1)}<\Lambda^{(2)}<\cdots<\Lambda^{(3N-3)}<\Lambda^{(3N-2)}<\Lambda^{(3N-1)}
\}\subset\mathcal{H}_{\gamma}[\Theta(\mathtt{Q})]$.
$b$ denotes an ordinal such that
$b=\mT+1$ if $\mT\in SSt\cup\{0\}$, 
and $b=\mT$ else.
Let 
$\vec{d}=(d_{1},\ldots,d_{N})$ be such that
$d_{1}\leq\cdots\leq d_{N}\leq\mT^{\dagger}$.

Let
$
(\mathcal{H}_{\gamma},\Theta, \mathtt{Q}
)
\vdash^{a}_{b,\vec{d},\mT^{\dagger},\Lambda^{(0)}}
\Phi
$ for $b_{0}+\omega^{\Lambda_{0}+1}\leq\gamma$ and
$a<\Lambda^{(0)}$.

For $1\leq i\leq N$, let $\mS_{i},\kappa_{i},\rho_{i}$ be ordinals, $\mathtt{Q}_{i}$ finite families,
$\Phi_{i}$ sequents and ordinals $\alpha_{i}(\mathtt{Q})=\alpha(\mathtt{Q};\mS_{i})$ 
defined from ordinals $b$, $\alpha_{0}=a$,
 the ordinals $\vec{d}$, and strongly critical numbers $\Lambda^{(n)}$
in Definition \ref{df:main_lemmas}.\ref{df:main_lemmas.2}.

Then 
\[
(\mathcal{H}_{\gamma},\Theta, \mathtt{Q}_{N}
)
\vdash^{\alpha_{N}}_{b,\vec{b},\mT^{\dagger},\Lambda^{(3N-1)}}
\Phi_{N}
\]
holds for $\vec{b}=(b,\ldots, b)$.

\ecor
\bprf
We see from Lemma \ref{mlem:singlemainl_I_1} that
$(\mathcal{H}_{\gamma},\Theta, \mathtt{Q}_{i}
)
\vdash^{\alpha_{i}}_{b,\vec{d}_{i},\mT^{\dagger},\Lambda^{(3i-1)}}
\Phi_{i}
$ by induction on $i\leq N$, where
$\vec{d}_{i}=(b,\ldots,b,d_{i+1},\ldots,d_{N})$.
\eprf

\subsection{Eliminations of inferences (rfl)}\label{subsec:elimrfl}

In this subsection, inferences $(i {\rm -rfl}_{\mS}(g,x,f,\Delta))$ 
are removed from operator controlled derivations of sequents of formulas in $\Sig(\Ome)\cup\Pi(\Ome)$.

\bprp\label{prp:main.2_L_I_1}
Let $\mathtt{Q}$ be a finite family with thresholds $\gamma^{\mathtt{Q}}_{\cdot}, \delta^{\mathtt{Q}}_{\cdot}$, and $\mT\in \mathrm{dom}(\mathtt{Q})$ be a stable ordinal not in $SSt$.
Suppose
$ 
(\mathcal{H}_{\gamma},\Theta, \mathtt{Q})\vdash^{a}_{\mT,\vec{\mT},\mT^{\dagger},\Lambda}\Gamma
$ for $\vec{\mT}=(\mT,\ldots,\mT)$ and $\mathrm{rk}(\Gamma)<\mT$.
Then
$ 
(\mathcal{H}_{\gamma},\Theta, \mathtt{Q}\restrict\mT)\vdash^{a}_{\mT,\vec{\mT},\mT,\Lambda}\Gamma
$.
\eprp
\bprf
By induction on $a$.
We have 
$\mathrm{srk}(\Gamma)\subset \mathrm{dom}(\mathtt{Q})\cap\mT= \mathrm{dom}(\mathtt{Q}\restrict\mT)$ by $\mathrm{rk}(\Gamma)<\mT$ for (\ref{eq:controlder_cap_cover}).
\eprf

\blem\label{lem:main.2_L_I_1}{\rm (Elimination of one limit stable ordinal)}\\
Let $\mathtt{Q}$ be a finite family with thresholds $\gamma^{\mathtt{Q}}_{\cdot},\delta^{\mathtt{Q}}_{\cdot}$, and $\mT\in\mathrm{dom}(\mathtt{Q})$ be a stable ordinal not in $SSt$.
Let
$\Lambda^{(n)}>\Lambda_{0}$ be strongly critical numbers such that
$\{
\Lambda^{(0)}<\Lambda^{(1)}<\Lambda^{(2)}<\cdots<\Lambda^{(3N-3)}<\Lambda^{(3N-2)}<\Lambda^{(3N-1)}
\}\subset\mathcal{H}_{\gamma}[\Theta(\mathtt{Q})]$.
Let $\Phi\subset\bigvee(\mT)$ and
$
(\mathcal{H}_{\gamma},\Theta, \mathtt{Q}
)
\vdash^{a}_{\mT^{\dagger},\vec{d},\mT^{\dagger},\Lambda^{(0)}}
\Phi
$ for 
$a<\Lambda^{(0)}$ and $\vec{d}=(\mT^{\dagger},\ldots,\mT^{\dagger})$.

For $1\leq i\leq N$, let $\mS_{i},\kappa_{i},\rho_{i}$ be ordinals, $\mathtt{Q}_{i}$ finite families,
$\Phi_{i}$ sequents and ordinals $\alpha_{i}(\mathtt{Q})=\alpha(\mathtt{Q};\mS_{i})$ 
defined from ordinals $\mT$, $\alpha_{0}=\varphi_{\mT^{\dagger}}(a)$,
 the ordinals $\vec{d}$, and strongly critical numbers $\Lambda^{(n)}$
in Definition \ref{df:main_lemmas}.\ref{df:main_lemmas.2}.

Then 
\[
(\mathcal{H}_{\gamma},\Theta, \mathtt{Q}_{N}\restrict\mT
)
\vdash^{\alpha_{N}}_{\mT,\vec{\mT},\mT,\Lambda^{(3N-1)}}
\Phi_{N}
\]
holds.

\elem
\bprf
By Cut-elimination \ref{lem:CE} we obtain
$(\mathcal{H}_{\gamma},\Theta, \mathtt{Q})\vdash^{\alpha_{0}}_{\mT,\vec{\mT^{\dagger}},\mT^{\dagger},\Lambda^{(0)}}\Gamma$ for $\alpha_{0}=\varphi_{\mT^{\dagger}}(a)$.
Corollary \ref{cor:singlemainl_I_1}  then yields
$(\mathcal{H}_{\gamma},\Theta, \mathtt{Q}_{N}
)
\vdash^{\alpha_{N}}_{\mT,\vec{\mT},\mT^{\dagger},\Lambda^{(3N-1)}}
\Gamma_{N}$.
By Proposition \ref{prp:main.2_L_I_1} we obtain the lemma.
\eprf

\bdf\label{df:restrict}
{\rm
Let $\mathtt{Q}$ be a finite family such that $\mS\in\mathrm{dom}(\mathtt{Q})\subset\mS^{\dagger}$.
$\mathtt{Q}\restrict\mS$ denotes the restriction of $\mathtt{Q}$ to $\mS$.

\benu
\item\label{df:restrict.1}
For a cap $\bar{\rho}\in\prod_{\mT\in\mathrm{dom}(\mathtt{Q})}\mathtt{Q}(\mT)$ over $\mathtt{Q}$,
$\bar{\rho}\restrict\mS\in\prod_{\mT\in\mathrm{dom}(\mathtt{Q}\, \restrict\, \mS)}\mathtt{Q}(\mT)$ denotes the restriction of $\bar{\rho}$
to $\mathrm{dom}(\mathtt{Q}\restrict\mS)$.

\item\label{df:restrict.2}
Let $\mS\in SSt$ be a successor stable ordinal, and assume $\mathrm{dom}(\mathtt{Q})\cap SSt\subset M_{\rho}$ for $\rho=\rho_{\mathtt{Q}(\mS)}$.
For  a capped formula $A^{(\bar{\rho})}$ over $\mathtt{Q}$ such that $\mathsf{k}(A^{(\bar{\rho})};\mathtt{Q})\subset M_{\mathtt{Q}}$,
$(A^{(\bar{\rho})})^{[\mS]}(\mathtt{Q})$ denotes a capped formula over
$\mathtt{Q}\restrict\mS$ defined as follows.

 \benu

 \item\label{df:restrict.21}
Let
$A\equiv B^{[\sigma/\mS]}\not\equiv B\equiv \mathrm{inv}^{(\mathtt{u})}(A^{(\bar{\rho})};\mathtt{Q})$ and $\mathsf{k}(B)\subset M_{\mathtt{Q}}\subset M_{\rho}$ with
$\rho=\rho_{\mathtt{Q}(\mS)}\leq\sigma\in(\mathtt{Q}(\mS))_{0}\cap(\bar{\rho}(\mS))_{0}$,
cf.\,Definition \ref{df:inverse}.\ref{df:inverse2}.
In this case let 
$(A^{(\bar{\rho})})^{[\mS]}(\mathtt{Q}):\equiv (B^{[\rho/\mS]})^{(\bar{\rho}\,\restrict\,\mS)}$.

 \item\label{df:restrict.23}
Otherwise let
$(A^{(\bar{\rho})})^{[\mS]}(\mathtt{Q}) :\equiv (A^{[(\bar{\rho}(\mS))_{0}/\mS]})^{(\bar{\rho}\,\restrict\,\mS)}$.
 \eenu

\item\label{df:restrict.3}
For a set $\Gamma$ of capped formulas over $\mathtt{Q}$ such that $\mathsf{k}(A^{(\bar{\rho})};\mathtt{Q})\subset M_{\mathtt{Q}}$ for each $A^{(\bar{\rho})}\in\Gamma$, let
$\Gamma^{[\mS]}(\mathtt{Q})=\{(A^{(\bar{\rho})})^{[\mS]}(\mathtt{Q}) : A^{(\bar{\rho})}\in\Gamma\}$.
\eenu
}
\edf

We have $\rho\in\bigcap_{\mU<\mS}M_{\mathtt{Q}(\mU)}$ for $\rho\in(\mathtt{Q}(\mS))_{0}$.
We obtain $\mathsf{k}(A^{[\rho/\mS]})\subset M_{\bar{\rho}\,\restrict\,\mS}$ in Definition \ref{df:restrict}.\ref{df:restrict.2}.

Note that $\mathsf{k}(A)\subset M_{(\bar{\rho}(\mS))_{0}}$ holds in Definition \ref{df:restrict}.\ref{df:restrict.23}.
Suppose
$A\equiv B^{[\sigma/\mT]}\not\equiv B\equiv \mathrm{inv}^{(\mathtt{u})}(A^{(\bar{\rho})};\mathtt{Q})$ and $\mathsf{k}(B)\subset M_{\mathtt{Q}}\subset M_{\rho}$, where
$\mT\in\mathrm{dom}(\mathtt{Q})\cap SSt\cap\mS$. By the assumption we obtain
$\mT\in\mathrm{dom}(\mathtt{Q})\cap SSt\cap\mS\subset M_{\rho_{\mathtt{Q}(\mS)}}\cap\mS\subset 
M_{(\bar{\rho}(\mS))_{0}}\cap \mS=(\bar{\rho}(\mS))_{0}$, and
$\mathsf{k}(B^{[\sigma/\mT]})\subset\mT\subset(\bar{\rho}(\mS))_{0}$.
Hence
$(A^{(\bar{\rho})})^{[\mS]}(\mathtt{Q}) \equiv ((B^{[\sigma/\mT]})^{(\bar{\rho}\,\restrict\,\mS)}$.

\bprp\label{prp:restrict_I_1.2}
Let $\mS\in\mathrm{dom}(\mathtt{Q})\cap SSt\subset\mS^{\dagger}$ and $\mathsf{k}(A^{(\bar{\rho})};\mathtt{Q})\subset M_{\mathtt{Q}}$.

Then $\mathrm{rk}((A^{(\bar{\rho})})^{[\mS]}(\mathtt{Q}))<\mS$,
$\mathrm{srk}((A^{(\bar{\rho})})^{[\mS]}(\mathtt{Q}))\in
(\{\mathrm{srk}(A^{(\bar{\rho})})\}\cup\mathrm{cl}(\mS)\cup(\{\rho_{\mathtt{Q}(\mS)}\}\cap St))\cap\mS$, and
if $\mathrm{srk}(A^{(\bar{\rho})})<\mS$, then $\mathrm{srk}(A^{(\bar{\rho})})=\mathrm{srk}((A^{(\bar{\rho})})^{[\mS]}(\mathtt{Q}))$.
\eprp
\bprf
We see readily $\mathrm{rk}((A^{(\bar{\rho})})^{[\mS]}(\mathtt{Q}))<\mS$ from Definition \ref{df:restrict}.
Let $\alpha\in M_{\rho}$ with $\rho\in\Psi_{\mS}$.
$\alpha[\rho/\mS]<\mS$ holds if $\mS\leq\alpha$, and $\alpha[\rho/\mS]=\alpha$ else.
We have
$\mathrm{srk}(\alpha[\rho/\mS])\subset \mathrm{srk}(\alpha)\cup
(\mathrm{cl}(\mS)\cap\rho)\cup(\{\rho\}\cap St)$
by Proposition \ref{prp:pd_closed_cover_I_1}.\ref{prp:pd_closed_cover_I_1.3}.

Let $u$ be an $RS$-term.
We obtain $\mathrm{srk}(u^{[\rho/\mS]})\subset\mathrm{srk}(u)\cup \mathrm{cl}(\mS)\cup(\{\rho\}\cap St)$,
and $\mathrm{rk}(u)<\mS \Rarw \mathrm{srk}(u)=\mathrm{srk}(u^{[\rho/\mS]})$
if $\mathsf{k}(u)\subset  M_{\rho}$.

Also
$\mathrm{srk}((u^{[\rho/\mS]})^{[\sigma/\mT]})=\mathrm{srk}(u^{[\sigma/\mT]})$ if $\mT<\mS$, $\sigma\in \Psi_{\mT}$, $\mathsf{k}(u)\subset M_{\sigma}\cap M_{\rho}$
and $\mathsf{k}(u^{[\rho/\mS]})\subset M_{\sigma}$.
\eprf

\blem\label{mlem:singlemainl_S_I_1}
Let $\mathtt{Q}$ be a finite family with thresholds $\gamma^{\mathtt{Q}}_{\cdot}, \delta^{\mathtt{Q}}_{\cdot}$  such that 
$\mS\in\mathrm{dom}(\mathtt{Q})\cap SSt$.
Let 
$\Gamma$ be a set of $\bigvee$-formulas such that $\mathrm{rk}(\Gamma)<\mS+1$.
Let
$
(\mathcal{H}_{\gamma},\Theta, \mathtt{Q}
)
\vdash^{a}_{\mS+1,\vec{\mS}+1,\mS^{\dagger},\Lambda}
\Gamma
$ for $a<\Lambda$ and $\vec{\mS}+1=(\mS+1,\ldots,\mS+1)$. 

Let $\mathtt{R}$ denote the restriction $\mathtt{Q}\restrict\mS$ of $\mathtt{Q}$ to $\mS$, 
$\Theta_{\mS}=\Theta\cup(\mathtt{Q}(\mS))_{0}$ and
$\gamma_{1}=\gamma_{\mathbb{S}}^{\mathtt{Q}}+\mS$.
Then 
\begin{equation}\label{eq:singlemainl_S1_I_1}
(\mathcal{H}_{\gamma_{1}},\Theta_{\mS}, \mathtt{R}
)
\vdash^{\omega a}_{\mS,\vec{\mS},\mS,\Lambda}
\Gamma^{[\mS]}(\mathtt{Q})
\end{equation}
holds.
\elem
\bprf
By induction on $a$.
We write $(A^{(\bar{\rho})})^{[\mS]}$ for $(A^{(\bar{\rho})})^{[\mS]}(\mathtt{Q})$ when $\mathtt{Q}$ is fixed.
We have
$\mathrm{cl}(\mS)\cup(\{\rho_{\mathtt{Q}(\mS)}\}\cap St)\subset\mathrm{dom}(\mathtt{Q})$
by $\mS\in\mathrm{dom}(\mathtt{Q})$ and Definition \ref{df:ffthreshold}.\ref{df:QJ}.
We obtain
$\mathrm{dom}(\mathtt{R})\cup\{\mathrm{srk}((A^{(\bar{\rho})})^{[\mS]})\}=
(\mathrm{dom}(\mathtt{Q})\cup\{\mathrm{srk}(A)\})\cap\mS$
by Proposition \ref{prp:restrict_I_1.2}.
Hence we may assume that
$\mathtt{R}[(\mathrm{srk}((A^{(\bar{\rho})})^{[\mS]})]= (\mathtt{Q}[\mathrm{srk}(A)])\restrict\mS$.

We have 
 $\{a,\mS\}
\subset
\mathcal{H}_{\gamma}[\Theta(\mathtt{Q})]$ by (\ref{eq:controlder_cap22}), and $\gamma\leq \gamma^{\mathtt{Q}}_{\mS}$
by (\ref{eq:controlder_cap0}).
We obtain
$\omega a\in\mathcal{H}_{\gamma_{1}}[\Theta_{\mS}(\mathtt{R})]$ for (\ref{eq:controlder_cap22}).

We have 
$\gamma_{\mathbb{S}}^{\mathtt{Q}}<\gamma_{1}=\gamma_{\mathbb{S}}^{\mathtt{Q}}+\mS<\gamma_{\mathbb{T}}^{\mathtt{Q}}$ for every 
$\{\mT<\mS\}\subset \mathrm{dom}(\mathtt{Q})$
by Definition \ref{df:caphat}.\ref{df:caphat.42}.
Hence 
(\ref{eq:controlder_cap0}) is enjoyed in (\ref{eq:singlemainl_S1_I_1}).

\bclm\label{clm:singlemainl_next_I_1}
Each of 
(\ref{eq:controlder_cap_cover}) and
(\ref{eq:controlder_cap2}) is enjoyed in (\ref{eq:singlemainl_S1_I_1}).
\eclm
\textbf{Proof} of Claim \ref{clm:singlemainl_next_I_1}.
Let $A^{(\bar{\rho})}\in\Gamma$.
We have $\mT<\mS$ for $\mT<\mathrm{srk}((A^{(\bar{\rho})})^{[\mS]}(\mathtt{Q}))$, and 
$\mathrm{srk}((A^{(\bar{\rho})})^{[\mS]}(\mathtt{Q}))\in\mathrm{dom}(\mathtt{R})=\mathrm{dom}(\mathtt{Q})\cap\mS$
by Proposition \ref{prp:restrict_I_1.2}.
Hence (\ref{eq:controlder_cap_cover}) follows.

Next we show $\mathsf{k}((A^{(\bar{\rho})})^{[\mS]};\mathtt{R})\subset
\mathcal{H}_{\gamma_{1}}[\Theta_{\mS}(\mathtt{R})]$ for (\ref{eq:controlder_cap2}).
We have $\mathsf{k}(A^{(\bar{\rho})};\mathtt{Q})\subset\mathcal{H}_{\gamma}[\Theta(\mathtt{Q})]\subset M_{\mathtt{Q}}$.
On the other hand we have $(\mathtt{Q}(\mS))_{0}\subset M_{\mathtt{R}}$ by 
Definition \ref{df:ffthreshold}.\ref{df:QJ.2}.
We obtain
$\Theta_{\mS}(\mathtt{R})=(\Theta\cup(\mathtt{Q}(\mS))_{0})\cap M_{\mathtt{R}}=
\Theta(\mathtt{R})\cup(\mathtt{Q}(\mS))_{0}$ with $\Theta(\mathtt{Q})\subset\Theta(\mathtt{R})$.
Let $\rho=\rho_{\mathtt{Q}(\mS)}$.

First consider the case in Definition \ref{df:restrict}.\ref{df:restrict.21}.
Let $\rho=\rho_{\mathtt{Q}(\mS)}\leq\sigma<(\bar{\rho}(\mS))_{0}$.
Then we have $A\equiv B^{[\sigma/\mS]}\not\equiv B$ and
$B^{[\rho/\mS]}\equiv \mathrm{inv}^{(\mathtt{u})}((A^{(\bar{\rho})})^{[\mS]};\mathtt{R})\not\equiv\mathrm{inv}^{(\mathtt{u})}(A^{(\bar{\rho})};\mathtt{Q})
\equiv B$.
By $\mathsf{k}(B)=\mathsf{k}(A^{(\bar{\rho})};\mathtt{Q})\subset\mathcal{H}_{\gamma}[\Theta(\mathtt{Q})]$,
we obtain $\mathsf{k}(B^{[\rho/\mS]})=\mathsf{k}((B^{[\rho/\mS]})^{(\bar{\rho}\, \restrict\, \mS)};\mathtt{R})=\mathsf{k}((A^{(\bar{\rho})})^{[\mS]};\mathtt{R})\subset\mathcal{H}_{\gamma_{1}}[\Theta_{\mS}(\mathtt{R})]$
by Proposition \ref{prp:EH_Pi11}.\ref{prp:EH_Pi11.4}, $\gamma_{1}\geq\mS$ and 
$\rho\in(\mathtt{Q}(\mS))_{0}\subset\mathcal{H}_{\gamma_{1}}[\Theta_{\mS}(\mathtt{R})]$.

Second let $(A^{(\bar{\rho})})^{[\mS]}\equiv (A^{[(\bar{\rho}(\mS))_{0}/\mS]})^{(\bar{\rho}\,\restrict\,\mS)}$ in Definition \ref{df:restrict}.\ref{df:restrict.23}.
Then $\mathsf{k}(A^{[(\bar{\rho}(\mS))_{0}/\mS]})=\mathsf{k}( \mathrm{inv}^{(\mathtt{u})}((A^{(\bar{\rho})})^{[\mS]};\mathtt{R}))\subset\mathcal{H}_{\gamma_{1}}[\Theta_{\mS}(\mathtt{R})]$ is seen as in the first case from 
$(\bar{\rho}(\mS))_{0}\in(\mathtt{Q}(\mS))_{0}\subset\mathcal{H}_{\gamma_{1}}[\Theta_{\mS}(\mathtt{R})]$.
\hspace*{\fill} $\Box$ of Claim \ref{clm:singlemainl_next_I_1}.
\\
\textbf{Case 1}. First the last inference is a $(\bigvee)$ introducing an $A^{(\bar{\rho})}\in\Gamma$:
Let 
$A\simeq\bigvee(B_{\iota})_{\iota\in J}$ and $B_{\iota}\simeq\bigwedge(A_{\nu})_{\nu\in J_{\iota}}$.
Let 
 $\iota\in[\bar{\rho}\mathtt{Q}]J$ be such that
$\mathsf{k}(B_{\iota}^{(\bar{\rho})};\mathtt{Q})\subset\mathcal{H}_{\gamma}[\Theta(\mathtt{Q})]$
for which there is an $a_{0}\in\mathcal{H}_{\gamma}[\Theta(\mathtt{Q})]\cap a$ such that
$
(\mathcal{H}_{\gamma},\Theta_{\nu}, \mathtt{Q}_{\nu}
)
\vdash^{a_{0}}_{\mS+1,\vec{\mS}+1,\mS^{\dagger},\Lambda}
(\Gamma\cup\{A_{\nu}^{(\bar{\rho})}\})(\mathtt{Q}_{\nu})
$ holds for every $\nu\in[\bar{\rho}\mathtt{Q}]J_{\iota}$, where
$\Theta_{\nu}=\Theta_{\mathtt{Q}_{\nu}}\cup \mathsf{k}(A_{\nu}^{(\bar{\rho})};\mathtt{Q})$,
$\mathtt{Q}_{\nu}= \mathtt{Q}[\mathrm{srk}(A_{\nu})]$ is a good extension, and
$\mathsf{k}(A_{\nu}^{(\bar{\rho})};\mathtt{Q})\subset \mathcal{H}_{\gamma}[\Theta(\mathtt{Q})]$.

For $A^{(\bar{\rho})}\in\Gamma$, we have
$\mathtt{R}[\mathrm{srk}((A_{\nu}^{(\bar{\rho})})^{[\mS]})]=\mathtt{Q}[\mathrm{srk}(A_{\nu})]\restrict\mS$.
Let $\rho=\rho_{\mathtt{Q}(\mS)}$ and $\bar{\rho}^{\mS}=\bar{\rho}\restrict\mS$.

There are two cases for the formula $(A^{(\bar{\rho})})^{[\mS]}$.
\\
\textbf{Case 1.1}. The case in Definition \ref{df:restrict}.\ref{df:restrict.21}:
Let $(\bar{\rho}(\mS))_{0}>\sigma\geq\rho_{\mathtt{Q}(\mS)}=\rho$.
We have $(A^{(\bar{\rho})})^{[\mS]}\equiv (C^{[\rho/\mS]})^{(\bar{\rho}^{\mS})}$ for a $C$ such that
$A\equiv C^{[\sigma/\mS]}$ and $\mathsf{k}(C)\subset M_{\mathtt{Q}}$.
Let $C\simeq\bigvee(C_{\mu})_{\mu\in I}$. Then $A\equiv C^{[\sigma/\mS]}\simeq\bigvee(C_{\mu}^{[\sigma/\mS]})_{\mu\in [\sigma]I}$ and
$C^{[\rho/\mS]}\simeq\bigvee(C_{\mu}^{[\rho/\mS]})_{\mu\in [\rho]I}$.
Let $\mu\in I$ be such that $B_{\iota}\equiv C_{\mu}^{[\sigma/\mS]}$ for $\iota=\mu^{[\sigma/\mS]}$.
We have
$\mathsf{k}(C_{\mu})=\mathsf{k}(B_{\iota}^{(\bar{\rho})};\mathtt{Q})\subset\mathcal{H}_{\gamma}[\Theta(\mathtt{Q})]\cap M_{\mathtt{Q}}\subset M_{\rho}$.
We see
$\mathsf{k}((C_{\mu}^{[\rho/\mS]})^{(\bar{\rho}^{\mS})};\mathtt{R})=\mathsf{k}(C_{\mu}^{[\rho/\mS]})
\subset\mathcal{H}_{\gamma_{1}}[\Theta_{\mS}(\mathtt{R})]\cap M_{\mathtt{R}}$ from $\rho\in M_{\mathtt{R}}$ as in Claim \ref{clm:singlemainl_next_I_1}.

Let $C_{\mu}\simeq\bigwedge(D_{\alpha})_{\alpha\in K}$.
We obtain $(B_{\iota}^{(\bar{\rho})})^{[\mS]}\equiv (C_{\mu}^{[\rho/\mS]})^{(\bar{\rho}^{\mS})}$ and
$C_{\mu}^{[\rho/\mS]}\simeq\bigwedge(D_{\alpha}^{[\rho/\mS]})_{\alpha\in [\rho]K}$,
while $B_{\iota}\equiv C_{\mu}^{[\sigma/\mS]}\simeq\bigwedge(D_{\alpha}^{[\sigma/\mS]})_{\alpha\in [\sigma]K}$.
On the other hand we have
$\mathsf{k}((D_{\alpha}^{[\rho/\mS]})^{(\bar{\rho}^{\mS})};\mathtt{R})=\mathsf{k}(D_{\alpha}^{[\rho/\mS]})$ and
$\mathsf{k}((D_{\alpha}^{[\sigma/\mS]})^{(\bar{\rho})};\mathtt{Q})=\mathsf{k}(D_{\alpha})$.

\bclm\label{clm:singlemainl_SA}
Let $\alpha\in[\rho]K$ be such that $\mathsf{k}(D_{\alpha}^{[\rho/\mS]})\subset M_{\mathtt{R}}$.
Let $\nu=\alpha^{[\sigma/\mS]}$, $A_{\nu}\equiv D_{\alpha}^{[\sigma/\mS]}$,
and $\Theta_{\mS,\nu}=\Theta\cup(\mathtt{Q}(\mS))_{0}\cup \mathsf{k}(D_{\alpha}^{[\rho/\mS]})$.
Then $\mathsf{k}(D_{\alpha})
\subset\mathcal{H}_{\gamma_{1}}[\Theta_{\mS.\nu}(\mathtt{R})]\cap M_{\mathtt{Q}}$
and
$\nu\in[\bar{\rho}\mathtt{Q}]J_{\iota}$.
\eclm
\textbf{Proof} of Claim \ref{clm:singlemainl_SA}.
We have $\mathsf{k}(D_{\alpha})\subset M_{\rho}\subset M_{\sigma}$.
We obtain $(A_{\nu}^{(\bar{\rho})})^{[\mS]}\equiv (D_{\alpha}^{[\rho/\mS]})^{(\bar{\rho}^{\mS})}$ and
$\mathsf{k}(A_{\nu}^{(\bar{\rho})};\mathtt{Q})=\mathsf{k}(D_{\alpha})\subset M_{\mathtt{R}}\cap M_{\rho}=M_{\mathtt{Q}}$
by $\rho=\rho_{\mathtt{Q}(\mS)}$ and Lemma \ref{lem:fTaut}.
Hence $\nu\in[\bar{\rho}\mathtt{Q}]J_{\iota}$.
On the other hand 
we have $\max\{\mS,\mathtt{p}_{0}(\rho)\}\leq \gamma_{1}$ and 
$\{\rho\}\cup \mathsf{k}(D_{\alpha}^{[\rho/\mS]}) \subset \mathcal{H}_{\gamma_{1}}[\Theta_{\mS,\nu}(\mathtt{R})]$.
Lemma \ref{lem:inverse} yields
$\mathsf{k}(A_{\nu}^{(\bar{\rho})};\mathtt{Q})=\mathsf{k}(D_{\alpha})\subset \mathcal{H}_{\gamma_{1}}[\Theta_{\mS,\nu}(\mathtt{R})]$.
\hspace*{\fill} $\Box$ of Claim \ref{clm:singlemainl_SA}.

Given a good extension $\mathtt{R}_{\nu}=\mathtt{R}[\mathrm{srk}((A_{\nu}^{(\bar{\rho})})^{[\mS]})]$ of $\mathtt{R}$ for $\gamma_{1}$ and 
$\Theta_{\mS,\nu}=\Theta\cup\mathtt{Q}(\mS)\cup \mathsf{k}(A_{\nu}^{(\bar{\rho})};\mathtt{Q})$, let 
$\mathtt{Q}_{\nu}= \mathtt{Q}[\mathrm{srk}(A_{\nu})]$ be an extension of $\mathtt{Q}$ such that
$\mathtt{Q}_{\nu}\restrict\mS=\mathtt{R}_{\nu}$, where $\mathtt{Q}_{\nu}=\mathtt{R}_{\nu}[\mS]$ is an extension of $\mathtt{R}_{\nu}$
by the ordinal $\gamma^{\mathtt{Q}}_{\mS}$, and a good extension of $\mathtt{Q}$
for $\gamma$ and $\Theta_{\nu}$. 
IH with Claim \ref{clm:singlemainl_SA} yields
\[
(\mathcal{H}_{\gamma_{1}},\Theta_{\mS,\nu}, \mathtt{R}_{\nu}
)
\vdash^{\omega a_{0}}_{\mS,\vec{\mS},\mS,\Lambda}
(\Gamma^{[\mS]}\cup\{((D_{\alpha}^{[\sigma/\mS]})^{(\bar{\rho})})^{[\mS]}\})(\mathtt{R}_{\nu})
\]
for every $\alpha\in[\rho]K$ such that $\mathsf{k}(D_{\alpha}^{[\rho/\mS]})\subset M_{\mathtt{R}}$, where 
$((D_{\alpha}^{[\sigma/\mS]})^{(\bar{\rho})})^{[\mS]}\equiv (D_{\alpha}^{[\rho/\mS]})^{(\bar{\rho}^{\mS})}$ and
$\mathtt{R}_{\nu}=\mathtt{R}[\mathrm{srk}((A_{\nu}^{(\bar{\rho})})^{[\mS]})]$.
We obtain 
$\{\mathrm{srk}(D_{\alpha}^{[\sigma/\mS]})\}\cap SSt=\{\mathrm{srk}(D_{\alpha}^{[\rho/\mS]})\}\cap SSt$.
On the other hand we have
$\mathsf{k}((C_{\mu}^{[\rho/\mS]})^{(\bar{\rho}^{\mS})};\mathtt{R})
\subset\mathcal{H}_{\gamma_{1}}[\Theta_{\mS}(\mathtt{R})]\cap M_{\mathtt{R}}$.
A $(\bigvee)$ yields (\ref{eq:singlemainl_S1_I_1}).
\\
\textbf{Case 1.2}. The case in Definition \ref{df:restrict}.\ref{df:restrict.23}:
Then 
$(A^{(\bar{\rho})})^{[\mS]}\equiv (A^{[(\bar{\rho}(\mS))_{0}/\mS]})^{(\bar{\rho}^{\mS})}$.
By the assumption we have $\mathrm{rk}(A)\leq\mS$.
If $\mathrm{rk}(A)<\mS$, then we have 
$\mathsf{k}(A)\subset M_{\rho_{\mathtt{Q}(\mS)}}\cap \mS\subset\rho_{\mathtt{Q}(\mS)}$.
\\
\textbf{Case 1.2.1}. $A\equiv(\exists x\in\mathsf{L}_{\mS}\, C(x))\in\Sigma_{1}(\mS)$ for a formula $C$:
Then $A\simeq\bigvee((v\dot{\in}\mathsf{L}_{\mS})^{\lor}\land C(v)^{\lor})_{v\in Tm(\mS)}$, and $\mathsf{k}(A)\subset M_{\mathtt{Q}}$,
where $(v\dot{\in}\mathsf{L}_{\mS})\equiv\bigwedge\emptyset$.
Let $\iota=v\in Tm(\mS)$ be such that $B_{\iota}\equiv C(v)$.
If $\mathrm{inv}^{(\mathtt{u})}(C(v)^{(\bar{\rho})};\mathtt{Q})\equiv C(v)$, then $\mathsf{k}(v)\subset M_{\rho_{\mathtt{Q}(\mS)}}\cap \mS$ and
$(C(v)^{(\bar{\rho})})^{[\mS]}\equiv C(v)^{(\bar{\rho}^{\mS})}$.
IH yields
$(\mathcal{H}_{\gamma_{1}},(\Theta_{\mS})_{\mathtt{R}_{v}}, \mathtt{R}_{v}
)
\vdash^{\omega a_{0}}_{\mS,\vec{\mS},\mS,\Lambda}
(\Gamma^{[\mS]}\cup\{((C(v)^{\lor})^{(\bar{\rho})})^{[\mS]}\})(\mathtt{R}_{v})$, where
$\mathtt{R}_{v}= \mathtt{R}[\mathrm{srk}(C(v))]$.
On the other hand we have
$(\mathcal{H}_{\gamma_{1}},\Theta_{\mS}, \mathtt{R}
)
\vdash^{\widehat{a_{0}}}_{\mS,\vec{\mS},\mS,\Lambda}\Gamma^{[\mS]},
(((v\dot{\in}\mathsf{L}_{\mS})^{\lor})^{(\bar{\rho})})^{[\mS]}$.
A $(\bigvee)$ yields (\ref{eq:singlemainl_S1_I_1}).

Let $\mathrm{inv}^{(\mathtt{u})}(C(v)^{(\bar{\rho})};\mathtt{Q})\equiv D\not\equiv C(v)$.
Then there exists a term $u$ such that $C(v)\equiv C(u^{[\sigma/\mS]})\equiv C(u)^{[\sigma/\mS]}$ 
for a $\sigma\in(\mathtt{Q}(\mS))_{0}\cap (\bar{\rho}(\mS))_{0}$, and $D\equiv C(u)$.
We obtain $\mathsf{k}(C(u))=\mathsf{k}(((C(u)^{[\sigma/\mS]})^{(\bar{\rho})};\mathtt{Q})\subset\mathcal{H}_{\gamma}[\Theta(\mathtt{Q})]\cap M_{\mathtt{Q}}$, and Definition \ref{df:restrict}.\ref{df:restrict.21} applies to $(C(u)^{[\sigma/\mS]})^{(\bar{\rho})}$ so that
$(((C(u)^{[\sigma/\mS]})^{(\bar{\rho})})^{[\mS]}\equiv (C(u^{[\rho/\mS]}))^{(\bar{\rho}^{\mS})}$.
By $\rho=\rho_{\mathtt{Q}(\mS)}\in (\mathtt{Q}(\mS))_{0}\subset M_{\mathtt{R}}$ we obtain
$\mathsf{k}( (C(u^{[\rho/\mS]}))^{(\bar{\rho}^{\mS})};\mathtt{R})=\mathsf{k}(C(u^{[\rho/\mS]}))
\subset\mathcal{H}_{\gamma_{1}}[\Theta_{\mS}(\mathtt{R})] \cap M_{\mathtt{R}}$.

IH yields
$(\mathcal{H}_{\gamma_{1}},(\Theta_{\mS})_{\mathtt{R}_{u}}, \mathtt{R}_{u}
)
\vdash^{\omega a_{0}}_{\mS,\vec{\mS},\mS,\Lambda}
(\Gamma^{[\mS]}\cup\{((C(u^{[\rho/\mS]})^{\lor})^{(\bar{\rho})})^{[\mS]}\})(\mathtt{R}_{u})$, where 
$\mathtt{R}_{u}=\mathtt{R}[\mathrm{srk}(C(u^{[\rho/\mS]}))]$ and
$\{\mathrm{srk}(C(u^{[\rho/\mS]}))\}\cap SSt=\{\mathrm{srk}(C(u^{[\sigma/\mS]}))\}\cap SSt$.
We have
$(\mathcal{H}_{\gamma_{1}},\Theta_{\mS}, \mathtt{R}
)
\vdash^{\widehat{a_{0}}}_{\mS,\vec{\mS},\mS,\Lambda}\Gamma^{[\mS]},
(((u^{[\rho/\mS]}\dot{\in}\mathsf{L}_{\mS})^{\lor})^{(\bar{\rho})})^{[\mS]}$.

On the other hand we have $u^{[\rho/\mS]}\in Tm((\bar{\rho}(\mS))_{0})$ by 
$(\bar{\rho}(\mS))_{0}>\sigma\geq\rho=\rho_{\mathtt{Q}(\mS)}$, 
$(A^{(\bar{\rho})})^{[\mS]}\equiv ((\exists x\in\mathsf{L}_{\mS}\, C(x))^{(\bar{\rho})})^{[\mS]}\equiv
( (\exists x\in\mathsf{L}_{\mS}\, C(x))^{[(\bar{\rho}(\mS))_{0}/\mS]})^{(\bar{\rho}^{\mS})}\equiv 
(\exists x\in\mathsf{L}_{(\bar{\rho}(\mS))_{0}}C(x))^{(\bar{\rho}^{\mS})}$, and $(u^{[\rho/\mS]}\dot{\in}\mathsf{L}_{(\bar{\rho}(\mS))_{0}})\equiv\bigwedge\emptyset$.
A $(\bigvee)$ yields (\ref{eq:singlemainl_S1_I_1}).
\\
\textbf{Case 1.2.2}.
Otherwise:
Then $\mathsf{k}(A)\subset M_{\rho_{\mathtt{Q}(\mS)}}\cap\mS=\rho_{\mathtt{Q}(\mS)}\subset
(\bar{\rho}(\mS))_{0}$,
$A^{[(\bar{\rho}(\mS))_{0}/\mS]}\equiv A$,
and $(A^{(\bar{\rho})})^{[\mS]}\equiv A^{(\bar{\rho}^{\mS})}$.
We obtain
$\mathsf{k}(\iota)\subset\rho_{\mathtt{Q}(\mS)}$,
$(B_{\iota}^{(\bar{\rho})})^{[\mS]}\equiv (B_{\iota}^{[(\bar{\rho}(\mS))_{0}/\mS]})^{(\bar{\rho}^{\mS})}$,
$\mathsf{k}(B_{\iota}^{(\bar{\rho}^{\mS})};\mathtt{R})=\mathsf{k}(B_{\iota}^{(\bar{\rho})};\mathtt{Q})
\subset\mathcal{H}_{\gamma}[\Theta(\mathtt{Q})]\subset\mathcal{H}_{\gamma_{1}}[\Theta_{\mS}(\mathtt{R})]$
and $\iota\in[\bar{\rho}^{\mS}\mathtt{R}]J$.
Let $\nu\in[\bar{\rho}^{\mS}\mathtt{R}]J_{\iota}$. Then we obtain
$\nu\in[\bar{\rho}\mathtt{Q}]J_{\iota}$.

IH yields
$(\mathcal{H}_{\gamma_{1}},\Theta_{\mS,\nu},\mathtt{R}_{\nu}
)
\vdash^{\omega a_{0}}_{\mS,\vec{\mS},\mS,\Lambda}
(\Gamma^{[\mS]}\cup\{((A_{\nu})^{(\bar{\rho})})^{[\mS]}\})(\mathtt{R}_{\nu})$ for each $\nu\in[\bar{\rho}^{\mS}\mathtt{R}]J_{\iota}$,
where
$\mathtt{R}_{\nu}= \mathtt{R}[\mathrm{srk}(A_{\nu}^{[\bar{\rho}(\mS)/\mS]})]$.
A $(\bigvee)$ yields (\ref{eq:singlemainl_S1_I_1}).
\\
\textbf{Case 2}. 
Second the last inference is an $(i {\rm -rfl}_{\mT}(g,x,f,\Delta))$ with 
$\mT\in \mathrm{dom}(\mathtt{Q})\subset\mS^{\dagger}$: Let $\bar{\rho}=\bar{\rho}_{\mathtt{Q}}$ and $\bar{\rho}^{\mS}=\bar{\rho}\restrict\mS$.
We have an ordinal $a_{0}<a$, 
a finite set $\Delta$ of formulas such that $\mathrm{rk}(\Delta)<\mS+1=d_{i}+1$,
where $d_{i}$ is the $i$-th element in the list $\vec{\mS}=(\mS,\ldots,\mS)$.

Let $\mU$ be a successor $j$-stable ordinal 
such that $j \geq i$, $\mT\leq\mU\leq\mathrm{rk}(\Delta)\leq\mS$ and $\min\{s(g), s(f)\}\geq\mU+1$.
Let $\mathtt{P}=\mathtt{Q}[\mU]$ be a good extension of $\mathtt{Q}$ for $\gamma,\Theta$ 
such that
$(\bar{\rho_{\mathtt{Q}}}(\mT))_{1}=g=(\bar{\rho_{\mathtt{P}}}(\mU))_{1}$.
Let $\rho_{\mU}=\rho_{\mathtt{P}(\mU)}$.

Then
$(\mathcal{H}_{\gamma},\Theta_{\iota}, \mathtt{P}_{\iota})\vdash^{a_{0}}_{\mS+1,\vec{\mS}+1,\mS^{\dagger},\Lambda}
(\Gamma(\mathtt{P})\cup\{\lnot \delta_{\iota}^{(\bar{\rho}_{\mathtt{P}})}\})(\mathtt{P}_{\iota})$ holds
for each $\delta\in\Delta$, $\iota\in[\bar{\rho}_{\mathtt{P}}\mathtt{P}]J_{\delta}$ with $\delta\simeq\bigvee(\delta_{\iota})_{\iota\in J_{\delta}}$
and each good extension $\mathtt{P}_{\iota}=\mathtt{P}[\mathrm{srk}(\delta_{\iota})]$ where
$\Theta_{\iota}=\Theta_{\mathtt{P}_{\iota}}\cup\mathsf{k}(\delta_{\iota}^{(\bar{\rho}_{\mathtt{P}})};\mathtt{P})$.
On the other,
$(\mathcal{H}_{\gamma},\Theta_{\mathtt{P}},
\mathtt{P}^{\sigma}
)
\vdash^{a_{0}}_{\mS+1,\vec{\mS}+1,\mS^{\dagger},\Lambda}
\Delta^{(\bar{\rho}_{\sigma})}, \Gamma(\mathtt{P})
$
for every $\sigma\in H_{\rho_{\mU}}(f,\mathtt{P},\Theta)$ with
$\mathtt{P}^{\sigma}=\mathtt{P}\cup\{(\mU,(\sigma,f))\}$, where $\bar{\rho}_{\sigma}=\bar{\rho}_{\mathtt{P}}[\sigma/\rho_{\mU}]$.

Let $\mathtt{S}=(\mathtt{P}\restrict\mS)[\mU]$.
We have $\mathtt{S}=\mathtt{R}[\mU]$ when $\mU<\mS$.
\\
\textbf{Case 2.1}. $\mT<\mS$: 
We have 
$(\delta_{\iota}^{(\bar{\rho}_{\mathtt{P}})})^{[\mS]}(\mathtt{P})\equiv 
(\delta_{\iota}^{[(\bar{\rho}(\mS))_{0}/\mS]})^{(\bar{\rho}_{\mathtt{P}}^{\mS})}$
and
$(\Delta^{(\bar{\rho}_{\sigma})})^{[\mS]}(\mathtt{P}^{\sigma})\equiv (\Delta^{[(\bar{\rho}(\mS))_{0}/\mS]})^{(\bar{\rho}_{\sigma}^{\mS})}$
for $\bar{\rho}=\bar{\rho}_{\mathtt{Q}}$ and  $\bar{\rho}_{\sigma}=\bar{\rho}_{\mathtt{P}^{\sigma}}$, where $\bar{\rho}_{\mathtt{P}}^{\mS}=\bar{\rho}_{\mathtt{P}}\restrict\mS$ and
$\bar{\rho}_{\sigma}^{\mS}=\bar{\rho}_{\sigma}\restrict\mS=\bar{\rho}_{\mathtt{P}}^{\mS}[\sigma/\rho]$.
IH yields
$(\mathcal{H}_{\gamma_{1}},\Theta_{\mS,\iota}, \mathtt{S}_{\iota}
) 
\vdash^{\omega a_{0}}_{\mS,\vec{\mS},\mS,\Lambda}
((\Gamma(\mathtt{P}))^{[\mS]}(\mathtt{P})\cup\{\lnot (\delta_{\iota}^{[(\bar{\rho}(\mS))_{0}/\mS]})^{(\bar{\rho}_{\mathtt{P}}^{\mS})}\})(\mathtt{S}_{\iota})$ 
for each $\iota\in[\bar{\rho}_{\mathtt{P}}^{\mS}\mathtt{P}]J_{\delta}$ and $\mathtt{S}_{\iota}=\mathtt{S}[\mathrm{srk}(\delta_{\iota}^{[\bar{\rho}(\mS)/\mS]})]$.
On the other hand we have
$(\mathcal{H}_{\gamma_{1}},(\Theta_{\mS})_{\mathtt{S}},\mathtt{S}^{\sigma})\vdash^{\omega a_{0}}_{\mS,\vec{\mS},\mS,\Lambda}(\Delta^{[(\bar{\rho}(\mS))_{0}/\mS]})^{(\bar{\rho}_{\sigma}^{\mS})} , (\Gamma(\mathtt{P})^{[\mS]}(\mathtt{P}^{\sigma}))(\mathtt{S})
$ for $\sigma\in H_{\rho_{\mU}}(f,\mathtt{S},\Theta_{\mS})\subset H_{\rho_{\mU}}(f,\mathtt{P},\Theta)$.
We have $\Theta\cap M_{\rho_{\mU}}\subset M_{\sigma}$ and $(\sigma,f)\not\in\mathtt{P}(\mU)$.

\bclm\label{clm:singlemainl_S_case2}
For $A^{(\bar{\tau})}\in\Gamma(\mathtt{P})$, let 
$(A^{(\bar{\tau})})^{[\mS]}(\mathtt{P})\equiv B^{(\bar{\tau}\,\restrict\,\mS)}$ and 
$(A^{(\bar{\tau})})^{[\mS]}(\mathtt{P}^{\sigma})\equiv C^{(\bar{\tau}\,\restrict\,\mS)}$.
Then $B\equiv C$.
\eclm
\textbf{Proof} of Claim \ref{clm:singlemainl_S_case2}.
We may assume $A\equiv D^{[\sigma/\mT]}\not\equiv D$ for a $D$.
Then $B\equiv (D^{[\sigma/\mT]})^{[\rho_{\mathtt{Q}(\mS)}/\mS]}\equiv D^{[\sigma/\mT]}\equiv C$
for $(A^{(\bar{\tau})})^{[\mS]}(\mathtt{P})\equiv B^{(\bar{\tau}\,\restrict\,\mS)}$
and $(A^{(\bar{\tau})})^{[\mS]}(\mathtt{P}^{\sigma})\equiv C^{(\bar{\tau}\,\restrict\,\mS)}$
by Definition \ref{df:restrict}.\ref{df:restrict.23}.
\hspace*{\fill} $\Box$ of Claim \ref{clm:singlemainl_S_case2}.
\\

We obtain $\Gamma(\mathtt{P})^{[\mS]}(\mathtt{P}^{\sigma})=\Gamma(\mathtt{P})^{[\mS]}(\mathtt{P})$
by Claim \ref{clm:singlemainl_S_case2}.
\\
(\ref{eq:singlemainl_S1_I_1}) follows by an $(i {\rm -rfl}_{\mT}(g,x,f,\Delta^{[(\bar{\rho}(\mS))_{0}/\mS]}))$ with 
$\mathrm{rk}(\Delta^{[(\bar{\rho}(\mS))_{0}/\mS]})<\mS=d_{i}$.
\\
\textbf{Case 2.2}. $\mT=\mS$:
Then $\mS\leq\mU\leq\mathrm{rk}(\Delta)<\mS+1$ means $\mU=\mS=\mathrm{rk}(\Delta)$.

Let $\rho=(\bar{\rho}(\mS))_{0}=\rho_{\mathtt{Q}(\mS)}$.
We have 
a finite set $\Delta\subset\bigvee(\mS+1)$
and an ordinal $a_{0}<a$ such that 
$(\mathcal{H}_{\gamma},\Theta_{\iota}, \mathtt{Q}_{\iota})\vdash^{a_{0}}_{\mS+1,\vec{\mS}+1,\mS^{\dagger},\Lambda}
(\Gamma\cup\{\lnot \delta_{\iota}^{(\bar{\rho})}\})(\mathtt{Q}_{\iota})$
for each $\delta\in\Delta$, $\iota\in[\bar{\rho}\mathtt{Q}]J_{\delta}$ with $\delta\simeq\bigvee(\delta_{\iota})_{\iota\in J_{\delta}}$,
and each good extension $\mathtt{Q}_{\iota}=\mathtt{Q}[\mathrm{srk}(\delta_{\iota})]$, where
$\Theta_{\iota}=\Theta\cup\mathsf{k}(\delta_{\iota}^{(\bar{\rho})};\mathtt{Q})$.
$f$ is a finite function such that 
$\mathrm{fld}(f)\subset\mathcal{H}_{\gamma}[\Theta(\mathtt{Q})]$.
On the other hand we have
\begin{equation}\label{eq:singlemainl_SS_I_1}
(\mathcal{H}_{\gamma},\Theta_{\mathtt{Q}},
\mathtt{ Q}^{\sigma}
)
\vdash^{a_{0}}_{\mS+1,\vec{\mS}+1,\mS^{\dagger},\Lambda}
\Delta^{(\bar{\rho}_{\sigma})}, \Gamma
\end{equation}
for every $\sigma\in H_{\rho}(f,\mathtt{Q},\Theta)$ and $\mathtt{Q}^{\sigma}=\mathtt{Q}\cup\{(\mS,(\sigma,f))\}$.

We obtain
$(\delta_{\iota}^{(\bar{\rho})})^{[\mS]}(\mathtt{Q})\equiv (\delta_{\iota}^{[\rho/\mS]})^{(\bar{\rho}^{\mS})}\equiv \delta_{\iota}^{(\bar{\rho}^{\mS})}$
by $\mathrm{rk}(\delta_{\iota})<\mS$. On the other hand we have
$(\Delta^{(\bar{\rho}_{\sigma})})^{[\mS]}(\mathtt{Q}^{\sigma})\equiv (\Delta^{[\sigma/\mS]})^{(\bar{\rho}^{\mS})}$ by Definition \ref{df:restrict}.\ref{df:restrict.23} for 
$\sigma=(\bar{\rho}_{\sigma}(\mS))_{0}=\rho_{\mathtt{Q}^{\sigma}(\mS)}$.

We have $\mathrm{fld}(f)\cup\{\mathtt{p}_{0}(\rho)\}\subset\mathcal{H}_{\gamma}[\Theta(\mathtt{Q})]$ 
for $\gamma\leq \gamma^{\mathtt{Q}}_{\mS}$ by (\ref{eq:controlder_cap2}) for
$\rho=\rho_{\mathtt{Q}(\mS)}$.
By Lemma  \ref{lem:exist_H}.\ref{lem:exist_H15} pick an ordinal $\sigma$ such that
$\mathtt{p}_{0}(\sigma)=\mathtt{p}_{0}(\rho)$,
$\{b(\sigma),\sigma\}\subset\mathcal{H}_{\gamma_{1}}[\Theta_{\mS}(\mathtt{R})]$, 
$m(\sigma)=f$, and
$\sigma\in H_{\rho}(f,\mathtt{Q},\Theta)$.

The following Claim \ref{clm:singlemainl_S_case22} is seen as in Claim \ref{clm:singlemainl_S_case2}.

\bclm\label{clm:singlemainl_S_case22}
For $A^{(\bar{\tau})}\in\Gamma$, let  $(A^{(\bar{\tau})})^{[\mS]}(\mathtt{Q})\equiv B^{(\bar{\tau}\,\restrict\,\mS)}$ and 
$(A^{(\bar{\tau})})^{[\mS]}(\mathtt{Q}^{\sigma})\equiv C^{(\bar{\tau}\,\restrict\,\mS)}$.
Then $B\equiv C$.
\eclm

IH with (\ref{eq:singlemainl_SS_I_1}) and Claim \ref{clm:singlemainl_S_case22} 
yields 
\[
(\mathcal{H}_{\gamma_{1}},(\Theta_{\mS})_{\mathtt{R}},
\mathtt{R}
)
\vdash^{\omega a_{0}}_{\mS,\vec{\mS},\mS,\Lambda}
(\Delta^{[\sigma/\mS]})^{(\bar{\rho}^{\mS})}, \Gamma^{[\mS]}(\mathtt{Q})
\]
where $\sigma\in\mathcal{H}_{\gamma_{1}}[\Theta_{\mS}(\mathtt{R})]$, $\widehat{a_{0}}<\hat{a}$, and
$\Gamma^{[\mS]}(\mathtt{Q}^{\sigma})=\Gamma^{[\mS]}(\mathtt{Q})$ by Claim \ref{clm:singlemainl_S_case22}.

Let $\delta\in\Delta$ with $\delta\simeq\bigvee(\delta_{\iota})_{\iota\in J_{\delta}}$.
We obtain $\delta^{[\sigma/\mS]}\simeq\bigvee(\delta_{\iota})_{\iota\in [\sigma]J_{\delta}}$.
Let $\iota\in [\sigma]J_{\delta}$ be such that $\mathsf{k}(\delta_{\iota})=\mathsf{k}(\delta_{\iota}^{(\bar{\rho}^{\mS})};\mathtt{R})\subset M_{\mathtt{R}}$.
IH then yields 
\[
(\mathcal{H}_{\gamma_{1}},\Theta_{\mS,\iota}, \mathtt{R}_{\iota})\vdash^{\omega a_{0}}_{\mS,\vec{\mS},\mS,\Lambda}
(\Gamma^{[\mS]}(\mathtt{Q})\cup\{\lnot \delta_{\iota}^{(\bar{\rho}^{\mS})}\})(\mathtt{R}_{\iota})
\]
where $\mathtt{R}_{\iota}=\mathtt{R}[\mathrm{srk}(\delta_{\iota})]$.
We have $\mathrm{rk}(\Delta^{[\sigma/\mS]})<\mS$ and
$\omega a_{0}+n<\omega a$ for any $n<\omega$.
A series of $(cut)$'s yields
(\ref{eq:singlemainl_S1_I_1}).
\\
\textbf{Case 3}. Third the last inference is a $(cut)$:
Let $C^{(\bar{\rho}_{\mathtt{Q}})}$ be a cut formula with $C\simeq\bigvee(C_{\iota})_{\iota\in J}$ and  $\mathrm{rk}(C)<\mS+1$.
We have an $a_{0}<a$ such that
$(\mathcal{H}_{\gamma},\Theta_{\iota},
\mathtt{Q}_{\iota})\vdash^{a_{0}}_{\mS+1,\vec{\mS}+1,\mS^{\dagger},\Lambda}
(\Gamma\cup\{\lnot C_{\iota}^{(\bar{\rho}_{\mathtt{Q}})}\})(\mathtt{Q}_{\iota})$ holds for each 
$\iota\in[\bar{\rho}_{\mathtt{Q}}\mathtt{Q}]J$ and each $\mathtt{Q}_{\iota}=\mathtt{Q}[\mathrm{srk}(C_{\iota})]$,
where
$\Theta_{\iota}=\Theta_{\mathtt{Q}_{\iota}}\cup\mathsf{k}(C_{\iota})$ with
$\mathsf{k}(C_{\iota})=\mathsf{k}(C_{\iota}^{(\bar{\rho}_{\mathtt{Q}})};\mathtt{Q})$.
On the other hand we have
$(\mathcal{H}_{\gamma},\Theta, \mathtt{Q})\vdash^{a_{0}}_{\mS+1,\vec{\mS}+1,\mS^{\dagger},\Lambda}C^{(\bar{\rho}_{\mathtt{Q}})},\Gamma$.
We obtain $(C^{(\bar{\rho}_{\mathtt{Q}})})^{[\mS]}\equiv 
(C^{[\rho_{\mathtt{Q}(\mS)}/\mS]})^{(\bar{\rho}_{\mathtt{R}})}$
and $(C_{\iota}^{(\bar{\rho}_{\mathtt{Q}})})^{[\mS]}\equiv 
(C_{\iota}^{[\rho_{\mathtt{Q}(\mS)}/\mS]})^{(\bar{\rho}_{\mathtt{R}})}\equiv
C_{\iota}^{(\bar{\rho}_{\mathtt{R}})}$.

Let $\iota\in[\rho_{\mathtt{Q}(\mS)}]J$ be such that 
$\mathsf{k}(C_{\iota})\subset M_{\mathtt{R}}$.
IH yields
$(\mathcal{H}_{\gamma_{1}},\Theta_{\mS,\iota},\mathtt{R}_{\iota}
)
\vdash^{\omega a_{0}}_{\mS,\vec{\mS},\mS,\Lambda}
(\Gamma^{[\mS]}\cup\{\lnot C_{\iota}^{(\bar{\rho}_{\mathtt{R}})}\})(\mathtt{R}_{\iota})$
for $\mathtt{R}_{\iota}=\mathtt{R}[\mathrm{srk}(C_{\iota})]$, and
$(\mathcal{H}_{\gamma_{1}},\Theta_{\mS}, \mathtt{R})\vdash^{\omega a_{0}}_{\mS,\vec{\mS},\mS,\Lambda}(C^{[\rho_{\mathtt{Q}(\mS)}/\mS]})^{(\bar{\rho}_{\mathtt{R}})},\Gamma^{[\mS]}$.
A $(cut)$ with $\mathrm{rk}(C^{[\rho_{\mathtt{Q}(\mS)}/\mS]})<\mS$ yields (\ref{eq:singlemainl_S1_I_1}).
\\
\textbf{Case 4}. Fourth  the last inference is a $(\mathrm{dom})$:
There are an ordinal $a_{0}<a$ and
a stable ordinal $\mT\not\in\mathrm{dom}(\mathtt{Q})$ such that 
$\mT\in\mathcal{H}_{\gamma}[\Theta(\mathtt{Q})]\cap\mS^{\dagger}$
for which
\[
\infer[(\mathrm{dom})]{(\mathcal{H}_{\gamma},\Theta,\mathtt{Q})\vdash^{a}_{\mS+1,\vec{\mS}+1,\mS^{\dagger},\Lambda}\Gamma}
{
(\mathcal{H}_{\gamma},\Theta_{\mathtt{P}},\mathtt{P})\vdash^{a_{0}}_{\mS+1,\vec{\mS}+1,\mS^{\dagger},\Lambda}\Gamma(\mathtt{P})
}
\]
holds for each good extension $\mathtt{P}=\mathtt{Q}[\mT]$,
where 
$\Gamma(\mathtt{P})=\{A^{(\bar{\rho}*(\rho_{\mathtt{P}(\mT)}))}: A^{(\bar{\rho})}\in\Gamma\}$.
We obtain $\mT<\mS$ by $\mT\not\in\mathrm{dom}(\mathtt{Q})\ni\mS$.

Given a good extension $\mathtt{R}_{1}=\mathtt{R}[\mT]$ of $\mathtt{R}$ for $\gamma_{1}$ and $\Theta_{\mS}$,
let $\mathtt{P}=\mathtt{R}_{1}[\mS]$ be the extension by the ordinal $\gamma^{\mathtt{Q}}_{\mS}$
and the set $\Theta$.
Then $\mathtt{R}_{1}=\mathtt{P}\restrict\mS$, and $\mathtt{P}=\mathtt{Q}[\mT]$ is a good extension of $\mathtt{Q}$ for $\gamma$ and $\Theta$ since
$\gamma<\gamma_{1}$ 
and $(\Theta\cup\mathtt{Q}(\mS))(\mathtt{R})=(\Theta_{\mS})(\mathtt{R})\subset M_{\rho_{\mathtt{R}_{1}(\mT)}}=M_{\rho_{\mathtt{P}(\mT)}}$.
In particular we obtain
$\mathtt{Q}(\mS)\subset M_{\rho_{\mathtt{R}_{1}(\mT)}}=M_{\mathtt{R}_{1}(\mT)}$
for $\mathtt{Q}(\mS)\subset M_{\mathtt{R}}$.

Let $A^{(\bar{\rho})}\in\Gamma$ and $(A^{(\bar{\rho})})^{[\mS]}(\mathtt{Q})\equiv B^{(\bar{\rho}\,\restrict\,\mS)}$.
We see $(A^{(\bar{\rho}*(\rho_{\mathtt{P}(\mT)}))})^{[\mS]}(\mathtt{P})\equiv B^{((\bar{\rho}*(\rho_{\mathtt{P}(\mT)}))\,\restrict\,\mS)}$, where
$(\bar{\rho}*(\rho_{\mathtt{P}(\mT)}))\restrict\mS=(\bar{\rho}\restrict\mS)*(\rho_{\mathtt{P}(\mT)})$ by  $\mT<\mS$.
IH followed by a $(\mathrm{dom})$ yields
\[
\infer[(\mathrm{dom})]{
(\mathcal{H}_{\gamma_{1}},\Theta_{\mS},\mathtt{R})\vdash^{\omega a}_{\mS,\vec{\mS},\mS,\Lambda}\Gamma^{[\mS]}(\mathtt{Q})
}
{
(\mathcal{H}_{\gamma_{1}},(\Theta_{\mS})_{\mathtt{R}_{1}},\mathtt{R}_{1})\vdash^{\omega a_{0}}_{\mS,\vec{\mS},\mS,\Lambda}
 (\Gamma(\mathtt{P})^{[\mS]})(\mathtt{P})
}
\]

The case when the last inference is a $(\Sigma(\Ome)\mbox{{\rm -rfl}})$ is seen from SIH.
\eprf

\blem\label{lem:main.2_S_I_1}{\rm (Elimination of one successor stable ordinal)}\\
Let $\mathtt{Q}$ be a finite family with thresholds $\gamma^{\mathtt{Q}}_{\cdot}, \delta^{\mathtt{Q}}_{\cdot}$
such that $\mT\in \mathrm{dom}(\mathtt{Q})\cap SSt$, and
$\Lambda^{(n)}>\Lambda_{0}$ strongly critical numbers such that
$\{
\Lambda^{(0)}<\Lambda^{(1)}<\Lambda^{(2)}<\cdots<\Lambda^{(3N-3)}<\Lambda^{(3N-2)}<\Lambda^{(3N-1)}
\}\subset\mathcal{H}_{\gamma}[\Theta(\mathtt{Q})]$.
Let 
Let $\vec{d}=(\mT^{\dagger},\ldots,\vec{\mT}^{\dagger})$ with $d=\mT^{\dagger}$.

Let
$
(\mathcal{H}_{\gamma},\Theta, \mathtt{Q}
)
\vdash^{a}_{\mT^{\dagger},\vec{d},\mT^{\dagger},\Lambda^{(0)}}
\Phi
$,
where $a<\Lambda^{(0)}$ and $\Phi$ is a set of $\bigvee$-formulas such that
$\mathrm{rk}(\Phi)<\mT+1$.

For $1\leq i\leq N$, let $\mS_{i},\kappa_{i},\rho_{i}$ be ordinals, $\mathtt{Q}_{i}$ finite families,
$\Phi_{i}$ sequents and ordinals $\alpha_{i}(\mathtt{Q})=\alpha(\mathtt{Q};\mS_{i})$ 
defined from ordinals $\mT$, $\alpha_{0}=\varphi_{\mT^{\dagger}}(a)$,
 the ordinals $\vec{d}$, and strongly critical numbers $\Lambda^{(n)}$
in Definition \ref{df:main_lemmas}.\ref{df:main_lemmas.2}.

Then 
\[
(\mathcal{H}_{\gamma_{1}},\Theta\cup\mathtt{Q}_{N}(\mT), \mathtt{Q}_{N}\restrict\mT
)
\vdash^{\alpha_{N}}_{\mT,\vec{\mT},\mT,\Lambda^{(3N-1)}}
(\Phi^{[\mT]}(\mathtt{Q}))^{(\vec{\kappa}/\vec{\rho})}
\]
holds for $\gamma_{1}=\gamma^{\mathtt{Q}}_{\mT}+\mT$.

\elem
\bprf 
By Cut-elimination \ref{lem:CE} we obtain
$ 
(\mathcal{H}_{\gamma},\Theta, \mathtt{Q})\vdash^{\alpha_{1}}_{\mT+1,\vec{d},\mT^{\dagger},\Lambda^{(0)}}\Phi
$ for $\alpha_{0}=\varphi_{\mT^{\dagger}}(a)$.
Corollary \ref{cor:singlemainl_I_1} then yields
$
(\mathcal{H}_{\gamma},\Theta, \mathtt{Q}_{N})\vdash^{\alpha_{N}}_{\mT+1,\vec{\mT}+1,\mT^{\dagger},\Lambda^{(3N-1)}}\Phi_{N}
$, where $\vec{\mT}+1=(\mT+1,\ldots,\mT+1)$.
We obtain
$
(\mathcal{H}_{\gamma_{1}},\Theta\cup\mathtt{Q}_{N}(\mT), \mathtt{Q}_{N}\restrict\mT)\vdash^{\alpha_{N}}_{\mT,\vec{\mT},\mT,\Lambda^{(3N-1)}}
(\Phi_{N})^{[\mT]}(\mathtt{Q}_{N})
$
by Lemma \ref{mlem:singlemainl_S_I_1}, where 
$\mathrm{rk}(\Phi_{N})\leq\mathrm{rk}(\Phi)$,
$\omega \alpha_{N}=\alpha_{N}$ and
$(\Phi_{N})^{[\mT]}(\mathtt{Q}_{N})=
(\Phi^{[\mT]}(\mathtt{Q}))^{(\vec{\kappa}/\vec{\rho})}$.
\eprf

\blem\label{lem:main.1_I_1}{\rm (Elimination of stable ordinals)}\\
Let $\xi\leq\Lambda_{0}$
be a stable ordinal, and
$\mathtt{Q}$
a finite family with thresholds $\gamma^{\mathtt{Q}}_{\cdot}, \delta^{\mathtt{Q}}_{\cdot}$.

Suppose
$
(\mathcal{H}_{\gamma},\Tht, \mathtt{Q})
\vdash^{a}_{\xi,\vec{\xi},\xi,\Lambda^{(k)}_{\xi}}
\Gamma
$,
where $\gamma=b_{0}+\omega^{\Lambda_{0}+1}$, $k<\omega$,
$a<\Lambda^{(k)}_{\xi}$ and
$\Gamma$ is a set of $\bigvee$-formulas.
Let $\mW\in St\cup\{0\}$ be an ordinal such that
$\mathrm{rk}(\Gamma)<\mW^{\dagger}<\xi$ and 
$\mW\in\mathrm{dom}(\mathtt{Q})\cap\mathcal{H}_{\gamma}[\Theta(\mathtt{Q})]$.

Let 
$f_{\mW}(\mathtt{Q},\xi,a)=\varphi_{\Lambda^{\mathtt{Q}}_{\mW}+\xi}(a)$ with
$\Lambda^{\mathtt{Q}}_{\mW}=\psi_{\mI_{N}}(b_{0}+\delta^{\mathtt{Q}}_{\mW})$ and
$\delta^{\mathtt{Q}}_{\mW}\leq\omega^{\Lambda_{0}+1}$, where
$\delta^{\mathtt{Q}}_{0}=\omega^{\Lambda_{0}+1}$,
$\Lambda^{\mathtt{Q}}_{0}=\Lambda^{(n)}_{0}=\psi_{\mI_{N}}(b_{0}+\delta^{\mathtt{Q}}_{0})$ and 
$\gamma^{\mathtt{Q}}_{0}=\max(\{\gamma\}\cup\mathsf{k}_{0}(\Theta))+\omega^{\Lambda_{0}\cdot 2}$ when $\mW=0$.
$\Lambda^{(n)}_{\mW}=\psi_{\mI_{N}}(b_{0}+\delta^{\mathtt{Q}}_{\mW}+n)$ for $\mW\neq 0$ and
$0<n<\omega$.

Let
$\{\mT\in \mathrm{dom}(\mathtt{Q}): \mT\geq\mW^{\dagger}\}=\{\mT_{0}>\mT_{1}>\cdots>\mT_{n-1}\}$,
and
$b_{m}=\mT_{m}+1$ if $\mT_{m}\in SSt$, and $b_{m}=\mT_{m}$ else.
Let $\mT_{n}=\mT_{n+1}=\mW$ and $b_{n}=\mW^{\dagger}$.

For $m\leq n$, let us introduce finite sets $\Theta_{m}$ of ordinals,
vectors $\vec{\kappa}_{m}$ of ordinals, and finite families $\mathtt{R}_{m}$ with
$\mathrm{dom}(\mathtt{R}_{m})\cap SSt_{i}=\mathrm{dom}(\mathtt{Q})\cap SSt_{i}\cap\mT_{m+1}^{\dagger}$
recursively as follows, 
cf.\,Definition \ref{df:main_lemmas}.\ref{df:main_lemmas.2}.
Let $\Theta_{-1}=\Theta$, $\vec{\kappa}_{-1}=\emptyset$, $\mathtt{R}_{-1}=\mathtt{Q}$
and $\Gamma_{-1}=\Gamma$
with $\mT_{-1}:=\mT_{0}^{\dagger}:=\xi^{\dagger}$.

Let $0\leq m\leq n$ and assume that $\Theta_{m-1}$, $\vec{\kappa}_{m-1}$ and $\mathtt{R}_{m-1}$
are given.
Let $\Theta_{m}=\Theta_{m-1}\cup (\mathtt{R}_{m-1}(\mT_{m}))_{0}
\cup\{\gamma^{\mathtt{Q}}_{\mT_{m}}\}$ for $m<n$, and $\Theta_{n}=\Theta_{n-1}$.
Let 
$1\leq i\leq N$ be such that
$\mathrm{dom}(\mathtt{Q})\cap SSt_{i}\cap\mT_{m}^{\dagger}\neq\emptyset$
and
$\mS_{i}=\max(\mathrm{dom}(\mathtt{Q})\cap SSt_{i}\cap\mT_{m}^{\dagger})$.
Let
$(\rho_{i},g_{i})\in\mathtt{R}_{m-1}(\mS_{i})$ with $\rho_{i}=\rho_{\mathtt{R}_{m-1}(\mS_{i})}$.
If $s(g_{i})\leq b_{m}$, then let $\kappa_{i}=\rho_{i}$ and $f_{i}=g_{i}$.
Assume $b_{m}<s(g_{i})$. 

For $1\leq j< i$, let $\alpha_{j}(\mathtt{R}_{m-1})=\Lambda_{\mT_{m}}^{(3j-2)}\cdot\omega^{\Lambda_{0}+1}$ if 
$\mathrm{dom}(\mathtt{Q})\cap SSt_{j}\cap\mT_{m}^{\dagger}=\emptyset$.
If
$\mathrm{dom}(\mathtt{Q})\cap SSt_{j}\cap\mT_{m}^{\dagger}\neq\emptyset$, then let
$\alpha_{j}(\mathtt{R}_{m-1})=\Lambda_{\mT_{m}}^{(3j-2)}\cdot\delta^{\mathtt{Q}}_{\mS_{j}}+on_{\Lambda_{\mT_{m}}^{(3j-3)}}(g_{j})$.
Define ordinals $\alpha_{j}$ recursively by
$\alpha_{0}=\varphi_{\mT_{m}^{\dagger}}(f_{\mT_{m}}(\mathtt{Q},\xi,a))$,
 and
$\alpha_{j}=\varphi_{\alpha_{j}(\mathtt{R}_{m-1})}(\alpha_{j-1})$.

Then let $f_{i}=h^{b_{m}}(g_{i};\varphi_{\mT_{m}^{\dagger}}(\alpha_{i-1}))$. and $\kappa_{i}\in H_{\rho_{i}}(f_{i},\mathtt{R}_{m-1},\Theta_{m-1})$.

$\vec{\rho}_{m}=(\rho_{i})_{i}$ and $\vec{\kappa}_{m}=(\kappa_{i})_{i}$ then denote
vectors of ordinals, where $i$ ranges over numbers such that
$1\leq i\leq N$ and $\mathrm{dom}(\mathtt{Q})\cap SSt_{i}\cap\mT_{m}^{\dagger}\neq\emptyset$.
Let $\mathtt{R}_{m}=(\mathtt{R}_{m-1}^{(\vec{\kappa}_{m}/\vec{\rho}_{m})})\restrict\mT_{m+1}^{\dagger}$.

Let
$A_{-1}\equiv A$ be a formula such that $\mathsf{k}(A^{(\bar{\rho})};\mathtt{Q})\subset M_{\mathtt{Q}}$,
and $\bar{\tau}_{-1}=\bar{\rho}$.
$(A^{(\bar{\rho})})_{m}\equiv A_{m}^{(\bar{\tau}_{m})}$ denotes a capped formula over $\mathtt{R}_{m}$ recursively defined as follows.
Let $\bar{\tau}_{m}=(\bar{\tau}_{m-1}^{(\vec{\kappa}_{m}/\vec{\rho}_{m})})\restrict\mT_{m+1}^{\dagger}$.

If $\mT_{m}\in SSt$ and $m<n$, then let 
$A_{m}^{(\bar{\tau}_{m})}\equiv \left(
(A_{m-1}^{(\bar{\tau}_{m-1})})^{[\mT_{m}]}(\mathtt{R}_{m-1})
\right)^{(\vec{\kappa}_{m}/\vec{\rho}_{m})}$.
Otherwise let $A_{m}^{(\bar{\tau}_{m})}\equiv (A_{m-1}^{(\bar{\tau}_{m-1})})^{(\vec{\kappa}_{m}/\vec{\rho}_{m})}$.
Let $\Gamma_{m}=\{(A^{(\bar{\rho})})_{m}  : A^{(\bar{\rho})}\in\Gamma\}$.

Finally let $\Theta_{\mW^{\dagger}*}=\Theta_{n}$,
$\mathtt{Q}_{\mW}=\mathtt{R}_{n}$, $A^{[\mW^{\dagger}*]}\equiv (A^{(\bar{\rho})})_{n}$ and 
 $\Gamma^{[\mW^{\dagger}*]}=\Gamma_{n}=\{A^{[\mW^{\dagger}*]}  : A^{(\bar{\rho})}\in\Gamma\}$.
Then 
\begin{equation}\label{eq:main.1_I_1}
(\mathcal{H}_{\gamma^{\mathtt{Q}}_{\mW}},\Theta_{\mW^{\dagger}*}, \mathtt{Q}_{\mW}
)
\vdash^{f_{\mW}(\mathtt{Q},\xi,a)}_{\mW^{\dagger},\vec{\eta}_{\mW},\mW^{\dagger},
\Lambda^{(3N-2)}_{\mW}}
\Gamma^{[\mW^{\dagger}*]}
\end{equation}
holds for $\vec{\eta}_{\mW}=(\mW^{\dagger},\ldots,\mW^{\dagger})$.

\elem
\bprf
By main induction on $\xi$ with subsidiary induction on $a$.
For $A\in\Gamma$ we have $\mathrm{rk}(A^{[\mW^{\dagger}*]})\leq\mathrm{rk}(A)<\mW^{\dagger}$.

We have
$\{\mW,a,\xi,\delta^{\mathtt{Q}}_{\mW}, \Lambda^{\mathtt{Q}}_{\mW}, \Lambda^{(3N-2)}_{\mW}
\}\subset\mathcal{H}_{\gamma}[\Theta(\mathtt{Q})]$
by the assumption, $b_{0}\in C_{0}(0)$, $\gamma\geq b_{0}+\delta^{\mathtt{Q}}_{\mW}$ and
(\ref{eq:controlder_cap22}),
and $\max\{\gamma,\gamma^{\mathtt{Q}}_{\mU}\}\leq\gamma^{\mathtt{Q}}_{\mW}$
for $\mW<\mU\in\mathrm{dom}(\mathtt{Q})$.
On the other hand we see inductively that
$\Theta(\mathtt{Q})\cup(\mathtt{R}_{m-1}(\mT_{m-1}))_{0}\subset\Theta_{m}(\mathtt{R}_{m})$ from
$\Theta(\mathtt{Q})\subset M_{\kappa_{i}}$.
Hence
$\{\mW^{\dagger},f_{\mW}(\mathtt{Q},\xi,a), \Lambda^{(3N-2)}_{\mW}\}\subset\mathcal{H}_{\gamma^{\mathtt{Q}}_{\mW}}[\Tht_{\mW^{\dagger}*}(\mathtt{Q}_{\mW})]$ follows.

We see $\gamma^{\mathtt{Q}}_{\mW}\leq\gam^{\mathtt{Q}}_{\mS}$
for every $\mS\in \mathrm{dom}(\mathtt{Q})\cap\mW^{\dagger}$, and
(\ref{eq:controlder_cap0}) is enjoyed.

Let $A^{(\bar{\rho})}\in\Gamma$.
$(A^{(\bar{\rho})})^{[\mW^{\dagger}*]}$ is a capped formula over $\mathtt{Q}_{\mW}$.
We see as in Claim \ref{clm:singlemainl_next_I_1} that
$\mathsf{k}((A^{(\bar{\rho})})^{[\mW^{\dagger}*]};\mathtt{Q}_{\mW})
\subset\mathcal{H}_{\gamma_{\mW}}[\Tht_{\mW^{\dagger}*}(\mathtt{Q}_{\mW})]$ from
 $\mT<\Lambda_{0}<\gamma^{\mathtt{Q}}_{\mW}$ and 
 $\bigcup\{(\mathtt{Q}(\mT))_{0}: \mW^{\dagger}\leq\mT\in \mathrm{dom}(\mathtt{Q})\}\subset\mathcal{H}_{\gamma^{\mathtt{Q}}_{\mW}}[\Tht_{\mW^{\dagger}*}(\mathtt{Q}_{\mW})]$.
 Hence (\ref{eq:controlder_cap2}) and (\ref{eq:controlder_cap22}) are enjoyed in (\ref{eq:main.1_I_1}).

(\ref{eq:controlder_cap_cover}) is enjoyed since
$\mathrm{srk}(\Gamma^{[\mW^{\dagger}*]})\subset\mathrm{dom}(\mathtt{Q}_{\mW})$
is seen as in Claim \ref{clm:singlemainl_next_I_1} by Proposition \ref{prp:restrict_I_1.2}.
We see 
$SC_{\mI_{N}}(\mathrm{fld}(f_{i}))<\Lambda_{\mW}^{(1)}$
from $\max\{\Lambda_{\xi}^{(k)},\Lambda_{\mT_{m}}^{(3N-4)}:m\leq n\}<\Lambda_{\mW}^{(3N-2)}$, cf.\,(\ref{eq:delta_Lambda}).
(\ref{eq:fld}) is enjoyed.

Let $b\in\mathcal{H}[\Theta(\mathtt{Q})]\cap a$, and
$\beta_{i}=\varphi_{\alpha_{i}(\mathtt{Q}_{m-1})}(\beta_{i-1})$
denote ordinals defined from
$\beta_{0}=\varphi_{\mT_{m}^{\dagger}}(f_{\mT_{m}}(\mathtt{Q},\xi,b))<\varphi_{\mT_{m}^{\dagger}}(f_{\mT_{m}}(\mathtt{Q},\xi,a))$.
$\mathtt{R}_{m}^{b}$ [$\Theta_{m}^{b}$] denotes finite families [sets of ordinals]
defined from $b$, resp.
$k_{i}=h^{b_{m}}(g_{i};\varphi_{\mT_{m}^{\dagger}}(\beta_{i-1}))\leq h^{b_{m}}(g_{i};\varphi_{\mT_{m}^{\dagger}}(\alpha_{i-1}))=h_{i}$ follows.
Hence if $\kappa_{i}\in H_{\rho_{i}}(f_{i},\mathtt{R}_{m-1},\Theta_{m-1})$, then
$\kappa_{i}\in H_{\rho_{i}}(k_{i},\mathtt{R}_{m-1},\Theta_{m-1})$, and we may assume 
$\mathtt{R}_{m}=\mathtt{R}_{m}^{b}$ and $\Theta_{m}=\Theta_{m}^{b}$.
\\
\textbf{Case 1}. 
First consider the case when the last inference is a $(\mathrm{dom})$:
We have an ordinal $a_{0}<a$, a number $k<\omega$ and
a stable ordinal $\mS\not\in \mathrm{dom}(\mathtt{Q})$ such that 
$\mS\in\mathcal{H}_{\gamma}[\Theta(\mathtt{Q})]\cap\xi$
for which
$(\mathcal{H}_{\gamma},\Theta_{\mathtt{P}}, \mathtt{P})\vdash^{a_{0}}_{\xi,\vec{\xi},\xi,\Lambda^{(k)}_{\xi}}\Gamma(\mathtt{P})$ 
holds for each good extension $\mathtt{P}=\mathtt{Q}[\mS]$ by 
ordinals $\gamma^{\mathtt{P}}_{\mS}$, $\delta^{\mathtt{P}}_{\mS}$ and 
$\mathtt{P}(\mS)=\{(\rho_{\mS},g_{\mS})\}$.
One of these extensions is very good.
We have $\mW\geq\mathrm{srk}(\mathrm{rk}(\Gamma))=\max(\{0\}\cup\{\mathrm{srk}(A^{(\bar{\rho})}): A^{(\bar{\rho})}\in\Gamma\})$.
Let $\mS\in SSt_{i}$.
\\
\textbf{Case 1.1}. $\mS\leq\mW$: We have $\mS\not\in\mathrm{dom}(\mathtt{Q})$ and 
$\mW\in\mathrm{dom}(\mathtt{Q})$.
We obtain $\mS<\mW$, $\gamma^{\mathtt{P}}_{\mW}=\gamma^{\mathtt{Q}}_{\mW}$  and
$\Lambda^{\mathtt{P}}_{\mW}=\Lambda^{\mathtt{Q}}_{\mW}$
since $\mT\in\mathrm{dom}(\mathtt{P})\Lrarw\mT\in\mathrm{dom}(\mathtt{Q})$ for any $\mT\geq\mW$.
Let $a_{1}=f_{\mW}(\mathtt{P},\xi,a_{0})=f_{\mW}(\mathtt{Q},\xi,a_{0})$.

Let $\Theta_{m}^{\mathtt{P}}$, $\vec{\kappa}_{m}^{\mathtt{P}}$, $\mathtt{R}_{m}^{\mathtt{P}}$
and $(A^{(\bar{\rho})})_{m}^{\mathtt{P}}\, (A^{(\bar{\rho})}\in\Gamma)$
denote sets, vectors of ordinals, finite families and formulas defined from $\mathtt{P}$, resp.
We may assume $\Theta_{m}^{\mathtt{P}}=\Theta_{m}$ by $\mT_{n-1}>\mS$, and
$(\Theta_{\mathtt{P}})_{\mW^{\dagger}*}=(\Theta_{\mW^{\dagger}*})_{\mathtt{P}}=
\Theta_{\mW^{\dagger}*}\cup\{\delta^{\mathtt{P}}_{\mS}\}\cup SC_{\mI_{N}}(\mathrm{fld}(g_{\mS}))$.

SIH yields
$(\mathcal{H}_{\gamma^{\mathtt{Q}}_{\mW}},(\Theta_{\mW^{\dagger}*})_{\mathtt{P}},
 \mathtt{P}_{\mW})\vdash^{a_{1}}_{\mW^{\dagger},\eta_{\mW},\mW^{\dagger},\Lambda^{(3N-2)}_{\mW}}
( \Gamma^{[\mW^{\dagger}*]}(\mathtt{P}))^{\mathtt{P}}$.

For $m<n$, let
$C(m):\Lrarw \left(
\mS=\max(\mathrm{dom}(\mathtt{P})\cap SSt_{i}\cap\mT_{m}^{\dagger}) \spand 
b_{m}<s(g_{\mS})
\right)$ and
$D(m):\Lrarw \left(
C(m)\spand \mathrm{dom}(\mathtt{Q})\cap SSt_{i}\cap\mT_{m}^{\dagger}\neq\emptyset
\spand b_{m}<s(g_{i})
\right)$, where
$(\rho_{i},g_{i})\in \mathtt{Q}(\mS_{i})$ with $\rho_{i}=\rho_{\mathtt{Q}(\mS_{i})}$ for
$\mS_{i}=\max(\mathrm{dom}(\mathtt{Q})\cap SSt_{i}\cap\mT_{m}^{\dagger})$.
\\
\textbf{Case 1.1.1}. $\lnot D(m)$ for any $m<n$:
If there is no $m$ such that $C(m)$, then we see 
$\mathtt{P}_{\mW}(\mS)=\mathtt{R}^{\mathtt{P}}_{n-1}(\mS)=\{(\rho_{\mS},g_{\mS})\}=\mathtt{P}(\mS)$ and 
$(\Gamma^{[\mW^{\dagger}*]}(\mathtt{P}))^{\mathtt{P}}=\Gamma^{[\mW^{\dagger}*]}(\mathtt{P})=\Gamma^{[\mW^{\dagger}*]}(\mathtt{P}_{\mW})$.

Suppose that there is an $m$ such that $C(m)$.
Then $\{(\rho_{\mS},g_{\mS})\}=\mathtt{P}(\mS)$ is replaced by
a $\{(\kappa_{\mS},f_{\mS})\}$ in such a way that
$\mathtt{P}_{\mW}(\mS)=\mathtt{R}^{\mathtt{P}}_{n-1}(\mS)=\{(\kappa_{\mS},f_{\mS})\}$.
If there is no $m$ such that $D(m)$, then we see $\mathtt{P}_{\mW}(\mT)=\mathtt{Q}_{\mW}(\mT)$
for any $\mT\neq\mS$, and
$(\Gamma^{[\mW^{\dagger}*]})^{\mathtt{P}}=\Gamma^{[\mW^{\dagger}*]}(\mathtt{P}_{\mW})$.

In each case, 
given a good extension $\mathtt{P}_{\mW}=\mathtt{Q}_{\mW}[\mS]$ of $\mathtt{Q}_{\mW}$ by 
ordinals $\gamma^{\mathtt{P}}_{\mS}$, $\delta^{\mathtt{P}}_{\mS}$ and 
$\mathtt{P}_{\mW}(\mS)$, let $\mathtt{P}=\mathtt{Q}[\mS]$ be the extension by
ordinals $\gamma^{\mathtt{P}}_{\mS}$, $\delta^{\mathtt{P}}_{\mS}$ and 
$\mathtt{P}_{\mW}(\mS)$.
Then $\mathtt{P}$ is a good extension of $\mathtt{Q}$.
Let $\mathtt{P}$ be a very good extension of $\mathtt{Q}$ for $\gamma,\Theta$.
Then we see from $\gamma\leq\gamma^{\mathtt{Q}}_{\mW}$ and
$\Theta(\mathtt{Q})\subset\Theta_{\mW^{\dagger}*}(\mathtt{Q}_{\mW})$
that $\mathtt{P}_{\mW}$ is a very good extension of $\mathtt{Q}_{\mW}$ for
$\gamma^{\mathtt{Q}}_{\mW},\Theta_{\mW^{\dagger}*}$.

A $(\mathrm{dom})$ yields
$(\mathcal{H}_{\gamma^{\mathtt{Q}}_{\mW}},\Theta_{\mW^{\dagger}*}, \mathtt{Q}_{\mW})\vdash^{a_{1}+1}_{\mW^{\dagger},\eta_{\mW},\mW^{\dagger},\Lambda^{(3N-2)}_{\mW}}\Gamma^{[\mW*]}$.
(\ref{eq:main.1_I_1}) follows.
\\
\textbf{Case 1.1.2}.
$D(m)$ for an $m<n$:
Then $\mS_{i}<\mS$, 
$\mathtt{P}_{\mW}(\mT)=\mathtt{Q_{\mW}}(\mT)$ holds for any $\mT\not\in \{\mS,\mS_{i}\}$, and
$\mathtt{P}_{\mW}(\mS_{i})=\mathtt{Q}(\mS_{i})\neq \mathtt{Q_{\mW}}(\mS_{i})$, where
$\mS_{i}=\max(\mathrm{dom}(\mathtt{Q})\cap SSt_{i}\cap\mT_{m}^{\dagger})$.

For $A^{(\bar{\rho})}\in\Gamma$, let
$(A^{[\mW^{\dagger}*]}(\mathtt{P}))^{\mathtt{P}}\equiv B^{(\bar{\lambda})}\in
(\Gamma^{[\mW^{\dagger}*]}(\mathtt{P}))^{\mathtt{P}}$ and
$A^{[\mW^{\dagger}*]}\equiv A_{n-1}^{(\bar{\tau})}\in \Gamma^{[\mW^{\dagger}*]}$.
We have $\bar{\lambda}(\mS)=(\kappa_{\mS}, f_{\mS})$, 
$\bar{\lambda}(\mS_{i})=\bar{\rho}(\mS_{i})\neq\bar{\tau}(\mS_{i})$ and
$\bar{\lambda}(\mT)=\bar{\tau}(\mT)$ for any $\mT\not\in\{\mS,\mS_{i}\}$.
We see $B\equiv A_{n-1}$ from $\mT_{m}\not\in\{\mS,\mS_{i}\}$ and
$\mathtt{R}_{m-1}^{\mathtt{P}}(\mT_{m})=\mathtt{R}_{m-1}(\mT_{m})$ for each $m<n$.
Let $m_{0}=\min\{m<n: D(m)\}$ and $m\geq m_{0}$.

Let $(\kappa_{m_{0}-1},f_{m_{0}-1})=(\rho_{i},g_{i})\in \mathtt{Q}(\mS_{i})$ and,
$(\kappa_{m},f_{m})\in\mathtt{R}_{m}(\mS_{i})$ with $\kappa_{m}=\rho_{\mathtt{R}_{m}(\mS_{i})}$.
Then $\mathtt{R}_{m}(\mS_{i})=(\mathtt{R}_{m}^{\mathtt{P}}(\mS_{i}))^{((\kappa_{m},f_{m})/(\kappa_{m_{0}-1},f_{m_{0}-1}))}$, where
$\mathtt{R}_{m}^{\mathtt{P}}(\mS_{i})=\mathtt{Q}(\mS_{i})$.
Let
$\mathtt{S}_{m}=\mathtt{Q}_{\mW}^{((\kappa_{m-1},f_{m-1})/(\kappa_{n-1},f_{n-1}))}$
and
$\Phi_{m}=(\Gamma^{[\mW*]})^{((\kappa_{m-1},f_{m-1})/(\kappa_{n-1},f_{n-1}))}$
for $m\geq m_{0}$.
$(\mathcal{H}_{\gamma^{\mathtt{Q}}_{\mW}},\Theta_{\mW^{\dagger}*}, 
\mathtt{S}_{m_{0}}
)\vdash^{a_{1}+1}_{\mW^{\dagger},\eta_{\mW},\mW^{\dagger},\Lambda^{(3N-2)}_{\mW}}
\Phi_{m_{0}}$
follows by a $(\mathrm{dom})$.

We have $s(f_{m_{0}-1})=s(g_{i})>b_{m_{0}}\geq\mW^{\dagger}$ and $s(f_{m-1})\geq b_{m-1}>b_{m}\geq\mW^{\dagger}$ for $m>m_{0}$.
Let $\alpha_{m_{0}}=a_{1}+1$ and $\alpha_{m}=\varphi_{\mW^{\dagger}}(\alpha_{m-1})$ for $m>m_{0}$.
We see inductively from
Recapping \ref{lem:recapping} that
$(\mathcal{H}_{\gamma^{\mathtt{Q}}_{\mW}},\Theta_{\mW^{\dagger}*}, 
\mathtt{S}_{m}
)\vdash^{\alpha_{m}}_{\mW^{\dagger},\eta_{\mW},\mW^{\dagger},\Lambda^{(3N-2)}_{\mW}}
\Phi_{m}$, and
$(\mathcal{H}_{\gamma^{\mathtt{Q}}_{\mW}},\Theta_{\mW^{\dagger}*}, 
\mathtt{S}_{n}
)\vdash^{\alpha_{n}}_{\mW^{\dagger},\eta_{\mW},\mW^{\dagger},\Lambda^{(3N-2)}_{\mW}}
\Phi_{n}$, where $\mathtt{S}_{n}=\mathtt{Q}_{\mW}$ and $\Phi_{n}=\Gamma^{[\mW*]}$.
We see $\alpha_{n}<f_{\mW}(\mathtt{Q},\xi,a)$ inductively from $a_{1}<f_{\mW}(\mathtt{Q},\xi,a)$ 
and $\mW^{\dagger}<\Lambda^{\mathtt{Q}}_{\mW}+\xi$.
(\ref{eq:main.1_I_1}) follows.
\\
\textbf{Case 1.2}. $\mW<\mS<\xi$: 
We have $\Lambda^{\mathtt{P}}_{\mS}<\Lambda^{\mathtt{P}}_{\mW}=\Lambda^{\mathtt{Q}}_{\mW}$
by $\{\mS>\mW\}\subset\mathrm{dom}(\mathtt{P})$.
We see $a_{2}=f_{\mS}(\mathtt{P},\xi,a_{0})=\varphi_{\Lambda^{\mathtt{P}}_{\mS}+\xi}(a_{0})<
\varphi_{\Lambda^{\mathtt{Q}}_{\mW}+\xi}(a)=f_{\mW}(\mathtt{Q},\xi,a)$
from $\Lambda^{\mathtt{P}}_{\mS}<\Lambda^{\mathtt{Q}}_{\mW}$ and
$a_{0}<a$.

\beqn\label{eq:main.1_I_1_case1_SIH}
(\mathcal{H}_{\gamma^{\mathtt{P}}_{\mS}},(\Theta_{\mS^{\dagger}*})_{\mathtt{P}},
 \mathtt{P}_{\mS})\vdash^{a_{2}}_{\mS^{\dagger},\eta_{\mS},\mS^{\dagger},\Lambda^{(3N-2)}_{\mS}}
( \Gamma^{[\mS^{\dagger}*]})^{\mathtt{P}}
\eeqn follows by SIH 
as in \textbf{Case 1.1}.

Let $\mT_{m-1}>\mS>\mT_{m}$ with 
$0\leq m\leq n$, $\mT_{-1}=\xi$ and
$\mT_{n}=\mW$.
Let us examine how $\mathtt{R}_{m}$ [$\mathtt{P}_{\mS}$] is obtained from $\mathtt{R}_{m-1}$
[from $\mathtt{R}_{m-1}^{\mathtt{P}}$], resp.

We obtain $(\Theta_{\mS^{\dagger}*})_{\mathtt{P}}=\Theta_{m-1}\cup\{\delta^{\mathtt{P}}_{\mS}\}
\cup SC_{\mI_{N}}(\mathrm{fld}(g_{\mS}))$.
Let $\mT_{m}\geq\mT\in\mathrm{dom}(\mathtt{Q})$ and
$\mS_{j}=\max(\mathrm{dom}(\mathtt{Q})\cap SSt_{j}\cap\mT_{m}^{\dagger})$ if
$\mathrm{dom}(\mathtt{Q})\cap SSt_{j}\cap\mT_{m}^{\dagger}\neq\emptyset$.
Let $\mU_{j}=\mS_{j}$ for $j\neq i$, and $\mU_{i}=\mS$.

We have $\mathtt{R}_{m-1}^{\mathtt{P}}(\mT)=\mathtt{R}_{m-1}(\mT)$ for $\mT\neq\mS_{i}$, and 
$\mathtt{P}_{\mS}(\mT)=\mathtt{R}_{m-1}^{\mathtt{P}}(\mT)=\mathtt{R}_{m-1}(\mT)$ for $\mT\not\in\{\mS_{j}:1\leq j\leq N\}$, while
$\mathtt{P}_{\mS}(\mS_{i})=\mathtt{Q}(\mS_{i})$.

Let
$(\rho_{j},g_{j})\in \mathtt{R}_{m-1}(\mS_{j})$ 
with $\rho_{j}=\rho_{\mathtt{R}_{m-1}(\mS_{j})}$, where $\mathtt{R}_{m-1}^{\mathtt{P}}(\mS_{j})=\mathtt{R}_{m-1}(\mS_{j})$ for $j\neq i$.
Let $(\kappa_{j},f_{j})=(\rho_{j},g_{j})$ if $s(g_{j})\leq b_{m}$.
Let $s(g_{j})> b_{m}$. Then let $\alpha_{0}=\varphi_{\mT_{m}^{\dagger}}(f_{\mT_{m}}(\mathtt{Q},\xi,a))$,
$\alpha_{k}(\mathtt{R}_{m-1})=\Lambda_{\mT_{m}}^{(3k-2)}\cdot\delta^{\mathtt{Q}}_{\mS_{k}}+on_{\Lambda_{\mT_{m}}^{(3k-3)}}(g_{k})$ and
$\alpha_{k}=\varphi_{\alpha_{k}(\mathtt{R}_{m-1})}(\alpha_{k-1})$.
Then $f_{j}=h^{b_{m}}(g_{j};\varphi_{\mT_{m}^{\dagger}}(\alpha_{j-1}))$ and 
$\kappa_{j}\in H_{\rho_{j}}(f_{j},\mathtt{R}_{m-1},\Theta_{m-1})$.
Let $\vec{\rho}=\{(\rho_{j},g_{j})\}_{j}$ and $\vec{\kappa}=\{(\kappa_{j},f_{j})\}_{j}$.

Let $j\neq i$.
If $s(g_{j})\leq \mS^{\dagger}$, then let
$(\kappa^{\mathtt{P}}_{j},f^{\mathtt{P}}_{j})=(\rho_{j},g_{j})$.
Let $s(g_{j})> \mS^{\dagger}$, $\beta_{0}=\varphi_{\mS^{\dagger}}(a_{2})$,
$\alpha_{k}(\mathtt{P}_{\mS})=\alpha_{k}(\mathtt{R}_{m-1}^{\mathtt{P}})=\Lambda_{\mS}^{(3k-2)}\cdot\delta^{\mathtt{P}}_{\mU_{k}}+on_{\Lambda_{\mS}^{(3k-3)}}(g_{k})$ and
$\beta_{k}=\varphi_{\alpha_{k}(\mathtt{R}_{m-1}^{\mathtt{P}})}(\beta_{k-1})$.
Then $f^{\mathtt{P}}_{j}=h^{\mS^{\dagger}}(g_{j};\varphi_{\mS^{\dagger}}(\beta_{j-1}))$ and 
$\kappa^{\mathtt{P}}_{j}\in H_{\rho_{j}}(f^{\mathtt{P}}_{j},\mathtt{R}^{\mathtt{P}}_{m-1},\Theta^{\mathtt{P}}_{m-1})$.
Let 
$\vec{\kappa}_{\mathtt{P}}=\{(\kappa^{\mathtt{P}}_{j},f^{\mathtt{P}}_{j})\}_{j\neq i}\cup\{(\rho_{i},g_{i})\}$.
We obtain 
\beqn\label{eq:main.1_I_1_case1_0}
\mathtt{P}_{\mS}\restrict\mS=\mathtt{R}_{m-1}^{(\vec{\kappa}_{\mathtt{P}}/\vec{\rho})}
\spand \mathtt{R}_{m-1}^{(\vec{\kappa}/\vec{\rho})}\restrict\mT_{m}=\mathtt{R}_{m}
\eeqn

Let $\Phi=(\Gamma^{[\mS^{\dagger}*]})^{\mathtt{P}}$.
We obtain
\beqn\label{eq:main.1_I_1_case1_1}
(\mathcal{H}_{\gamma^{\mathtt{P}}_{\mS}+\mS},
\Theta(\mS),
 \mathtt{P}_{\mS}\restrict\mS)\vdash^{\beta_{N}}_{\mS,\vec{\mS},\mS,\Lambda^{(6N-3)}_{\mS}}
 \Phi_{\mS}
 \eeqn
 by (\ref{eq:main.1_I_1_case1_SIH}), Lemmas \ref{lem:main.2_L_I_1} and \ref{lem:main.2_S_I_1}, where
 $\Theta(\mS)=(\Theta_{\mS^{\dagger}*})_{\mathtt{P}}\cup (\mathtt{P}_{\mS}(\mS))_{0}$, and
 $\Phi_{\mS}=(\Phi^{[\mS]}(\mathtt{P}_{\mS}))^{(\vec{\kappa}_{\mathtt{P}}/\vec{\rho})}$
 if $\mS\in SSt$, and
$\Phi_{\mS}=\Phi^{(\vec{\kappa}_{\mathtt{P}}/\vec{\rho})}$ else.

\bclm\label{clm:main.1_I_1_1}
$\Phi^{[\mS]}(\mathtt{P}_{\mS})=\Phi$ when $\mS\in SSt$.
\eclm
\textbf{Proof} of Claim \ref{clm:main.1_I_1_1}.
For $A^{(\bar{\tau})}\in\Phi$, suppose $(A^{(\bar{\tau})})^{[\mS]}(\mathtt{P}_{\mS})
\not\equiv A^{(\bar{\tau})}$, cf.\,Definition \ref{df:restrict}.\ref{df:restrict.2}.
First let 
$A\equiv B^{[\sigma/\mS]}\not\equiv B\equiv\mathrm{inv}^{(\mathtt{u})}(A^{(\bar{\tau})};\mathtt{P}_{\mS})$ with $\rho_{\mathtt{P}_{\mS}(\mS)}\leq\sigma< (\bar{\tau}(\mS))_{0}$.
Then we would have $\{\rho_{\mathtt{P}_{\mS}(\mS)}< (\bar{\tau}(\mS))_{0}\}\subset(\mathtt{P}_{\mS}(\mS))_{0}$.
On the other, since $\mathtt{P}(\mS)=\{(\rho_{\mS},g_{\mS})\}$ is a singleton,
so is the set $\mathtt{P}_{\mS}(\mS)$.

Otherwise we have
$(A^{(\bar{\tau})})^{[\mS]}(\mathtt{P}_{\mS})\equiv 
(A^{[\rho/\mS]})^{(\bar{\tau} \, \restrict \,\mS)}$ for $\rho=\rho_{\mathtt{P}_{\mS}(\mS)}$.
On the other hand we have
$\mathrm{rk}(A)\leq\mathrm{rk}((\Gamma^{[\mS^{\dagger}*]})^{\mathtt{P}})\leq\mathrm{rk}(\Gamma)<\mW^{\dagger}\leq\mS$ and $\mathsf{k}(A)\subset M_{\rho}$ for
$A^{(\bar{\tau})}\in\Phi=(\Gamma^{[\mS^{\dagger}*]})^{\mathtt{P}}$.
Hence $A^{[\rho/\mS]}\equiv A$, and
$\Phi_{\mS}=\Phi^{(\vec{\kappa}_{\mathtt{P}}/\vec{\rho})}=
((\Gamma^{[\mS^{\dagger}*]})^{\mathtt{P}})^{(\vec{\kappa}_{\mathtt{P}}/\vec{\rho})}
=\Gamma_{m-1}^{(\vec{\kappa}_{\mathtt{P}}/\vec{\rho})}$.
\hspace*{\fill} $\Box$ of Claim \ref{clm:main.1_I_1_1}.
\\

$(\mathcal{H}_{\gamma^{\mathtt{P}}_{\mS}+\mS},
\Theta(\mS),
 \mathtt{R}_{m-1}^{(\vec{\kappa}_{\mathtt{P}}/\vec{\rho})}
 )\vdash^{\beta_{N}}_{\mS,\vec{\mS},\mS,\Lambda^{(6N-3)}_{\mS}}
\Gamma_{m-1}^{(\vec{\kappa}_{\mathtt{P}}/\vec{\rho})}$ 
follows by (\ref{eq:main.1_I_1_case1_0}),
(\ref{eq:main.1_I_1_case1_1}) and
Claim \ref{clm:main.1_I_1_1}.

Now let $\mathtt{P}=\mathtt{Q}[\mS]$ be a very good extension.
Let $\mathtt{Q}_{1}=\mathtt{Q}\restrict\mS$.
We have 
$\gamma^{\mathtt{P}}_{\mS}\in\mathcal{H}_{\gamma}[\Theta_{\gamma,\mS}(\mathtt{Q}_{1})]$, 
$\{\delta^{\mathtt{P}}_{\mS}\}\cup SC_{\mI_{N}}(\mathrm{fld}(g_{\mS}))\subset
\mathcal{H}_{\gamma}[\Theta_{\delta,\mS}(\mathtt{Q})]$ and
$\rho_{\mS}\in\mathcal{H}_{\gamma^{\mathtt{P}}_{\mS}+\mS}[\Theta_{\gamma,\mS}(\mathtt{Q}_{1})]$ with $\{(\rho_{\mS},g_{\mS})\}=\mathtt{P}(\mS)$,
where 
$\Theta_{\gamma,\mS}=\Theta\cup
\{\gamma^{\mathtt{Q}}_{\mU}: \mS<\mU\in\mathrm{dom}(\mathtt{Q})\}$ and 
$\Theta_{\delta,\mS}=\Theta\cup\{\delta^{\mathtt{Q}}_{\mU}: \mS<\mU\in\mathrm{dom}(\mathtt{Q})\}$.
We obtain $\{\delta^{\mathtt{Q}}_{\mU}: \mS<\mU\in\mathrm{dom}(\mathtt{Q})\}\subset\mathcal{H}_{\gamma}[\Theta(\mathtt{Q})]$ by (\ref{eq:controlder_cap22}), and
$\{\delta^{\mathtt{P}}_{\mS}\}\cup SC_{\mI_{N}}(\mathrm{fld}(g_{\mS}))\subset
\mathcal{H}_{\gamma}[\Theta(\mathtt{Q})]$.

We have $\{\gamma^{\mathtt{Q}}_{\mU}: \mS<\mU\in\mathrm{dom}(\mathtt{Q})\}\subset\Theta_{m-1}$, and
$\gamma^{\mathtt{Q}}_{\mU}\in M_{\mathtt{Q}_{1}}$
for $\mS<\mU\in\mathrm{dom}(\mathtt{Q})$.
Hence
$\Theta_{\gamma,\mS}(\mathtt{Q}_{1})\subset\Theta_{m-1}(\mathtt{Q}_{1})$, and
$\rho_{\mS}\in
\mathcal{H}_{\gamma^{\mathtt{P}}_{\mS}+\mS}[\Theta_{m-1}(\mathtt{Q}_{1})]$.
Let $\mathtt{P}_{\mS}(\mS)=\{(\kappa_{\mS},f_{\mS}\}\}$.
As in \textbf{Case 1.1.1} we see that $(\kappa_{\mS},f_{\mS})=(\rho_{\mS},g_{\mS})$ if 
there is no $m$ such that $C(m)$.
Otherwise $(\kappa_{\mS},f_{\mS})$ is obtained from $(\rho_{\mS},g_{\mS})$ by
applying steps $(\rho,g)\mapsto(\kappa,f)$ finitely many times,
where $f=h^{b_{m}}(g;\varphi_{\mT_{m}^{\dagger}}(b))$ and $\kappa\in H_{\rho}(f,\mathtt{R}_{m-1},\Theta_{m-1})$.
By Lemma \ref{lem:exist_H}.\ref{lem:exist_H2} and $\rho_{\mS}\in
\mathcal{H}_{\gamma^{\mathtt{P}}_{\mS}+\mS}[\Theta_{m-1}(\mathtt{Q}_{1})]$,
we may assume that
$\kappa_{\mS}\in \mathcal{H}_{\gamma^{\mathtt{P}}_{\mS}+\mS}[\Theta_{m-1}(\mathtt{Q}_{1})]$.
Therefore 
$\mathcal{H}_{\gamma^{\mathtt{P}}_{\mS}+\mS}[(\Theta(\mS))(\mathtt{P}_{\mS}\restrict\mS)]=
\mathcal{H}_{\gamma^{\mathtt{P}}_{\mS}+\mS}[\Theta_{m-1}(\mathtt{P}_{\mS}\restrict\mS)]$
for $\Theta(\mS)=\Theta\cup\{\delta^{\mathtt{P}}_{\mS}\}\cup SC_{\mI_{N}}(\mathrm{fld}(g_{\mS}))\cup (\mathtt{P}_{\mS}(\mS))_{0}$, and we obtain
\[
(\mathcal{H}_{\gamma^{\mathtt{P}}_{\mS}+\mS},
\Theta_{m-1},
 \mathtt{R}_{m-1}^{(\vec{\kappa}_{\mathtt{P}}/\vec{\rho})}
 )\vdash^{\beta_{N}}_{\mS,\vec{\mS},\mS,\Lambda^{(6N-3)}_{\mS}}
\Gamma_{m-1}^{(\vec{\kappa}_{\mathtt{P}}/\vec{\rho})}
\]
We see
$\beta_{N}<\varphi_{\Lambda^{\mathtt{Q}}_{\mW}+\xi}(a)=f_{\mW}(\mathtt{Q},\xi,a)$ from
$a_{2}<f_{\mW}(\mathtt{Q},\xi,a)$ and $\Lambda_{\mS}^{(k)}<\Lambda^{\mathtt{Q}}_{\mW}$.
We obtain
$a_{3}=f_{\mT_{m}}( \mathtt{R}_{m-1}^{(\vec{\kappa}_{\mathtt{P}}/\vec{\rho})},\mS,\beta_{N})
=\varphi_{\Lambda^{\mathtt{Q}}_{\mT_{m}}+\mS}(\beta_{N})<
\varphi_{\Lambda^{\mathtt{Q}}_{\mW}+\xi}(a)=f_{\mW}(\mathtt{Q},\xi,a)$
from $\Lambda^{\mathtt{R}_{m-1}^{(\vec{\kappa}_{\mathtt{P}}/\vec{\rho})}}_{\mT_{m}}=\Lambda^{\mathtt{Q}}_{\mT_{m}}\leq\Lambda^{\mathtt{Q}}_{\mW}$ and $\mS<\xi$.
MIH with $\mS<\xi$ yields

\beqn\label{eq:main.1_I_1_case1_2}
(\mathcal{H}_{\gamma^{\mathtt{Q}}_{\mT_{m}}},
\Theta_{m-1},
 \mathtt{R}_{m-1}^{(\vec{\kappa}_{\mathtt{P}}/\vec{\rho})}
 )\vdash^{a_{3}}_{\mT_{m}^{\dagger},\vec{\eta}_{\mT_{m}},\mT_{m}^{\dagger},\Lambda^{(3N-2)}_{\mT_{m}}}
\Gamma_{m-1}^{(\vec{\kappa}_{\mathtt{P}}/\vec{\rho})}
\eeqn

Next let us compare $\vec{\kappa}$ with $\vec{\kappa}_{\mathtt{P}}$.
\bclm\label{clm:main.1_I_1_2}
If $(\kappa_{j},f_{j})\neq(\kappa^{\mathtt{P}}_{j},f^{\mathtt{P}}_{j})$, then
$\kappa_{j}\in H_{\kappa^{\mathtt{P}}_{j}}(k_{j}, \mathtt{R}_{m}^{(\vec{\kappa}_{\mathtt{P}}/\vec{\kappa})},\Theta_{m-1})$, where
$k_{j}=h^{b_{m}}(f^{\mathtt{P}}_{j};\varphi_{\mT_{m}^{\dagger}}(\delta))$
for an ordinal $\delta<\Lambda^{(3N-2)}_{\mT_{m}}$.
\eclm
\textbf{Proof} of Claim \ref{clm:main.1_I_1_2}.
If $s(g_{j})\leq b_{m}\leq\mT_{m}^{\dagger}<\mS^{\dagger}$, then 
$(\kappa_{j},f_{j})=(\kappa^{\mathtt{P}}_{j},f^{\mathtt{P}}_{j})=(\rho_{j},g_{j})$.
Assume $b_{m}<s(g_{j})$.

If $b_{m}<s(g_{j})\leq\mS^{\dagger}$, then 
$f_{j}=h^{b_{m}}(f^{\mathtt{P}}_{j};\varphi_{\mT_{m}^{\dagger}}(\alpha_{j-1}))$ and
$\kappa_{j}\in H_{\kappa^{\mathtt{P}}_{j}}(f_{j},\mathtt{R}_{m-1},\Theta_{m-1})$
with $(\kappa^{\mathtt{P}}_{j},f^{\mathtt{P}}_{j})=(\rho_{j},g_{j})$,
where $\Theta_{m-1}(\mathtt{R}_{m}^{(\vec{\kappa}_{\mathtt{P}}/\vec{\kappa})})=\Theta_{m-1}(\mathtt{R}_{m-1})$ and
$\alpha_{j-1}<\Lambda^{(3N-2)}_{\mT_{m}}$.

Let $\mS^{\dagger}<s(g_{j})$. Then
$f_{j}=h^{b_{m}}(g_{j};\varphi_{\mT_{m}^{\dagger}}(\alpha_{j-1}))$
and $\kappa_{j}\in H_{\rho_{j}}(f_{j},\mathtt{R}_{m-1},\Theta_{m-1})$.
On the other hand we have $f^{\mathtt{P}}_{j}=h^{\mS^{\dagger}}(g_{j};\varphi_{\mS^{\dagger}}(\beta_{j-1}))$ and $\kappa^{\mathtt{P}}_{j}\in H_{\rho_{j}}(f^{\mathtt{P}}_{j},\mathtt{R}^{\mathtt{P}}_{m-1},\Theta^{\mathtt{P}}_{m-1})$.
We may assume $\mathtt{p}_{0}(\kappa_{j})\leq\mathtt{p}_{0}(\rho_{j})=\mathtt{p}_{0}(\kappa^{\mathtt{P}}_{j})$ and $\kappa_{j}<\kappa_{j}^{\mathtt{P}}$.
Let $\delta_{1}=\beta_{N}$ and $\delta_{j+1}=\varphi_{\mT_{m}^{\dagger}}(\delta_{j})$.
Then let
$k_{j}=h^{b_{m}}(f^{\mathtt{P}}_{j};\varphi_{\mT_{m}^{\dagger}}(\delta_{j}))$.
We obtain
$\mT_{m}^{\dagger}<\mS^{\dagger}\leq\xi\leq\Lambda_{0}<\Lambda_{\mS}^{\mathtt{P}}\leq\Lambda_{\mathtt{P}}^{(n)}<\Lambda_{\mT_{m}}^{\mathtt{P}}=\Lambda_{\mT_{m}}^{\mathtt{Q}}$
by $\{\mS>\mT_{m}\}\subset\mathrm{dom}(\mathtt{P})$.
$a_{2}=f_{\mS}(\mathtt{P},\xi,a_{0})=
\varphi_{\Lambda_{\mS}^{\mathtt{P}}+\xi}(a_{0})<
\varphi_{\Lambda_{\mT_{m}}^{\mathtt{Q}}+\xi}(a)=f_{\mT_{m}}(\mathtt{Q},\xi,a)$ follows by $a_{0}<a$,
and $\varphi_{\mS^{\dagger}}(\beta_{j-1})<\beta_{N}<\varphi_{\mT_{m}^{\dagger}}(\delta_{j})<f_{\mT_{m}}(\mathtt{Q},\xi,a)\leq\varphi_{\mT_{m}^{\dagger}}(\alpha_{j-1})$.

We obtain $k_{j}\leq f_{j}$.
$\kappa_{j}\in H_{\kappa^{\mathtt{P}}_{j}}(k_{j}, \mathtt{R}_{m}^{(\vec{\kappa}_{\mathtt{P}}/\vec{\kappa})},\Theta_{m-1})$ follows
by $\mathtt{p}_{0}(\kappa_{j})\leq\mathtt{p}_{0}(\kappa^{\mathtt{P}}_{j})$, $\kappa_{j}<\kappa_{j}^{\mathtt{P}}$, $\kappa_{j}\in H_{\rho_{j}}(f_{j},\mathtt{R}_{m-1},\Theta_{m-1})$ and
$\Theta_{m-1}(\mathtt{R}_{m}^{(\vec{\kappa}_{\mathtt{P}}/\vec{\kappa})})=\Theta_{m-1}(\mathtt{R}_{m-1})$.

\hspace*{\fill} $\Box$ of Claim \ref{clm:main.1_I_1_2}.
\\

Let $c_{0}=a_{3}$ and $c_{k+1}=\varphi_{\mS}(c_{k})$.
By Claim \ref{clm:main.1_I_1_2} and (\ref{eq:main.1_I_1_case1_2}), Recapping \ref{lem:recapping} yields
\beqn\label{eq:main.1_I_1_case1_3}
(\mathcal{H}_{\gamma^{\mathtt{Q}}_{\mT_{m}}},
\Theta_{m-1},
  \mathtt{R}_{m-1}^{(\vec{\kappa}/\vec{\rho})}
 )\vdash^{c_{k}}_{\mT_{m}^{\dagger},\vec{\eta}_{\mT_{m}},\mT_{m}^{\dagger},\Lambda^{(3N-2)}_{\mT_{m}}
 }
\Gamma_{m-1}^{(\vec{\kappa}/\vec{\rho})}
\eeqn
where $k=\#\{j\leq N : (\kappa_{j},f_{j})\neq(\kappa^{\mathtt{P}}_{j},f^{\mathtt{P}}_{j})\}$, and
$c_{k}<f_{\mW}(\mathtt{Q},\xi,a)$ by $a_{3}<f_{\mW}(\mathtt{Q},\xi,a)$.
If $m=n$, i.e., if $\mT_{m}=\mW$, then we are done. 

Let $\mW<\mT_{m}$.
By (\ref{eq:main.1_I_1_case1_0}) we have
$\mathtt{R}_{m-1}^{(\vec{\kappa}/\vec{\rho})}\restrict\mT_{m}=\mathtt{R}_{m}$ and
$\Theta_{m-1}\cup(\mathtt{Q}(\mT_{m}))_{0}\subset\Theta_{m}$.
By Lemmas \ref{lem:main.2_L_I_1} and \ref{lem:main.2_S_I_1} with  (\ref{eq:main.1_I_1_case1_3}) 
we obtain
\beqn\label{eq:main.1_I_1_case1_4}
(\mathcal{H}_{\gamma^{\mathtt{Q}}_{\mT_{m}}+\mT_{m}},
\Theta_{m},
  \mathtt{R}_{m}
 )\vdash^{c}_{\mT_{m},\vec{d},\mT_{m},\Lambda^{(6N-3)}_{\mT_{m}}
 }
\Gamma_{m}
\eeqn
for $\vec{d}=(\mT_{m},\ldots,\mT_{m})$ and an ordinal $c<f_{\mW}(\mathtt{Q},\xi,a)$ by $\Lambda_{\mT_{m}}^{(k)}<\Lambda^{\mathtt{Q}}_{\mW}$.

Let
$a_{4}=f_{\mW}(\mathtt{R}_{m},\mT_{m},c)<
f_{\mW}(\mathtt{Q},\xi,a)$ by $\mT_{m}<\xi$ and $c<f_{\mW}(\mathtt{Q},\xi,a)$.
MIH with (\ref{eq:main.1_I_1_case1_4}) yields (\ref{eq:main.1_I_1}).
\\
\textbf{Case 2}. 
Second the case when the last inference is an $(i {\rm -rfl}_{\mS}(g,x,f,\Delta))$
for a successor $i$-stable ordinal $\mS\in \mathrm{dom}(\mathtt{Q})\cap SSt_{i}\subset\xi$:
Then $\bar{\rho}=\bar{\rho}_{\mathtt{Q}}$ and $(\rho,g)\in\mathtt{Q}(\mS)$ 
with $\rho=\rho_{\mathtt{Q}(\mS)}$.
We have an ordinal $a_{0}<a$, and
a finite set $\Delta$ of formulas such that $\mS\leq\mathrm{rk}(\Delta)<\xi$.

Let $\mV$ be a successor $j$-stable ordinal 
such that $j\geq i$, $\mS\leq\mV\leq\mathrm{rk}(\Delta)$ and $\min\{s(g), s(f)\}\geq\mV+1$.
Let $\mathtt{P}=\mathtt{Q}[\mV]$ be a good extension of $\mathtt{Q}$ for $\gamma,\Theta$
such that
$(\bar{\rho_{\mathtt{Q}}}(\mS))_{1}=g=(\bar{\rho_{\mathtt{P}}}(\mV))_{1}$.

Then we have
$(\mathcal{H}_{\gamma},\Theta_{\iota}, \mathtt{P}_{\iota})
\vdash^{a_{0}}_{\xi,\vec{\xi},\xi,\Lambda^{(k)}_{\xi}}
(\Gamma(\mathtt{P})\cup\{\lnot \delta_{\iota}^{(\bar{\rho}_{\mathtt{P}})}\})(\mathtt{P}_{\iota})$ for each $\delta\in\Delta$,
 $\iota\in[\bar{\rho}_{\mathtt{P}}\mathtt{P}]J$ and each good extension 
 $\mathtt{P}_{\iota}=\mathtt{P}[\mathrm{srk}(\delta_{\iota})]$, where
$\mathsf{k}(\delta_{\iota})=\mathsf{k}(\lnot \delta_{\iota}^{(\bar{\rho}_{\mathtt{P}})};\mathtt{P})\subset M_{\mathtt{P}}$,
$\Theta_{\iota}=\Theta\cup\mathsf{k}(\delta_{\iota})$, and
$\delta\simeq\bigvee(\delta_{\iota})_{\iota\in J}$.

On the other hand we have 
$
(\mathcal{H}_{\gamma},\Theta,
\mathtt{P}^{\sigma}
)\vdash^{a_{0}}_{\xi,\vec{\xi},\xi,\Lambda^{(k)}_{\xi}}
\Gamma(\mathtt{P}),\Delta^{(\bar{\rho}_{\sigma})}
$ for every $\sigma\in H_{\rho_{\mV}}(f,\mathtt{P},\Theta)$ with $\rho_{\mV}=\rho_{\mathtt{P}(\mV)}$, where $\bar{\rho}_{\sigma}=\bar{\rho}_{\mathtt{P}}[(\sigma,f)/(\rho_{\mV},g)]$.

Let $\mT=\mathrm{srk}(\mathrm{rk}(\Delta))\geq\mS$ and $\mU=\max\{\mW,\mT\}<\xi$.
Then $\mS\leq\mV\leq\mathrm{rk}(\Delta)$ is equivalent to $\mS\leq\mV\leq\mT$.

We have
$\{\mS,\mathrm{rk}(\Delta)\}\subset\mathcal{H}_{\gamma}[\Theta(\mathtt{Q})]$ by (\ref{eq:controlder_cap22}), and
$\mT\in\mathrm{dom}(\mathtt{Q})$ by (\ref{eq:controlder_cap_cover}).
We obtain
$\mT=\mathrm{srk}(\mathrm{rk}(\Delta))\in\mathcal{H}_{\gamma}[\Theta(\mathtt{Q})]$ by Proposition \ref{prp:pd_closed_cover_I_1}.\ref{prp:pd_closed_cover_I_1.1}, and
$\mT^{\dagger}\in\mathcal{H}_{\gamma}[\Theta(\mathtt{Q})]$ follows.
\\
\textbf{Case 2.1}. $\mT\leq\mW$: Then $\mU=\mW\geq\mS$ and 
$\mathrm{srk}(\mathrm{rk}(\Gamma\cup\{\lnot\delta_{\iota}\}))=
\mathrm{srk}(\mathrm{rk}(\Gamma\cup\Delta))=\mW$.
Let $\mS\leq\mV\leq\mathrm{rk}((\Delta^{(\bar{\rho}_{\sigma})})^{[\mW^{\dagger}*]})\leq\mathrm{rk}(\Delta)\leq\mW$ and consider good extensions $\mathtt{P}=\mathtt{Q}[\mV]$.

For $\delta\in\Delta$ we have 
$\mathrm{rk}(\delta^{[\mW^{\dagger}*]})\leq\mathrm{rk}(\delta)<\mT^{\dagger}\leq\mW^{\dagger}$.
As in \textbf{Case 1}, SIH yields
$(\mathcal{H}_{\gamma^{\mathtt{Q}}_{\mW}},\Theta_{\mW^{\dagger}*,\iota}, (\mathtt{P}_{\iota})_{\mW}
)
\vdash^{f_{\mW}(\mathtt{Q},\xi,a_{0})}_{\mW^{\dagger},\vec{\eta}_{\mW}, \mW^{\dagger}, \Lambda^{(3N-2)}_{\mW}}
((\Gamma(\mathtt{P})\cup\{\lnot (\delta_{\iota}^{(\bar{\rho}_{\mathtt{P}})})\})(\mathtt{P}_{\iota}))^{[\mW^{\dagger}*]}$, where
$\Theta_{\mW^{\dagger}*,\iota}$ is defined from $\Theta_{\iota}$,
$\gamma^{\mathtt{Q}}_{\mW}=\gamma^{\mathtt{P}}_{\mW}$ and
$f_{\mW}(\mathtt{Q},\xi,a_{0})=f_{\mW}(\mathtt{P},\xi,a_{0})$ by
$\delta^{\mathtt{Q}}_{\mW}=\delta^{\mathtt{P}}_{\mW}$.

On the other hand we have
$
(\mathcal{H}_{\gamma},\Theta,
\mathtt{P}^{\sigma}
)\vdash^{a_{0}}_{\xi,\vec{\xi},\xi,\Lambda^{(k)}_{\xi}}
\Gamma(\mathtt{P}),\Delta^{(\bar{\rho}_{\sigma})}
$. 
For each $\sigma\in H_{\rho_{\mV}}(f,\mathtt{P}_{\mW},\Theta_{\mW^{\dagger}*})\subset H_{\rho_{\mV}}(f,\mathtt{P},\Theta)$,
$
(\mathcal{H}_{\gamma^{\mathtt{Q}}_{\mW}},\Theta_{\mW^{\dagger}*},
(\mathtt{P}^{\sigma})_{\mW}
)\vdash^{f_{\mW}(\mathtt{Q},\xi,a_{0})}_{\mW^{\dagger},\vec{\eta}_{\mW},\mW^{\dagger},\Lambda^{(3N-2)}_{\mW}}
\Gamma(\mathtt{P})^{[\mW^{\dagger}*]}, (\Delta^{(\bar{\rho}_{\sigma})})^{[\mW^{\dagger}*]}
$ follows by SIH.
$
(\mathcal{H}_{\gamma^{\mathtt{Q}}_{\mW}},\Theta_{\mW^{\dagger}*},
\mathtt{Q}_{\mW}
)\vdash^{f_{\mW}(\mathtt{Q},\xi,a_{0})+1}_{\mW^{\dagger},\vec{\eta}_{\mW},\mW^{\dagger},\Lambda^{(3N-2)}_{\mW}}
\Gamma^{[\mW^{\dagger}*]}
$ follows by an
$(i {\rm -rfl}_{\mS}(g,x,f,(\Delta^{(\bar{\rho}_{\sigma})})^{[\mW^{\dagger}*]}))$ with 
$\mS\leq\mathrm{rk}((\Delta^{(\bar{\rho}_{\sigma})})^{[\mW^{\dagger}*]})<\mW^{\dagger}$.
We obtain (\ref{eq:main.1_I_1}).
\\
\textbf{Case 2.2}. $\mW<\mT$: Then $\mU=\mT\geq\mS$ and
$\mathrm{srk}(\mathrm{rk}(\Gamma\cup\{\lnot\delta_{\iota}\}\}))\leq
\mathrm{srk}(\mathrm{rk}(\Gamma\cup\Delta))=\mT$.
Let $\mT=\mT_{m}$ for an $m<n$. Then 
$\mathtt{Q}_{\mT}=\mathtt{R}_{m-1}^{(\vec{\kappa}/\vec{\rho})}$ for some $\vec{\kappa}$ and $\vec{\rho}$ of ordinals, cf.\,(\ref{eq:main.1_I_1_case1_0}).
Let $\mS\leq\mV\leq\mathrm{rk}((\Delta^{(\bar{\rho}_{\sigma})})^{[\mT^{\dagger}*]})\leq\mathrm{rk}(\Delta)$ and consider good extensions $\mathtt{P}=\mathtt{Q}[\mV]$.

Let $a_{1}=f_{\mT}(\mathtt{Q},\xi,a_{0})=
\varphi_{\Lambda^{\mathtt{Q}}_{\mT}+\xi}(a_{0})<f_{\mT}(\mathtt{Q},\xi,a)$, where
$a_{0},\xi<\Lambda^{(k)}_{\xi}\leq\Lambda^{(1)}_{\mT}$ and
$\Lambda^{\mathtt{Q}}_{\mT}<\Lambda^{(1)}_{\mT}$.
We obtain $a_{1}<\Lambda^{(1)}_{\mT}$, and
$\Lambda^{\mathtt{Q}}_{\mT}<\Lambda^{\mathtt{Q}}_{\mW}$ by $\mT\in\mathrm{dom}(\mathtt{Q})$.
Hence
$a_{1}<\varphi_{\Lambda^{\mathtt{Q}}_{\mW}+\xi}(a)=f_{\mW}(\mathtt{Q},\xi,a)$.
On the other hand we have $\gamma^{\mathtt{Q}}_{\mT}=\gamma^{\mathtt{P}}_{\mT}$
and $\delta^{\mathtt{Q}}_{\mT}=\delta^{\mathtt{P}}_{\mT}$ by $\mT\in\mathrm{dom}(\mathtt{Q})$.
We obtain $a_{1}=f_{\mT}(\mathtt{Q},\xi,a_{0})=f_{\mT}(\mathtt{P},\xi,a_{0})$.
Let $\vec{\eta}_{\mT}=(\mT^{\dagger},\ldots,\mT^{\dagger})$.

We obtain
$(\mathcal{H}_{\gamma^{\mathtt{Q}}_{\mT}},\Theta_{\mT^{\dagger}*,\nu}, (\mathtt{P}_{\iota})_{\mT}
)\vdash^{a_{1}}_{\mT^{\dagger},\vec{\eta}_{\mT}, \mT^{\dagger},\Lambda^{(3N-2)}_{\mT}}
((\Gamma(\mathtt{P})\cup\{ \lnot (\delta_{\iota}^{(\bar{\rho}_{\mathtt{P}})}\})(\mathtt{P}_{\iota}))^{[\mT^{\dagger}*]}$ by SIH.
On the other,
$
(\mathcal{H}_{\gamma^{\mathtt{Q}}_{\mT}},\Theta_{\mT^{\dagger}*},
(\mathtt{P}^{\sigma})_{\mT}
)\vdash^{a_{1}}_{\mT^{\dagger},\vec{\eta}_{\mT},\mT^{\dagger},\Lambda^{(3N-2)}_{\mT}}
\Gamma(\mathtt{P})^{[\mT^{\dagger}*]}, (\Delta^{(\bar{\rho}_{\sigma})})^{[\mT^{\dagger}*]}
$ follows by SIH for every $\sigma\in H_{\rho_{\mV}}(f,\mathtt{P}_{\mT},\Theta_{\mT^{\dagger}*})\subset H_{\rho_{\mV}}(f,\mathtt{P},\Theta)$.

An $(i {\rm -rfl}_{\mS}(g,x,f, (\Delta^{(\bar{\rho}_{\sigma})})^{[\mT^{\dagger}*]}))$ with 
$\mS\leq\mathrm{rk}((\Delta^{(\bar{\rho}_{\sigma})})^{[\mT^{\dagger}*]})\leq \mathrm{rk}(\Delta)<\mT^{\dagger}$
yields 
$
(\mathcal{H}_{\gamma^{\mathtt{Q}}_{\mT}},\Theta_{\mT^{\dagger}*},
\mathtt{Q}_{\mT}
)\vdash^{a_{1}+1}_{\mT^{\dagger},\vec{\eta}_{\mT},\mT^{\dagger},\Lambda^{(3N-2)}_{\mT}}
\Gamma^{[\mT^{\dagger}*]}
$.

We have $\mathrm{rk}(\Gamma)<\mW^{\dagger}\leq\mT$, and $\mT\geq\mS\in\mathrm{dom}(\mathtt{Q})\cap SSt_{i}$.
We claim that $\mS=\max(\mathrm{dom}(\mathtt{Q}\restrict\mT^{\dagger})\cap \bigcup_{j\geq i}SSt_{j})$.
We may assume that there exists a $\mU$ such that $\mS<\mU\in\mathrm{dom}(\mathtt{Q})\cap SSt_{j}$ with $j\geq i$.
Then
$\mT\leq\mathrm{rk}(\Delta)<\mU$ by (r1), and the claim follows.
On the other hand we have 
$a_{1}
<\Lambda^{(1)}_{\mT}$.
\\
\textbf{Case 2.2.1}. $\mT\in SSt$: If $\mT\in SSt_{j}$ with $j\geq i$, then
 $\mS=\mT$ and $(\rho,g)=(\rho_{\mathtt{Q}(\mT)},g)\in\mathtt{Q}(\mT)$.
 If $\mT\in SSt_{j}$ with $j< i$, then 
 $\mS<\mT$ and $(\rho_{\mathtt{Q}(\mT)},g_{\mT})\in\mathtt{Q}(\mT)$.

We have $(\mathcal{H}_{\gamma^{\mathtt{Q}}_{\mT}},\Theta_{\mT^{\dagger}*},
\mathtt{Q}_{\mT}
)\vdash^{a_{1}+1}_{\mT^{\dagger},\vec{\eta}_{\mT},\mT^{\dagger},\Lambda^{(3N-2)}_{\mT}}
\Gamma^{[\mT^{\dagger}*]}
$.
On the other hand we have
$a_{1}<\Lambda^{(1)}_{\mT}$ and 
$\mathrm{rk}(\Gamma^{[\mT^{\dagger}*]})\leq\mathrm{rk}(\Gamma)<\mT$.

For $1\leq j\leq N$, let $\mS_{j}\in\{0, \max(\mathrm{dom}(\mathtt{Q}\restrict\mT^{\dagger})\cap SSt_{j})\}$,
$\kappa_{j},\rho_{j}$ be ordinals, $\mathtt{Q}_{j}$ finite families,
$\Phi_{j}$ sequents and ordinals $\alpha_{j}(\mathtt{Q}_{\mT})=\alpha(\mathtt{Q}_{\mT};\mS_{j})$ 
defined from the sequent $\Phi_{0}=\Gamma^{[\mT^{\dagger}*]}$, 
ordinals $\mT$, $\alpha_{0}=\varphi_{\mT^{\dagger}}(a_{1}+1)$,
 the ordinals $\vec{\eta}_{\mT}$, and strongly critical numbers $\Lambda_{\mT}^{(n)}\, (n>0)$
in Definition \ref{df:main_lemmas}.\ref{df:main_lemmas.2}.

We have
$\varphi_{\mT^{\dagger}}(a_{1}+1)<\varphi_{\mT^{\dagger}}(f_{\mT}(\mathtt{Q},\xi,a))$, and
we may assume that each $\kappa_{j}$ is defined from the ordinal
$\varphi_{\mT^{\dagger}}(f_{\mT}(\mathtt{Q},\xi,a))$.
Hence 
$\mathtt{Q}_{N}\restrict\mT=\mathtt{R}_{m}$ and 
$\Phi_{N}=\{A_{m}^{(\vec{\tau}_{m})}: A^{(\bar{\rho})}\in\Gamma\}=\Gamma_{m}$
and
$\Gamma^{[\mT^{\dagger}*]}=\{A_{m-1}^{(\vec{\tau}_{m-1})}: A^{(\bar{\rho})}\in\Gamma\}=\Gamma_{m-1}$, where
$A_{m}^{(\vec{\tau}_{m})}\equiv (A_{m-1}^{(\vec{\tau}_{m-1})})^{[\mT]}(\mathtt{Q}_{\mT})$ and
$\mathtt{Q}_{\mT}=\mathtt{R}_{m-1}^{(\vec{\kappa}/\vec{\rho})}$ with $\mT=\mT_{m}$.
Then 
\[
(\mathcal{H}_{\gamma_{1}},\Theta_{\mT*}, \mathtt{R}_{m}
)
\vdash^{\alpha_{N}}_{\mT,\vec{\mT},\mT,\Lambda_{\mT}^{(6N-3)}}
\Gamma_{m}
\]
holds by Lemma \ref{lem:main.2_S_I_1},
where 
$\gamma_{1}=\gamma^{\mathtt{Q}}_{\mT}+\mT$,
$\Lambda^{\mathtt{Q}_{N}}_{\mT}=\Lambda^{\mathtt{Q}}_{\mT}$ by 
$\delta^{\mathtt{Q}_{N}}_{\cdot}=\delta^{\mathtt{Q}}_{\cdot}$, and
$\Theta_{\mT^{\dagger}*}\cup(\mathtt{Q}_{N}(\mT))_{0}=\Theta_{\mT *}$.
We have
$\mathrm{rk}(\Gamma_{m})\leq\mathrm{rk}(\Gamma_{m-1})\leq\mathrm{rk}(\Gamma)$.

Let $\mT_{n}=\mW$ and $m<k\leq n$.
We show 
$\alpha_{j}<f_{\mT_{k}}(\mathtt{Q},\xi,a)=\vphi_{\Lambda^{\mathtt{Q}}_{\mT_{k}}+\xi}(a)
\leq f_{\mW}(\mathtt{Q},\xi,a)$ by induction on $j$.
We have
$a_{1}=f_{\mT_{m}}(\mathtt{Q},\xi,a_{0})<f_{\mT_{k}}(\mathtt{Q},\xi,a)$ by
$\{\mT_{m}>\mT_{k}\}\subset\mathrm{dom}(\mathtt{Q})$ and $a_{0}<a$.
$\alpha_{0}=\varphi_{\mT^{\dagger}}(a_{1}+1)<f_{\mT_{k}}(\mathtt{Q},\xi,a)$ follows from
$\mT^{\dagger}\leq\xi\leq\Lambda_{0}<\Lambda^{\mathtt{Q}}_{\mT_{k}}$.

Assuming $\alpha_{j-1}<f_{\mT_{k}}(\mathtt{Q},\xi,a)$, we obtain
$\alpha_{j}=\varphi_{\alpha}(\alpha_{j-1})<f_{\mT_{k}}(\mathtt{Q},\xi,a)$ for
$\alpha=\alpha_{j}(\mathtt{Q}_{\mT})=\alpha(\mathtt{Q}_{\mT};\mS_{j})=
\Lambda_{\mT}^{(3j-2)}\cdot\delta^{\mathtt{Q}}_{\mS_{j}}+on_{\Lambda_{\mT}^{(3j-3)}}(g_{j})$,
where $\Lambda_{\mT}^{(p)}=\psi_{\mI_{N}}(b_{0}+\delta^{\mathtt{Q}}_{\mT}+p)<\psi_{\mI_{N}}(b_{0}+\delta^{\mathtt{Q}}_{\mT_{k}})=\Lambda^{\mathtt{Q}}_{\mT_{k}}$ for $p<\omega$ by
$\delta^{\mathtt{Q}}_{\mT}+p<\delta^{\mathtt{Q}}_{\mT_{k}}$ with 
$\{\mT=\mT_{m}>\mT_{k}\}\subset\mathrm{dom}(\mathtt{Q})$.

In particular we obtain
$\varphi_{\mT_{k}^{\dagger}}(f_{\mT_{k}}(\mathtt{Q}_{m},\mT,\alpha_{N}))<
\varphi_{\mT_{k}^{\dagger}}(f_{\mT_{k}}(\mathtt{Q},\xi,a))$ by $\mT<\xi$ and 
$\Lambda^{\mathtt{Q}_{m}}_{\mT_{k}}=\Lambda^{\mathtt{Q}}_{\mT_{k}}$.
We may assume that
each $\kappa_{i}\in H_{\rho_{i}}(f_{i},\mathtt{Q}_{k-1},\Theta_{k-1})$ is in the resolvent class defined from $\varphi_{\mT_{k}^{\dagger}}(f_{\mT_{k}}(\mathtt{Q}_{m},\mT,\alpha_{N}))$.

We have $a_{2}=f_{\mW}(\mathtt{Q},\mT,\alpha_{N})=\vphi_{\Lambda^{\mathtt{Q}}_{\mW}+\mT}(\alpha_{N})<
\vphi_{\Lambda^{\mathtt{Q}}_{\mW}+\xi}(a)=f_{\mW}(\mathtt{Q},\xi,a)$ by $\mT<\xi$ and $\alpha_{N}<f_{\mW}(\mathtt{Q},\xi,a)$.
By MIH with $\mT<\xi$ we obtain (\ref{eq:main.1_I_1}).
\\
\textbf{Case 2.2.2}. $\mT\not\in SSt$:
We have $(\mathcal{H}_{\gamma^{\mathtt{Q}}_{\mT}},\Theta_{\mT^{\dagger}*},
\mathtt{Q}_{\mT}
)\vdash^{a_{1}+1}_{\mT^{\dagger},\vec{\eta}_{\mT},\mT^{\dagger},\Lambda^{(3N-2)}_{\mT}}
\Gamma^{[\mT^{\dagger}*]}
$.
On the other hand we have
$a_{1}<\Lambda^{(1)}_{\mT}$ and $\mathrm{rk}(\Gamma^{[\mT^{\dagger}*]})\leq\mathrm{rk}(\Gamma)<\mT$.

As in \textbf{Case 2.2.1} we see from Lemma \ref{lem:main.2_L_I_1}
\[
(\mathcal{H}_{\gamma^{\mathtt{Q}}_{\mT}},\Theta_{\mT*}, \mathtt{R}_{m}
)
\vdash^{\alpha_{N}}_{\mT,\vec{\mT},\mT,\Lambda_{\mT}^{(6N-3)}}
\Gamma_{m}
\]
where 
$\Theta_{\mT^{\dagger}*}=\Theta_{\mT *}$ by
$\mathtt{Q}(\mT)=\emptyset$.
Let $\mT_{n}=\mW$ and $m<k\leq n$.
We see 
$\alpha_{N}<f_{\mT_{k}}(\mathtt{Q},\xi,a)=\vphi_{\Lambda^{\mathtt{Q}}_{\mT_{k}}+\xi}(a)
\leq f_{\mW}(\mathtt{Q},\xi,a)$, and
$\varphi_{\mT_{k}^{\dagger}}(f_{\mT_{k}}(\mathtt{Q}_{m},\mT,\alpha_{N}))<
\varphi_{\mT_{k}^{\dagger}}(f_{\mT_{k}}(\mathtt{Q},\xi,a))$.
We may assume that
each $\kappa_{i}\in H_{\rho_{i}}(f_{i},\mathtt{Q}_{k-1},\Theta_{k-1})$ is in the resolvent class defined from $\varphi_{\mT_{k}^{\dagger}}(f_{\mT_{k}}(\mathtt{Q}_{m},\mT,\alpha_{N}))$.
We obtain $a_{2}=f_{\mW}(\mathtt{Q},\mT,\alpha_{N})<f_{\mW}(\mathtt{Q},\xi,a)$.
By MIH with $\mT<\xi$ we obtain (\ref{eq:main.1_I_1}).
\\
\textbf{Case 3}. Third the last inference is a $(cut)$:
Let $\bar{\rho}=\bar{\rho}_{\mathtt{Q}}$.
We have an $a_{0}<a$ and a $C\simeq\bigvee(C_{\iota})_{\iota\in J}$ such that $\mathrm{rk}(C)<\xi$ for which
$(\mathcal{H}_{\gamma},\Theta_{\iota}, \mathtt{Q}_{\iota})
\vdash^{a_{0}}_{\xi,\vec{\xi},\xi,\Lambda^{(k)}_{\xi}}
(\Gamma\cup\{\lnot C_{\iota}^{(\bar{\rho})}\})(\mathtt{Q}_{\iota})$ holds for each $\iota\in[\bar{\rho}\mathtt{Q}]J$ and
each good extension $\mathtt{Q}_{\iota}=\mathtt{Q}[\mathrm{srk}(C_{\iota})]$, where
$\Theta_{\iota}=\Theta\cup\mathsf{k}(C_{\iota})$ with $\mathsf{k}(C_{\iota})=\mathsf{k}(C_{\iota}^{(\bar{\rho})};\mathtt{Q})$.
On the other hand we have $(\mathcal{H}_{\gamma},\Theta, \mathtt{Q})
\vdash^{a_{0}}_{\xi,\vec{\xi},\xi,\Lambda^{(k)}_{\xi}}
\Gamma, C^{(\bar{\rho})}$.
Let $d=\mathrm{rk}(C)$ and $\mT=\mathrm{srk}(d)$.
\\
\textbf{Case 3.1}. $\mT\leq\mW$:
We obtain by SIH,
$(\mathcal{H}_{\gamma^{\mathtt{Q}}_{\mW}},\Theta_{\mW^{\dagger}*,\iota}, 
(\mathtt{Q}_{\iota})_{\mW})
\vdash^{a_{1}}_{\mW^{\dagger},\eta_{\mW},\mW^{\dagger},\Lambda^{(3N-2)}_{\mW}}
((\Gamma\cup\{\lnot C_{\iota}^{(\bar{\rho})}\})(\mathtt{Q}_{\iota}))^{[\mW^{\dagger}*]}$, where
$a_{1}=f_{\mW}(\mathtt{Q},\xi,a_{0})$.
On the other hand we have
$(\mathcal{H}_{\gamma^{\mathtt{Q}}_{\mW}},\Theta_{\mW^{\dagger}*}, \mathtt{Q}_{\mW})
\vdash^{a_{1}}_{\mW^{\dagger},\eta_{\mW},\mW^{\dagger},\Lambda^{(3N-2)}_{\mW}}
\Gamma^{[\mW^{\dagger}*]}, (C^{(\bar{\rho})})^{[\mW^{\dagger}*]}$.

We have $\mathrm{rk}((C^{(\bar{\rho})})^{[\mW^{\dagger}*]})\leq \mathrm{rk}(C^{(\bar{\rho})})=d<\mW^{\dagger}$.
A $(cut)$ yields (\ref{eq:main.1_I_1}).
\\
\textbf{Case 3.2}. $\mW<\mT$:
Let 
$a_{1}=f_{\mT}(\mathtt{Q},\xi,a_{0})=\varphi_{\Lambda^{\mathtt{Q}}_{\mT}+\xi}(a_{0})<\Lambda^{(1)}_{\mT}$.
We obtain $\mW<\mT\leq d<\mT^{\dagger}$, and $\{\mT,\mT^{\dagger}\}\subset\mathcal{H}_{\gamma}[\Theta(\mathtt{Q})]$
by $d\in\mathcal{H}_{\gamma}[\Theta(\mathtt{Q})]$ and Proposition \ref{prp:pd_closed_cover_I_1}.\ref{prp:pd_closed_cover_I_1.1}.
SIH yields
$(\mathcal{H}_{\gamma^{\mathtt{Q}}_{\mT}},\Theta_{\mT^{\dagger}*,\iota}, (\mathtt{Q}_{\iota})_{\mT})
\vdash^{a_{1}}_{\mT^{\dagger},\eta_{\mT},\mT^{\dagger},\Lambda^{(3N-2)}_{\mT}}
((\Gamma\cup\{\lnot C_{\iota}^{(\bar{\rho})}\})(\mathtt{Q}_{\iota}))^{[\mT^{\dagger}*]}$, and
$(\mathcal{H}_{\gamma^{\mathtt{Q}}_{\mT}},\Theta_{\mT^{\dagger}*}, \mathtt{Q}_{\mT})
\vdash^{a_{1}}_{\mT^{\dagger},\eta_{\mT},\mT^{\dagger},\Lambda^{(3N-2)}_{\mT}}
\Gamma^{[\mT^{\dagger}*]}, (C^{(\bar{\rho})})^{[\mT^{\dagger}*]}$.
We have
$\mathrm{rk}((C^{(\bar{\rho})})^{[\mT^{\dagger}*]})\leq \mathrm{rk}(C^{(\bar{\rho})})=d<\mT^{\dagger}$.
A $(cut)$ yields 
\[
(\mathcal{H}_{\gamma^{\mathtt{Q}}_{\mT}},\Theta_{\mT^{\dagger}*}, \mathtt{Q}_{\mT})
\vdash^{a_{1}+1}_{\mT^{\dagger},\eta_{\mT},\mT^{\dagger},\Lambda^{(3N-2)}_{\mT}}
\Gamma^{[\mT^{\dagger}*]}
\]
where $a_{1}<\Lambda^{(1)}_{\mT}$ and $\mT=\mathrm{srk}(C)\in\mathrm{dom}(\mathtt{Q})$ by (\ref{eq:controlder_cap_cover}).

As in \textbf{Case 2.2},
we obtain (\ref{eq:main.1_I_1}) by Lemmas \ref{lem:main.2_S_I_1} and \ref{lem:main.2_L_I_1}
followed by MIH with $\mT<\xi$.
\\
\textbf{Case 4}. The last inference is a $(\bigvee)$ introducing an $A^{(\bar{\rho})}\in\Gamma$:
The case is seen from SIH as in \textbf{Case 1} of the proof of Lemma \ref{mlem:singlemainl_S_I_1}.
Note that
$\varphi_{\mT_{m}^{\dagger}}(f_{\mT_{m}}(\mathtt{Q},\xi,a_{0}))<\varphi_{\mT_{m}^{\dagger}}(f_{\mT_{m}}(\mathtt{Q},\xi,a))$, and each $\kappa_{i}\in H_{\rho_{i}}(f_{i},\mathtt{Q}_{m-1},\Theta_{m-1})$ is in the resolvent class defined from $\varphi_{\mT_{m}^{\dagger}}(f_{\mT_{m}}(\mathtt{Q},\xi,a_{0}))$.

The case when the last inference is a 
$(\Sigma(\Omega){\rm -rfl})$
is seen from SIH.
\eprf

\subsection{Uncapped calculus}\label{subsec:uncapped}

Let $\mS_{0}=0^{\dagger}$ be the least stable ordinal.
Since there is no stable ordinal below $\mS_{0}$, $\mathrm{dom}(\mathtt{Q})\subset\mS_{0}$ means 
$\mathrm{dom}(\mathtt{Q})\subset\{0\}$.
Let $\mathtt{Q}_{0}$ denote the family such that $\mathrm{dom}(\mathtt{Q}_{0})=\{0\}$.
Then $\mathtt{Q}_{0}(\mS)=\emptyset$ for any (successor) stable ordinal $\mS$.
In a witnessed derivation of the fact $(\mathcal{H}_{\gamma},\Theta,\mathtt{Q}_{0})\vdash^{a}_{\mS_{0},\vec{\eta}_{0},\mS_{0},\Lambda}\Gamma$ with $\vec{\eta}_{0}=(\mS_{0},\ldots,\mS_{0})$,
there occurs no
inference $(i {\rm -rfl}_{\mS}(g,x,f,\Delta))$ 
nor $(\mathrm{dom})$,
cf.\,Definition \ref{df:controldercap}.
Hence each formula $A^{(\bar{\rho})}$ occurring in it is uncapped with the empty list $\bar{\rho}=\emptyset$.
Let us write $A^{(\mathtt{u})}$ for $A^{(\emptyset)}$ in this subsection.

We see from Definition \ref{df:Qmin}.\ref{df:Qmin.1} that $\mathrm{inv}(A^{(\mathtt{u})}; \mathtt{Q}_{0})\equiv A^{(\mathtt{u})}$,
$\mathsf{k}(A^{(\mathtt{u})};\mathtt{Q}_{0})=\mathsf{k}(A)$ and
$\iota\in[\mathtt{u}\mathtt{Q}_{0}]J$ iff $\mathsf{k}(A_{\iota})\subset M_{\mathtt{Q}_{0}}$ for
$A\simeq\bigvee(A_{\iota})_{\iota\in J}$.
By (\ref{eq:controlder_cap2})
we have $
\mathsf{k}(\Gamma)
 \subset
\mathcal{H}_{\gamma}[
\Theta
]$, where $\Theta=\Theta\cap M_{\mathtt{Q}_{0}}=\Theta(\mathtt{Q}_{0})$.

In this subsection we write
$(\mathcal{H}_{\gamma},\Theta)\vdash^{a}_{c}\Gamma^{(\mathtt{u})}$
for $(\mathcal{H}_{\gamma},\Theta,\mathtt{Q}_{0})\vdash^{a}_{c,\vec{\eta}_{0},\mS_{0},\Lambda}\Gamma^{(\mathtt{u})}$, where
$\{a,c\}
\subset
\mathcal{H}_{\gamma}[\Theta]
$
by (\ref{eq:controlder_cap22}), and
$\gamma^{\mathtt{Q}_{0}}_{0}\in\mathcal{H}_{\gamma^{\mathtt{Q}_{0}}_{0}}[\Theta]$
by (\ref{eq:controlder_cap0}).

\blem\label{lem:CEu}{\rm (Cut-elimination)}

\benu
\item\label{lem:CEu1}
Assume $(\mathcal{H}_{\gamma},\Theta)\vdash^{a}_{\Omega+1+b}\Gamma^{(\mathtt{u})}$.
Then
$(\mathcal{H}_{\gamma},\Theta)\vdash^{\varphi_{b}(a)}_{\Omega+1}\Gamma^{(\mathtt{u})}$.

\item\label{lem:CEu2}
Assume $(\mathcal{H}_{\gamma},\Theta)\vdash^{a}_{b}\Gamma^{(\mathtt{u})}$ with $b<\Omega$.
Then
$(\mathcal{H}_{\gamma},\Theta) \vdash^{\varphi_{b}(a)}_{0}\Gamma^{(\mathtt{u})}$
holds.

\eenu

\elem
\bprf
Each is seen by main induction on $b$ with subsidiary induction on $a$.
\eprf

\blem\label{lem:collapseOmega}{\rm (Collapsing)}
Let
$\Gamma\subset\Sigma(\mathcal{L}_{0}:\Omega)$ be a set of formulas.
Suppose
$\Theta\subset
C_{\gamma+1}(\psi_{\Omega}(\gamma+1))$
 and
$
(\mathcal{H}_{\gamma},\Theta
)\vdash^{a}_{\Omega+1}\Gamma
$ with $\{a,\gamma\}\subset\mathcal{H}_{\gamma}[\Theta(\mathtt{Q})]$.
Let
$\beta=\psi_{\Omega}(\hat{a})$ with $\hat{a}=\gamma+\omega^{a}$.
Then 
$(\mathcal{H}_{\hat{a}},\Theta)
\vdash^{\beta}_{\beta}
\Gamma^{(\beta,\Omega)}$ holds.
\elem
\bprf
By induction on $a$ as in \cite{Buchholz}.
\eprf
\\

Let us prove Theorem \ref{thm:2}.
Let 
$S_{\mI_{N}}\vdash\theta^{\mathsf{L}_{\Omega}}$ for a $\Sigma$-sentence $\theta$.
By Embedding \ref{th:embedregthm} pick an $m>0$ so that 
 $(\mathcal{H}_{0},\emptyset; \emptyset)
 \vdash_{\mI_{N}+1+m}^{* \mI_{N}\cdot 2+m}
\theta^{L_{\Omega}}$.
Cut-elimination \ref{lem:predcereg*} yields 
$(\mathcal{H}_{0},\emptyset;\emptyset)
\vdash_{\mI_{N}+1}^{* a}\theta^{\mathsf{L}_{\Omega}}$
for $a=\omega_{m}(\mI_{N}\cdot 2+m)<\omega_{m+1}(\mI_{N}+1)$.
Then Collapsing \ref{lem:Kcollpase.1} yields
$(\mathcal{H}_{\hat{a}},\emptyset;\emptyset)
\vdash_{\beta}^{* \beta}\theta^{\mathsf{L}_{\Omega}}$
for $\beta=\psi_{\mI_{N}}(\hat{a})\in LSt_{N}$ with 
$b_{0}:=\hat{a}=\ome^{a}=\omega_{m+1}(\mI_{N}\cdot 2+m)>\bet$.

Now let 
$\Lambda_{0}:=\beta=\psi_{\mI_{N}}(b_{0})<\Lambda_{1}=\psi_{\mI_{N}}(b_{0}+1)<\mathbb{I}_{N}$ be the next strongly critical number, and
$\gamma=b_{0}+\omega^{\Lambda_{0}+1}$.
Capping \ref{lem:capping} then yields
$(\mathcal{H}_{\gamma},\emptyset,\mathtt{Q}_{0})
\vdash_{\Lambda_{0},\vec{d},\Lambda_{0},\Lambda_{1}}^{\Lambda_{0}\cdot 2} (\theta^{\mathsf{L}_{\Omega}})^{(\mathtt{u})}$
where  $c=\xi=\Lambda_{0}=2\Lambda_{0}$, $\vec{d}=(\Lambda_{0},\ldots,\Lambda_{0})$,
$\mathtt{u}$ denotes the empty cap,
$\mathtt{Q}_{0}$ is a finite family such that
$\mathrm{dom}(\mathtt{Q}_{0})=\{0\}$ with thresholds 
$\gamma_{0}=\gamma^{\mathtt{Q}_{0}}_{0}=\gamma+\omega^{\Lambda_{0}\cdot 2}=b_{0}+\omega^{\Lambda_{0}\cdot 2}\geq\gamma+\omega$ and
$\delta^{\mathtt{Q}_{0}}_{0}=\omega^{\Lambda_{0}+1}$, and
$\Lambda_{\Lambda_{0}}^{(k)}=\Lambda_{1}$.

In what follows each finite function is an $f:\mI_{N}\to\Gamma(\mI_{N})$
with the base $\mI_{N}$ of the $\tilde{\theta}$-function.
Let
$\mS_{0}=0^{\dagger}$ be the least stable ordinal.

By Lemma \ref{lem:main.1_I_1} we obtain
$(\mathcal{H}_{\gamma_{0}},\emptyset,\mathtt{Q}_{0})
\vdash^{\alpha}_{0^{\dagger},\vec{\eta}_{0},0^{\dagger},\Lambda}
(\theta^{\mathsf{L}_{\Omega}})^{(\mathtt{u})}$, where
\\
$\alpha=f_{0}(\mathtt{Q}_{0},\Lambda_{0},\Lambda_{0}\cdot 2)=
\varphi_{\Lambda^{\mathtt{Q}_{0}}_{0}+\Lambda_{0}}(\Lambda_{0}\cdot 2)$ with
$\Lambda^{\mathtt{Q}_{0}}_{0}=\psi_{\mI_{N}}(b_{0}+\delta^{\mathtt{Q}_{0}}_{0})$,
$\Lambda=\psi_{\mI_{N}}(b_{0}+\delta^{\mathtt{Q}_{0}}_{0}+3N-2)$ and
$\gamma_{0}=\gamma^{\mathtt{Q}_{0}}_{0}\in\mathcal{H}_{\gamma_{0}}[\emptyset]$
by (\ref{eq:controlder_cap0}).
We obtain
$(\mathcal{H}_{\gamma_{0}},\emptyset)
\vdash_{\mS_{0}}^{\alp}
(\theta^{\mathsf{L}_{\Omega}})^{(\mathtt{u})}$.

Cut-elimination \ref{lem:CEu}.\ref{lem:CEu1} yields
$(\mathcal{H}_{\gamma_{0}},\emptyset)
\vdash_{\Omega+1}^{\alpha}
(\theta^{\mathsf{L}_{\Omega}})^{(\mathtt{u})}$ for $\alpha=\varphi_{\mS_{0}}(\alp)$.

$(\mathcal{H}_{\gamma_{1}},\emptyset)
\vdash_{\delta}^{\del}(\theta^{\mathsf{L}_{\del}})^{(\mathtt{u})}$ follows from 
Collapsing \ref{lem:collapseOmega}
for
$\gamma_{1}=\gamma_{0}+\alpha\in\mathcal{H}_{\gamma}[\emptyset]$ and
$\del=\psi_{\Ome}(\gamma_{0}+\alpha)$ with the epsilon
number $\alpha$.
We obtain
$(\mathcal{H}_{\gamma_{1}},\emptyset)
\vdash_{0}^{\vphi_{\del}(\del)}(\theta^{\mathsf{L}_{\del}})^{(\mathtt{u})}$
by Cut-elimination \ref{lem:CEu}.\ref{lem:CEu2}.

We see that $\theta^{L_{\del}}$ is true by induction up to $\vphi_{\del}(\del)$,
where $\del<\psi_{\Ome}(\ome_{m+2}(\mI_{N}+1))<\psi_{\Ome}(\veps_{\mI_{N}+1})$.

This completes a proof of Theorem \ref{thm:2} assuming the well-foundedness of $OT(\mI_{N})$.

\section{Some ordinals in well-foundedness proof}\label{sect:prop}
In this section we introduce some ordinals needed in our well-foundedness proof.

In \cite{singlewfprf} the following Lemmas \ref{lem:of} and \ref{lem:oflx} are shown.
Lemma \ref{lem:of} is used in showing the finiteness of the sequence
$\rho_{0}\succ \rho_{1}\succ \rho_{2}\succ \cdots$, 
cf.\,Definition \ref{df:prec} and Lemma \ref{lem:LSwinding}.
Lemma \ref{lem:oflx} is needed in showing
Corollary \ref{cor:psiwf}.

\begin{definition}\label{df:Lam.of}
{\rm
Let $\Lam\leq\mI_{N}$ be a strongly critical number.
\benu
\item\label{df:Lam.of.1}
For $\xi<\Gamma(\Lam)$,
$a_{\Lam}(\xi)$ denotes an ordinal defined recursively by
$a_{\Lam}(0)=0$, and
$a_{\Lam}(\xi)=\sum_{i\leq m}\tilde{\theta}_{b_{i}}(\omega\cdot a_{\Lam}(\xi_{i});\Lam)$
when $\xi=_{NF}\sum_{i\leq m}\tilde{\theta}_{b_{i}}(\xi_{i};\Lam)$ in (\ref{eq:CantornfLam}).

\item\label{df:Lam.of.2}
For irreducible functions $f:\Lam\to\Gamma(\Lam)$ with base $\Lam$
 let us associate ordinals $o_{\Lam}(f)<\Gamma(\Lam)$ as follows.
$o_{\Lam}(\emptyset)=0$ for the empty function $f=\emptyset$.
Let $\{0\}\cup \supp(f)=\{0=c_{0}<c_{1}<\cdots<c_{n}\}$,
$f(c_{i})=\xi_{i}<\Gamma(\Lambda)$
for $i>0$, and $\xi_{0}=0$.
Define ordinals $\zeta_{i}=o_{\Lam}(f;c_{i})$ by
$\zeta_{n}=\omega\cdot a_{\Lam}(\xi_{n})$, and 
$\zeta_{i}=\omega\cdot a_{\Lam}(\xi_{i})+\tilde{\theta}_{c_{i+1}-c_{i}}(\zeta_{i+1}+1;\Lam)$.
Finally let $o_{\Lam}(f)=\zeta_{0}=o_{\Lam}(f;c_{0})$.

\item\label{df:Lam.of.3}
For $d\not\in \{0\}\cup \supp(f)$, let
$o_{\Lam}(f;d)=0$ if $f^{d}=\emptyset$.
Otherwise $o_{\Lam}(f;d)=\tilde{\theta}_{c-d}(o_{\Lam}(f;c)+1;\Lam)$ for
$c=\min( {\rm supp}(f^{d}))$.
\eenu
}
\end{definition}

\begin{lemma}\label{lem:of}
Let $f:\Lam \to\Gamma(\Lam)$ be an irreducible finite function with base $\Lam$
defined from an irreducible function $g:\Lam\to\Gamma(\Lam)$ and ordinals $c,d$
as follows.
$f_{c}\leq g_{c}$, 
$c<d\in \supp(g)$ with
$(c,d)\cap \supp(g)=
(c,d)\cap \supp(f)=
\emptyset$, 
$f(c)<g(c)+\tilde{\theta}_{d-c}(g(d);\Lam)\cdot\omega$, and
$f<_{\Lam}^{d}g(d)$, cf.\,Definition \ref{df:notationsystem}.\ref{df:notationsystem.6}.
Then $o_{\Lam}(f)<o_{\Lam}(g)$ holds.
\end{lemma}

\begin{lemma}\label{lem:oflx}
For irreducible finite functions $f,g:\Lam\to\Gamma(\Lam)$ with base $\Lam$, 
assume
$f<_{lx}^{0}g$.
Then $o_{\Lam}(f)<o_{\Lam}(g)$ holds.
\end{lemma}

\subsection{A preview of well-foundedness proof}\label{subsec:preview_wf}

To prove the well-foundedness of a computable notation system,
we utilize
the distinguished class introduced by W. Buchholz\cite{Buchholz75}.
Also cf.\,\cite{J2} for a well-foundedness in terms of a maximal distinguished class.

Let $OT$ be a computable notation system of ordinals with an ordinal term $\Omega_{1}$.
$\Omega_{1}$ denotes the least recursively regular ordinal $\omega_{1}^{CK}$.
Assume that
the well-founded part $W(OT)$ of $OT$ exists as a set.
A parameter-free $\Pi^{1-}_{1}\mathrm{+CA}$ suffices to show the existence.
Then the well-foundedness of such a notation system $OT$ is provable.
When the next recursively regular ordinal $\Omega_{2}$ is in $OT$,
we further assume that 
a well-founded part $W(\mathcal{C}^{\Omega_{1}}(W_{0}))$ of a set $\mathcal{C}^{\Omega_{1}}(W_{0})$ exists,
where $W_{0}=W(OT)\cap\Omega_{1}$, and
$\alpha\in \mathcal{C}^{\Omega_{1}}(W_{0})$ iff
each component$<\!\!\Omega_{1}$ of $\alpha$ is in $W_{0}$.
Likewise when $OT$ contains $\alpha$-many terms denoting increasing sequence of
recursively regular ordinals,
we need to iterate the process of defining the well-founded parts $\alpha$-times.

Let us consider a notation system $OT$ for recursively inaccessible universes.
There are $\alpha$-many ordinal terms denoting recursively regular ordinals in $OT$
with the order type $\alpha$ of $OT$.
The whole process then should be internalized.
We need to specify a feature of sets arising in the process.
Then \textit{distinguished sets} emerge. $D[P]$ denotes the fact that $P$ is a distinguished class and
defined by
\[
D[P]:\Lrarw
\forall\alpha
\left(
\alpha\leq P \to W(C^{\alpha}(P))\cap\alpha^{+}= P\cap\alpha^{+}
\right)
\]
where $\alp\leq P \Lrarw \exi \bet\in P(\alp\leq\bet)$ and
$\alp^{+}$ denotes the next recursively regular ordinal above $\alp$ if such an ordinal exists.

$W_{0}=W(OT)\cap\Ome_{1}$ is the smallest distinguished set, and $W_{1}=W(\mathcal{C}^{\Omega_{1}}(W_{0}))\cap\Omega_{2}$
 is the next one.
Given two distinguished sets, it turns out that one is an initial segment of the other,
and the union $\mathcal{W}_{0}=\bigcup\{P\subset OT: D[P]\}$ of all distinguished sets is distinguished, the \textit{maximal distinguished class}.
The maximal distinguished class $\mathcal{W}_{0}$ is $\Sigma^{1-}_{2}$-definable, and a proper class
without assuming $\Sigma^{1-}_{2}\mathrm{-CA}$.

Assuming the maximal distinguished class $\mathcal{W}_{0}$ exists as a set,
the well-foundedness of $OT$ for a single stable ordinal is provable in \cite{singlewfprf}.
Consider now a notation system $OT$ for several stable ordinals $\mS_{0},\mS_{1},\ldots$.
We then need several maximal distinguished sets 
$\mathcal{W}_{0},\mathcal{W}_{1},\ldots$ to prove the well-foundedness.
$\mathcal{W}_{0}$ is the maximal distinguished set in an absolute sense.

A moment reflection on the emergence of distinguished sets shows that
$\mathcal{W}_{1}$ could be a \textit{maximal distinguished set relative to} $\mathcal{W}_{0}$
and $\mS_{0}$.
Specifically cf.\,(\ref{eq:Di}),
a set $P$ is said to be a $0$-\textit{distinguished set} for $\gam$ and $X$, denoted by
$D^{\gam}[P;X]$, iff $P$ is well-founded and 
\[
P\cap\gam^{-\dagger}=X\cap\gam^{-\dagger} \spand 
\forall\alpha\geq \gam^{-\dagger}
\left(
\alpha\leq P \to W(C^{\alpha}(P))\cap\alpha^{+}= P\cap\alpha^{+}
\right)
\]
where 
$\gam^{-\dagger}=\max\{\mS\in St\cup\{0\}:\mS\leq\gam\}$.
Then let, cf.\,(\ref{eq:Wi+1})
\[
W_{1}^{\gam}(X)  :=  \bigcup\{P\subset OT  :D^{\gam}[P;X]\}.
\]
Observe that $W_{1}^{\gam}(X)$ is a $\Sig^{1}_{2}$-definable class, and hence a set
assuming $\Sig^{1}_{2}\mathrm{-CA}$.
We see in Lemma \ref{lem:3wf6}.\ref{lem:3wf6.3} that $W_{1}^{\gam}(X)$ is the maximal
$0$-distinguished class for $\gam$ and $X$ provided that $X\cap\gam^{-\dagger}$ is well-founded.

Assume that there are $\alp$-many stable ordinals with the order type $\alp$ of
a notation system $OT$ of ordinals.
Then we have to introduce distinguished sets in the next level.
In the higher level the recursive regularity is replaced by the stability, and
the $\Pi^{1}_{1}$-sets $W(\mathcal{C}^{\alp}(P))$ by $\Sig^{1}_{2}$-sets $\mathcal{W}_{1}^{\gam}(X)$.

A set $X$ is a $1$-\textit{distinguished set}, denoted by $D_{1}[X]$ iff $X$ is well-founded and
\[
\forall\gam
\left(
\gam\leq X\to \mathcal{W}_{1}^{\gam}(X)\cap\gam^{\dagger}= X\cap\gam^{\dagger}
\right).
\]
where $\alp^{\dagger}=\min\{\mS\in St: \alp<\mS\}$ if such a stable ordinal $\mS$ exists.
We see that $\mathcal{W}_{0}=W_{1}^{0}(\emptyset)$ is the smallest $1$-distinguished set, and
$\mathcal{W}_{1}=W_{1}^{\mS_{0}}(\mathcal{W}_{0})$ is the next $1$-distinguished set, and so forth.
In Lemma \ref{lem:6.15} it is shown that if $D_{1}[X]$ and $\gam\in X$, then
$X\cap\gamma^{\dagger}$ is a $0$-distinguished set for $\gam$ and $X$, i.e., $D^{\gam}[X\cap\gamma^{\dagger};X]$, 
and $\gam\in W(\mathcal{C}^{\gam}(X))\cap\gam^{+}=X\cap\gam^{+}$, where
$\gam\in W_{1}^{\gam}(X)\cap\gam^{\dagger}=X\cap\gam^{\dagger}$.
This crucial lemma allows us to prove facts by going down to the lowest level, i.e., to the well-foundedness.

$\mathcal{W} :=  \bigcup\{X\subset OT:D_{1}[X]\}$ is the $1$-maximal distinguished class, which is a 
$\Sig^{1-}_{3}$-definable class.
Although $\mathcal{W}$ is a proper class in a set theory with $\Pi_{1}$-Collection or
equivalently in $\Sig^{1}_{3}\mathrm{-DC+BI}$, the theories proves that
if $\mS\in\mathcal{W}$ for $\mS\in St\cup\{0\}$, then
$\mS^{\dagger}\in\mathcal{W}$, cf.\,Lemma \ref{cor:6.21}.
In showing that a limit of stable ordinals is in $\mathcal{W}$, we invoke 
$\Sig^{1}_{3}\mathrm{-DC}$ in Lemma \ref{lem:6.43}:
if $\alp\in\mathcal{G}^{\mathcal{W}}$, then 
there exists a $1$-distinguished \textit{set} $Z$ such that $Z$ is closed under
$\mS\mapsto\mS^{\dagger}$ and $\alp\in\mathcal{G}^{Z}$ for a $\Pi^{1}_{0}$-set $\mathcal{G}^{Z}$ 
in Definition \ref{df:calg} of subsection \ref{subsec:C(X)}.

By iterating this `jump' operators, we arrive at a $\Sig_{N+1}$-formula $D_{N}[X]$ for $N$-distinguished sets
for positive integers $N$, cf.\,Definition \ref{df:3wfdtg32}.
The maximal $N$-distinguished class $\bigcup\{X\subset OT:D_{N}[X]\}$ is $\Sig^{1-}_{N+2}$-definable proper class
in $\Pi_{N}$-Collection or in $\Sig^{1}_{N+2}\mathrm{-DC+BI}$.
\\

Up to this, everything seems to go well. 
But as long as we have an infinite increasing sequence $\{\mS_{n}\}_{n}$ of successor stable ordinals,
a technical difficulty is hidden as follows.
Above a successor stable ordinal $\mS_{0}$, there are increasing sequence
$\mS_{1}=\mS_{0}^{\dagger}, \mS_{2}=\mS_{1}^{\dagger},\ldots$ of successor stable ordinals.
Let $\rho_{n}\prec\mS_{n}$.
Let us define ordinals $\kap_{n,i}$ and $\sig_{n,i}$ for $i\leq n$ recursively by
$\kap_{n,n}=\mS_{n}$, $\kap_{n,i}=\kap_{n,i+1}[\rho_{i}/\mS_{i}]$,
$\sig_{n,n}=\rho_{n}$ and
$\sig_{i}=\sig_{n,i+1}[\rho_{i}/\mS_{i}]$.
Let $\kap_{n}=\kap_{n,0}$ and $\sig_{n}=\sig_{n,0}$.
Then we see that 
$\sig_{0}<\sig_{1}<\sig_{2}<\cdots<\kap_{2}<\kap_{1}<\kap_{0}$.
This might yield an infinite decreasing chain $\{\kap_{n}\}_{n}$ of collapsed ordinals.

For simplicity let $\rho_{i}=\psi_{\mS_{i}}^{f_{i}}(\alp_{i})$.
Then $M_{\rho_{i}}=C_{\alp_{i}}(\rho_{i})$.
In order to collapse $\kap_{n,i+1}$ by $\rho_{i}$, $\rho_{j}\in M_{\rho_{i}}$ has to be enjoyed for $j>i$.
Since $\rho_{j}>\rho_{i}$, this means that $\alp_{j}<\alp_{i}$.
Namely there must exist an infinite decreasing chain $\alp_{0}>\alp_{1}>\alp_{2}>\cdots$ in advance
to have another chain $\kap_{0}>\kap_{1}>\kap_{2}>\cdots$.
Here $\alp_{i}$ is the ordinal $\mathtt{p}_{0}(\rho_{i})$ in Definition \ref{df:p0}.\ref{df:p0.1}.
Let $\eta\in L(\mS)$ be an ordinal in the layer $L(\mS)$ of a successor stable ordinal $\mS$, 
cf.\,Definition \ref{df:precR}.
A pair $(\mathtt{g}_{1}(\eta),\mathtt{g}_{2}(\eta))$ of ordinals is associated with such an ordinal $\eta$
in Definitions \ref{df:gh.02} and \ref{df:gh.1}, and we show in
Lemma \ref{lem:LSwinding} that 
$(\mathtt{g}_{1}(\gam),\mathtt{g}_{2}(\gam))<_{lx}(\mathtt{g}_{1}(\eta),\mathtt{g}_{2}(\eta))$
when $\gam\in R(\eta)$ for the set $R(\eta)$ in Definition \ref{df:R0}.
It turns out that this suffices to prove the well-foundedness in Lemma \ref{lem:psimS}.

\subsection{Props}\label{subsec:prop}
In this subsection an ordinal $\mathtt{ p}_{\mS}(\alp)$ and
a pair $\mathtt{ g}(\alp)=(\mathtt{ g}_{1}(\alp),\mathtt{ g}_{2}(\alp))$
are introduced for ordinal terms $\alp$.
These are needed to show that there is no infinite sequence $\{\rho_{n},\kap_{n}\}_{n}$
such that $\rho_{0}\prec\mS_{0}$,
$\kap_{n}\in\{\mI_{N}[\rho_{n}]\}\cup\{\rho_{n}^{\dagger\vec{i}_{n}}, \mS_{n}^{\dagger \vec{i}_{n}}[\rho_{n}/\mS_{n}]\}$ and
either $\rho_{n+1}\prec \mS_{n}^{\dagger \vec{i}_{n}}[\rho_{n}/\mS_{n}]=\kap_{n}$ or $\rho_{n+1}\prec\tau^{\dagger \vec{i}_{n}}$
for $\tau\prec \mI_{N}[\rho_{n}]=\kap_{n}$,
cf.\,Proposition \ref{prp:LSwinding.2}, 
Lemmas \ref{lem:LSwinding} and \ref{lem:psimS}.

Recall that $\alp\in SSt^{M}$ iff either $\alp$ is a successor stable ordinal in $SSt$ or
$\alp=\bet[\rho/\mS]$ for a $\bet\in SSt^{M}$ and a successor stable ordinal $\mS$, 
cf.\,Definition \ref{df:notationsystem}.\ref{df:notationsystem.8}.

\bdf\label{df:Next}
{\rm
For $\rho\in\Psi_{\mS}$ with $\mS\in SSt^{M}$, let
$N(\rho)=\{\mI_{N}[\rho]\}\cup\{\rho^{\dagger\vec{i}}, \mS^{\dagger \vec{i}}[\rho/\mS]: \vec{i}\neq\emptyset\}\cap OT(\mI_{N})$ if $\mS\not\in SSt$.
Otherwise
$N(\rho)=\{\mI_{N}[\rho]\}\cup\{ \mS^{\dagger \vec{i}}[\rho/\mS]: \vec{i}\neq\emptyset\}\cap OT(\mI_{N})$.
}
\edf

Note that $\rho^{\dagger\vec{i}}\in SSt$ when $\mS\in SSt$, and $N(\rho)\cap\Psi=\emptyset$.
Recall that, cf.\,Definition \ref{df:precR}, $L(\mS)$ denotes the layer of $\mS$, and 
$\alp\in L(\mS)$ iff $\alp\prec^{R}\mS$ iff there are ordinals $\{\rho_{i},\kap_{i}\}_{i}$ such that
$\kap_{0}=\mS$, $\rho_{i}\prec\kap_{i}$ , $\kap_{i+1}\in N(\rho_{i})$, and 
$\alp\in\{\rho_{0}\}\cup\{\rho_{i},\kap_{i}\}_{i>0}$.

\bdf\label{df:ps}
{\rm
Let $\mS\in SSt_{i}$ and $\mT\in St\cup\{0\}$ be the least such that
$\mS=\mT^{\dagger i}$.
For $a\in OT(\mI_{N})$,
the \textit{prop} $\mathtt{ p}_{\mS}(a)$ of $a$ 
denotes an ordinal term defined recursively as follows.

 \benu
 \item
 $\mathtt{ p}_{\mS}(\mI_{N})=\mathtt{ p}_{\mS}(a)=0$ if $a\leq\mT$
 
 In what follows assume $\mI_{N}\neq a>\mT$.
 \item
 $\mathtt{ p}_{\mS}(a)=\max_{i\leq m} \mathtt{ p}_{\mS}(a_{i})$ if
 $a=a_{0}+\cdots+a_{m}$.
 
 $\mathtt{ p}_{\mS}(a)=\max\{\mathtt{ p}_{\mS}(b),\mathtt{ p}_{\mS}(c)\}$ if $a=\vphi bc$.
 
 \item
 $\mathtt{ p}_{\mS}(a)=\mathtt{ p}_{\mS}(\kap)$ if 
 $a\in N(\kap)$
 for a $\kap\in L(\mU)\cap\Psi$ with a $\mU>\mT$.

 \item
 $\mathtt{ p}_{\mS}(\mU^{\dagger k})=\mathtt{ p}_{\mS}(\mU)$ for $\mT\leq\mU\in St$.

 \item
 $\mathtt{ p}_{\mS}(\psi_{\mI_{N}}(a))=\mathtt{ p}_{\mS}(a)$. 
 
 \item
 For $\mathtt{ p}_{\mS}(\mathrm{fld}(f))=\max\{\mathtt{ p}_{\mS}(b): b\in \mathrm{fld}(f)\}$, let
 \[
 \mathtt{ p}_{\mS}(\psi_{\kap}^{f}(a))=
 \left\{
 \begin{array}{ll}
 \max\{\mathtt{ p}_{\mS}(\kap),\mathtt{ p}_{\mS}(a),\mathtt{ p}_{\mS}(\mathrm{fld}(f))\}
 & \mbox{if  } \kap>\mS
 \\
 \max\{a,\mathtt{ p}_{\mS}(a)\}
 & \mbox{if  } \kap=\mS
 \\
 \mathtt{ p}_{\mS}(\kap)
 & \mbox{if  } 
 \kap<\mS
 \end{array}
 \right.
 \]
 
 \eenu
}
\edf

\bprp\label{prp:ttb}
Let $\mS\in SSt$ and $\alp=\psi_{\mS}^{f}(a)$,
$\bet=\psi_{\mS}^{g}(b)$ with
$\{\alp,\bet\}\subset OT(\mI_{N})$,

\benu
\item\label{prp:ttb.0}
Let $c\in C_{b}(\bet)$ with
$\mathtt{ p}_{\mS}(c)\neq 0$.
Then there exists a subterm 
$\gam\in C_{b}(\bet)$ of $c$ such that
$\gam\prec\mS$ and $\mathtt{ p}_{\mS}(\gam)=\mathtt{ p}_{\mS}(c)$.

\item\label{prp:ttb.01}
$\mathtt{ p}_{\mS}(\mathrm{fld}(f))\leq\max\{a,\mathtt{ p}_{\mS}(a)\}=\mathtt{ p}_{\mS}(\alp)$ holds.

\item\label{prp:ttb.1}
$\mathtt{ p}_{\mS}(\bet)\leq\mathtt{ p}_{\mS}(\alp)$
if $\bet<\alp$.

\item\label{prp:ttb.2}
Let $\del<\alp<\bet$ with $\del\prec\bet$.
Then
$\mathtt{ p}_{\mS}(\bet)\leq\mathtt{ p}_{\mS}(\alp)$.

\item\label{prp:ttb.3}
Let $\{\gam,\del\}\subset OT(\mI_{N})$.
Then $\mathtt{ p}_{\mS}(\gam)\leq\mathtt{ p}_{\mS}(\del)$ if $\gam<\del$.
\eenu
\eprp
\bprf
\ref{prp:ttb}.\ref{prp:ttb.0}.
By induction on $\ell c$.
\\
\ref{prp:ttb}.\ref{prp:ttb.01}.
By (\ref{eq:notationsystem.5}) in 
Definition \ref{df:notationsystem}.\ref{df:notationsystem.5} we obtain
$\mathrm{fld}(f)\subset C_{0}(SC(a))$.

By induction on $\ell b$,
we see
$b\in C_{0}(SC(a)) \Rarw 
\mathtt{ p}_{\mS}(b)\leq \max\{a,\mathtt{ p}_{\mS}(a)\}$.
\\

\noindent
We show Propositions \ref{prp:ttb}.\ref{prp:ttb.1} and \ref{prp:ttb}.\ref{prp:ttb.2}
simultaneously by induction on $\ell\bet+\ell\alp$.
\\
\ref{prp:ttb}.\ref{prp:ttb.1}.
If $a=b$, then $\mathtt{ p}_{\mS}(\bet)=\mathtt{ p}_{\mS}(\alp)$.
Let $b<a$. We can assume $a<c=\mathtt{ p}_{\mS}(b)$.
By Proposition \ref{prp:ttb}.\ref{prp:ttb.0} pick a shortest subterm 
$\gam\in C_{b}(\bet)\cap\mS\subset\bet$ of $b$
such that $\gam\prec\mS$ and $\mathtt{ p}_{\mS}(\gam)=\mathtt{ p}_{\mS}(b)=c$
for $b\in C_{b}(\bet)$.
Then $\gam\preceq\del=\psi_{\mS}^{h}(c)$ for some $h$ and $\gam<\bet$.
If $\del\leq\alp$, then IH yields $c=\mathtt{ p}_{\mS}(\del)\leq\mathtt{ p}_{\mS}(\alp)$.
Assume $\gam<\alp<\del$ with $\gam\prec\del$.
IH for Proposition \ref{prp:ttb}.\ref{prp:ttb.2} then yields 
$c\leq\mathtt{ p}_{\mS}(\alp)$.

Next let $a<b$. 
Pick a subterm $\eta$ of a term in $\{a\}\cup \mathrm{fld}(f)$ such that 
$\bet\leq\eta\in C_{a}(\alp)$ and
$\eta\prec\mS$.
Let $\eta\preceq\psi_{\mS}^{h}(d)=\sig$ for some $h$ and $d$.
Then we obtain $\bet\leq\sig$, and IH yields
$\mathtt{ p}_{\mS}(\bet)\leq\mathtt{ p}_{\mS}(\sig)=\mathtt{ p}_{\mS}(\eta)$.
On the other hand we have
$\mathtt{ p}_{\mS}(\eta)\leq\max\{a,\mathtt{ p}_{\mS}(a)\}$
by Proposition \ref{prp:ttb}.\ref{prp:ttb.01}.
Hence $\mathtt{ p}_{\mS}(\bet)\leq\mathtt{ p}_{\mS}(\alp)$.
\\
\ref{prp:ttb}.\ref{prp:ttb.2}.
Pick a subterm $\eta$ of a term in $\{a\}\cup \mathrm{fld}(f)$ such that 
$\del\leq\eta\in C_{a}(\alp)$,
$\eta\prec\mS$ and
$\mathtt{ p}_{\mS}(\eta)\leq\max\{a,\mathtt{ p}_{\mS}(a)\}$ by Proposition \ref{prp:ttb}.\ref{prp:ttb.01}.
Let $\eta\preceq\psi_{\mS}^{h}(d)=\sig$ for some $h$ and $d$.
Then we obtain $\del\leq \sig$.
If $\bet\leq\sig$, then IH for Proposition \ref{prp:ttb}.\ref{prp:ttb.1} yields 
$\mathtt{ p}_{\mS}(\bet)\leq\mathtt{ p}_{\mS}(\sig)=\mathtt{ p}_{\mS}(\eta)\leq\mathtt{ p}_{\mS}(\alp)$.
Otherwise we obtain $\del\leq\sig<\bet$ with $\del\prec\bet$.
IH yields $\mathtt{ p}_{\mS}(\bet)\leq\mathtt{ p}_{\mS}(\sig)\leq\mathtt{ p}_{\mS}(\alp)$.
\\

\noindent
\ref{prp:ttb}.\ref{prp:ttb.3}.
This is seen by induction on $\ell\gam+\ell\del$ using Definition \ref{df:lessthanpsi}, and Propositions
\ref{prp:ttb}.\ref{prp:ttb.1} and \ref{prp:ttb}.\ref{prp:ttb.2}.
\eprf
\\

The set $Cr$ of strongly critical numbers in $OT(\mI_{N})$ is divided to
$Cr=St\cup \bigcup\{L(\mS):\mS\in SSt\}\cup(Cr\cap(\Ome+1))$,
cf.\,Definition \ref{df:precR}.

\bdf\label{df:gh.02}
{\rm
Let $\mS\in SSt$ and $\alp\in L(\mS)$.
Let us define ordinals $\mathtt{ g}_{0}(\alp), \mathtt{ g}^{*}_{0}(\alp)$ and $\mathtt{ g}_{2}(\alp)$
 as follows.

\benu
\item
$\mathtt{ g}_{0}(\alp)=\mathtt{ g}_{2}(\alp)=0$ for $\alp\not\in\Psi$.

\item
If $\rho\prec\mS$,
then let
$\mathtt{ g}_{0}(\rho)=\mathtt{ g}^{*}_{0}(\rho)=\mathtt{ g}_{0}(\psi_{\mI_{N}[\rho]}(b))=
\mathtt{ g}^{*}_{0}(\psi_{\mI_{N}[\rho]}(b))=\mathtt{ p}_{\mS}(\rho)$
for every $b$.
Also $\mathtt{ g}_{2}(\rho)=o_{\mI_{N}}(m(\rho))+1$ for
$m(\rho):\mI_{N} \to \Gamma(\mI_{N})$ with base $\mI_{N}$, and
$\mathtt{ g}_{2}(\psi_{\mI_{N}[\rho]}(b))=0$.

\item
Let $\rho\prec\mS$ and 
$\alp\prec^{R}\tau\in N(\rho)$,
where $\alpha\neq\psi_{\mI_{N}[\rho]}(a)$ for any $b$
if $\tau=\mI_{N}[\rho]$.
Let $\mathtt{ g}^{*}_{0}(\alp)=\mathtt{ g}_{0}(\rho)=\mathtt{ p}_{\mS}(\rho)$.
Let $\bet\in M_{\rho}$ be such that $\alp=\bet[\rho/\mS]$.
If $\alp\in\Psi$, let
$\mathtt{ g}_{i}(\alp)=\mathtt{ g}_{i}(\bet)$ for $i=0,2$.

\eenu

}
\edf

\bprp\label{prp:g01}
Let $b=\mathtt{ p}_{0}(\alp)$ for $\alp\in L(\mS)\cap\Psi$ with $\mS\in SSt$.
Then
$SC(\mathtt{ g}_{2}(\alp))\subset \psi_{\mI_{N}}(b)$.
Moreover
$\mathtt{p}_{0}(\alp)\leq\mathtt{g}^{*}_{0}(\alp)$.
\eprp
\bprf
By induction on $\ell\alp$.
Cf.\,Definition \ref{df:p0}.\ref{df:p0.1} for $\mathtt{p}_{0}(\alp)$.
\\
\textbf{Case 1}.
First let $\alp\preceq\psi_{\mS}^{g}(b)$ with an $\mS\in SSt$ and $f=m(\alp)$.
By Proposition \ref{prp:Stclass}.\ref{prp:Stclass.2} let $\mT\in LSt\cup\{0\}$ be such that $\mS=\mT^{\dagger \vec{i}}$ for a sequence $\vec{i}$.
We obtain $SC(\mathtt{ g}_{2}(\alp))\subset \mathrm{fld}(f)$
for $\mathtt{ g}_{2}(\alp)=o_{\mI_{N}}(f)+1$.
By 
(\ref{eq:notationsystem.6}) in Definition \ref{df:notationsystem}
we obtain
$\mathrm{fld}(f)\subset M_{\alp}\cap\mI_{N}=C_{b}(\alp)\cap\mI_{N}$.
 On the other hand we have $\mathtt{p}_{0}(\alp)=b\leq\mathtt{p}_{\mS}(\alp)=\mathtt{g}^{*}_{0}(\alp)$.

We claim that $\alp<\psi_{\mI_{N}}(b)$.
$SC(\mathtt{ g}_{2}(\alp))\subset 
 C_{b}(\psi_{\mI_{N}}(b))\cap\mI_{N}\subset\psi_{\mI_{N}}(b)$
 follows from the claim.
For the claim it suffices to show $\mS<\psi_{\mI_{N}}(b)$.
Let $\{(\mT_{m},\mS_{m},\vec{i}_{m})\}_{m\leq n}$ be the sequence
such that
$\mT_{0}\in LSt_{N}\cup\{0\}$,
$\mS_{m}=\mT_{m}^{\dagger\vec{i}_{m}}$ and $\mT_{m+1}\prec\mS_{m}\,(m<n)$, and
$\mS=\mS_{n}$, cf.\,the trail to $\mS$ in Proposition \ref{prp:trail}.
If $\mT_{0}=0$, then $\mS\leq\mS_{0}<\psi_{\mI_{N}}(b)\in LSt_{N}$.
Let $\mT_{0}=\psi_{\mI_{N}}(c)$.
Proposition \ref{prp:Stclass}.\ref{prp:Stclass.3} yields $c<b$, and
$\mT_{0}=\psi_{\mI_{N}}(c)<\psi_{\mI_{N}}(b)\in LSt_{N}$.
Hence $\mS\leq\mS_{0}<\psi_{\mI_{N}}(b)\in LSt_{N}$.
\\
\textbf{Case 2}.
Next let $LSt_{i}\ni\rho\prec\mS\in SSt$,
$\alp\prec^{R}\tau\in N(\rho)$ and
$\alp=\bet[\rho/\mS]$ for a 
$\bet\in M_{\rho}$.
Then $b=\mathtt{ p}_{0}(\alp)=\mathtt{ p}_{0}(\bet)$, and 
IH yields $SC(\mathtt{ g}_{2}(\alp))=SC(\mathtt{ g}_{2}(\bet))\subset \psi_{\mI_{N}}(b)$
for
$\mathtt{ g}_{2}(\alp)=\mathtt{ g}_{2}(\bet)$, and
$\mathtt{ p}_{0}(\bet)\leq\mathtt{g}^{*}_{0}(\bet)$.

On the other hand we have $\mathtt{g}^{*}_{0}(\alp)=\mathtt{g}_{0}(\rho)=\mathtt{p}_{\mS}(\rho)\geq\mathtt{p}_{0}(\rho)=c$ with $M_{\rho}=C_{c}(\rho)$.
Thus it suffices to show $\mathtt{g}^{*}_{0}(\del)\leq\mathtt{p}_{0}(\rho)$ for
$\del\in C_{c}(\rho)$ by induction on $\ell\del$.
If $\del\preceq\psi_{\mT}^{f}(d)$ with a $\mT\in SSt$, then
$\mathtt{g}^{*}_{0}(\del)=\mathtt{g}_{0}(\del)=\mathtt{p}_{\mT}(\del)=\max\{d,\mathtt{p}_{T}(d)\}$.
We obtain $d<c$ and $d\in C_{c}(\rho)$.
IH yields $\mathtt{p}_{\mT}(d)<c$.

Next let $\del=\gam[\tau/\mT]$ with a $\gam\in M_{\tau}$.
Then $\mathtt{g}^{*}_{0}(\del)=\mathtt{g}^{*}_{0}(\tau)$ and $\tau\in M_{\rho}$.
IH yields $\mathtt{g}^{*}_{0}(\tau)\leq\mathtt{p}_{0}(\rho)$.
\eprf

\bprp\label{prp:LSwinding.1}
Let $\tau\in L(\mS)\cup\{\mS\}$
and $\mS\in SSt$.

For $\rho, \eta\prec\tau$,
if $\rho<\eta$, then
$\mathtt{ g}_{0}(\rho)\leq\mathtt{ g}_{0}(\eta)$.
\eprp
\bprf 
By induction on $\ell\rho$.
\\
\textbf{Case 1}. $\tau=\mS$:
Let $\eta\preceq\alp=\psi_{\mS}^{f}(a)$ and $\rho\preceq\bet=\psi_{\mS}^{g}(b)$.
Then $\mathtt{ g}_{0}(\eta)=\mathtt{ p}_{\mS}(\eta)=\mathtt{ p}_{\mS}(\alp)$ and
$\mathtt{ g}_{0}(\rho)=\mathtt{ p}_{\mS}(\rho)=\mathtt{ p}_{\mS}(\bet)$.
If $\bet\leq\alp$, then Proposition \ref{prp:ttb}.\ref{prp:ttb.1} yields
$\mathtt{ p}_{\mS}(\bet)\leq\mathtt{ p}_{\mS}(\alp)$.
Suppose $\rho<\alp<\bet$ with $\rho\prec\bet$.
We obtain $\mathtt{ p}_{\mS}(\bet)\leq\mathtt{ p}_{\mS}(\alp)$
by Proposition \ref{prp:ttb}.\ref{prp:ttb.2}.
\\
\textbf{Case 2}. $\tau\neq\mS$:
Let $\kap\prec\mS$ be such that either $\tau\preceq^{R}\mS^{\dagger \vec{i}}[\kap/\mS]$ or
$\tau\prec^{R}\mI_{N}[\kap]$.
Then $\mathtt{ g}_{0}(\rho)=\mathtt{ g}_{0}(\rho_{1})$ and
$\mathtt{ g}_{0}(\eta)=\mathtt{ g}_{0}(\eta_{1})$ for
$\rho_{1}=\rho[\kap/\mS]^{-1}$ and $\eta_{1}=\eta[\kap/\mS]^{-1}$,
cf.\,Definition \ref{df:divide} for uncollapsing.
We obtain $\rho_{1}<\eta_{1}$, $\rho_{1}\prec\tau_{1}$
and $\eta_{1}\prec\tau_{1}$ for
$\tau_{1}=\tau[\kap/\mS]^{-1}$.
IH with $\ell\rho_{1}<\ell\rho$ yields
$\mathtt{ g}_{0}(\rho_{1})\leq\mathtt{ g}_{0}(\eta_{1})$.
\eprf

\bprp\label{prp:LSwinding.2}
Let $\mS\in SSt$,
$\rho\prec\tau\in (L(\mS)\cup\{\mS\})\cap SSt^{M}$, and
$\alp\prec\sig\in SSt^{M}$, where 
$\sig\preceq^{R}\kap\in N(\rho)$.
Then
$\mathtt{ g}_{0}(\alp)<\mathtt{ g}_{0}(\rho)$.
\eprp
\bprf
We may assume that either $\sig=\kap=\mS^{\dagger \vec{i}}[\rho/\mS]$ or
$\kap=\mI_{N}[\rho] \spand \sig=(\psi_{\mI_{N}[\rho]}(\gam))^{\dagger \vec{i}}$ for a $\gam$ and
an $\vec{i}$.
By induction on $\ell\alp$ we show $\mathtt{ g}_{0}(\alp)<\mathtt{ g}_{0}(\rho)$.
\\
\textbf{Case 1}. $\rho\prec\mS$:
Let $\rho\preceq\bet=\psi_{\mS}^{g}(b)$.
Then $\mathtt{ g}_{0}(\rho)=\mathtt{ p}_{\mS}(\psi_{\mS}^{g}(b))$.
From $b\in C_{b}(\psi_{\mS}^{g}(b))$
we see $\mathtt{ p}_{\mT}(b)<\mathtt{ p}_{\mT}(\psi_{\mT}(b))=b\leq\mathtt{ p}_{\mS}(\psi_{\mS}^{g}(b))=\mathtt{ g}_{0}(\rho)$
for any $\mS<\mT\in SSt$.
\\
\textbf{Case 1.1}. $\sig=\mS^{\dagger \vec{i}}[\rho/\mS]$ :
Let
$\alp\preceq\psi_{\sig}^{h_{1}}(c_{1})=
\left(
\psi_{\mS^{\dagger \vec{i}}}^{h}(c)
\right)[\rho/\mS]$, where
$h_{1}=h[\rho/\mS]\neq\emptyset$, $c_{1}=c[\rho/\mS]$ and 
$\sig=\mS^{\dagger \vec{i}}[\rho/\mS]=(\mS^{\dagger \vec{i}})[\rho/\mS]$.
Then $\mathtt{ g}_{0}(\alp)=\mathtt{ p}_{\mS^{\dagger \vec{i}}}(\psi_{\mS^{\dagger \vec{i}}}^{h}(c))$.
We have 
$\rho<\psi_{\mS^{\dagger \vec{i}}}^{h}(c)\in M_{\rho}=C_{b}(\rho)$,
and hence $c<b$.
We obtain
$\mathtt{ p}_{\mS^{\dagger \vec{i}}}(c)\leq\mathtt{ p}_{\mS^{\dagger \vec{i}}}(b)$
by Proposition \ref{prp:ttb}.\ref{prp:ttb.3}.
\\
\textbf{Case 1.2}. $\sig=(\psi_{\mI_{N}[\rho]}(\gam_{1}))^{\dagger \vec{i}}$ for a $\gam_{1}$:
Let $\alp\preceq\psi_{\sig}^{h_{1}}(c_{1})=
\left(
\psi_{\mT^{\dagger \vec{i}}}^{h}(c)
\right)[\rho/\mS]$, where
$h_{1}=h[\rho/\mS]\neq\emptyset$, $c_{1}=c[\rho/\mS]$ and 
$\sig=\mT[\rho/\mS]$ with $\mT=\psi_{\mI_{N}}(\gam)$ and 
$\gam_{1}=\gam[\rho/\mS]$.
Then $\mathtt{ g}_{0}(\alp)=\mathtt{ p}_{\mT^{\dagger \vec{i}}}(\psi_{\mT^{\dagger \vec{i}}}^{h}(c))$.
We have 
$\psi_{\mT^{\dagger \vec{i}}}^{h}(c)\in M_{\rho}$.
As in \textbf{Case 1.1} we see 
$c<b$ and
$\mathtt{ p}_{\mT^{\dagger \vec{i}}}(c)\leq\mathtt{ p}_{\mT^{\dagger \vec{i}}}(b)$ from $\mS<\mT^{\dagger \vec{i}}$, i.e., from
$\mS<\mT=\psi_{\mI_{N}}(\gam)\in LSt_{N}$.
\\
\textbf{Case 2}. $\rho\prec\tau\neq\mS$:
Let $\lam\prec\mS$ be such that either $\alp\prec^{R}\mS^{\dagger \vec{j}}[\lam/\mS]$ or
$\alp\prec^{R}\mI_{N}[\lam]$.
Then $\mathtt{ g}_{0}(\alp)=\mathtt{ g}_{0}(\alp_{1})$ with $\alp=\alp_{1}[\lam/\mS]$
and $\mathtt{ g}_{0}(\rho)=\mathtt{ g}_{0}(\rho_{1})$ with $\rho=\rho_{1}[\lam/\mS]$.
We have $\rho_{1}\prec \tau[\lam/\mS]^{-1}$ and
$\alp_{1}\prec \sig[\lam/\mS]^{-1}$
with $\sig=\mS^{\dagger \vec{i}}[\rho/\mS]$ or $\sig=(\psi_{\mI_{N}[\rho]}(\gam))^{\dagger\vec{i}}$ for a $\gam$.
If $\tau[\lam/\mS]^{-1}\in SSt$, then
we obtain $\mathtt{ g}_{0}(\alp_{1})<\mathtt{ g}_{0}(\rho_{1})$ by \textbf{Case 1}.
Otherwise IH with $\ell\alp_{1}<\ell\alp$ yields the proposition.
\eprf

\bprp\label{prp:LSwinding.3}
Let $\{\alp, \bet\}\subset L(\mS)$ with an $\mS\in SSt$.
If $\alp<\bet$, then
$\mathtt{ g}^{*}_{0}(\alp)\leq\mathtt{ g}^{*}_{0}(\bet)$.
\eprp
\bprf
Let $\rho,\eta\prec\mS$ be such that
either $\alp=\rho$ or $\alp\preceq^{R}\kap\in N(\rho)$,
and either $\bet=\eta$ or 
$\bet\preceq^{R}\sig\in N(\eta)$.
Then $\rho\leq\eta$ by $\alp<\bet$.
Proposition \ref{prp:LSwinding.1} yields
$\mathtt{ g}^{*}_{0}(\alp)=\mathtt{ g}_{0}(\rho)\leq\mathtt{ g}_{0}(\eta)=\mathtt{ g}^{*}_{0}(\bet)$.
\eprf

\bdf\label{df:R0}
{\rm
A 
set $R(\eta)\subset\Psi$ is defined.
\benu
\item\label{df:R00}
Let $\eta\prec\mI_{N}$.
$\gam\in R(\eta)$ holds iff
there exists an $SSt\ni\mS<\eta$ such that $\gam\in L(\mS)\cap\Psi$.

\item\label{df:R01}
Let 
$\eta\in L(\mS)$ with an $\mS\in SSt$.
$\gam\in R(\eta)$ holds iff $\gam\in\Psi$ and one of the following holds:

\benu
\item\label{df:R0.0}
$\gam\prec^{R}\eta$.

\item\label{df:R0.100}
There exist $\tau\in L(\mS)\cap\Psi$ and $\{\kappa<\rho\}\subset N(\tau)$ such that
$\eta\preceq^{R}\rho$ and $\gamma\prec^{R}\kappa$.

\item\label{df:R0.2}
There exist $\tau\in L(\mS)$, $\rho$ and $\vec{i}$ such that
$\eta,\rho\prec\mI_{N}[\tau]$, $\rho<\eta$ and
$\gam\prec^{R}\rho^{\dagger \vec{i}}$.

\item\label{df:R0.3}
There exist  $\tau\in (L(\mS)\cup\{\mS\})\cap SSt^{M}$, $\rho$ and $\kap$ such that
$\eta,\rho\prec\tau$, $\rho<\eta$, 
$\gam\prec^{R}\kap\in N(\rho)$.

\eenu
\eenu
}
\edf

\bprp\label{prp:lessthan_layer}
Let $\mS\in SSt$.
\benu
\item\label{prp:lessthan_layer.0}
Let $\rho\prec\kappa\in SSt^{M}$ and $\tau\in N(\rho)$.
Then $\tau<\kappa$.

\item\label{prp:lessthan_layer.1}
If $\alpha\prec^{R}\kappa$, then $\alpha<\kappa$.

\item\label{prp:lessthan_layer.25}
Let $\kappa\in L(\mS)\cap\Psi$, $\{\rho<\tau\}\subset N(\kappa)$ and $\alpha\prec\tau$.
Then $\rho<\alpha$.

\item\label{prp:lessthan_layer.3}
Let $\rho\in L(\mS)\cap\Psi$, $\kappa\in N(\rho)$ and $\alpha\prec\kappa$.
Then $\rho<\alpha$.

\item\label{prp:lessthan_layer.4}
Let $\rho,\tau\in L(\mS)\cap\Psi$ and $\kappa,\sigma_{0}\in (L(\mS)\cup\{\mS\})\cap SSt^{M}$ be
such that $\rho<\tau\prec\sigma_{0}$ and $N(\rho)\ni\kappa<\sigma_{0}$. 
Then $\kappa<\tau$.

\item\label{prp:lessthan_layer.7}
$\gamma<\eta$ for $\gamma\in R(\eta)$.

\item\label{prp:lessthan_layer.8}
Let $\rho\in R(\eta)\cap\Psi$ and $\gamma\prec\kappa\in N(\rho)$. Then $\gamma\in R(\eta)$.

\eenu
\eprp
\bprf
These are seen by an inspection to Definition \ref{df:lessthanpsi}, and follow from
the facts that
$\alpha\prec\beta\Rarw\alpha<\beta$ and
 $\gamma<\delta\Rarw \gamma[\rho/\mS]<\delta[\rho/\mS]<\mS$ for $\gamma,\delta\in M_{\rho}$ with $\rho\in\Psi_{\mS}$
 as follows.
 \\
 \ref{prp:lessthan_layer}.\ref{prp:lessthan_layer.0} by induction on $\ell\tau$.
If $\kappa=\mS\in SSt$, then $N(\rho)\ni\tau<\mS$.
Otherwise let $\kappa[\sigma/\mS]^{-1}\in SSt^{M}$ for a $\sigma\in\Psi_{\mS}$ and an $\mS\in SSt$.
Then $\rho[\sigma/\mS]^{-1}\prec\kappa[\sigma/\mS]^{-1}$ and $\tau[\sigma/\mS]^{-1}\in N(\rho[\sigma/\mS]^{-1})$.
IH yields $\tau[\sigma/\mS]^{-1}<\kappa[\sigma/\mS]^{-1}$. We obtain $\tau<\kappa$.
\\
\ref{prp:lessthan_layer}.\ref{prp:lessthan_layer.1}. By induction on $\ell\alpha$ using Proposition \ref{prp:lessthan_layer}.\ref{prp:lessthan_layer.0}.
\\
\ref{prp:lessthan_layer}.\ref{prp:lessthan_layer.25} by induction on $\ell\kappa$.
We may assume that $\kappa\prec\mS$.
Consider the case when $\tau=\mI_{N}[\kappa]$ and $\rho=\mS^{\dagger \vec{j}}[\kappa/\mS]$.
Let $\alpha\prec\mI_{N}[\kappa]$.
Then $\alpha[\kappa/\mS]^{-1}\prec(\mI_{N}[\kappa])[\kappa/\mS]^{-1}=\mI_{N}$ and $\mS<\alpha[\kappa/\mS]^{-1}\in LSt$.
We obtain $(\mS^{\dagger \vec{j}}[\kappa/\mS])[\kappa/\mS]^{-1}=\mS^{\dagger j}<\alpha[\kappa/\mS]^{-1}$, and $\mS^{\dagger \vec{j}}[\kappa/\mS]<\alpha$.
\\
\ref{prp:lessthan_layer}.\ref{prp:lessthan_layer.3} by induction on $\ell\alpha$.
Let $\rho\in L(\mS)\cap\Psi$, $\kappa\in N(\rho)$ and $\alpha\prec\kappa$.
We show $\rho<\alpha$. Let $\kappa_{0}\in (L(\mS)\cup\{\mS\})\cap SSt^{M}$ be such that $\rho\prec\kappa_{0}$.
If $\kappa_{0}=\mS$, then we obtain $\mS<\alpha[\rho/\mS]^{-1}\prec\kappa[\rho/\mS]^{-1}$, and
$\rho=\mS[\rho/\mS]<\alpha$.
Otherwise let $\rho_{1}\prec\mS$ be such that $\rho\prec^{R}\kappa_{1}$ for a $\kappa_{1}\in N(\rho_{1})$.
We obtain $\alpha[\rho_{1}/\mS]^{-1}\prec\kappa[\rho_{1}/\mS]^{-1}\in N(\rho[\rho_{1}/\mS]^{-1})$.
IH yields $\rho[\rho_{1}/\mS]^{-1}<\alpha[\rho_{1}/\mS]^{-1}$, and $\rho<\alpha$.
\\
\ref{prp:lessthan_layer}.\ref{prp:lessthan_layer.4} by induction on $\ell\rho+\ell\tau$.
Let $\rho,\tau\in L(\mS)\cap\Psi$ and $\kappa,\sigma_{0}\in (L(\mS)\cup\{\mS\})\cap SSt^{M}$ be
such that $\rho<\tau\prec\sigma_{0}$ and $N(\rho)\ni\kappa<\sigma_{0}$. 
We show $\kappa<\tau$.
Let $\rho\prec\kappa_{0}\in (L(\mS)\cup\{\mS\})\cap SSt^{M}$.
\\
\textbf{Case 1}. $\kappa_{0}=\sigma_{0}$: Let $\tau=\psi_{\pi}^{f}(a)$. If $\pi=\sigma_{0}$, then we obtain $\kappa\in C_{a}(\tau)\cap\sigma_{0}=\tau$ by
$\rho\in C_{a}(\tau)$. Otherwise IH yields $\kappa<\pi$, and $\kappa\in C_{a}(\tau)\cap\pi=\tau$.
\\
\textbf{Case 2}. $\kappa_{0}<\sigma_{0}$: Let $\kappa_{0}\in N(\rho_{1})$. We obtain
$\rho_{1}<\rho<\tau$ by Proposition \ref{prp:lessthan_layer}.\ref{prp:lessthan_layer.3}.
On the other hand we have $\kappa<\kappa_{0}$ by Proposition \ref{prp:lessthan_layer}.\ref{prp:lessthan_layer.0}.
 IH yields $\kappa<\kappa_{0}<\tau$.
\\
\textbf{Case 3}. $\kappa_{0}>\sigma_{0}$:
Let $\sigma_{0}\in N(\tau_{1})$.
\\
\textbf{Case 3.1}. $\rho>\tau_{1}$: IH with $\sigma_{0}<\kappa_{0}$ yields $\tau<\sigma_{0}<\rho$. 
This is not the case
\\
\textbf{Case 3.2}. $\rho=\tau_{1}$: 
Proposition \ref{prp:lessthan_layer}.\ref{prp:lessthan_layer.25} with $\kappa<\sigma_{0}$ yields $\kappa<\tau$.
\\
\textbf{Case 3.3}. $\rho<\tau_{1}$: Let $\sigma_{1}\in (L(\mS)\cup\{\mS\})\cap SSt^{M}$ be such that $\tau_{1}\prec\sigma_{1}$.
We have $\kappa<\sigma_{0}<\sigma_{1}$ by Proposition \ref{prp:lessthan_layer}.\ref{prp:lessthan_layer.0}.
On the other hand we have $\tau_{1}<\tau$ by Proposition \ref{prp:lessthan_layer}.\ref{prp:lessthan_layer.3}.
IH yields $\kappa<\tau_{1}<\tau$.
\\
\ref{prp:lessthan_layer}.\ref{prp:lessthan_layer.7}. 
By Propositions \ref{prp:lessthan_layer}.\ref{prp:lessthan_layer.0}, \ref{prp:lessthan_layer}.\ref{prp:lessthan_layer.1},
\ref{prp:lessthan_layer}.\ref{prp:lessthan_layer.25},
\ref{prp:lessthan_layer}.\ref{prp:lessthan_layer.3} and \ref{prp:lessthan_layer}.\ref{prp:lessthan_layer.4},
we may assume that either $\eta\prec\mI_{N}$ by Definition \ref{df:R0}.\ref{df:R00} or $\eta\prec\mI_{N}[\tau]$ for a $\tau\in L(\mS)\cap\Psi$ with $\mS\in SSt$ by Definition \ref{df:R0}.\ref{df:R0.2}.
If $\eta\prec\mI_{N}$, then $\gamma\in L(\mS)\cap\Psi$ with $\mS<\eta$.
Then $\gamma<\mS^{\dagger}<\eta\in LSt$.
Let $\eta\prec\mI_{N}[\tau]$. Then there is a $\rho\prec\mI_{N}[\tau]$ such that $\rho<\eta$ and $\gamma\prec^{R}\rho^{\dagger \vec{i}}$ for a $\vec{i}$.
We obtain $\gamma<\rho^{\dagger \vec{i}}<\eta$.
\\
\ref{prp:lessthan_layer}.\ref{prp:lessthan_layer.8}.
Let $\rho\in R(\eta)$ and $\gamma\prec\kappa\in N(\rho)$. We show $\gamma\in R(\eta)$.
We may assume that one of Definitions \ref{df:R0}.\ref{df:R00} and  Definition \ref{df:R0}.\ref{df:R0.2} is the case.
If $\eta\prec\mI_{N}$, then $\rho\in L(\mS)\cap\Psi$ with $\mS<\eta$.
We obtain $\gamma\in L(\mS)$.
Let $\eta\prec\mI_{N}[\tau]$. Then there is a $\rho_{0}\prec\mI_{N}[\tau]$ such that $\rho_{0}<\eta$ and $\rho\prec^{R}\rho^{\dagger\vec{i}}$ for a $\vec{i}$.
We obtain $\gamma\prec^{R}\rho^{\dagger\vec{i}}$.
\eprf

\blem\label{lem:R}
\benu
\item\label{lem:R0}
Let $\eta\prec\mI_{N}$ and $\Psi\ni\gamma<\eta$ with $\Omega<\gamma\not\in\Psi_{\mI_{N}}$.
Then 
$\gamma\in R(\eta)$ holds.

\item\label{lem:R1}
Let 
$\eta\in L(\mS)$ for an $\mS\in SSt$, $\eta>\gamma\in L(\mS)\cap\Psi$ and $\gamma\prec\tau\in (L(\mS)\cup\{\mS\})\cap SSt^{M}$.
Assume the following two:
\benu
\item
If $\eta\prec\mS$ and $\gamma\prec\mS$, then $\gamma\prec\eta$.
\item
If $\eta\not\prec\mS$, then $\delta<\gamma<\tau\leq\eta$, where $\delta\in L(\mS)\cap\Psi$ denotes
the ordinal such that $\eta\preceq\kappa\in N(\delta)$.
\eenu
Then 
$\gamma\in R(\eta)$ holds
\eenu
\elem
\bprf
Each is seen by induction on $\ell\gamma$.
\\
\ref{lem:R}.\ref{lem:R0}.
Let $\Psi\ni\gamma<\eta\prec\mI_{N}$ with $\Omega<\gamma\not\in\Psi_{\mI_{N}}$.
Let $\gamma\prec\tau\in St^{M}$. If $\tau\in N(\rho)$ for a $\rho\in\Psi$, then $\rho<\gamma$ by Proposition \ref{prp:lessthan_layer}.\ref{prp:lessthan_layer.25}.
IH yields $\rho\in R(\eta)$, and $\gamma\in R(\eta)$ by Proposition \ref{prp:lessthan_layer}.\ref{prp:lessthan_layer.8}.
Let $\gamma\prec\mS\in SSt$. Suppose $LSt\ni\eta<\mS$. Then $\eta\leq\mT$ for $\mS=\mT^{\dagger i}$ and an $i$.
Then we would have $\mT<\gamma$.
Hence $\mS<\eta$, and $\gamma\in R(\eta)$.
\\
\ref{lem:R}.\ref{lem:R1}.
First let $\eta\prec\mS$. If $\gamma\prec\mS$, then $\gamma\prec\eta$ by the first assumption,
and $\gamma\in R(\eta)$ by Definition \ref{df:R0}.\ref{df:R0.0}.
Let $\gamma\not\prec\mS$, and let $\rho$ denote the ordinal such that $\rho\prec\mS$ and
$\gamma\prec^{R}\sigma\in N(\rho)$.
We obtain $\rho<\gamma<\eta$ by Proposition \ref{prp:lessthan_layer}.\ref{prp:lessthan_layer.3}.
$\gamma\in R(\eta)$ follows from Definition \ref{df:R0}.\ref{df:R0.3}.

In what follows let $\eta\not\prec\mS$.
By the second assumption we have $\delta<\gamma<\tau\leq\eta$ for the ordinal $\delta$
such that $\eta\preceq\kappa\in N(\delta)$.
If $\gamma\prec\tau=\eta$, then $\gamma\in R(\eta)$. Let $\tau<\eta$, and $\gamma_{0}\in L(\mS)\cap\Psi$ denote the ordinal such that
$\tau\in N(\gamma_{0})$.

Let $\gamma_{0}<\delta$. 
By Proposition \ref{prp:lessthan_layer}.\ref{prp:lessthan_layer.0}
 we obtain $N(\gamma_{0})\ni\tau<\kappa<\kappa_{0}$ for the ordinal $\kappa_{0}$ such that
$\delta\prec\kappa_{0}\in(L(\mS)\cup\{\mS\})\cap SSt^{M}$.
We would have $\gamma<\tau<\delta$ by Proposition \ref{prp:lessthan_layer}.\ref{prp:lessthan_layer.4}.
We obtain $\delta\leq\gamma_{0}$.

If $\gamma_{0}=\delta$, then $\{\tau<\kappa\}\subset N(\delta)$ with $\gamma\prec\tau$ and $\eta\preceq\kappa$.
We obtain $\gamma\in R(\eta)$ by Definition \ref{df:R0}.\ref{df:R0.100}.

Let $\delta<\gamma_{0}\prec\tau_{0}\in(L(\mS)\cup\{\mS\})\cap SSt^{M}$.
If $\tau_{0}\leq\eta$, then IH yields $\gamma_{0}\in R(\eta)$.
We obtain $\gamma\in R(\eta)$ by Proposition \ref{prp:lessthan_layer}.\ref{prp:lessthan_layer.8}.
Suppose $\eta<\tau_{0}$. If $\kappa<\tau_{0}$, then we would have $\kappa<\gamma_{0}<\gamma<\eta\leq\kappa$ by
$\delta<\gamma_{0}$ and Propositions \ref{prp:lessthan_layer}.\ref{prp:lessthan_layer.4} and \ref{prp:lessthan_layer}.\ref{prp:lessthan_layer.0}.
We obtain $\gamma_{0}<\gamma<\eta<\tau_{0}\leq\kappa$, where $\gamma_{0}\prec\tau_{0}$, $\eta\prec\kappa$
and $\{\tau_{0},\kappa\}\subset(L(\mS)\cup\{\mS\})\cap SSt^{M}$.
If $\gamma\prec\eta$, then $\gamma\in R(\eta)$.
Let $\gamma\not\prec\eta$. We claim that $\tau_{0}=\kappa$.
Let $\gamma\preceq\psi_{\pi}^{f}(a)<\eta<\pi\preceq\tau_{0}$. If $\eta\prec\pi$, then $\tau_{0}=\kappa$.
Let $\eta\not\prec\pi$, and $\eta\preceq\psi_{\sigma}^{g}(b)<\pi<\sigma\preceq\kappa$.
Then we obtain $\tau_{0}=\kappa$ by Proposition \ref{prp:jumpover}.
We have shown $\tau_{0}=\kappa$.
Then $\gamma\in R(\eta)$ follows from Definition \ref{df:R0}.\ref{df:R0.3}.
\eprf

\bdf\label{df:gh.1}
{\rm
Let $\mS\in SSt$ and $\alp\in L(\mS)$.

Let $\vec{i}=(i_{0}\geq i_{1}\geq\cdots\geq i_{m})$ be a weakly descending chain of positive
integers with $i_{0}\leq N$.
Then let
$o(\vec{i}\,):=\ome^{i_{0}-1}+\ome^{i_{1}-1}+\cdots+\ome^{i_{m}-1}<\ome^{N}$.

Let us define ordinals $\mathtt{ g}^{\prime}_{1}(\alp)$ and
$\mathtt{ g}_{1}(\alp)$ as follows.
Let $\lam=\ome^{N+1}$.

\benu
\item
Let $\rho\prec\mS$. 
Then 
$\mathtt{ g}^{\prime}_{1}(\rho)=\lam^{\mathtt{g}_{0}(\rho)}$
and
$\mathtt{ g}_{1}(\rho)=\lam^{\mathtt{g}_{0}(\rho)+1}$.

\item
Let $\rho\in L(\mS)$ be such that
$\rho\prec\mT\in SSt^{M}\cap(L(\mS)\cup\{\mS\})$,
$\alp\prec\kap\in N(\rho)\cup
\{(\psi_{\mI_{N}[\rho]}(a))^{\dagger \vec{i}}:\vec{i}\neq\emptyset \}$,
where $\alp\neq\psi_{\mI_{N}[\rho]}(b)$ for any $b$ if $\kap=\mI_{N}[\rho]$.
Let $\mathtt{ g}_{1}(\alp)=\mathtt{ g}^{\prime}_{1}(\alp)+\lam^{\mathtt{g}_{0}(\rho)}$.

 \benu
  \item
$\mathtt{ g}^{\prime}_{1}(\mT^{\dagger \vec{i}}[\rho/\mT])=\mathtt{ g}^{\prime}_{1}(\rho)+\lam^{\mathtt{g}_{0}(\rho)}\cdot (o(\vec{i}\, )+1)$.

 \item
 $\alp\prec \mT^{\dagger \vec{i}}[\rho/\mT]$:
$\mathtt{ g}^{\prime}_{1}(\alp)=\mathtt{ g}^{\prime}_{1}(\rho)+\lam^{\mathtt{g}_{0}(\rho)}\cdot o(\vec{i}\,)$.

\item
$\mathtt{ g}^{\prime}_{1}(\mI_{N}[\rho])=\mathtt{ g}^{\prime}_{1}(\rho)+\lam^{\mathtt{g}_{0}(\rho)}\cdot (\ome^{N}+1)$.

\item
$\alp\prec\mI_{N}[\rho]$: 
$\mathtt{ g}^{\prime}_{1}(\alp)=\mathtt{ g}^{\prime}_{1}(\rho)+\lam^{\mathtt{g}_{0}(\rho)}\cdot \ome^{N}$.

\item
$\mathtt{ g}^{\prime}_{1}((\psi_{\mI_{N}[\rho]}(a))^{\dagger \vec{i}})=\mathtt{ g}^{\prime}_{1}(\rho)+\lam^{\mathtt{g}_{0}(\rho)}\cdot (\ome^{N}+o(\vec{i}\,)+1)$.

\item
$\alp\prec(\psi_{\mI_{N}[\rho]}(a))^{\dagger \vec{i}}$:
$\mathtt{ g}^{\prime}_{1}(\alp)=\mathtt{ g}^{\prime}_{1}(\rho)+\lam^{\mathtt{g}_{0}(\rho)}\cdot 
(\ome^{N}+o(\vec{i}\,))$.

\item
$\mathtt{ g}^{\prime}_{1}(\rho^{\dagger\vec{i}})=\mathtt{ g}^{\prime}_{1}(\rho)+\lam^{\mathtt{g}_{0}(\rho)}\cdot (\ome^{N}+\ome^{N}+o(\vec{i}\,)+1)$.

\item
$\alp\prec\rho^{\dagger\vec{i}}$:
$\mathtt{ g}^{\prime}_{1}(\alp)=\mathtt{ g}^{\prime}_{1}(\rho)+\lam^{\mathtt{g}_{0}(\rho)}\cdot (\ome^{N}+\ome^{N}+o(\vec{i}\,))$.

\eenu

\eenu

Let $\mathtt{ g}(\alp)=(\mathtt{ g}_{1}(\alp),\mathtt{ g}_{2}(\alp))$.
}
\edf

\blem\label{lem:LSwinding}
Let $\eta\in L(\mS)$
with $\mS\in SSt$.
Then
$\mathtt{g}^{*}_{0}(\gam)\leq\mathtt{g}^{*}_{0}(\eta)$,
$\mathtt{ g}(\gam)<_{lx}\mathtt{ g}(\eta)$ and
$SC(\mathtt{ g}_{2}(\gam))\subset \psi_{\mI_{N}}(b)$ 
for $\gam\in R(\eta)$ and $b=\mathtt{ g}^{*}_{0}(\eta)$.
\elem
\bprf
\\
\textbf{Case 1}. $\gam\prec\eta$: We have $\mathtt{g}^{*}_{0}(\gam)=\mathtt{g}^{*}_{0}(\eta)$.
If $\eta\in\Psi$, then 
$\mathtt{ g}_{1}(\eta)=\mathtt{ g}_{1}(\gam)$ and
$\mathtt{ g}_{2}(\gam)<\mathtt{ g}_{2}(\eta)$ by
Lemma \ref{lem:of}.
Otherwise $\mathtt{ g}_{1}(\gam)<\mathtt{ g}_{1}(\eta)$.
In what follows assume $\gam\not\prec\eta$.
We claim that $\mathtt{ g}_{1}(\gam)<\mathtt{ g}_{1}(\eta)$.
\\
\textbf{Case 2}. $\eta\preceq\tau_{1}$, $\gam\prec^{R}\tau_{2}$ with
$\{\tau_{2}\leq\tau_{1}\}\subset N(\tau)$
for a $\tau\in L(\mS)$,
cf.\,Definitions \ref{df:R0}.\ref{df:R0.0} and \ref{df:R0}.\ref{df:R0.100}: We have
$\mathtt{ g}_{1}(\eta)=\mathtt{ g}^{\prime}_{1}(\tau)+\lam^{\mathtt{ g}_{0}(\tau)}\cdot (\alp+1)$
for an $\alp<\ome^{N+1}$.
If $\gam\prec\tau_{2}$, then
$\mathtt{ g}_{1}(\gam)=\mathtt{ g}^{\prime}_{1}(\tau)+\lam^{\mathtt{ g}_{0}(\tau)}\cdot(\bet+1)$
with $\bet<\alp$.
Otherwise let $\sig\prec\tau_{2}$ be such that
$\gam\prec\sig_{1}\in SSt^{M}$,
$\sig_{1}\preceq^{R}\kap_{1}\in N(\sig)$.
We obtain $\mathtt{g}_{0}(\sig)<\mathtt{g}_{0}(\tau)$ by
Proposition \ref{prp:LSwinding.2}, and
$\mathtt{ g}_{1}(\gam)=\mathtt{ g}^{\prime}_{1}(\tau)+\lam^{\mathtt{ g}_{0}(\tau)}\cdot\bet
+\del$ with $\del<\lam^{\mathtt{g}_{0}(\sig)+1}\leq\lam^{\mathtt{g}_{0}(\tau)}$.
\\
\textbf{Case 3}. 
$\rho,\eta\prec\tau_{1}
\in N(\tau)$,
$\rho<\eta$ and $\gam\prec^{R}\kap\in N(\rho)$,
cf.\,Definitions \ref{df:R0}.\ref{df:R0.2} and \ref{df:R0}.\ref{df:R0.3}:
We have 
$\mathtt{g}^{\prime}_{1}(\rho)=\mathtt{ g}^{\prime}_{1}(\tau)+\lam^{\mathtt{ g}_{0}(\tau)}\cdot \alp$
for an $\alp<\ome^{N+1}$,
$\mathtt{ g}_{1}(\eta)=\mathtt{ g}^{\prime}_{1}(\rho)+\lam^{\mathtt{ g}_{0}(\tau)}$, and
$\mathtt{ g}_{1}(\gam)=\mathtt{g}^{\prime}_{1}(\rho)+\del$ for $\del<\lam^{\mathtt{ g}_{0}(\tau)}$
by Proposition \ref{prp:LSwinding.2}.
\\
\textbf{Case 4}. $\rho,\eta\prec\mS$, $\rho<\eta$ and $\gam\prec^{R}\kap\in N(\rho)$,
cf.\,Definition \ref{df:R0}.\ref{df:R0.3}:
We have 
$\mathtt{g}_{1}(\eta)=\lam^{\mathtt{g}_{0}(\eta)+1}$,
and
$\mathtt{g}^{\prime}_{1}(\rho)=\lam^{\mathtt{g}_{0}(\rho)}$,
where
$\mathtt{g}_{0}(\rho)\leq\mathtt{g}_{0}(\eta)$ by Proposition \ref{prp:LSwinding.1}.
On the other hand we have
$\mathtt{g}_{1}(\gam)=\mathtt{g}^{\prime}_{1}(\rho)+\del$ with $\del<\lam^{\mathtt{g}_{0}(\rho)+1}$.

Thus $\mathtt{g}(\gam)<_{lx}\mathtt{g}(\eta)$ is shown.
In each case  $c=\mathtt{ p}_{0}(\gam)\leq \mathtt{ g}^{*}_{0}(\gam)\leq\mathtt{ g}^{*}_{0}(\eta)=b$ holds by Proposition \ref{prp:g01}.
We obtain $\psi_{\mI}(c)\leq\psi_{\mI_{N}}(b)$ by
$c\in C_{c}(\psi_{\mI_{N}}(c))$ and $b\in C_{b}(\psi_{\mI_{N}}(b))$.
On the other hand we have
$SC(\mathtt{ g}_{2}(\gam))\subset \psi_{\mI_{N}}(c)$ by Proposition \ref{prp:g01}.
Hence $SC(\mathtt{ g}_{2}(\gam))\subset \psi_{\mI_{N}}(b)$.
\eprf

\bprp\label{prp:hg}
Let $\{\alp_{1},\bet\}\subset L(\mS)$ for an $\mS\in SSt$,
 $\alp_{1}=\psi_{\kap}^{f}(a)\leq\psi_{\sigma}^{h}(c)=\bet$ and 
$\bet\in C_{a}(\alp_{1})$. Then
$c<a$ and 
$\mathtt{g}^{*}_{0}(\bet)\leq\mathtt{g}^{*}_{0}(\alp_{1})$.
\eprp
\bprf 
By induction on $\ell\bet$.
We have $c\in K_{\alp_{1}}(\bet)<a$, and $\{\sig, c\}\subset C_{a}(\alp_{1})$. 
We show $\mathtt{ g}^{*}_{0}(\bet)\leq\mathtt{g}_{0}^{*}(\alp_{1})$.
First let $\bet\prec\mS$.
We show $\mathtt{ g}^{*}_{0}(\bet)=\mathtt{ p}_{\mS}(\bet)\leq\mathtt{g}^{*}_{0}(\alp_{1})$.
We can assume $\sig=\mS$ by IH.
Let $\gam=\psi_{\mS}^{g}(b)$ be a proper subterm of $\bet$.
If $\gam\in K_{\alp_{1}}(\bet)$, then $b<a$.
If $\gam<\alp_{1}$, then $\mathtt{ g}^{*}_{0}(\gam)\leq\mathtt{ g}^{*}_{0}(\alp_{1})$ by
Proposition \ref{prp:LSwinding.3}.

Second let $\rho\prec\mS$ and $\bet\prec^{R}\kap\in N(\rho)$.
Then  $\mathtt{ g}^{*}_{0}(\bet)=\mathtt{ p}_{\mS}(\rho)$.
If $\alp_{1}\leq\rho$, then $\rho\in C_{a}(\alp_{1})$ and
$\mathtt{ p}_{\mS}(\rho)\leq\mathtt{g}^{*}_{0}(\alp_{1})$ by the first case.
Let $\rho<\alp_{1}<\bet$.
Then we obtain $\mathtt{ g}^{*}_{0}(\alp)=\mathtt{ p}_{\mS}(\rho)=\mathtt{ g}^{*}_{0}(\bet)$.
\eprf

\subsection{Coefficients}
In this subsection we introduce coefficient sets $\mathcal{E}(\alpha),G_{\delta}(\alpha), F_{X}(\alpha),k_{X}(\alpha)$ of 
$\alpha\in OT(\mI_{N})$ for $X\subset OT(\mI_{N})$,
each of which is a finite set of subterms of $\alpha$.
These are utilized in our well-foundedness proof.
Roughly $\mathcal{E}(\alpha)$ is the set of subterms of the form $\psi_{\pi}^{f}(a)$, and
$F_{X}(\alpha)$ [$k_{X}(\alpha)$] the set of subterms in $X$ [subterms not in $X$], resp.

Let us write for $\alp<\mI_{N}$,
$\alp^{\dagger 0}=\min\{\sig\in Reg: \sig>\alp\}$ for the next regular ordinal $\alp^{+}$ above $\alp$.
Let $\alp^{\dagger i}:=\infty$ if $\alp\geq\mI_{N}$.
For $0\leq i\leq N$, let 
$\alp^{- i}:=\max\{\sig\in St_{i}\cup\{0\}: \sig\leq\alp\}$ when $\alp<\mI_{N}$, and
$\alp^{- i}:=\mI_{N}$ if $\alp\geq\mI_{N}$.

Although $\alp^{- 1}$ is similar to the Mostowski uncollapsing $\alp[\rho/\mS]^{-1}$ in
Definition \ref{df:divide}, no confusion likely occurs.

Since $St_{i+1} \subset St_{i}$, we obtain $\alp^{\dagger i}\leq\alp^{\dagger(i+1)}$ and
$\bet^{\dagger 0}<\sig$ if $\bet<\sig\in St\cap \mI_{N}$ since each $\sig\in St$ is a limit of regular ordinals.

Note that $R(\eta)\subset L(\mS)$ if $\eta\in L(\mS)$, and
$\gam^{- N}=\eta^{- N}$ for every $\gam,\eta\in L(\mS)$.

\begin{definition}\label{df:EGFk}

{\rm 
For terms $\alpha,\delta\in OT(\mI_{N})$ and $X\subset OT(\mI_{N})$, finite sets 
$\mathcal{E}(\alpha)$, $G_{\delta}(\alpha)$, $F_{X}(\alpha)$, $k_{X}(\alpha)$ of terms 
 are defined recursively as follows.
\benu

\item\label{df:EGFk.1}
$\mathcal{E}(\alpha)=\emptyset$ for $\alpha\in\{0,\Omega,\mI_{N}\}$.
$\mathcal{E}(\alpha_{m}+\cdots+\alpha_{0})=\bigcup_{i\leq m}\mathcal{E}(\alpha_{i})$.
$\mathcal{E}(\varphi\beta\gamma)=\mathcal{E}(\beta)\cup\mathcal{E}(\gamma)$.
$\mathcal{E}(\mI_{N}[\rho])=\mathcal{E}(\rho^{\dagger\vec{i}})=\mathcal{E}(\mS^{\dagger\vec{i}}[\rho/\mS])=\mathcal{E}(\rho)$.
$\mathcal{E}(\psi_{\pi}^{f}(a))=\{\psi_{\pi}^{f}(a)\}$.
$\mathcal{E}(\psi_{\mI_{N}}(a))=\{\psi_{\mI_{N}}(a)\}$.

\item\label{df:EGFk.2}
$\mathcal{A}(\alpha)=\bigcup\{\mathcal{A}(\beta): \beta\in\mathcal{E}(\alpha)\}$
 for $\mathcal{A}\in\{G_{\delta},F_{X},k_{X}\}$.

\item\label{df:EGFk.3}
$G_{\delta}(\psi_{\mI_{N}}(a))=G_{\delta}(a)$.
$F_{X}(\psi_{\mI_{N}}(a))=F_{X}(a)$ if $\psi_{\mI_{N}}(a)\not\in X$, and
$F_{X}(\psi_{\mI_{N}}(a))=\{\psi_{\mI_{N}}(a)\}$ if $\psi_{\mI_{N}}(a)\in X$.
$k_{X}(\psi_{\mI_{N}}(a))=\{\psi_{\mI_{N}}(a)\}\cup k_{X}(a)$ if $\psi_{\mI_{N}}(a)\not\in X$, and
$k_{X}(\psi_{\mI_{N}}(a))=\emptyset$ if $\psi_{\mI_{N}}(a)\in X$.

\[
G_{\delta}(\psi_{\pi}^{f}(a))=\left\{
\begin{array}{ll}
G_{\delta}(\{\pi,a\}\cup \mathrm{fld}(f)) & \delta<\pi
\\
\{\psi_{\pi}^{f}(a)\} & \pi\leq\delta
\end{array}
\right.
\]
\[
F_{X}(\psi_{\pi}^{f}(a))=\left\{
\begin{array}{ll}
F_{X}(\{\pi,a\}\cup \mathrm{fld}(f)) & \psi_{\pi}^{f}(a)\not\in X
\\
\{\psi_{\pi}^{f}(a)\} & \psi_{\pi}^{f}(a)\in X
\end{array}
\right.
\]
\[
k_{X}(\psi_{\pi}^{f}(a))=\left\{
\begin{array}{ll}
\{\psi_{\pi}^{f}(a)\}\cup k_{X}(\{\pi,a\}\cup \mathrm{fld}(f)) & \psi_{\pi}^{f}(a)\not\in X
\\
\emptyset & \psi_{\pi}^{f}(a)\in X
\end{array}
\right.
\]

\item\label{df:EGFk.4}
For 
$\alp\in N(\rho)$
\[
G_{\delta}(\alp)=\left\{
\begin{array}{ll}
\{\alp\} & \alp<\delta
\\
G_{\delta}(\rho) & \delta\leq\alp
\end{array}
\right.
\]
$F_{X}(\alp)=F_{X}(\rho)$ and
$k_{X}(\alp)=k_{X}(\rho)$.

\eenu

For $\mathcal{A}\in\{K_{\delta},G_{\delta},F_{X},k_{X}\}$ and sets $Y\subset OT(\mI_{N})$,
$\mathcal{A}(Y):=\bigcup\{\mathcal{A}(\alpha): \alpha\in Y\}$.
}
\end{definition}

\begin{definition}\label{df:Simm}
{\rm $S(\eta)$ denotes a set of immediate subterms of $\eta$.
For example $S(\varphi\beta\gamma)=\{\beta,\gamma\}$.
$S(\eta):=\emptyset$ when $\eta\in\{0,\Omega,\mathbb{I}_{N}\}$,
$S(\alp)=\{\rho\}$ for $\alp\in N(\rho)$,
$S(\eta)=\{\eta\}$ when $\eta\in\Psi$.
}
\end{definition}

\begin{proposition}\label{prp:G}
For $\{\alpha,\delta, a,b,\rho\}\subset OT(\mI_{N})$,
\benu
\item\label{prp:G1}
$G_{\delta}(\alpha)\leq\alpha$.

\item\label{prp:G2}
$\alpha\in C_{a}(b) 
\Rarw G_{\delta}(\alpha)\subset C_{a}(b)$.

\eenu
\end{proposition}
\bprf
These are shown simultaneously by induction on $\ell\alpha$.
It is easy to see that
\begin{equation}\label{eq:G}
G_{\delta}(\alpha)\ni\beta \Rarw \beta<\delta \spand \ell\beta\leq \ell\alpha
\end{equation}
\ref{prp:G}.\ref{prp:G1}.
Consider the case $\alpha=\psi_{\pi}^{f}(a)$ with $\delta<\pi$.
Then $G_{\delta}(\alpha)=G_{\delta}(\mathrm{fld}(f)\cup\{\pi,a\})$.
On the other hand we have  
$\mathrm{fld}(f)\cup\{\pi,a\}\subset C_{a}(\alpha)$.
Proposition \ref{prp:G}.\ref{prp:G2} with (\ref{eq:G}) 
yields 
$G_{\delta}(\mathrm{fld}(f)\cup\{\pi,a\})\subset C_{a}(\alpha)\cap\pi\subset\alpha$.
Hence $G_{\delta}(\alpha)<\alpha$.

Next let $\alp\in N(\rho)$
with $\delta\leq\alp$.
Then $G_{\delta}(\alpha)=G_{\delta}(\rho)$.
By IH we have $G_{\delta}(\rho)\leq\rho<\alp$.
Hence $G_{\delta}(\alpha)<\alpha$.
\\
\ref{prp:G}.\ref{prp:G2}.
Since $G_{\delta}(\alpha)\leq\alpha$ by Proposition \ref{prp:G}.\ref{prp:G1}, we can assume $\alpha\geq b$.

Consider the case $\alpha=\psi_{\pi}^{f}(a)$ with $\delta<\pi$.
Then $\mathrm{fld}(f)\cup\{\pi,a\}\subset C_{a}(b)$ and 
$G_{\delta}(\alpha)= G_{\delta}(\mathrm{fld}(f)\cup\{\pi,a\})$.
IH yields the lemma.

Next let $\alp\in N(\rho)$
with $\delta\leq\alp$.
Then $G_{\delta}(\alpha)=G_{\delta}(\rho)$ and $\rho<\alp$.
$b\leq\alpha\in C_{a}(b)$ yields 
$\rho
\in C_{a}(b)$. 
IH yields the lemma.
\hspace*{\fill} $\Box$

\begin{proposition}\label{prp:G3}
If
$\beta\not\in C_{a}(\alp)$ and $K_{\del}(\beta)<a$, then there exists a $\gamma\in F_{\del}(\beta)$
such that $C_{a}(\alp)\not\ni\gamma<\del$.
\end{proposition}
\bprf
By induction on $\ell\beta$. Assume $\beta\not\in C_{a}(\alpha)$ and 
$K_{\del}(\beta)<a$.
By IH we can assume that $\beta=\psi_{\kappa}^{f}(b)$.
If $\beta<\del$, then $\beta\in F_{\del}(\beta)$, and $\gamma=\beta$ is a desired one.
Assume $\beta\geq\del$. Then we obtain 
$K_{\del}(\beta)=\{b\}\cup K_{\del}(\{b,\kappa\}\cup \mathrm{fld}(f))<a$.
In particular $b<a$, and hence $\{b,\kappa\}\cup \mathrm{fld}(f)\not\subset C_{a}(\alp)$.
By IH there exists a $\gamma\in F_{\del}(\{b,\kappa\}\cup \mathrm{fld}(f))=F_{\del}(\beta)$
such that $C_{a}(\alp)\not\ni\gamma<\del$.
\hspace*{\fill} $\Box$

\section{Well-foundedness proof with the maximal distinguished sets}\label{sec:distinguished}
In this section working in the second order arithmetic 
$\Sigma^{1}_{N+2}\mbox{{\rm -DC}}+\mbox{{\rm BI}}$, 
we show the well-foundedness of the notation system $OT(\mI_{N})$ up to
\textit{each} $\alp<\Ome$.
The proof is based on distinguished classes, which was first
introduced by Buchholz\cite{Buchholz75}.
Each ordinal term $\alp\in OT(\mI_{N})$ is identified with its code $\lc\alp\rc\in\mathbb{N}$, cf.\,Lemma \ref{lem:compT}.

\subsection{Distinguished sets}

In this subsection
we establish elementary facts on
distinguished classes.

 $X,Y,Z,\ldots$ range over \textit{subsets} of $OT(\mI_{N})$, 
 while $\mathcal{X},\mathcal{Y},\ldots$ range over \textit{classes}, which are
 definable by second-order formulas in the language of arithmetic.
Following \cite{BuchholzBSL},
we define sets $\mathcal{C}^{\alpha}(X)\subset OT(\mI_{N})$ for $\alpha\in OT(\mI_{N})$ and
$X\subset OT(\mI_{N})$ as follows.

\begin{definition}\label{df:CX}
{\rm 
For $\alpha,\beta\in OT(\mI_{N})$ and $X\subset OT(\mI_{N})$, let us define
a set
$\mathcal{C}^{\alpha}(X)$ recursively as follows.
\benu
\item\label{df:CX.1}
$\{0,\Omega,\mathbb{I}_{N}\}\cup(X\cap\alpha)\subset \mathcal{C}^{\alpha}(X)$.

\item\label{df:CX.2}
Let $(\alp_{1}+\cdots+\alp_{n})\in OT(\mI_{N})$ with
$\{\alp_{1},\ldots,\alp_{n}\}\subset \mathcal{C}^{\alpha}(X)$.
Then
$(\alp_{1}+\cdots+\alp_{n})\in \mathcal{C}^{\alpha}(X)$.

\item\label{df:CX.3}
Let $\varphi\beta\gamma\in OT(\mI_{N})$ with
$\{\beta,\gamma\}\subset \mathcal{C}^{\alpha}(X)$.
Then $\varphi\beta\gamma\in \mathcal{C}^{\alpha}(X)$.

\item\label{df:CX.4}

Let $\psi_{\mI_{N}}(\beta)\in OT(\mI_{N})$ with
$\beta\in \mathcal{C}^{\alpha}(X)$.
Then
$\psi_{\mI_{N}}(\beta)\in \mathcal{C}^{\alpha}(X)$ if $\mI_{N}>\alpha$.

\item\label{df:CX.5}

Let $\psi_{\sigma}^{f}(\beta)\in OT(\mI_{N})$ with
$\{\sigma,\beta\}\cup \mathrm{fld}(f)\subset \mathcal{C}^{\alpha}(X)$.
Then
$\psi_{\sigma}^{f}(\beta)\in \mathcal{C}^{\alpha}(X)$ if $\sigma>\alpha$.

\item\label{df:CX.6}
Let 
$\bet\in N(\rho)$
with $\rho\in  \mathcal{C}^{\alpha}(X)$.
Then
$\bet\in  \mathcal{C}^{\alpha}(X)$ if $\bet\geq\alp$.

\eenu

}
\end{definition}

\begin{proposition}\label{lem:CX2}
Assume $\forall\gamma\geq\alp[\gam\in P \Rarw  \gamma\in \mathcal{C}^{\gamma}(P)]$ for a set $P\subset OT(\mI_{N})$.

\benu
\item\label{lem:CX2.3} 
$\alpha\leq\beta \Rarw  \mathcal{C}^{\beta}(P)\subset  \mathcal{C}^{\alpha}(P)$.

\item\label{lem:CX2.4} 
$\alpha\leq\beta<\alpha^{\dagger 0} \Rarw  \mathcal{C}^{\beta}(P)= \mathcal{C}^{\alpha}(P)$.
\eenu
\end{proposition}
\bprf
\ref{lem:CX2}.\ref{lem:CX2.3}.
We see by induction on $\ell\gamma\,(\gamma\in OT(\mI_{N}))$ that
\begin{equation}\label{eq:CX2.3}
\forall\beta\geq\alpha[\gamma\in  \mathcal{C}^{\beta}(P) \Rarw \gamma\in  \mathcal{C}^{\alpha}(P)\cup(P\cap\beta)]
\end{equation}
For example, if $\psi_{\pi}^{f}(\delta)\in  \mathcal{C}^{\beta}(P)$ with $\pi>\beta\geq\alpha$
and $\{\pi,\delta\}\cup \mathrm{fld}(f)\subset \mathcal{C}^{\alpha}(P)\cup(P\cap\beta)$, then 
$\pi\in \mathcal{C}^{\alpha}(P)$, and 
for any $\gamma\in\{\delta\}\cup \mathrm{fld}(f)$, either $\gamma\in \mathcal{C}^{\alpha}(P)$ or
$\gamma\in P\cap\beta$. If $\gamma<\alpha$, then 
$\gamma\in P\cap\alpha\subset \mathcal{C}^{\alpha}(P)$.
If $\alpha\leq\gamma\in P\cap\beta$, then $\gamma\in \mathcal{C}^{\gamma}(P)$ by the assumption, and
by IH we have $\gamma\in \mathcal{C}^{\alpha}(P)\cup(P\cap\gamma)$, i.e., $\gamma\in \mathcal{C}^{\alpha}(P)$.
Therefore $\{\pi,\delta\}\cup \mathrm{fld}(f)\subset \mathcal{C}^{\alpha}(P)$, and 
$\psi_{\pi}^{f}(\delta)\in \mathcal{C}^{\alpha}(P)$.

Using (\ref{eq:CX2.3}) we see from the assumption that
$
\forall\beta\geq\alpha[ \gamma\in \mathcal{C}^{\beta}(P) \Rarw \gamma\in \mathcal{C}^{\alpha}(P)]
$.
\\
\ref{lem:CX2}.\ref{lem:CX2.4}.
Assume $\alpha\leq\beta<\alpha^{\dagger 0}$. Then by Proposition \ref{lem:CX2}.\ref{lem:CX2.3} we have
$ \mathcal{C}^{\beta}(P)\subset \mathcal{C}^{\alpha}(P)$.
$\gam\in \mathcal{C}^{\alpha}(P) \Rarw \gam\in \mathcal{C}^{\beta}(P)$ is seen 
by induction on $\ell\gam$ using the facts $\bet^{- 0}=\alp^{- 0}$ and
$\beta^{\dagger 0}=\alpha^{\dagger 0}$.
\eprf

\begin{definition}\label{df:wftg}
{\rm
\benu
\item 
$\mathrm{Prg[}X,Y] :\Lrarw \forall\alpha\in X(X\cap\alpha\subset Y \to \alpha\in Y)$.

\item 
For a definable class $\mathcal{X}$, $\mathrm{TI}[\mathcal{X}]$ denotes the schema:\\
$\mathrm{TI}[\mathcal{X}] :\Lrarw \mathrm{Prg}[\mathcal{X},\mathcal{Y}]\to \mathcal{X}\subset\mathcal{Y} \mbox{ {\rm holds for} any definable classes } \mathcal{Y}$.
\item
For $X\subset OT(\mI_{N})$, $W(X)$ denotes the \textit{well-founded part} of $X$. 
\item 
$Wo[X] : \Lrarw X\subset W(X)$.
\eenu
}
\end{definition}
Note that for $\alpha\in OT(\mI_{N})$,
$W(X)\cap\alpha=W(X\cap\alpha)$.

\begin{definition} \label{df:3wfdtg32}
{\rm 
For $P,X\subset OT(\mI_{N})\cap \mI_{N}$ and
$\gam\in OT(\mI_{N})\cap\mI_{N}$, define $W_{i}^{\alp}(P)\, (0\leq i\leq N)$ and 
$D_{i}^{\gam}[P;X]\,(0\leq i\leq N)$
recursively on $i\leq N$ as follows.
\beqnarr
W_{0}^{\alp}(P) & := & W( \mathcal{C}^{\alpha}(P))
\label{eq:W0}
\\
D_{i}^{\gam}[P;X] & :\Lrarw &
Wo[P] \spand
P\cap\gam^{-(i+1)}=X\cap\gam^{-(i+1)} \spand 
\label{eq:Di}
\\
&&
\forall\alpha<\mI_{N}
\left(
\gam^{-(i+1)}\leq\alpha\leq P \to W_{i}^{\alp}(P)\cap\alpha^{\dagger i}= P\cap\alpha^{\dagger i}
\right)
\nonumber
\\
W_{i+1}^{\gam}(X) & := & \bigcup\{P\subset OT(\mI_{N})\cap\mI_{N} :D_{i}^{\gam}[P;X]\}\,
(i<N)
\label{eq:Wi+1}
\eeqnarr
where $\gam^{-(N+1)}:=0$.
Obviously $D_{N}[X;Y] \Lrarw D_{N}[X;Z]$ for every $Y$ and $Z$.
From $ \mathcal{W}_{N}^{\gam}(X)$ define
\beqnarr
D_{N}[X] & :\Lrarw  & D_{N}[X;X] 
\nonumber
\\
& \Lrarw &
Wo[X]  \spand
\forall\gam
\left(
\gam\leq X\to \mathcal{W}_{N}^{\gam}(X)\cap\gam^{\dagger N}= X\cap\gam^{\dagger N}
\right)
\label{eq:DN}
\\
\mathcal{W}_{N+1} & := & \bigcup\{X\subset OT(\mI_{N})\cap \mI_{N}  :D_{N}[X]\}
\label{eq:WN+1}
\eeqnarr
A set $P$ is said to be an $i$-\textit{distinguished set} for $\gam$ and $X$ if $D_{i}^{\gam}[P;X]$,
and a set $X$ is an $N$-\textit{distinguished set} if $D_{N}[X]$.

}
\end{definition}

Observe that in $S_{\mI_{N}}$,
$W^{\alp}_{0}(P)$ as well as $D^{\gam}_{0}[P;X]$ are $\Del_{1}$.
Assuming that $D^{\gam}_{i}[P;X]$ is $\Del_{i+1}$,
$\mathcal{W}_{i+1}^{\gam}(X)$ is $\Sig_{i+1}$, and 
$D^{\gam}_{i+1}[P;X]$ is $\Del_{i+2}$.
Hence $D_{N}[X]$ is $\Del_{N+1}$, and
$\mathcal{W}=\mathcal{W}_{N+1}$ is a $\Sig_{N+1}$-class.
In $S_{\mI_{N}}$,
each $\mathcal{W}_{i}^{\gam}(X)$ is a set, i.e.,
$\fal\gam\in OT(\mI_{N})\cap\mI_{N}\fal X\subset OT(\mI_{N})\exi Y[Y=\mathcal{W}_{i}^{\gam}(X)]$
for $0\leq i\leq N$,
and $\mathcal{W}_{N+1}$ is a proper class.

\begin{proposition}\label{lem:3.11.632}
Let $D_{0}^{\gam}[P;X]$ and $\gam^{- 1}\leq\alpha\in P$. 
Then $\forall\beta\geq\gam^{- 1}[\alpha\in \mathcal{C}^{\beta}(P)]$.
\end{proposition}
\bprf
Let $D_{0}^{\gam}[P;X]$ and $\gam^{- 1}\leq\alpha\in P$. 
We obtain $\alpha\in P\cap\alpha^{\dagger 0}=W( \mathcal{C}^{\alpha}(P))\cap\alpha^{\dagger 0}\subset  \mathcal{C}^{\alpha}(P)$ by (\ref{eq:W0}) and (\ref{eq:Di}).
Hence $\forall\del\geq\gam^{- 1}(\del\in P\Rarw\del\in \mathcal{C}^{\del}(P))$, and
$\alpha\in \mathcal{C}^{\beta}(P)$ for any $\gam^{- 1}\leq\beta\leq\alpha$ by Proposition \ref{lem:CX2}.\ref{lem:CX2.3}. 
Moreover for $\beta>\alpha$ we have $\alpha\in P\cap\beta\subset \mathcal{C}^{\beta}(P)$.
\eprf

\bprp\label{prp:3wf6.2-1}
If $P\cap\alp=Q\cap\alp$, then
$W_{i}^{\alp}(P)=W_{i}^{\alp}(Q)$.

\eprp
\bprf
For $i>0$, this follows from (\ref{eq:Di}) and $\alp^{- i}\leq\alp$.
For $i=0$, we obtain $ \mathcal{C}^{\alp}(P)= \mathcal{C}^{\alp}(Q)$ by 
$P\cap\alp=Q\cap\alp$.
Hence
$W_{0}^{\alp}(P)=W( \mathcal{C}^{\alp}(P))=W( \mathcal{C}^{\alp}(Q))=W_{0}^{\alp}(Q)$ 
by (\ref{eq:W0}).
\eprf

\blem\label{lem:3wf6.2}
$\alpha\leq P \spand \alpha\leq Q \Rarw P\cap\alpha^{\dagger i}=Q\cap\alpha^{\dagger i}$
if $D_{i}^{\gam}[P;X]$ and $D_{i}^{\gam}[Q;X]$.
\elem
\bprf
Suppose $\alpha\leq P$, $\alpha\leq Q$, $D_{i}^{\gam}[P;X]$ and $D_{i}^{\gam}[Q;X]$.
We have $P\cap\gam^{-(i+1)}=X\cap\gam^{-(i+1)}=Q\cap\gam^{-(i+1)}$.
We may assume that $\gam^{-(i+1)}\leq\alp$ since
$\alp^{\dagger i}\leq\gam^{-(i+1)}$ when $\alp<\gam^{-(i+1)}$.

By (\ref{eq:Di}) we obtain
$\mathcal{W}^{\alp}_{i}(P)\cap\alp^{\dagger i}=P\cap\alp^{\dagger i}$ and
$\mathcal{W}^{\alp}_{i}(Q)\cap\alp^{\dagger i}=Q\cap\alp^{\dagger i}$.
We obtain $Wo[P\cup Q]$ by $Wo[P]$ and $Wo[Q]$.
We show $\bet\in P\cap Q$ by induction on $\bet\in (P\cup Q)\cap\alp^{\dagger i}$.
Let $\bet\in(P\cup Q)\cap\alp^{\dagger i}$ and $P\cap\bet=Q\cap\bet$.
If $\bet<\gam^{-(i+1)}$, then $\bet\in P\cap Q$ by $P\cap\gam^{-(i+1)}=Q\cap\gam^{-(i+1)}$.
Let $\gam^{-(i+1)}\leq\bet$. 

If $\alp\leq\bet$, then $P\cap\alp=Q\cap\alp$, and
$\mathcal{W}^{\alp}_{i}(P)\cap\alp^{\dagger i}=\mathcal{W}^{\alp}_{i}(Q)\cap\alp^{\dagger i}$
by Proposition \ref{prp:3wf6.2-1}.
Hence
$\bet\in P\cap Q$.

Let $\gam^{-(i+1)}\leq\bet<\alp$.
We obtain $\mathcal{W}^{\bet}_{i}(P)\cap\bet^{\dagger i}=\mathcal{W}^{\bet}_{i}(Q)\cap\bet^{\dagger i}$ by $P\cap\bet=Q\cap\bet$ and Proposition \ref{prp:3wf6.2-1}.
By (\ref{eq:Di}),
$\bet\leq P$ and $\bet\leq Q$,
we obtain 
$P\cap\bet^{\dagger i}=\mathcal{W}^{\bet}_{i}(P)\cap\bet^{\dagger i}=\mathcal{W}^{\bet}_{i}(Q)\cap\bet^{\dagger i}=Q\cap\bet^{\dagger i}$.
Hence $\bet\in P\cap Q$.
\eprf

\blem\label{lem:3wf6} $(\Sigma^{1}_{N+1}\mbox{{\rm -CA}})$\\
For each $i\leq N$, $\fal\gam<\mI_{N}\fal X\exists Y(Y=W^{\gam}_{i}(X))$.
Let $\gam<\mI_{N}$.
\benu
\item\label{lem:3wf6.1}
For $i\leq N$, $W^{\gam}_{i}(X)$ is a well order:
$Wo[W^{\gam}_{i}(X)]$.

\item\label{lem:3wf6.3}
For $i<N$,
$W^{\gam}_{i+1}(X)$ is the maximal $i$-distinguished {\rm set}
for $\gam$ and $X$ if $X\cap\gam^{-(i+1)}$ is a well order:
$Wo[X\cap\gam^{-(i+1)}] \Rarw D^{\gam}_{i}[W^{\gam}_{i+1}(X);X]$.
In particular $W^{\gam}_{i+1}(X)\cap\gam^{-(i+1)}=X\cap\gam^{-(i+1)}$ holds.
\eenu

\elem
\bprf
\ref{lem:3wf6}.\ref{lem:3wf6.1}.
Clearly $W^{\gam}_{0}(X)=W( \mathcal{C}^{\gam}(X))$ is a well order.
We show $Wo[W^{\gam}_{i+1}(X)]$.
Let 
$\{\bet<\alp\}\subset W^{\gam}_{i+1}(X)$.
Pick a $P$ and a $Q$ such that
$D^{\gam}_{i}[P;X]$, $\alp\in P$, $D^{\gam}_{i}[Q;X]$ and $\bet\in Q$ by (\ref{eq:Wi+1}).
 Lemma \ref{lem:3wf6.2} yields 
$\bet\in Q\cap\bet^{\dagger i}\subset P$.
We obtain $Wo[W^{\gam}_{i+1}(X)\cap\alp]$ by $Wo[P]$.
\\
\ref{lem:3wf6}.\ref{lem:3wf6.3}.
Assuming that $X\cap\gam^{-(i+1)}$ is a well order, we see that
$X\cap\gam^{-(i+1)}$ is the minimal $i$-distinguished set for $\gam$ and $X$:
$D_{i}^{\gam}[X\cap\gam^{-(i+1)};X]$.
We obtain $W^{\gam}_{i+1}(X)\cap\gam^{-(i+1)}=X\cap\gam^{-(i+1)}$.
Lemma \ref{lem:3wf6}.\ref{lem:3wf6.1} yields $Wo[W^{\gam}_{i+1}(X)]$.

Let $\gam^{-(i+1)}\leq\alp\leq W^{\gam}_{i+1}(X)$.
We show 
$W^{\alp}_{i}(W^{\gam}_{i+1}(X))\cap\alp^{\dagger i}=W^{\gam}_{i+1}(X)\cap\alp^{\dagger i}$.
Pick a $P$ such that $D^{\gam}_{i}[P;X]$ and $\alp\leq P$.
We obtain 
$W^{\alp}_{i}(P)\cap\alp^{\dagger i}=P\cap\alp^{\dagger i}
\subset W^{\gam}_{i+1}(X)\cap\alp^{\dagger i}$
by  (\ref{eq:Di}).
Let $D^{\gam}_{i}[Q;X]$ and $\bet\in Q\cap\alp^{\dagger i}$. 
Lemma \ref{lem:3wf6.2} yields 
$\bet\in Q\cap\bet^{\dagger i}=P\cap\bet^{\dagger i}$ for $\bet^{\dagger i}\leq\alp^{\dagger i}$.

Therefore we obtain
$W^{\alp}_{i}(P)\cap\alp^{\dagger i}=P\cap\alp^{\dagger i}
=W^{\gam}_{i+1}(X)\cap\alp^{\dagger i}$, a fortiori
$P\cap\alp=\mathcal{W}^{\gam}_{i+1}(X)\cap\alp$.
Hence
$W^{\gam}_{i+1}(X)\cap\alp^{\dagger i}=P\cap\alp^{\dagger i}=W^{\alp}_{i}(P)\cap\alp^{\dagger i}=W^{\alp}_{i}(W^{\gam}_{i+1}(X))\cap\alp^{\dagger i}$
by Proposition \ref{prp:3wf6.2-1}.
\eprf

\blem\label{lem:8wf2}
\benu
\item\label{lem:8wf2.1}
Let $X$ and $Y$ be $N$-distinguished sets, and $\gam<\mI_{N}$.
Then $\gam\leq X\spand \gam\leq Y
\Rarw 
X\cap\gam^{\dagger N}=Y\cap\gam^{\dagger N}$.

\item\label{lem:8wf2.4}
$\calw_{N+1}$ is the $N$-maximal distinguished {\rm class}, i.e.,
$D_{N}[\calw_{N+1}]$.

\item\label{lem:8wf2.5}
For a family $\{Y_{j}\}_{j\in J}$ of $N$-distinguished sets,
the union $Y=\bigcup_{j\in J} Y_{j}$ is also an $N$-distinguished set.

\eenu
\elem
\bprf 
\ref{lem:8wf2}.\ref{lem:8wf2.1} is seen as in Lemma \ref{lem:3wf6.2}.
\ref{lem:8wf2}.\ref{lem:8wf2.4} and \ref{lem:8wf2}.\ref{lem:8wf2.5}
follow from Lemma \ref{lem:8wf2}.\ref{lem:8wf2.1} as in Lemma \ref{lem:3wf6}.
\eprf

\blem\label{lem:6.15}
Let $D_{N}[X]$ and $\gam\in X\subset \mI_{N}$. Then for each 
$0<i\leq N$,
$\gam\in W_{i}^{\gam}(X)\cap\gam^{\dagger i}=X\cap\gam^{\dagger i}$ and
$D_{i-1}^{\gam}[X\cap\gamma^{\dagger i};X]$ 
holds.
In particular 
 $\gam\in \mathcal{C}^{\gam}(X)$ and $\mathcal{C}^{\gamma}(X)\cap\gamma\subset X$.
\elem
\bprf
By induction on $N-i$.
We obtain $\gam\in W_{N}^{\gam}(X)\cap\gam^{\dagger N}=X\cap\gam^{\dagger N}$
by $D_{N}[X]$ and $\gam\in X$.
Lemma \ref{lem:3wf6} with $Wo[X]$ yields 
$D_{N-1}^{\gam}[W_{N}^{\gam}(X);X]$, and 
$D_{N-1}^{\gam}[X\cap\gam^{\dagger N};X]$ 
follows.

Assuming 
$D_{i+1}^{\gam}[X\cap\gamma^{\dagger (i+2)};X]$, 
we obtain
$W_{i+1}^{\gam}(X)\cap\gam^{\dagger(i+1)}=X\cap\gam^{\dagger(i+1)}$ by 
$\gam^{-(i+1)}\leq\gam\in X$, and
$D_{i}^{\gam}[W_{i+1}^{\gam}(X);X]$ by Lemma \ref{lem:3wf6}.
Hence 
$D_{i}^{\gam}[X\cap\gamma^{\dagger (i+1)};X]$ 
and $\gam\in W_{i}^{\gam}(X)\cap\gam^{\dagger i}=X\cap\gam^{\dagger i}$.

$\gam\in W_{0}^{\gam}(X)\cap\gamma^{\dagger 0}=W(\mathcal{C}^{\gamma}(X))\cap\gamma^{\dagger 0}=X\cap\gamma^{\dagger 0}$ 
yields  $\gam\in \mathcal{C}^{\gam}(X)$ and $\mathcal{C}^{\gamma}(X)\cap\gamma\subset X$.
\eprf

\bprp\label{lem:6.16}
Let $D_{N}[X]$, $\alp\leq \gam\in X$ and $\alp\in  \mathcal{C}^{\gam}(X)$.
Then $\alp\in X$.
\eprp
\bprf
Lemma \ref{lem:6.15} yields
$\gam\in W( \mathcal{C}^{\gam}(X))\cap\gam^{\dagger 0}=W_{0}^{\gam}(X)\cap\gam^{\dagger 0}=X\cap\gam^{\dagger 0}$.
$\gam\geq\alp\in  \mathcal{C}^{\gam}(X)$ yields 
$\alp\in W_{0}^{\gam}(X)\cap\gam^{\dagger 0}=X\cap\gam^{\dagger 0}$.
\eprf

\bprp\label{prp:noncritical}
Let $D_{N}[X]$.
\benu
\item\label{prp:noncritical.1}
Let $\{\alp,\bet\}\subset X$ with $\alp+\bet=\alp\#\bet$ and $\alp>0$.
Then $\gam=\alp+\bet\in X$.

\item\label{prp:noncritical.2}
If $\{\alp,\bet\}\subset X$, then $\vphi\alp\bet\in X$.
\eenu
\eprp
\bprf
Proposition \ref{prp:noncritical}.\ref{prp:noncritical.2} is seen 
by main induction on $\alp\in X$ with
subsidiary induction on $\bet\in X$ using 
Proposition \ref{prp:noncritical}.\ref{prp:noncritical.1}.
We show Proposition \ref{prp:noncritical}.\ref{prp:noncritical.1}.
By Lemma \ref{lem:6.15}
we obtain
$\alp\in X\cap\alp^{\dagger 0}=W_{0}^{\alp}(X)\cap\alp^{\dagger 0}$.
We see that $\alp+\bet\in W_{0}^{\alp}(X)=W( \mathcal{C}^{\alp}(X))$ by induction on 
$\bet\in X\cap(\alp+1)\subset \mathcal{C}^{\alp}(X)$.
\eprf

\blem\label{lem:6.34}
\benu
\item\label{lem:6.34.2}
$ \mathcal{C}^{\mI_{N}}(\mathcal{W}_{N+1})\cap\mI_{N}=\mathcal{W}_{N+1}\cap\mI_{N}=
W( \mathcal{C}^{\mI_{N}}(\mathcal{W}_{N+1}))\cap\mI_{N}$.

\item\label{lem:6.34.3}
{\rm (BI)}
For {\rm each} $n<\ome$, 
$\mathrm{TI}[ \mathcal{C}^{\mI_{N}}(\mathcal{W}_{N+1})
\cap\ome_{n}(\mI_{N}+1)
]$, i.e.,
for each class $\mathcal{X}$,
$\mathrm{Prg}[ \mathcal{C}^{\mI_{N}}(\mathcal{W}_{N+1}),\mathcal{X}] \to 
 \mathcal{C}^{\mI_{N}}(\mathcal{W}_{N+1})
\cap\ome_{n}(\mI_{N}+1)
\subset\mathcal{X}$.

\item\label{lem:6.34.4}
For {\rm each} $n<\ome$, 
$ \mathcal{C}^{\mI_{N}}(\mathcal{W}_{N+1})
\cap\ome_{n}(\mI_{N}+1)
\subset 
W( \mathcal{C}^{\mI_{N}}(\mathcal{W}_{N+1}))$.
In particular 
$\{\mI_{N},\ome_{n}(\mI_{N}+1)\}\subset W( \mathcal{C}^{\mI_{N}}(\mathcal{W}_{N+1}))$.

\eenu
\elem
\bprf
\ref{lem:6.34}.\ref{lem:6.34.2}.
$\alp\in \mathcal{C}^{\mI_{N}}(\mathcal{W}_{N+1})\cap\mI_{N} \Rarw \alp\in\mathcal{W}_{N+1}$
is seen by induction on $\ell\alp$ using Proposition \ref{prp:noncritical} and
Lemma \ref{lem:8wf2}.\ref{lem:8wf2.4}.
Since $\mathcal{W}_{N+1}$ is well-founded, we obtain
$ \mathcal{C}^{\mI_{N}}(\mathcal{W}_{N+1})\cap\mI_{N} =W( \mathcal{C}^{\mI_{N}}(\mathcal{W}_{N+1}))\cap\mI_{N}$.
\\
\ref{lem:6.34}.\ref{lem:6.34.3}.
We show $\mathrm{TI}[ \mathcal{C}^{\mI_{N}}(\mathcal{W}_{N+1})\cap\ome_{n}(\mI_{N}+1)]$
by metainduction on $n<\ome$.
Let $D_{N}[Y]$. We obtain $Wo[Y]$, and $\mathrm{TI}[Y]$ follows from (BI).
We have $ \mathcal{C}^{\mI_{N}}(\mathcal{W}_{N+1})\cap\mI_{N}=\mathcal{W}_{N+1}\cap\mI_{N}$,
and $\mathcal{W}_{N+1}\cap\gam^{\dagger N}=Y\cap\gam^{\dagger N}$
for $\gam\in Y\cap\mI_{N}$ by Lemma \ref{lem:8wf2}.\ref{lem:8wf2.1}.
We obtain $\mathrm{TI}[\mathcal{W}_{N+1}\cap\mI_{N}]$, from which
$\mathrm{TI}[ \mathcal{C}^{\mI_{N}}(\mathcal{W}_{N+1})\cap(\mI_{N}+1)]$ follows.

Assuming $\mathrm{TI}[ \mathcal{C}^{\mI_{N}}(\mathcal{W}_{N+1})\cap\ome_{n}(\mI_{N}+1)]$,
$\mathrm{TI}[ \mathcal{C}^{\mI_{N}}(\mathcal{W}_{N+1})\cap\ome_{n+1}(\mI_{N}+1)]$ is seen from
the fact that $\mathrm{Prg}[ \mathcal{C}^{\mI_{N}}(\mathcal{W}_{N+1}),A] \to 
\mathrm{Prg}[ \mathcal{C}^{\mI_{N}}(\mathcal{W}_{N+1}),\mathtt{ j}[A]]$, where
 for a given formula $A$, $\mathtt{ j}[A](\alp)$ denotes the formula
 \\
{\small
$
\fal\bet\in  \mathcal{C}^{\mI_{N}}(\mathcal{W}_{N+1})
\left[
\fal\gam\in  \mathcal{C}^{\mI_{N}}(\mathcal{W}_{N+1})\cap\bet\, A(\gam)
\to
\fal\gam\in  \mathcal{C}^{\mI_{N}}(\mathcal{W}_{N+1})\cap(\bet+\ome^{\alp}) A(\gam)
\right]
$.
}
\eprf

\subsection{Sets $\mathcal{G}^{X}$}\label{subsec:C(X)}

In this subsection we establish a key fact, Lemma \ref{th:3wf16} on distinguished sets.

\begin{definition}\label{df:calg}
$\mathcal{G}^{X}:=\{\alpha\in OT(\mI_{N})\cap\mI_{N} :\alpha\in  \mathcal{C}^{\alpha}(X)
\spand  \mathcal{C}^{\alpha}(X)\cap\alpha\subset X\}$.
\end{definition}

\bprp\label{prp:calg}
Let $D_{N}[X]$ and $\alp\in X$. Then $\alp\in\mathcal{G}^{X}$.
\eprp
\bprf
By Lemma \ref{lem:6.15} we obtain 
$\alp\in W_{0}^{\alp}(X)=W( \mathcal{C}^{\alp}(X))$.
Hence $\alpha\in  \mathcal{C}^{\alpha}(X)$.
On the other side Proposition \ref{lem:6.16} yields
$ \mathcal{C}^{\alpha}(X)\cap\alpha\subset X$.
\eprf

\blem\label{lem:6.21}$(\Sigma^{1}_{N+1}\mbox{{\rm -CA}})$\\
Suppose $D_{N}[Y]$ and $\alp\in\mathcal{G}^{Y}$.
Let 
$P_{N}=W_{N}^{\alp}(Y)\cap\alp^{\dagger N}$.
Assume that the following condition (\ref{eq:6.21.550})
is fulfilled.
Then $\alp\in P_{N}$ and $D_{N}[P_{N}]$.
In particular $\alp\in\mathcal{W}_{N+1}$ holds.

Moreover if there exists a set $Z$ and an ordinal $\gam$ such that
$Y=W_{N}^{\gam}(Z)$ and $\alp^{- N}=\gam^{- N}$,
then $\alp\in Y$ holds.

\beqn\label{eq:6.21.550}
\fal\bet\geq\alp^{- 1}
\left(Y\cap\alp^{\dagger 1}<\bet \spand \bet^{\dagger 0}<\alp^{\dagger 0} \to
W_{0}^{\bet}(Y)\cap\bet^{\dagger 0}\subset Y
\right)
\eeqn

\elem
\bprf
If $Y=W_{N}^{\gam}(Z)$ with $\alp^{- N}=\gam^{- N}$, then
$Y\cap\alp^{- N}=Z\cap\alp^{- N}$ and $W_{N}^{\gam}(Z)=W_{N}^{\alp}(Y)$.
Hence if $\alp\in W_{N}^{\alp}(Y)$, then $\alp\in Y$.

Lemma \ref{lem:3wf6}.\ref{lem:3wf6.3} yields 
\beqn\label{eq:6.21.55i}
\fal i<N
\left[
W_{i+1}^{\bet}(Y)\cap\bet^{\dagger(i+1)}= Y\cap\bet^{\dagger(i+1)}
\right]
\eeqn

Let $P_{i}=W_{i}^{\alp}(Y)\cap\alp^{\dagger i}$ for $0\leq i\leq N$.
By $ \mathcal{C}^{\alp}(Y)\cap\alp\subset Y$ and $Wo[Y]$
we obtain for $P_{0}=W( \mathcal{C}^{\alp}(Y))\cap\alp^{\dagger 0}$
\beqn\label{eq:6.21.57}
P_{0}\cap\alp=Y\cap\alp= \mathcal{C}^{\alp}(Y)\cap\alp
\eeqn
Hence $\alp\in P_{0}$.
We obtain $D_{i-1}^{\alp}[W_{i}^{\alp}(Y);Y]$ for $i>0$.
This together with (\ref{eq:6.21.57}) yields for $0\leq i\leq N$
\beqn\label{eq:6.21.57+i}
P_{i}\cap\alp^{-i}=Y\cap\alp^{-i}
\eeqn

\bclm\label{clm:6.22-1}
$\alp^{\dagger 0}=\gam^{\dagger 0} \spand \gam\in P_{0} \Rarw \gam\in \mathcal{C}^{\gam}(P_{0})$.
\eclm
\textbf{Proof} of Claim \ref{clm:6.22-1}.
Let $\alp^{\dagger 0}=\gam^{\dagger 0}$ and
$\gam\in P_{0}=W( \mathcal{C}^{\alp}(Y))\cap\alp^{\dagger 0}$. 
We obtain $\gam\in \mathcal{C}^{\alp}(Y)= \mathcal{C}^{\gam}(Y)$ by 
Propositions \ref{prp:calg} and \ref{lem:CX2}.
Hence
 $Y\cap\gam\subset \mathcal{C}^{\gam}(Y)\cap\gam= \mathcal{C}^{\alp}(Y)\cap\gam$.
$\gam\in W( \mathcal{C}^{\alp}(Y))$ yields
$Y\cap\gam\subset P_{0}$.
Therefore
we obtain
$\gam\in \mathcal{C}^{\gam}(Y)\subset \mathcal{C}^{\gam}(P_{0})$.
\hspace*{\fill} $\Box$ of Claim \ref{clm:6.22-1}.

\bclm\label{clm:6.22}
$D_{i}^{\alp}[P_{i};Y]$ and $\alp\in P_{i+1}$ for each $0\leq i< N$.
\eclm
\textbf{Proof} of Claim \ref{clm:6.22}.
Obviously $Wo[P_{i}]$.
(\ref{eq:6.21.57+i}) yields $P_{i}\cap\alp^{- (i+1)}=Y\cap\alp^{- (i+1)}$.
Let $\alp^{- (i+1)}\leq\bet\leq P_{i}$. 
We show $W_{i}^{\bet}(P_{i})\cap\bet^{\dagger i}=P_{i}\cap\bet^{\dagger i}$.
\\
\textbf{Case 1}. $\bet^{\dagger i}=\alp^{\dagger i}$:
First let $i=0$.
We obtain $ \mathcal{C}^{\bet}(P_{0})= \mathcal{C}^{\alp}(P_{0})$ by
Proposition \ref{lem:CX2} and Claim \ref{clm:6.22-1}.
Hence the assertion follows from (\ref{eq:6.21.57}).

Next let $i>0$. 
(\ref{eq:6.21.57+i}) with $\bet^{-i}=\alp^{-i}$ yields 
$W_{i}^{\bet}(P_{i})=W_{i}^{\alp}(P_{i})=W_{i}^{\alp}(Y)$.
\\
\textbf{Case 2}. $\bet^{\dagger i}<\alp^{\dagger i}$: 
For $i>0$, (\ref{eq:6.21.55i}) 
yields
$W_{i}^{\bet}(Y)\cap\bet^{\dagger i}= Y\cap\bet^{\dagger i}$.
We obtain
$W_{i}^{\bet}(P_{i})\cap\bet^{\dagger i}=W_{i}^{\bet}(Y)\cap\bet^{\dagger i}= Y\cap\bet^{\dagger i}=P_{i}\cap\bet^{\dagger i}$
by (\ref{eq:6.21.57+i}).

Let $i=0$. We have $\bet^{\dagger 0}\leq\alp^{-0}$.
First let $Y\cap\alp^{\dagger 1}<\bet$. 
Then the assumption (\ref{eq:6.21.550}) with $\alp^{- 1}\leq\bet$
yields
$W_{0}^{\bet}(Y)\cap\bet^{\dagger 0}\subset Y$.
We obtain
$W_{0}^{\bet}(P_{0})\cap\bet^{\dagger 0}=W_{0}^{\bet}(Y)\cap\bet^{\dagger 0}\subset Y\cap\bet^{\dagger 0}=P_{0}\cap\bet^{\dagger 0}$
by (\ref{eq:6.21.57}).
It remains to show $Y\cap\bet^{\dagger 0}\subset W_{0}^{\bet}(Y)$.
Let $\gam\in Y\cap\bet^{\dagger 0}$.
We obtain $\gam\in W_{0}^{\gam}(Y)$ by Lemma \ref{lem:6.15}.
On the other hand we have $\mathcal{C}^{\bet}(Y)\subset\mathcal{C}^{\gam}(Y)$
by Propositions \ref{prp:calg} and \ref{lem:CX2}.
Moreover (\ref{eq:6.21.57}) with Propositions \ref{prp:calg} and \ref{lem:CX2} yields 
$\gam\in\mathcal{C}^{\alp}(Y)\subset\mathcal{C}^{\bet}(Y)$.
Hence $\gam\in W_{0}^{\bet}(Y)$.

Next let $\bet\leq Y\cap\alp^{\dagger 1}$.
We obtain $Y\cap\bet^{\dagger 1}=\mathcal{W}_{1}^{\bet}(Y)\cap\bet^{\dagger 1}$, and
$\bet^{- 1}=\alp^{- 1}\leq\bet<\alp^{\dagger 1}=\bet^{\dagger 1}$ with 
$\bet<\bet^{\dagger 0}\leq\alp<\bet^{\dagger 1}$.
On the other hand we have $D_{0}^{\bet}[\mathcal{W}_{1}^{\bet}(Y);Y]$ by
Lemma \ref{lem:3wf6}.
Therefore
$P_{0}\cap\bet^{\dagger 0}=Y\cap\bet^{\dagger 0}=\mathcal{W}_{1}^{\bet}(Y)\cap\bet^{\dagger 0}=
W_{0}^{\bet}(\mathcal{W}_{1}^{\bet}(Y))\cap\bet^{\dagger 0}=W_{0}^{\bet}(P_{0})\cap\bet^{\dagger 0}$ by (\ref{eq:6.21.57}).

Thus $D_{i}^{\alp}[P_{i};Y]$ is shown.
From $\alp\in P_{0}$ we see by induction on $i<N$ that
$\alp\in P_{i}\cap\alp^{\dagger (i+1)}\subset 
W_{i+1}^{\alp}(Y)\cap\alp^{\dagger (i+1)}=P_{i+1}$
for the maximal $i$-distinguished set $W_{i+1}^{\alp}(Y)$ for $\alp$ and $Y$.
\hspace*{\fill} $\Box$ of Claim \ref{clm:6.22}.

\bclm\label{clm:6.23}
$D_{N}[P_{N}]$.
\eclm
\textbf{Proof} of Claim \ref{clm:6.23}.
Let $\bet\leq P_{N}=W_{N}^{\alp}(Y)\cap\alp^{\dagger N}$. 
Then $\bet<\alp^{\dagger N}$, and $\bet^{- N}\leq\alp^{- N}<\alp^{\dagger N}$.
We show 
$W_{N}^{\bet}(W_{N}^{\alp}(Y))\cap\bet^{\dagger N}=W_{N}^{\alp}(Y)\cap\bet^{\dagger N}$.
\\
\textbf{Case 1}. $\alp^{- N}\leq\bet$: 
By $W_{N}^{\alp}(Y)\cap\alp^{- N}=Y\cap\alp^{- N}$ with $Wo[Y]$, 
and $\alp^{- N}=\bet^{- N}$ we obtain
$W_{N}^{\alp}(Y)=W_{N}^{\alp}(W_{N}^{\alp}(Y))
=W_{N}^{\bet}(W_{N}^{\alp}(Y))$.
\\
\textbf{Case 2}. $\bet<\alp^{- N}$ and $\bet^{- N}\leq Y$:
We obtain $\bet^{\dagger N}\leq\alp^{- N}$.
Hence
$W_{N}^{\alp}(Y)\cap\bet^{\dagger N}=Y\cap\bet^{\dagger N}=
W_{N}^{\bet}(Y)\cap\bet^{\dagger N}$
by $D_{N}[Y]$.
Therefore
$W_{N}^{\bet}(Y)=W_{N}^{\bet}(W_{N}^{\alp}(Y))$.
We obtain 
$W_{N}^{\alp}(Y)\cap\bet^{\dagger N}=
W_{N}^{\bet}(W_{N}^{\alp}(Y))\cap\bet^{\dagger N}$.
\\
\textbf{Case 3}. $\bet<\alp^{- N}$ and $Y<\bet^{- N}$:
Then $\bet^{\dagger N}\leq\alp^{- N}$.
(\ref{eq:6.21.55i}) yields 
$Y\cap\bet^{\dagger N}=W_{N}^{\bet}(Y)\cap\bet^{\dagger N}$.
On the other hand we have
$Y\cap\bet^{\dagger N}=W_{N}^{\alp}(Y)\cap\bet^{\dagger N}$ and 
$W_{N}^{\bet}(Y)\cap\bet^{\dagger N}=
W_{N}^{\bet}(W_{N}^{\alp}(Y))\cap\bet^{\dagger N}$.
Therefore
$W_{N}^{\bet}(W_{N}^{\alp}(Y))\cap\bet^{\dagger N}=
W_{N}^{\alp}(Y)\cap\bet^{\dagger N}$.

\hspace*{\fill} $\Box$ of Claim \ref{clm:6.23}.

This completes a proof of Lemma \ref{lem:6.21}.
\eprf

\blem\label{cor:6.21} 
Assume $D_{N}[Y]$, $\mI_{N}>\mS\in Y\cap(St_{k}\cup\{0\})$ and 
$\{0,\Ome\}\subset Y$ for $0<k\leq N$.
Then $\mS^{\dagger k}\in\mathcal{W}_{N+1}$.
\elem
\bprf
Let us verify the condition (\ref{eq:6.21.550}) in Lemma \ref{lem:6.21}
for $\alp=\mS^{\dagger k}$.
Let $\alp^{- 1}\leq\bet$.
We have $\alp=\alp^{- 1}\leq\bet$. Hence $\alp^{\dagger 0}\leq\bet^{\dagger 0}$,
and (\ref{eq:6.21.550}) is vacuously fulfilled.

Thus
it suffices to show that $\alp=\mS^{\dagger k}\in\mathcal{G}^{Y}$.
$\alp\in \mathcal{C}^{\alp}(Y)$ follows from
$\mS\in Y\cap\alp$, cf.\,Definition \ref{df:CX}.\ref{df:CX.6}.
We show
$\gam\in \mathcal{C}^{\alp}(Y)\cap\alp \Rarw \gam\in Y$ by induction on 
$\ell\gam$.
By Proposition \ref{prp:noncritical} and the assumption $\{0,\Ome\}\subset Y$, 
we can assume $\mS\neq\gam=\psi_{\sig}^{f}(a)<\alp=\mS^{\dagger k}<\sig$,
cf.\,Definition \ref{df:CX}.\ref{df:CX.6}.
Suppose $\mS<\gam$. Then $\mS\in C_{a}(\gam)$, and 
$\alp=\mS^{\dagger k}\in C_{a}(\gam)\cap\sig\subset\gam$.
We obtain $\gam<\mS$. 
Lemma \ref{lem:6.15} with $\mS\in Y$ and $D_{N}[Y]$ yields
$\mS\in W_{0}^{\mS}(Y)\cap\mS^{\dagger 0}=Y\cap\mS^{\dagger 0}$
for $W_{0}^{\mS}(Y)=W(\mathcal{C}^{\mS}(Y))$, where
$\forall\del[\del\in Y \Rarw  \del\in \mathcal{C}^{\del}(Y)]$.
We obtain $\gam\in \mathcal{C}^{\mS}(Y)$ by $\gam\in \mathcal{C}^{\alp}(Y)$, 
$\mS<\alp$ and Proposition \ref{lem:CX2}.\ref{lem:CX2.3} 
Hence $\gam\in W_{0}^{\mS}(Y)\cap\mS^{\dagger 0}\subset Y$ follows.
Therefore $\alp\in\mathcal{G}^{Y}$.
\eprf

\bprp\label{lem:6.29}
$\{0,\Ome\}
\subset\mathcal{W}_{N+1}$.
\eprp
\bprf
For each $\alp\in\{0,\Ome\}$ and any set $Y\subset OT(\mI_{N})$ we have
$\alp\in \mathcal{C}^{\alp}(Y)$.
First let $\alp=0$.
We obtain $ \mathcal{C}^{0}(\emptyset)\cap\alp\subset \emptyset$, and 
$0\in\mathcal{G}^{\emptyset}$.
Moreover $D_{N}[\emptyset]$, and 
there is no $\bet$ such that 
$\bet^{\dagger 0}<\alp^{\dagger 0}$
since $\alp^{\dagger 0}=\Ome$ is the least in $SSt_{0}$.
Hence the condition
(\ref{eq:6.21.550}) 
is fulfilled, and we obtain
$0\in X=W_{N}^{0}(\emptyset)\cap 0^{\dagger N}$ with $D_{N}[X]$ by 
Lemma \ref{lem:6.21}.

Next let $\alp=\Ome$. 
Let $\gam\in \mathcal{C}^{\alp}(X)\cap\alp$.
We see $\gam\in X$ by induction on $\ell\gam$ using Lemma \ref{lem:6.15} and 
$0\in X$ since
there is no strongly critical number $\gam\in \mathcal{C}^{\alp}(X)\cap\alp$.
Therefore we obtain
$\alp\in\mathcal{G}^{X}$.

Let $\bet^{\dagger 0}<\alp^{\dagger 0}$. Then $\bet^{\dagger 0}=\Ome$ and $\bet<\Ome$.
Let $\gam\in W_{0}^{\bet}(X)\cap\Ome$.
We show $\gam\in X$.
We obtain 
$D_{0}^{0}[X\cap 0^{\dagger 1};X]$
 by Lemma \ref{lem:6.15}, 
 and 
$\gam\in W_{0}^{\bet}(X)\cap\Ome=
W_{0}^{0}(X)\cap\Ome=X\cap\Ome$. 
Hence the condition (\ref{eq:6.21.550}) is fulfilled, and we obtain
$\Ome\in\mathcal{W}_{N+1}$ by
Lemma \ref{lem:6.21}.
\eprf

\bprp\label{th:id5wf21_Omega}
For $\alpha=\psi_{\Omega}(a)\in OT(\mI_{N})$,
assume $\alpha\in\mathcal{G}^{\mathcal{W}_{N+1}}$.
Then 
$\alpha\in \mathcal{W}_{N+1}$.
\eprp
\bprf
Let $\alpha=\psi_{\Omega}(a)\in\mathcal{G}^{\mathcal{W}_{N+1}}$ and $X=W_{N}^{0}(\emptyset)\cap 0^{\dagger N}$.
In the proof of Proposition \ref{lem:6.29}, $0\in X$ and $D_{N}[X]$ are shown.
We obtain $\mathcal{W}_{N+1}\cap\Omega=X\cap\Omega$ by Lemma \ref{lem:8wf2}, and 
$\alpha\in\mathcal{G}^{X}$.
Since there is no $\beta$ such that $\beta^{\dagger 0}<\alpha^{\dagger 0}=\Omega$, 
(\ref{eq:6.21.550}) holds vacuously for $X$.
By Lemma \ref{lem:6.21} we obtain $D_{N}[X]$ and $\alpha\in X\subset\mathcal{W}_{N+1}$.
\eprf

\blem\label{lem:6.43}$(\Sig^{1}_{N+2}\mbox{{\rm -DC}})$\\
If 
$\alp\in\mathcal{G}^{\mathcal{W}_{N+1}}$, 
then there exists an $N$-distinguished {\rm set}
$Z$ such that $\{0,\Ome\}\subset Z$, 
$\alp\in\mathcal{G}^{Z}$ and
$\fal k\fal\mS\in Z\cap(St_{k}\cup\{\Ome\})[\mS^{\dagger k}\in Z]$.
\elem
\bprf
Let $\alp\in\mathcal{G}^{\mathcal{W}_{N+1}}$.
We have $\alp\in \mathcal{C}^{\alp}(\mathcal{W}_{N+1})$.
Pick an $N$-distinguished set $X_{0}$ such that $\alp\in \mathcal{C}^{\alp}(X_{0})$.
We can assume $\{0,\Ome\}\subset X_{0}$ by Proposition \ref{lem:6.29}.
On the other hand we have
$ \mathcal{C}^{\alp}(\mathcal{W}_{N+1})\cap\alp\subset\mathcal{W}_{N+1}$
and
$\fal k\fal\mS\in \mathcal{W}_{N+1}\cap(St_{k}\cup\{\Ome\})[\mS^{\dagger k}\in \mathcal{W}_{N+1}]$
by Lemma \ref{cor:6.21}.
We obtain 
\beqnarrs
&&
\fal n\fal X\exi Y
\{ 
D_{N}[X] \to D_{N}[Y]
\\
& \land &
\fal\bet\in OT(\mI_{N})
\left(
\ell(\bet)\leq n \land \bet\in \mathcal{C}^{\alp}(X)\cap\alp \to
\bet\in Y
\right)
\\
& \land &
\fal k \fal \mS\in (St_{k}\cup\{\Ome\})
\left( 
\ell (\mS)\leq n \land \mS\in X \to \mS^{\dagger k}\in Y
\right)
\}
\eeqnarrs

Since $D_{N}[X]$ is $\Del^{1}_{N+2}$, 
$\Sig^{1}_{N+2}\mbox{-DC}$ yields a set $Z$ such that $Z_{0}=X_{0}$ and
\beqnarrs
&&
\fal n
\{ 
D_{N}[Z_{n}] \to D_{N}[Z_{n+1}]
\\
& \land &
\fal\bet\in OT(\mI_{N})
\left(
\ell(\bet)\leq n \land \bet\in \mathcal{C}^{\alp}(Z_{n})\cap\alp \to
\bet\in Z_{n+1}
\right)
\\
& \land &
\fal k \fal \mS\in (St_{k}\cup\{\Ome\})
\left( 
\ell(\mS)\leq n \land \mS\in Z_{n} \to \mS^{\dagger k}\in Z_{n+1}
\right)
\}
\eeqnarrs
Let $Z=\bigcup_{n}Z_{n}$.
We see by induction on $n$ that $D_{N}[Z_{n}]$ for every $n$.
Lemma \ref{lem:8wf2}.\ref{lem:8wf2.5} yields $D_{N}[Z]$.
Let $\bet\in \mathcal{C}^{\alp}(Z)\cap\alp$.
Pick an $n$ such that $\bet\in \mathcal{C}^{\alp}(Z_{n})$ and $\ell\bet\leq n$.
We obtain $\bet\in Z_{n+1}\subset Z$.
Therefore $\alp\in\mathcal{G}^{Z}$.
Furthermore let $\mS\in Z\cap(St_{k}\cup\{\Ome\})$.
Pick an $n$ such that $\mS\in Z_{n}$ and $\ell(\mS)\leq n$.
We obtain $\mS^{\dagger k}\in Z_{n+1}\subset Z$.
\eprf

\begin{proposition}\label{lem:CMsmpl}
Let $D_{N}[Y]$ and
$\alpha\in \mathcal{C}^{\beta}(Y)$.
Assume $Y\cap\beta<\delta$.
Then $F_{\delta}(\alpha)\subset \mathcal{C}^{\beta}(Y)$.
\end{proposition}
\bprf
By induction on $\ell\alpha$.
Let $\{0,\Omega,\mathbb{I}_{N}\}\not\ni\alpha\in \mathcal{C}^{\beta}(Y)$.
We have $\mathcal{E}(\alp)\leq\alp$
First consider the case $\alpha\not\in\mathcal{E}(\alpha)$.
If $\alpha\in Y\cap\beta\subset\mathcal{G}^{Y}$ by Proposition \ref{prp:calg}, then 
$\mathcal{E}(\alpha)\subset \mathcal{C}^{\alpha}(Y)\cap\alpha\subset Y\subset \mathcal{C}^{\beta}(Y)$ by Proposition \ref{lem:3.11.632}.
Otherwise we have $\alpha\not\in\mathcal{E}(\alpha)\subset \mathcal{C}^{\beta}(Y)$.
In each case IH yields $F_{\delta}(\alpha)=F_{\delta}(\mathcal{E}(\alpha))\subset \mathcal{C}^{\beta}(Y)$.

Let $\alpha=\psi_{\pi}^{f}(a)$ for some $\pi,f,a$. 
If $\alpha<\delta$, then $F_{\delta}(\alpha)=\{\alpha\}$, and there is nothing to prove.
Let $\alpha\geq\delta$. Then $F_{\delta}(\alpha)=F_{\delta}(\{\pi,a\}\cup \mathrm{fld}(f))$.
On the other side we see $\{\pi,a\}\cup \mathrm{fld}(f)\subset \mathcal{C}^{\beta}(Y)$ from $\alpha\in \mathcal{C}^{\beta}(Y)$ 
and the assumption.
IH yields $F_{\delta}(\alpha)\subset \mathcal{C}^{\beta}(Y)$.

Finally let $\alp\in N(\rho)$.
Then $F_{\del}(\alp)=F_{\del}(\rho)$.
If $\rho\in \mathcal{C}^{\bet}(Y)$, then IH yields $F_{\del}(\rho)\subset \mathcal{C}^{\bet}(Y)$.
Otherwise we have $\alp\in Y$, and $\alp\in \mathcal{C}^{\alp}(Y)$.
Hence $\rho\in \mathcal{C}^{\alp}(Y)\cap\alp\subset Y\subset \mathcal{C}^{\bet}(Y)$.
\hspace*{\fill} $\Box$

\begin{proposition}\label{lem:CX3}
Let $\gamma<\beta$.
Assume $\alpha\in \mathcal{C}^{\gamma}(Y)$ and 
$G_{\bet}(\alpha)<\gamma$.
Moreover assume 
$\forall\delta[\ell\delta\leq\ell\alpha\spand\delta\in \mathcal{C}^{\gamma}(Y)\cap\gamma\Rarw\delta\in  \mathcal{C}^{\beta}(Y)]$.
Then $\alpha\in \mathcal{C}^{\beta}(Y)$.

\end{proposition}
\bprf
By induction on $\ell\alpha$. 
If $\alpha<\gamma$, then $\alpha\in  \mathcal{C}^{\gamma}(Y)\cap\gamma$.
The third assumption yields $\alpha\in  \mathcal{C}^{\beta}(Y)$. 
Assume $\alpha\geq\gamma$. 
Consider the case $\alpha=\psi_{\pi}^{f}(a)$ for some 
$\{\pi,a\}\cup \mathrm{fld}(f)\subset \mathcal{C}^{\gamma}(Y)$ and $\pi>\gamma$.
If $\pi\leq\beta$, then $\{\alpha\}=G_{\bet}(\alpha)<\gamma$ by the second assumption. Hence this is not the case, and we obtain 
$\pi>\beta$.
 Then $G_{\bet}(\{\pi,a\}\cup \mathrm{fld}(f))=G_{\bet}(\alpha)<\gamma$. 
IH yields $\{\pi,a\}\cup \mathrm{fld}(f)\subset \mathcal{C}^{\beta}(Y)$. We conclude $\alpha\in \mathcal{C}^{\beta}(Y)$ from $\pi>\beta$.

Next let $\gam\leq\alp\in N(\rho)$
with $\rho\in \mathcal{C}^{\gam}(Y)$.
If $\alp<\bet$, then $\{\alp\}=G_{\bet}(\alp)<\gam$, and this is not the case.
Let $\alp\geq\bet$. Then $G_{\bet}(\alp)=G_{\bet}(\rho)$.
IH yields $\rho\in \mathcal{C}^{\bet}(Y)$, and
$\alp\in \mathcal{C}^{\bet}(Y)$ by $\alp\geq\bet$.
\hspace*{\fill} $\Box$
\\

The following Lemma \ref{th:3wf16} is a key result on distinguished classes.

\begin{lemma}\label{th:3wf16}
Suppose $D_{N}[Y]$ 
with $\{0,\Ome\}\subset Y$, and for $\eta\in \Psi_{\mI_{N}}\cup\bigcup_{\mS\in SSt}L(\mS)$, cf.\,Definition \ref{df:R0},
\begin{equation}\label{eq:3wf16hyp.132}
\eta\in\mathcal{G}^{Y}
\end{equation}

\begin{equation}\label{eq:3wf16hyp.232}
R(\eta)\cap \{\gam\in OT(\mI_{N})\cap\mI_{N} : Y\cap\eta^{\dagger 1}<\gamma\}
\cap\mathcal{G}^{Y}\subset Y
\end{equation}
and
\begin{equation}\label{eq:3wf16hyp.332}
\eta\in L(\mS) \spand \mS=\mT^{\dagger k} \Rarw
\fal i<k(\mT^{\dagger i}\in  Y)
\end{equation}
Then
$\eta\in\mathcal{W}_{N+1}$.
Moreover if there exists a set $Z$ and an ordinal $\gam$ such that
$Y=W_{N}^{\gam}(Z)$ and $\eta^{-N}=\gam^{-N}$,
then $\eta\in Y$ holds.
\end{lemma}
\bprf
By Lemma \ref{lem:6.21} 
and the hypothesis (\ref{eq:3wf16hyp.132}) it suffices to show 
(\ref{eq:6.21.550})
\[
\fal\bet\geq\eta^{- 1}\left(
Y\cap\eta^{\dagger 1}<\bet \spand \bet^{\dagger 0}<\eta^{\dagger 0} \to
W_{0}^{\bet}(Y)\cap\bet^{\dagger 0}\subset Y
\right)
.\]
Assume $Y\cap\eta^{\dagger 1}<\bet$ and $\bet^{\dagger 0}<\eta^{\dagger 0}$. 
We have to show  
$W_{0}^{\bet}(Y)\cap\bet^{\dagger 0}\subset Y$. 
We prove this by induction on $\gamma\in W_{0}^{\bet}(Y)\cap\bet^{\dagger 0}$. 
Suppose $\gamma\in  \mathcal{C}^{\beta}(Y)\cap\bet^{\dagger 0}$ and 
\[
\mbox{{\rm MIH : }}  \mathcal{C}^{\beta}(Y)\cap\gamma\subset Y.
\]
We show $\gamma\in Y$. 
We can assume that
\begin{equation}
\label{eq:3wf9hyp.232X}
Y\cap\eta^{\dagger 1}<\gamma
\end{equation}
since if $\gamma\leq \delta$ for some $\delta\in Y\cap\eta^{\dagger 1}$, then by 
$Y\cap\eta^{\dagger 1}<\beta$ and $\gamma\in  \mathcal{C}^{\beta}(Y)$ 
we obtain $\delta<\beta$, 
$\gamma\in  \mathcal{C}^{\delta}(Y)$ and $\delta\in W( \mathcal{C}^{\delta}(Y))\cap\delta^{\dagger 0}=Y\cap\delta^{\dagger 0}$ by Lemma \ref{lem:6.15}.
Hence $\gamma\in W( \mathcal{C}^{\delta}(Y))\cap\delta^{\dagger 0}\subset Y$.

Moreover we can assume $\gam
\not\in (Reg_{0}\setm\{\Ome,\mI_{N}\})
\cap\bet$.
For otherwise $\gam\in Y$
by Definition \ref{df:CX}.\ref{df:CX.6} and $\gam\in \mathcal{C}^{\bet}(Y)\cap\bet$.

We show first 
\begin{equation}
\label{eq:3wf9hyp.232}
\gamma\in\mathcal{G}^{Y}
\end{equation}
First $\gamma\in  \mathcal{C}^{\gamma}(Y)$ by 
$\gamma\in  \mathcal{C}^{\beta}(Y)\cap\beta^{\dagger  0}$ and Proposition \ref{lem:CX2}. 
Second we show the following claim by induction on $\ell\alpha$:

\begin{equation}\label{clm:3wf1632}
\alpha\in \mathcal{C}^{\gamma}(Y)\cap\gamma \Rarw  \alpha\in Y
\end{equation}
\textbf{Proof} of (\ref{clm:3wf1632}). 
Assume $\alpha\in \mathcal{C}^{\gamma}(Y)\cap\gamma$.
We can assume $\gamma^{\dagger  0}\leq\beta$ for otherwise we have 
$\alpha\in  \mathcal{C}^{\gamma}(Y)\cap\gamma= \mathcal{C}^{\beta}(Y)\cap\gamma\subset Y$ by MIH.

By induction hypothesis on lengths, 
Proposition \ref{prp:noncritical}, and $\{0,\Ome\}\subset Y$, 
 we can assume that
$\alpha=\psi_{\pi}^{f}(a)$ for some $\pi>\gamma$ such that
 $\{\pi,a\}\cup \mathrm{fld}(f)\subset \mathcal{C}^{\gamma}(Y)$.
\\
\textbf{Case 1}. $\beta<\pi$: 
Then $G_{\bet}(\{\pi,a\}\cup \mathrm{fld}(f))=G_{\bet}(\alpha)<\alpha<\gamma$ by Proposition \ref{prp:G}.\ref{prp:G1}.
 Proposition \ref{lem:CX3} with induction hypothesis on lengths yields 
$\{\pi,a\}\cup \mathrm{fld}(f)\subset \mathcal{C}^{\beta}(Y)$.
Hence $\alpha\in  \mathcal{C}^{\beta}(Y)\cap\gamma$ by $\pi>\beta$.
MIH yields $\alpha\in Y$.
\\
\textbf{Case 2}. $\beta\geq\pi$: 
We have $\alpha<\gamma<\pi\leq\beta$. 
It suffices to show that $\alpha\leq Y\cap\eta^{\dagger 1}$.
Then by (\ref{eq:3wf9hyp.232X}) we have $\alpha\leq\delta\in Y\cap\eta^{\dagger  1}$ for some $\delta<\gamma$.
$ \mathcal{C}^{\delta}(Y)\ni\alpha\leq\delta\in Y\cap\delta^{\dagger  0}=W( \mathcal{C}^{\delta}(Y))\cap\delta^{\dagger  0}$
yields $\alpha\in W( \mathcal{C}^{\delta}(Y))\cap\delta^{\dagger  0}\subset Y$.

Consider first the case $\gamma\not\in\mathcal{E}(\gamma)$.
By $\alpha=\psi_{\pi}^{f}(a)<\gamma<\pi$,
we can assume that 
$\gamma\not\in\{0,\Omega,\mathbb{I}_{N}\}$.
Then let $\delta=\max S(\gamma)$ denote the largest immediate subterm of $\gamma$.
Then 
$\delta\in \mathcal{C}^{\gamma}(Y)\cap\gamma$, and 
by (\ref{eq:3wf9hyp.232X}), 
$Y\cap\eta^{\dagger  1}<\gamma\in \mathcal{C}^{\beta}(Y)$ we have 
$\delta\in  \mathcal{C}^{\beta}(Y)\cap\gamma$.
Hence $\delta\in Y\cap\eta^{\dagger 1}$ by MIH.
Also by $\alpha<\gamma$, we obtain $\alpha\leq\delta$, i.e., $\alpha\leq Y\cap\eta^{\dagger 1}$, and we are done.

Next let $\gam\not\in(Reg_{0}\setm\{\Ome,\mI_{N}\})$ and $\gamma\in\mathcal{E}(\gamma)$.
This means that $\gam\in\Psi$. Let
$\gamma=\psi_{\kappa}^{g}(b)$ for some $b,g$ and $\kappa>\beta$ by (\ref{eq:3wf9hyp.232X}) and $\gamma\in  \mathcal{C}^{\beta}(Y)$.
We have $\alpha<\gamma<\pi\leq\beta<\kappa$.
Let $\pi\preceq\rho$ and $\kappa\preceq\tau$
with $\{\rho,\tau\}\subset Reg_{0}$.
We obtain $\rho=\tau$ by Proposition \ref{prp:jumpover}.

$\pi\not\in C_{b}(\gamma)$ since otherwise by $\pi<\kappa$ we would have $\pi<\gamma$.
Then by Proposition \ref{prp:psicomparison} we have $a\geq b$ and 
$\mathrm{fld}(g)\cup\{\kappa,b\}\not\subset C_{a}(\alpha)$.
On the other hand we have 
$K_{\gamma}(\mathrm{fld}(g)\cup\{\kappa,b\})<b\leq a$, i.e.,
$\mathrm{fld}(g)\cup\{\kappa,b\}\subset C_{a}(\gamma)$.
By Proposition \ref{prp:G3} pick a 
$\delta\in F_{\gamma}(\mathrm{fld}(g)\cup\{\kappa,b\})$
such that $C_{a}(\alpha)\not\ni\delta\in\gamma$.
In particular 
$\delta<\gamma$.
Also we have $\mathrm{fld}(g)\cup\{\kappa,b\}\subset \mathcal{C}^{\beta}(Y)$,
$Y\subset\mathcal{G}^{Y}$ by Proposition \ref{prp:calg},
and $Y\cap\eta^{\dagger 1}<\gamma$ by (\ref{eq:3wf9hyp.232X}).
Therefore by Proposition \ref{lem:CMsmpl} with MIH
we obtain 
$\alpha\leq\delta\in \mathcal{C}^{\beta}(Y)\cap\gamma\subset Y$.

\hspace*{\fill} $\Box$ of (\ref{clm:3wf1632}) and (\ref{eq:3wf9hyp.232}).
\\

\noindent
Hence we obtain $\gamma\in\mathcal{G}^{Y}$.
We have $\gamma<\beta^{\dagger 0}\leq\eta$ and 
$\gamma\in \mathcal{C}^{\gamma}(Y)$.
If $\gamma\in R(\eta)$, then the hypothesis 
(\ref{eq:3wf16hyp.232}) yields $\gamma\in Y$.
In what follows assume $\gamma\not\in R(\eta)$.

If $G_{\eta}(\gamma)<\gamma$, then  Proposition \ref{lem:CX3} yields 
$\gamma\in \mathcal{C}^{\eta}(Y)\cap\eta\subset Y$ by $\eta\in\mathcal{G}^{Y}$.

In what follows suppose $G_{\eta}(\gamma)=\{\gamma\}$.
This means $\gam\in\Psi$ by $\gamma\not\in (Reg_{0}\setm\{\Ome,\mI_{N}\})$,
 and 
$\gamma\prec\tau$ for a $\tau<\eta$ by $\gamma\not\prec\eta$ and Definition \ref{df:EGFk}.\ref{df:EGFk.3}.
If $\eta\prec\mI_{N}$, then $\gam\prec\mI_{N}$ by $\gam\not\in R(\eta)$.
Hence this is not the case.

Let $\eta\in L(\mS)$ with $\mS=\mT^{\dagger k}\in SSt$.
By (\ref{eq:3wf16hyp.332}) and (\ref{eq:3wf9hyp.232X}) we obtain 
$\fal i<k(\mT^{\dagger i}\in Y\cap\eta^{\dagger 1}<\gam)$.
This yields $\gam\in L(\mS)$.
Let $\tau$ be maximal such that $\gamma\prec\tau<\eta$.

First assume $\eta\not\prec\mS$, and
let $\delta$ denote the ordinal such that $\delta\in N(\kappa_{0})$ with $\eta\preceq\kappa_{0}\in (L(\mS)\cup\{\mS\})\cap SSt^{M}$.
We see $\delta<\eta$.
On the other hand we have $\eta\leq\kappa\in \mathcal{C}^{\eta}(Y)$ by $\eta\in\mathcal{C}^{\eta}(Y)$, and 
$\delta\in \mathcal{C}^{\eta}(Y)$.
$\delta\in \mathcal{C}^{\eta}(Y)\cap\eta\subset Y$ follows.
We obtain $\delta<\gamma$ by (\ref{eq:3wf9hyp.232X}).
Then 
we obtain $\tau\in\Psi$ by 
$\gamma\not\in R(\eta)$ and Lemma \ref{lem:R}.
$\tau\in\Psi$ holds even if $\eta\prec\mS$.
From $\gamma\in \mathcal{C}^{\gamma}(Y)$ we see 
$\tau\in \mathcal{C}^{\gamma}(Y)$.
Next we show that
\begin{equation}\label{eq:3wf1632last}
G_{\eta}(\tau)<\gamma
\end{equation}
Let $\tau=\psi_{\kap}^{f}(b)$.
Then $\eta<\kap$ by the maximality of $\tau$, and
$G_{\eta}(\tau)=G_{\eta}(\{\kap,b\}\cup \mathrm{fld}(f))<\tau$ by Proposition \ref{prp:G}.\ref{prp:G1}.
On the other hand we have $\tau\in C_{a}(\gam)$.
Proposition \ref{prp:G}.\ref{prp:G2} yields
$G_{\eta}(\tau)\subset C_{a}(\gam)$.
We see $G_{\eta}(\tau)<\gam$ inductively.

(\ref{eq:3wf1632last}) is shown.
Proposition \ref{lem:CX3} yields 
$\tau\in \mathcal{C}^{\eta}(Y)$, and
$\tau\in \mathcal{C}^{\eta}(Y)\cap\eta\subset Y$ by $\eta\in\mathcal{G}^{Y}$.
Therefore $Y\cap\eta^{\dagger 1}<\gamma<\tau\in Y$.
This is not the case by (\ref{eq:3wf9hyp.232X}).
We are done.
\hspace*{\fill} $\Box$

\bprp\label{prp:wfpsiI}
For $\alp_{1}=\psi_{\mI_{N}}(a)$,
$\alp_{1}\in\mathcal{G}^{\mathcal{W}_{N+1}} \Rarw \alp_{1}\in\mathcal{W}_{N+1}$.
\eprp
\bprf
By Lemma \ref{lem:6.43} pick an $N$-distinguished set
$Z$ such that $\{0,\Ome\}\subset Z$, 
$\alp_{1}\in\mathcal{G}^{Z}$ and
$\fal k\fal\mS\in Z\cap(St_{k}\cup\{\Ome\})[\mS^{\dagger k}\in Z]$.

\bclm\label{clm:wfpsiI.2}
Let $SSt\ni\mT<\alp_{1}$ and $\mathcal{G}^{Z}\ni\gam\in L(\mT)\cap\Psi$.
Then $\gam<Z\cap\alp_{1}$.
\eclm
\textbf{Proof} of Claim \ref{clm:wfpsiI.2}.
Let $\rho\prec\mS^{\dagger \vec{k}}=\mT<\alp_{1}$ for an $\mS\in LSt\cup\{\Ome\}$ and a $\vec{k}\neq\emptyset$.
Also $\gam\prec^{R}\kap\in N(\rho)$.
We show $\rho\in Z$ by induction on $\ell\gam$.
First let $\gam=\psi_{\mI_{N}[\sig]}(b)$ for some $b$ and $\sig\preceq^{R}\kap$.
Then  we obtain $\mI_{N}[\sig]\in \mathcal{C}^{\gam}(Z)$ by $\gam\in \mathcal{C}^{\gam}(Z)$,
and $\sig\in \mathcal{C}^{\gam}(Z)\cap\gam\subset Z$.
Proposition \ref{prp:calg} yields $\sig\in\mathcal{G}^{Z}$.
If $\sig=\kap=\mI_{N}[\rho]$, then $\sig\in \mathcal{C}^{\sig}(Z)$ yields
$\rho\in \mathcal{C}^{\sig}(Z)\cap\sig\subset Z$.
Otherwise IH yields $\rho\in Z$.
Second let $\gam=\psi_{\sig^{\dagger \vec{i}}}^{f}(b)\in \mathcal{C}^{\gam}(Z)$ 
for some $f$, $b$ and 
$\sig^{\dagger\vec{i}}\preceq^{R}\kap$.
We obtain $\sig\in \mathcal{C}^{\gam}(Z)\cap\gam\subset Z$, and 
$\sig\in\mathcal{G}^{Z}$.
We obtain $\sig\prec^{R}\kap$.
IH yields $\rho\in Z$.
Third let $\gam=\psi_{\tau}^{f}(a)$ with
$\tau=\mW^{\dagger\vec{j}}[\sig/\mW]$.
We obtain $\gam<\tau\in \mathcal{C}^{\gam}(Z)$, and $\sig\in \mathcal{C}^{\gam}(Z)\cap\gam$.
Hence $\sig\in\mathcal{G}^{Z}$.
If $\tau=\mW^{\dagger\vec{j}}[\sig/\mW]=\mU^{\dagger\vec{i}}[\rho/\mS]$, then
$\sig\in \mathcal{C}^{\sig}(Z)$ yields
$\rho\in \mathcal{C}^{\sig}(Z)\cap\sig\subset Z$.
Otherwise IH yields $\rho\in Z$.

Now $\rho\in Z$ yields $\rho\in \mathcal{C}^{\rho}(Z)$, and this yields
$\mS\in \mathcal{C}^{\rho}(Z)\cap\rho\subset Z$.
Since $Z$ is closed under $\mT\mapsto\mT^{\dagger i}$, we obtain
$\gam<\mS^{\dagger \vec{k}}\in Z\cap\alp_{1}$.
\eprf\,of Claim \ref{clm:wfpsiI.2}.

Since there is no $\gam\prec\alp_{1}$, if $\gam\in R(\alp_{1})$,
then $\gam\in L(\mT)\cap\Psi$ for a $SSt\ni\mT<\alp_{1}$ by Definition \ref{df:R0}.\ref{df:R00}.
We conclude $\alp_{1}\in\mathcal{W}_{N+1}$
by Lemma \ref{th:3wf16}.
\eprf

\begin{lemma}\label{lem:wfpsiI}
For {\rm each} $n<\ome$, the following holds:

Let
$a\in\mathcal{C}^{\mI_{N}}(\mathcal{W}_{N+1})\cap\ome_{n}(\mI_{N}+1)$. 
Then
$
\psi_{\mI_{N}}(a)\in\mathcal{W}_{N+1}$ holds.

\end{lemma}
\bprf
For each $n<\ome$, we have $\mathrm{TI}[\mathcal{C}^{\mI_{N}}(\mathcal{W}_{N+1})\cap(\ome_{n}(\mI_{N}+1))]$ by 
Lemma \ref{lem:6.34}.\ref{lem:6.34.3}.
We show the lemma by induction on $a\in\mathcal{C}^{\mI_{N}}(\mathcal{W}_{N+1})\cap\ome_{n}(\mI_{N}+1)$.
Assume
\[
\mbox{IH} :\Lrarw \fal b\in\mathcal{C}^{\mI_{N}}(\mathcal{W}_{N+1})\cap a\left(
\psi_{\mI_{N}}(b)\in OT(\mI_{N}) \Rarw \psi_{\mI_{N}}(b)\in\mathcal{W}_{N+1}
\right)
\]
Let $\alp_{1}=\psi_{\mI_{N}}(a)\in OT(\mI_{N})$ with $a\in\mathcal{C}^{\mI_{N}}(\mathcal{W}_{N+1})\cap\ome_{n}(\mI_{N}+1)$.
By Proposition \ref{prp:wfpsiI}
it suffices to show 
$\alpha_{1}\in\mathcal{G}^{\mathcal{W}_{N+1}}$.

From $a\in\mathcal{C}^{\mI_{N}}(\mathcal{W}_{N+1})$ with $\alp_{1}<\mI_{N}$ we see
$\alpha_{1}\in\mathcal{C}^{\alpha_{1}}(\mathcal{W}_{N+1})$.
It suffices to show the following claim by induction on $\ell\beta_1$.
\begin{equation}\label{clm:wfpsiI.1}
\forall\beta_{1}\in\mathcal{C}^{\alpha_{1}}(\mathcal{W}_{N+1})\cap\alpha_{1}[\beta_1\in\mathcal{W}_{N+1}].
\end{equation}
\textbf{Proof} of (\ref{clm:wfpsiI.1}). 
Assume $\beta_{1}\in\mathcal{C}^{\alpha_{1}}(\mathcal{W}_{N+1})\cap\alpha_{1}$ and let
\[
\mbox{LIH} :\Lrarw
\forall\gamma\in\mathcal{C}^{\alpha_{1}}(\mathcal{W}_{N+1})\cap\alpha_{1}[\ell\gamma<\ell\beta_{1} \Rarw \gamma\in\mathcal{W}_{N+1}].
\]

We show $\beta_1\in\mathcal{W}_{N+1}$. 
We can assume $\bet_{1}\not\in\{0,\Ome\}$ by Proposition \ref{lem:6.29}.
\\
\textbf{Case 1}. 
$\beta_1\not\in\mathcal{E}(\beta_{1})$:
Assume $\beta_{1}\not\in\mathcal{W}_{N+1}$.
Then 
$\bet_{1}\not\in N(\rho)$
for any $\rho$ by 
$\beta_{1}\in\mathcal{C}^{\alpha_{1}}(\mathcal{W}_{N+1})\cap\alpha_{1}$ and
Definition \ref{df:CX}. We obtain
$S(\beta_{1})\subset\mathcal{C}^{\alpha_{1}}(\mathcal{W}_{N+1})\cap\alpha_{1}$. 
LIH yields $S(\beta_{1})\subset\mathcal{W}_{N+1}$. 
Hence we conclude $\beta_{1}\in\mathcal{W}_{N+1}$ from Proposition \ref{prp:noncritical}.
\\
\textbf{Case 2}.
In what follows consider the cases when $\beta_{1}=\psi_{\pi}^{g}(b)$ for some $\pi,b,g$.
We can assume $\pi>\alp_{1}$. Then we see $\pi=\mI_{N}$ and $\beta_{1}=\psi_{\mI_{N}}(b)$ with
$b\in\mathcal{C}^{\alpha_{1}}(\mathcal{W}_{N+1})$.
We obtain
$b<a$ by Proposition \ref{prp:comparisonrud}.\ref{prp:comparisonrud.1},
and $b\in C_{b}(\bet_{1})$.
By IH it suffices to show $b\in\mathcal{C}^{\mI_{N}}(\mathcal{W}_{N+1})$.

By induction on $\ell c$ we see that
$c\in C_{b}(\bet_{1}) \Rarw G_{\mI_{N}}(c)<\bet_{1}$.
For example let $c=\gam_{1}^{\dagger \vec{i}}$ with $\gam_{1}\in LSt_{N}\cup\{\Ome\}$ and $\vec{i}\neq\emptyset$.
Suppose $c>\bet_{1}$. Then $\gam_{1}\in C_{b}(\bet_{1})$.
IH yields $\{\gam_{1}\}=G_{\mI_{N}}(\gam_{1})<\bet_{1}\in LSt_{N}$, and hence $\{c\}=G_{\mI_{N}}(c)<\bet_{1}$.

In particular we obtain $G_{\mI_{N}}(b)<\bet_{1}$.
Proposition \ref{lem:CX3} with LIH yields $b\in\mathcal{C}^{\mI_{N}}(\mathcal{W}_{N+1})$.
This shows (\ref{clm:wfpsiI.1}).
\eprf

\subsection{Layers of stable ordinals}

In this subsection we examine ordinals in layers 
$L(\mS)=\{\alp\in OT(\mI_{N}): \alp\prec^{R}\mS\}$ for $\mS\in SSt$.
We show that there is no infinite descending chain in $L(\mS)$,
cf.\,Lemma \ref{lem:psimS}.
Here we need the condition (\ref{eq:notationsystem.6}) and the fact that
 $\alp\in M_{\rho}$ if $\alp$ is in
the domain of the Mostowski collapsing 
$\alp\mapsto \alp[\rho/\mS]$, cf.\,Definition \ref{df:notationsystem.2M}
and Proposition \ref{prp:psimS.2}.

Let $k(\psi_{\kap}^{f}(a))=\{\kap,a\}\cup \mathrm{fld}(f)$ and
$\mathtt{ h}(\psi_{\kap}^{f}(a))=\{a,\mathtt{g}_{0}^{*}(\psi_{\kap}^{f}(a))\}$.
For a set $\Tht\subset OT(\mI_{N})$,
When $\kappa\in\{\Omega,\mI_{N}\}$, let $k(\psi_{\kappa}(a))=\mathtt{h}(\psi_{\kappa}(a))=\{a\}$ and
$N(\psi_{\kappa}(a))=\emptyset$.

\bprp\label{prp:psimS}

Let $Z$ be an $N$-distinguished set such that $\{0,\Ome\}\subset Z$
and
$\fal k\fal\mS\in Z\cap(St_{k}\cup\{\Ome\})[\mS^{\dagger k}\in Z]$.
Assume $\psi_{\mI_{N}}(b)\in Z$, and let
\beqnarrs
&&
\mbox{{\rm MIH}}(b;Z) :\Lrarw 
\fal\gam\in \Psi
\\
&&
\left[
k(\gam)
\subset\mathcal{C}^{\mI_{N}}(Z)\spand 
\mathtt{h}(\gam)\subset \mathcal{C}^{\mI_{N}}(Z)\cap b
 \Rarw 
\{\gam\}\cup
N(\gam)
\subset Z
\right].
\eeqnarrs

Then for any $\Tht\subset Z$,
$
C_{b}(\Tht)\subset
\mathcal{C}^{\mI_{N}}(Z)
$ holds.

\eprp
\bprf 
Let $\Tht\subset Z$.
Assuming $\gam\in C_{b}(\Tht)$,
we show $\gam\in\mathcal{C}^{\mI_{N}}(Z)$ by induction on $\ell\gam$.

Let $\gam\not\in \Tht$.
By IH and Proposition \ref{prp:noncritical}, we can assume 
$\gam\in\Psi\cup(Reg_{0}\setm\{\Ome,\mI_{N}\})$.
If $\gam\in N(\gam_{1})$
for a $\gam_{1}\in\Psi$, then
 $\gam_{1}\in C_{b}(\Tht)$.
 Hence we may assume that
$\gam=\psi_{\kap}^{f}(a)$
with $k(\gam)\subset C_{b}(\Tht)$.
IH yields $\{\kap,a\}\subset k(\gam)\subset \mathcal{C}^{\mI_{N}}(Z)$.
If $\kap=\mI_{N}$ and $f=\emptyset$, then we obtain 
$\gam=\psi_{\mI_{N}}(a)<\psi_{\mI_{N}}(b)=\del\in Z$.
$a\in\mathcal{C}^{\mI_{N}}(Z)\subset\mathcal{C}^{\del}(Z)$ yields 
$\gam\in\mathcal{C}^{\del}(Z)\cap\del\subset Z$.
If $\gamma=\psi_{\Omega}(a)$, then we obtain $a\in \mathcal{C}^{\mI_{N}}(Z)\cap b$.
$\mbox{{\rm MIH}}(b;Z)$ yields 
$\gamma\in Z$.

Let $\kap<\mI_{N}$ and $\gam\in L(\mS)$ with $\mS=\mT^{\dagger k}$ and $\mT\in St\cup\{\Ome\}$.
We claim that $\mT\in Z$ and
$\mathtt{h}(\gam)\subset\mathcal{C}^{\mI_{N}}(Z)\cap b$.
We have $\kap\in\mathcal{C}^{\mI_{N}}(Z)\cap\mI_{N}\subset Z$.
We obtain $\kap\in\mathcal{G}^{Z}$.
Let $\rho\prec\mS$ be such that either $\rho=\kap$ or
$\kap\prec^{R}\sig\in N(\rho)$.
We obtain $\rho\in Z$ and $\rho\in\mathcal{G}^{Z}$, from which we see
$\mS\in\mathcal{C}^{\rho}(Z)$ and
$\mT\in\mathcal{C}^{\rho}(Z)\cap\rho\subset Z$.

On the other,
IH yields $a\in \mathcal{C}^{\mI_{N}}(Z)\cap b$.
We show $\mathtt{ g}^{*}_{0}(\gam)\in\mathcal{C}^{\mI_{N}}(Z)\cap b$.

If $\gam\prec\mS$, then
$\mathtt{ g}^{*}_{0}(\gam)=\mathtt{ p}_{\mS}(\gam_{0})$ for
$\gam\preceq\gam_{0}=\psi_{\mS}^{g}(c)$.
IH with $\mathtt{ p}_{\mS}(\gam_{0})\in C_{b}(\Tht)$
yields $\mathtt{ p}_{\mS}(\gam_{0})\in\mathcal{C}^{\mI_{N}}(Z)$.
On the other hand we have
$\mathtt{ p}_{\mS}(\gam_{0})<b$ by 
$\gam_{0}\in C_{b}(\Tht)$.

Next let $\rho\prec\mS$ and $\gam\prec^{R}\sig\in N(\rho)$.
Then $\mathtt{ g}^{*}_{0}(\gam)=\mathtt{ g}^{*}_{0}(\rho)$.
We obtain $\rho\in C_{b}(\Tht)$.
IH yields $\mathtt{ g}^{*}_{0}(\rho)\in\mathcal{C}^{\mI_{N}}(Z)\cap b$.

Therefore $\mbox{{\rm MIH}}(b;Z)$ yields 
$\{\gam\}\cup N(\gam)
\subset Z$.
\eprf

\bprp\label{prp:psim-1}
\benu
\item\label{prp:psim-1.1}
Let $\gam_{1}=\gam[\rho/\mS]^{-1}$ be the Mostowski uncollapsing,
and $\{\mS,\gam\}\subset \mathcal{C}^{\rho}(Z)$.
Then $\gam_{1}\in \mathcal{C}^{\rho}(Z)$.

\item\label{prp:psim-1.2}
$\gam\in C_{b}(\rho)\cap \mathcal{C}^{\rho}(Z) \Rarw \gam\in C_{b}( \mathcal{C}^{\rho}(Z)\cap\rho)$.

\eenu
\eprp
\bprf
Both are seen by induction on $\ell\gam$.
For Proposition \ref{prp:psim-1}.\ref{prp:psim-1.1},
use the fact $\gam_{1}=\gam[\rho/\mS]^{-1}\geq\gam$.
\eprf

\bprp\label{prp:psimS.2}

Let $\mS\in SSt$,
$\eta\in L(\mS)$
and
$Z$ be an $N$-distinguished set such that $\{0,\Ome\}\subset Z$,
$\fal k\fal\mS\in Z\cap(St_{k}\cup\{\Ome\})[\mS^{\dagger k}\in Z]$.
Assume $\eta\in\mathcal{G}^{Z}$, $\psi_{\mI_{N}}(b)\in Z$ and
$\mbox{{\rm MIH}}(b;Z)$
in Proposition \ref{prp:psimS} for a $b\geq\mathtt{ g}^{*}_{0}(\eta)$.
Then the following holds.

\benu
\item\label{prp:psimS.20}
$\mathtt{g}_{0}(\eta)\in \mathcal{C}^{\mI_{N}}(Z)$.

\item\label{prp:psimS.21}
$\mathtt{g}_{1}(\eta)\in \mathcal{C}^{\mI_{N}}(Z)$.

\item\label{prp:psimS.22}
$\mathtt{g}_{2}(\eta)\in \mathcal{C}^{\mI_{N}}(Z)$.
\eenu
\eprp
\bprf
Proposition \ref{prp:psimS.2}.\ref{prp:psimS.21} is seen from Proposition \ref{prp:psimS.2}.\ref{prp:psimS.20} by induction on $\ell\eta$.
\\
\ref{prp:psimS.2}.\ref{prp:psimS.20}. Let $\eta\in\Psi$.
\\
\textbf{Case 1}. 
$\eta=\rho$ or $\eta=\psi_{\mI_{N}[\rho]}(c)$ for a $\rho\prec\mS$ and a $c$: 
 Then 
 $\mathtt{p}_{0}(\rho)\leq\mathtt{ g}_{0}(\rho)=\mathtt{ g}_{0}(\eta)=\mathtt{ g}^{*}_{0}(\eta)\leq b$.
We show $\mathtt{p}_{\mS}(\rho)=\mathtt{g}_{0}(\rho)\in\mathcal{C}^{\mI_{N}}(Z)$.
By (\ref{eq:notationsystem.6}) in Definition \ref{df:notationsystem}.\ref{df:notationsystem.6}
we have $\mathtt{p}_{0}(\rho)\in C_{b}(\rho)$,
and $\mathtt{p}_{\mS}(\rho)\in C_{b}(\rho)$.
On the other hand we have $\rho\in\mathcal{G}^{Z}$.
We obtain $\mathtt{p}_{\mS}(\rho)\in\mathcal{C}^{\rho}(Z)$ by $\rho\in\mathcal{C}^{\rho}(Z)$, and
$\mathtt{p}_{\mS}(\rho)\in C_{b}(\mathcal{C}^{\rho}(Z)\cap\rho)$
by Proposition \ref{prp:psim-1}.\ref{prp:psim-1.2}.
Moreover we have
 $\calc^{\rho}(Z)\cap\rho\subset Z$.
Proposition \ref{prp:psimS} with $\mbox{{\rm MIH}}(b;Z)$ yields 
$\mathtt{p}_{\mS}(\rho)\in C_{b}(\mathcal{C}^{\rho}(Z)\cap\rho)\subset\mathcal{C}^{\mI_{N}}(Z)$.
\\
\textbf{Case 2}. Otherwise:
Let $\rho\prec\mS$ be such that
$\eta\prec^{R}\tau\in N(\rho)$.
Let $\eta_{1}\in M_{\rho}$ be such that $\eta=\eta_{1}[\rho/\mS]$.
Then $\mathtt{ g}_{0}(\eta)=\mathtt{ g}_{0}(\eta_{1})$ and $\mathtt{p}_{0}(\rho)\leq\mathtt{g}_{0}(\rho)=\mathtt{g}^{*}_{0}(\eta)\leq b$.

On the other hand we have $\eta\in\mathcal{G}^{Z}$.
$\eta\in\mathcal{C}^{\eta}(Z)$ yields $\rho\in\mathcal{C}^{\eta}(Z)\cap\eta\subset Z$.
Hence $\rho\in Z$.
We obtain $\rho\in\mathcal{G}^{Z}$.
We see  $\mS\in\calc^{\rho}(Z)$ from $\rho\in\calc^{\rho}(Z)$.
Hence $\{\mS,\eta\}\subset\mathcal{C}^{\rho}(Z)$.
Proposition \ref{prp:psim-1}.\ref{prp:psim-1.1} yields $\eta_{1}\in \mathcal{C}^{\rho}(Z)$,
and $\mathtt{g}_{0}(\eta_{1})\in \mathcal{C}^{\rho}(Z)$ by $\eta_{1}>\rho$.
$\eta_{1}\in M_{\rho}\subset  C_{b}(\rho)$ yields $\mathtt{g}_{0}(\eta_{1})\in  C_{b}(\rho)$.
We obtain $\mathtt{g}_{0}(\eta_{1})\in  C_{a}(\mathcal{C}^{\rho}(Z)\cap\rho)$
by Proposition \ref{prp:psim-1}.\ref{prp:psim-1.2}.
$\rho\in\mathcal{G}^{Z}$ yields $\mathcal{C}^{\rho}(Z)\cap\rho\subset Z$.
Hence Proposition \ref{prp:psimS} yields $\mathtt{g}_{0}(\eta_{1})\in\mathcal{C}^{\mI_{N}}(Z)$.
\\
\ref{prp:psimS.2}.\ref{prp:psimS.22}.
By induction on $\ell\eta$. Let $\eta\in\Psi$.
\\
\textbf{Case 1}. $\eta=\rho\prec\mS$:
Then $\mathtt{ p}_{0}(\rho)=\mathtt{ p}_{0}(\eta)$.
Let 
$f=m(\rho)$
and
$\zeta
\geq b \geq \mathtt{ p}_{0}(\eta)$.
Then $\mathtt{ g}_{2}(\rho)=o_{\mS}(f)+1$ and 
$\mathrm{fld}(f)\subset  C_{b}(\rho)$
by (\ref{eq:notationsystem.6}) in 
Definition \ref{df:notationsystem}.\ref{df:notationsystem.6}.
Moreover $\mathrm{fld}(f)\subset \calc^{\rho}(Z)$ by $\rho\in\calc^{\rho}(Z)$.
Hence we obtain $\mathrm{fld}(f)\subset  C_{b}(\calc^{\rho}(Z)\cap\rho)$
by Proposition \ref{prp:psim-1}.\ref{prp:psim-1.2},
where $\calc^{\rho}(Z)\cap\rho\subset Z$.
Proposition \ref{prp:psimS} with $\mbox{{\rm MIH}}(b;Z)$ yields 
$\mathrm{fld}(f)\subset\mathcal{C}^{\mI_{N}}(Z)$, and 
$o_{\mS}(f)\in\mathcal{C}^{\mI_{N}}(Z)$.
\\
\textbf{Case 2}. Otherwise: 
Let $\rho\prec\mS$ be such that
$\eta=\eta_{1}[\rho/\mS]$ 
with $\eta_{1}\in M_{\rho}\cap L(\mS_{1})$ and $\mathtt{ g}_{0}(\eta_{1})=\mathtt{ g}_{0}(\eta)$,
where $\eta\prec^{R} (\mS_{1}[\rho/\mS])$ and $M_{\rho}\subset C_{b}(\rho)$ for 
$\mathtt{p}_{0}(\rho)\leq\mathtt{g}_{0}(\rho)=\mathtt{g}^{*}_{0}(\eta)\leq b$.
Then $\ell\eta_{1}<\ell\eta$.
$\eta\in\calc^{\eta}(Z)$ with $\rho<\eta$ yields $\eta\in\calc^{\rho}(Z)$ and
$\rho\in\calc^{\eta}(Z)\cap\eta\subset Z$.
We see  $\mS\in\calc^{\rho}(Z)$ from $\rho\in\calc^{\rho}(Z)$.
We obtain $\eta_{1}\in\calc^{\rho}(Z)$ by Proposition \ref{prp:psim-1}.\ref{prp:psim-1.1}
 and $\eta\in\calc^{\rho}(Z)$, and
$\eta_{1}\in C_{b}(\calc^{\rho}(Z)\cap\rho)$
by Proposition \ref{prp:psim-1}.\ref{prp:psim-1.2}.
On the other hand we have $\calc^{\rho}(Z)\cap\rho\subset Z$.
By Proposition \ref{prp:psimS} we obtain 
$\eta_{1}\in\calc^{\mI_{N}}(Z)\cap\mI_{N}\subset Z$.
Hence $\eta_{1}\in\calg^{Z}$.
Moreover we see $\mathtt{g}^{*}_{0}(\eta_{1})<b=\mathtt{p}_{0}(\rho)\leq\mathtt{g}^{*}_{0}(\eta)$
from $\eta_{1}\in C_{b}(\rho)$.
IH yields
$\mathtt{ g}_{2}(\eta)=\mathtt{ g}_{2}(\eta_{1})\in\mathcal{C}^{\mI_{N}}(Z)$.
\eprf

\blem\label{lem:psimS}
Let $\mS=\mT^{\dagger k}\in SSt$ with $\mT\in \{\Ome\}\cup(LSt\cap\Psi)$, $\mT<\eta\in L(\mS)$,
and
$Z$ be an $N$-distinguished set such that $\{0,\Ome, \mT\}\subset Z$,
$\fal k\fal\mU\in Z\cap(St_{k}\cup\{\Ome\})[\mU^{\dagger k}\in Z]$.
Assume $\eta\in\mathcal{G}^{Z}$,  
$b\geq\mathtt{ g}^{*}_{0}(\eta)$,
$\Lam=\psi_{\mI_{N}}(b)\in Z$, and
$\mbox{{\rm MIH}}(b;Z)$
in Proposition \ref{prp:psimS}.
Then 
$\eta\in Z$.

\elem
\bprf
By Lemma \ref{lem:LSwinding} we obtain 
$SC(\mathtt{ g}_{2}(\eta))\subset \Lam=\psi_{\mI_{N}}(b)$.
An ordinal $\mathtt{ g}_{2}^{\Lam}(\eta)=o_{\Lam}(f)+1<\mI_{N}$ is obtained from
$\mathtt{ g}_{2}(\eta)=o_{\mI_{N}}(f)+1$ in Definition \ref{df:Lam.of}.\ref{df:Lam.of.2}
by changing the base $\mI_{N}$ to $\Lam$.
Then for $SC(\mathtt{ g}_{2}(\gam))\cup SC(\mathtt{ g}_{2}(\del))\subset \Lam$,
$\mathtt{ g}_{2}(\del)<\mathtt{ g}_{2}(\gam) \Lrarw \mathtt{ g}_{2}^{\Lam}(\del)<\mathtt{ g}_{2}^{\Lam}(\gam)$ by
Proposition \ref{prp:thtcollapse}, and
$\mathtt{ g}_{2}(\gam)\in\mathcal{C}^{\mI_{N}}(Z) \Lrarw 
\mathtt{ g}_{2}^{\Lam}(\gam)\in Z$ by the assumption $\Lam\in Z$.

On the other side, we see
$W_{N}^{\mT}(Z)\cap\mS=Z\cap\mS$ from $\mT\in Z$ and $D_{N}[Z]$.
Hence $\mathcal{G}^{Y}\cap\mS=\mathcal{G}^{Z}\cap\mS$ for $Y=W_{N}^{\mT}(Z)\cap\mS$.

We see $Wo[\mathcal{C}^{\mI_{N}}(Z)]$ from $\mathcal{C}^{\mI_{N}}(Z)\cap\mI_{N}=Z\cap\mI_{N}$
as in Lemma \ref{lem:6.34}.
We show $\eta\in Z\cap\mS=W_{N}^{\mT}(Z)\cap\mS$ 
 by induction on 
$\mathtt{ g}^{\Lam}(\eta)=(\mathtt{ g}_{1}(\eta),\mathtt{ g}_{2}^{\Lam}(\eta))$
with respect to the lexicographic order $<_{lx}$ on $\mathcal{C}^{\mI_{N}}(Z)\times Z$.

Let $\gam\in R(\eta)$ be such that 
$\gam\in\mathcal{G}^{Z}$. 
Then $\gam\in R(\eta)\subset L(\mS)$, $\mT^{-N}=\gam^{- N}=\eta^{- N}$ and
$\mT<\gam<\eta<\mS$.
By Lemma \ref{lem:LSwinding} we obtain $\mathtt{ g}^{*}_{0}(\gam)\leq\mathtt{ g}^{*}_{0}(\eta)$,
$\mathtt{ g}(\gam)<_{lx}\mathtt{ g}(\eta)$ and
$SC(\mathtt{ g}_{2}(\gam))\subset \Lam=\psi_{\mI_{N}}(b)$.
Proposition \ref{prp:psimS.2} yields 
$\{\mathtt{g}_{1}(\gam),\mathtt{g}_{2}(\gam),\mathtt{g}_{1}(\eta),\mathtt{g}_{2}(\eta)\}\subset \mathcal{C}^{\mI_{N}}(Z)$.
We obtain $\mathtt{ g}^{\Lam}(\gam)<_{lx}\mathtt{ g}^{\Lam}(\eta)$.
IH yields $\gam\in Z$, and (\ref{eq:3wf16hyp.232}) is shown.
Lemma \ref{th:3wf16}  yields $\eta\in Z$.
\eprf

\bprp\label{prp:MG}
Let $D_{N}[Z]$ and $\rho\in L(\mS)\cap Z\cap\Psi$ with an $\mS\in SSt$.
Then $N(\rho)
\subset\mathcal{G}^{Z}$.
\eprp
\bprf
Let $\alp\in N(\rho)$.
We obtain $\alp\in\mathcal{C}^{\alp}(Z)$ by $\rho\in Z\cap\alp$.
We show
$\bet\in\mathcal{C}^{\alp}(Z)\cap\alp \Rarw \alp\in Z$
by induction on $\ell\bet$.
Let $\rho\neq\bet\in\mathcal{C}^{\alp}(Z)\cap\alp$. 
If $\bet<\rho$, then $\bet\in\mathcal{C}^{\rho}(Z)\cap\rho\subset Z$ by
Propositions \ref{lem:CX2}.\ref{lem:CX2.3} and \ref{prp:calg}.
Let $\rho<\bet<\alp$.
By IH, Proposition \ref{prp:noncritical} and Definition \ref{df:CX}
we may assume that $\bet=\psi_{\sig}^{f}(c)$ with $\sig>\alp$.
Then $\bet<\rho$ by Proposition \ref{prp:jumpover}.
\eprf

\bcor\label{cor:MG}
For {\rm each} $\zeta\in \mathcal{C}^{\mI_{N}}(\mathcal{W}_{N+1})$, the following holds:

Let $\mS=\mT^{\dagger k}\in SSt$ with $\mT\in St\cup\{\Ome\}$, 
$\eta\in N(\rho)$
with $\rho\in L(\mS)$,
and
$\{ \mT,\rho\}\subset \mathcal{W}_{N+1}$.
Assume
$\zeta\geq \mathtt{ g}^{*}_{0}(\eta)$
and
$\mbox{{\rm MIH}}(\zeta;\mathcal{W}_{N+1})$
in Proposition \ref{prp:psimS}.
Then 
$\eta\in \mathcal{W}_{N+1}$.
\ecor
\bprf
By $\mathtt{ g}^{*}_{0}(\eta)\leq \zeta\in\mathcal{C}^{\mI_{N}}(\mathcal{W}_{N+1})$
and Lemma \ref{lem:wfpsiI}
we obtain $\psi_{\mI_{N}}(\zeta)\in \mathcal{W}_{N+1}$.
As in the proof of Lemma \ref{lem:6.43}
we see that there exists 
an $N$-distinguished set $Z$ such that $\{0,\Ome, \mT,\rho\}\subset Z$,
$\fal k\fal\mU\in Z\cap(St_{k}\cup\{\Ome\})[\mU^{\dagger k}\in Z]$,
$\psi_{\mI_{N}}(\zeta)\in Z$, and
$\mbox{{\rm MIH}}(\zeta;Z)$.
Then $\eta\in Z\subset\mathcal{W}_{N+1}$ follows from Lemma \ref{lem:psimS}
and Proposition \ref{prp:MG}.
\eprf

\begin{definition}\label{df:id4wfA}
{\rm 
For irreducible functions $f:\Lambda\to\Gamma(\Lambda)$ let
\[
f\in J:\Lrarw \mathrm{fld}(f)\subset\mathcal{W}_{N+1}.
\]
Let $\mathtt{g}^{*}_{0}(\psi_{\Omega}(a))=0$.
For $a\in OT(\mI_{N})$ and irreducible functions $f$, define:
\[
 A(\zeta,a,f)  :\Lrarw 
 \forall\sigma\in \mathcal{W}_{N+1}\cap\mI_{N}
 [
 \psi_{\sigma}^{f}(a)\in OT(\mI_{N}) \spand
 \mathtt{g}_{0}^{*}(\psi_{\sigma}^{f}(a))\leq\zeta
 \Rarw \psi_{\sigma}^{f}(a)\in\mathcal{W}_{N+1}]
\]
\beqnarrs
\mbox{{\rm SIH}}(\zeta,a) & :\Lrarw &
 \forall b\in \mathcal{C}^{\mI_{N}}(\mathcal{W}_{N+1})\cap a\forall f\in J \, A(\zeta,b,f).
\\
\mbox{{\rm SSIH}}(\zeta,a,f) & :\Lrarw &
 \forall g\in J [g<^{0}_{lx}f  \Rarw A(\zeta,a,g)].
\eeqnarrs
}
\end{definition}

\begin{lemma}\label{th:id5wf21}
For {\rm each} $\zeta\in\mathcal{C}^{\mI_{N}}(\mathcal{W}_{N+1})$, the following holds:

Assume 
$a\in\mathcal{C}^{\mI_{N}}(\mathcal{W}_{N+1})\cap(\zeta+1)$,
$f\in J$, 
$\mbox{{\rm SIH}}(\zeta,a)$, $\mbox{{\rm SSIH}}(\zeta,a,f)$ in Definition \ref{df:id4wfA}.
Moreover assume $\mbox{{\rm MIH}}(\zeta;\mathcal{W}_{N+1})$
in Proposition \ref{prp:psimS}.
Then
\[
\psi_{\Omega}(a)\in OT(\mI_{N}) 
 \Rarw 
 \psi_{\Omega}(a)\in\mathcal{W}_{N+1}
 \]
 and
for any $\mS=\mT^{\dagger k}\in SSt$ with $\mT\in\mathcal{W}_{N+1}$ and any
$\kappa\in\mathcal{W}_{N+1}\cap (L(\mS)\cup\{\mS\})$ the following holds:
\[
\mathtt{g}^{*}_{0}(\psi_{\kappa}^{f}(a)) \leq \zeta
 \Rarw 
 \psi_{\kappa}^{f}(a)\in\mathcal{W}_{N+1}.
\]

\end{lemma}
\bprf
Let
$\alpha_{1}=\psi_{\kappa}^{f}(a)\in OT(\mI_{N})$ with 
$a\in\mathcal{C}^{\mI_{N}}(\mathcal{W}_{N+1})\cap(\zeta+1)$,
$\kappa\in\mathcal{W}_{N+1}\cap((L(\mS)\cup\{\mS\})$ and $f\in J$ such that 
$\mS=\mT^{\dagger k}$ with $\mT\in\mathcal{W}_{N+1}$, and
$\mathtt{g}^{*}_{0}(\alp_{1})\leq\zeta$.
By Lemma \ref{lem:wfpsiI} we have $\psi_{\mI_{N}}(\zeta)\in\mathcal{W}_{N+1}$.
By 
Proposition \ref{th:id5wf21_Omega}, Lemma \ref{lem:psimS} and the assumption $\mbox{{\rm MIH}}(\zeta;\mathcal{W}_{N+1})$
it suffices to show 
$\alpha_{1}\in\mathcal{G}^{\mathcal{W}_{N+1}}$.

By Lemma \ref{lem:6.15} we have
$\{\kappa,a\}\cup \mathrm{fld}(f)\subset \mathcal{C}^{\alpha_{1}}(\mathcal{W}_{N+1})$, and hence $\alpha_{1}\in \mathcal{C}^{\alpha_{1}}(\mathcal{W}_{N+1})$.
It suffices to show the following claim by induction on $\ell\beta_1$.
\bclm\label{clm:id5wf21.1}
$\forall\beta_{1}\in \mathcal{C}^{\alpha_{1}}(\mathcal{W}_{N+1})\cap\alpha_{1}[\beta_1\in\mathcal{W}_{N+1}]$.
\eclm
\textbf{Proof} of Claim \ref{clm:id5wf21.1}.
Assume $\beta_{1}\in \mathcal{C}^{\alpha_{1}}(\mathcal{W}_{N+1})\cap\alpha_{1}$ and let
\[
\mbox{LIH} :\Lrarw
\forall\gamma\in \mathcal{C}^{\alpha_{1}}(\mathcal{W}_{N+1})\cap\alpha_{1}[\ell\gamma<\ell\beta_{1} \Rarw \gamma\in\mathcal{W}_{N+1}].
\]

We show $\beta_1\in\mathcal{W}_{N+1}$. 
We can assume $\bet_{1}\not\in\{0,\Ome\}$ by Proposition \ref{lem:6.29}.

\noindent
\textbf{Case 1}. 
$\beta_1\not\in\mathcal{E}(\beta_{1})$:
Assume $\beta_{1}\not\in\mathcal{W}_{N+1}$.
Then 
$\bet_{1}\not\in N(\rho)$
for any $\rho$ by 
$\beta_{1}\in\mathcal{C}^{\alpha_{1}}(\mathcal{W}_{N+1})\cap\alpha_{1}$ and
Definition \ref{df:CX}.
We obtain 
$S(\beta_{1})\subset \mathcal{C}^{\alpha_{1}}(\mathcal{W}_{N+1})\cap\alpha_{1}$. 
LIH yields $S(\beta_{1})\subset\mathcal{W}_{N+1}$. 
Hence we conclude $\beta_{1}\in\mathcal{W}_{N+1}$ from Proposition \ref{prp:noncritical}.
\\

In what follows consider the cases when $\beta_{1}=\psi_{\pi}^{g}(b)$ for some $\pi,b,g$.
We can assume $\pi>\alp_{1}$ and
$\{\pi,b\}\cup \mathrm{fld}(g)\subset \mathcal{C}^{\alpha_{1}}(\mathcal{W}_{N+1})$.
Then either $\pi=\mI_{N}$ or $\bet_{1}\in L(\mS)$ for $\alp_{1}\in L(\mS)$.
\\
\textbf{Case 2}. $\pi=\mI_{N}$ and $b<a$:
As in the proof of Lemma \ref{lem:wfpsiI} we see $b\in\mathcal{C}^{\mI_{N}}(\mathcal{W}_{N+1})$.
We obtain $\bet_{1}=\psi_{\mI_{N}}(b)\in\mathcal{W}_{N+1}$
by $b<a\leq\zeta$ and Lemma \ref{lem:wfpsiI}.
\\
\textbf{Case 3}. $\pi<\mI_{N}$,
$b<a$, $\beta_{1}<\kappa$ and $\{\pi,b\}\cup \mathrm{fld}(g)\subset C_{a}(\alp_{1})$:
Then $\bet_{1}\in L(\mS)$.

Let $B$ denote a set of subterms of $\beta_{1}$ defined recursively as follows.
First $\{\pi,b\}\cup \mathrm{fld}(g)\subset B$.
Let $\alpha_{1}\leq\beta\in B$. 
If $\beta=_{NF}\gamma_{m}+\cdots+\gamma_{0}$, then $\{\gamma_{i}:i\leq m\}\subset B$.
If $\beta=_{NF}\varphi\gamma\delta$, then $\{\gamma,\delta\}\subset B$.
If $\beta=\psi_{\sigma}^{h}(c)$, then $\{\sigma,c\}\cup \mathrm{fld}(h)\subset B$.
If $\bet\in N(\tau)$,
then $\tau\in B$.

Then from $\{\pi,b\}\cup \mathrm{fld}(g)\subset \mathcal{C}^{\alpha_{1}}(\mathcal{W}_{N+1})$ we see inductively that
$B\subset \mathcal{C}^{\alpha_{1}}(\mathcal{W}_{N+1})$.
Hence by LIH we obtain $B\cap\alpha_{1}\subset\mathcal{W}_{N+1}$.
Moreover if $\alpha_{1}\leq\psi_{\sigma}^{h}(c)\in B$, then 
$c\in K_{\alpha_{1}}(\{\pi,b\}\cup \mathrm{fld}(g))<a$.

We claim that
\begin{equation}\label{eq:case2A}
\forall\beta\in B(\beta\in \mathcal{C}^{\mI_{N}}(\mathcal{W}_{N+1}))
\end{equation}
\textbf{Proof} of (\ref{eq:case2A}) by induction on $\ell\beta$.
Let $\beta\in B$. 
We may assume that $\alp_{1}\leq\bet$ is a strongly critical number such that $\bet\not\in\{\Ome,\mI_{N}\}\cup SSt$
by induction hypothesis on the lengths.
First consider the case when $\alpha_{1}\leq\beta=\psi_{\sigma}^{h}(c)$.
By induction hypothesis we have 
 $\{\sigma,c\}\cup \mathrm{fld}(h)\subset\mathcal{C}^{\mI_{N}}(\mathcal{W}_{N+1})$.
On the other hand we have $c<a$ and 
$\mathtt{g}^{*}_{0}(\bet)\leq\mathtt{g}^{*}_{0}(\alp_{1})$ by 
Proposition \ref{prp:hg}.
$\mbox{SIH}(\zeta,a)$ yields $\beta\in\mathcal{W}_{N+1}$.

Second let $\alp_{1}\leq\bet\in N(\tau)$
for a $\tau\in L(\mS)$.
By IH we obtain $\tau\in\mathcal{W}_{N+1}$.
We claim that $\mathtt{g}^{*}_{0}(\tau)\leq\mathtt{g}^{*}_{0}(\alp_{1})\leq\zeta$.
If $\tau\leq\alp_{1}$, then we obtain 
$\mathtt{g}_{0}^{*}(\tau)=\mathtt{g}_{0}^{*}(\alp_{1})$.
Otherwise $\alp_{1}<\tau=\psi_{\sig}^{h}(c)\in B\subset C_{a}(\alp_{1})$ for some $\sig,h,c$.
We obtain $\mathtt{g}^{*}_{0}(\tau)\leq\mathtt{g}^{*}_{0}(\alp_{1})$ by 
Proposition \ref{prp:hg}.
On the other hand we have 
$\mT\in\mathcal{W}_{N+1}$ by one of the assumptions.
Corollary \ref{cor:MG} yields $\bet\in\mathcal{W}_{N+1}$.
\hspace*{\fill} $\Box$ of (\ref{eq:case2A})
\\

In particular we obtain $\{\pi,b\}\cup \mathrm{fld}(g)\subset  \mathcal{C}^{\mI_{N}}(\mathcal{W}_{N+1})$.
Moreover we have $b<a$ and $\mathtt{ g}^{*}_{0}(\bet_{1})\leq\mathtt{ g}^{*}_{0}(\alp_{1})$
by Proposition \ref{prp:LSwinding.1}.
Therefore once again $\mbox{SIH}(\zeta, a)$ yields $\beta_{1}\in\mathcal{W}_{N+1}$.
\\
\textbf{Case 4}. 
$b=a$, $\pi=\kappa$, $\forall\delta\in \mathrm{fld}(g)(K_{\alpha_{1}}(\delta)<a)$ and $g<^{0}_{lx}f$: 
Obviously $\mathtt{g}_{0}^{*}(\bet_{1})=\mathtt{g}_{0}^{*}(\alp_{1})$.
As in (\ref{eq:case2A}) we see that $\mathrm{fld}(g)\subset\mathcal{W}_{N+1}$ from 
$\mbox{SIH}(\zeta,a)$.
$\mbox{SSIH}(\zeta,a,f)$ yields $\beta_{1}\in\mathcal{W}_{N+1}$.
\\
\textbf{Case 5}.
$a\leq b\leq K_{\beta_{1}}(\delta)$ for some 
$\delta\in \mathrm{fld}(f)\cup\{\kappa,a\}$:
It suffices to find a $\gamma$ such that 
$\beta_{1}\leq\gamma\in\mathcal{W}_{N+1}\cap\alpha_{1}$.
Then $\beta_{1}\in\mathcal{W}_{N+1}$ follows from $\beta_{1}\in \mathcal{C}^{\alpha_{1}}(\mathcal{W}_{N+1})$ and Propositions \ref{lem:CX2}.\ref{lem:CX2.3} and \ref{lem:6.16}.

$k_{X}(\alpha)$ denotes the set in Definition \ref{df:EGFk}.
In general we see that $a\in K_{X}(\alpha)$ iff $\psi_{\sigma}^{h}(a)\in k_{X}(\alpha)$ for some $\sigma,h$,
and for each $\psi_{\sigma}^{h}(a)\in k_{X}(\psi_{\sigma_{0}}^{h_{0}}(a_{0}))$ there exists a sequence
$\{\alpha_{i}\}_{i\leq m}$ of subterms of $\alpha_{0}=\psi_{\sigma_{0}}^{h_{0}}(a_{0})$ such that 
$\alpha_{m}=\psi_{\sigma}^{h}(a)$, 
$\alpha_{i}=\psi_{\sigma_{i}}^{h_{i}}(a_{i})$ for some $\sigma_{i},a_{i},h_{i}$,
and for each $i<m$,
$X\not\ni \alpha_{i+1}\in\mathcal{E}(C_{i})$ for $C_{i}=\{\sigma_{i},a_{i}\}\cup \mathrm{fld}(h_{i})$.

Let $\delta\in \mathrm{fld}(f)\cup\{\kappa,a\}$ such that $b\leq \gamma$
for a $\gamma\in K_{\beta_{1}}(\delta)$.
Pick an $\alpha_{2}=\psi_{\sigma_{2}}^{h_{2}}(a_{2})\in \mathcal{E}(\delta)$ 
such that $\gamma\in K_{\beta_{1}}(\alpha_{2})$, and 
an $\alpha_{m}=\psi_{\sigma_{m}}^{h_{m}}(a_{m})\in k_{\beta_{1}}(\alpha_{2})$ for some $\sigma_{m},h_{m}$ and
$a_{m}\geq b\geq a$.
We have $\alpha_{2}\in\mathcal{W}_{N+1}$ by $\delta\in\mathcal{W}_{N+1}$.
If $\alpha_{2}<\alpha_{1}$, then
$\beta_{1}\leq\alpha_{2}\in\mathcal{W}_{N+1}\cap\alpha_{1}$, and we are done.
Assume $\alpha_{2}\geq\alpha_{1}$, i.e., $\alp_{2}\not\in\alp_{1}$.
Then $a_{2}\in K_{\alpha_{1}}(\alpha_{2})<a\leq b$, and $m>2$.

Let $\{\alpha_{i}\}_{2\leq i\leq m}$ be the sequence of subterms of $\alpha_{2}$ such that
$\alpha_{i}=\psi_{\sigma_{i}}^{h_{i}}(a_{i})$ for some $\sigma_{i},a_{i},h_{i}$,
and for each $i<m$,
$\beta_{1}\leq\alpha_{i+1}\in\mathcal{E}(C_{i})$ for $C_{i}=\{\sigma_{i},a_{i}\}\cup \mathrm{fld}(h_{i})$.

Let $\{n_{j}\}_{0\leq j\leq k}\, (0<k\leq m-2)$ be the increasing sequence $n_{0}<n_{1}<\cdots<n_{k}\leq m$ 
defined recursively by $n_{0}=2$, and assuming $n_{j}$ has been defined so that
$n_{j}<m$ and $\alpha_{n_{j}}\geq\alpha_{1}$, $n_{j+1}$ is defined by
$n_{j+1}=\min(\{i: n_{j}< i<m: \alpha_{i}<\alpha_{n_{j}}\}\cup\{m\})$.
If either $n_{j}=m$ or $\alpha_{n_{j}}<\alpha_{1}$, then $k=j$ and $n_{j+1}$ is undefined.
Then we claim that
\begin{equation}\label{eq:case4A}
\forall j\leq k(\alpha_{n_{j}}\in\mathcal{W}_{N+1}) \spand \alpha_{n_{k}}<\alpha_{1}
\end{equation}
\textbf{Proof} of (\ref{eq:case4A}).
By induction on $j\leq k$ we show first that 
$\forall j\leq k(\alpha_{n_{j}}\in\mathcal{W}_{N+1})$. 
We have $\alpha_{n_{0}}=\alpha_{2}\in\mathcal{W}_{N+1}$.
Assume $\alpha_{n_{j}}\in\mathcal{W}_{N+1}$ and $j<k$.
Then $n_{j}<m$, i.e., $\alpha_{n_{j+1}}<\alpha_{n_{j}}$, and 
by $\alpha_{n_{j}}\in \mathcal{C}^{\alpha_{n_{j}}}(\mathcal{W})$, we have $C_{n_{j}}\subset C^{\alpha_{n_{j}}}(\mathcal{W}_{N+1})$,
and hence $\alpha_{n_{j}+1}\in\mathcal{E}(C_{n_{j}})\subset C^{\alpha_{n_{j}}}(\mathcal{W}_{N+1})$.
We see inductively that
$\alpha_{i}\in  \mathcal{C}^{\alpha_{n_{j}}}(\mathcal{W}_{N+1})$ for any $i$ with $n_{j}\leq i\leq n_{j+1}$.
Therefore 
$\alpha_{n_{j+1}}\in  \mathcal{C}^{\alpha_{n_{j}}}(\mathcal{W}_{N+1})\cap\alpha_{n_{j}}\subset\mathcal{W}_{N+1}$ by 
Propositions \ref{lem:CX2}.\ref{lem:CX2.3} and \ref{lem:6.16}.

Next we show that $\alpha_{n_{k}}<\alpha_{1}$.
We can assume that $n_{k}=m$.
This means that $\forall i(n_{k-1}\leq i<m \Rarw \alpha_{i}\geq\alpha_{n_{k-1}})$.
We have
$\alpha_{2}=\alpha_{n_{0}}>\alpha_{n_{1}}>\cdots>\alpha_{n_{k-1}}\geq\alpha_{1}$, and
$\forall i<m(\alpha_{i}\geq\alpha_{1})$.
Therefore $\alpha_{m}\in k_{\alpha_{1}}(\alpha_{2})\subset k_{\alpha_{1}}(\{\kappa,a\}\cup \mathrm{fld}(h))$, i.e.,
$a_{m}\in K_{\alpha_{1}}(\{\kappa,a\}\cup \mathrm{fld}(h))$ for $\alpha_{m}=\psi_{\sigma_{m}}^{h_{m}}(a_{m})$.
On the other hand we have $K_{\alpha_{1}}(\{\kappa,a\}\cup \mathrm{fld}(h))<a$ for 
$\alpha_{1}=\psi_{\sigma}^{h}(a)$.
Thus $a\leq a_{m}<a$, a contradiction.

(\ref{eq:case4A}) is shown, and we obtain $\beta_{1}\leq\alpha_{n_{k}}\in\mathcal{W}_{N+1}\cap\alpha_{1}$.

This completes a proof of Claim \ref{clm:id5wf21.1} and of the lemma.
\hspace*{\fill} $\Box$

\bcor\label{cor:psiwf}
For {\rm each} $\zeta\in\mathcal{C}^{\mI_{N}}(\mathcal{W}_{N+1})$, 
$\mbox{{\rm MIH}}(\zeta;\mathcal{W}_{N+1})$ holds.
\ecor
\bprf
For each $n<\ome$, we have $\mathrm{TI}[\mathcal{C}^{\mI_{N}}(\mathcal{W}_{N+1})\cap\ome_{n}(\mI_{N}+1)]$ by Lemma \ref{lem:6.34}.\ref{lem:6.34.3}.
We show $\mbox{{\rm MIH}}(\zeta;\mathcal{W}_{N+1})$ by induction on 
$\zeta\in\mathcal{C}^{\mI_{N}}(\mathcal{W}_{N+1})$.
Assume $\fal\xi\in\mathcal{C}^{\mI_{N}}(\mathcal{W}_{N+1})\cap\zeta\, \mbox{{\rm MIH}}(\xi;\mathcal{W}_{N+1})$.

Let $\mS=\mT^{\dagger k}$ with $\mT\in\mathcal{W}_{N+1}$, and
$\gam=\psi_{\kap}^{f}(a)\in L(\mS)$ be such that
$k(\gam)=\{\kap,a\}\cup \mathrm{fld}(f)\subset\mathcal{C}^{\mI_{N}}(\mathcal{W}_{N+1})$ and 
$\xi=\mathtt{h}(\gam)=\{a,\mathtt{g}_{0}^{*}(\gam)\}\subset\mathcal{C}^{\mI_{N}}(\mathcal{W}_{N+1})\cap\zeta$.
We obtain $\mbox{{\rm MIH}}(\xi;\mathcal{W}_{N+1})$ by IH.

We obtain
$
\gam\in\mathcal{W}_{N+1}
$ by Lemma \ref{th:id5wf21} and $\mbox{{\rm MIH}}(\xi;\mathcal{W}_{N+1})$ 
with subsidiary induction on $a\in\mathcal{C}^{\mI_{N}}(\mathcal{W}_{N+1})\cap(\xi+1)$ and 
sub-subsidiary induction 
on $f\in J$.
Then Corollary \ref{cor:MG} yields 
$N(\gam)\subset\mathcal{W}_{N+1}$.

Here by induction on $f\in J$ we mean by induction along $g<_{lx}^{0}f$.
In the proof of Lemma \ref{th:id5wf21}, $\mbox{SSIH}(\zeta,a,f)$ is invoked 
in \textbf{Case 4}, i.e., only when
$\psi^{g}_{\kap}(a)<\psi_{\kap}^{f}(a)$ with $\kap<\mI_{N}$.
Then Lemma \ref{lem:oflx} yields
$o_{\mI_{N}}(g)<o_{\mI_{N}}(f)\in\mathcal{C}^{\mI_{N}}(\mathcal{W}_{N+1})$ for 
$\mathrm{fld}(f)\subset\mathcal{C}^{\mI_{N}}(\mathcal{W}_{N+1})\cap\Lam$,
where $\Lam=\psi_{\mI_{N}}(b)$ and
$b=\mathtt{g}^{*}_{0}(\psi_{\kap}^{f}(a))\geq\mathtt{ p}_{0}(\psi_{\kap}^{f}(a))$.
Hence
$o_{\Lam}(g)<o_{\Lam}(f)\in\mathcal{W}_{N+1}$ by $\Lam\in\mathcal{W}_{N+1}$.
\eprf

\begin{lemma}\label{lem:psiw}
For {\rm each} $n<\ome$, the following holds:

If one of the followings holds, then 
$\alp\in\mathcal{W}_{N+1}$ for $\alp\in OT(\mI_{N})$.

\benu
\item\label{lem:psiw.2}
$\alp=\mS^{\dagger k}$ with $\mS\in\mathcal{W}_{N+1}\cap (St_{k}\cup\{\Ome\})$.

\item\label{lem:psiw.5}
$\alp=\psi_{\mI_{N}}(a)$ with
$a\in \mathcal{C}^{\mI_{N}}(\mathcal{W}_{N+1})\cap\ome_{n}(\mI_{N}+1)$.

\item\label{lem:psiw.1}
$\alp=\psi_{\Omega}(a)$ with
$a\in \mathcal{C}^{\mI_{1}}(\mathcal{W}_{N+1})\cap\ome_{n}(\mI_{N}+1)$.

\item\label{lem:psiw.6}
$\alp=\psi_{\kappa}^{f}(a)\in L(\mS)$ 
for $\mS=\mT^{\dagger k}$ with $\mT\in\mathcal{W}_{N+1}$ 
and 
$k(\alp)\cup\mathtt{h}(\alp)=
\{
\kappa,a,\mathtt{ g}^{*}_{0}(\alp)\}\cup \mathrm{fld}(f)\subset\mathcal{C}^{\mI_{N}}(\mathcal{W}_{N+1})\cap\ome_{n}(\mI_{N}+1)$.

\item\label{lem:psiw.3}
$\alp\in N(\rho)$
for $\rho\in\mathcal{W}_{N+1}\cap L(\mS)$ with
$\mS=\mT^{\dagger k}$ such that $\mT\in\mathcal{W}_{N+1}$
 and $\mathtt{ g}^{*}_{0}(\rho)<\ome_{n}(\mI_{N}+1)$.

\eenu

\end{lemma}
\bprf
\ref{lem:psiw}.\ref{lem:psiw.2} is seen from Lemma \ref{cor:6.21}.
\\
\ref{lem:psiw}.\ref{lem:psiw.5} follows from Lemma \ref{lem:wfpsiI}.
\\
\ref{lem:psiw}.\ref{lem:psiw.1} and \ref{lem:psiw}.\ref{lem:psiw.6} are seen from Lemma \ref{th:id5wf21} and Corollary \ref{cor:psiwf}.
\\
\ref{lem:psiw}.\ref{lem:psiw.3} follows from Corollaries 
\ref{cor:MG} and \ref{cor:psiwf}.
\eprf
\\

\noindent
Let us conclude Theorem \ref{th:wf}.
For each $\alp\in OT(\mI_{N})$,
$\alp\in\mathcal{C}^{\mI_{N}}(\mathcal{W}_{N+1})$
is seen
by metainduction on the lengths $\ell\alp$ using
Propositions \ref{prp:noncritical},
\ref{lem:6.29}
and Lemma \ref{lem:psiw}.
Note that $\ell(\mathtt{ g}^{*}_{0}(\psi_{\kap}^{f}(a)))<\ell(\psi_{\kap}^{f}(a))$
and $\ell(\mT)<\ell(\rho)$ for $\rho\in L(\mS)$ and $\mS=\mT^{\dagger k}$.
Therefore we obtain
$\Sig^{1}_{N+2}\mbox{{\rm -DC+BI}}\vdash \alp\in\mathcal{C}^{\mI_{N}}(\mathcal{W}_{N+1})\cap\Ome=\mathcal{W}_{N+1}\cap\Ome=W(\mathcal{C}^{0}(\mathcal{W}_{N+1}))\cap\Ome=W(OT(\mI_{N}))\cap\Ome$, and 
$\Sig^{1}_{N+2}\mbox{{\rm -DC+BI}}\vdash Wo[\alp]$
for each $\alp<\psi_{\Ome}(\veps_{\mI_{N}+1})$.

\section{Outcomes on $\mathbf{Z}_{2}$
}\label{sect:consis}
In this final section let us conclude some standard outcomes of
an ordinal analysis of the theory $\mathbf{Z}_{2}$.

Let $\mathrm{TI}[\Pi^{1-}_{0},\psi_{\Ome}(\veps_{\mI_{N}+1})]$ denote a schema of transfinite induction 
$\fal\alp\in OT(\mI_{N})\cap\Ome\,( \mathrm{Prg}[OT(\mI_{N}),A]\to OT(\mI_{N})\cap\alp\subset A)$ up to $\psi_{\Ome}(\veps_{\mI_{N}+1})$ in $OT(\mI_{N})$
applied to arithmetic formulas $A\in\Pi^{1-}_{0}$ in the language of the first-order arithmetic ${\sf PA}$.
Let $T_{0}={\sf PA}+\bigcup\{\mathrm{TI}[\Pi^{1-}_{0},\psi_{\Ome}(\veps_{\mI_{N}+1})]: N<\ome\}$, and 
$T_{1}=\mathrm{FiX}^{i}(T_{0})$ denote the intuitionistic fixed point theory over $T_{0}$.
The language of the theory $T_{1}$ is expanded by unary predicate symbols $I$ for each operator $\Phi(X,x)$, in which
every occurrence of a unary predicate symbol $X$ is strictly positive.
The axioms in $T_{1}$ are obtained from $T_{0}$ by adding the axioms $\fal x[I(x)\lrarw \Phi(I,x)]$ for a fixed point $I$.
The axiom schema $\mathrm{TI}[\Pi^{1-}_{0},\psi_{\Ome}(\veps_{\mI_{N}+1})]$ of transfinite induction as well as
schema of complete induction may be applied to arbitrary first-order formulas in the expanded language with the predicates $I$.
The underlying logic in $T_{1}$ is the intuitionistic first-order logic with the axiom $\fal x,y(x=y\to I(x)\to I(y))$.
The excluded middle $\fal x(\lnot I(x)\lor I(x))$ for the predicate $I$ is not available in $T_{1}$.

\blem\label{lem:FiX}
$\mathrm{FiX}^{i}(T_{0})$ is a conservative extension of $T_{0}$. 
Moreover the fact is provable in the fragment $I\Sig^{0}_{1}$ of the first-order arithmetic:
$I\Sig^{0}_{1}\vdash \mathrm{Pr}_{T_{1}}(\lc\vphi\rc)\to\mathrm{Pr}_{T_{0}}(\lc\vphi\rc)$, 
where $\mathrm{Pr}_{T}(x)$ is a standard provability predicate for a theory $T$.
\elem
\bprf
The fact is seen as in \cite{IntFix, OA}.
To formalize a proof of the fact in $I\Sig^{0}_{1}$, follow a finitary analysis in section 4.4 of \cite{OA}.
\eprf

\begin{theorem}\label{th:final}
$\mathbf{Z}_{2}$ is a conservative extension of ${\sf PA}+\bigcup\{{\rm TI}[\Pi^{1-}_{0},\psi_{\Ome}(\veps_{\mI_{N}+1})]: N<\ome\}$.
Moreover the fact is provable in the fragment $I\Sig^{0}_{1}$.
\end{theorem}
\bprf
Assume that $\mathbf{Z}_{2}\vdash A$ for an arithmetic sentence $A\in\Pi^{1-}_{0}$.
Pick an $N<\ome$ such that $\Sig^{1}_{N+2}\mbox{{\rm -DC+BI}}\vdash A$.
By Lemma \ref{lem:9scdset} we obtain ${\sf KP}\ome+\Pi_{N}\mbox{{\rm -Collection}}+(V=L)\vdash A^{set}$, and hence
${\sf KP}\ome+\Pi_{N}\mbox{{\rm -Collection}}\vdash A^{set}$.
Then by Lemma \ref{lem:9setI} we obtain $S_{\mI_{N}}\vdash A^{set}$.

Now we see that the proof of Theorem \ref{thm:2} in sections \ref{sec:operatorcont} and \ref{sec:proofonestep} is formalizable in
the intuitionistic fixed point theory $T_{1}=\mathrm{FiX}^{i}(T_{0})$ over $T_{0}$.
Let us regard each of the relations
$( \mathcal{H}_{\gamma},\Theta;\mathtt{D})\vdash^{* a}_{c, \gam_{0}} \Gamma$ and 
$( \mathcal{H}_{\gamma},\Tht, \mathtt{ Q})\vdash^{a}_{c,\xi,\Lambda_{0},\gam_{0}} \Gamma$ as a fixed point of a strictly positive operator.
Then by applying transfinite induction to first-order formulas with the fixed point predicates, Theorem \ref{thm:2} is proved.
Therefore we obtain $\mathrm{FiX}^{i}(T_{0})\vdash A$, and $T_{0}\vdash A$ by Lemma \ref{lem:FiX}.
\eprf
\\

We see readily that 
the transfinite induction  ${\rm TI}(\psi_{\Ome}(\mI_{\ome}))$ up to $\psi_{\Ome}(\mI_{\ome})$
is equivalent to the $\Pi^{1}_{1}$-soundness ${\rm RFN}_{\Pi^{1}_{1}}(\mathbf{Z}_{2})$ of $\mathbf{Z}_{2}$
over $\mathrm{RCA}_{0}$, where
$\mathrm{TI}(\psi_{\Ome}(\mI_{\ome}))$ denotes a $\Pi^{1}_{1}$-sentence
$\fal N\fal\alp\in OT(\mI_{N})\cap\Ome\fal Y\,( \mathrm{Prg}[OT(\mI_{N}),Y]\to OT(\mI_{N})\cap\alp\subset Y)$.

\bdf\label{df:DS}
{\rm
Let $\alp\in OT(\mI_{N})$ be an ordinal  term.
\benu
\item
$DS_{\alp}$ denotes a $\Pi^{0}_{2}$-sentence saying that `there is no primitive recursive and descending sequence $\{f(n)\}_{n}$ of ordinals
with $f(0)<\alp$'.
This means that $f(0)<\alp\Rarw\exi n(f(n+1)\not<f(n))$.

\item
$WDS_{\alp}$ denotes a $\Pi^{0}_{3}$-sentence saying that `for every primitive recursive and weakly descending sequence $\{f(n)\}_{n}$ of ordinals
with $f(0)<\alp$, there exists an $n$ such that $\fal m\geq n(f(m)=f(n))$'.
This is equivalent to the principle that `for every primitive recursive sequence $\{f(n)\}_{n}$ of ordinals, there exists 
an $n$ such that $\fal m(f(n)\leq f(m))$.

\item
$DS_{\psi_{\Ome}(\veps_{\mI_{N}+1})}:\Lrarw \fal\alp\in OT(\mI_{N})\cap\Ome\, DS_{\alp}$ and
$DS_{\psi_{\Ome}(\mI_{\ome})}:\Lrarw \fal N>0\, DS_{\psi_{\Ome}(\veps_{\mI_{N}+1})}$.

Also
 $WDS_{\psi_{\Ome}(\veps_{\mI_{N}+1})}:\Lrarw \fal\alp\in OT(\mI_{N})\cap\Ome\, WDS_{\alp}$ and
 $WDS_{\psi_{\Ome}(\mI_{\ome})}:\Lrarw \fal N>0\, WDS_{\psi_{\Ome}(\veps_{\mI_{N}+1})}$.

\item A computable (total) function $f$ on integers is said to be \textit{$\psi_{\Ome}(\veps_{\mI_{N}+1})$-recursive}
 if $f$ is defined from $\psi_{\Ome}(\veps_{\mI_{N}+1})$-recursive functions $g,r,h$ by $\psi_{\Ome}(\veps_{\mI_{N}+1})$-recursion:
\[
f(y,x)=
\left\{
\begin{array}{ll}
g(y,x,f(y,r(y,x))) & \mbox{{\rm if } } r(y,x)<x<\Ome \mbox{ {\rm in }} OT(\mI_{N})
\\
h(y,x) & \mbox{otherwise}
\end{array}
\right.
\]

\item
$\mathrm{RFN}_{\Sig^{0}_{n}}(\mathbf{Z}_{2})$ denotes the uniform reflection principle of $\mathbf{Z}_{2}$ for $\Sig^{0}_{n}$-formulas.
\eenu

}
\edf

\bcor\label{cor:4main}

\benu

\item\label{cor:4main.22}
The 2-consistency $\mathrm{RFN}_{\Sig^{0}_{2}}(\mathbf{Z}_{2})$ of $\mathbf{Z}_{2}$ is equivalent to 
$WDS_{\psi_{\Ome}(\mI_{\ome})}$ over $I\Sig^{0}_{1}$.

\item\label{cor:4main.23}
$\mathbf{Z}_{2}$ is $\Pi^{0}_{3}$-conservative over $I\Sig^{0}_{1}+\{WDS_{\psi_{\Ome}(\veps_{\mI_{N}+1})}:0<N<\ome\}$.

\item\label{cor:4main.25}
The 1-consistency $\mathrm{ RFN}_{\Sig^{0}_{1}}(\mathbf{Z}_{2})$ of $\mathbf{Z}_{2}$ is equivalent to 
$DS_{\psi_{\Ome}(\mI_{\ome})}$ over $I\Sig^{0}_{1}$.

\item\label{cor:4main.26}
$\mathbf{Z}_{2}$ is $\Pi^{0}_{2}$-conservative over $I\Sig^{0}_{1}+\{DS_{\psi_{\Ome}(\veps_{\mI_{N}+1})}:0<N<\ome\}$.

\item\label{cor:4main.27}
For computable total function $f$ on $\Natural$,
$f$ is provably computable in $\mathbf{Z}_{2}$
iff $f$ is $\psi_{\Ome}(\veps_{\mI_{N}+1})$-recursive for an $N<\ome$.

\eenu

\ecor
\bprf
Each follows from Theorem \ref{th:final} as in chapter 4 of \cite{OA}.
\eprf
\\

For the consistency $\mathrm{Con}(\mathbf{Z}_{2})$ of $\mathbf{Z}_{2}$ we obtain the following.

\bcor\label{cor:4main30}
There are primitive recursive predicate $B$ and primitive recursive function  $f$ such that
both of $\fal N>0\fal\alp\in OT(\mI_{N})\cap\Ome(f(N,\alp)<\alp\to B(N,f(N,\alp)) \to B(N,\alp))$
and $\fal N>0\fal\alp\in OT(\mI_{N})\cap\Ome\, B(N,\alp) \to \mathrm{Con}(\mathbf{Z}_{2})$
is provable in $I\Sig^{0}_{1}$.
\ecor
\bprf
This is seen from Theorem \ref{th:final} as in section 4.3 of \cite{OA}.
\eprf

\end{document}

%% file: header.tex
\newcommand{\brem}{\begin{remark}}
\newcommand{\erem}{\end{remark}}
\newcommand{\blem}{\begin{lemma}}
\newcommand{\elem}{\end{lemma}}
\newcommand{\bth}{\begin{theorem}}
\newcommand{\ethm}{\end{theorem}}
\newcommand{\benu}{\begin{enumerate}}
\newcommand{\eenu}{\end{enumerate}}
\newcommand{\bdes}{\begin{description}}
\newcommand{\edes}{\end{description}}
\newcommand{\bdf}{\begin{definition}}
\newcommand{\edf}{\end{definition}}
\newcommand{\bcor}{\begin{cor}}
\newcommand{\ecor}{\end{cor}}
\newcommand{\bprp}{\begin{proposition}}
\newcommand{\eprp}{\end{proposition}}
\newcommand{\bmlem}{\begin{mlemma}}
\newcommand{\emlem}{\end{mlemma}}
\newcommand{\bclm}{\begin{claim}}
\newcommand{\eclm}{\end{claim}}
\newcommand{\bprf}{{\bf Proof}.\hspace{2mm}}

\newcommand{\eprf}{\hspace*{\fill} $\Box$}

\newcommand{\beqn}{\begin{equation}}
\newcommand{\eeqn}{\end{equation}}
\newcommand{\beqnarr}{\begin{eqnarray}}
\newcommand{\eeqnarr}{\end{eqnarray}}
\newcommand{\beqnarrs}{\begin{eqnarray*}}
\newcommand{\eeqnarrs}{\end{eqnarray*}}

\newcommand{\spand}{\,\&\,}

\newcommand{\Natural}{\mathbb{N}}

\newcommand{\restrict}{\!\upharpoonright\!}

\newtheorem{theorem}{Theorem}[section]
\newtheorem{definition}[theorem]{Definition}
\newtheorem{proposition}[theorem]{Proposition}
\newtheorem{lemma}[theorem]{Lemma}
\newtheorem{cor}[theorem]{Corollary}
\newtheorem{remark}[theorem]{Remark}
\newtheorem{mlemma}[theorem]{Main Lemma}
\newtheorem{claim}[theorem]{Claim}

\newcommand{\alp}{\alpha}

\newcommand{\veps}{\varepsilon}
\newcommand{\del}{\delta}
\newcommand{\Del}{\Delta}
\newcommand{\ome}{\omega}
\newcommand{\Ome}{\Omega}
\newcommand{\bet}{\beta}
\newcommand{\gam}{\gamma}
\newcommand{\Gam}{\Gamma}
\newcommand{\kap}{\kappa}
\newcommand{\sig}{\sigma}
\newcommand{\Sig}{\Sigma}
\newcommand{\tht}{\theta}
\newcommand{\Tht}{\Theta}
\newcommand{\lam}{\lambda}
\newcommand{\Lam}{\Lambda}
\newcommand{\vphi}{\varphi}

\newcommand{\fal}{\forall}
\newcommand{\exi}{\exists}
\newcommand{\rarw }{\rightarrow}
\newcommand{\Rarw }{\Rightarrow}

\newcommand{\lrarw}{\leftrightarrow}
\newcommand{\Lrarw}{\Leftrightarrow}

\newcommand{\calc}{{\cal C}}

\newcommand{\calg}{{\cal G}}
\newcommand{\calh}{{\cal H}}

\newcommand{\calw}{{\cal W}}

\newcommand{\la}{\langle}
\newcommand{\ra}{\rangle}
\newcommand{\lc}{\lceil}
\newcommand{\rc}{\rceil}
\newcommand{\rk}{\mbox{{\rm rk}}}

\newcommand{\sfk}{{\sf k}}

\newcommand{\setm}{\setminus}

%% file: piN_Collection_arxiv_v3.bbl
\begin{thebibliography}{99}




\bibitem{IntFix}T. Arai,
Quick cut-elimination for strictly positive cuts, Ann. Pure Appl. Logic 162 (2011), 807-815.

\bibitem{KPPiN}
T. Arai, A simplified ordinal analysis of first-order reflection,
Jour. Symb. Logic 85 (2020) 1163-1185.

\bibitem{OA}
T. Arai,
Ordinal Analysis with an Introduction to Proof Theory,
(Springer, Singapore, 2020)

\bibitem{singlewfprf}
T. Arai, Wellfoundedness proof with the maximal distinguished set,
Arch. Math. Logic 62 (2023) 333-357.

\bibitem{singlestable}
T. Arai,
An ordinal analysis of a single stable ordinal,
submitted.


\bibitem{Ba}  J. Barwise, 
Admissible Sets and Structures (Springer, Berlin, 1975)

\bibitem{Buchholz75}
W. Buchholz, 
Normalfunktionen und konstruktive Systeme von Ordinalzahlen.
in: J. Diller, G. H. M\"uller, eds.
Proof Theory Symposion Kiel 1974, Lect. Notes Math. vol. 500, pp. 4-25 (Springer, Berlin, 1975)

\bibitem{Buchholz86}
W. Buchholz, 
A new system of proof-theoretic ordinal functions,
Ann. Pure Appl. Logic 32 (1986), 195-207.

\bibitem{Buchholz} W. Buchholz, 
A simplified version of local predicativity, 
in: P. H. G. Aczel, H. Simmons and S. S. Wainer, eds. Proof Theory. pp. 115-147
(Cambridge UP, Cambridge, 1992)

\bibitem{BuchholzBSL}
W. Buchholz, Review of the paper: A. Setzer, Well-ordering proofs for Martin-L\"of type theory, 
Bull. Symb. Logic 6 (2000) 478-479.


\bibitem{J2}G. J\"ager, 
A well-ordering proof for Feferman's theory $T$, 
Archiv f. math. Logik u. Grundl. 23 (1983) 65-77.

\bibitem{J3} G. J\"ager, 
Theories for admissible sets, A unifying approach to proof theory, Studies in Proof Theory Lecture Notes 2 (Bibliopolis, Napoli, 1986)

\bibitem{Rathjen94} 
M. Rathjen, 
Proof theory of reflection, 
Ann. Pure Appl. Logic 68 (1994) 181-224.

%\bibitem{RathjenAFML1}M. Rathjen, 
%An ordinal analysis of stability, 
%\textit{Arch. Math. Logic} {\bf 44} (2005) 1-62.

\bibitem{RathjenAFML2}
M. Rathjen, 
An ordinal analysis of parameter free $\Pi^{1}_{2}$-comprehension,
Arch. Math. Logic 44 (2005) 263-362.


\bibitem{Simpson}S. Simpson, 
Subsystems of Second Order Arithmetic, second edition (Cambridge UP, Cambridge, 2009)

\end{thebibliography}
